\documentclass{book} 

\usepackage{etex} 
\reserveinserts{28} 

\usepackage{hyperref}
\usepackage{framed,graphicx}
\usepackage{enumerate}
\usepackage{enumitem}
\usepackage{mathtools, skull}
\usepackage{faktor,stmaryrd}
\makeatletter
\DeclareRobustCommand*{\mfaktor}[3][]
{
 { \mathpalette{\mfaktor@impl@}{{#1}{#2}{#3}} }
}
\newcommand*{\mfaktor@impl@}[2]{\mfaktor@impl#1#2}
\newcommand*{\mfaktor@impl}[4]{
 \settoheight{\faktor@zaehlerhoehe}{\ensuremath{#1#2{#3}}}%
 \settoheight{\faktor@nennerhoehe}{\ensuremath{#1#2{#4}}}%
 \raisebox{-0.5\faktor@zaehlerhoehe}{\ensuremath{#1#2{#3}}}%
 \mkern-4mu\diagdown\mkern-5mu%
 \raisebox{0.5\faktor@nennerhoehe}{\ensuremath{#1#2{#4}}}%
}
\makeatother

\usepackage{tikz}
\usetikzlibrary{cd}
\usetikzlibrary{decorations.markings}
			\tikzset{->-/.style={decoration={
 					 markings,
 					 mark=at position #1 with {\arrow{>}}},postaction={decorate}}}
			\tikzset{-<-/.style={decoration={
 					markings,
 					mark=at position #1 with {\arrow{<}}},postaction={decorate}}}

\input{epsf.sty}
\usepackage{pgfplots}

\usepackage{makeidx}
\makeindex
 
\usepackage{fancyhdr}

\newcommand{\arccosh}{{\rm arccosh}}
\newcommand{\QQ}[1]{\langle #1 \rangle}
\renewcommand{\Re}{{\mathfrak{Re}}} 
\newcommand{\imaginary}{{\mathfrak{Im}}} 

\newcommand{\Hyp}[1]{\mathbb{#1}\mathbf{H}}
\newcommand{\Bhyp}[1]{\partial_\infty\mathbb{#1}\mathbf{H}}

\newcommand*\conj[1]{\overline{#1}}

\newlength{\originalbase}
\newcommand{\spacing}[1]{\setlength{\baselineskip}{#1\originalbase}}

\usepackage{amscd,amssymb,amsthm,verbatim,
amsmath,amsfonts,mathrsfs,xspace}
\usepackage{subfigure}

\usepackage{color}
\usepackage[all]{xy}

\newcommand{\acts}{\curvearrowright}
\newcommand{\Ad}{\operatorname{Ad}}
\newcommand{\ad}{\operatorname{ad}}

\newcommand{\Area}{\operatorname{Area}}
\newcommand{\Aut}{\operatorname{Aut}}

\newcommand{\bbraket}[1]{\langle\!\langle #1 \rangle\!\rangle}

\newcommand{\Bx}{\mathtt{Box}} 
\newcommand{\C}{{\mathbb C}}

\renewcommand{\d}{\operatorname{d}}

\newcommand{\dd}{\,\mathrm{d}}	

\newcommand{\dcc}{d_{\rm cc}} 

\newcommand{\diam}{\operatorname{diam}}

\newcommand{\End}{\operatorname{End}}

\newcommand{\g}{\mathfrak{g}}

\newcommand{\GL}{\operatorname{GL}}
\newcommand{\gl}{\mathfrak{gl}}	

\newcommand{\h}{\mathfrak{h}}

\newcommand{\id}{\operatorname{id}}
\newcommand{\Id}{\operatorname{Id}}
\newcommand{\IFF}{\Leftrightarrow}	

\newcommand{\into}{\hookrightarrow}		

\newcommand{\Isom}{\operatorname{Isom}}
\renewcommand{\k}{\mathfrak{k}}
\newcommand{\Ker}{\operatorname{Ker}}

\newcommand{\Lip}{\mathtt{Lip}}

\renewcommand{\H}{{\mathcal H}}
\renewcommand{\S}{{\mathscr S}}
\newcommand{\Jac}{\mathtt{Jac}} 
\newcommand{\Length}{\operatorname{Length}}

\newcommand{\LL}{{\mathcal L}}

\newcommand{\Mat}{\rm{Mat}}	

\newcommand{\Mult}{\rm{Mult}}
\newcommand{\N}{{\mathbb N}}

\newcommand{\norm}[1]{\left\Vert#1\right\Vert}
\newcommand{\norma}[1]{\|#1\|}

\newcommand{\R}{{\mathbb R}}
\newcommand{\G}{{\mathbb G}}

\newcommand{\scr}[1]{\mathscr{#1}}

\newcommand{\Span}{\operatorname{span}}

\newcommand{\THEN}{\Rightarrow}	

\renewcommand{\Vec}{\mathrm{Vec}}	

\newcommand{\vol}{\operatorname{vol}}
\newcommand{\Z}{{\mathbb Z}}

\newcommand{\eps}{\epsilon}

\newcommand{\ra}{\rightarrow}

\newcommand{\Lie}{\operatorname{Lie}}

\newcommand{\abs}[1]{\left\lvert#1\right\rvert}

\def\XXint#1#2#3{{\setbox0=\hbox{$#1{#2#3}{\int}$}
\vcenter{\hbox{$#2#3$}}\kern-.5\wd0}}

\numberwithin{equation}{section}

\theoremstyle{definition}
\newtheorem{example}[equation]{Example}

\newtheorem{remark}[equation]{Remark}
\newtheorem{Rem}[equation]{Remark} 
\newtheorem{definition}[equation]{Definition}
\newtheorem{exe}[equation]{\textsc{Exercise}}
\newenvironment{exercise}{\begin{exe}}{\hfill $\blacksquare$\end{exe}}
\newtheorem{formula}[equation]{Formula}

\theoremstyle{plain}

\newtheorem{lemma}[equation]{Lemma}
\newtheorem{theorem}[equation]{Theorem}
\newtheorem{proposition}[equation]{Proposition}

\newtheorem{corollary}[equation]{Corollary}
\newtheorem{conjecture}[equation]{Conjecture}

\newtheorem{question}[equation]{Question}
\newtheorem{problem}[equation]{Problem}

\newtheorem{fact}[equation]{Fact}

\newtheorem{customtheorem}{Theorem}

\makeindex

\begin{document}						

\spacing{1.62} 

\begin{titlepage}

	\begin{centering} 
	
	{-- book project --}
	
		\begin{framed}
\begin{framed}
	\vspace*{.1in}
 	{\Huge \textsf{Metric	Lie Groups	} 
}\\
	\vspace*{.2in}
		{\large \textsf{Carnot-Carath\'eodory spaces\\
		from the homogeneous viewpoint	
		}}
		\vspace*{.4in}
			\end{framed}
		\end{framed}
			\vspace*{1in}
			\includegraphics[height=2.8in]{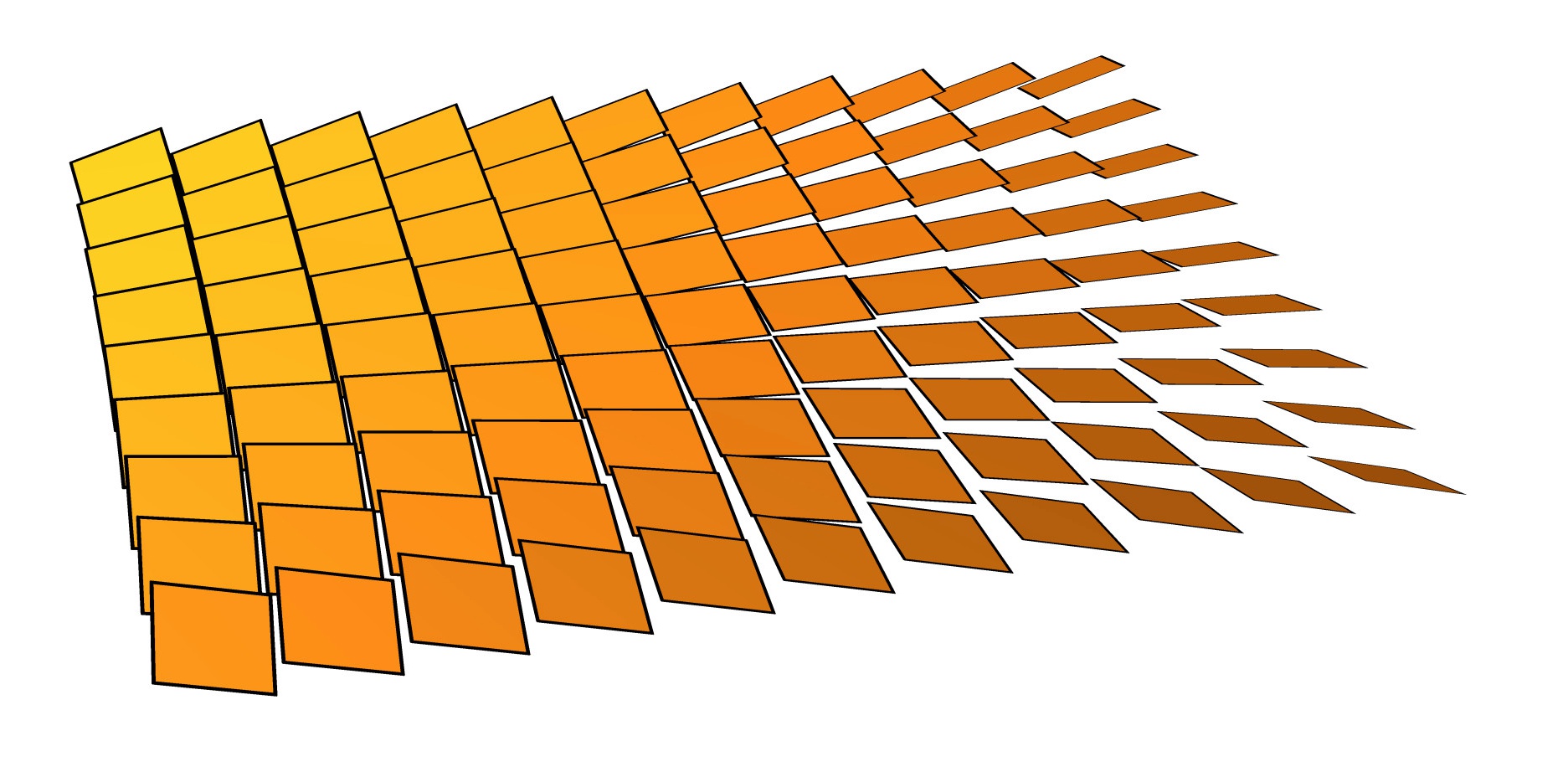}

	\vspace*{1cm}
		by Enrico Le Donne.\\
https://sites.google.com/view/enricoledonne/	\\
	Version of October 9, 2024
		
	\end{centering}
\end{titlepage}
\newpage
\pagestyle{empty}
\frontmatter
\pagestyle{plain}
\pagestyle{fancy}
\setcounter{tocdepth}{2}
\tableofcontents 

%

\mainmatter
\setcounter{chapter}{-1}



\section*{Preface}
The present book treats sub-Riemannian geometries and their generalizations, which go under the names of sub-Finsler geometries or Carnot-Carathéodory spaces. We will discuss these non-smooth geometries, focusing on the cases where there is the additional presence of a group structure. The techniques from Lie group theory will then be mixed with metric geometry, giving a new viewpoint.

This preface describes the origins of this book, outlines its purpose, and identifies its intended audience. 
 This text draws heavily on the following references: 
 \cite{
Pontryagin_topgrp, 
Warner, 
 Corwin-Greenleaf, bellaiche, Gromov, Ambook, Burago:book,
Helgason, 
 Montgomery, 
 Knapp,
 Hilgert_Neeb:book}, 
 as well as on various articles due to the author, his collaborators, and students, such as \cite{AKL,
 Breuillard-LeDonne1,
 LeDonne_characterization, ledonne_primer,LMPOV, LeDonne_NicolussiGolo_metric,CKLNO,Antonelli_LeDonne_MR4645068,
 Pilastro_master_thesis,
 Purho_master,
 Sopio_master}.
 Additionally, it incorporates insights gained from numerous conversations that the author had with his collaborators and mentors.
He would like to acknowledge, in particular, Bruce Kleiner, Urs Lang, Emmanuel Breuillard, Alessandro Ottazzi, Pierre Pansu, and Yves de Cornulier, who have provided invaluable guidance and support throughout this work.

 This text has its origin in lecture notes for a course titled `Sub-Riemannian Geometry' taught at ETH Z\"urich during the fall of 2009 and later at the University of Jyv\"askyl\"a in the spring of 2014. Additional sections were included after the author delivered courses on `Carnot Groups' at a summer school in Levico Terme (Trento, Italy) in 2015 and `Riemannian and Sub-Riemannian Geometry on Lie Groups' at the Neurogeometry summer school in Cortona (Italy) in 2017. The notes were further expanded for the course `Sub-Riemannian Geometry' taught at the University of Fribourg (Switzerland) in spring 2021.
 

This book's primary audience consists of young researchers seeking an introduction to sub-Riemannian geometry. It can serve as background reading material for a master's thesis or as an initial reference for those beginning a PhD program focusing on subjects at the crossroads of geometry, analysis, and group theory.

In contrast to other sources, such as \cite{Montgomery, Agrachev_Barilari_Boscain:book, jeancontrol, Rifford:book, Bonfiglioli:et:al}, this book employs the language and formalism of Lie groups and treats the general category of Carnot-Carath\'eodory spaces, considering norms that not necessarily come from scalar products. In fact, one of the aims of this book is to demonstrate how sub-Finsler geometries manifest in other mathematical domains, including hyperbolic geometry and geometric group theory, through the perspective of Lie groups.

Prerequisite topics from differential geometry, measure theory, and group theory will be discussed within the chapters' main flow. Given the positive feedback received from many students regarding this approach, the author has finally decided to publish this text.

\chapter{Introduction}

%
%

This book is an exploration of 
Carnot-Carathéodory spaces, through the perspective of Lie groups. 
It is intended to study these non-smooth geometries, focusing on the prototypical examples called Carnot groups. Carnot groups
are a fundamental class of nilpotent Lie groups equipped with sub-Riemannian or, more generally, sub-Finsler structures.
They play an essential role in various mathematical domains, including metric geometry and geometric group theory -- as we shall see.

Carnot groups are particular examples of metric Lie groups, i.e., Lie groups equipped with left-invariant distances inducing the manifold topology. 
Moreover, metric Lie groups for which the distance is geodesic are precisely the sub-Finsler Lie groups.
In several problems, it is more natural to consider the abstract setting of metric groups. 

In this introductory chapter, we begin by providing a glimpse into the core concepts of Carnot-Carathéodory spaces: contact distributions and sub-Riemannian distances. 
The second section provides an outline of the book’s content and structure, offering a roadmap for our readers. 
While experts in the field can directly proceed to the subsequent chapters, readers new to the topic will find guidance within these pages. 
In Section \ref{sec_example_models}, we present a series of applications in mathematics, physics, and various other scientific disciplines where Carnot-Carathéodory spaces appear, underscoring the versatility and significance of these spaces in real-world scenarios.

\usetikzlibrary{shapes}
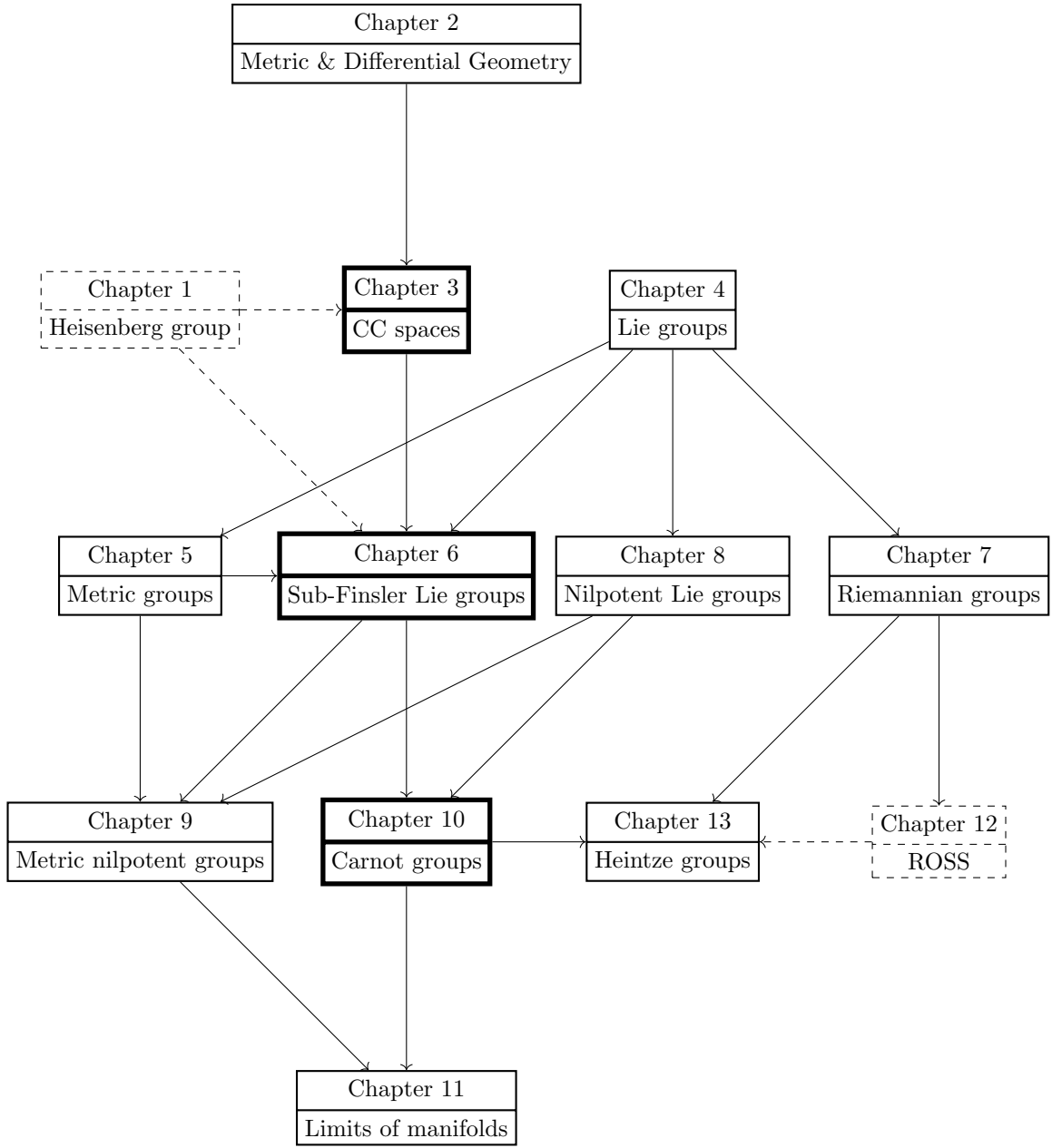
\begin{figure}
\centering
\begin{tikzpicture}[node distance=3.9cm, auto]
 \node[draw, rectangle split, rectangle split parts=2, dashed ] (chapter1) {Chapter~\ref{ch_Heisenberg}\nodepart{second}Heisenberg group};
 \node[draw, rectangle split, rectangle split parts=2, line width=2pt, right of=chapter1] (chapter3) {Chapter~\ref{ch_CCspaces}\nodepart{second}CC spaces};
 \node[draw, rectangle split, rectangle split parts=2, thick, right of=chapter3] (chapterLie) {Chapter~\ref{ch_LieGroups}\nodepart{second}Lie groups};
 \node[draw, rectangle split, rectangle split parts=2, thick, above of=chapter3] (chapter2) {Chapter~\ref{ch_MGgeometry}\nodepart{second}Metric \& Differential Geometry};
 \node[draw, rectangle split, rectangle split parts=2, line width=2pt, below of=chapter3] (chapter5) {Chapter~\ref{ch_SFLie}\nodepart{second}Sub-Finsler Lie groups}; 
 \node[draw, rectangle split, rectangle split parts=2, thick, left of=chapter5] (chapterMetricGroups) {Chapter~\ref{ch_MetricGroups}\nodepart{second}Metric groups}; 
 \node[draw, rectangle split, rectangle split parts=2, thick, below of=chapterLie] (chapterNilpotent) {Chapter~\ref{ch_Nilpotent}\nodepart{second}Nilpotent Lie groups}; 
 \node[draw, rectangle split, rectangle split parts=2, thick, right of=chapterNilpotent] (chapter7) {Chapter~\ref{ch_RiemaLie}\nodepart{second}Riemannian groups};

 \node[draw, rectangle split, rectangle split parts=2, line width=2pt, below of=chapter5] (chapterCarnot) {Chapter~\ref{ch_Carnot}\nodepart{second}Carnot groups};
 \node[draw, rectangle split, rectangle split parts=2, thick, left of=chapterCarnot] (chapterMetricNilpotent) {Chapter~\ref{ch_MetricNilpotent}\nodepart{second}Metric nilpotent groups};
 
 \node[draw, rectangle split, rectangle split parts=2, thick, below of=chapterCarnot] (chapter9) {Chapter~\ref{ch_limits}\nodepart{second}Limits of manifolds}; 
 \node[draw, rectangle split, rectangle split parts=2, thick, right of=chapterCarnot] (chapter13) {Chapter~\ref{ch_Heintze}\nodepart{second}Heintze groups
}; 
 \node[draw, rectangle split, rectangle split parts=2, dashed, below of=chapter7] (chapter12) {Chapter~\ref{ch_ROSS}\nodepart{second}ROSS};

 \draw[->, dashed] (chapter1) -- (chapter3);
 \draw[->, dashed] (chapter1) -- (chapter5);
 \draw[->] (chapter2) -- (chapter3);
 \draw[-> ] (chapter3) -- (chapter5);
 \draw[-> ] (chapterLie) -- (chapter5);
 \draw[-> ] (chapterLie) -- (chapterNilpotent);
 \draw[-> ] (chapterLie) -- (chapterMetricGroups);
 \draw[-> ] (chapterMetricGroups) -- (chapterMetricNilpotent);
 \draw[-> ] (chapterMetricGroups) -- (chapter5);
 \draw[-> ] (chapterNilpotent) -- (chapterMetricNilpotent);

 \draw[-> ] (chapterLie) -- (chapter7);

 \draw[-> ] (chapter5) -- (chapterCarnot);
 \draw[-> ] (chapter5) -- (chapterMetricNilpotent);

 \draw[-> ] (chapter7) -- (chapter13);
 \draw[-> ] (chapter7) -- (chapter12);
 
 \draw[-> ] (chapterNilpotent) -- (chapterCarnot);
 
 \draw[-> ] (chapterMetricNilpotent) -- (chapter9);
 
 \draw[-> ] (chapterCarnot) -- (chapter9);
 \draw[-> ] (chapterCarnot) -- (chapter13);
 
 \draw[->, dashed ] (chapter12) -- (chapter13); 
\end{tikzpicture}
\caption{Interdependence of chapters. Thick boxes represent the main chapters. Dashed boxes represent chapters devoted to examples.}
\end{figure}

\section{What sub-Riemannian geometry is}

 Sub-Riemannian geometry is also known as non-holonomic Riemannian geometry in Russia and in France got the name Carnot geometry, or Carnot-Carath\'eodory, when it is appropriately generalized.
 Since the 1980s, it has been a full research domain, with motivations and ramifications in several parts of pure and applied mathematics.
However, historically, it was not clear that such theories were heading into the same notions. Thus, each source provided its own jargon in the field.
Accordingly, some concepts have multiple terminologies: a contact structure is a particular distribution of hyper-planes in an odd-dimensional manifold, and the concept of Carnot-Carath\'eodory metric is a generalization of a sub-Riemannian distance.

Sub-Riemannian geometry is a generalization of Riemannian geometry. 
In addition to a Riemannian structure, in each sub-Riemannian manifold, there is a constraint on admissible velocities for curves. 
Geometrically, in Riemannian geometry, every smooth curve has locally finite length. 
In sub-Riemannian geometry, curves that fail to satisfy the constraint have infinite length.

Typical examples to keep in mind come from mechanics.
Indeed, the state of a moving object is determined by its position in space and the velocities of its parts: the momenta.
Thus, in the space `positions and speeds' of configurations, the possible evolutions of the object should satisfy the fact that the derivatives of the first coordinates are equal to the second coordinates.
For example, the movement of a single particle in space is described by a curve $t\mapsto (x(t), v(t))$ in $\R^3\times \R^3$. However, not every curve in $\R^3\times \R^3$ is possible. 
Some trajectories are not allowed by the dynamical constraint. 
 As trivial examples, you cannot vary your speed without changing your position or, similarly, 
you cannot move into another place at speed zero!

\begin{figure}
 \centering
 \includegraphics[width=5in]{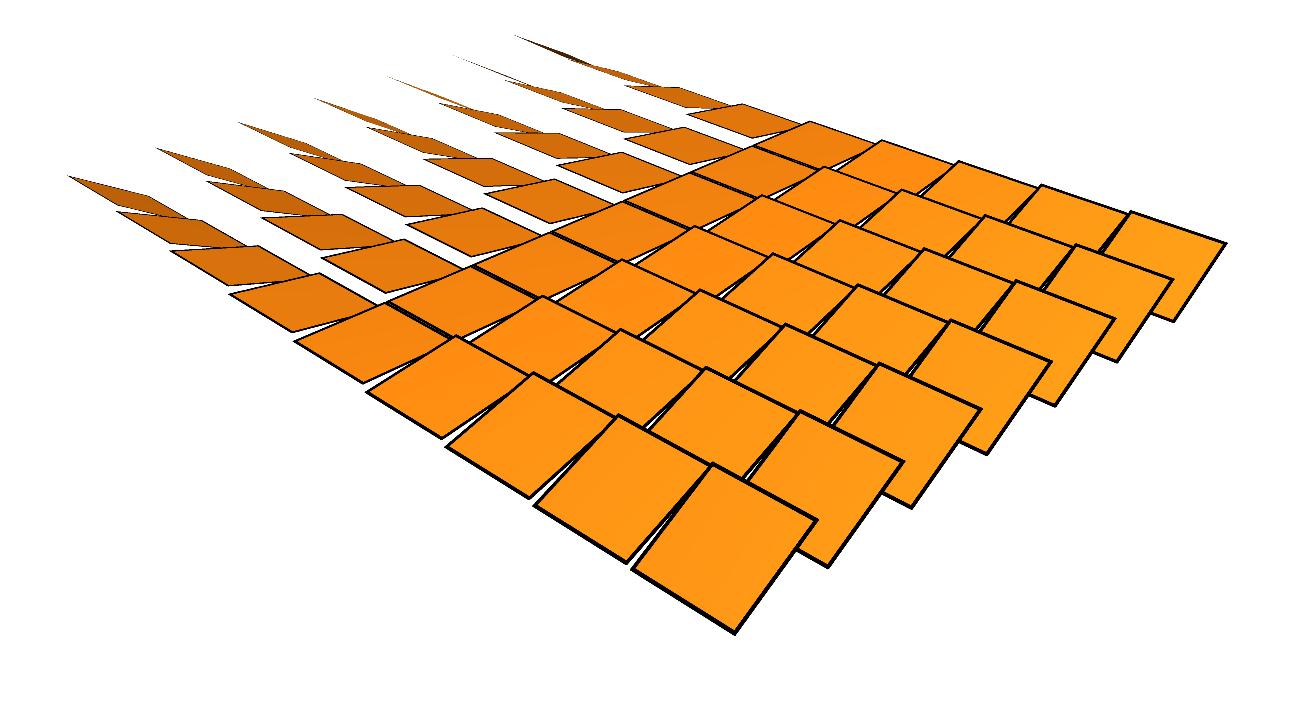}
	\caption{A contact distribution on $\R^3$ is a polarization by planes.}
	\label{fig:contact_structure}
\end{figure}

\index{Heisenberg! -- group}
The 3-dimensional (3D, for short) Heisenberg group equipped with its contact geometry is one of the most essential examples in sub-Riemannian geometry among those that actually are not Riemannian manifolds. 
Visualizing some of its features is relatively easy. 
As a set and topological space, the 3D Heisenberg group $\mathbb{H}$ is equivalent to \(\mathbb{R}^3\). 
The constraint on curves in this space is determined by what is called a `distribution of planes', or a `rank-2 polarization'. 
Similar to how a smooth vector field $X$ smoothly assigns a tangent vector $X_p$ at each point $p$ of a manifold, 
a distribution $\Delta$ of planes in $\R^3$ smoothly assigns to each point $p\in \R^3$ a plane $\Delta_p$ within the 3D tangent space at $p$. 
The curves that we call `admissible' are those that are tangent to one such a distribution, in the sense that a smooth curve $\gamma $ is {\em admissible with respect to a distribution} $\Delta$ if, for every $t$ in the domain of the curve, the velocity vector $\dot\gamma(t)$ belongs to the plane $\Delta_{\gamma(t)}$. Refer to Figure~\ref{fig:contact_structure} for a visual representation of a distribution.

The particular feature of the Heisenberg group $\mathbb{H}$ is that it comes with a distribution that is curly enough in a way that each pair of points can be connected by at least one admissible curve. 
Therefore, one can define a finite-valued distance similarly to the Riemannian case: 
the distance between two points $p$ and $q$ in $\mathbb{H}$ is given by the infimum of the length of all admissible curves from $p$ to $q$,
\begin{equation}\label{intro_distance}\tag{$\star$}
 d(p,q)=\inf\{ \Length (\gamma)\;:\; \gamma \text{ admissible curve from } p \text{ to } q\}. 
\end{equation}

%

More generally, a sub-Riemannian manifold is a triple $(M,g,\Delta)$ consisting of a manifold $M$, a Riemannian tensor $g$ inducing a length structure, and a distribution $\Delta$ of subspaces in the tangent bundle of $M$. These data define admissible curves and a distance as described in \eqref{intro_distance}. These geometries provide a broad generalization of Riemannian geometries in several aspects. Those sub-Riemannian spaces that are non-Riemannian, 
in the sense that we are dealing with a proper distribution,
are rather different from Riemannian spaces:
 They exhibit fractal properties as their Hausdorff dimension exceeds the topological dimension. Additionally, there exist smooth curves with locally infinite length, as well as other smooth admissible curves that are isolated in the topological space of smooth admissible curves with the same endpoints. Consequently, sub-Riemannian geometry requires techniques that are different from those used in Riemannian geometry.

According to Mikhail Gromov, \cite{Gromov1}, the concept of sub-Riemannian distance can be traced back to the ideas of Nicolas Carnot in 1824 and Constantin Carathéodory in 1909; the early references are \cite{Carnot, Caratheodory}. For this reason, sub-Riemannian manifolds are also referred to as Carnot-Carathéodory spaces. Since the 1980s, this geometry has emerged as a vibrant research field with applications and connections to various areas of pure and applied mathematics, including classical mechanics, control theory, metric geometry, group theory, and the analysis of hypoelliptic differential operators.

\section{Content and structure of this text}

%
We shall explore Carnot-Carathéodory spaces from the perspective of Lie groups. The main objective is to illustrate how these possibly non-Riemannian geometries manifest in other mathematical domains, such as metric geometry and geometric group theory, through Carnot groups, which are key examples to keep in mind.
Carnot groups are a class of nilpotent Lie groups equipped with sub-Riemannian or, more generally, sub-Finsler structures.
This text aims to explain the role of Carnot groups as asymptotic cones of finitely generated nilpotent groups.
We also explore their presence as parabolic boundaries of rank-one symmetric spaces and their involvement as limits of Riemannian manifolds and tangents of Carnot-Carath\'eodory spaces.

In Chapter~\ref{ch_Heisenberg}, we will focus on the plane distribution in the 3D Heisenberg group $\mathbb{H}$. We will consider the induced distance \eqref{intro_distance}. Specifically, we will discuss the following facts:
\begin{enumerate}
\item This distance $d$ turns the space $\R^3$ into a metric space. The topology associated with $d$ is the standard topology of $\R^3$. In particular, nearby points can be connected by short admissible curves.

\item 
The distance between every two points is equal to the length of some admissible curve connecting them. If a curve is admissible, its length is comparable to its Euclidean length. However, non-admissible curves have infinite length.

\item This metric space is distinct and different from Riemannian spaces. It is not biLipschitz equivalent to any Riemannian manifold. This is because the Heisenberg geometry exhibits characteristics of fractal geometry. Indeed, the metric on this topologically 3-dimensional space has a metric dimension equal to 4, as determined by its Hausdorff 4-measure.

\item The geometry is homogeneous. Actually, it is invariant under a smooth group structure. 
\end{enumerate}
The general definition of a Carnot-Carathéodory space arises when we formally define the concept of a distribution being `curly enough'. This property should ensure that every pair of points can be connected by an admissible curve. To explore this topic, we require basic notions of both differential geometry and metric geometry, which we review in Chapter~\ref{ch_MGgeometry}.

Then, in Chapter~\ref{ch_CCspaces}, we delve into the sub-Finsler geometry of Carnot-Carathéodory spaces, focusing on their distributions, distances, and dimensions. 
A {\em distribution} on $M$ refers to a sub-bundle of the tangent bundle $TM$ or, more generally, a subset of $TM$ that, locally on the manifold, can be expressed as the span of a collection of vector fields. Constant-rank distributions are also referred to as {\em polarizations}.\index{polarization}
A distribution $\Delta \subseteq TM$ is said to be {\em bracket generating}\index{bracket-generating}
 if, for every $p\in M$, the Lie algebra generated by the sections of $\Delta$ evaluated at $p$ is the entire tangent space $T_p M$. In other words, a distribution $\Delta$ is bracket generating if every tangent vector $v\in T M$ can be represented as a linear combination of vectors of the following form: the evaluation at $p$ of vector fields $X_1$, $[X_2, X_3]$, $[[X_4,[X_5, X_6]]]$, and so on, where all the vector fields $X_1, X_2, X_3, \ldots$ are tangent to $\Delta$, and $v\in T_pM$.
 
A {\em sub-Riemannian manifold} 
\index{sub-Riemannian! -- manifold}
is a triple $(M,\Delta,g)$, where $M$ is a differentiable manifold, $\Delta$ is a bracket-generating distribution, and $g$ is a smooth section of positive-definite quadratic forms on $\Delta$. In fact, 
the map $g$ can be considered as the restriction to $\Delta$ of a Riemannian metric tensor on the manifold $M$.
A curve $\gamma$ on $M$ is called {\em admissible}, or {\em horizontal}, with respect to $\Delta$ if it is absolutely continuous and $\dot\gamma(t)\in\Delta_{\gamma(t)}$ for almost every $t$.\index{admissible! -- curve}\index{horizontal! -- curve} 
Then the {\em sub-Riemannian distance}
 is defined by the same formula (\ref{intro_distance}).
More generally, if the length comes from a smoothly varying norm, then the distance is called {\em Carnot-Carath\'eodory metric} or {\em sub-Finsler}.
Most of the previously mentioned results on the Heisenberg group will be valid for every Carnot-Carath\'eodory space. 
\index{sub-Finsler! -- metric}\index{sub-Riemannian! -- metric}\index{sub-Finsler! -- distance}\index{sub-Riemannian! -- distance}\index{Carnot-Carath\'eodory! -- metric}
\\

The understanding of many Riemannian geometric properties comes from the fact that `metric' tangent spaces $T_pM$ of each Riemannian manifold $M$ are Euclidean spaces, and Euclidean geometry is sufficiently understood. Such a notion of tangent space is defined in terms of limits of metric spaces, and we call them {\em tangent cones} or {\em metric tangents}. What are the metric tangents in sub-Riemannian geometry? The answer is not immediate. Under further assumptions of equiregularity, as we will discuss in Section~\ref{sec:equiregular},  
for each $3$-dimensional (non-Riemannian) sub-Riemannian manifold $M$, each $T_p M$ is isomorphic to the Heisenberg group -- another reason for it to be important. In general, alas, fixed a topological dimension greater or equal to 7, the possible tangents are infinitely many. They may not be the same one at every point, even for a given sub-Riemannian manifold. The good news is that, analogously as the Heisenberg structure has a group structure, the metric tangent of a Carnot-Carath\'eodory space $M$ has a Lie group structure for most points $p\in M$, and at the remaining points, it is still a quotient of some Lie group modulo a closed subgroup.
 The metric tangent at `regular' points has even more structure: it has a dilation property, and consequently, it is a nilpotent Lie group. Such metric Lie groups are called Carnot groups.
 
 We shall review the necessary theory of Lie groups in Chapter~\ref{ch_LieGroups}. 
 Chapter~\ref{ch_MetricGroups} introduces the general objects that combine a group structure with metric geometry: metric Lie groups and, more broadly, isometrically homogeneous spaces.
 In Chapter~\ref{ch_SFLie}, we then discuss the geometry of sub-Riemannian Lie groups and, more generally, sub-Finsler Lie groups.
 This chapter is fundamental since we present Chow's theorem about the fact that Carnot-Carath\'eodory distances induce the manifold topologies, and we discuss the notion of endpoint map, which leads to a first-order study of geodesics in sub-Riemannian manifolds.
 Chapter~\ref{ch_RiemaLie} is a side excursion into the classical basic theory of Riemannian Lie groups.
 Before getting into Carnot groups, we review many properties of nilpotent Lie groups in some detail. Respectively, in Chapter~\ref{ch_Nilpotent}, we discuss the classical differential geometry of nilpotent Lie groups, while Chapter~\ref{ch_MetricNilpotent} is devoted to their metric geometry. 
 In Chapter~\ref{ch_Carnot}, we finally define and study Carnot groups.
 
 The aim of the rest of the book is twofold:
 
 {\bf [Aim 1]} We explore the role of Carnot groups with their Carnot-Carath\'eodory distances in other mathematical areas. Namely, they appear as 
 
 (A) limits of Riemannian manifolds, asymptotic cones of nilpotent groups, and tangents of Carnot-Carath\'eodory spaces; see Chapter~\ref{ch_limits}. 
 
 
 (B) parabolic boundaries of rank-one symmetric spaces and other negatively curved spaces; see Chapter~\ref{ch_Heintze}.
 
 In harmonic analysis on stratified Lie groups, and more generally on graded groups, Carnot-Carath\'eodory distances enter in the study of hypoelliptic differential operators. In complex analysis, Carnot-Carath\'eodory spaces occur as boundaries of strictly pseudo-convex complex domains.
 We do not treat these last two settings in this monograph but refer to the books \cite{Stein:book, Capogna-et-al} as initial references.
 
 As we will explain in Section~\ref{sec_largescale}, Carnot groups with Carnot-Carath\'eodory distances appear in geometric group theory as asymptotic cones of nilpotent finitely generated groups; see
\cite{Gromov1, Pansu}. 
Part of this text is devoted to the study of the coarse geometry of nilpotent groups.
We will see how a geometric notion such as the polynomial growth of balls in the Cayley graph of a discrete group relates with the geometry of the tangent cone at infinity of this graph, which in this case turns out to be a Carnot group endowed with a Finsler-Carnot-Carath\'eodory metric, and eventually gives an algebraic consequence: the group is (virtually) nilpotent. 
 
 There is a general natural explanation for why Carnot groups appear in the above situations. 
 The reason is that Carnot groups are the analogs of finite-dimensional normed vector spaces in the non-commutative case. Indeed, Carnot groups admit the following axiomatic characterization:
 \begin{customtheorem}\label{characterization_Carnot_intro}
The sub-Finsler Carnot groups are the only metric spaces that are locally compact, geodesic, isometrically homogeneous, and admitting metric dilations.
 \end{customtheorem}
 In Chapter~\ref{ch_Carnot}, we will prove this result. We will also previously discuss the setting where the geodesic assumption is replaced with connectedness; see Sections~\ref{sec: self-similar} and~\ref{sec: self-similar2}.
 
 {\bf [Aim 2]} 
 With the use of Lie group theory, one may develop calculus, analysis, geometric measure theory, calculus of variations, and geometric analysis on Carnot groups and, more generally, on sub-Finsler Lie groups. In Chapter~\ref{ch_Carnot}, we shall prove some crucial results in this regard. 
 
This book's approach is to focus on studying specific examples of tangent spaces of sub-Riemannian manifolds (i.e., Carnot groups) to shed light on the general case of Carnot-Carathéodory spaces. 
The reason for this perspective is that Carnot groups
provide tools for developing 
 calculus in such settings, thanks to the availability of translations by group elements and the dilation property. 
The classical definition of the derivative of a real function $f:\R\to \R$ relies on addition, multiplication, and limits:
$$f'(x)=\lim_{h\to0}\dfrac{f(x+h)-f(x)}{h}.$$
All these operations are present in Carnot groups, where addition is replaced with a possibly non-commutative group operation. Consequently, we can define a metric notion of derivative known as the {\em blow-up differential}, or {\em Pansu derivative}, named after Pierre Pansu, who made pioneering contributions to the field in the late 1980s; see \cite{Pansu}.
The main differentiability result, obtained by Pansu and then generalized in \cite{Margulis-Mostow, Vodopyanov}, is the following:
 
 \begin{customtheorem}[Pansu's Rademacher Theorem]\label{Thm_Pansu_intro}
Given a Lipschitz map between sub-Riemannian manifolds, at almost all points, its blow-up differential exists, is a group homomorphism of the tangent cones, and is equivariant with respect to their dilations.
\end{customtheorem}
This theorem will be proved in Chapter~\ref{ch_Carnot} for Carnot groups.
In fact, the theorem also holds for quasi-conformal maps between Carnot-Carath\'eodory spaces; see \cite{Vodopyanov}.
The theory of quasi-conformal mappings has been used in \cite{Pansu} to prove rigidity theorems on hyperbolic spaces over the division algebras of real, complex, or quaternionic numbers.
Indeed, as we shall see in Chapter~\ref{ch_Heintze}, the `parabolic visual boundaries' of rank-one symmetric spaces are Carnot groups. More generally, all negatively curved homogenous Riemannian manifolds have graded groups as boundaries. This last fact is mainly based on the work of Heintze. 
In Chapter~\ref{ch_ROSS}, we discuss rank-one symmetric spaces, which can be seen as semi-direct products of Lie groups. One of the factors of the semi-direct product is a Carnot group of Heisenberg type, and the other one is the one-dimensional simply connected Lie group acting on the Carnot group by its dilations. In Chapter~\ref{ch_Heintze}, we review the notion of boundaries of CAT$(-1)$ spaces and their visual boundaries. First, we observe that boundaries of rank-one symmetric spaces are the particular Carnot group of Heisenberg type. Second, we discuss the general viewpoint of Heintze groups and show that their boundaries are metric Lie groups admitting dilations.

\section{Sub-Riemannian geometries as models}

In this section, we highlight several contexts in which Carnot-Carath\'eodory metrics find applications: either as specific examples, generalizations, or tools. While this overview is not comprehensive, it provides specific references for readers interested in delving deeper into these areas.
 
 \subsection{Examples from mathematics}\label{sec_example_models}

\subsubsection{Control theory and metric geometry} 
Control theory is an interdisciplinary branch of engineering and mathematics that deals with the behavior of dynamical systems. 
The primary objective is to manipulate a differential system by choosing control variables to guide trajectories toward a desired state, possibly optimally.
Sub-Riemannian geometry specifically focuses on control systems that are linear in the controls. Many theorems in sub-Riemannian geometry have broader validity in control theory. For instance, foundational sub-Riemannian theorems such as the ones by Chow, Pontryagin, and Goh have more general statements within the framework of geometric control theory. Readers interested in this perspective should refer to the book \cite{Agrachev_Sachkov}.

 
 However, the distinctive characteristic of sub-Riemannian geometry is the presence of the induced metric, which transforms these manifolds into metric spaces. In this way, methods from metric geometry become applicable in the study of sub-Riemannian geometry, and sub-Riemannian geometries, in turn, offer intriguing examples of metric spaces, sometimes exhibiting strange or even pathological behavior. This geometric perspective on sub-Riemannian geometry was pioneered by Gromov in his seminal work \cite{Gromov1}.
 
 \subsubsection{Geometric group theory and asymptotic geometry}
 Sub-Riemannian geometry has a significant presence in geometric group theory, which is the study of groups from a geometric perspective. 
 Sub-Riemannian structures naturally arise as asymptotic cones of groups with polynomial growth. A group $\Gamma$, generated by a finite set $S$, is said to {\em growth polynomially} if the cardinality of the product set $S^n$ is polynomially bounded in $n\in \N$. One verifies easily that this property does not depend on the choice of the finite generating set $S$.
The asymptotic properties of these groups are intricately connected to the geometry of their asymptotic cones.

For instance, according to a theorem by Pansu, as presented in Section~\ref{sec_largescale} of this book, each nilpotent group $\Gamma$ generated by a finite set $S$ has polynomial growth, and there exist constants $Q\in \N$ and $V>0$ such that 
$$\#(S^n)/n^Q\to V, \qquad \text{ as } n\to \infty.$$
Here, $Q$ and $V$ possess clear geometric significance: $Q$ is the Hausdorff dimension of the asymptotic cone of $\Gamma$, and $V$ is the volume of the unit ball in the asymptotic cone. Furthermore, this asymptotic cone is an example of a Carnot group, equipped with a Carnot-Carath\'eodory distance. These are the spaces on which this book focuses. See Chapters~\ref{ch_Carnot} and \ref{ch_limits}. 
 
 \begin{figure}
 \centering
 \includegraphics[width=4.1in]{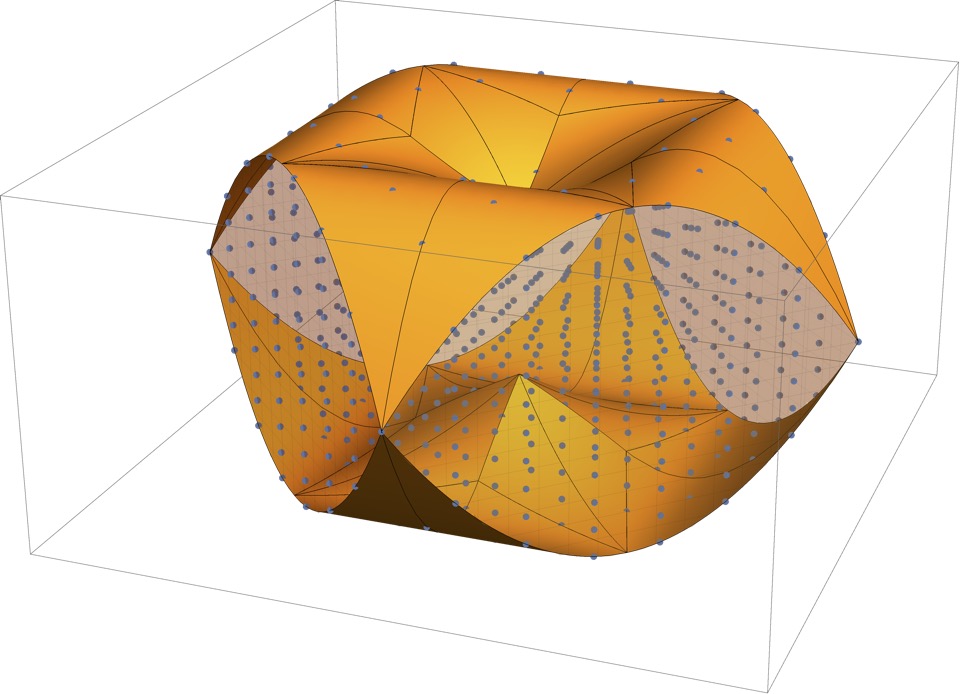}
	\caption{A Carnot-Carath\'eodory ball, which is the limit of large balls for a word distance on a finitely generated nilpotent group. See Example~\ref{example:standard generating set in Heisenberg}.}
\end{figure} 
 
\subsubsection{Complex analysis and Cauchy-Riemann geometry} 
Sub-Riemannian geometry arises when studying the geometry of Cauchy-Riemann (CR) manifolds. Typical examples are domains in complex Euclidean space, $\mathbb{C}^n$.
The boundaries of strictly pseudo-convex domains, of great importance in analysis of several complex variables, are naturally equipped with visual distances that are of Carnot-Carath\'eodory type. 

A domain in $\mathbb{C}^n$ is called {\em strictly pseudo-convex} if, near every point on its boundary, there exists a defining function whose Levi form is positive definite. The Levi form encodes information about the local geometry of the boundary and is related to the complex Hessian of the defining function.

The boundaries of strictly pseudo-convex domains in $\mathbb{C}^n$ exhibit rich geometric and analytic properties. Understanding them is crucial in the study of several complex variables, where it serves as a foundation for topics like plurisubharmonic functions, the Cauchy-Riemann equations, and the theory of complex manifolds.
In this book, we will not discuss strictly pseudo-convex domains, but we refer to \cite{Gromov1, Ma_1991, Balogh_Bonk, Pilastro_master_thesis}.





The main example that the reader should have in mind is the unit ball $\mathbb B$ in $\C^n$. On the one hand, it is an example of a strictly pseudo-convex domain. On the other hand, it is a rank-one symmetry space: Namely, $\mathbb B$ equipped with its Kobayashi-Carath\'eodory distance is the complex hyperbolic space. 
The real, complex, and quaternionic hyperbolic spaces,
together with the octonionic plane, exhaust the list of the rank-one symmetry spaces of non-compact type. All of these spaces carry a natural metric on the boundary, which is a Carnot-Carath\'eodory distance.
We present rank-one symmetry spaces from the Lie group viewpoint as semi-direct products in
Chapter~\ref{ch_ROSS}.
 Then, in Chapter~\ref{ch_Heintze}, we discuss the visual boundaries of these spaces and, more generally, of all the negatively curved homogeneous spaces. 

\subsubsection{Analysis of hypoelliptic operators and singular integral operators} 

The theory of partial differential equations (PDEs) deals with a wide variety of equations, each requiring different approaches.
One of the most important operators is the Laplacian, which, in Euclidean space $\R^n$, is defined as:
$$f\longmapsto \Delta f:= \frac{\partial^2}{\partial x_1^2} f +\ldots + \frac{\partial^2}{\partial x_n^2} f .$$
This operator is called {\em hypoelliptic} because it has the property that 
$$\Delta f\in C^\infty\qquad\Longrightarrow \qquad f\in C^\infty.$$
There is a connection between the Laplace operator in Euclidean space and the Euclidean distance.
Indeed, the fundamental solution to the Laplace equation, called the {\em Green's function}, is written as a function of the Euclidean distance. 

When considering differential equations defined by bracket-generating vector fields, 
the Laplacian generalizes to the {\em sub-Laplacian}.
An example is the operator on functions on $\R^3$ given by 
$$f\longmapsto \Delta_{\rm sub} f:= X(X f) +Y(Y f ) ,$$
where 
$$X:=
\frac{\partial }{\partial x_1} -x_2\frac{\partial }{\partial x_3} \quad \text{ and }\quad 
Y:=\frac{\partial }{\partial x_2} +x_1\frac{\partial }{\partial x_3} . $$
The sub-Laplacian takes into account the non-holonomic constraints imposed by the vector fields.
In fact, it is intrinsically linked to the geometry of the sub-Riemannian metric obtained from the bracket-generating vector fields. Namely,  certain regularity bounds for the sub-Laplacian, and many other subelliptic PDEs, can be obtained
 in terms of the Sobolev spaces with respect to the associated sub-Riemannian distances.

This relation is used to study diffusion processes, heat equations, singular integrals, and other differential equations on such manifolds or on spaces of homogeneous type.
For more on the analysis of hypoelliptic operators, one can consult the following references:
 \cite{Folland73, Roth:Stein, Goodman, nagelstwe, Capogna-cpam}.
 

\subsubsection{Classical mechanics and non-holonomic dynamical systems}
Sub-Riemannian geometry plays a significant role in classical mathematical mechanics, particularly when studying mechanical systems with constraints. These systems, referred to as {\em non-holonomic systems}, involve specific limitations on velocities or accelerations, constraining the possible motions of the system according to defined rules. Sub-Riemannian geometry offers a comprehensive framework for analyzing and understanding the geometric properties of these constrained mechanical systems.
Some classical references on non-holonomic dynamical systems are
\cite{MR1562963, MR0003103, MR0043608, MR0661282}. A recent review is \cite{MR3531848}; see also \cite{Bloch,MR1999579}.
 
\subsubsection{Contact geometry, Engel structures, and exterior differential systems}

Contact structures, fundamental in various mathematical contexts, have their roots in Carathéodory's formalization of thermodynamics \cite{Caratheodory}. They are instrumental in understanding processes like the {\em Carnot cycle}, where work is transformed into heat due to the bracket-generating property of the distribution that is constraining the dynamics.
 
Contact manifolds are particular polarized manifolds that involve certain distributions on odd-dimensional manifolds.
The distributions are hyperplane fields defined by certain differential one-forms, called {\em contact forms}. 
 The Heisenberg group, as discussed in Chapter~\ref{ch_Heisenberg}, is the standard example of a contact structure.
Introductions to the subject can be found in \cite{Geiges_MR2397738, Etnyre_intro}.

In the realm of low-dimensional topology and geometry, equiregular bracket-generating distributions of dimension 2 on 4-dimensional manifolds are of particular interest. These are known as {\em Engel structures}; see Section~\ref{examples_Carnot} for the model space: the {\em Engel group}.\index{Engel! -- structure}\index{Engel! -- group}
For further insights into these structures and their relation to contact structures, refer to \cite{Montgomery_MR1736216}.
 
The theoretical foundations of contact structures are rooted in {\em exterior differential systems}, systems of equations on a manifold defined by exterior differential forms. This theory, pioneered by Élie Cartan \cite{MR0016174}, serves as a general framework for understanding contact manifolds and many other polarized manifolds. For a comprehensive exploration of this subject, we refer to \cite{Bryant_MR1083148}.\index{exterior differential system}

In some sense, sub-Riemannian geometry metrizes contact manifolds and, more generally, polarized manifolds, similar to how Riemannian geometry metrizes differentiable manifolds. In the presence of a distance function, problems concerning the contact equivalence of contact structures transform and extend to queries like bi-Lipschitz equivalence or quasi-conformal equivalence.

\subsubsection{Riemannian geometry}
%
 
 Sub-Riemannian geometry is a natural generalization of Riemannian geometry. Sub-Riemannian metrics often emerge as limit cases of Riemannian metrics.

A crucial instance of this phenomenon occurs in the asymptotic cones of Riemannian Lie groups with polynomial growth. Specifically, these Lie groups are of type (R), and their asymptotic cones are Carnot groups. 
Further details on these results will be explored in Sections~\ref{sec type (R)} and~\ref{sec_largescale}. 

\subsubsection{Univalent function theory}
There is a remarkable recent application of sub-Riemannian
geometry to univalent function theory. 
The following quick summary is based on 
the article \cite{Markina_Prokhorov_Vasilev}
and on kind 
help from 
Jeremy Tyson.

Classical univalent function theory considers the class $\mathcal S$ of injective analytic functions $f:\mathbb D\to \C$ defined in the unit disc $\mathbb D$ in the complex plane $\C$, normalized by $f(0)=0$ and
$f'(0)=1$. 
A major subject of investigation is the so-called {\em coefficient body}
$$M := \left\{(c_k)\in \C^\N:
f\in \mathcal S, f(z) = z \left(1+ \sum_{k=1}^\infty c_k z^k \right)\right\}.$$
This set can be seen as the limit of its {\em finite-dimensional slices}:
\begin{eqnarray*}
M_n := \left\{(c_1,\ldots,c_{n}): f\in \mathcal S, f|_{\mathbb S^1}\in C^\infty(\mathbb S^1), f(z) = z \left(1+ \sum_{k=1}^\infty c_k z^k \right)
 \right\}.
\end{eqnarray*}

By a famous result by de Branges \cite{MR0772434} (formerly known as the Bieberbach conjecture), we have $|c_n| \le n+1$, for all $n$. This gives information on the size of $M_n$ and $M$. However, there is no explicit description of $M_n$ except for the trivial cases $n=2$ and the case $n=3$; see
\cite{MR0037908}.

One of the basic tools in the subject is the Loewner (or Loewner-Kufarev) parametric representation, which embeds each function $f \in \mathcal S$ into an ODE flow within the class $\mathcal S$. 
Loewner parametrizations were used by de Branges in his proof. 
Nowadays, there is a stochastic version of the Loewner flow, called SLE, which stands for either the {\em stochastic Loewner equation} or for {\em Schramm--Loewner evolution}, a topic at the intersection of probability, complex analysis, stochastic PDE, mathematical physics; see \cite{
MR2087784, 
 MR1776084, 
 MR2334202,
MR2768999}. 

I. Markina, D. Prokhorov, and A. Vasilev used the Loewner flow on $\mathcal S$ to define a natural (partially integrable) Hamiltonian system on each set $M_n$. They found certain first integrals of the flow and calculated all the relevant commutators. Consequently, they constructed a complex sub-Riemannian structure on $M_n$ that is naturally adapted to the underlying univalent function theory. In fact, the Loewner trajectories become horizontal curves with respect to this sub-Riemannian structure.

An interesting problem in the field is to extend Markina-Prokhorov-Vasilev's setup to include SLE as well as the classical (deterministic) Loewner equation.
We refer to \cite{MR3203100} for more on classical and stochastic Loewner-Kufarev equations.

\subsubsection{Diffusion processes and financial mathematics} 

Sub-Riemannian geometry has also found applications in the field of 
diffusion processes and stochastic calculus, particularly in the modeling and analysis of financial systems. Here are some ways in which sub-Riemannian geometry has been used in this context, as it has been in part explained to the author by Josef Teichmann and by Andrea Pascucci; see, for example \cite{Friz_MR2604669, TeichmannANDco, Pascucci_MR2791231, Pascucci_MR2176727}. 


In mathematical finance one often deals with controlled ordinary differential equations and with stochastic processes, in particular, with Ito diffusions. It is well known that Ito diffusions can be related to Riemannian and sub-Riemannian structures. When high-dimensional markets are modeled, it is quite natural to assume that the instantaneous covariance matrix is degenerate but that densities still exist, which leads to sub-Riemannian structures.

Another significant stream of literature involving sub-Riemannian structures is rough path theory and its applications ranging from mathematical finance to machine learning. In its very foundations, free-nilpotent Lie groups play a role, for instance, when strategies are approximated by signatures, which are natural maps with values in nilpotent Lie groups.

To learn more about signature-kernel methods in machine learning, shuffle product, and SiG stochastic differential equations in mathematical finance, the reader might consult the following references: 
\cite{Bayer_MR4551549, Chevyrev_MR3572331, Chevyrev_MR4577129, Cuchiero_MR4628930, Cuchiero_MR4506280, Friz_MR4174393, Friz_MR2604669, Hambly_MR2630037, Kiraly_MR3948071, Lyons_MR4070272, Lyons_MR4260459, Salvi_MR4309861}.


Sub-Riemannian geometry has been employed to model the evolution of asset prices and the dynamics of some financial markets.
For example, in the Asian market, the price of the options depends on the evolution. 
Namely, asset-price equations are a class of path-dependent options characterized by a payoff that is a function of the history of the underlying asset price. The equations are strongly degenerate partial differential equations in three dimensions; see
\cite{Anceschi_MR4273088, Kim_MR2472945, Barucci_MR1830951}.

In a simplified exposition of Asian option pricing, we may say that some variables are the integrals of other variables. 
 Formally, we have a differential constraint implying that the dynamics are tangent to a proper subbundle of the manifold of the variables. The subbundle is bracket-generating. The equation that gives the dynamics is known as the Black--Scholes equation. It is an equation of Kolmogorov-Fokker-Planck type, which is a type of sub-Riemannian Laplacian with respect to Hörmander vector fields. When looking at the fundamental solution of these equations, the sub-Riemannian metric comes naturally.
 
The Black \& Scholes model for the pricing problem goes back to the early 1970s; see the reprint \cite{Black_MR3235225}.  
In 1973, Myron Scholes, in close collaboration with Fischer Black, developed what we call nowadays the {\em Black-Scholes model}. For his contributions, he earned the 1997 Nobel Prize in economics.
The model is generally used to determine the fair price or theoretical value for a call or a put option. 
In the case of path-dependent options, the model leads to Kolmogorov operators associated with some stochastic processes. 




General references, also in connection with stochastic equations, are  
\cite[Chapter 9.5]{Pascucci_MR2791231} 
and the review article \cite{Pascucci_MR2176727}.
In addition, 
we suggest the following topics and references. For applications to Asian options, there is
\cite{Barucci_MR1830951}, 
while for American options, there is \cite{Pascucci_MR2433509}. 
For path-dependent volatility models of asset prices,
stochastic models of stock prices,
 originally underlying models, and the geometric Brownian motion with constant volatility, 
 see
 \cite{Foschi_MR2385732, Hobson_MR1613366}.
For an introduction to the probabilistic theory of arbitrage pricing of financial derivatives, including stochastic optimal control theory, optimal stopping theory, and arbitrage theory in continuous time, we suggest \cite{Bjork_book} as many other books by Tomas Björk. 
The paper \cite{Pascucci_MR2176727} contains a survey of results about partial differential equations of Kolmogorov type arising in physics and mathematical finance.

\subsection{Examples from physics}
Sub-Riemannian geometry models various structures beyond pure mathematics, from finance to mechanics, from bio-medicine to quantum
phases, from robots to falling cats! 
Lacking space as well as competence, in this section, we do not enter in details. Instead, we direct the interested reader to suitable references.

\subsubsection{Geometry of principal bundles with connections: falling cats}
Theoretical physics describes most mechanical systems by a kinetic energy and a potential energy.
Gauge theory, also known as the geometry of principal bundles with connections, studies systems with physical symmetries, i.e., when there is a group $G$ acting by isometries on a configuration space $M$. 
Most of the time, it is easier to understand the {\em dynamics up to isometries}. Namely, one first studies the system trajectories in the quotient space $G/M$. Subsequently, one has to study the `lift' of the dynamics into the initial configuration space to obtain the trajectories in $M$. Such lifts will be subject to a sub-Riemannian restriction. This viewpoint is elaborated in \cite[Part II]{Montgomery}.

The formalism of principal bundles with connections is well presented by the example of the fall of a cat.
A cat dropped upside down will land on its legs. The reason for this ability is the good flexibility of the cat in changing its shape.
 \begin{figure}
 \centering
\subfigure[A photo.] 
{
 \includegraphics[height=3.5in]{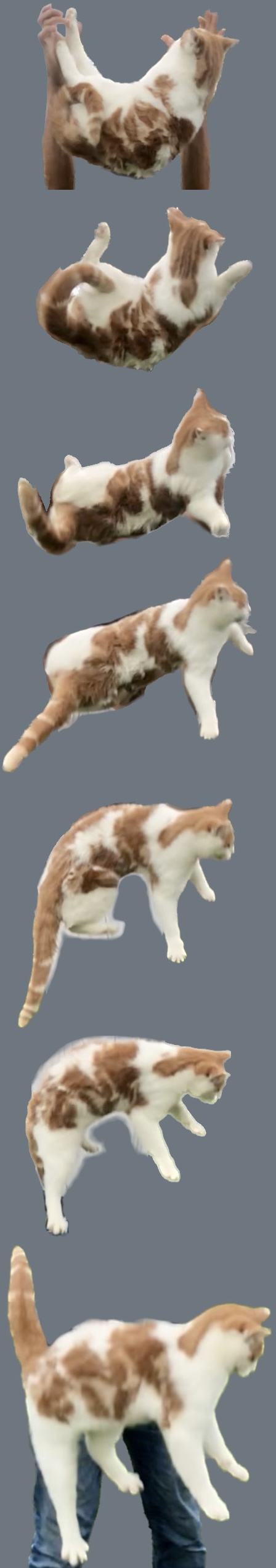} 
}
 \hspace{1cm}
\subfigure[A cartoon.] 
{
 \includegraphics[height=3.5in]{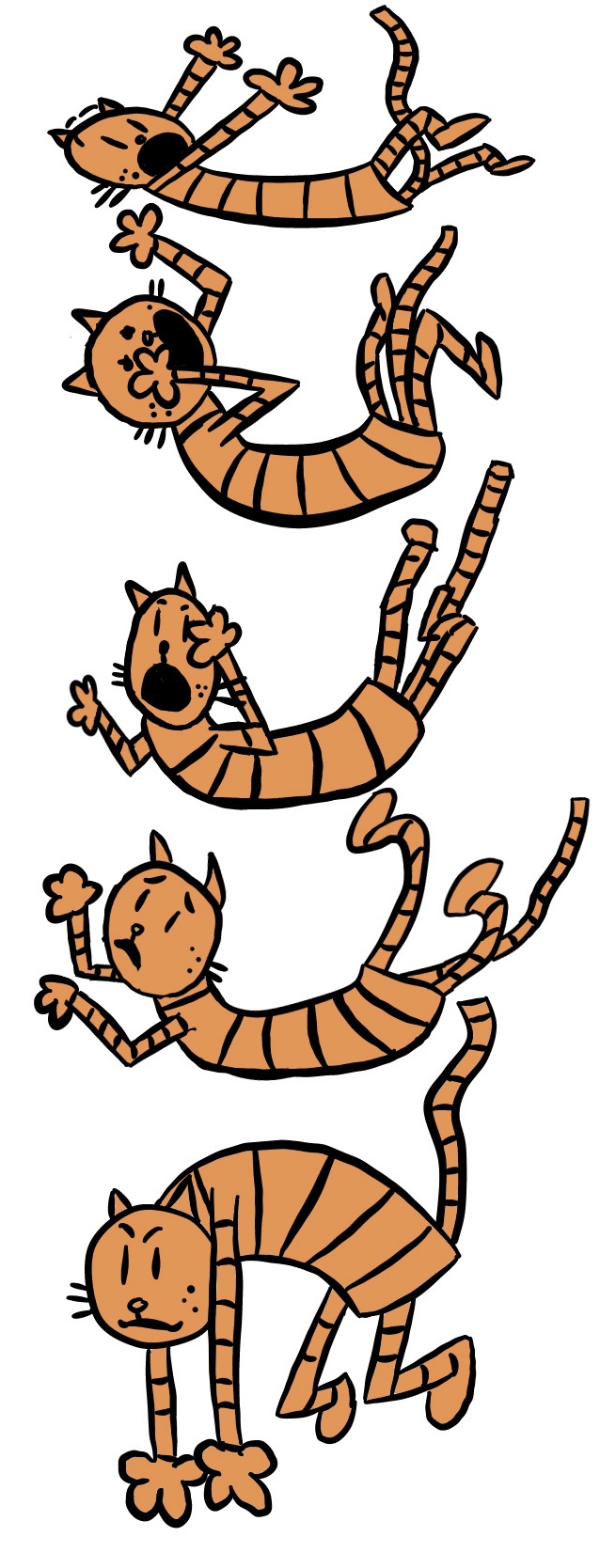}
 }
 \hspace{1cm}
\subfigure[A sketch.] 
{
 \includegraphics[height=3.5in]{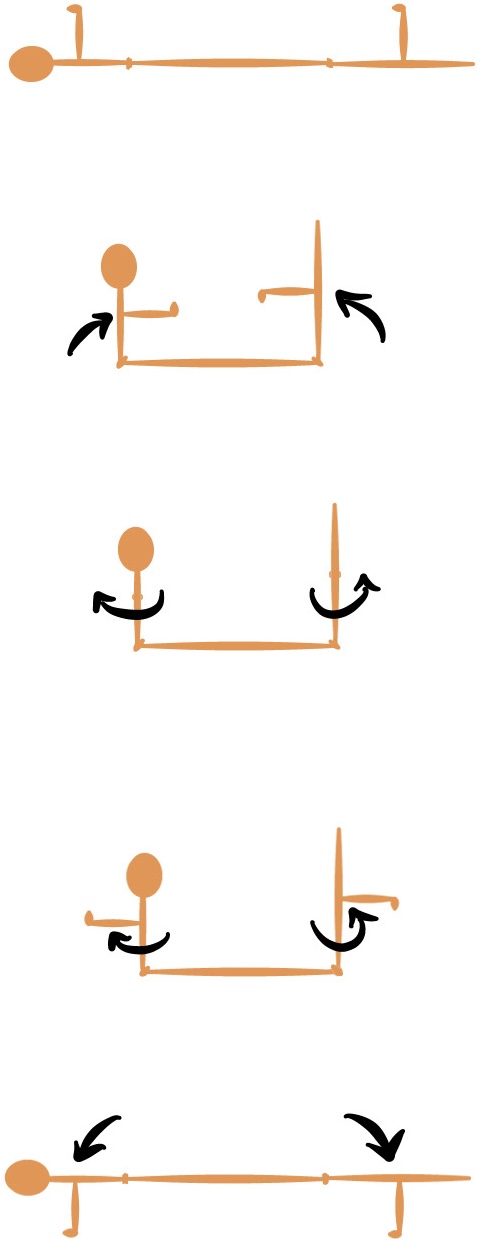}
 
}
\caption{The cat spins itself around and right itself.}
\label{cat} 
\end{figure} 

Let us fix some formalism. Let $M$ be the set of all the possible configurations in the 3D space of a given cat. Let $S$ be the set of all the shapes that a cat can assume. Both $M$ and $S$ are manifolds of large dimensions. 
The position of a cat is determined by its shape plus its orientation in space. Otherwise said, the group $G:=\Isom(\R^3)$ of isometries of the Euclidean 3D space acts on the configuration space $M$, and the shape space is just the quotient of the action:
$$\pi : M\to M/G=S.$$
In fancy words, we say that $M$ is a principal $G$-bundle.

The key fact is that the cat has considerable freedom in deciding its shape $\sigma(t)\in S$ at each time $t$.
However, during the fall, each strategy $\sigma(t)$ of changing shapes will give, as a result, a change in configurations $\tilde\sigma(t)\in M$. The curve $\tilde\sigma(t)$ satisfies 
$$\pi(\tilde\sigma)=\sigma.$$
Moreover, the lifted curve is unique: it has to satisfy the constraint given by what is called the {\em natural mechanical connection}. 
What the cat is demonstrating is that such a connection has non-trivial holonomy. In other words, the cat can choose to vary its shape from the standard normal shape into the same shape, resulting in a change in configuration: the legs initially point towards the sky, but in the end, they are toward the floor.
We refer to \cite[Section 14.2]{Montgomery} for an explanation of an engineers' model for the cat and to Figure~\ref{cat} for an action shot, a cartoon, and a schematic drawing of the drop. 

\subsubsection{From mechanics: parking cars, rolling balls, moving robots, and satellites}

Sub-Riemannian geometry has been extensively used in mechanics and robotics. The study of sub-Riemannian structures provides a mathematical foundation for analyzing the motion planning and control of underactuated mechanical systems. These systems have fewer control inputs than the degrees of freedom, leading to nontrivial constraints on the achievable motions. By understanding the sub-Riemannian geometry associated with such systems, researchers can develop efficient control strategies for navigating robots in complex environments.

{\it Parking a car or riding a bike.} The configuration space is $3$-dimensional: the position in the $2$-dimensional street plus the angle with respect to some fixed line. However, the driver has only two degrees of freedom: turning and pushing. Again, using non-trivial holonomy, we can move the car to any position we like. In this model, we only parametrize the back wheel. Taking into account also the front wheel, we get a 4D manifold with still only two controls.
\\

 \begin{figure}
 \centering
 \includegraphics[width=4.5in]{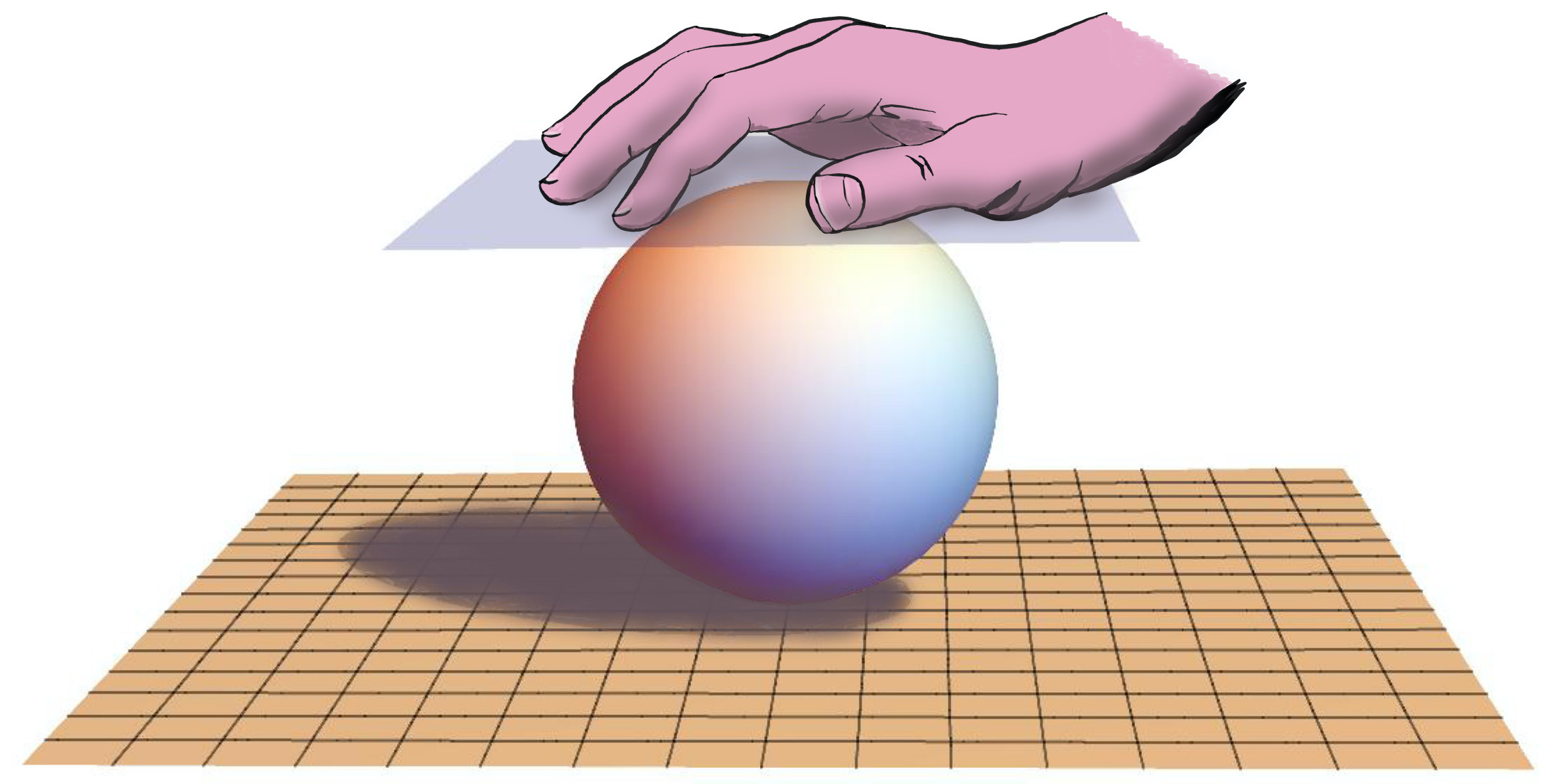}
	\caption{A ball rolling on the plane without sliding, nor spinning.}
\end{figure}
{\it Rolling a ball on the plane.} The position of a ball lying on a plane is described by five coordinates: two real numbers to characterize the point on the plane that touches the ball, another two spherical coordinates describing the point of the ball that touches the plane, and finally, an angular coordinate for spinning the ball around its vertical axis.
When one rolls the ball without sliding, there are only three admissible control directions: two to choose a direction in which to roll the ball, and the third one for spinning it. 
One can also omit the possibility of spinning and have only a 2-dimensional space of control variables.
In both cases, the ball can arrive at every position regardless of its initial position.
This problem has been studied intensively. 
One is interested in minimizing the length of the curve traced on the plane.
The locally optimal trajectories are known and are Gauss' elastic curves. However, it is not known up to which length each such curve is optimal. To more information see \cite{MR1240578, MR2907816,MR2390213, MR1791692}.
\\

In \textit{robotics}, mechanisms are subject to constraints of movement that do not restrict the manifold of possible positions, like for the arm of a robot. 
We refer to \cite{Bellaiche_Jean_Risler, jeancontrol, MR1300410, MR2233026} for the motion planning problem.
The situation of \textit{satellites} is similar.
One should really think about a satellite as a falling cat: it should choose its strategy of modifying the shape properly to have the necessary change in configuration. Another similar example is the case of an astronaut in outer space. See \cite{MR2031058, MR2957013} for applications to aerospace engineering.

\subsubsection{Quantum control and quantum information} 
Sub-Riemannian geometry has been utilized in the control and manipulation of quantum systems. Quantum control aims to steer quantum systems to desired states by applying suitable control fields. Sub-Riemannian structures naturally arise when considering the controllability of quantum systems subject to constraints on the available control resources. By applying techniques from sub-Riemannian geometry, researchers can design control protocols for quantum systems, which find applications in quantum computing, quantum information processing, and quantum sensing.

The following application comes from conversations with Ugo Boscain and from reading his `Habilitation \`a diriger des recherches', \cite{Boscain_Habilitation}.

Let $\H$ be a separable complex Hilbert space. 
Let us denote by $\mathbb S$ the unit sphere in $\H$.
In quantum mechanics
the time evolution of a quantum mechanical system (e.g., an atom, a molecule, or a system of particles with
spin) is described by a map $\psi: \R \to \mathbb S$, called {\em wave function}. The vector $\psi(t)$ is called the {\em state of the system} at time $t$.
The equation of evolution of the state is the so-called {\em Schr\"odinger equation}. If the system is isolated, the equation has the form:
$$i\dfrac{\dd \psi}{\dd t} (t)= H_0 \psi(t),$$
where $H_0$ is a self-adjoint operator acting on $\H$ called the {\em free Hamiltonian}.

For simplicity of notation, let us assume that the spectrum of $H_0$ is discrete
and non-degenerate, with eigenvalues $E_1, E_2,\ldots$ (called {\em energy levels}) and eigenvectors $\psi_1, \psi_2,\ldots\in\mathbb S$.

Assume now to act on the system with some external fields (e.g., an electromagnetic field) whose amplitude is represented by some functions $u_1,\ldots, u_m\in L^\infty(\R;\R)$. 
In this case, the Schr\"odinger equation becomes
$$i\dfrac{\dd \psi}{\dd t} (t)= H(t) \psi(t), \quad \text{ where } H(t) := H_0 +\sum_{j=1}^m u_j(t)H_j,$$
and $H_j$ are self-adjoint operators representing the coupling
between the system and the external fields. 
The time-dependent
operators $H(t)$ and $\sum_{j=1}^m u_j(t)H_j$ are called the {\em Hamiltonian} and the {\em control Hamiltonian}, respectively. 
The typical problem of quantum control is the so-called {\em population transfer problem}:

{\it Assume that at time zero, the system is in an eigenstate $\phi_j$ of the free Hamiltonian $H_0$. Design controls $u_1,\ldots, u_m$ such that at a fixed time $T$, the system is in another prescribed eigenstate $\phi_l$ of $H_0$.}

Nowadays, quantum control has many applications in chemical physics, in nuclear magnetic resonance (also
in medicine) and it is central in the implementation of the so-called {\em quantum gates} (the basic blocks of a
quantum computer); see \cite{MR2808612, MR2357229, MR2761009, MR2808612}. 

For a finite-dimensional quantum mechanical system, if $n$ is the number of energy levels, then we have $\H = \C^n$ and the state space $\mathbb S$ is the unit sphere $\mathbb S^{2n-1} \subset \C^n$.
In this setting, problems of quantum
mechanics can be formulated with matrices. The solution is of the form
$$\psi(t) = g(t) \psi(0), \quad \text{ 
with } g(t) \in SU(n).$$ The Schr\"odinger equation becomes
$\dfrac{\dd}{\dd t}g(t) = -iH(t)g(t)$, 
and now $-iH(t)$ is a skew trace-zero Hermitian matrix, i.e., belongs to the Lie algebra $\mathfrak{su}(n)$.

The controllability problem (i.e., proving that for every couple of points in
$SU(n)$ one can find controls steering the system from one point to the other) is nowadays well understood.
Indeed, the system is controllable if and only if the H\"ormander's condition holds:
$${\rm Lie} \{iH_0, iH_1,\ldots, iH_m\} = \mathfrak{su}(n).$$

Once that controllability is proved, one would like to steer the system between two fixed points in the state
space in the most efficient way. In applications, typical costs that are to be minimized are:
either the energy transferred by the controls to the system or the time of transfer. 
 
\subsubsection{Quantum Berry's phases}
Berry's phase is a phase factor that accumulates during the adiabatic evolution of quantum systems. It arises when a quantum system undergoes slow changes while staying in its instantaneous eigenstate. The connection with sub-Riemannian geometry arises in systems with degenerate energy levels, where the parameter space exhibits a geometric structure with the structure of a fiber bundle. 
The constraints imposed by the structure affect the system's evolution and lead to the accumulation of Berry's phase.
The interested reader should start with the introduction and Chapter 13 in \cite{Montgomery} and the references therein.

\subsection{Scientific applications}

\subsubsection{Neurobiology}


Neurobiological research over the past few decades has enhanced our understanding of the functional mechanisms of the first layer of the visual cortex. This layer comprises various types of cells, including the so-called `simple cells', which are sensitive to specific orientation-based brightness gradients. Recently, this cortical structure has been modeled using a sub-Riemannian manifold. The following summary of this sub-Riemannian application is based on \cite{Sarti-Citti-Petitot} and \cite{MR3729403}, as well as discussions with S. Pauls and G. Citti.


The modeling space is $\mathbb{R}^2 \times \mathbb{S}^1$, where each point $(x, y, \theta)$ represents a cell associated with the point of retinal data $(x, y) \in \mathbb{R}^2$ and oriented according to the angle $\theta \in \mathbb{S}^1$. Namely, the vector $(\cos\theta,\sin\theta)$ indicates the direction of the maximum rate of change of brightness at the point $(x, y)$ in the image perceived by the eye. In pictures with high contrast, this vector is normal to the image contours. 

The moral is that when an image stimulates the cortex cells, the border of the image gives a curve inside this 3D space. Such curves are constrained to be tangent to the distribution spanned by the vector fields\label{roto_transl_neuro}
 $$X_1 = \cos (\theta) \partial_x +\sin (\theta) \partial_y \quad \text{ and } \quad X_2 =\partial_\theta.$$
These vector fields are left-invariant under a Lie group structure. The Lie group is also known as {\em rototranslation group}.\index{rototranslation group}
Researchers think that if a piece of the contour of a picture is missing to the eye vision (possibly by been covered by an object), then the brain tends to `complete' the curve by minimizing some kind of energy. In other words, there is some sub-Riemannian structure in the space of visual cells, and the brain constructs a sub-Riemannian geodesic between the endpoints of the missing data. For more about the neurogeometry of vision, we refer to the book \cite{MR3729403}.


\subsubsection{Image processing and computer vision tracking}
 Sub-Riemannian geometry has been applied in image processing and computer vision. By representing images as curves in a suitable space, sub-Riemannian techniques allow for the extraction of intrinsic geometric features that are invariant under certain transformations. This enables robust object recognition, shape analysis, and image-matching algorithms, which can be applied in fields like pattern recognition, medical imaging, and image-based navigation.

Image processing has benefited from sub-Riemannian geometric methods with the use of geometric flows \cite{MR3463049, 10.1007/978-3-030-00928-1_50}, hypoelliptic diffusion \cite{MR3196939}, and energy minimization \cite{MR3011633}.
For anthropomorphic image reconstruction, see
\cite{MR2968057}.

A very applied problem that is currently tackled with sub-Riemannian methods is the tracking of blood vessels on spherical images of the retina.
Namely, given (several) images of eyes, we seek algorithms that recognize the blood vessels as curves that possibly pass one over the others. The problem can be solved via sub-Riemannian geodesics, 
 see \cite{
MR3726782,
MR3737794}.
Using nilpotent approximations (namely, Carnot groups as studied in this book), researchers have obtained 
fast perceptual grouping of blood vessels in 2D and 3D, 
see \cite{MR3812952}.


\subsubsection{Neuroscience: brain modeling and deep learning} 
Sub-Riemannian geometry has been applied to model and analyze the connectivity and activity patterns in the brain's neural networks. The brain's white matter, which consists of axonal bundles, can be viewed as a sub-Riemannian manifold, where the propagation of nerve impulses is subject to constraints imposed by the underlying anatomy. By studying the sub-Riemannian geometry of neural networks, researchers can analyze brain connectivity and information processing.

The theory of non-isotropic propagation describes the brain propagation of the signals that occur along the dense network of axons that constitutes the neural connectivity. Visual perception phenomena are expressed as anisotropic partial differential equations. The first differential models of the cortex go back to \cite{MR0277268, MR1272050, MR1697185}. 
These geometric brain models have been developed in \cite{Citti_Sarti, MR3318290, MR3865031}, where each family of cells is described via a Lie group with a sub-Riemannian metric. See also \cite{MR2663001, MR2663002, Duits2023} for results in image analysis. 
 
The geometric analysis of sub-Riemannian Lie groups led to group-invariant brain-inspired architectures of deep learning: group-equivariant convolutional neural networks. The considerable advances are in terms of parameters to be trained; see 
\cite{pmlr-v48-cohenc16, 10.1007/978-3-030-00928-1_50, Sifre2014RigidMotionSF}. 
 A more recent brain-inspired architecture is PDE-group equivariant convolution neural networks \cite{10.1007/978-3-030-75549-2_3},
 which have shown to have major advantages in terms of efficiency.

 
\chapter{The main example: the Heisenberg group}\label{ch_Heisenberg}

%
%

The sub-Riemannian Heisenberg group is the first prominent example of sub-Riemannian geometry that deviates from the Riemannian framework. Such a geometry is connected to the solution of the isoperimetric problem on the plane and has a formulation in terms of contact geometry.

In this chapter, we present the geometric models of the sub-Riemannian Heisenberg group and explore certain properties that will be further examined in Carnot groups. Given that the topological dimension of the Heisenberg group is 3, visualizing its sub-Riemannian geodesics and spheres becomes relatively simple.

\section{An isoperimetric problem on the plane}

The {\em isoperimetric problem} is a mathematical challenge where the goal is to find the maximum area among domains with a fixed length as perimeter. In our study, we will focus on a specific variation of the standard isoperimetric problem known as the problem of Dido.


Dido, as described in ancient Greek and Roman sources, is renowned as the founder and first queen of Carthage, located in modern-day Tunisia. Her story is famously depicted in the epic poem Aeneid by the Roman poet Virgil. According to this account, King Jarbas was convinced by Dido to grant her a parcel of land along the African coast for settlement. The condition set forth was that Queen Dido could claim as much land as she could enclose with a leather string, utilizing the coastline as part of the boundary. The optimal solution to maximize the area in this scenario involves a half-circle, assuming the coastline can be treated as a straight line. Actually, whatever the shape of the coastline is, the leather string will take the shape of an arc of a circle. 
\\

We next provide a mathematical model of such a problem. In $\R^2$ with coordinates $(x,y)$, the area form is denoted by $ \vol:=\dd x\wedge \dd y$, which is the differential of the differential one-form 
$$\alpha:=\dfrac{1}{2} (x \dd y-y \dd x)=\dfrac{1}{2}r^2 \dd\theta,$$
where the latter is the expression in polar coordinates $(r, \theta)$.
By applying Stokes' Theorem, we deduce that if a closed, smooth, counterclockwise-oriented curve $\gamma$ in $\R^2$ encloses a domain $D_\gamma$, then the area of $D_\gamma$ is equivalent to the line integral of $\alpha$ along $\gamma$:
$$\Area(D_\gamma):=\iint_{D_\gamma} \vol=\int_\gamma\alpha.$$

\begin{figure}[ht]
 \centering
\begin{tikzpicture}[scale=1.8]
 \fill[blue!20, domain=0:2, variable=\x]
 (0,0) .. controls (0.5,2.7) and (3.5,0.3) .. (2,2)
 -- (2.5,2.5) -- (-1,-1) -- cycle;
 
 \draw[->] (0,-0.2) -- (0,2.5) node[left]{$y$};
 
 \draw[->] (-0.5,0) -- (3.5,0) node[above]{$x$};
 
 \draw[thick, blue] (-0.1,-0.1) -- (2.5,2.5);
 
 \draw[red, thick, smooth, tension=0.6, postaction={decorate},
 decoration={markings, mark=at position 0.2 with {\arrow{>}},
 mark=at position 0.7 with {\arrow{>}}}]
 (0,0) .. controls (0.5,2.7) and (3.5,0.3) .. (2,2);
 \draw[red, thick, smooth, tension=0.6 ]
 (2,2) .. controls (1.8,2.2) .. (2,2.4);
 
 \node at (3,2) {\textbf{$(x(t), y(t))$}}; 
 \node at (-0.2,0.15) {\textbf{$0$}};
 \node at (2,1.99) {{$\Large\cdot$}};

 \node at (2,1.6) {\textbf{$+$}};
 
 \node at (0.5,0.8) {\textbf{$-$}};
\end{tikzpicture}
	\caption{The lift of the curve is performed by defining the third coordinate $z(t)$ as the oriented area of the region between the arc of the curve up to the point $(x(t),y(t))$ and the straight segment from $(0,0)$ to $(x(t),y(t))$.}
	\label{orient-area}
\end{figure}
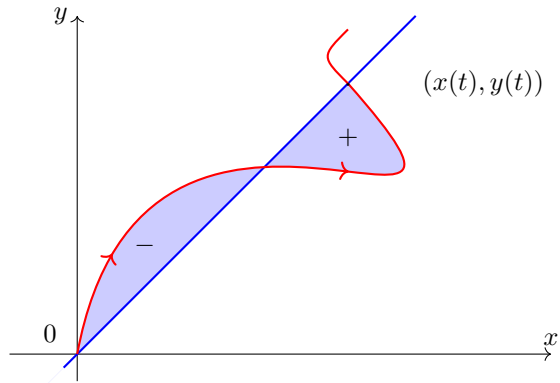

Observe that at each point $(x,y)\in\R^2$, the vector $(x,y)$ is in the kernel of $\alpha$.
Consequently, if $L$ is a line passing through the origin, we have that $\int_L\alpha=0$.
This observation leads us to the conclusion that for a smooth curve $\gamma$, starting from the origin and not necessarily closed, the integral $\int_\gamma \alpha$ represents the signed area enclosed by $\gamma$ and the line segment connecting the origin to the final point of $\gamma$. Refer to Figure~\ref{orient-area} for a visual representation. 

Therefore, Dido's problem can be reformulated as the task of maximizing the integral 
$\int_\gamma\alpha$ 
while fixing the integral 
$ \int_\gamma \dd s$,
 which expresses the length of the curve obtained by integrating it with respect to the element of arc length $\dd s$.


\section{Contact-geometry formulation of the problem} \label{sec:contact}
One of the models of the Heisenberg geometry is constructed as follows, and it has the property that the projection
$\pi:\R^3\to\R^2$ onto the first two coordinates sends geodesics into the solutions of Dido's isoperimetric problem.

If we begin with a curve $\sigma(t)=(x(t),y(t))$ in $\R^2$, with $x(0)=y(0)=0$, we can lift it to a curve in 3D space, where the third coordinate $z(t)$ is the signed area enclosed by the arc $ \sigma_{[0,t]}$, obtained by restricting $\sigma$ to the interval $[0,t]$, and the segment connecting $0$ to $(x(t),y(t))$, as in Figure~\ref{orient-area}. 
Namely, we have
\begin{equation}
\label{third}z(t) :=\int_{\sigma_{[0,t]}}\alpha=\int_{\sigma_{[0,t]}}\dfrac{1}{2} (x\dd y-y\dd x)
=\int_0^t \dfrac{1}{2}\left( x(s) \dot y(s) - y(s) \dot x(s) \right) \dd s.
\end{equation}
Differentiating in $t$ we get
\begin{equation}
 \label{contact1a}
\dot z=\dfrac{1}{2} (x\dot y-y\dot x).
\end{equation}
Set $\xi:=\dd z-\dfrac{1}{2} (x \dd y-y \dd x)$.
Consider a curve $\gamma=(\gamma_1, \gamma_2, \gamma_3):[0,1]\to\R^3$ starting at $0$. 
Then we have that such lifted curves are exactly those satisfying $\dot\gamma\in\ker(\xi)$, which is $\xi(\dot\gamma_1, \dot\gamma_2, \dot\gamma_3)\equiv 0$.


The differential one-form $\xi$ can be written in cylindrical coordinates $(r, \theta, z)$ as $\dd z-\dfrac{1}{2}r^2\dd\theta$.
\begin{definition}[Standard contact form]
We refer to the differential one-form
\begin{equation}
 \label{contact1}\xi:=\dd z-\dfrac{1}{2} (x\dd y-y \dd x)=\dd z-\dfrac{1}{2}r^2\dd \theta\end{equation} 
as the {\em standard contact form} in $\R^3$.
More generally, a {\em contact form} on a $(2 n + 1)$-dimensional differentiable manifold $M$ is a differential $1$-form $\alpha$, with the property that
$$
 \alpha \wedge (\dd\alpha)^n \ne 0, \, 
\quad\text{ 
where
}\quad
 (\dd\alpha)^n := \underbrace {\dd\alpha \wedge \dots \wedge \dd\alpha}_n.\, 
$$
Sometimes, the contact forms
$
\dd z-x \dd y+y \dd x=\dd z- r^2 \dd\theta
$ and $\dd z+x \dd y$ are also called {\em standard}.\index{contact! form}\index{standard contact form} 
\end{definition}


 \begin{figure}
\centering
 \includegraphics[width=5in]{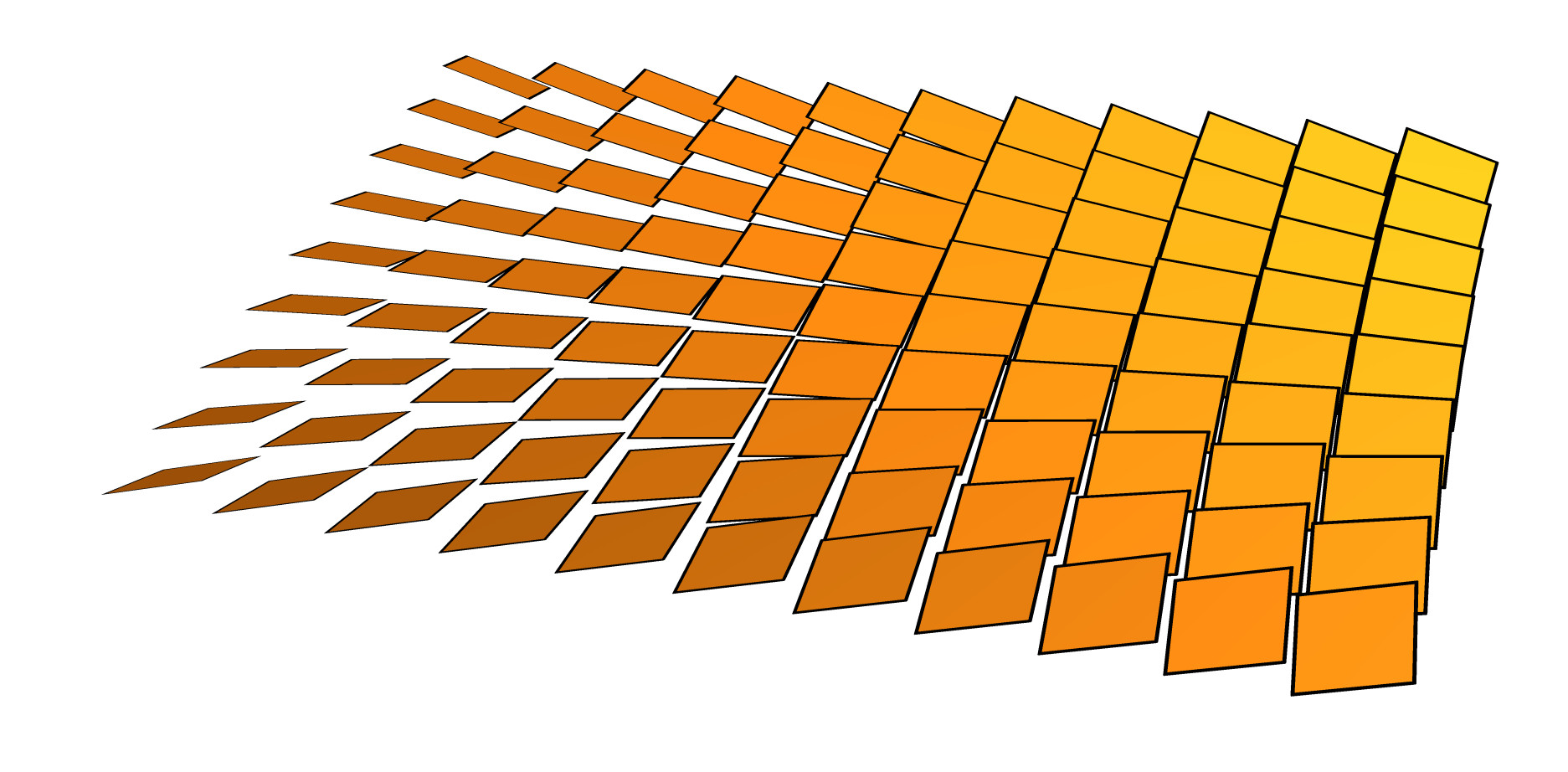} 
\caption{Standard contact distribution on $\R^3.$}
\label{Heisenberg_distribution_Large}
\end{figure}

As with every never-vanishing differential one-form on $\R^3$, the standard contact form \eqref{contact1} gives, at each point $(x,y,z)\in\R^3$, a 2D kernel inside the tangent space $T_{(x,y,z)}\R^3\cong\R^3$ at $(x,y,z)$:
\begin{equation}\label{first_def_Delta}
\Delta_{(x,y,z)}:=\ker (\xi_{(x,y,z)})=\left\{(v_1,v_2,v_3)\in\R^3 \;:\;v_3=\dfrac{1}{2} (xv_2-y v_1)\right\}.
\end{equation}
Geometrically, the set $\Delta$ forms a field of 2D planes in the 3D space, known as {\em distribution of planes}, or {\em contact distribution}.\index{contact! -- distribution}
Now, given vectors $v = (v_1,v_2,v_3)$ and $w = (w_1,w_2,w_3)$, we consider the linear product given by
\begin{equation}
 \langle v,w\rangle := v_1w_1+v_2w_2. \label{SC} 
\end{equation}

Notice that, since each plane $\Delta_{(x,y,z)}$ never includes the $z$-axis, then the restriction of $ \langle \cdot, \cdot\rangle$ on $\Delta_{(x,y,z)}$ is a positive-definite inner product. 
Alternatively, one can view this restriction as a restriction of a Riemannian tensor on $\R^3$, i.e., a positive-definite inner product on the entire tangent bundle of $\R^3$.
In fact, we can consider the following frame\footnote{A {\em frame} is a set of vector fields on a differentiable manifold $M$ that, at each point $p\in M$, gives a basis for the tangent space $T_p M$.\index{frame}} of $\R^3$:
\begin{equation}
 \left\{\begin{array}{ccl}
X&:=&\frac{\partial}{\partial x} - \frac{1}{2} y\frac{\partial}{\partial z}, \\
 Y&:=&\frac{\partial}{\partial y} + \frac{1}{2} x\frac{\partial}{\partial z}, \\ \label{ONF}
 Z&:=&\frac{\partial}{\partial z}, \end{array}
\right. \end{equation}
and declare it orthonormal.
Let us verify that such a Riemannian metric gives the linear product (\ref{SC}) when restricted to the plane $\Delta_{(x,y,z)}$. Indeed, since $\frac{\partial}{\partial x}= X+\frac{1}{2} y Z$ and $\frac{\partial}{\partial y}= Y-\frac{1}{2} x Z$, then we can write
$$v=v_1 X+ v_2 Y + \left(\frac{v_1}{2}y-\frac{v_2}{2}x+v_3\right) Z.$$
So, if $v\in \Delta_{(x,y,z)}$, we have $v=v_1 X+ v_2 Y$ and thus (\ref{SC}) holds.

In contact geometry, a curve $\gamma$ is called {\em Legendrian}\index{Legendrian} with respect to the differential 1-form $\xi$ if $\xi(\dot\gamma)\equiv 0$. In other words, the tangent vector $\dot\gamma(t)$ must be in the plane $\Delta_{\gamma(t)}$, as defined in \eqref{first_def_Delta}.
For a given Legendrian curve $\gamma$, we define its length $L(\gamma)$ as the integral of the norm of $\dot\gamma$ using the scalar product (\ref{SC}). In simpler terms, the value $L(\gamma)$ is precisely the Euclidean length of the projection of $\gamma$ onto the first two components of $\R^3$.

At this point, we introduce a new distance on $\R^3$ to which we refer as the {\em contact distance}. For every pair of points $p$ and $q$ in $\R^3$, we define it as follows:
\begin{equation}
\label{contactdist}
d_c(p,q):=\inf\{ L(\gamma)\;:\; \gamma \text{ is a Legendrian curve between } p \text{ and } q\}.
\end{equation}

An important fact to note is that because $\xi$ was derived from Dido's problem, for every pair of points in $\R^3$, there are multiple Legendrian curves connecting them. In other words, it is always possible to find a Legendrian curve between any two points. Let's underline this crucial point:
\\
{\bf{A crucial fact:}} Every pair of points in $\R^3$ is connected by a curve that is Legendrian with respect to the differential form $\xi$ defined in \eqref{contact1}.
\\
In practice, to connect, for instance, $(0,0,0)$ to $(x,y,z)$, it is enough to take a curve $\sigma$ in $\R^2$ from $(0,0)$ to $(x,y)$ with the property that the signed area enclosed by $\sigma$ and the segment from $(0,0)$ to $(x,y)$ is exactly $z$.
Then, the lifted curve $\tilde\sigma$ will connect $(0,0,0)$ to $(x,y,z)$.

Furthermore, it's important to stress that the length of $\tilde\sigma$ is equal to the planar Euclidean length of $\sigma$. Consequently, there exists a correspondence between geodesics concerning the metric $d_c$ (or, more accurately, curves that realize the infimum in \eqref{contactdist}) and solutions to the {\em dual} Dido's isoperimetric problem:
Fixed a value for the area, minimize the perimeter.
Since it is relatively straightforward to find solutions to Dido's problem, we will be able to explicitly determine the geodesics of the metric space $(\R^3,d_c)$. We will delve into this topic further in Section~\ref{geod-Heis}.

 \begin{figure}
 \centering
 \includegraphics[height=1.4in]{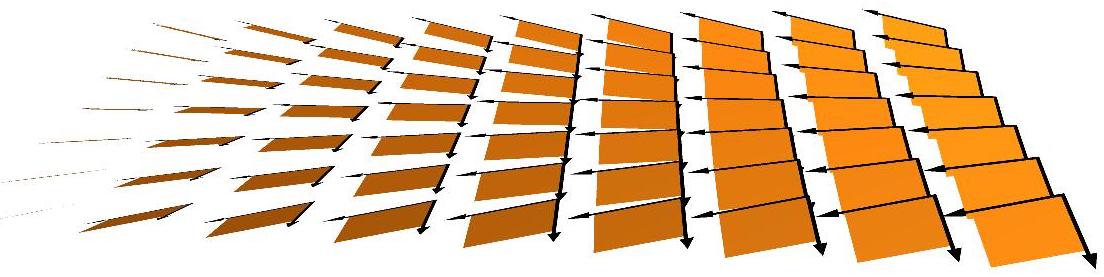}
	\caption{The horizontal bundle spanned by the vector fields $X$ and $Y$.}
	\label{frame_field_arrows}
\end{figure}

\section{The Heisenberg group}\label{sec:Heisenberg}\index{Heisenberg! -- group} 
\subsection{Heisenberg-group invariance of the standard contact structure}
At this point, we have introduced a geometry, which we will call {\em contact geometry}. Specifically, we are considering the plane distribution that at every point $(x,y,z)\in \R^3$ is spanned by the vectors: 
\begin{eqnarray}\label{vector_field_heisenberg}
X(x,y,z)&:=&\frac{\partial}{\partial x} - \frac{y}{2} \frac{\partial}{\partial z} = (1, 0, - \frac{y}{2} ), \\ \nonumber
 Y(x,y,z)&:=&\frac{\partial}{\partial y} + \frac{x}{2} \frac{\partial}{\partial z} = ( 0, 1, \frac{x}{2} ) ;\
\end{eqnarray}
at each point $(x,y,z)$ we are considering $X(x,y,z)$ and $ Y(x,y,z)$ as an orthonormal basis on their span $\Delta(z,y,z)$;
 for each smooth curve $\gamma:[a,b] \to \R^3$ for which $\dot \gamma(t)$ is in $\Delta(\gamma(t))$ we define its length.
 Namely, if $u_1(t), u_2(t)$ are such that 
 $\dot\gamma(t) = u_1(t) X_{\gamma(t)} + u_2(t) Y_{\gamma(t)}$, then the length of $\gamma $ is defined as $
\int_a^b \sqrt{ u_1(t)^2 + u_2(t)^2} \dd t.
$ Such a length structure defines the contact distance \eqref{contactdist}.

A crucial property of the contact geometry is that the space is metrically homogeneous. In fact, the space $\R^3$ can be endowed with a group structure (different from the Euclidean one) in such a way that all of the above constructions are preserved by the action of the group onto itself.

This group structure is named after Werner Heisenberg (1901--1976), a German theoretical physicist and pioneer of quantum mechanics. The Heisenberg group provides a geometric framework to represent Heisenberg's famous uncertainty principle. In our coordinates of $\R^3$, the group product of this structure is defined as follows:
\begin{equation}\label{Hprod}
(x,y,z)\cdot (x',y',z'):=\left(x+x',y+y',z+z'+\dfrac{1}{2} (xy'-yx') \right).
\end{equation}
One can easily check that \eqref{Hprod} gives a group structure, and it turns $\R^3$ into a Lie group, i.e., multiplication and inversion are smooth maps. We will revisit the general theory of Lie groups in Chapter~\ref{ch_LieGroups}. We shall refer to the group $\R^3$ equipped with group law \eqref{Hprod} as the {\em Heisenberg group}.

We claim that the left translations preserve the distribution $\Delta$ and, in fact, preserve the orthonormal frame $X, Y, Z$ defined by (\ref{ONF}).
Let us verify this claim for $X$.
 Fix a left translation $f$, say $f= L_{(s,t,u)}$, for $ {(s,t,u)}\in\R^3$, i.e., 
\begin{equation}
\label{left-trans}f(x,y,z):=L_{(s,t,u)}(x,y,z):=(s,t,u)\cdot (x,y,z)\stackrel{\eqref{Hprod}}{=}\left(x+s,y+t,z+u+\dfrac{1}{2} (sy-tx) \right).
\end{equation}
The differential is 
\begin{equation}\label{detLT}\dd f=\begin{bmatrix}
	1&0&0\\
	0&1&0\\
	-t/2&s/2&1
 \end{bmatrix}.\end{equation}
So, on the one hand, we have that $\dd f (X)$ is given by:
 $$\dd f (X)
=\begin{bmatrix}
	1&0&0\\
	0&1&0\\
	-t/2&s/2&1
 \end{bmatrix}\begin{bmatrix}
	1 \\
	0 \\
	-y/2 
 \end{bmatrix}
= \begin{bmatrix}
	1 \\
	0 \\
	-t/2-y/2 
 \end{bmatrix}
= \frac{\partial}{\partial x} + \left(-\frac{t}{2} -\frac{y}{2} \right) \frac{\partial}{\partial z}.$$ On the other hand, we have $X\circ f= \frac{\partial}{\partial x} -\frac{1}{2} \left(t+y \right) \frac{\partial}{\partial z}.$
Therefore $\dd f (X)=X\circ f$, i.e., $X$ is left-invariant.
Analogously, $\dd f (Y)= \frac{\partial}{\partial y} +\frac{1}{2} \left(s+x \right) \frac{\partial}{\partial z} =Y\circ f$
and $\dd f (Z)= \frac{\partial}{\partial z} =Z\circ f$.

As a consequence of the fact that each left translation by the product \eqref{Hprod} preserves the orthonormal frame $\{X, Y\}$ we deduce that each of these translations preserves the length of Legendrian curves and, consequently, preserves the contact distance as defined in \eqref{contactdist}.

The next proposition summarizes the above discussion.
\begin{proposition}
 The Heisenberg geometry is metrically homogeneous: the space has a Lie group structure in which each left translation acts as an isometry with respect to the contact distance $d_c$. 
\end{proposition}
%
The above model of the Heisenberg group offers the advantage that its one-dimensional subgroups are easily computable and visually understandable. Specifically, the one-parameter subgroups of this group structure correspond to the standard Euclidean lines passing through the origin:
$$t\in\R\longmapsto\gamma_v(t)=\exp \left( t(v_1, v_2,v_3)\right)=\left(t v_1, t v_2,t v_3\right), \qquad \text{ for } (v_1, v_2,v_3)\in \R^3.$$
Furthermore, it is worth noting that all the lines through $0$ in the $x y$-plane are curves that minimize the contact distance. This last statement is recommended as an exercise for novice readers.

\subsection{The 3D nilpotent non-abelian matrix group}
The Heisenberg group can also be represented using matrices. It is a subgroup of the group of invertible matrices and can be defined as the set of $3\times3$ upper-triangular matrices equipped with the row-by-column matrix product:
\begin{equation}\label{Heis_as_matrixes} \mathbb G :=\left\{ \begin{bmatrix} 1 & a & c\\ 0 & 1 & b\\ 0 & 0 & 1\\ \end{bmatrix} \;:\; a,b,c \in \R\right\} < \GL(3, \R).
\end{equation}

This matrix representation is particularly useful because, first, it simplifies the comprehension of the group structure. Second, it allows the Lie algebra associated with this group to be viewed as a matrix Lie algebra. Furthermore, the exponential of the Lie group corresponds to the classical matrix exponential.

The Lie algebra corresponding to this matrix group $\mathbb G$ is:
$$ \g :=\left\{ \begin{bmatrix} 0 & a & c\\ 0 & 0 & b\\ 0 & 0 & 0\\ \end{bmatrix} \;:\; a,b,c \in \R\right\} .$$
A basis of the Lie algebra is 
\begin{equation} \label{MONF}
 X:=\begin{bmatrix}
 0 & 1 & 0\\ 0 & 0 & 0\\ 0 & 0 & 0
 \end{bmatrix}, \qquad 
Y:=\begin{bmatrix}
 0 & 0 & 0\\ 0 & 0 & 1\\ 0 & 0 & 0
 \end{bmatrix}, \qquad 
Z:=\begin{bmatrix}
 0 & 0 & 1\\ 0 & 0 & 0\\ 0 & 0 & 0
 \end{bmatrix} .\end{equation}
One parameter subgroups, as $a,b,c$ vary in $\R$, are represented as:
\begin{eqnarray*}
t\in \R\longmapsto \gamma_{(a,b,c)}(t)&:=&\exp \left( t \begin{bmatrix} 0 & a & c\\ 0 & 0 & b\\ 0 & 0 & 0\\ \end{bmatrix}\right)\\
	&=& I + t \begin{bmatrix} 0 & a & c\\ 0 & 0 & b\\ 0 & 0 & 0\\ \end{bmatrix}+\frac{ t^2}{2!} \begin{bmatrix} 0 & a & c\\ 0 & 0 & b\\ 0 & 0 & 0\\ \end{bmatrix}^2+\ldots \\
	&=& I + t \begin{bmatrix} 0 & a & c\\ 0 & 0 & b\\ 0 & 0 & 0\\ \end{bmatrix}+\frac{ t^2}{2} \begin{bmatrix} 0 & 0 & ab\\ 0 & 0 & 0\\ 0 & 0 & 0\\ \end{bmatrix}+0 \\
&=& \begin{bmatrix} 1 & at & ct +a b t^2/2\\ 0 & 1 & b t\\ 0 & 0 & 1\\ \end{bmatrix}.
\end{eqnarray*}

We claim that the map 
$$\varphi: (x,y,z) \longmapsto \begin{bmatrix} 1 & x & z+\dfrac{1}{2} x y\\ 0 & 1 & y\\ 0 & 0 & 1\\ \end{bmatrix}
$$
is a Lie group isomorphism from the Lie group $\R^3$ with the product (\ref{Hprod}) to the Lie group $\mathbb G$ from \eqref{Heis_as_matrixes} 
with the matrix product.
Indeed, the map $\varphi$ is a group homomorphism (a result of straightforward calculation), and its differential at the identity maps the left-invariant vector fields $X, Y, Z$ from \eqref{ONF} to $X, Y, Z$ from \eqref{MONF}, respectively. In fact, the subsequent section will elaborate on this concept further.

\subsection{Characterization of the Heisenberg algebra}\index{Heisenberg! -- algebra} 
The Lie algebra of the Heisenberg group is spanned by three vectors, denoted as \(X\), \(Y\), and \(Z\), where the only non-trivial Lie bracket relation is given by \([X, Y]=Z\). Notably, for every three vectors \(X_1\), \(X_2\), and \(X_3\) in this Lie algebra, the Lie bracket operation satisfies the property \([X_1, [X_2, X_3]]=0\). This characteristic feature establishes the Heisenberg Lie algebra as a nilpotent Lie algebra with a step of \(2\). To clarify, a Lie algebra is said to be {\em nilpotent} with a {\em nilpotency step} \(s\) if, for every selection of more than \(s\) vectors within it, the iterated bracket of these vectors results in \(0\).

In every vector space, we can define the zero bracket operation, thus forming a Lie algebra. Such a Lie algebra is referred to as {\em commutative} and is nilpotent with a step of 1.

We assert that there are only two examples of 3D simply connected nilpotent Lie groups: the 3D vector space \((\R^3,+)\) and the Heisenberg group.
To prove this, let \(\mathfrak{g}\) be the Lie algebra of one such group. Since \(\mathfrak{g}\) is nilpotent, we can take a non-zero element \(Z\) in the center of \(\mathfrak{g}\). We extend \(Z\) to form a basis \(\{X, Y, Z\}\) for \(\mathfrak{g}\). Now, either \(X\) and \(Y\) commute, making the algebra commutative, or the vector \(W := [X, Y] \) is not 0. In this second case, we write \(W = aX + bY + cZ\), for some $a,b,c\in \R$. Then, we have that \([W, Y] = aW\). Due to the nilpotency of \(\mathfrak{g}\), we have \(a = 0\). Similarly, \(b = 0\). Thus, we infer that \(c \neq 0\), and by replacing \(Z\) with \(cZ\), we can define the algebra structure of \(\mathfrak{g}\) by the relations:
\[ [X, Y] = Z \quad\text{and}\quad [X, Z] = [Y, Z] = 0. \]
We can conclude the proof by invoking the uniqueness of a simply connected Lie group with a fixed Lie algebra, as stated in Theorem~\ref{Lie's Third Theorem}.

\section{The sub-Riemannian Heisenberg group}\label{sec:subRiem_Heis}

Our preferred model for the Heisenberg group is $\R^3$ with the product law \eqref{Hprod}, which we observed yield the following left-invariant vector fields: $ \partial_x - \frac y2\partial_z $, $ \partial_y + \frac x2 \partial_z$, $\partial_z$.
The reason why this model is advantageous lies in its canonical identification of the group with its Lie algebra. In other words, we are using exponential coordinates, a perspective that will be elucidated in Section~\ref{sec1120}.
It's important to note that, given the uniqueness of the Heisenberg structure, all other models can be considered equivalent through a smooth group morphism.

In $\R^3$, we begin by selecting three vector fields, denoted as $X$, $Y$, and $Z$, that as for \eqref{vector_field_heisenberg} are linearly independent at every point and satisfy the following commutation relations:
\[ [X, Y] = Z \quad\text{and}\quad [X,Z] = [Y,Z] = 0. \]
It is a fact that the space can be endowed with a group law that renders them left-invariant.

We consider the subbundle $\Delta$ of the tangent bundle $ T(\R^3)$ over $\R^3$ such that for every $p\in\R^3$
\[
\Delta_p := \Delta\cap T_p\R^3:= \Span\{X_p,Y_p\} .
\] 
A smooth (or, more generally, absolutely continuous) curve $\gamma:[0,1]\to \R^3$ such that $\dot \gamma\in\Delta$ is called {\em horizontal curve}. 
\index{horizontal! -- curve}
In this case, if we express $\dot\gamma(t) = u_1(t) X_{\gamma(t)} + u_2(t) Y_{\gamma(t)}$, for almost all $t\in[0,1]$, for some integrable functions $u_1$ and $ u_2$ on $[0,1]$, then the {\em length} of $\gamma$ is defined as
\[
L(\gamma) := \int \sqrt{ u_1(t)^2 + u_2(t)^2} \dd t.
\]
We define the \emph{Carnot-Carath\'eodory distance} 
\index{Carnot-Carath\'eodory! -- distance}
between two points $p$ and $q$ in $\R^3$ as
\begin{equation}\label{first_def_dcc}
\dcc(p,q) := \inf\left\{ L(\gamma) : \gamma\text{ is a horizontal curve from $p$ to $q$}\right\}.
\end{equation}
Therefore, we have broadened the notion of Legendrian curve with the one of horizontal curve, and the concept of a contact distance has been expanded to a Carnot-Carath\'eodory distance.
This variation in terminology arises because sub-Riemannian geometry draws from various mathematical domains, resulting in multiple jargon.
 
We refer to the space $(\R^3, \dcc)$ as {\em(a model for) the sub-Riemannian Heisenberg group}. \index{sub-Riemannian!-- Heisenberg group}
In the rest of this section, we will exclusively operate within our preferred model: $\R^3$ with the product law \eqref{Hprod} and the orthonormal frame \eqref{vector_field_heisenberg}.

\subsection{Geodesics and spheres in the Heisenberg group} \label{geod-Heis}
From Sections \ref{sec:contact} and~\ref{sec:Heisenberg}, we can deduce the following properties of a curve $\gamma(t)=(x(t),y(t),z(t))$:
\begin{itemize}
\item The curve \(\gamma\) is horizontal, meaning \(\dot{\gamma} \in \Delta\), if and only if 
\[
\dot{z} = \frac{1}{2}(x\dot{y} - y\dot{x}),
\]
which is equivalent to stating that \(z(t)\) represents the area spanned by the curve \((x(\cdot),y(\cdot))\) up to the point \(t\), as in Figure~\ref{orient-area}.

\item The condition \(\dot{\gamma} \in \Delta\) holds if and only if \(\dot{\gamma} = u_1X + u_2Y\), where \(u_1 = \dot{x}\) and \(u_2 = \dot{y}\). To see this, on the one hand, we have \(\pi(\dot{\gamma}) =\pi((\dot{x}, \dot{y}, \dot{z}))= (\dot{x}, \dot{y})\) and, on the other hand, we have
\[
\pi(\dot{\gamma}) = \pi(u_1X + u_2Y) = u_1\partial_x + u_2\partial_y = (u_1,u_2).
\]

\item If \(\dot{\gamma} \in \Delta\), then the length of \(\gamma\) is given by 
\[
L(\gamma) = \int \sqrt{\dot{x}^2 + \dot{y}^2} = {\rm Length_{Eucl}}(\pi\circ\gamma).
\]
\end{itemize}
%
Due to the preceding discussion, we can derive explicit formulas for the geodesics in the sub-Riemannian Heisenberg group. This is made possible by our understanding of the solutions to the isoperimetric problem, as detailed in Section~\ref{sec:isoper}. Specifically, we have discovered that given the way the geometry in the Heisenberg group has been constructed, the shortest curves concerning the length structure are the lifts of solutions to a variant of the isoperimetric problem.
Namely, we search for the shortest curves in the plane that enclose a fixed area and join two specified points. Such curves turn out to be arcs of circles or straight-line segments. Consequently, the geodesics in the Heisenberg group are essentially the lifted versions of circles in the plane.
%
\begin{fact}
For a fixed endpoint $(x(1),y(1),z(1))$, every curve $(x(t),y(t))$ that encloses an area equal to $z(1)$, is such that $(x(0),y(0))=(0,0)$, and among such curves minimizes ${\rm Length_{Eucl}}(x(\cdot),y(\cdot))$ is either a segment of a circle or a straight line in the plane.

Hence, length-minimizing curves originating from $(0,0,0)$ are lifts of circles if $z(1)\neq0$ and straight lines if $z(1)=0$.
\end{fact}

We aim to parametrize the curves that solve Dido's problem.
A circle of length $\frac{2\pi}{|k|}$, with $k\neq0$, passing through $(0,0)$ at time $0$ is described as
\[
(x_0(t),y_0(t)) = 
\left( \frac{\cos(kt)-1}{k}, \frac{\sin(kt)}{k} \right)
\]
for $0\le t\le \frac{2\pi}{|k|}$. Such a circle is parametrization by arc length and has its center on the $x$-axis.
If $k>0$, the center is on the negative $x$-axis; if $k<0$, it is on the positive $x$-axis.
Moreover, if $k>0$, then the circle $(x_0,y_0)$ encloses positive area; if $k<0$ it encloses negative area. For $k=0$, we can still consider the formula in the limit sense: the circles degenerate into the line $(0,t)$, defined for all $t\in\R$.

We obtain every other circle by rotating these initial ones by an angle $\theta\in\R/2\pi\Z$:
\[
R_\theta(x_0(t),y_0(t)) 
:= \begin{bmatrix}
 	\cos\theta & -\sin\theta \\
	\sin\theta & \cos\theta
\end{bmatrix}
\cdot
\begin{bmatrix}
 	\frac{\cos(kt)-1}{k} \\ \frac{\sin(kt)}{k}
\end{bmatrix}
= \begin{bmatrix}
 	\cos\theta \frac{\cos(kt)-1}{k} - \sin\theta \frac{\sin(kt)}{k} \\
	\sin\theta \frac{\cos(kt)-1}{k} + \cos\theta \frac{\sin(kt)}{k}
\end{bmatrix}.
\]

We can calculate the third coordinate of the lift of such circles via (\ref{third}).
\begin{eqnarray*}
z(T) &=& \int_{0}^T\dfrac{1}{2} (x \dd y-y \dd x )= \frac12\int_{0}^T x\dot y-y\dot x
\\
 &=&\dfrac{1}{2}\int_{0}^T 
\left(\cos\theta\dfrac{ \cos(k t)-1}{k}-\sin\theta\dfrac{ \sin(k t) }{k} \right)
\left(-\sin\theta\sin(k t)+\cos\theta \cos(k t) \right)\\
&&\qquad- \left(\sin\theta\dfrac{ \cos(k t)-1}{k}+\cos\theta\dfrac{ \sin(k t) }{k}\right)\left(-\cos\theta\sin(kt)-\sin\theta\cos(k t) \right)\dd t\\
&=&\dfrac{1}{2 k}\int_{0}^T 
-\cos\theta(\cos(k t)-1)\sin\theta\sin(k t)+(\cos\theta)^2(\cos(k t)-1) \cos(k t)\\
&&\qquad\qquad+(\sin\theta)^2 \left(\sin(k t)\right)^2-\sin\theta \sin(k t)\cos\theta \cos(k t)
\\
&&\qquad\qquad
+\sin\theta(\cos(k t )-1)\cos\theta\sin(k t)+(\sin\theta)^2(\cos(k t)-1)\cos(k t) \\
&&\qquad\qquad+(\cos\theta)^2 \left( \sin(k t)\right)^2+\cos\theta \sin(k t)\sin\theta\cos(k t)\dd t\\
&=&\dfrac{1}{2 k}\int_{0}^T 
(\cos(kt)-1) \cos(k t)+\left( \sin(k t)\right) ^2\dd t\\
&=&\dfrac{1}{2 k}\int_{0}^T 
 1-\cos(k t)\dd t=\dfrac{1}{2 k^2}(T k- 
 \sin(k T)).
\end{eqnarray*}


 \begin{figure}
\centering
\subfigure[The top view] 
{
 \includegraphics[height=2.5in]{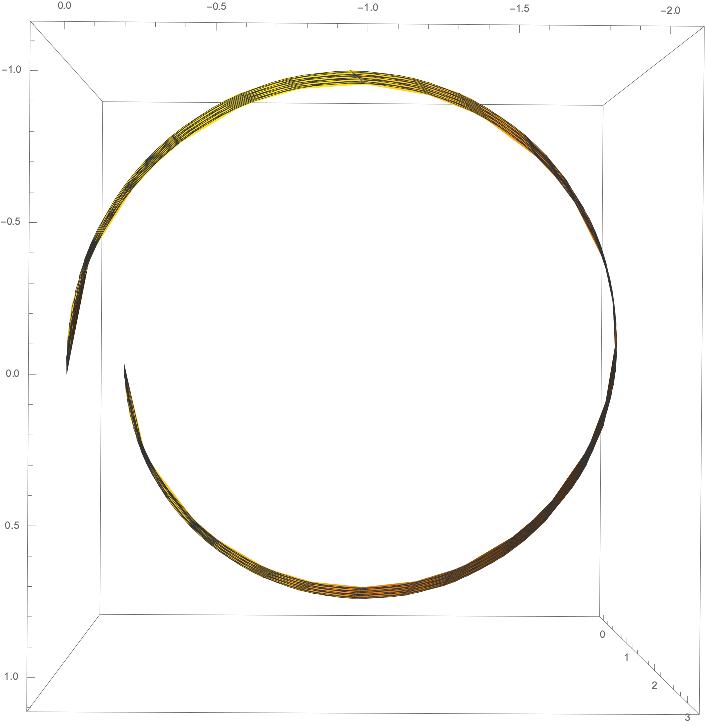}}
 \hspace{0.1cm}
\subfigure[A front view] 
{
 \includegraphics[height=2.5in]{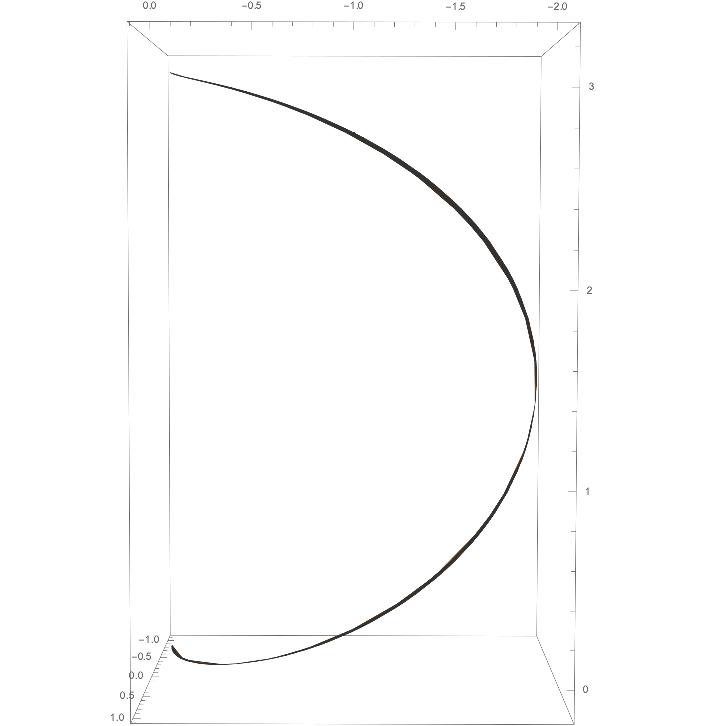}} 
 \hspace{0.1cm}
\subfigure[A side view] 
{
 \includegraphics[height=2.5in]{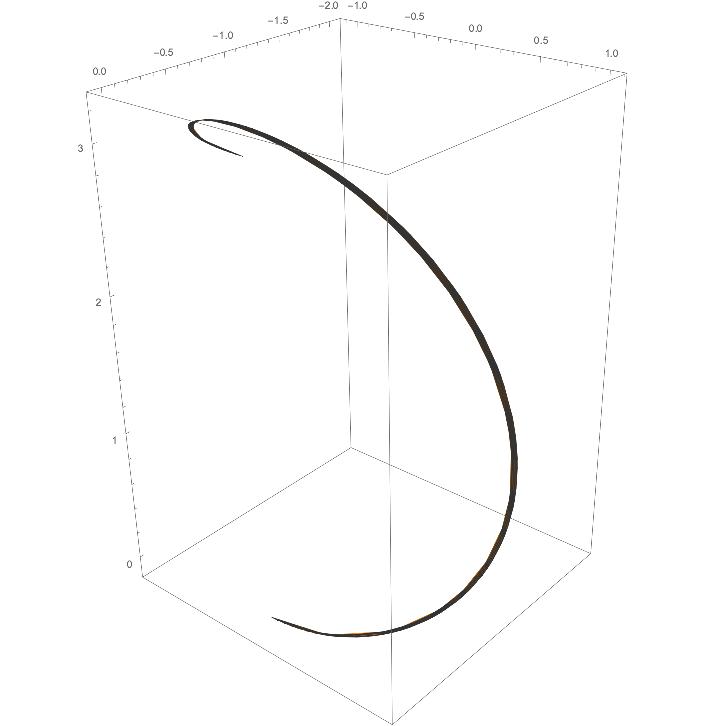}} 
 \hspace{0.1cm}
\subfigure[A side view] 
{
 \includegraphics[height=2.5in]{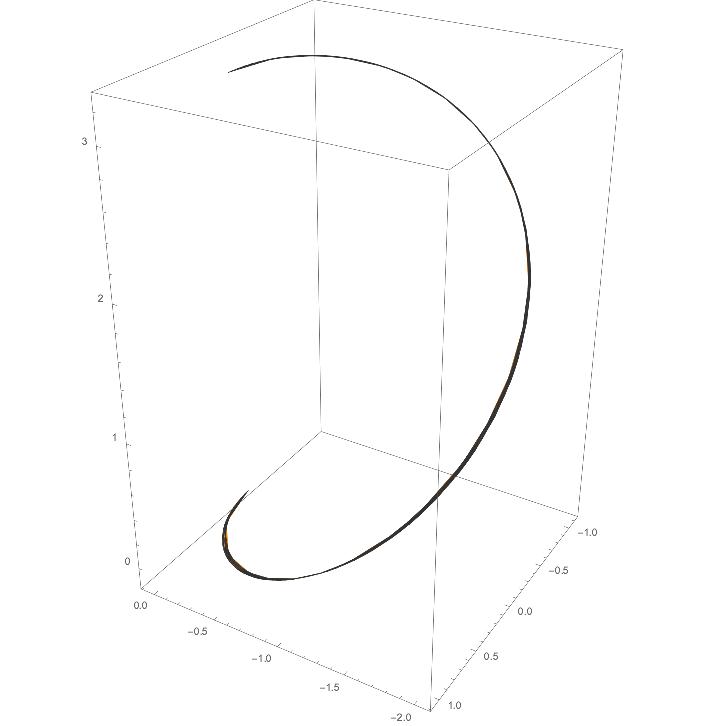}} 
 \caption{Various views of a geodesic with non-zero curvature in the sub-Riemannian Heisenberg geometry}
\end{figure} 

 \begin{figure}	
\centering
\subfigure[A geodesic with zero curvature
] 
{
 \includegraphics[width=4in]{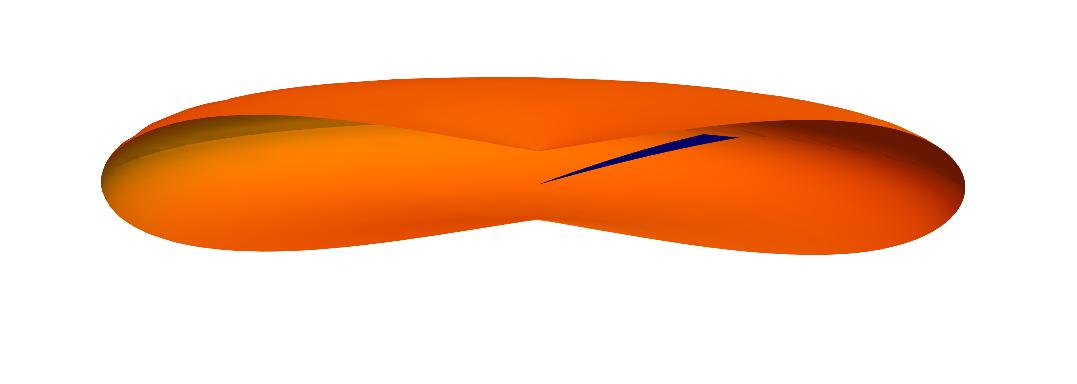}} 
 \hspace{0.1cm}
\subfigure[A geodesic with small curvature
] 
{
 \includegraphics[width=4in]{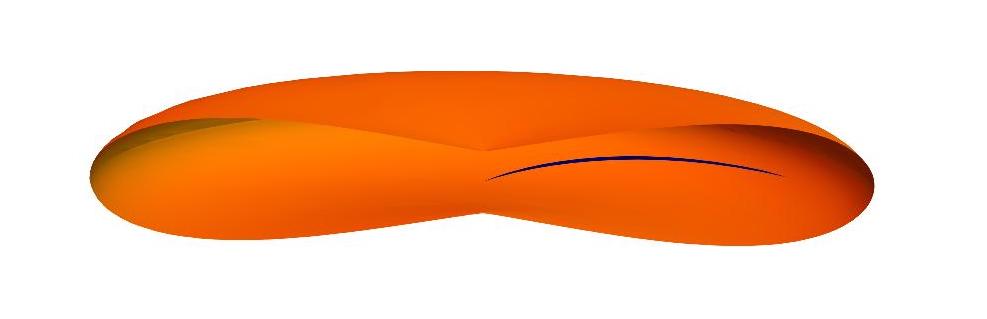}} 
 \hspace{0.1cm}
\subfigure[A geodesic with some curvature less than $\frac{1}{2\pi}$.
] 
{
 \includegraphics[width=4in]{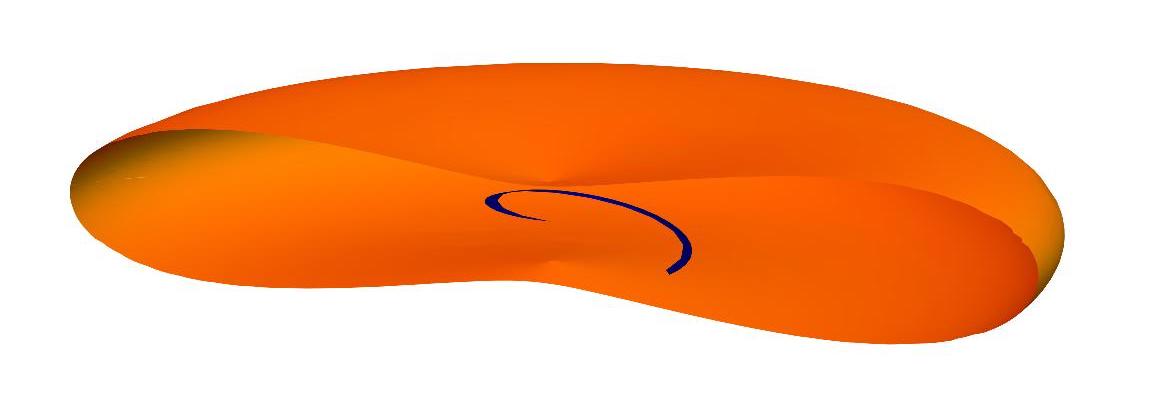}} 
 \hspace{0.1cm}
 \subfigure[A geodesic some curvature equal to $\frac{1}{2\pi}$. It joins points that can be connected with infinitely many geodesics.] 
 {
 \includegraphics[width=4in]{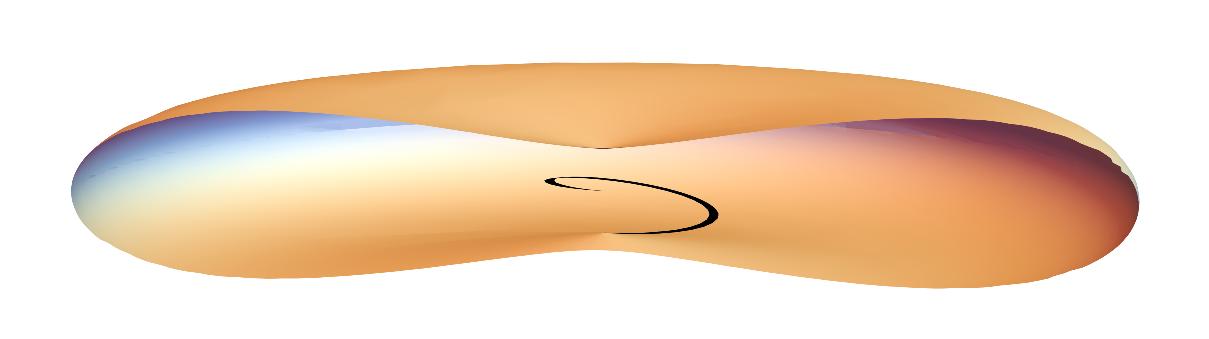}}
 \caption{Geodesics within the unit sphere in the sub-Riemannian Heisenberg geometry.}
\end{figure}


We conclude that the length-minimizing curves originating from the point $(0,0,0)$ in $\R^3$ are smooth curves $\gamma = (\gamma_1, \gamma_2, \gamma_3)$
given by the equations:
\begin{equation}\label{geod-in-H}
\begin{cases}
 \gamma_1(t) = 	\cos\theta \frac{\cos(kt)-1}{k} - \sin\theta \frac{\sin(kt)}{k} \\
 \gamma_2(t) = \sin\theta \frac{\cos(kt)-1}{k} + \cos\theta \frac{\sin(kt)}{k} \\
 \gamma_3(t) = \frac{ kt - \sin(kt) }{2k^2}
\end{cases},
\end{equation}
for some $\theta\in\R/2\pi\Z$ and $k\in\R$.

These curves are defined for $t\in[0, \frac{2\pi}{|k|}]$ and have a length of $\frac{2\pi}{|k|}$.
When $k=0$, these curves degenerate into lines:
\[
\begin{cases}
 	\gamma_1(t) = 	-t\sin\theta \\
 	\gamma_2(t) = t\cos\theta \\
 	\gamma_3(t) = 0
\end{cases}.
\]
Indeed, lines passing through the origin in the $x y$-plane are geodesics.

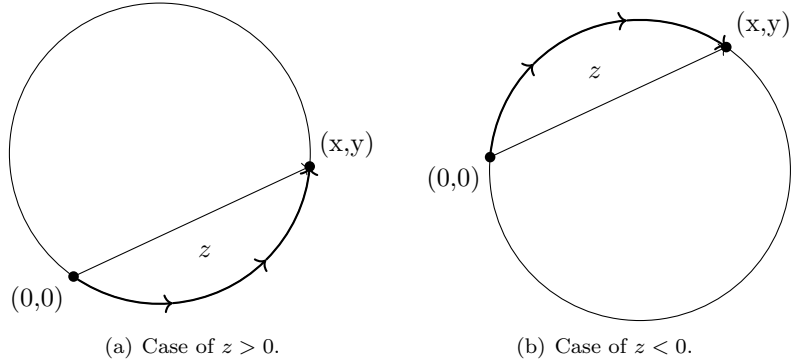
\begin{figure}
\centering
\subfigure[Case of $z>0$.] 
{\begin{tikzpicture}
 \draw (0,0) circle (2cm);
 \node at (0.6cm,-1.3cm) {$z$};

 \draw[thick,->] (-125:2cm) arc (-125:-85:2cm);
 \draw[thick,->] (-85:2cm) arc (-85:-45:2cm); 
 \draw[thick,->] (-45:2cm) arc (-45:-5:2cm); 
 \coordinate (A) at (-125:2cm);
 \coordinate (B) at (-5:2cm);
 
 \fill (A) circle (2pt) node[below left] {(0,0)};
 \fill (B) circle (2pt) node[above right] {(x,y)};
 
 \draw[->] (A) -- (B);
 
\end{tikzpicture}}
 \hspace{0.1cm}
\subfigure[Case of $z<0$.] 
{
\begin{tikzpicture}
 \draw (0,0) circle (2cm);
 \node at (-0.6cm,1.3cm) {$z$};
 
 \draw[thick,->] (175:2cm) arc (175:135:2cm);
 \draw[thick,->] (135:2cm) arc (135:95:2cm); 
 \draw[thick,->] (95:2cm) arc (95:55:2cm); 
 \coordinate (A) at (175:2cm);
 \coordinate (B) at (55:2cm);
 
 \fill (A) circle (2pt) node[below left] {(0,0)};
 \fill (B) circle (2pt) node[above right] {(x,y)};
 
 \draw[->] (A) -- (B);
 
\end{tikzpicture}} 
 
 \caption{Projections of minimizing curves from $(0,0,0)$ to $(x,y,z)$ in the Heisenberg model. When the third coordinate $z$ is positive, the curve follows a circle counterclockwise. If $z$ is negative, it follows clockwise. In both cases, the area enclosed by the curve and the circle equals $|z|$.}
 \label{cerchi_area_z}
\end{figure}

 
We have identified \emph{all} length-minimizing curves in the sub-Riemannian Heisenberg group.
This characterization of the geodesics leads to several interesting conclusions:
\begin{enumerate}
\item 	If a point $(x,y,z)\in\R^3$ lies on the $z$-axis, i.e., $(x,y)=(0,0)$, then there exist infinitely many length-minimizing curves between this point and the origin $(0,0,0)$. 
In fact, these curves form a one-parameter family and are the lifts of circles with area $z$, all containing the point $(0,0)$.
\item 	If $(x,y)\neq(0,0)$, then there exists a unique length-minimizing curve from $(x,y,z)$ to $(0,0,0)$. 
This curve is the lift of a circular arc enclosing area $z$ together with the segment connecting $(0,0)$ to $(x,y)$.
Please refer to Figure~\ref{cerchi_area_z} for a visual representation. 

\end{enumerate}

Since the distance $\dcc$ is left-invariant and also $Z=\partial_z$ is left-invariant, we get that for all $p,q\in\R^3$ there exist infinitely many length-minimizing curves between $p$ and $q$ if and only if $\pi(p)=\pi(q)$, i.e., the points lie on the same vertical line. 
Conversely, if $\pi(p)\neq\pi(q)$, then there is only a single length-minimizing curve between them.

We deduce that this sub-Riemannian geometry differs fundamentally from Riemannian geometry. Although all metric balls and metric spheres in the Heisenberg group remain topological balls and spheres, respectively (see Exercise~\ref{ex: spheres: balls}), this geometry is not Riemannian and, furthermore, it cannot be biLipschitz equivalent to any Riemannian geometry, as it will be shown in Corollary~\ref{corollary_HD_Heis4}.

\begin{figure}
\centering 
\subfigure[The unit sphere has a singularity at the intersection with the $z$-axis.]
 {\includegraphics[width=5in]{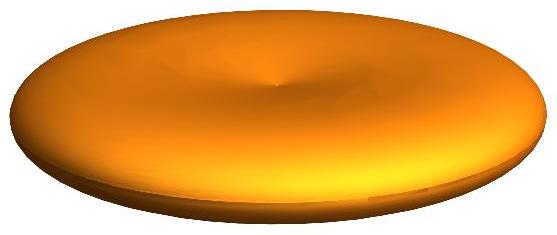}}
 \subfigure[The portion of the unit sphere in the half-space $\{y>0\}$.]
 {\includegraphics[width=5in]{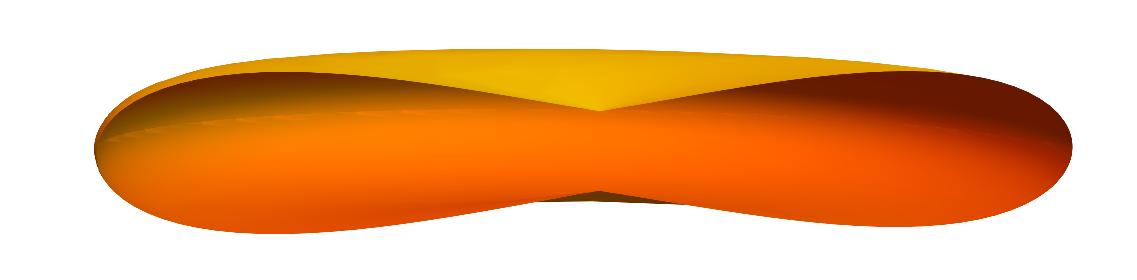}}
 \subfigure[A section of the sphere as intersection with the $xz$-plane.]
{ \includegraphics[width=5in]{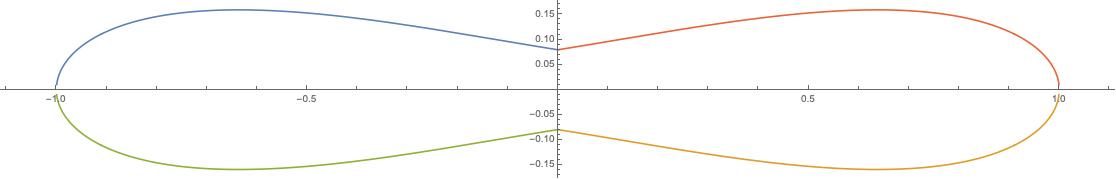}
}
\caption{Balls in the sub-Riemannian Heisenberg group are not smooth surfaces. At the two ``poles'' the sphere is not $C^1$; there is no cusp, but there is a corner. For a parametrization, see Exercise~\ref{ex: spheres: balls}.} 
\end{figure}

\subsection{Dilations on the Heisenberg group}

For every $\lambda\in \R$, we define the map
\begin{equation}\label{def_delta_Heis}
\begin{array}{rccc}
 	\delta_\lambda : &\R^3 &\longrightarrow & \R^3 \\
		& (x,y,z) & \longmapsto &(\lambda x, \lambda y, \lambda^2 z) .
\end{array}
\end{equation}
Notice the squared $\lambda$ in the third component. For $\lambda=0$, such a map is constantly equal to the origin ${\bf{0}}:=(0,0,0)$, which is the identity element for the group law \eqref{Hprod}.

\begin{lemma}\label{Prop_dilations_Heis}
The dilations \eqref{def_delta_Heis} on the Heisenberg group satisfy the following properties: 	For all $\lambda, \mu\in \R$ and all $p,q\in\R^3$:
	\begin{description}
\item[\ref{Prop_dilations_Heis}.i.] 
 	$\delta_\lambda(p\cdot q) = \delta_\lambda(p)\cdot\delta_\lambda(q)$;
\item[\ref{Prop_dilations_Heis}.ii.] 		$\delta_\lambda \circ \delta_\mu = \delta_{\lambda\mu}$;
\item[\ref{Prop_dilations_Heis}.iii.]	$\delta_\lambda$ is a Lie group isomorphism with inverse $\delta_{\frac1\lambda}$, if $\lambda\neq 0$;
	\item[\ref{Prop_dilations_Heis}.iv.] 	$\dcc(\delta_\lambda(p), \delta_\lambda(q)) = |\lambda| \dcc(p,q)$, where $\dcc$ is defined in \eqref{first_def_dcc}.
	\end{description}
\end{lemma}
\proof[Proof of \ref{Prop_dilations_Heis}.i.]
		From the group law~\ref{Hprod}, we have:
\[
\begin{aligned}
\delta_\lambda(p\cdot q) &= \delta_\lambda \left( p_1+q_1, p_2+q_2, p_3+q_3+\frac12(p_1q_2-p_2q_1) \right) \\
&= \left( \lambda p_1+\lambda q_1, \lambda p_2+\lambda q_2, \lambda^2p_3+\lambda^2q_3+\frac12(\lambda p_1 \lambda q_2-\lambda p_2\lambda q_1) \right) \\
&= (\lambda p_1, \lambda p_2, \lambda^2 p_3)\cdot(\lambda q_1, \lambda q_2, \lambda^2 q_3) \\
&= \delta_\lambda(p)\cdot\delta_\lambda(q).
\end{aligned}
\]
	\proof[Proof of \ref{Prop_dilations_Heis}.ii.] 	This is obvious from the definition \eqref{def_delta_Heis}:\\
\[
\begin{aligned}
(\delta_\lambda \circ \delta_\mu) (x,y,z) &= \delta_\lambda(\mu x, \mu y, \mu^2 z) \\
&= (\lambda \mu x, \lambda \mu y, \lambda^2 \mu^2 z) \\
&= (\lambda \mu x, \lambda \mu y,(\lambda \mu)^2 z) \\
&= \delta_{\lambda\mu}(x,y,z).
\end{aligned}
\]	
	\proof[Proof of \ref{Prop_dilations_Heis}.iii.] From the previous points, we conclude that each $\delta_\lambda$ is a group homomorphism and $(\delta_\lambda)^{-1}=\delta_{\frac1\lambda}$, if $\lambda\neq 0$. 
Moreover, it is evident that each map is smooth.

\proof[Proof of \ref{Prop_dilations_Heis}.iv.]	Regarding the last point, we shall give three methods of proof for educational reasons.
	\begin{description}
	\item[Method 1]
		We claim that the map $\delta_\lambda$ is such that $(\delta_\lambda)_*X = \lambda X$ and $(\delta_\lambda)_*Y=\lambda Y$, where $X, Y$ are the vector fields defining the subbundle $\Delta$. (This last statement is suggested as an exercise)
		Hence, the map $\delta_\lambda$ preserves horizontal curves and multiplies their length by $\lambda$.
	\item[Method 2]
		By the left-invariance of $\dcc$ and by (\ref{Prop_dilations_Heis}.ii), we have 
		\[
		\dcc(\delta_\lambda(p), \delta_\lambda(q)) 
		= \dcc((\delta_\lambda(p))^{-1}\cdot\delta_\lambda(q),{\bf{0}})
		= \dcc(\delta_\lambda(p^{-1}q),{\bf{0}}) .
		\]
		Hence, it is enough to show that 
		\begin{equation}\label{eq1231}
		 	\dcc(\delta_\lambda(p),{\bf{0}}) = \lambda\dcc(p,{\bf{0}}) .
		\end{equation}
		
		Let $\gamma$ be a length-minimizing curve from ${\bf{0}}$ to an arbitrary $p$.
		Recall that we have an explicit formula for such curves.
		An easy calculation (see Exercise~\ref{dilation_of_geodesics}) shows that $\delta_\lambda\circ\gamma$ is still of the same form, up to a linear reparametrization by $\lambda$.
		Hence, its length got multiplied by the factor $\lambda$.
	\item[Method 3]
		Reasoning as at the beginning of Method 2, proving \eqref{eq1231} is enough.
		Take a horizontal curve $\gamma=(x,y,z)$ from ${\bf{0}}$ to $p$.
		Notice that the linear map of $\R^2$ represented by the matrix $\begin{bmatrix} \lambda & 0 \\ 0 & \lambda \end{bmatrix}$ multiplies length by $\lambda$ and area by $\lambda^2$. 
		Therefore, the curve $(\lambda x, \lambda y)$ spans areas that are $\lambda^2$ times the areas of $(x,y)$ and has length $\lambda$ times the length of $(x,y)$.
		Thus the curve $(\lambda x, \lambda y, \lambda^2 z)$ is horizontal and has length $\lambda L(\gamma)$.
		Hence $\dcc(\delta_\lambda(p),{\bf{0}}) \le \lambda\dcc(p,{\bf{0}})$.
		
		We conclude by arguing similarly with each curve $\sigma$ joining $\delta_\lambda(p)$ to ${\bf{0}}$ and considering the curve $\delta_{\frac1\lambda}\circ\sigma$.
		\qedhere
	\end{description}
 
As a consequence of Lemma~\ref{Prop_dilations_Heis}, we have the following properties for the metric balls in the Heisenberg group.
\begin{corollary}\label{cor_dil_Heis}
 	In the sub-Riemannian Heisenberg group, we have
 		\begin{description}
	\item[\ref{cor_dil_Heis}.i.] 	$B_{\dcc}({\bf{0}},r) = \delta_r(B_{\dcc}({\bf{0}},1))$;
	\item[\ref{cor_dil_Heis}.ii.] $B_{\dcc}(p,r) = L_p(\delta_r(B_{\dcc}({\bf{0}},1)))$,
	\end{description}
	 for all points $p $ and all radii $r$. 
\end{corollary}

In other words, we deduce that that if $B_{\dcc}({\bf{0}},r)$ is the ball of center ${\bf{0}}$ and radius $r$, then 
\begin{equation}\label{growH}
(x,y,z)\in B_{\dcc}(({\bf{0}},1)\Longleftrightarrow(r x,r y,r^2 z)\in B_{\dcc}({\bf{0}},r).
\end{equation}
Notice that we did not use the homogeneous dilation ${{\bf{v}}}\mapsto r{\bf{v}}$; the third coordinate has been multiplied by $r^2$. 
Thus, such a map $(x,y,z)\mapsto (r x,r y,r^2 z)$ multiplies the volume by a factor of $r^4$, and not $r^3$ as the usual Euclidean dilation of factor $r$ does!

We can now deduce the growth of the balls in the Heisenberg geometry. 

\begin{corollary} 
Let $\vol$ be the 3D Lebesgue volume in $\R^3$.
The Heisenberg sub-Riemannian distance $\dcc$ satisfies 
		\begin{equation}\label{vol_Heis}
		\vol(B_{\dcc}(p,r)) = r^4 \vol(B_{\dcc}({{\bf0}},1)),
		\qquad\forall p\in\R^3, \ \forall r>0 .
		\end{equation}
\end{corollary}

\begin{proof}
From (\ref{growH}), we know that $\vol(B({\bf{0}},r))=r^4\vol(B({\bf{0}},1))$. Now, we can conclude the proof using the fact that left translations (\ref{left-trans}) in the Heisenberg group are isometries, and they preserve the volume. This last fact can be verified by noticing that the determinant of the differential of a left translation is $1$, as seen in (\ref{detLT}).
	 Namely, every left translation $L_p$ is such that $\dd L_p = \begin{bmatrix} 1 & 0 & 0 \\ * & 1 & 0 \\ * & * & 1 \end{bmatrix} $, and thus 
		$\Jac(L_p) = \det(\dd L_p) = 1$.
		Notice that $\Jac(\delta_\lambda) = \det(\dd \delta_\lambda) = \det\begin{bmatrix}
		 \lambda & 0 & 0 \\ 0 & \lambda & 0 \\ 0 & 0 & \lambda^2 
		\end{bmatrix} = \lambda^4$. 
		Then
$$ 
		 	\vol(B(p,r)) = \vol(L_p(B({\bf{0}},r))) = \vol (B({\bf{0}},r)) 
			= \vol(\delta_r(B({\bf{0}},1))) = r^4 \vol(B({\bf{0}},1) ,
	$$ 
	where we used \ref{cor_dil_Heis}.
\end{proof}


\subsubsection{The dimension of the Heisenberg group}

\begin{figure}
\centering 
 \includegraphics[height=2.5in]{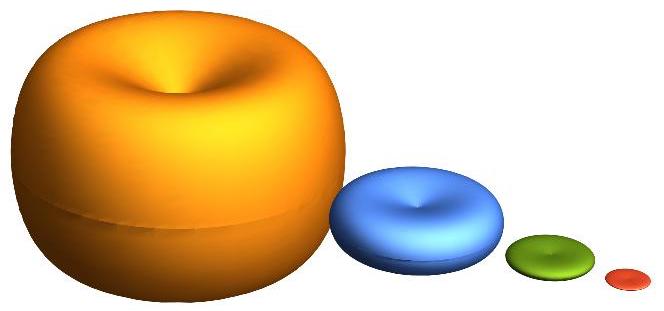}
\caption{Balls of different sizes in the Heisenberg geometry. All the balls are with the origin as the center. From the left, these are the balls of radius 2, 1, 1/2, 1/4. } 
\end{figure} 

The following theorem will demonstrate that the sub-Riemannian Heisenberg group is not locally biLipschitz equivalent to any Riemannian manifold. For the notions of biLipschitz and Hausdorff dimension, please refer to Section~\ref{MetricGeometry}.
\begin{corollary}\label{corollary_HD_Heis4}
 The Heisenberg group endowed with the standard Carnot-Carath\'eodory distance has Hausdorff dimension equal to $4$. In particular, locally, this metric space is not biLipschitz equivalent to the Euclidean space.
\end{corollary}

\proof From the general metric geometry theory, which we will review in Section~\ref{subsec_Hausdorff}, it is sufficient to prove the existence of positive constants $k_1$ and $k_2$ such that the minimal number $N_\eps$ of balls of radius $\eps$, with $\eps\in(0,1)$, needed to cover the unit ball satisfies:
\begin{equation}\label{bounds for Hdim Heisenberg}
k_1\eps^{-4}<N_\eps<k_2\eps^{-4}.
\end{equation}
For the lower bound, consider $B_1, \ldots,B_{N_\eps}$ as such balls. Using \eqref{vol_Heis}, we have:
$$\vol(B(0,1))\leq \sum_{j=1}^{N_\eps}\vol(B_j)=N_\eps \eps^4\vol(B(0,1)).$$
For the upper bound, let 
$x_1, \ldots, x_{N}$ be a maximal set (which exists by Zorn's Lemma) of points in the unit ball such that the distance between each pair is at least $\eps$.
Hence, the balls $B(x_1, \eps/2), \ldots,B(x_{N}, \eps/2)$ are disjoint balls of radius $\eps/2$ contained in the ball of radius $1+\eps/2$.
Then from \eqref{vol_Heis} we infer:
$$(1+\eps/2)^4\vol(B(0,1)) = \vol(B(0,1+\eps/2))\geq \sum_{j=1}^{N }\vol(B(x_j, \eps/2)=N \left( \dfrac{\eps}{2}\right)^4\vol(B(0,1)).$$
Therefore, using that $\eps<1$, we get:
$$ 6>(1+\eps/2)^4\geq N \dfrac{\eps^4}{16}.$$
Now, since the set $\{x_j\}_j$ is maximal, the balls $B(x_j, \eps)$, which have the same centers but radius $\eps$, make up a cover of the unit ball. Thus, we conclude:
$$N_\eps\leq N \leq 96 
\eps^{-4}.
$$
Hence, both bounds in \eqref{bounds for Hdim Heisenberg} are proven.
\qed

\subsubsection{A ball-box theorem}
In this section, we provide an elementary explanation of why the balls in the sub-Riemannian Heisenberg geometry behave like boxes with inhomogeneous sides. Specifically, let us define:
 \begin{equation}
\Bx(r) := [-r,r]\times[-r,r]\times[-r^2,r^2]\subseteq\R^3. 
\end{equation}

\begin{proposition}\label{BB_for_Heis}
 In the sub-Riemannian Heisenberg group (in the standard coordinates as above), the balls at the origin satisfy
 \begin{equation}\Bx(c_1 r)) \subset B_{cc}(1,r) \subset \Bx(c_2 r)), \end{equation}
for some universal constants $c_1, c_2>0$ and for all $r>0$.
\end{proposition}

\begin{proof}
In the following argument, we do not aim to find the best possible choices for $c_1$ and $ c_2$. Moreover, using the dilations $\delta_r$ from the previous section, one can prove the result for the unit ball and then dilate to obtain the general case. The existence of the two boxes (inside and outside) comes from the fact that the unit ball is a bounded open set. Nonetheless, we provide a direct proof without relying on the solution of the isoperimetric problem.

First, observe that for all $(x,y,z)\in B_{cc}(1,r) $ we have $|x|, |y|<r $ since the length of a horizontal curve is equal to its projection on the $xy$-plane, so effectively $\norm{(x, y)}<r $. Moreover, we claim that we have an upper bound on $z$ as a function of $r$. Indeed, we should bound the oriented area enclosed by a curve of length $r$. We stress that the curve is not closed, and the area is a signed area. In other words, the coordinate $z(t)$ satisfies \eqref{contact1a}. Hence, for the considered curve (which we might think is parametrized on the interval $[0,r]$ at unit speed, so that $\dot y, \dot x\leq 1$), we bound
\begin{equation*}
 \label{contact1a2021}
 |z(r)|=\left|\int_0^r\dfrac{1}{2} (x\dot y-y\dot x)\right|\leq \int_0^r\dfrac{1}{2} (|x| \,|\dot y|+|y|\,|\dot x|) \leq \int_0^r\dfrac{1}{2} (r 1 + r 1 )=r^2.
\end{equation*}
We then get 
$$ 
 B_{cc}(1,r) \subset [-r,r]\times[-r,r]\times[-r^2,r^2], \qquad \forall r>0. $$
 
Second, to show that the $r$-ball contains a specific box, we claim that 
\begin{equation}\label{BB_bound2021}
\left[-\frac{r}{3}, \frac{r}{3}\right]\times\left[-\frac{r}{3}, \frac{r}{3}\right]\times\left[-\frac{r^2}{100}, \frac{r^2}{100}\right]\subset B_{cc}(1,r), \qquad \forall r>0. 
 \end{equation}
Indeed, take a point $(x,y,z)$ such that 
$|x|, |y|\leq r/3 $ and $|z|\leq {r^2}/{100}$. Then, construct the following planar curve: start from $(0,0)$ and follow a square of area $z$ (clockwise if $z<0$, counterclockwise otherwise), then follow the segment from $(0,0)$ to $(x,y)$. This curve encloses area $z$, so its lift is an admissible curve reaching $(x,y,z)$. The length of the curve is four times the side length of the square plus the length of the segment.
The square has area at most $ \frac{r^2}{100}$, so its side length is at most $ \frac{r}{10}$. The segment has a length of at most $\frac{\sqrt{2}r}{3}$. From these bounds, we have
$4\frac{r}{10} +\frac{\sqrt{2}r}{3} < r .$
Therefore the point $(x,y,z)$ is inside the $r$-ball, confirming \eqref{BB_bound2021}.
\end{proof}

 \begin{figure}\index{Pansu! -- sphere}	
\centering
\subfigure[The so-called Pansu sphere is $C^\infty$ outside of the poles, and $C^2$ around them. In the above picture, the $z$-axis has been rescaled for aesthetics] 
{
 \includegraphics[height=3in]{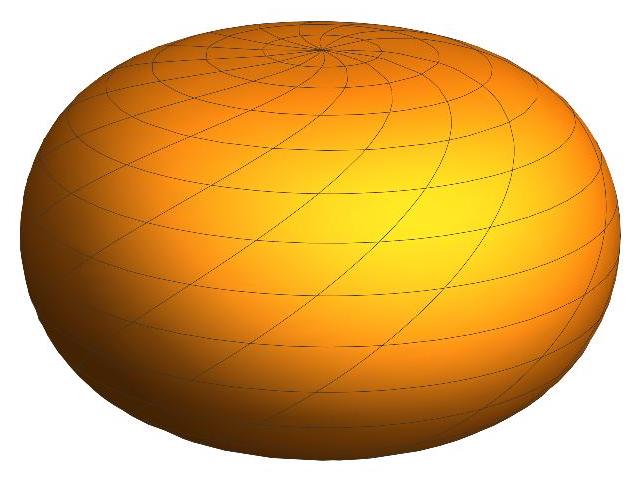}
 }
 \hspace{0.1cm}
\subfigure[Another picture of the Pansu sphere with its true axis.
] 
{
 \includegraphics[height=1.25in]{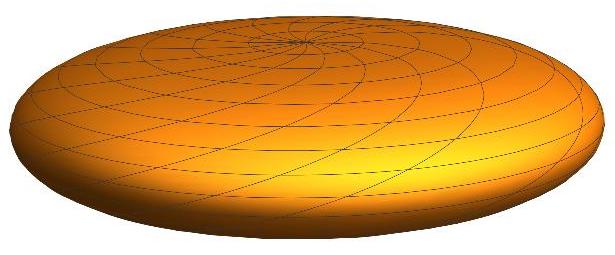}}
 \hspace{0.1cm}
\subfigure[The Pansu sphere is obtained by rotating a complete geodesic around the $z$-axis.
] 
{
 \includegraphics[height=1.25in]{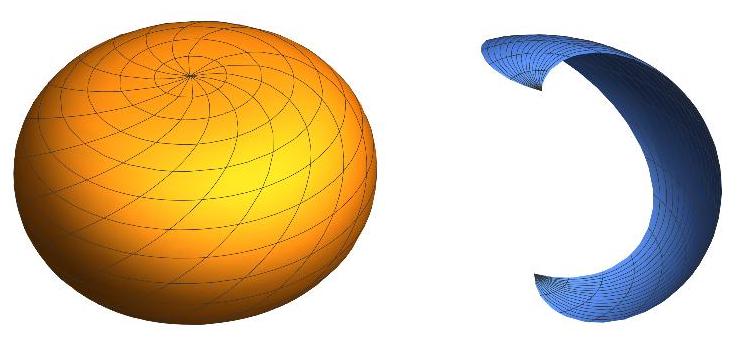}}
\caption{The (conjectured) isoperimetric sphere in the sub-Riemannian Heisenberg geometry}
\end{figure}

\section{\it Supplementary material}
 
\subsection{Dido's problem: A proof of the isoperimetric problem}\label{sec:isoper}

For a better understanding of how in Section~\ref{geod-Heis} we obtained formulas for the geodesics in the sub-Riemannian Heisenberg group, we discuss in this section the solutions of the isoperimetric problem. We then solve Dido's problem. The proof will be done under the nontrivial assumption that the minimizers of the problems are curves that are smooth enough. For the general case, we refer the reader to \cite{Morgan, Maggi}.

We shall use the formalism of Calculus of Variations to prove that each of the shortest closed curves in the plane that encloses a fixed amount of area is a circle. We will not need to show any preliminary on the curve such as the fact that it is locally a graph or that the enclosed domain is convex.
We prove that the only critical points of the variational integral functional
$$\mathcal{L}(\sigma):= \Length( \sigma), $$
subjected to the bond
$$\mathcal{A}(\sigma):= \text{Area enclosed by }\sigma=A_0, \text{ for some } A_0>0,$$
are circles.
However, we shall assume that such a $\sigma$ is a $C^1$ curve with Lipschitz derivative.

\subsubsection{Variation of length}
A necessary condition for $\sigma$ being a critical point is the vanishing of the first variation of $\mathcal{L}$.
Let $\sigma:[0,l]\to \R^2$ be a Lipschitz curve with coordinates $(\sigma_1, \sigma_2)$. Its length is given by
$$\mathcal{L}(\sigma)=\int_0^l \sqrt{\dot\sigma_1^2(t)+\dot\sigma_2^2(t)} \d t.$$
The fact that $\sigma$ is critical with respect to a variation $h$ is expressed in Calculus of Variations as follows.
\begin{definition}\index{variation}\index{variation! -- of length}
Given a Lipschitz curve $\sigma:[0,l]\to \R^2$, a {\em variation} on the interval $[0,l]$ is a smooth curve $h:[0,l]\to \R^2$ with $h(0)=h(l)=0$, and we define the associated {\em variation of length} as
$$\delta \mathcal{L}(\sigma,h):=\left.\dfrac{\dd}{\dd \eps} \mathcal{L}(\sigma+\eps h)\right|_{\eps=0}.$$
We say that a 
$\sigma$ is {\em critical in the direction of $h$} if $\delta \mathcal{L}(\sigma,h)=0.$
\end{definition}
Let us calculate the variation $\delta\mathcal{L}$ of length in the case when $\sigma$ is parametrized by arc length. So $|\dot\sigma|=1$ and $l=\text{Length}(\sigma)$. The variation in this case is
\begin{eqnarray*}
\delta \mathcal{L}(\sigma,h)&:=&\left.\dfrac{\dd}{\dd \eps} \mathcal{L}(\sigma+ \eps h)\right|_{\eps=0}\\
&=&\left.\dfrac{\dd}{\dd \eps} \int_0^l \sqrt{\left(\dot\sigma_1(t)+ \eps\dot h_1(t)\right)^2 + \left(\dot\sigma_2(t)+ \eps\dot h_2(t)\right)^2} \d t\right|_{\eps=0}\\
&=&\int_0^l\left.\dfrac{\dd}{\dd \eps} \sqrt{\dot\sigma_1(t)^2+ 2\eps \dot\sigma_1(t) \dot h_1(t)+\eps^2\dot h_1(t)^2 +\dot\sigma_2(t)^2+ 2\eps \dot\sigma_2(t) \dot h_2(t)+\eps^2\dot h_2(t)^2 } \right|_{\eps=0} \d t\\
&=&\int_0^l \left. \dfrac{ 2 \dot\sigma_1(t) \dot h_1(t)+2\eps\dot h_1(t)^2 + 2 \dot\sigma_2(t) \dot h_2(t)+2\eps\dot h_2(t)^2 }{2\sqrt{\dot\sigma_1(t)^2+ 2\eps \dot\sigma_1(t) \dot h_1(t)+\eps^2\dot h_1(t)^2 +\dot\sigma_2(t)^2+ 2\eps \dot\sigma_2(t) \dot h_2(t)+\eps^2\dot h_2(t)^2 } }\right|_{\eps=0} \d t\\
&=&\int_0^l \dfrac{ \dot\sigma_1(t) \dot h_1(t) + \dot\sigma_2(t) \dot h_2(t) }{\sqrt{\dot\sigma_1(t)^2 +\dot\sigma_2(t)^2 } } \d t\\
&=&\int_0^l \dfrac{ \langle\dot\sigma(t), \dot h(t) \rangle }{|\dot\sigma(t)|} \d t\\
&=&\int_0^l \langle\dot\sigma(t), \dot h (t)\rangle \d t.
\end{eqnarray*}
We conclude the following:
\begin{lemma}\label{lemmaA1_1}
Let $\sigma:[0,l]\to \R^2$ be a planar curve parametrized by arc length.
For every variation $h$, we have
$$\delta \mathcal{L}(\sigma,h) = \int_0^l \langle\dot\sigma, \dot h \rangle \d t.$$
\end{lemma}
\subsubsection{Area functional and its variation}
Because of Stokes Theorem, the area enclosed by a Lipschitz curve $\sigma$ can be computed by the formula
$$\mathcal{A}(\sigma)=\dfrac{1}{2}\int_0^l\sigma_1(t)\dot\sigma_2(t)-\sigma_2(t)\dot\sigma_1(t) \d t.$$
For convenience of notation, we use the {\em cross product} on $\R^2$, that is, the real number
$$v\times w:= v_1w_2-w_1v_2=\left\langle 
\begin{bmatrix}
0\\
0\\
1
\end{bmatrix},
\begin{bmatrix}
v_1\\
v_2\\
0
\end{bmatrix}
\times 
\begin{bmatrix}
w_1\\
w_2\\
0
\end{bmatrix}
\right\rangle,
\text{ for } v=(v_1,v_2), w=(w_1,w_2)\in\R^2.$$
We have linearity in $v$ and $w$ and $w\times v=-v\times w$. The area enclosed by $\sigma$ is
$$\mathcal{A}(\sigma)=\dfrac{1}{2}\int_0^l\sigma\times\dot\sigma \d t.$$
Let $h$ be a variation on $[0,l]$. The new area is
\begin{eqnarray*}
\mathcal{A}(\sigma+h)&=&\dfrac{1}{2}\int_0^l(\sigma+h)\times(\dot\sigma+\dot h) \d t\\
&=&\dfrac{1}{2}\int_0^l \sigma\times\dot\sigma + \sigma\times\dot h+h\times\dot \sigma+ h\times\dot h \d t\\
&=&\mathcal{A}(\sigma)+ \dfrac{1}{2}\left. h\times \sigma\right|_0^l+ \dfrac{1}{2}\int_0^l - \dot\sigma\times h+ h\times\dot \sigma+ h\times\dot h \d t\\
&=&\mathcal{A}(\sigma) +\int_0^l h\times\dot \sigma \d t+\dfrac{1}{2}\int_0^l h\times\dot h \d t.
\end{eqnarray*}
We deduce that
a variation $h$ is area-preserving for a curve $\sigma$ if and only if 
$$\int_0^l h\times\dot \sigma +\dfrac{ h\times\dot h}{2} \d t=0.$$
We need a weaker notion than the area preservation. We require $\mathcal{A}(\sigma+\eps h)=\mathcal{A}(\sigma) +o(\eps),$ as $\epsilon\to0$.
\begin{definition}\index{infinitesimal! -- preservation of area}
We say that a variation $h$ {\em infinitesimally preserves the area} of a curve $\sigma$
if
%
$$ \left.\dfrac{\dd}{\dd \eps} \mathcal{A}(\sigma+\eps h)\right|_{\eps=0}=0.$$
\end{definition}
By the above calculation, a variation $h$ infinitesimally preserves the area if and only if
$$ 0=\left.\dfrac{\dd}{\dd \eps}\int_0^l \eps h\times\dot \sigma +\dfrac{ \eps h\times\eps\dot h}{2} \d t\right|_{\eps=0}
=\int_0^l h\times\dot \sigma \d t.$$

\begin{proposition}
Let $\sigma:[0,l]\to \R^2$ be a curve parametrized by arc length.
If $\sigma$ is a critical curve for the length functional under an area constrain, then 
$\sigma$ has zero first variation of length with respect to all infinitesimally area-preserving variations. In particular,
$$\int_0^l \langle\dot\sigma,\dot h \rangle \d t=0,$$
for all $h:[0,l]\to \R^2$ with $h(0)=h(l) =0$ and 
$$\int_0^l h\times\dot \sigma \d t=0.$$
\end{proposition} 
\proof
Set
$a_\eps:= \mathcal{A}(\sigma+\eps h)$, which we may assume positive. Since $h$ infinitesimally preserves the area, we have
$ \left.\dfrac{\dd}{\dd \eps} a_\eps \right|_{\eps=0}=0.$
Consider the curves 
$$\sigma_\eps := \sqrt{\dfrac{a_0}{a_\eps}} (\sigma+\eps h).$$
Then $\sigma_0=\sigma$ and the area enclosed by $\sigma_\eps$ is independent of $\eps$.
Since $\sigma$ is critical for the length functional under the area constraint, we have that 
$ \left.\dfrac{\dd}{\dd \eps} \mathcal{L}(\sigma_\eps) \right|_{\eps=0}=0.$
Therefore,
\begin{eqnarray*}
0 &=& \left.\dfrac{\dd}{\dd \eps} \mathcal{L}(\sigma_\eps) \right|_{\eps=0}\\
&=& \left.\dfrac{\dd}{\dd \eps} \sqrt{\dfrac{a_0}{a_\eps}} \mathcal{L}(\sigma+\eps h) \right|_{\eps=0}\\
&=& \left.\dfrac{\dd}{\dd \eps} \sqrt{\dfrac{a_0}{a_\eps}} \right|_{\eps=0} \mathcal{L}(\sigma) 
+ \sqrt{\dfrac{a_0}{a_0}}\left.\dfrac{\dd}{\dd \eps} \mathcal{L}(\sigma+\eps h) \right|_{\eps=0}\\
&=& -\dfrac{1}{2}\sqrt{ a_0} a_\eps^{-3/2}
 \left.\dfrac{\dd}{\dd \eps} {a_\eps} \right|_{\eps=0} \mathcal{L}(\sigma) 
+ 1\cdot \delta \mathcal{L}(\sigma, h) \\
&=&0 +\int_0^l \langle\dot\sigma, \dot h \rangle \d t,
\end{eqnarray*}
where we used Lemma~\ref{lemmaA1_1}.
\qed

\subsubsection{Conclusion}
\begin{proposition}
If $\sigma$ is a $C^{1,1}$ closed curve in the plane that is one of the shortest among all Lipschitz curves that enclose the same amount of area, then $\sigma$ is a circle.
\end{proposition}
\proof Assume, without loss of generality, that $\sigma$ has unit speed.
Let $\phi:[0,l]\to \R$ be a $C^\infty$ function with $\phi(0)= \phi(l)=0$ and $\int_0^l \phi(t) \d t=0$.
Take $h(t):=\phi(t)(\dot\sigma_2(t),-\dot\sigma_1(t))$, which, since $\sigma$ is $C^{1,1}$, is Lipschitz.
Such an $h$ is an admissible variation since clearly $h(0)=h(l)=0$ and also
\begin{eqnarray*}
\int_0^l h\times\dot \sigma \, \d t&=& 
\int_0^l 
\phi(t)\left(\dot\sigma_2(t),-\dot\sigma_1(t)) \times (\dot\sigma_1(t), \dot\sigma_2(t)\right) \, \d t \\&=&
\int_0^l 
\phi(t) ( \dot\sigma_2(t)^2 +\dot\sigma_1(t)^2) \, \d t \\
&=&\int_0^l 
\phi(t) |\dot\sigma|^2 \, \d t
\\&=&\int_0^l 
\phi(t)\cdot 1 \, \d t
\\&=&\int_0^l 
\phi(t) \d t=0.\end{eqnarray*}
Then, since
$\dot h(t)=\dot\phi(t)(\dot\sigma_2(t),-\dot\sigma_1(t))+\phi(t)(\ddot\sigma_2(t),-\ddot\sigma_1(t))$,
 the vanishing of the first variation of length becomes
\begin{eqnarray*}
0&=&\int_0^l \langle\dot\sigma, \dot h \rangle \d t\\
&=&\int_0^l \left\langle(\dot\sigma_1, \dot\sigma_2), \dot\phi(t)(\dot\sigma_2(t),-\dot\sigma_1(t))+\phi(t)(\ddot\sigma_2(t),-\ddot\sigma_1(t)) \right\rangle \d t \\
&=&\int_0^l \dot\phi(t)\left\langle(\dot\sigma_1, \dot\sigma_2), (\dot\sigma_2(t),-\dot\sigma_1(t))\right\rangle +\phi(t)\left\langle (\dot\sigma_1, \dot\sigma_2), (\ddot\sigma_2(t),-\ddot\sigma_1(t)) \right\rangle \d t \\
&=&\int_0^l \dot\phi(t)(\dot\sigma_1(t)\dot\sigma_2(t)-\dot\sigma_2(t)\dot\sigma_1(t))+\phi(t)(\dot\sigma_1(t)\ddot\sigma_2(t)-\dot\sigma_2(t)\ddot\sigma_1(t)) \d t\\
&=&\int_0^l \phi(t)(\dot\sigma_1(t)\ddot\sigma_2(t)-\dot\sigma_2(t)\ddot\sigma_1(t)) \d t.
\end{eqnarray*}
The conclusion is that the function $\kappa(t):=\dot\sigma_1(t)\ddot\sigma_2(t)-\dot\sigma_2(t)\ddot\sigma_1(t)$, which is in fact the curvature of the curve $\sigma$, is such that
$$\int_0^l \phi(t)\kappa(t) \d t=0 \text{ for all } 
\phi\in C^\infty([0,l]) \text{ such that } \phi(0)= \phi(1)\text{ and }\int_0^l \phi(t) \d t=0.$$
By the (second) Fundamental Lemma of Calculus of Variations (due to DuBois and Reymond), we deduce that $\kappa$ is constant. The only planar curves of constant curvature are circles (and lines).\qed

The assumption that the curve is $C^{1,1}$ can be dropped, but the proof of the result would not be as brief. We refer to other texts for the more general result. For example, a complete proof, based on Poincare-Wirtinger inequality, can be found in \cite[pp. 1183-1185]{Osserman1978}.
The following general statement of the isoperimetric solution is for curves that are absolutely continuous.

\begin{theorem}[Isoperimetric solution]\label{isopersolution}
If $\sigma$ is an absolutely continuous closed curve in the plane that is one of the shortest among all absolutely continuous curves that enclose the same amount of area, then $\sigma$ is a parametrization of a circle.
\end{theorem}

From the solution of the isoperimetric problem, Dido's problem has an immediate solution.
\begin{theorem}[Dido's solution]	
Given two points $p$ and $q$ on the plane and a value $A>0$, the shortest curve from $p$ to $q$ that, together with the segment from $p$ to $q$, encloses area $A$ is an arc of a circle.
\end{theorem}
\proof
Assume by contradiction that there is a shortest curve $\sigma$ that is not an arc of a circle.
Let $\gamma$ be an arc of a circle joining the points and enclosing area $A$. (Notice that such an arc is unique).
Let $\hat \gamma$ be the circle of which $\gamma$ is an arc.
Let $\tilde \gamma$ be the complementary arc of $\gamma$, i.e., $\gamma$ followed by $\tilde \gamma$ is $\hat \gamma$.
Observe that the curve $\hat \sigma$ obtained following $\tilde \gamma$ after $\sigma$ is such that
$$\mathcal A(\hat \sigma)= \mathcal A(\hat \gamma)\qquad \text{ and }\qquad
\mathcal L(\hat \sigma)< \mathcal L(\hat \gamma).$$
Hence, we get a contradiction with Theorem~\ref{isopersolution}.
\qed

\section{Exercises}
Some of the following exercises require a basic understanding of differential geometry, ranging from calculus to Lie group theory. Novice readers could choose to skip some of these exercises and return to them after having read the relevant material in Chapters~\ref{ch_MGgeometry} and~\ref{ch_LieGroups}.
 
\begin{exercise}[Dido's solution] The maximal area enclosed by a curve of length $l$ on the plane together with a fixed line is $ \frac{l^2}{2\pi}$, and it is only obtained as a half disk.
\end{exercise} 
 
\begin{exercise}\label{ex prop vol contact form} The differential forms $ \vol:=\dd x\wedge \dd y$ and $\alpha:=\dfrac{1}{2} (x \dd y-y \dd x )$ in $\R^2$ have the following properties:
\\\ref{ex prop vol contact form}.i. $\dd(\alpha)= \vol;$
\\\ref{ex prop vol contact form}.ii.   In polar coordinate $(r, \theta)$, we have 
$\alpha=\dfrac{1}{2}r^2 \dd\theta;$
 \\\ref{ex prop vol contact form}.iii.  
 If $L$ is a line through the origin, then the line integral $\int_L\alpha$ is $0$.
 \end{exercise} 
 
\begin{exercise} Let $t\in [0,1]\mapsto\sigma(t)=(x(t),y(t))\in \R^2$ be a Lipschitz curve on the plane. Let $ \sigma_{[0,t]}$ be the arc restricted to the interval $[0,t]$. Let $f:\R^2\to\R$ be a smooth function. Then,
$$\dfrac{\dd}{\dd t}\left(\int_{\sigma_{[0,t]}}f(x,y)\dd x\right)=f(x(t),y(t))\dfrac{\dd x}{\dd t}(t), \qquad\text{ for almost every } t\in [0,1],$$
where $\int_{\sigma_{[0,t]}}f(x,y)\dd x$ denotes the line integral of a one-form along a curve.
\end{exercise} 
 
\begin{exercise} We have the relations 
$[X, Y]=Z$ and $[X,Z]=[Y,Z]=0$
in the following cases: (a) for the vector fields in (\ref{ONF});
			(b) for the matrices (\ref{MONF}).
\end{exercise} 
 
\begin{exercise} The inverse of an element 
$(x,y,z)$, with respect to the group structure given by (\ref{Hprod}), is $(-x,-y,-z)$.
\end{exercise} 
 
\begin{exercise} In the group structure on $\R^3$ given by (\ref{Hprod}), the lines 
$t\in\R\longmapsto\gamma_v(t) =\left(t v_1, t v_2,t v_3\right)$
are one-parameter subgroups.
\end{exercise} 
 
\begin{exercise} 
Let $L$ be a line through $0$ in the $x y$-plane of $\R^3$.
Then, the curve $L$ is a geodesic with respect to the contact distance $d_c$ defined in (\ref{contactdist}).
\end{exercise} 
 
\begin{exercise} Consider the map
$$\varphi: (x,y,z) \longmapsto \begin{bmatrix} 1 & x & z+\dfrac{1}{2} x y\\ 0 & 1 & y\\ 0 & 0 & 1\\ \end{bmatrix}$$
from $\R^3$ with the product (\ref{Hprod}) to the space of $3\times3$ upper-triangular matrices with the usual matrix product.
Then
\\ 
(a) the map is a Lie group isomorphism,
 \\ 
 (b) the map sends the standard basis $X$, $Y$, and $Z$ (defined in (\ref{ONF})) of the first Lie algebra to the standard basis $X$, $Y$, and $Z$ (defined in (\ref{MONF})) of the second Lie algebra.
\end{exercise} 
 
\begin{exercise} On the vertical $z$-axis, the distance $d_c$ defined in (\ref{contactdist}) is a multiple of the square root of the Euclidean one. Find this multiple.
 \end{exercise}

\begin{exercise}\label{dilation_of_geodesics}
Each map $\delta_\lambda$ as in \eqref{def_delta_Heis} sends geodesics of the Heisenberg group to reparametrized geodesics.
\\{\it Solution.}
 If $(\gamma_1, \gamma_2, \gamma_3)$
is a geodesic arc of length $1$ starting from the origin, 
then it is of the form (\ref{geod-in-H}) for some $k\in \R$ with
$2\pi/|k|\geq 1$, 
and the parameter $t$ of \eqref{geod-in-H} ranges from 0 to 1.
Now, consider the curve 
$(r\gamma_1,r\gamma_2,r^2\gamma_3)$ which is 
$$\left(\dfrac{ \cos\theta(\cos(kt)-1)-\sin\theta\sin(k t) }{k/r},
\dfrac{\sin\theta (\cos(kt)-1)+\cos\theta \sin(k t) }{k/r},
\dfrac{ k t -\sin(k t)}{2(k/r)^2}\right), \quad \text{ for } t\in[0,1].$$
This curve is a geodesic, albeit not parametrized by arc length but by a multiple of it, specifically $r$. Consequently, its length is $r$. 
\end{exercise} 
 
 \begin{exercise} We have Corollary~\ref{cor_dil_Heis}, from Lemma~\ref{Prop_dilations_Heis}.
\end{exercise}
 
\begin{exercise}\label{formula:geodesics}
The map $$\Phi: \left\{(\theta,k,t) \;:\; \theta\in\R/2\pi\Z, \; k\in\R, \; t \in \left(0, \frac{2\pi}{|k|}\right) \right\}\to \R^3 \setminus \left(\{0\}^2\times\R\right)$$ given
by
\begin{equation*}
\Phi(\theta,k,t)
=\left(\dfrac{ \cos\theta(\cos(kt)-1)-\sin\theta\sin(k t) }{k},
\dfrac{\sin\theta (\cos(kt)-1)+\cos\theta \sin(k t) }{k},
\dfrac{ k t -\sin(k t)}{2k^2}\right)
\end{equation*}
is a homeomorphism.
\end{exercise}

\begin{exercise} \label{ex: spheres: balls}
(i) Let $\Phi$ be the map defined in Exercise~\ref{formula:geodesics}.
The unit ball in the Heisenberg geometry is given by
\begin{eqnarray*}
 B(0,1)&=&\{\Phi(\theta,k,t) \,:\, \theta\in\R/2\pi\Z, k\in\R, t \in (0, 1)\} \\
&=&\{\Phi(\theta,k,t) \,:\, \theta\in\R/2\pi\Z, k\in[-2\pi,2\pi], t \in (0, 1)\}, 
\end{eqnarray*}
and the unit sphere is 
$$S(0,1)=\{\Phi(\theta,k,1) \,:\, \theta\in\R/2\pi\Z, k\in[-2\pi,2\pi]\}. $$ 
(ii) All the metric balls and metric spheres in the sub-Riemannian Heisenberg group are topological balls and topological spheres, respectively.
\end{exercise}

\chapter{A review of metric and differential geometry}\label{ch_MGgeometry}


Metric and differential geometry are the main tools for studying sub-Riemannian geometries. Metric geometry provides the foundation for understanding distances, geodesics, and intrinsic geometric properties in sub-Riemannian manifolds, including the broader context of Carnot-Carathéodory spaces.
Differential geometry is central and indispensable in sub-Riemannian geometry. It provides the mathematical framework for studying fundamental geometric objects such as tangent bundles and vector fields. It allows for analyzing the geometric interpretation of the sub-Riemannian distance as the minimization of a cost functional. This geometric cost functional can be viewed in metric and differential geometry as a length functional defined on curves.
To provide a clearer understanding of the setting and terminology, it is essential to have an overview of these key concepts. While there are several excellent books, such as \cite{Federer, Gromov, Ambook, Heinonenbook, Burago:book, AmbTil}, that offer clear and detailed expositions of the material, this discussion aims to provide some insights for non-experts.

 \section[Metric geometry]{Metric geometry: lengths, geodesics spaces, and Hausdorff measures}\label{MetricGeometry}
\subsection{Metric spaces}

Let $M$ be a set. A function 
\begin{equation}\label{eq_distance}
d:M\times M \to [0,+\infty]
\end{equation}
is called a {\em distance function} 
\index{distance! -- function}
(or just a {\em distance}, 
\index{distance|see {distance function}} or a {\em metric})\index{metric|see {distance function}} on $M$ if, for all $x,y,z\in M$, it satisfies
\begin{description}
\item[(\ref{eq_distance}.i)] {\em positiveness}: $d(x,y)=0 \iff x=y$,
\item[(\ref{eq_distance}.ii)] {\em symmetry}: $d(x,y)=d(y,x)$,
\item[(\ref{eq_distance}.iii)] {\em triangle inequality}: $d(x,y)\leq d(x,z)+d(z,y)$.
\end{description}
The pair $(M,d)$ is called {\em metric space}.\index{metric! -- space}
If it is clear what metric we are considering or do not want to specify the notation for the distance, 
we shall write just $M$ as an abbreviation for $(M,d)$.
We will use the term `metric' as a synonym of distance function, and never as a shortening of `Riemannian metric', which will be revised in Section~\ref{Riemannian and Finsler geometry}. 

A metric space has a natural topology which is generated by the {\em open balls} \index{balls! open --}
$$B(p,r):=\{q\in M:d(p,q)<r\}, \qquad \forall p\in M, \forall r>0.$$

We also consider distance functions that may have value $\infty$. However, on each connected component of the metric space, the distance is finite (see Exercise~\ref{d-finite}).

 For a subset $E$ of a metric space $(M,d)$, the {\rm diameter} of $E$ is defined as 
\begin{equation}
\label{def: diameter}\index{diameter}
\diam(E) :=\sup\{d(p,q)\,:\,p,q\in E\}.
\end{equation}

A {\em curve} (or {\em path}, or {\em trajectory}\index{curve}\index{path}\index{trajectory}) in
a metric space $M$ is 
 a continuous map $\gamma : I \rightarrow M$, where $I\subset\R$ is an interval. 
 The interval $I$ may be open, closed, half open, bounded, or unbounded.
 When $\gamma$ is injective, the map
 might 
 be conflated with its image $\gamma(I)$.
 We will 
 say that
 the curve $\gamma : [a,b] \rightarrow M$, with $a,b\in \R$, is a {\em curve from $p$ to $q$} (or 
 that {\em joins $p$ to $q$}) if $\gamma(a)=p$ and $\gamma(b)=q$.
 
\subsection{Length of curves in metric spaces}
\begin{definition}[Length of a curve]\label{def_length}
Let $M$ be a metric space with distance function $d$.
The {\em length (with respect to $d$)} of a curve
\index{length! -- of a curve}$\gamma : [a, b] \rightarrow M$ is
\begin{equation}\label{def_length_eq}
 L(\gamma):= \Length_d(\gamma):=\sup \left\{ \sum_{i=1}^{k} d(\gamma(t_{i-1}),\gamma(t_{i})) \;:\; k \in \mathbb{N}, a = t_0 < t_1 < \cdots < t_k = b \right\}. 
 \index{rectifiable curve}\index{curve! rectifiable --}
 \end{equation}
\end{definition} 
A {\em rectifiable curve} is a curve with finite length. Verifying that the length does not depend on the parametrization is easy, as shown in Exercise~\ref{parametrization: length}.
A curve $\gamma: [a, b] \rightarrow M$ is said to be {\em parametrized by arc length} 
\index{parametrization by arc length}
if 
$$
 \Length(\gamma|_{[t_1,t_2]})=|t_2-t_1|, \qquad \forall t_1, t_2\in[a,b]. 
$$
Every rectifiable curve admits a reparametrization by arc length; see Exercise~\ref{rect: PBAL}. With this parametrization, the curve is a $1$-Lipschitz map; see Section~\ref{sec_Lip} for the classical definition of Lipschitz map.

We shall rephrase the definition of length in terms of partitions. A {\em partition}
\index{partition}
 $\mathcal P$ of an interval $[a,b]$ is a tuple $(t_0, t_1, \ldots, t_k)\in [a,b]^{k+1}$
with $k \in \mathbb{N}$ such that 
$a = t_1 \leq t_2 \leq \cdots \leq t_k = b$.
We define
$$L(\gamma, \mathcal P):=\sum_{i=1}^{k} d(\gamma(t_{i-1}),\gamma(t_{i})).$$
Hence, we have
$$L(\gamma) = \sup \{ L(\gamma, \mathcal P) \;:\; \mathcal P \text{ is a partition of } [a,b]\}.$$

Next, we recall the lower semicontinuity of length for sequences of curves that converge pointwise. A sequence of curves $\gamma_j : [a, b] \rightarrow M$ in a metric space $M$ 
{\em converges pointwise} 
\index{convergence! -- pointwise}
to a curve $\gamma : [a, b] \rightarrow M$ in the same metric space (note that all such curves have the same interval of definition), if,
for all $ t\in [a, b]$, we have $\gamma_j(t) \to \gamma (t)$.
Furthermore, we say that $\gamma_j$ {\em converges uniformly} 
\index{convergence! -- uniform}
to $\gamma$ if
$\sup_{t\in [a, b]}d(\gamma_j(t), \gamma (t)) \to 0$, as $j\to \infty$.

\begin{theorem}[Semicontinuity of length]\label{semicontinuity:length}\index{semicontinuity of lentgh}
Let $\gamma, \gamma_1, \gamma_2, \ldots$ be curves in a metric space defined on the same interval.
If $\gamma_j \to \gamma $ pointwise, then 
$L(\gamma )\leq \liminf_{j\to\infty} L(\gamma_j ) $.
\end{theorem}
\begin{proof}
We could make the result follow from the 
fact that for each partition $ \mathcal P$ the function $\gamma\mapsto L(\gamma, \mathcal P)$ is sequentially continuous (see Exercise~\ref{Gioacchino1}) and the
general fact that the supremum of sequentially continuous functions is a sequentially lower semicontinuous function (see Exercise~\ref{Gioacchino2}). A proof of the latter fact is given by a straightforward adaptation of the following argument.

Let $\mathcal P$ be a partition of $ [a,b]$. 
Say $\mathcal P= (t_0, t_1, \ldots, t_k)$, for some $k\in \N$. Let $\eps>0$. 
Hence, there exists $N\in \N$ such that, for all $j>N$, we have
$d(\gamma_j(t_i), \gamma (t_i)) <\eps/k$, for all $i\in\{0,1, \ldots, k\}$.
By the triangle inequality, for all $j>N$, we have 
$$d(\gamma(t_{i-1}),\gamma(t_{i}))\leq d(\gamma_j(t_{i-1}),\gamma_j(t_{i})) +2\eps/k, \qquad \forall i\in\{1, \ldots, k\}.$$
See Figure~\ref{fig:triang_ineq_curves_semicontinuity_length} for a visualization. 
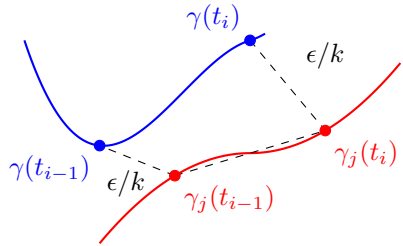
\begin{figure}[ht]
\centering
\begin{tikzpicture}
\draw[red, thick, domain=-2:2, samples=100] plot (\x, {0.3*\x^2});
\draw[blue, thick, domain=-3:0.2, samples=100] plot (\x, {-0.3*\x^3+0.8*\x^2 +0.3*\x + 1.5});
\draw[ dashed, domain=-1:1, samples=100] plot (\x, {0.3*\x});
\draw[ dashed, domain=0:1, samples=100] plot (\x, {-1.2*\x + 1.5});
\draw[ dashed, domain=-1:-2, samples=100] plot (\x, { -0.4*\x - 0.7});
 \node[above] at (1, {1}) {$\eps/k$};
 \node[below] at (-1.65, {-0.1}) {$\eps/k$};
 
\filldraw[red] (1, {0.3*1^2}) circle (2pt) node[below right] {$\gamma_j(t_{i})$};
\filldraw[red] (-1, {-0.3*1^2}) circle (2pt) node[below right] {$\gamma_j(t_{i-1})$};
\filldraw[blue] (0, { 1.5}) circle (2pt) node[above left] {$\gamma(t_{i})$};
\filldraw[blue] (-2, {0.1
}) circle (2pt) node[below left] {$\gamma(t_{i-1})$};
\end{tikzpicture}
\caption{The triangle inequality in the proof of Theorem~\ref{semicontinuity:length} for bounding $d(\gamma(t_{i-1}),\gamma(t_{i}))$.}
\label{fig:triang_ineq_curves_semicontinuity_length}
\end{figure}

Consequently, for all $j>N$, we have
$$
L(\gamma, \mathcal P)\leq L(\gamma_j, \mathcal P)+2(\eps/k)\cdot k
\leq L(\gamma_j) +2\eps.$$
Taking the liminf as $j$ goes to $\infty$ and considering that $\eps$ is arbitrarily, we obtain 
$$
L(\gamma, \mathcal P)\leq \liminf_{j\to \infty} L(\gamma_j).$$
By taking the supremum over all partitions $\mathcal P$, we conclude the result. 
\end{proof}
To show the existence of length-minimizing curves, we now recall a version of the Ascoli--Arzel\`a Compactness Theorem. 
\begin{theorem}[Ascoli--Arzel\`a]\label{AATHM}\index{Ascoli--Arzel\`a Theorem}\index{Theorem! Ascoli--Arzel\`a --}
In compact metric spaces, sequences of curves with uniformly bounded lengths contain subsequences that converge uniformly, up to reparameterization. 
\end{theorem}

\begin{proof}
Let \( (M,d) \) be a compact metric space and $(\gamma_n)_{n\in \N}$ a sequence of curves in $M$ with uniformly bounded length. Because of the bound on the lengths, 
the curves can be reparametrized with uniformly bounded constant speed to be curves $\gamma_n:[0,1] \to M$ that are uniformly Lipschitz, say $L$-Lipschitz; see Exercise \ref{rect: PBAL} and Exercise \ref{rect:constant_speed}.
The key fact of the argument of Ascoli--Arzel\`a is that the family $\mathcal F:=\{\gamma_n:n\in \N\}$ is equi-uniformly continuous (see later) and is equi-uniformly bounded (in our case this is trivial since $M$ is bounded, being compact).

Our aim is to show that $\mathcal F$ is precompact within the space $C^0([0,1]; M )$ equipped with the uniform convergence, which, when considered with the sup-distance is a complete space; see Exercise~\ref{sup: dist: complete}.
It is an exercise in topology \cite[Theorem 45.1]{Munkres} to show that in a complete metric space, a subset is precompact if and only if it is totally bounded.
Namely, by definition of totally bounded,\index{totally bounded} we need to show that for all $\eps>0$ there exists a finite set $\Lambda$ and, for all $\lambda \in \Lambda$, there exists $\mathcal F_\lambda\subset \mathcal F$ such that 
$\mathcal F= \cup_{\lambda\in \Lambda} \mathcal F_\lambda$
and $\diam \mathcal F_\lambda\leq \eps$, for all $\lambda \in \Lambda$, where the diameter is defined by \eqref{def: diameter}. 

We start from the fact that, because of the uniform Lipschitz property, the family $\mathcal F $ is {\em equi-uniformly continuous}, i.e., for every $\varepsilon>0$ there is \( \delta > 0 \) such that if \( |s-t| < \delta \),  then \( d(\gamma(t),\gamma(s)) < \varepsilon \), for all \( \gamma \in \mathcal{F} \). 
In our case, it is enough to take $\delta:=\varepsilon/L$. Given this $\delta=\delta_\varepsilon$, cover \( [0,1] \) with $k_\varepsilon$ intervals of radius \( \delta \) and center \( x_i\in [0,1] \), that is,
$  
[0,1] \subset \bigcup_{i=1}^{k_\varepsilon} B(x_i,\delta).
$
In addition, since $M$ is compact, there exists \( h_\varepsilon \in \mathbb{N}\) and points \( p_1, \ldots, p_{h_\varepsilon} \in M \) such that
$$  
M \subset \bigcup_{i=1}^{h_\varepsilon} B(p_i,\varepsilon).
$$
Next, given such $k_\varepsilon$ and $h_\varepsilon$, we define 
\[ \Lambda := \{ 
\lambda \colon \{1,\ldots,k_\varepsilon\} \ra \{1,\ldots,h_\varepsilon\}\} \] 
This set is finite, having \( h_\varepsilon^{k_\varepsilon} \) elements. We will use it as index-set. 
For $\lambda\in \Lambda$, define
$$ 
\mathcal{F}_\lambda := \{\gamma \in \mathcal{F} \,:\, d(\gamma(x_i) , p_{\lambda(i)} ) < \varepsilon \quad \forall i \in \{1,\ldots,k_\varepsilon\} \}
,$$
which is the set of those curves for which the centers of the intervals get mapped into the balls according to \(\lambda\).
Clearly, we have $\mathcal{F} = \bigcup_{\lambda \in \Lambda} \mathcal{F}_\lambda,
$
for how we choose the points $p_j$.
We need to bound the diameter of
 \( \mathcal{F}_\lambda \). 
 Pick \( \alpha,\beta \in \mathcal{F}_\lambda \) and consider their distance, given by the sup-norm. 
 For each \( t \in [0,1] \) take $i$ so that $t\in B(x_i,\delta)$.
Then
\begin{eqnarray*} 
d(\alpha(t),\beta(t)) 
&\le& d(\alpha(t),\alpha(x_i))
+ d(\alpha(x_i),p_{\lambda(i)} )+ d(p_{\lambda(i)}, \beta(x_i)) + d(\beta(x_i),\beta(t))
\\
&<& 4\varepsilon,
\end{eqnarray*}
where we used the equi-uniform continuity of \( \alpha,\beta\) and that \( \alpha,\beta \in \mathcal{F}_\lambda \).
Therefore, the family $\mathcal{F} $ is precompact.
\end{proof}

The above argument actually proves various more general statements; see, for example, Exercise~\ref{ex_AA}.
The following result is an essential consequence of Ascoli--Arzel\`a Theorem~\ref{AATHM}.

\begin{proposition}[Existence of shortest paths]\label{exists: short}\index{existence! -- of geodesics}
Let $M$ be a compact metric space.
For all $p,q\in M$ that a curve can join, there exists a curve $\gamma$ from $p$ to $q$ such that
\begin{equation}\label{Elle:inf}
L(\gamma)= \inf\{L(\sigma)\;:\; \sigma \text{ curve from } p \text{ to } q\}.
\end{equation}
\end{proposition}
\begin{proof}
Set $L$ as the right-hand side of \eqref{Elle:inf}.
If $L=\infty$, then there is nothing to prove since we can take any curve $\gamma$ joining $p$ to $q$.
We next assume that $L<\infty$. 
Let $\gamma_j$ curves from $p$ to $q$ with $L(\gamma_j)\to L$.
By Ascoli--Arzel\`a Theorem~\ref{AATHM}, up to passing to a subsequence, we may assume that 
$\gamma_j$ converges (uniformly and, hence, pointwise) to a curve $\gamma$ joining $p$ to $q$. 
By semicontinuity of length (Theorem~\ref{semicontinuity:length}), we get $L(\gamma )\leq \liminf_{j\to\infty} L(\gamma_j ) =L$.
Hence, we conclude that $L(\gamma )=L$.
\end{proof}

\subsection{Length spaces, intrinsic metrics, and geodesic spaces}\label{sec_length_sp}

If a metric space $(M,d)$ has the property that, for all $p,q\in M$, the value $d(p,q)$ is finite and 
\begin{equation}\label{def_length_sp}
d(p,q)=\inf\{\Length_d(\gamma)\;:\; \gamma \text{ curve from } p \text{ to } q \}, \end{equation}
then $(M,d)$ is called {\em length space} (or {\em path metric space}) and 
$d$ is called an {\em intrinsic metric}.
\index{length! -- space}\index{path! -- metric space}\index{intrinsic! -- metric}\index{metric! intrinsic --}
 Notice that we have chosen to require intrinsic metrics to be finite, although this decision may not be shared by all authors in the field.

If a metric space $(M,d)$ is such that, for all $p,q\in M$, there exists a curve $\gamma$ from $p$ to $q$ with the property that 
$d(p,q)= \Length_d(\gamma)$, then $(M,d)$ is called {\em geodesic space}\index{geodesic! -- space},
$d$ is called a {\em geodesic metric}\index{geodesic! -- metric}\index{metric! geodesic --}, and every such $\gamma$ is called a {\em length-minimizing curve}\index{curve! length minimizing --}
joining $p$ to $q$.
Length-minimizing curves are also called {\em length minimizers}.\index{length! -- minimizer}
Some authors use the term `geodesic' to denote length minimizers or curves that are locally length-minimizing, in agreement with Riemannian geometry.
In this text, we will call {\em geodesics} those length-minimizing curves that, in addition, are parametrized by arc length.
In other words, a curve $\gamma:[a,b]\to M$ is a {geodesics} if\index{geodesic}
\begin{equation}\label{def geodesics} d(\gamma(s), \gamma(t)) = |s-t|,\qquad \forall s,t\in [a,b].
\end{equation}

Every geodesic space is a length space (Exercise~\ref{ex: geodesic: length}).
Not all length spaces are geodesic spaces; one reason can be a lack of completeness, as, for example, $\R^2\setminus\{(0,0)\}$. As we will recall shortly, this is the only obstruction for locally compact spaces.
Recall that a topological space is called {\em locally compact}\index{locally compact} if every point of the space has a local basis of compact neighborhoods. For other topological notions, 
we refer to any standard book in topology, such as \cite{Munkres}, and we suggest playing with
the database of topological counterexamples: \url{http://topology.jdabbs.com}.\index{database! -- of topological counterexamples}

We prefer to assume the following stronger property, which gives compactness at every scale.

\begin{definition}[Boundedly compact]
A metric space is said to be {\em boundedly compact} (or {\em proper}) 
\index{boundedly compact}
\index{proper (metric space)|see {boundedly compact}}
if its bounded subsets are precompact. Equivalently, a space is 
boundedly compact
if its {\em closed balls}\index{balls! closed --}
$$\overline B(p,r):=\{q\in M:d(p,q)\leq r\}$$
are compact for all $p\in M$ and all $r>0$. Equivalently, for every $p\in M$, the distance function $q\in M\mapsto d(p,q)\in \R$ is a proper function, in the sense that preimages of compact sets are compact.
\end{definition}
\begin{proposition}\label{prop:length:geod}
Every boundedly compact length space
is a geodesic space.
\end{proposition}

\begin{proof}
Let $(M,d)$ be a 
 boundedly compact length space.
Fix $p,q\in M$. Since the metric $d$ is intrinsic, there is a curve $\gamma$ from
 $p$ and $q$ with
 $L(\gamma)< d(p,q)+1.$
 Notice that every other curve $\sigma$ from
 $p$ and $q$ with $L(\sigma)\leq L(\gamma)$
is inside 
$\overline B(p,d(p,q)+1)$, which is compact by assumption.
By Proposition
~\ref{exists: short}, we have the existence of a shortest path and hence of a geodesic joining
 $p$ to $q$ since the distance is intrinsic.
\end{proof}


With a bit more topological arguments, one can prove the following stronger result.
An explicit proof can be found in {\cite[Proposition 2.5.22 and Theorem 2.5.23]{Burago:book}}.
\begin{theorem}[Hopf-Rinow-Cohn-Vossen]\label{HRCV}\index{Hopf-Rinow Theorem}
 If a length space $(M,d)$ is complete and locally compact, then $(M,d)$ is boundedly compact and, hence, a geodesic space.
\end{theorem}

\subsection{Length as integral of metric derivative} 

Throughout the section, we will denote by $d$ the distance function of a metric space $M=(M,d)$.
\begin{definition}[Metric derivative]\label{definition metric derivative}\index{metric! -- derivative}
Given a curve $\gamma\colon [a,b]\to M$ on a metric space, we define the \emph{metric derivative} of $\gamma$ at the point $t\in (a,b)$ as the limit
\begin{equation*}
\lim_{h\to 0} \frac{d(\gamma(t+h),\gamma(t))}{\abs{h}},
\end{equation*}
whenever it exists and, in this case, we denote it by $\abs{\dot{\gamma}}(t)$.
\end{definition}

The following is the main result in this subsection:
\begin{theorem}\label{thm:length}
For each Lipschitz curve $\gamma\colon [a,b]\to M$ on a metric space, we have that
\begin{description}
\item[\ref{thm:length}.i.]\label{existence derivative} the metric derivative $\abs{\dot{\gamma}}$ exists almost everywhere;
\item[\ref{thm:length}.ii.]\label{formula variation} $\Length_d(\gamma)=\int_a^b \abs{\dot{\gamma}}(t)\, \dd t$.
\end{description}
\end{theorem}

\proof
For \ref{thm:length}.i, we start by 
noticing that by the triangle inequality 
\begin{equation}\label{triangl_gamma_2021}
 \abs{d(\gamma(s),y)-d(\gamma(t),y)} \leq d(\gamma(s),\gamma(t)), \qquad \forall s,t\in [a,b],\,\forall y\in M,
\end{equation}
with equality if $y= \gamma(t)$.
Fix a countable dense set $\{x_n\}_{n\in\N}$ in $\gamma([a,b])$ and define $$\varphi_n(t):=d(\gamma(t),x_n)
. $$Consequently, from \eqref{triangl_gamma_2021} (and its equality when $x_n\to \gamma(t))$, we have 
\begin{equation}\label{sup_phi_2021}
\sup_{n\in \N} \abs{\varphi_n(s)-\varphi_n(t)} = {d(\gamma(s),\gamma(t))} .
\end{equation}
Notice that each $\varphi_n:[a,b]\to \R$ is Lipschitz with the same Lipschitz constant as $\gamma$, and therefore differentiable almost everywhere and absolutely continuous, by the one-dimensional version of Rademacher Theorem; see \cite[Section 3.5]{Folland_book}.\index{Rademacher Theorem}
Let 
$$m(t):=\sup_n\abs{\dot{\varphi}_n(t)}.$$



We claim that \begin{equation}\label{metric_Deriv2021}
\abs{\dot{\gamma}}(t)=m(t),\qquad \text{ for almost all } t.
\end{equation}
For a first inequality, note that for each point $t$ of differentiability for ${\varphi}_n$, we have from \eqref{sup_phi_2021} that
\begin{equation*}
\abs{\dot{\varphi}_n}(t)\stackrel{\rm def}{=}\lim_{h\to 0} \frac{\abs{\varphi_n(t+h)-\varphi_n(t)}}{\abs{h}} \stackrel{\eqref{sup_phi_2021}}{\leq} \liminf_{h\to 0} \frac{d(\gamma(t+h),\gamma(t))}{\abs{h}}.
\end{equation*}
Hence
\begin{equation*}
m(t)\leq \liminf_{h\to 0} \frac{d(\gamma(t+h),\gamma(t))}{\abs{h}}.
\end{equation*}
Regarding the other inequality, using the Fundamental Theorem of Calculus, we have for $s\leq t$ that
\begin{eqnarray}\label{estimate part of variation}
d(\gamma(t),\gamma(s))&\stackrel{\eqref{sup_phi_2021} }{=}&\sup_n \abs{{\varphi}_n(t) -{\varphi}_n(s)} \nonumber\\
&=& \sup_n \abs{\int_s^t \dot{\varphi}_n(\tau)\, \dd \tau}\nonumber\\
&\leq &\sup_n \int_s^t \abs{\dot{\varphi}_n(\tau)}\, \dd \tau\nonumber\\
&\leq &\int_s^t m(\tau)\, \dd \tau.
\end{eqnarray}
Let us argue why the integral of $m$ is finite. It is because the derivative of each $\varphi_n$ is bounded from above by the Lipschitz constant of $\varphi_n$, which in turn is bounded from above by the one of $\gamma$. From Lebesgue's Differentiation Theorem, \cite[p.98]
{Folland_book}, at each Lebesgue point $t$ for $ m$ we have that 
\begin{equation*}
\limsup_{h\to 0} \frac{d(\gamma(t+h),\gamma(t))}{\abs{h}}
\stackrel{\eqref{estimate part of variation} }{\leq} \limsup_{h\to 0} \abs{\frac{1}{h}\int_t^{t+h} m(\tau)\, \dd \tau}=m(t).
\end{equation*}
So \eqref{metric_Deriv2021} holds, 
and in particular $\abs{\dot{\gamma}}$ 
exists almost everywhere.
 The first part of Theorem~\ref{thm:length} is proven.

Regarding \ref{thm:length}.ii, we first prove one inequality.
For every partition $(t_0, t_1, \ldots, t_k)$ of $ [a,b]$, for some $k\in \N$, 
we have
\begin{equation*}
\sum_{i=1}^{k} d(\gamma(t_{i-1}),\gamma(t_i))
\stackrel{\eqref{estimate part of variation}}{\leq} \sum_{i=1}^{k} \int_{t_{i-1}}^{t_{i}} m(\tau)\, \dd \tau
\stackrel{\eqref{metric_Deriv2021}}{=}
\sum_{i=1}^{k} \int_{t_{i-1}}^{t_{i}} \abs{\dot{\gamma}}(\tau)\, \dd \tau
= \int_{a}^{b} \abs{\dot{\gamma}}(\tau)\, \dd \tau.
\end{equation*}
Taking the supremum over all partitions gives $\Length(\gamma)\leq \int_a^b \abs{\dot{\gamma}}(t)\, \dd t$.

Regarding the other inequality, let $\varepsilon>0$ and $n\geq 2$ such that $h:=(b-a)/n\leq \varepsilon$. We set $t_i:=a+ih$ for $i\in \{0,1,\ldots,n\}$; so that we have $t_n=b$ and $b-\varepsilon<t_{n-1}$. Then
\begin{eqnarray}\label{estimate part of variation2}
\nonumber
 \int_a^{b-\varepsilon} d(\gamma(t ),\gamma(t+h))\, \dd t
&\leq & \sum_{i=1}^{n-1} \int_{t_{i-1}}^{t_{i}} d(\gamma(t),\gamma(t+h))\, \dd t\\
\nonumber
&=& \int_0^h \sum_{i=1}^{n-1} d(\gamma(\tau+t_{i-1}),\gamma(\tau+t_i))\, \dd \tau\\
&\leq &\int_0^h \Length(\gamma)\, \dd \tau=h \Length(\gamma).
\end{eqnarray}
Using Fatou's lemma \cite[p.52]{Folland_book}:
\begin{equation*}
\begin{split}
\int_a^{b-\varepsilon} \abs{\dot{\gamma}}(t)\, \dd t&\stackrel{\rm def}{=}\int_a^{b-\varepsilon}\liminf_{h\to 0^+} \frac{d(\gamma(t+h),\gamma(t))}{h}\, \dd t\\
&\stackrel{\rm Fatou}{\leq} \liminf_{h\to 0^+} \frac{1}{h} \int_a^{b-\varepsilon} d(\gamma(t+h),\gamma(t))\, \dd t \stackrel{\eqref{estimate part of variation2}}{\leq} \Length(\gamma).
\end{split}
\end{equation*}
Letting $\varepsilon\to 0^+$ gives the missing inequality.
\qed

\begin{example}\label{constant vector space}\index{finite-dimensional normed space}\index{normed! -- space}\index{vector space! normed --}
Some first interesting examples of geodesic metric spaces are given by finite-dimensional normed spaces. 
Let $(V,\norm{\cdot})$ be a finite-dimensional normed space. Equip $V$  with the metric $d$ induced by $\norm{\cdot}$, i.e., $d(p,q):=\norm{p-q} $, for all $p,q\in V$.
Let $\gamma:[a,b]\to V$ be 
a Lipschitz curve (either with respect to the distance $d$ or, equivalently, with respect to any other Euclidean distance). Hence, since $V$ has finite dimension, by Rademacher Theorem, the curve $\gamma$ is differentiable almost everywhere and absolutely continuous.
For every point $t$ of differentiability for $\gamma$, we have
\begin{equation*}
\norm{\gamma'(t)}=
\norm{\lim_{h\to 0}\frac{ \gamma(t+h)-\gamma(t)}{ {h}}}=
\lim_{h\to 0}\frac{\norm{\gamma(t+h)-\gamma(t)}}{\abs{h}}
\stackrel{\rm def}{=} 
\abs{\dot{\gamma}}(t),
\end{equation*}
where $\abs{\dot{\gamma}}(t)$ is the metric derivative and $\gamma'(t)$ instead denotes the (classical) derivative.
Consequently, from Theorem~\ref{thm:length} we infer 
\begin{equation}\label{same_length_str}
\Length_d(\gamma)=\int_a^b \norm{\gamma'(t)} \dd t.
\end{equation}
We deduce that for every two points $p,q\in V$ and every rectifiable curve $\gamma$ between $p$ and $q$ we have
\begin{equation}\label{norm_less_than_length}\norm{p-q} \stackrel{\rm def}{=} d(p,q)\leq \Length_d (\gamma ) = \int_a^b \norm{\gamma'(t)} \dd t.
\end{equation}
We stress that with the curve $t\in [0,1]\mapsto tp +(1-t)q$, we get equality in \eqref{norm_less_than_length}. In conclusion, we have proved that {\em every finite-dimensional normed space is a geodesic space with straight lines being length minimizering}.
\end{example}

\subsubsection{Energy functional} 
In geometric analysis, it is often more appropriate to consider the energy of curves rather than their length.
The reason is that the energy functional often possesses better analytic and geometric properties than the length functional. It may be smoother and more amenable to analysis, allowing for the application of variational techniques and optimization methods. 

Let $\gamma:[a,b]\to M$ be a Lipschitz curve on a metric space $(M,d)$. Hence, by Theorem~\ref{thm:length} its metric derivative $\abs{\dot{\gamma}}$ exists almost everywhere. The {\em energy} of $\gamma$ (with respect to the distance $d$) is defined as
\begin{equation}\label{eq:def_energy}\index{energy! -- of a curve}
 {\rm Energy}_ {d}(\gamma):=
\frac12
\int_a^b \left(\abs{\dot{\gamma}}(t)\right)^2\, \dd t
\end{equation}

Contrary to length, energy depends on the parametrization of the curve. 
However, we shall now see that parametrizations with constant speed minimize the energy among all of the reparametrizations of the curve, and in that case, the energy is a precise function of the length.

\begin{proposition}\label{prop: energy}
Let $\gamma:[a,b]\to M$ be a Lipschitz curve on a metric space $(M,d)$ and $p,q\in M$. Then, the energy satisfies the following properties:
\begin{description}
\item[\ref{prop: energy}.i.] $\Length_d(\gamma)\leq \sqrt{ 2\,(b-a)\, {\rm Energy}_ {d}(\gamma)}$.
\item[\ref{prop: energy}.ii.] If $\gamma$ is parametrized by a multiple of the arc length, then $\Length_d(\gamma)= \sqrt{ 2\,(b-a)\, {\rm Energy}_ {d}(\gamma)}$.
\item[\ref{prop: energy}.iii.] 
$\inf\{\Length_d(\gamma)\;:\;\gamma \text{ from } p \text{ to } q \} 
=
\inf\{ \sqrt{ 2\cdot{\rm Energy}_ {d}(\gamma)}\;:\;\gamma\text{ Lipschitz, on } [0,1] \text{ from } p \text{ to } q\}.$
\item[\ref{prop: energy}.iv.] A curve defined on an interval $I$, parametrized by a multiple of arc length, and going from $p$ to $q$ is length-minimizing among curves from $p$ to $q$ if and only if it is energy-minimizing among curves defined on $I$ going from $p$ to $q$.
\end{description}
\end{proposition}
\begin{proof}
The statements are straightforward consequences of Jensen's inequality or Cauchy–Shwarz's inequality.
\end{proof}

\begin{remark}\label{anche_in_Lp}
More generally, one can minimize the $p$-energies. 
For $p\in [1,\infty)$, the {\em $p$-energy} is defined as the $L^p$-norm of the metric derivative, up to possibly a normalizing multiplicative constant. 
For $p=\infty$, the $\infty$-energy is the essential supremum of the metric derivative. The 1-energy is the length, by Theorem~\ref{thm:length}.ii. 
Because of H\"older's inequality, for  $p\in (1,\infty]$, minimizing the $p$-energy among curves parametrized on $[0,1]$, exactly determines the length-minimizing curves parametrized by a multiple of arc length.
\index{$p$-energy}
\end{remark}

\subsection{Isometries, Lipschitz maps, and quasi-isometries} \label{sec_Lip}

Given two metric spaces $(X, d_X)$ and $(Y, d_Y)$, a map
$
 f: X \to Y
$
is called {\em Lipschitz} 
\index{Lipschitz! -- map}
if there exists a real constant $K \geq 0$ such that
$$
 d_Y(f(x_1), f(x_2)) \le K d_X(x_1, x_2),\qquad \forall x_1,x_2\in X.
$$
Every such a value $K$ (or, many times, the smallest value of such $K$'s) is called a (or the) {\em Lipschitz constant}\index{Lipschitz! -- constant} of the function $f$. 
A function is called {\em locally Lipschitz} if, for every $x \in X$, there exists a neighborhood $U$ of $x$ such that $f$ restricted to $U$ is Lipschitz.

If there exists a $K \geq 1$ with
$$
 \frac{1}{K}d_X(x_1,x_2) \le d_Y(f(x_1), f(x_2)) \le K d_X(x_1, x_2),
 \qquad
 \forall x_1,x_2 \in X,
$$
then $f$ is called 
{\em biLipschitz embedding}\index{bi-Lipschitz! -- embedding} (also written bi-Lipschitz or bilipschitz). 
Surjective 
biLipschitz embeddings are called {\em biLipschitz homeomorphisms} (or 
biLipschitz maps).
\index{bi-Lipschitz! -- map} 
\index{bi-Lipschitz! -- homeomorphisms} 
BiLipschitz homeomorphisms are the isomorphisms in the category 
of Lipschitz maps. 
To be more explicit about the value of the constant $K$, we would say that $f$ is $K$-biLipschitz. BiLipschitz embeddings are injective and, in fact, embeddings, i.e., they are homeomorphisms onto their image.
We call $1$-biLipschitz homeomorphisms {\em isometries}; while $1$-biLipschitz embeddings are {\em isometric embeddings}.\index{isometry}\index{isometric embeddings}

\index{bi-Lipschitz! -- equivalent functions}
Two functions $\alpha, \beta$ defined on the same set $X$ are {\em biLipschitz equivalent} if there exists $K>1$ such that
$$
 \frac{1}{K}\alpha(x) \le \beta(x) \le K \alpha(x), 
 \qquad
 \forall x \in X.
$$
Two important examples of functions for which we will consider biLipschitz equivalence will be 
distances and measures.
Notice that in particular, two distances $d_1, d_2$ on the same set $M$ 
are {\em biLipschitz equivalent} if and only if 
the identity map 
$(M,d_1)$ to $(M,d_2)$ is biLipschitz.
 \index{bi-Lipschitz! -- equivalent distances}

The notion of biLipschitz map can be further generalized to quasi-isometry. 
For the definition of the latter, we first introduce the notion of net.
If $X,Y\subset M$ are subsets of a metric space $M$ and $r>0$, we say that $X$ is an \emph{$r$-net} for $Y$ if\index{net}
\[
X\subset B(A,r):=\{ p\in M:d(p,A) < r \}.
\]
\begin{definition}[Quasi-isometry]\label{def quasi-isometry}\index{quasi-isometry}
Suppose $(M_1,d_1)$ and $(M_2,d_2)$ are metric spaces, and $f : M_1\to M_2$ is a function (not necessarily continuous). Then $f$ is called an $(L,C)$-{\em quasi-isometric embedding}, with $L\geq 1$ and $C\geq 0$, if\index{quasi-isometric embedding}
$$ \frac{1}{L}\; d_2(f(x),f(y))-C\leq d_1(x,y)\leq L\; d_2(f(x),f(y))+C\quad\mbox{for all}\quad x,y\in M_1.$$
Moreover, an $(L,C)$-quasi-isometric embedding is called an $(L,C)$-{\em quasi-isometry}
if $f(M_1) $ is a $C$-net for $M_2$.
Two metric spaces $M_1$ and $M_2$ are called {\em quasi-isometric} if a quasi-isometry exists between them.
\end{definition}

\subsection{Hausdorff measures and dimension}\label{subsec_Hausdorff}

Recall that a collection $\scr F$ of subsets of an arbitrary set $X$ is called \emph{$\sigma$-algebra} \index{algebra! $\sigma$-algebra} for $X$ if
\begin{itemize}
\item[(i)] 	$\emptyset,X\in\scr F$;
\item[(ii)] 	$A,B\in\scr F\THEN A\setminus B\in\scr F$;
\item[(iii)] 	$\{A_n\}_{n\in\N}\subset\scr F\THEN\bigcup_{n\in\N} A_n\in\scr F$.
\end{itemize}

If $X$ is a topological space, the smallest $\sigma$-algebra containing the open sets is called \emph{Borel $\sigma$-algebra}. \index{Borel! -- $\sigma$-algebra}
\begin{definition}[Measure]\label{def: measure}
 	A \emph{measure} \index{measure} on a $\sigma$-algebra $\scr F$ is a function $\mu:\scr F\to[0,+\infty]$ such that
	\begin{description}
\item[\ref{def: measure}.i.] 	$\mu(\emptyset)=0$;
\item[\ref{def: measure}.ii.] 	$\{A_n\}_{n\in\N}\subset\scr F$, pairwise disjoint $\THEN\mu(\bigcup_{n\in\N}A_n)= \sum_{n=0}^\infty\mu(A_n)$.
	\end{description}
	The latter condition is called {\em $\sigma$-additivity}.
\end{definition}

	Every measure has the property of being \emph{countably subadditive} \index{countably subadditive} on arbitrary elements of $\scr F$, i.e., 
	\begin{equation}\label{subadditivity}
	\{A_n\}_{n\in\N}\subset\scr F
	\Longrightarrow \mu\left(\bigcup_{n\in\N}A_n\right) \le \sum_{n=0}^\infty\mu(A_n),
	\end{equation}
	see Exercise~\ref{ex: countably subadditive}.

A measure on a topological space is called a \emph{Borel measure} \index{Borel! -- measure} if $\mu$ is defined on the Borel $\sigma$-algebra.
Hence, if $\mu$ is a Borel measure on a metric space $M$, then $\mu(B_M(p,r))$ is defined for all $p\in M$ and all $r>0$.

For the following definition, we use the notion of diameter from \eqref{def: diameter}.
\begin{definition}[Hausdorff measures]
Let $M$ be a metric space.
Let $S\subset M$ be a subset, $Q\in[ 0,\infty)$, and $\delta>0$. The \emph{$Q$-dimensional Hausdorff $\delta$-content} \index{Hausdorff! -- content} is defined as
\begin{equation}\label{defHau_cont}
\H^Q_\delta(S) := \inf\left\{
\sum_{i=1}^\infty \left(\diam(E_i)\right)^Q : E_i \subseteq S,\ \diam E_i<\delta,\ S\subseteq \bigcup_{i=1}^\infty E_i
\right\}, 
\end{equation}
with convention that $0^0=1$.
Notice that the function $\delta\mapsto \H^Q_\delta(S)$ is non-increasing.
The \emph{$Q$-dimensional Hausdorff measure}\index{Hausdorff! -- measure} of $S$ is defined as
\[
\H^Q(S) := \sup_{\delta>0}\H^Q_\delta(S) = \lim_{\delta\to0^+}\H^Q_\delta(S) .
\]
\end{definition}
Each measure $\H^Q$ is an outer measure, as explained in \cite{Folland_book}, which gives a measure when restricted to the Borel $\sigma$-algebra.


\begin{proposition}\label{Prop: Haus_dim_def}
 	Let $M$ be a metric space. Then there exists $Q_0\in[0,+\infty]$ such that
	\[
	\H^Q(M) = 0
	\quad\forall Q>Q_0
	\quad\text{and}\quad
	\H^Q(M)=\infty
	\quad\forall Q<Q_0 .
	\]
\end{proposition}
\begin{proof}
 	Set 
	\[
	Q_0 := \inf\{Q\ge 0:\H^Q(M)\neq\infty\} .
	\]
	Hence $\H^Q(M)=\infty$ for all $Q<Q_0$.
	
	If $Q_0=\infty$, then there is nothing else to prove.
	If $Q_0<\infty$, then take $Q>Q_0$.
	Then there is $Q'\in[Q_0,Q)$ with $\H^{Q'}(M)=:K<\infty$.
	Hence for all $\delta\in(0,1)$ we have $\H^{Q'}_\delta(M)\le K$, i.e., there are $E_i\subset M$ with $M=\bigcup_iE_i$, $\diam(E_i)<\delta$ and $\sum_i \diam(E_i)^{Q'} < K+1$.
	Notice that
	\[
	\sum \diam(E_i)^Q \le \delta^{Q-Q'} \sum_i\diam(E_i)^{Q'} < (K+1) \delta^{Q-Q'} .
	\]
	Thus $\H^Q_\delta(M)\le (K+1) \delta^{Q-Q'}$.
	Since $\delta^{Q-Q'}\to 0$ as $\delta\to0^+$, we get $\H^Q(M) = 0$.
\end{proof}

\begin{definition}[Hausdorff dimension]
 	The \emph{Hausdorff dimension} \index{Hausdorff! -- dimension} of a metric space $M$ is denoted by $\dim_H(M)$ and is equivalently defined as
	\begin{align*}
	 	\dim_H(M) 
		&:= \inf\{Q\ge0:\H^Q(M)=0\} \\
		&= \inf\{Q\ge0:\H^Q(M) \neq \infty\}\\
		&= \sup(\{Q\ge0:\H^Q(M) = \infty\}\cup\{0\}).
	\end{align*}
\end{definition}
The above definitions are equivalent because of Proposition~\ref{Prop: Haus_dim_def}.
We notice that Lipschitz maps increase in a controlled way the Hausdorff measures; see Exercise~\ref{ex_Hausdorff_Lipschitz}. 
Consequently, the Hausdorff dimension is preserved by biLipschitz homeomorphisms.


\begin{theorem}\label{thm: Ahlfors: regular: Hausdorf: dim}
 	Let $M$ be a metric space and $\mu$ a Borel measure on $M$.
	Assume that there are $Q>0$, $C>1$, and $R>0$ such that 
	\begin{equation}\label{eq2234}
	\frac1C r^Q \le \mu(B(p,r)) \le C r^Q,\qquad \forall p\in M,\ \forall r\in(0,R] .
	\end{equation}
	Then for all $p\in M$ 
	\begin{enumerate}
 \item[(i)] 	$\H^Q(B(p,R)) \in (0,\infty)$,
 \item[(ii)] 	$\dim_H B(p,R) = Q$,
	\end{enumerate}
	and, if in addition $M$ admits a countable cover of balls of radius $R$, then $\dim_HM=Q$.
\end{theorem}
\begin{proof}
	Fix $p\in M$.
	We first show that $\H^Q(B(p,R))<\infty$.
		Fix $r\in(0,R)$ and let $0<\delta< R-r$.
	We claim that we can take a finite maximal family of points $p_1,\dots,p_N\in B(p,r)$ such that $d(p_i,p_j)>\delta$ for all $i\neq j$. Indeed, such a finite set of points exists because if $p_1,\dots,p_k\in B(p,r)$ are such that $d(p_i,p_j)>\delta$, then the balls $B(p_i,\frac\delta2)$ are disjoint and contained in $ B(p, R)$, hence
	\begin{multline*}
		k \frac{\delta^Q}{2^QC} =
		\frac1C \sum_{i=1}^k \left(\frac\delta2\right)^Q
		\le \sum_{i=1}^k \mu\left(B(p_i,\frac\delta2)\right) 
	 	= \mu\left(\bigcup_{i=1}^k B(p_i,\frac\delta2)\right) 
		\le \mu(B(p,R)) \le CR^Q .
	\end{multline*}
	Therefore, the integer $k$ has to be bounded, and such a maximal set of points is finite.

	Maximality implies that $B(p_1,\delta),\dots,B(p_N,\delta)$ cover $B(p,r)$.
	Hence, we bound
	\begin{align*}
	 	\H^Q_{2\delta}(B(p,r)) 
		&\le \sum_{j=1}^N \Large(\diam(B(p_j,\delta))\Large)^Q \\
		&\le N(2\delta)^Q =
		4^Q CN \frac 1C\left(\frac\delta2\right)^Q \\
		&\le 4^QC \sum_{j=1}^N \mu\left(B(p_j,\frac\delta2)\right) \\
		&\le 4^Q C \mu (B(p,R)), 
	\end{align*}
	where in the second inequality we used that the diameter of a ball is at most twice its radius. 
	We stress that the last term is finite and independent of $\delta$. 
	Finally, for the ball of radius $R$, we have 
	$\H^Q(B(p,R)) =\H^Q( \bigcup_{r<R} B(p,r)) \le 4^Q C \mu (B(p,R))<\infty$,
	where we have used that the measure is continuous with respect to the increasing union of sets; see Exercise~\ref{ex_continuity_below}.
	
	We then show that $\H^Q(B(p,R)) > 0$.
	Let $\delta\in(0,R)$. To bound from below the $\delta$-Hausdorff content take $\epsilon>0$ and countably many sets $E_1, E_2,\ldots \subset M$ such that $\diam(E_i)<\delta$, $B(p,R)\subset\bigcup_iE_i$, and
	\[
	\H^Q_\delta(B(p,R)) \ge \sum_i(\diam E_i)^Q -\epsilon .
	\]
	Such a cover exists because $\H^Q(B(p,R))<\infty$.
	Take some $p_i\in E_i$, so $E_i\subset B(p_i,\diam (E_i))$ and 
	\[
	\mu(B(p_i,\diam(E_i))) \le C\diam(E_i)^Q .
	\]
	Thus, by the countably subadditivity \eqref{subadditivity} of $\mu$, we have, since $\bigcup_iB(p_i,\diam(E_i)) \supset \bigcup_i E_i \supset B(p,R)$,
	\begin{align*}
	\H^Q_\delta(B(p,R)) 
	&\ge \frac1C \sum_i \mu(B(p_i,\diam E_i))-\epsilon \\
	&\stackrel{\eqref{subadditivity}}{\ge} \frac1C \mu\left(\bigcup_iB\left(p_i,\diam(E_i)\right)\right) -\epsilon \\
	&\ge \frac1C \mu(B(p,R)) - \epsilon \\
	&\ge \frac1{C^2} R^Q - \epsilon .
	\end{align*}
	Since $\epsilon$ was arbitrary, we get that $\H^Q_\delta(B(p,R))$ 
	is greater than a positive constant independent of $\delta$.
	
	So $(i)$ is proved, and $(ii)$ is an immediate consequence.
	By countable subadditivity \eqref{subadditivity} of the Hausdorff measure, also the last statement of the theorem follows.
\end{proof}

\begin{remark}\label{rmk: Ahlfors: regular: Hausdorf: dim}
The above proof actually shows that the $Q$-dimensional Hausdorff measure $\H^Q$ is biLipschitz equivalent to the measure $\mu$. In particular, the measure $\mathcal H^Q$ satisfies equation \eqref{eq2234}, with possibly some other choice for the constant $C$. We shall rephrase the last theorem using the following definition.
\end{remark} 

\begin{definition}[Ahlfors regularity for measures]\index{Ahlfors regularity}
A measure $\mu$ a on 
 	 a metric space that is Borel and for which 
	there are $Q\in(0,\infty)$, $C>1$, and $R>0$ such that 
	\begin{equation}\label{Ahlfors_condition_2021}	
	\frac1C r^Q \le \mu(B(p,r)) \le C r^Q, 	\qquad \forall p\in M,\ \forall r\in(0,R],
	\end{equation}
	is said to be {\em Ahlfors $Q$-regular up to scale $R$}.
\end{definition}	

The following result is a consequence of Theorem~\ref{thm: Ahlfors: regular: Hausdorf: dim}.
\begin{corollary}\label{Corol_Ahlfors_Hausdorff}
If a metric space supports a measure that is Ahlfors $Q$-regular up to scale $R$, then the $Q$-dimensional Hausdorff measure $\H^Q$ of the metric space is Ahlfors $Q$-regular up to scale $R$, and the $R$-balls have Hausdorff dimension $Q$.
\end{corollary}

Using the Hausdorff measure, we rephrase the notion of length for injective curves.


\begin{proposition}\label{prop_H1_L}\index{length! -- of an injective curve}
If $\gamma:I\to M$ is an injective curve on a metric space $M$, then we have
\begin{equation}
\label{length_Hausdorff}
 \H^1 ( \gamma(I) )= \Length( \gamma ). 
 \end{equation}
\end{proposition}

 \begin{proof}
 We shall focus on the case when $ \Length( \gamma )<\infty$ and leave the other case to the reader; see Exercise~\ref{H1_L_infty_case}.
Thus, we reparametrize $\gamma:[0,\ell]\to M$ by arc length. 
 For proving \eqref{length_Hausdorff}, we shall consider one inequality at a time.
 
 For the inequality $\leq$, for each $\delta>0$ divide the interval $[0,\ell]$ into $n$ disjoint intervals $J_1,\ldots, J_n$ of diameter less than $\delta$. 
 Since $\gamma$ is parametrized by arc length, then it is 1-Lipschitz, and therefore, we have $\diam \gamma(J_j) < \delta$, for $j\in\{1,\ldots, n\}$.
 Hence 
 \begin{eqnarray*}
 \H^1_\delta ( \gamma([0,\ell]) )&\leq& \sum_{j=1}^n \diam \gamma(J_j)\\
 & \leq & \sum_{j=1}^n \diam J_j = \ell, 
 \end{eqnarray*}
 where we have used in the first inequality that $( \gamma(J_j) )_j$ is a admissible cover for \eqref{defHau_cont} and in the second inequality that $\gamma$ is 1-Lipschitz. Taking the limit for $\delta\to 0$, we infer the desired inequality in \eqref{length_Hausdorff}.
 
 For the inequality $\geq$, we shall use the general bound $d(\gamma(s),\gamma(t )) \leq  \H^1 (\gamma([ s,t ]))$; see  Exercise~\ref{ex_H1of_curve}. In fact, take a partition $t_0<t_1<\ldots <t_k$ of the interval $I$. Then we bound
 \begin{eqnarray*}
 \sum_{i=1}^{k} d(\gamma(t_{i-1}),\gamma(t_{i})) &\leq &
 \sum_{i=1}^{k} \H^1 (\gamma([ t_{i-1},t_{i}]))\\
 &\leq & \H^1 ( \gamma(I) ),
 \end{eqnarray*}
 where in the last inequality we have used that $\H^1$ is additive and that $\gamma$ is injective.
 \end{proof}

\subsection{Submetries}
 \begin{definition}[Submetry]\label{def submetry}
A map $\pi:X\rightarrow Y$ between metric spaces is a {\em submetry} if\index{submetry}
\begin{equation} \label{def:submetry}
\pi(\bar B(p,r)) = \bar B(\pi(p),r), \qquad \forall p\in X, \forall r>0.
\end{equation}
\end{definition}

We stress that in \eqref{def:submetry}, we consider closed balls. For boundedly compact metric spaces, it is equivalent to consider open balls; see Exercises~\ref{ex submetry open} and~\ref{ex submetry equivalent}. Also, notice that submetries are 1-Lipschitz; see Exercise~\ref{rmk:1lip}. In addition, they are open maps. They are surjective if the distance function on the target is finite valued.

We can equivalently restate the condition of submetry with the notion of parallel sets, in the following sense: two subsets $A$ and $B$ of a metric space $X$ are {\em parallel},\label{page_parallel} if for all $a\in A$ and all $b\in B$, there are $a'\in A$ and $b'\in B$ such that
 $d(A,B) = d(a, b ')=d( a ',b) $; see Exercise~\ref{ex_parallel_def} for other viewpoints.\index{parallel sets}
Then, a map $\pi: X \to Y$ is a submetry if and only if its fibers are parallel and the distance in $Y$ is exactly the distances of the fibers; see Exercises~\ref{charact-submetry} and~\ref{prop: parallelfibers-submetry}.

\subsubsection{Lifting of curves via submetries}\index{lift of curve}

\begin{figure}[h]
 \centering
 \begin{tikzcd}
 &X\ar[d,->>]{u}{\pi } \\
 {[0,T]} \ar[r]{u}{\gamma}\ar[ru]{u}{\tilde\gamma}& Y
 \end{tikzcd}
	\caption{In the presence of a submetry $\pi$, 
	each geodesic $\gamma$ in the target space can be lifted to a geodesic $\tilde \gamma$ in the source space.
Whereas each geodesic $\tilde\gamma$ realizing the distance between fibers projects to a geodesic $\gamma$ in the target.
	}
	\label{fig_diag_geod}
\end{figure}
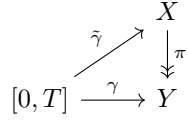

Submetries have the property that geodesics in the target space can be lifted to geodesics in the source space; see the diagram in Figure~\ref{fig_diag_geod}.
This lifting property will be important in Section~\ref{sec:sFLie_proj} when we lift curves from the abelianization of nilpotent simply connected sub-Finsler Lie groups.
Recall that a topological space $X$ is called {\em simply connected} if it is path-connected and every loop in $X$ is homotopic to a constant.\index{simply connected}

In some settings, this lifting property is equivalent to the submetry property; see
Proposition~\ref{prop_lifting_gives_submetry}.

\begin{proposition}\label{lem:pathlifting}
Let $\pi: X\rightarrow Y$ be a submetry between metric spaces.
Assume $X$ to be boundedly compact.
Then, for 
every $1$-Lipschitz curve $\gamma:[0,T]\rightarrow Y$ 
and every $x\in \pi^{-1}(\gamma(0))$,
there is a $1$-Lipschitz curve $\tilde{\gamma}:[0,T]\rightarrow X$ such that $\tilde{\gamma}(0)=x$ and $\pi\circ \tilde{\gamma} = \gamma$. 
\end{proposition}


\begin{proof}
Fix $m\in\mathbb{N}$. We define
$$\tilde{\gamma}_m: [0,T] \cap \frac{1}{m}\mathbb{Z} \rightarrow X$$
as follows.
Set $\tilde{\gamma}_m\left(0\right)= x$. By induction on $k\in \Z \cap[0,mT]$, assume that the value $\tilde{\gamma}_m\left(\frac{k}{m}\right)$ has been defined with $\pi\left(\tilde{\gamma}_m\left(\frac{k}{m}\right)\right) = \gamma\left(\frac{k}{m}\right)$.
By the submetry assumption, we know that
$$ \pi\left(\overline{B}\left( \tilde{\gamma}_m\left(\tfrac{k}{m}\right), \frac{1}{m}\right)\right) = \overline{B}\left( \gamma\left(\tfrac{k}{m}\right), \frac{1}{m}\right) \ni \gamma\left(\tfrac{k+1}{m}\right), $$
where the belonging follows from the fact that $\gamma$ is $1$-Lipschitz.
We therefore define $\tilde{\gamma}_m\left(\frac{k+1}{m}\right)$ to be any point of $\overline{B}\left( \tilde{\gamma}_m\left(\frac{k}{m}\right), 1/m\right)$ that maps onto $\gamma\left(\frac{k+1}{m}\right)$ under $\pi$.

It follows that $\tilde{\gamma}_m$ is $1$-Lipschitz for each $m\in\mathbb{N}$. By a standard Ascoli-Arzel\`a type argument, as in Theorem~\ref{AATHM}, a sub-sequence of $\{\tilde{\gamma}_m\}$ converges as $m\rightarrow \infty$ to a curve $\tilde{\gamma}$ as desired.
\end{proof}

As an immediate consequence, we can lift geodesics, which are isometric embeddings of intervals as in \eqref{def geodesics}.
As in the diagram of Figure~\ref{fig_diag_geod}, every geodesic $\gamma$ in the base is the projection of a geodesic $\tilde \gamma$ on the top.
\begin{corollary}\label{lift pf geodesics via submetries}
Let $\pi: X\rightarrow Y$ be a submetry between metric spaces.
Assume $X$ boundedly compact.
Then for every geodesic $\gamma:[0,T]\rightarrow Y$ and every $x\in \pi^{-1}(\gamma(0))$
there exists a geodesic $\tilde{\gamma}:[0,T]\rightarrow X$ such that $\tilde{\gamma}(0)=x$ and $\pi\circ \tilde{\gamma} = \gamma$. 
\end{corollary}

\begin{proof}
We apply Proposition~\ref{lem:pathlifting} to such a $\gamma$ in order to get some $\tilde{\gamma}:[0, T]\rightarrow X$ with the desired properties: to additionally see that it is a geodesic, we bound
$$|t-s| = d(\gamma(s), \gamma(t))\leq d(\tilde\gamma(s), \tilde\gamma(t))\leq |t-s|, \qquad s,t\in[0,T],$$
where we used that $\gamma$ is a geodesic and that $\pi$ and $\tilde\gamma$ are 1-Lipschitz.
\end{proof}

\begin{remark}
Proposition~\ref{lem:pathlifting} holds in a bigger generality. In fact, the notion of submetry generalizes to Lipschitz quotients. 
One can use Lipschitz quotients to lift rectifiable curves, as above. This result was first proven in \cite[Lemma~4.4]{BJLPS}
and \cite[Lemma 2.2]{MR1755254}. Though stated there for $\R^n$-targets, the proof
works the same way in the setting above, as was observed in \cite[Lemma~4.3]{MR3811770} and \cite[Lemma~3.3]{MR3964412}, where one can find a proof written in this
generality. We also point out the article \cite{MR4151334}.
\end{remark}

By considering the diagram in Figure~\ref{fig_diag_geod}, one shows that each geodesic $\tilde\gamma$ realizing the distance between fibers projects to a geodesic $\gamma$ in the target; see Exercise~\ref{ex-submetry project geodesic}.
A consequence of one such a small argument is the following fact.
\begin{proposition}\label{proj_geod_geod}
Let $\pi: X_1\to X_2$ be submetry between metric spaces.
If $X_1$ is a geodesic space, then so is $X_2$.
\end{proposition}

\begin{proof}
Take $p,q\in X_2$. Pick any $\tilde p\in \pi^{-1}(p)$, recalling that submetries are surjective.
Because $\pi$ is a submetry, there is $\tilde q\in \pi^{-1}(q)$ such that $d(p,q)=d(\tilde p, \tilde q)$. 
Since $X_1$ is a geodesic space, there is an isometric embedding $\tilde \gamma:[0,T] \to X_1$ with $\tilde \gamma(0)=\tilde p$, $\tilde \gamma(T)=\tilde q$, and $T=d(\tilde p, \tilde q)$.
Since both $\pi$ and $\tilde \gamma $ are 1-Lipschitz, then so is $\gamma:=\pi\circ \tilde \gamma$, which is a curve from $p$ to $q$. We in fact have that $\gamma$ is a geodesic, since
$d(p,q)\leq L(\gamma)\leq T=L(\tilde \gamma) = d(\tilde p, \tilde q)=d(p,q).$
\end{proof}
 
\begin{proposition}\label{prop_lifting_gives_submetry}
Let $\pi: X\to Y$ be a map between geodesic metric spaces.
Assume that
\begin{description}

\item{\ref{prop_lifting_gives_submetry}.1.} for every curve $\tilde\gamma$ in $X$ we have $L(\pi\circ\tilde\gamma)\leq L( \tilde\gamma)$ and
\item{\ref{prop_lifting_gives_submetry}.2.} for every rectifiable curve $\gamma:[0,T] \to Y$ and every $\tilde p \in \pi^{-1}(\gamma(0))$ there is a curve
$\tilde\gamma :[0,T] \to X$ such that
\begin{description}
	\item{\ref{prop_lifting_gives_submetry}.2.i.} $\tilde\gamma(0)=\tilde p$,
	\item{\ref{prop_lifting_gives_submetry}.2.ii.} $ \pi\circ\tilde\gamma=\gamma$, and
	\item{\ref{prop_lifting_gives_submetry}.2.iii.}	 $L( \tilde\gamma)=L(\gamma)$. 
\end{description} 
\end{description} Then $\pi$ is a submetry.	 
\end{proposition}

\begin{proof}
We first check that $\pi$ is 1-Lipschitz. Since the space $X$ is assumed geodesic, for every $\tilde p_1$ and $\tilde p_2$ there is a curve $\tilde \gamma $ in $X$ joining $\tilde p_1$ to $\tilde p_2$ such that 
$d( \tilde p_1,\tilde p_2)= L( \tilde \gamma)$. Then, since $\pi\circ \tilde \gamma$ joins $\pi(\tilde p_1)$ to $\pi(\tilde p_2) $ and is shorter than $\tilde \gamma$, we have
$$d( \pi(\tilde p_1),\pi(\tilde p_2)) \leq L( \pi\circ\tilde \gamma) \stackrel{\ref{prop_lifting_gives_submetry}.1}{\leq} L( \tilde \gamma)= d( \tilde p_1,\tilde p_2).$$
Thus, the map $\pi$ is 1-Lipschitz, i.e., 
$ \pi(\bar B(\tilde p,r)) \subseteq \bar B(\pi(\tilde p),r), $ for all $\tilde p\in X,$ and all $r>0$.
To prove the opposite inclusion, take $q\in \bar B(\pi(\tilde p),r)$. Since the space $Y$ is assumed geodesic, there is a curve $ \gamma $ in $Y$ joining $\pi(\tilde p)$ to $q$ with
$ L( \gamma)\leq r$.
By \ref{prop_lifting_gives_submetry}.2, there is a curve $\tilde\gamma :[0,T] \to X$ with the three properties
 \ref{prop_lifting_gives_submetry}.2.i--iii, with our choice of $\tilde p$.
 Then, the endpoint $\tilde q$ of $\tilde \gamma$ is such that $\pi(\tilde q) = q$ and
 $$ d(\tilde p, \tilde q) \leq L(\tilde \gamma) \stackrel{\ref{prop_lifting_gives_submetry}.2.iii}{=}L(\gamma)\leq r.$$
 Thus $\tilde q\in \ \bar B(\tilde p,r) $, and therefore
$ q\in \pi(\bar B(\tilde p,r)) $. Since $q$ was arbitrary in $ \bar B(\pi(\tilde p),r)$ we infer 
$\bar B(\pi(\tilde p),r)
\subseteq
\pi(\bar B(\tilde p,r)) $, for all $\tilde p\in X,$ and all $r>0$.
Altogether, the map $\pi$ satisfies \eqref{def:submetry}.
\end{proof}

\begin{proposition}\label{prop: same regularity2}
 Let
 $\pi: M_2 \to M_2$ be a smooth submetry between smooth manifolds equipped with admissible distances.
If each pair of points in $M_1$ can be joined by a $C^k$ geodesic (or a piecewise $C^k$ geodesic),
 then the same is valid on $M_2$. 
\end{proposition}
\begin{proof}
Take $p,q\in X_2$. As in the previous proof, pick any $\tilde p\in \pi^{-1}(p)$ and $\tilde q\in \pi^{-1}(q)$ such that $d(p,q)=d(\tilde p, \tilde q)$. 
Pick a `good' geodesic $\tilde \gamma:[0,T]\to G$ from $\tilde p$ to $\tilde q$, which exists by assumption. Define $\gamma:=\pi\circ \tilde \gamma $.
We stress that $\gamma$ connects $p$ to $q$ and has the same regularity of $\tilde \gamma$ since $\pi$ is smooth. Moreover, since $\pi$ is 1-Lipschitz, then $\gamma$ has length at most $d(p,q)$. Therefore, it is a geodesic.
\end{proof}


\section{Differential geometry}\label{sec_diff_geom}
\subsection{Vector fields and Lie brackets}\label{Sec:vect_fields_brackets}
In this section, we will denote by $M$ a smooth differentiable manifold. 
We will not review here the definition of a manifold nor the concept of a smooth map between manifolds, to which we refer to any introductory book, such as \cite{Lee_MR2954043}.
We denote by $C^\infty(M)$ the space of $C^\infty$ functions from $M$ to $\R.$
We shall prefer the following viewpoint for the space of smooth vector fields on $M$: A linear function $X: C^\infty(M) \to C^\infty(M) $ is a {\em smooth vector field} on $M$ if it satisfies the {\em Leibniz rule}:\index{vector field}\index{Leibniz rule}
$$ X(f g)= X(f) g + f X(g), \qquad \forall f,g\in C^\infty(M) .$$
We denote by $\Vec(M)$ or by $\Gamma(TM)$ the linear space of smooth vector fields;\index{$\Vec(M)$}\index{$\Gamma(TM)$}
we will typically use the letter $X,\,Y,\,Z$ to denote elements in $\Vec(M)$.

\begin{definition}[Vector fields in charts]\label{corresponding:base}\index{coordinate vector field} 
 	Let $\varphi: U\to\R^n$ be a coordinate chart for an $n$-manifold $M$. For 
	$j\in\{1,\dots,n\}$ we define the {\em $j$-th coordinate vector field} $\partial_j\in \Vec(U)$ by
	\[
	\partial_j ( f)(p) := \frac{\partial(f\circ\varphi^{-1})}{\partial x_j} (\varphi(p))
	= \left. \frac{\dd}{\dd t} f(\varphi^{-1}( \varphi(p) +t e_j))\right|_{t=0}, \quad \forall f\in C^\infty (U), \forall p\in U,
	\]
	where $e_j$ denotes the $j$-th element of the canonical basis of $\R^n$. 
	\index{corresponding basis vectors}\index{basis element induced by a chart} \index{$\partial_i$}
	\end{definition}
Given a chart $(U,\varphi) $ for $M$, every vector field	$X$ on $M$ restricted to $U$ can be written
 using the coordinate vector fields as
	\[
	X=\sum_{j=1}^n X^j {\partial_j}, \qquad \text{on } U, 
	\]
	 for some smooth functions $X^j\in C^\infty (U)$. Namely, we have $X(f)(p)=\sum_{j=1}^n X^j(p) {\partial_j}(f)(p)$, for all $f\in C^\infty (U)$ and all $p\in U$.


To also consider tangent vectors and form the tangent bundle of $M$, we use the notion of germs of functions: For every $p\in M$, a {\em germ of $C^\infty$ function}\index{germ of function} at $p$ is the equivalence class of smooth functions from $M$ to $\R$ with respect to the equivalence relation of being equal in some neighborhood of $p$. 
We denote by $C^\infty(p)$ the space of germs of $C^\infty$ functions at $p$.
The tangent bundle over $M$ is a set, denoted by $TM $, together with a map
 $\pi:TM\rightarrow M$ called \emph{(tangent bundle) projection map}\index{tangent! -- bundle} with the following property:
The {\em fiber space}\index{fiber! -- space}\index{tangent! -- space} $T_p M:=\pi^{-1}(p)$ of the tangent bundle $TM$ is the linear space formed by all the derivations on the space $C^\infty(p)$. 
In other words, the elements of $T_p M$, called {\em tangent vectors} at $p$, are those $\R$-linear applications
$v:C^\infty(p)\to\R$ that satisfy the Leibniz rule: $ v(f g)= v(f) g + f v(g), $ for all $f,g\in C^\infty(p) .$
 Therefore, if $X$ is a vector field on $M$ and $p$ is in $M$, then $X_p$, defined as
 $$X_p (f):= (X(f))(p), \qquad \forall f\in C^\infty(p),$$
 gives a tangent vector at $p$. Hence, vector fields on $M$ are sections of the tangent bundle $TM$, and moreover, one puts on $TM$ the structure of a manifold such that $X\in \Vec(M)$ if and only if 
 $X:M\to T(M)$ is smooth and $\pi\circ X$ is the identity on $M$. For this reason, we write $\Gamma(TM)$ for $\Vec(M)$.

If $F:M\to N$ is a smooth map between smooth manifolds and $p\in M$,
we shall denote by $\dd F_p:T_p M\to T_{F(p)}N$
its differential, defined as follows. The pull-back operator
$u\mapsto F_p^*(u):=u\circ F$ maps $C^\infty\left(F(p)\right)$
into $C^\infty(p)$; thus, for $v\in T_p M$ we have that
$$
\dd F_p(v)(f):=v( F_p^*(f)) =v(f\circ F)
,\qquad \forall f\in C^\infty(F(p)),
$$
defines an element of $T_{F(p)}N$.

Every smooth curve $\sigma : I \to M$ gives a derivation at $\sigma(t)$ for each $t\in I$ by 
$$\sigma'(t) (f):= \lim_{h\to 0} \dfrac{f(\sigma(t+h)) - f(\sigma(t))}{h}, \qquad \forall f \in C^\infty (\sigma(t) ).$$
If $F:M\to N$ is smooth and $\sigma$ is a smooth curve on $M$, then we have the formula
\begin{equation}
\dd F_{\sigma(t)}(\sigma'(t))=(F\circ\sigma)'(t),\end{equation} where
$\sigma'(t)\in T_{\sigma(t)}M$ and $(F\circ\sigma)'(t)\in
T_{F(\sigma(t))}N$ are the tangent vectors along the two
curves, in $M$ and $N$, respectively. If $f\in C^\infty(M)$ and $p\in M$, identifying
$T_{f(p)}\R$ with $\R$ itself, given $X\in \Gamma(TM)$, we have
$$\dd f_p(X_p)=X_p(f).$$

 	For a vector field $X\in \Gamma(TM)$, a smooth curve $\sigma:(a,b)\to M$ is an \emph{integral curve}, or a \emph{flow line, of $X$} if\index{integral curve}
\index{flow! -- line}
	\[
	\sigma'(t) = X_{\sigma(t)},\qquad\forall t\in(a,b) .
	\]
 	For all $X\in\Gamma(TM)$ and all $p\in M$ there are $a<0$, $b>0$, and $\sigma:(a,b)\to M$ such that $\sigma$ is an integral curve of $X$ and $\sigma(0)=p$.
		Moreover, such a $\sigma$ is unique and has a unique maximal extension.
We denote by
\index{flow}
 $t\mapsto \Phi_X^t(p)$ the integral curve of $X$ starting at $p$.
We call $\Phi_X^t(p)$ the \emph{flow at $p$ at time $t$ with respect to $X$}. Namely, we have
\begin{equation}
\left\{\begin{array}{l}
\quad\Phi_X^0(p)=p,\\
 \dfrac{\dd}{\dd t}\Phi_X^t(p) = X_{\Phi_X^t(p)}. 
\end{array}\right.
\end{equation}

One of the fundamental notions that we will utilize in our study is the Lie bracket of vector fields. The Lie bracket of vector fields has several equivalent definitions, and we will employ them all based on the viewpoint that we consider.

\begin{definition}[Lie bracket]\label{def:Lie_bracket_vector_fields}\index{Lie! -- bracket}
	The \emph{Lie bracket of vector fields} on a manifold $M$ is the map 
	\begin{align*}
	[ \cdot, \cdot ] : \mathrm{Vec}(M)\times\mathrm{Vec}(M)&\rightarrow \mathrm{Vec}(M)\\
	(X, Y)&\mapsto [X, Y]
	\end{align*}
	defined with any of the equivalent viewpoints a--d:
	\begin{description}
	\item[\ref{def:Lie_bracket_vector_fields}.a.] Viewpoint of \textbf{derivations}: For $f\in C^\infty(M)$, 
	\[
	[X,Y](f )=X(Yf)-Y(Xf).
	\]

	\item[\ref{def:Lie_bracket_vector_fields}.b.] Viewpoint in \textbf{coordinates}:
	In local coordinates, if two vector fields are given by $X=\sum_{k=1}^n X^k\partial_k$ and $Y=\sum_{k=1}^n Y^k\partial_k$ for some smooth functions $X^1,\ldots, X^n$ and $Y^1,\ldots, Y^n$, then 
	\[
	[X,Y]=\sum_{h,k=1}^n\left(X^h \partial_h Y^k-Y^h \partial_h X^k\right)\partial_k.
	\]

	\item[\ref{def:Lie_bracket_vector_fields}.c.] Viewpoint of \textbf{Lie derivative}: For $p\in M$,\index{Lie! -- derivative} 
	\[
	[X,Y]_p=\left.\frac{\dd}{\dd t}\left((\dd\Phi_X^{t})^{-1} Y_{\Phi_X^t(p)}\right)\right|_{t=0} =: \left(\mathcal L_X Y\right)_p.
	\]

	\item[\ref{def:Lie_bracket_vector_fields}.d.] Viewpoint of \textbf{commutation of flows}: For $p\in M$, \begin{eqnarray*}
	[X,Y]_p&=&\left.\frac{1}{2}\frac{\dd^2}{\dd t^2}\left(\Phi^{-t}_Y\circ\Phi^{-t}_X\circ \Phi^t_Y \circ \Phi^t_X\right)(p)\right|_{t=0}\\
	 &=& \left.\frac{\dd}{\dd t} ( \phi _Y^{-\sqrt t} \circ \phi _X^{-\sqrt t} \circ \phi _Y^{\sqrt t} \circ \phi _X^{\sqrt t}) (p)\right|_{t=0^+} .
	\end{eqnarray*}
 
	\end{description}
\end{definition}
	The Lie bracket induces on $\Vec(M)$ an infinite-dimensional Lie algebra structure; see Definition~\ref{def: Lie_algebra}.
Clearly, the push-forward via a diffeomorphism commutes with the Lie bracket operation; see Exercise~\ref{ex: push-forward_bracket}. Here, if $F: M\to N$ is a diffeomorphism and $X\in\Gamma(TM)$, the {\em push forward vector field} $F_*X\in \Gamma(TN)$ is defined by the identity\index{push forward vector field}
$(F_*X)_{F(p)}:=\dd F_p(X_p)$, for $p\in M$. Equivalently,
\begin{equation} (F_*X)f:=[X(f\circ F)]\circ F^{-1},\qquad\forall
f\in C^\infty(N).\end{equation}

\subsection{Vector bundles}\label{Vector bundles}

A simple example of a vector bundle of rank $r$ over a manifold $M$ is the product space $M\times\R^r$ with the projection on the first component $\pi_1:M\times\R^r\to M$.
The next important example of a vector bundle of rank $\dim(M)$ over a manifold $M$ is the tangent bundle $TM$ of $M$. The abstract definition is the following.

\begin{definition}[Vector bundle]
	\index{vector bundle}\index{rank! -- of a vector bundle} 
 	A \emph{vector bundle of rank $r$} over a manifold $M$ is a manifold $E$ together with a smooth surjective map $\pi:E\to M$ such that, for all $p\in M$, the following properties hold:
	\begin{enumerate}
	\item 	The fiber $E_p:=\pi^{-1}(p)$ is equipped with the structure of a vector space of dimension $r$.
	\item 	There is a neighborhood $U$ of $p$ in $M$ and a diffeomorphism $\chi:\pi^{-1}(U)\to U\times \R^r$ such that
	\begin{enumerate}
		\item 	$\pi_1\circ\chi = \pi$
		\item 	for all $q\in U$, the restricted map $\chi|_{E_q}:E_q\to \{q\}\times\R^r$ is an isomorphism of vector spaces.
		\end{enumerate}
	\end{enumerate}
	\begin{center}
 \begin{tikzcd}
\bigcup_{q\in U} E_q= \pi^{-1}(U)\ar[rr]{u}{\chi} \ar[rd,->>]{u}{\pi}&&U\times \R^r = \bigcup_{q\in U} \left( \{q\}\times\R^r\right) \ar[dl,->>]{u}{\pi_1}\\
&U& 
 \end{tikzcd}
\end{center}
	The space $E$ is called \emph{total space}, the manifold $M$ is the \emph{base}, the vector space $E_p$ is the \emph{fiber over $p$}, and every such a map $\chi$ is called a \emph{local trivialization}. 
	\index{total space of a bundle}
	\index{base of a bundle}
	\index{fiber! -- of a bundle}
	\index{local! -- trivialization of a bundle}
\end{definition}

\begin{definition}[Section]
	\index{section}
 	A \emph{section} of a vector bundle $\pi:E\to M$ is a smooth map $\sigma:M\to E$ such that $\pi\circ\sigma = \Id_M$. We will denote by $\Gamma(E)$ the set of all sections of $E$.
\end{definition}
\vspace{-1cm}
\[
\xymatrix{
E \ar[d]_{\pi}\\
M \ar@/^-2pc/[u]_{\sigma}
		}
\]

\begin{definition}[Frames and local frames]
	\index{frame}\index{local! -- frame}
	A \emph{frame} of a bundle $\pi: E\to M$ is a set $\{X_1,\dots,X_n\}\subset\Gamma(E)$ of sections on $M$ such that, for all $p\in M$, the $n$-tuple $(X_1(p),\dots,X_n(p))$ is a basis of the fiber $E_p$.
	A \emph{local frame} for $\pi: E\to M$ at a point $p\in M$ is a frame for the bundle $\pi|_{\pi^{-1}(U)} : \pi^{-1}(U) \to U$, where $U$ is some open neighborhood of $p$.
\end{definition}

\subsection{Riemannian and Finsler geometry}\label{Riemannian and Finsler geometry}
Let $M $ be a differentiable manifold of dimension $n$. A {\em Riemannian metric} \index{Riemannian! -- metric} on $M$ is a family of (positive-definite) inner products
$$
 \rho_p : T_p M\times T_p M\longrightarrow \R,\qquad p\in M,
$$
such that, for all smooth vector fields $X,Y$ on $M$, we have that
$$
 p\longmapsto \rho_p(X_p, Y_p)
$$
defines a smooth function from $M $ to $\R$. This smooth assignment of an inner product $\rho_p$ to each tangent space $T_pM$
 is called a {\em metric tensor}, or {\em Riemannian metric tensor}. \index{metric! -- tensor}\index{Riemannian! -- metric tensor}
 A metric tensor will also be denoted by $\langle \cdot,\cdot \rangle$.
Endowed with one such metric tensor, the pair $(M,\langle \cdot,\cdot \rangle)$ is called a {\em Riemannian manifold}. \index{Riemannian! -- manifold}

Given a chart $(U,\varphi) $ for the manifold $M$ we have the coordinate vector fields
$ {\partial_1},\dots, {\partial_n}$ from Definition~\ref{corresponding:base}, and we consider the {\em components of the metric tensor
relative to the coordinate system} as 
$$
 \rho_{i j}(p):=
 \rho_p\left(\left. {\partial_i}\right|_p,\left. {\partial_j}\right|_p\right),
 \qquad \forall p\in U.
$$
It is easy to verify that the functions $ ( \rho_{i j})_{ij}$ are smooth and contain all the information about $\rho$. 
\\

Finsler manifolds generalize Riemannian manifolds by no longer assuming that they are infinitesimally Euclidean. Namely,  on each tangent space, we have a norm but it is not necessarily induced by a scalar product. Two good references on Finsler geometry are \cite{Bao-Chern} and \cite{Abate-Patrizio}.

Classically, a Finsler structure on a differentiable manifold $M$ is given by a function $\norm{\cdot}: T M\to \R$ that is smooth on the complement of the zero section of $T M$ and
such that 
 the restriction of $\norm{\cdot}$ to every tangent space $T_p M$ is a (symmetric) norm (see Remark~\ref{remark_def_Finsler}). 
We will consider a more general definition for Finsler structures: as regularity, we only assume the continuity in the point and the convexity in the vector. 

Every Riemannian manifold $(M,\langle \cdot,\cdot \rangle)$ has an associated function $TM\to[0,\infty)$, $X\mapsto \|X\|:=\sqrt{\langle X,X \rangle}$.
This function is an example of a continuously varying norm.

\begin{definition}\label{def: continuously varying norm}
 	A \emph{continuously varying norm}\index{norm! continuously varying --} 
	 \index{continuously varying norm}
	 on a differentiable manifold $M$ is a continuous function from $TM$ to $[0,\infty)$ usually denoted by $\|\cdot\|$ 
	 $$ \|\cdot\|:TM \to [0,\infty), \quad X\in TM \mapsto \|X\|,$$
	 with the property that for all $p\in M$ the restriction of $\|\cdot\|$ to $T_pM$ is a symmetric norm, i.e.,\index{symmetric norm}\index{norm! symmetric --}
	\begin{enumerate}
	\item 	$\|\lambda X\| = |\lambda|\|X\|$, $\forall X\in TM$, $\forall \lambda\in\R$;
	\item 	$\|X+Y\| \leq \|X\|+\|Y\|$, $\forall p\in M$ and $\forall X,Y\in T_pM$;
	\item 	$\|X\|=0\THEN X=0$.
	\end{enumerate}
\end{definition}

\begin{definition}\label{Finsler manifold}
 	In this text, we say that a \emph{Finsler manifold} \index{Finsler! -- manifold} is a pair $(M,\|\cdot\|)$ where $M$ is a differentiable manifold and $\|\cdot\|$ is a continuously varying norm on $M$, in the sense of Definition~\ref{def: continuously varying norm}. In this case, the function $\|\cdot\|$ is also called \emph{Finsler structure}.\index{Finsler! -- structure}
\end{definition}
\begin{example}\label{esempi_Finsler}
There are at least two situations that we want the reader to keep in mind:
\begin{description}
\item[\ref{esempi_Finsler}.i.] 	Every Riemannian manifold $(M,\langle \cdot,\cdot \rangle)$ naturally has the structure of a Finsler manifold.
\item[\ref{esempi_Finsler}.ii.] Every finite-dimensional normed vector space naturally has the structure of a Finsler manifold. 
\end{description}
\end{example}
\begin{Rem}\label{remark_def_Finsler}
 	The notion of Finsler manifold is present in the literature with different meanings. On the one hand, the norm is classically required to be smooth (away from the zero section) and with a positive Hessian. Namely, some authors assume that norms for Finsler structures have strongly convex smooth unit spheres, while we do not in Definition~\ref{def: continuously varying norm}.
	On the other hand, some authors considered other weak notions of norms.
 For example, they allow asymmetric norms, i.e., the first condition in Definition~\ref{def: continuously varying norm} is assumed only for $\lambda>0$. 
\end{Rem}


\section{Length structures for Finsler manifolds}
Connected Riemannian and Finsler manifolds carry the structure of length metric spaces. Let us recall the notion of absolutely continuous curve and 
its length with respect to a Finsler structure.

\begin{definition}\label{def:AC}
 	A curve $\gamma:[a,b]\to\R^n$ is \emph{absolutely continuous} \index{absolutely continuous} if there exists a Lebesgue integrable $\R^n$-valued function $g:[a,b]\to \R^n$ such that 
	\[
	\gamma(t) - \gamma(a) = \int_a^t g(s) \dd s,
	\qquad\forall t\in[a,b].
	\]
	The function $g$ is sometimes denoted by $\dot\gamma$; however, it is only defined almost everywhere with respect to the Lebesgue measure on $[a,b]$.
	A curve $\gamma:[a,b]\to M$ into a differentiable manifold is said \emph{absolutely continuous} (or, {\em AC}, for short), if it is so when read in local coordinates, i.e., for all local coordinate map $\phi: U\to\R^n$ and for all $a',b'\in[a,b]$ such that $\gamma([a',b'])\subset U$, then $\phi\circ\gamma|_{[a',b']}$ is absolutely continuous.
For every absolutely continuous curve $\gamma:[a,b]\to M$, one can also define a {\em derivative} $\dot\gamma :[a,b]\to TM$ using local coordinates, which is defined almost everywhere as a measurable map (see Exercise~\ref{ex: 36584842}).\index{derivative! -- of AC curve, $\dot \gamma$}\index{$\dot \gamma$}	
\end{definition}
As it is usual for other notions in differential geometry, to check that a curve $\gamma:[a,b]\to TM$ is
absolutely continuous
it is sufficient that the image of the curve admits a covering of coordinate systems for $M$ on which $\gamma$ is 
absolutely continuous (see Exercise~\ref{ex: 1243124}).

\begin{definition}[Length of a curve in a Finsler manifold]\index{Finsler! -- length}\index{length}
 	Let $(M,\|\cdot\|)$ be a Finsler manifold in the sense of Definition~\ref{Finsler manifold}.
	Let $\gamma:[a,b]\to M$ be an absolutely continuous curve.
	We define 
\begin{equation}\label{def_Finsler_length}
	\Length_{\|\cdot\|}(\gamma) := \int_a^b\|\dot\gamma(t)\|\dd t .
\end{equation}
\end{definition}
We remark that the {\em Finsler length} \eqref{def_Finsler_length} of an absolutely continuous curve is finite because the derivative is integrable by assumption and $\|\cdot \|$ is continuous.

The arc length is independent of the chosen parametrization, as can be shown using the change-of-variables formula.
 In particular, a curve $\gamma:[a,b]\rightarrow M$ can be parametrized by its arc length, i.e., in such a way that
\index{parametrization by arc length}
$$
 \Length_{\|\cdot\|}(\gamma|_{[t_1,t_2]})=|t_2-t_1|, \qquad\forall t_1, t_2 \text{ with } a\leq t_1\leq t_2 \leq b. 
$$
 A curve is parametrized by arc length if and only if $\norm{\dot\gamma(t)}=1$, for almost all $t\in[a,b]$.

The {\em Finsler distance function} $d_{\|\cdot\|} : M\times M \to [0,+\infty)$ is defined by
\begin{equation}\label{def_Finsler_distance}\index{Finsler! -- distance function}
 d_{\|\cdot\|} (p,q) = \inf \Length_{\|\cdot\|}(\gamma),
\end{equation}
where the infimum is taken among all absolutely continuous curves $\gamma$ in $M$ joining $p $ to $q$.

The function $d_{\|\cdot\|} $ satisfies the properties of a distance function for a metric space. The only property that is not completely straightforward is that $d_{\|\cdot\|} (p,q)=0$ implies $p=q$. 
To prove this property, we claim that, locally in a coordinate system, every Finsler structure (as every Riemannian structure) is biLipschitz equivalent to the Euclidean structure, i.e., for some $c>0$, we have
\begin{equation}\label{equiv_to_Euclid}c^{-1}\|\cdot\|\leq \|\cdot\|_{\mathbb E} \leq c \|\cdot\|,\end{equation}
 where $\|\cdot\|_{\mathbb E} $ is the Euclidean norm.
Indeed, let $U\subseteq \R^n$ be an open set parametrizing the manifold and fix a compact set $K\subseteq U$, which we think has a nonempty interior. 
Consider $T^1K:= \{ (p,v): p\in K, v\in T_pU,\|v\|_{\mathbb E} =1 \}$ the bundle of unit vectors on $K$. Notice that $T^1K$ is compact. Hence, the continuous function $\|\cdot\|$ on $T^1K$ admits maximum and minimum; moreover, the minimum cannot be 0 since, otherwise, we would have a non-zero vector with norm 0. We deduce that there exists a constant $c>0$ such that if $p\in K$ and $v\in T_pK$ is such that $ \|v\|_{\mathbb E} =1$ then 
$c^{-1} \leq \|v\| \leq c$. By homogeneity, we have \eqref{equiv_to_Euclid} on $K$. 

Consequently, based on \eqref{equiv_to_Euclid}, we can establish the biLipschitz equivalence between distance functions. Specifically, we have proven that every two Finsler distance functions on the same manifold are biLipschitz equivalent on compact sets. We summarize our findings in the following proposition.
\begin{proposition}\label{prop_Finsler_Riem_top}
On every Finsler manifold 
in local coordinates, on compact sets, the Finsler distance function is biLipschitz equivalent to the Euclidean distance function.
Consequently, on every compact set of every manifold, every Finsler structure is biLipschitz equivalent to every Riemannian structure. 
In particular, Finsler distance functions induce the same topology as the manifold topology. 
\end{proposition}

On each Finsler manifold to every continuously varying norm, as defined in Definition~\ref{def: continuously varying norm}, we associated a length structure as in \eqref{def_Finsler_length} and a distance function as in \eqref{def_Finsler_distance}.
The distance function then induces another length structure, as in Definition~\ref{def_length}.
Next, we show that these two length structures coincide. 


\begin{proposition}\label{thm: lengthscoincide}
Assume $M$ is a differentiable manifold equipped with a continuously varying norm $\norm{\cdot}\colon TM\to \R$ with induced length structure $\Length_{{\norm{\cdot}}}$ and distance function $d_{\norm{\cdot}}$. 
If $\gamma\colon [a,b]\to M$ is an absolutely continuous 
 curve, then
\begin{equation}\label{same_length_str_on_mfds}
\Length_{d_{\norm{\cdot}}}(\gamma)=\Length_{{\norm{\cdot}}}(\gamma).
\end{equation}
\end{proposition}

\proof 
To prove the $\leq$ inequality in \eqref{same_length_str_on_mfds}, notice that for all $t, s\in [a,b]$ we have 
%
%
\begin{equation*}
d_{\norm{\cdot}}(\gamma(s),\gamma(t)) \stackrel{\rm def}{=}\inf_\sigma \int_{s}^{t} \norm{\dot{\sigma}(\tau)}\, \dd \tau \leq \int_{s}^{t} \norm{\dot{\gamma}(\tau)}\, \dd \tau
\stackrel{\rm def}{=}\Length_{{\norm{\cdot}}}(\gamma|_{[s,t]})
,\end{equation*}
where the infimum is taken over all AC curves $\sigma$ from $\gamma(s)$ to $\gamma(t)$.
Using the definition of length, \eqref{def_length_eq}, and the additivity of integrals, we deduce that $\Length_{d_{\norm{\cdot}}}
\leq \Length_{{\norm{\cdot}}}$.

Regarding the other inequality, we shall use the fact that the norm changes continuously. It is convenient to work in coordinates, and it is enough to prove our claim locally. Parametrizing $M$ with an open subset $U$ of $\R^n$, we write the norm as $\norm{v}_x=:F(x,v)$, for $x\in U$ and $v\in T_x U \simeq \R^n$.
Fix some $K>1$. Since $F$ is continuous and homogeneous in the second variable, then at each point $p\in U$, there exists a neighborhood $U_p$ of $p$ such that
\begin{equation}\label{compare metric}
\frac{1}{K} F(q,v)\leq F(p,v)\leq KF(q,v), \qquad \forall q\in U_p, \forall v\in \R^n .
\end{equation}
We find a partition $a=a_0<a_1<\cdots <a_n=b$ and points $p_1, \ldots, p_n \in M$ such that 
\begin{equation}\label{dentro_U}
\gamma ([a_{i-1},a_{i}])\subseteq U_{p_i}, \qquad \forall i\in\{1,\ldots, n\}.\end{equation}
Let us denote by $d_i$ 
the distance induced by the (constant) norm $F(p_i,\cdot)$.
Since we are in the case of a normed vector space (see Example~\ref{constant vector space}), we have
\begin{equation}\label{ccaso_banach}
\Length_{F(p_i,\cdot)} = \Length_{d_i}.
\end{equation} 
Moreover, 
as a consequence of \eqref{compare metric}, we have
\begin{equation}\label{Lipschitz_control}
{d_i} \leq K {d_{\norm{\cdot}}}.
\end{equation} 
Thus, using \eqref{compare metric}, \eqref{ccaso_banach}, and \eqref{Lipschitz_control}, together with \eqref{dentro_U}, we obtain that
\begin{eqnarray*}
\Length_{{\norm{\cdot}}}(\gamma) &\stackrel{\text{def}}{=}&
\int_a^b F(\gamma(t),\dot{\gamma}(t))
\\
&=&\sum_{i=1}^{n} \int_{a_{i-1}}^{a_{i}} F(\gamma(t),\dot{\gamma}(t))\\
&\stackrel{\eqref{compare metric}}{\leq}&
K\sum_{i=1}^{n} \int_{a_{i-1}}^{a_{i}} F(p_i,\dot{\gamma}(t))
\\
&\stackrel{\eqref{ccaso_banach}}{=}&
 K\sum_{i=1}^{n} \Length_{d_i}(\gamma |_{[a_{i-1},a_{i}]})
 \\
&\stackrel{\eqref{Lipschitz_control}}{\leq}&
K^2\sum_{i=1}^{n} \Length_{d_{\norm{\cdot}}}(\gamma|_{[a_{i-1},a_{i}]})=K^2\Length_{d_{\norm{\cdot}}}(\gamma).
\end{eqnarray*}
As $K$ can be chosen arbitrarily close to $1$, we also deduce that $
\Length_{{\norm{\cdot}}}
\leq
\Length_{d_{\norm{\cdot}}}
 $.
\qed

\begin{remark}\label{Un_altro_remark_curve} 
Let $\gamma:[a,b]\to M$ be a curve on a manifold that is equipped with a continuously varying norm $\norm{\cdot}$. With the following points, we shall clarify the relationship between absolute continuity (AC) and having finite length:
\begin{description} 
\item[\ref{Un_altro_remark_curve}.i.] 
If $\gamma$ is AC, then $\Length_{{\norm{\cdot}}} (\gamma)
=
\Length_{d_{\norm{\cdot}}} (\gamma)$ and both these quantities are finite; see Proposition~\ref{thm: lengthscoincide}.

\item[\ref{Un_altro_remark_curve}.ii.] If $\gamma$ is not AC, then $\Length_{{\norm{\cdot}}} (\gamma)$ is not defined.

\item[\ref{Un_altro_remark_curve}.iii.] If $\Length_{d_{\norm{\cdot}}} (\gamma)$ is finite, then up to reparametrization $\gamma$ is Lipschitz with respect to ${d_{\norm{\cdot}}}$, and thus with respect to any Euclidean distance, in coordinates. Therefore, by Rademacher Theorem, the curve $\gamma$ is AC.
\end{description}
\end{remark}

\section{Exercises}

\begin{exercise}\label{d-finite}
Let $(M,d)$ be a metric space equipped with its natural topology.

(i) If $M$ is connected, then $d$ is finite.

(ii) In general, the function $d$ is finite on each connected component of $M$. 
\end{exercise}

 \begin{figure}[h]
\centering
 \includegraphics[width=5in]{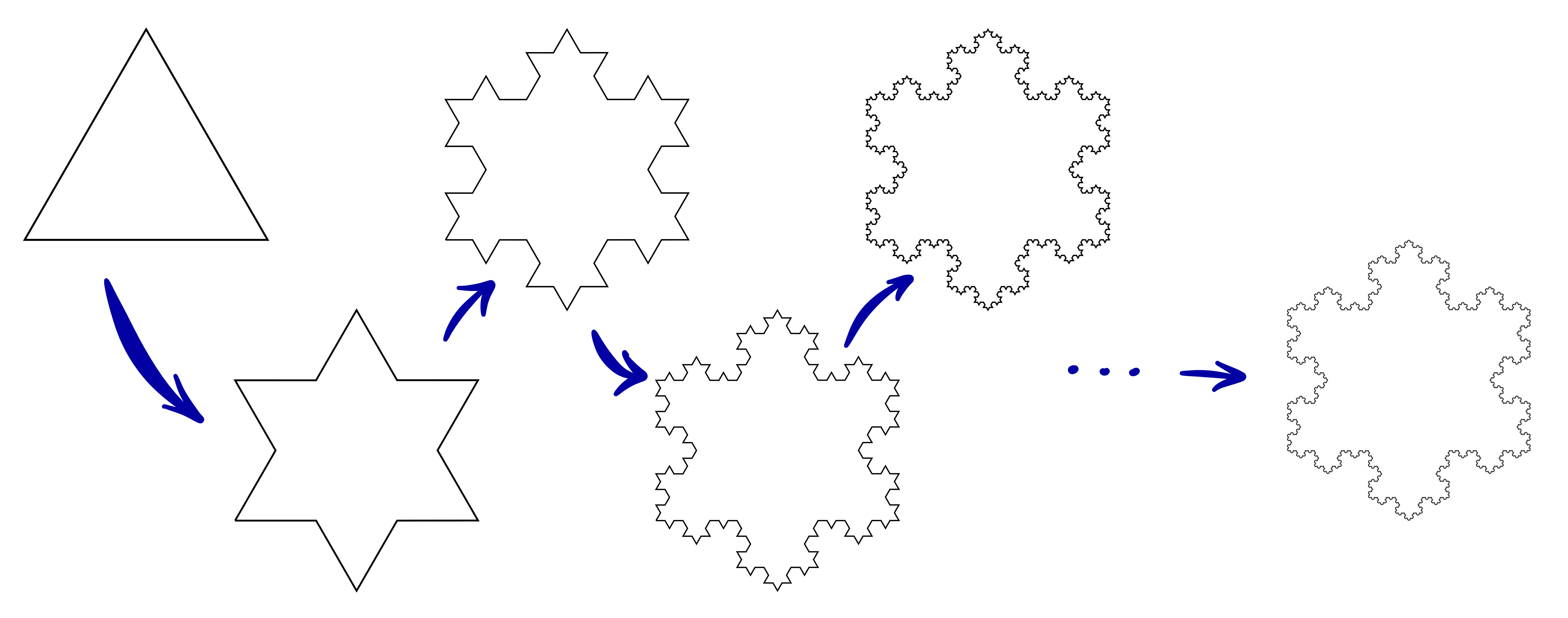}
\caption{Construction of the von Koch snowflake.}
\label{fig von Koch}
\end{figure} 
\begin{exercise}[Von Koch showflake]\label{ex_Koch}\index{snowflake! Von Koch --}
Start with a curve that is the boundary of an equilateral triangle and change iteratively each line
segment as follows to get a new curve: Remove the middle third of each segment and
replace it with the outward pointing side of an equilateral triangle without the base so that you get four line segments of equal length. Repeat now this procedure infinitely and add at each iteration smaller and smaller equilateral triangles on the previous line segments. The resulting limit curve $K$ is called the {\em von Koch snowflake}.
Equip $K\subset \R^2$ with the distance function $d_{\R^2}$ of the plane $\R^2$.
Then, the metric space $K$ is bi-Lipschitz equivalent to the circle $\mathbb S^1\subset \R^2$ equipped with the distance
$(d_{\R^2})^{\log(3)/\log(4)}$.  
\end{exercise} 
\begin{exercise}[Showflake of a metric space]\index{snowflake! -- of a metric space}
Let $(M,d)$ be a metric space and $\alpha\in (0,1)$. Then, the pair $(M,d^\alpha)$ is a metric space, called the  $\alpha$-{\em snowflake} of $(M,d)$. The Hausdorff dimension of $(M,d^\alpha)$ is $\alpha$ times the Hausdorff dimension of $(M,d)$.
The von Koch snowflake from Exercise~\ref{ex_Koch} is the $\log(3)/\log(4)$-snowflake of $\S^1$ and has dimension $\log(4)/\log(3)>1$.
\end{exercise}

\begin{exercise}
\index{mesh of a partition}
The {\em mesh} of a partition $\mathcal P=(t_1, \ldots, t_k)$ is defined as
$$\norm{\mathcal P}:= \max_{j\in\{1, \ldots, k-1\}} |t_{i+1} - t_i|.$$
If $\mathcal P_j$ are partitions such that $\norm{\mathcal P_j}\to 0$ as $j\to \infty$, then 
$L(\gamma)=\lim_{j\to \infty} L(\gamma, \mathcal P_j).$
\end{exercise}
\begin{exercise}
If $\mathcal P_1$ and $\mathcal P_2$ are partition of the same interval with 
$\mathcal P_1\subset\mathcal P_2$, then 
$L(\gamma, \mathcal P_1)\leq L(\gamma, \mathcal P_2).$
\end{exercise}
\begin{exercise} \label{parametrization: length}
The length of a curve is independent of its parameterization. Namely,
If $\gamma:I\to M$ is a curve in a metric space and $h:J\to I$ is a 
continuous weakly monotone surjection between intervals, then
$L(\gamma) = L(\gamma\circ h)$.
\end{exercise}

\begin{exercise}\label{rect: PBAL}
	If $\gamma:[a,b]\to (M,d)$ is rectifiable, then it can be reparametrized by arc length.
	{\it Hint:} consider the change of parametrization given by $s\mapsto \Length(\gamma|_{[a,s]})$.
\end{exercise}

\begin{exercise}\label{rect:constant_speed}
If $\gamma:[a,b]\to (M,d)$ is 
{\em 
parametrized with constant speed} $s$, with $s\in [0,\infty)$, i.e., 
$$
 \Length(\gamma|_{[t_1,t_2]})=s |t_2-t_1|, \qquad \forall t_1, t_2\in[a,b],
$$
then $L(\gamma) = s |a-b|$ and $\gamma $ is $s$-Lipschitz.
	\end{exercise}

\begin{exercise}\label{Gioacchino1}
For each partition $ \mathcal P$, if a sequence $(\gamma_n)_{n\in \N}$ of curves pointwise converges to $\gamma$ then 
$L(\gamma_n, \mathcal P)$ converges to $L(\gamma, \mathcal P)$. \end{exercise}

\begin{exercise}\label{Gioacchino2} Let $f_n: X\to \R$ be a sequence of continuous functions on a topological space. Then, the function $\sup_n f_n$ is lower semicontinuous. {\it Hint:} adapt the proof of Theorem \eqref{semicontinuity:length}. 
 \end{exercise}
 
 \begin{exercise}\label{sup: dist: complete}
Let $ (M,d) $ be a complete metric space, and let
$ 
\mathcal{F} :=C^0(I; M )$ be the family of all curves from a fixed interval $I$ into $M$. 
Endow \( \mathcal{F} \) with the metric
\[ d_{\text{sup}}(\sigma,\gamma) = \sup_{t\in I}\{d_M(\sigma(t),\gamma(t))\}, \qquad \forall \sigma, \gamma \in \mathcal F. \]
Then, the pair \( (\mathcal{F},d_{\text{sup}}) \) is a complete metric space.
\end{exercise}

\begin{exercise}[Ascoli--Arzel\`a]\label{ex_AA}\index{Ascoli--Arzel\`a Theorem}\index{Theorem! Ascoli--Arzel\`a --}
Let $X$ and $Y$ metric spaces. Assume $X$ is boundedly compact compact. 
If a family $\mathcal F $ of maps from $X$ to $Y$ is equicontinuous 
 and, for each $x\in X$, the set $\{f(x):f\in \mathcal F\}$ is precompact,
 then $\mathcal F $ is precompact in the uniform topology on compact sets.
\end{exercise}

\begin{exercise}[The space of $L$-Lipschitz functions]\label{sec01081010}
Let $X$ and $Y$ be metric spaces, with $X$ boundedly compact.
Let $\Lip^L(X; Y)$ be the set of Lipschitz functions $X\to Y$ of Lipschitz constant at most $L$. Fix a base point $o\in X$
Consider the function
\[
d_L(f,g) := \sup
	\left\{\frac{d(f(x),g(x))}{n^2} : \ n\in\N,\ x\in B(o,n) \right\} .
\]
\begin{description}
\item[\ref{sec01081010}.i.] 	The function $d_L$ is a distance function on $\Lip^L(X;Y)$.

\item[\ref{sec01081010}.ii.] 
	Convergence with respect to $d_L$ is equivalent to uniform convergence on compact sets.

\item[\ref{sec01081010}.iii.] The space $(\Lip^L(X;Y),d_L)$ is separable.	
\end{description}
{\it Solution of \ref{sec01081010}.ii.}
		Let $\{f_k\}_k\subset\Lip^L(X;Y)$ and $f\in\Lip^L(X;Y)$.
	Suppose that $\lim_{k\to\infty} d_L(f_k,f) = 0$.
	If $E\subset X$ is compact, then there is $N\in\N$ such that $E\subset B(o,N)$.
	Since
	\[
	\sup\{d(f_k(x),f(x)):x\in B(o,N)\} \le N^2 d_L(f_k,f) \to 0,
	\]
	then $f_k\to f$ uniformly on $E$.
	Since $E$ is an arbitrary compact set, $f_k\to f$ uniformly on compact sets.
	Suppose now that $f_k\to f$ uniformly on compact sets and let $\epsilon>0$.
	Since $\{o\}$ is compact, there is $C>0$ such that $d(f_k({o}),f({o}))\le C$ for all $k\in\N$.
	Notice that, for all $n\in\N$ and $x\in B(o,n)$, we have
	\begin{align*}
	\frac{d(f_k(x),f(x))}{n^2} 
	&\le \frac{d(f_k(x),f_k(o))+d(f_k(o),f(o))+d(f(o),f(x))}{n^2} 
	\le \frac{2L}{n} + \frac{C}{n^2}.
	\end{align*}
	Therefore, there is $N\in\N$ such that $\frac{d(f_k(x),f(x))}{n^2}<\epsilon$ for all $n\ge N$ and $x\in B(o,n)$.
Since $X$ is boundedly compact, then $B(o, N)$ is precompact. Hence, there is $K\in\N$ be such that 
	\[
	\sup\{d(f_k(x),f(x)):x\in B(o,N)\} \le \epsilon, \qquad\forall k>K.
	\]
	Then, for $k>K$, we have $d_L(f_k,f) \le \epsilon$.
	We conclude that $\lim_{k\to\infty} d_L(f_k,f) = 0$.
\\{\it Solution of \ref{sec01081010}.iii.}
	Fixing $o'\in Y$, 
	by Ascoli--Arzelà result (see Exercise~\ref{ex_AA}), for every $n\in\N$ the set
	\[
	\scr K(n) := \{f\in \Lip^L(X;Y): f(o)\in \bar B(o',n)\}
	\]
	is compact, hence separable.
	Since $\Lip^L(X; Y)=\bigcup_{n\in\N} \scr K(n)$ is a countable union of separable sets, then it is separable.
\end{exercise}

\begin{exercise}
Let $F: M_1\to M_2$ be a map between two metric spaces that is $K$-Lipschitz.
If $\gamma $ is a curve in $M_1$, then 
$L(F\circ \gamma) \leq K \cdot L(\gamma).$
\end{exercise}
\begin{exercise}
\label{ex: geodesic: length}
Every geodesic space is a length space -- what is not automatic is that the distance is finite.
\end{exercise}
\begin{exercise}
Every complete locally compact length space is boundedly compact. 
\\{\it Hint:} see \cite[Proposition 2.5.22]{Burago:book}.
\end{exercise}
\begin{exercise}
Let $L\geq 0$ and $F:M_1\to M_2$ a map that is {\em locally $L$-Lipschitz}\index{locally! -- Lipschitz}, i.e., for all $p\in M_1$ there is $r>0$ such that $F_{B(p,r)}$ is $L$-Lipschitz.
If $M_1$ is a length space, then $F$ is $L$-Lipschitz.
\end{exercise}
\begin{exercise}
Find a homeomorphism $F:M_1\to M_2$ 
between two metric spaces with the property that
$L(F \circ \gamma) = L(\gamma)$ for all curves $\gamma$ in $M_1$,
but $F$ is not an isometry.
\end{exercise}

\begin{exercise}\label{ex: countably subadditive}
	Each measure is {countably subadditive} as in \eqref{subadditivity}. 
	\\
	{\it Hint:} Given countably many sets, 
	split them into disjoint sets and apply Assumption \ref{def: measure}.ii.
\end{exercise}

\begin{exercise}
Let $\gamma:[a,b]\to\R^n$ be {absolutely continuous}.
The vector-valued function $\dot\gamma$ is unique up to a measure-zero subset of $[a,b]$.
\end{exercise}

\begin{exercise}\label{ex_continuity_below}
We have the {\em continuity from below} for measures, i.e., for every measure $\mu$ on a space $X$, if $E_1\subseteq E_2\subseteq \ldots\subseteq X$ are in the domain of $\mu$ then
$\mu(\cup_{i=1}^\infty E_i) =\lim_{i\to\infty} \mu( E_i) .$
\end{exercise}

\begin{exercise}\label{ex_Hausdorff_Lipschitz}
	If $F:M_1\to M_2$ is an $L$-Lipschitz map between metric spaces, $Q\ge 0$, and $S\subset M_1$, then
	\[
	\H^Q(F(S)) \le L^Q\H^Q(S) .
	\]
Consequently,
	if $F:M_1\to M_2$ is a biLipschitz homeomorphism, then $\dim_HM_1=\dim_HM_2$.
\end{exercise}

\begin{exercise}[to be generalized in Exercise~\ref{ex_H1of_connected}]\label{ex_H1of_curve}
For every continuous curve $\gamma:[a,b]\to M$ on a metric space $M$, we have
$$\H^1 ( \gamma([a,b]) ) \geq d(\gamma(a), \gamma(b) ). $$
{\it Solution}. Consider $\phi(x):= d(x, \gamma(a))$, which is 1-Lipschitz.
Then, using that on $\R$ the measure $\H^1$ coincides with Lebesgue measure, bound 
$\H^1 ( \gamma([a,b]) ) 
\geq \H^1 (\phi( \gamma([a,b]) ) )
\geq 
\diam(\phi( \gamma([a,b]) ) )
\geq d(\gamma(a), \gamma(b) ). $ 
\end{exercise}

\begin{exercise}\label{H1_L_infty_case}
Complete the proof of Proposition~\ref{prop_H1_L} by showing that for every curve $\gamma:I\to M$ on a metric space
 if $ \Length( \gamma )=\infty$, then $\H^1 ( \gamma(I) )=\infty$.
{\it Hint}. Use Exercise~\ref{ex_H1of_curve}.
\end{exercise}

\begin{exercise}\label{ex_H1of_connected}
For every connected subset $X$ of a metric space, we have
$\H^1 ( X ) \geq \diam(X ). $
\end{exercise}

\begin{exercise}[Maximal separated nets]\index{net}\index{separated set}
Given a metric space $M$ and $\delta>0$, a set $E\subseteq M$ is called $\delta$-{\em net} if 
$M\subseteq B(E,\delta):=\{p\in M:d(p,E)<\delta\}$.
While a set $E\subseteq M$ is called $\delta$-{\em separated} if for all distinct $p,q\in E$, we have $d(p,q)\geq \delta$.
Then, every separated set that is maximal (with respect to the inclusion) is a net.
Moreover, for each $\delta>0$, every compact metric space has a finite maximal $\delta$-separated set.
\end{exercise}

\begin{exercise}\index{quasi-isometry}
Let $X$ and $Y$ be metric spaces.
Then, there exists a quasi-isometry from $X$ to $Y$ if and only if there is a net in $X$ that is biLipschitz homeomorphic to a net in $Y$.
\end{exercise}

\begin{exercise}[Doubling distance]\index{doubling! -- metric space}\index{metrically doubling}
A metric space $M$ is called {\em doubling}, or, also, {\em metrically doubling}, 
if there is a constant $K\in \N$ such that, for $p\in M$ and $r>0$, there are $p_1, \ldots, p_K$ such that 
$B(p,r)\subseteq \bigcup_{i=1}^K B(p_i,r/2).$
Then, a metric space $M$ is doubling if and only if there is a constant $K\in \N$ such that, for every $d>0$, every set in $M$ of diameter $d$ can be covered by $K$ sets of diameter at most $d/2$. 
\end{exercise}

\begin{exercise}[Assouad dimension]\index{doubling! -- metric space}
A metric space $M$ is doubling if and only if there exist $C>1$ and $\beta>0$ 
such that, for every $p\in M$ and $r>0$, every $\eps r$- separated set in the ball $B(p,r)$ in $M$ has cardinality at most $C \eps^{-\beta}$. The infimum of such $\beta$'s is called {\em Assuoad dimension}.\index{Assouad dimension}
\end{exercise}

\begin{exercise} 
A metric space $M$ is doubling if and only if there exist $C>1$ and $\beta>0$ 
such that, for every $d, \eps>0$, every set in $M$ of diameter $d$ can be covered by at most $C \eps^{-\beta}$ sets of diameter at most $\eps d$. 
\end{exercise}

\begin{exercise}[Doubling measure]\index{doubling! -- measure}\label{ex_doubling_measure}
A measure $\mu$ on a
metric space $M$ is called {\em doubling}
 if 
 there is a constant $C>1$ such that 
\begin{equation}
0<\mu (B(p,2r)) \leq C \mu (B(p, r)) <\infty, \qquad \forall p\in M,\forall r>0.
\end{equation}

(i) If a metric space $(M,d)$ admits a doubling measure $\mu$, then the space is metrically doubling. Hence, we call the triple $(M,d,\mu)$ a {\em doubling metric measure space}.\index{metric! -- measure space}

(ii) On every metrically doubling metric space that is complete, there is a doubling measure. 
See \cite[Theorem~13.3]{Heinonenbook}.
\end{exercise}

\begin{exercise}\index{doubling! -- measure}\label{ex_doubling_measure2}\index{homogeneous! -- measure}
Every doubling measure $\mu$ on a metric space $M$ is $\alpha$-{\em homogeneous}, for some $\alpha>0$ (and some $C>0$), in the sense:
\begin{equation}\frac{\mu (B(p, r)) }{\mu (B(p, R))} \leq C\left(\frac{r}{R}\right)^\alpha,\qquad \forall p\in M, \forall 0< r \leq R <\diam (M).
\end{equation}
 See \cite[Equation~(4.16)]{Heinonenbook}.
\end{exercise}


\begin{exercise}\label{rmk:1lip}
For every map $\pi: X\rightarrow Y$ between metric spaces, the following are equivalent: 
\\i) $\quad \pi( B(p,r)) \subseteq B(\pi(p),r), \qquad \forall p\in X, \forall r>0$;
\\ii) $\quad\pi(\bar B(p,r)) \subseteq \bar B(\pi(p),r), \qquad \forall p\in X, \forall r>0$;
\\iii) $\quad\pi$ is 1-Lipschitz.
\end{exercise}

\begin{exercise}\label{ex submetry open}
If $\pi:X\rightarrow Y$ is a submetry as in Definition~\ref{def submetry}, then 
\begin{equation} \label{def:submetry open}
\pi( B(p,r)) = B(\pi(p),r), \qquad \forall p\in X, \forall r>0,
\end{equation}
where now we are considering open balls.
\end{exercise}

\begin{exercise}\label{ex submetry equivalent}
Let $X$ and $Y$ be metric spaces, with $X$ assumed to be boundedly compact. 
Then $\pi: X\rightarrow Y$ is a submetry as in Definition~\ref{def submetry} if and only if \eqref{def:submetry open} holds. Give a counterexample for when $X$ is not assumed boundedly compact. 
\end{exercise}


\begin{exercise} \label{charact-submetry}
Let $\pi: X\rightarrow Y$ be a surjective map between metric spaces. Then $\pi$ is a submetry if and only if
 for all $\hat{p}, \hat{q} \in Y$ and all $p\in \pi^{-1}(\hat p)$, there exists $q\in \pi^{-1}(\hat q)$ such that 
$d(p,q) = d(\hat{p},\hat{q}) = d(\pi^{-1}(\hat p),\pi^{-1}(\hat q)) $. 
\end{exercise}

\begin{exercise}[Hausdorff distance]\label{ex_parallel_def}\index{Hausdorff! -- distance}
Consider two subsets $A$ and $B$ of a metric space $X$.
The {\em Hausdorff distance} between $A$ and $B$ is 
$$d_H(A, B ) := \max \left\{
\sup_{a\in A} 
\{d(a, B )\},\sup_{b\in B} 
\{d(A, b)\}\right\}.$$
If $A$ and $B$ are parallel, as in the sense defined at page \pageref{page_parallel}, then $d(A,B)= d_H(A, B )$.
If $X$ is boundedly compact, $A$ and $B$ are compact, and $d(A,B)= d_H(A, B )$, then $A$ and $B$ are parallel.
\end{exercise}

\begin{exercise}\label{prop: parallelfibers-submetry}
Let $X$ be a metric space and $Y$ be a nonempty set. Let $\pi: X\rightarrow Y$ be a surjective map. Assume that the fibers of $\pi$ are parallel, i.e., for all $\hat{p} \in Y$, all $\hat{q} \in Y$, and all $p\in \pi^{-1}(\hat p)$, one can find $q\in \pi^{-1}(\hat q)$ such that
 $d(\pi^{-1}( \hat{p}),\pi^{-1}(\hat q)) = d(p,q)\, .$ 
Then 
\begin{equation*} 
d(\hat{p},\hat{q}):= d(\pi^{-1}(\hat{p}),\pi^{-1}(\hat q)), \qquad \forall \hat{p},\hat{q}\in Y,
\end{equation*}
defines a distance function on $Y$ and $\pi$ is a submetry from $X$ onto $Y$.
\end{exercise}

\begin{exercise}\label{ex-submetry project geodesic}
Let $\pi: X\rightarrow Y$ be a submetry between metric spaces. Let $p,q\in X$ be such that 
$d(p,q) = d(\pi^{-1}(\pi(p)), \pi^{-1}(\pi(q)) )$. If $\tilde\gamma$ is a geodesic between $p$ and $q$, then $\pi\circ\tilde\gamma$ is a geodesic between $\pi(p)$ and $\pi(q)$.
\\{\it Hint.} Check the proof of Proposition~\ref{proj_geod_geod}.
\end{exercise}



\begin{exercise}
	If $E$ is a vector bundle of rank $r$ over a manifold $M$, then $\dim(E)=\dim(M)+r$.
\end{exercise}
\begin{exercise}
If $\pi:E\to M$ is a vector bundle and $U\subset M$ is an open set, then
$\pi|_{\pi^{-1}(U)} : \pi^{-1}(U) \to U$ is a vector bundle.
\end{exercise}

\begin{exercise}\label{ex: 1243124}
Let $\gamma: I\to M$ be a curve into a differentiable manifold.
Then, the curve $\gamma$ is absolutely continuous if
for all $t\in I$ there exist $\eps >0$
and
a local coordinate map $\varphi:U\subset M\to\R^n$ 
with $\gamma([t-\eps,t+\eps])\subset U$ and 
such that 
 $\varphi\circ\gamma|_{[t-\eps,t+\eps]}$ is absolutely continuous.
\end{exercise}

\begin{exercise}\label{ex: 36584842}
Let $\gamma: I\to M$ be an absolutely continuous curve into a differentiable manifold. Let $\varphi_1, \varphi_2:U\subset M\to\R^n$ be two coordinate maps.
Then, the derivative of $\varphi_1\circ\gamma$ is related to the derivative of $\varphi_2\circ\gamma$ by the differential of $\varphi_1\circ \varphi_2^{-1}$ and hence one can define the derivative $\dot\gamma$ up to measure-zero sets.
\end{exercise}

\begin{exercise}
Every absolutely continuous curve in $\R^n$ can be re-parametrized to be a Lipschitz curve with respect to the Euclidean distance.
\end{exercise}


	\begin{exercise}\label{ex: push-forward_bracket}
The push-forward commutes with the Lie bracket: if $F: M\to N$ is a diffeomorphism of differentiable manifolds, then
\begin{equation}\label{bir}
[F_*X,F_*Y]=F_*[X,Y],\qquad\forall X,\,Y\in\Gamma(TM).
\end{equation}
\end{exercise}

	\begin{exercise}\label{Seba 3 oct 2024} 
	If $X$ and $Y$ are vector fields tangent to a smooth submanifold $N\subseteq M$ of a differentiable manifold $M$, then also $[X, Y]$ is tangent to $N$.
\end{exercise}

\begin{exercise}
 	Let $M$ be a differentiable manifold. For all $X,Y \in\mathrm{Vec}(M)$ and for all $f,g\in C^\infty(M)$
	\[
	[fX,gY] = fg[X,Y] + f(Xg) Y - g(Yf) X .
	\]
\end{exercise}
 
\chapter{General theory of Carnot-Carath\'eodory spaces}\label{ch_CCspaces}


We have reached the point where we are ready to introduce the main object of our investigation: sub-Riemannian manifolds and, more generally, sub-Finsler manifolds, also known as Carnot-Carathéodory spaces. 
These spaces will be equipped with Carnot-Carathéodory distances.
Our first significant result is the Chow-Rashevsky Theorem, which states that on every sub-Finsler manifold, the Carnot-Carathéodory distance induces the same topology as the manifold structure itself. 
It is important to emphasize that this result relies on the crucial assumption that the horizontal subbundle generates brackets.

\section{Definition of Carnot-Carath\'eodory spaces}

In this chapter, we denote by \(M\) a differentiable manifold, with its dimension primarily indicated as \(n\). The tangent bundle of \(M\), denoted as \(TM\), is a \(2n\)-dimensional manifold with the following local parametrization: if \(\varphi: U \subset \mathbb{R}^n \to M\) is a local parametrization for \(M\), then it induces vector fields \(\partial_{x_1}, \ldots, \partial_{x_n}\). The map \(U \times \mathbb{R}^n \to TM\), \((x, v) \mapsto v_1 \partial_{x_1}|_{\varphi(x)} + \cdots + v_n \partial_{x_n}|_{\varphi(x)}\), is a local parametrization for \(TM\). In other words, the vector fields \(\partial_{x_1}, \ldots, \partial_{x_n}\) form a local frame for \(TM\).

\subsection{Bracket-generating distributions}
\begin{definition}[Polarization, a.k.a. distribution or tangent subbundle]\label{def_distribution}\index{distribution}\index{polarization}\index{polarized! -- manifold}
A {\em distribution of tangent subspaces} on a manifold $M$ is a subset $\Delta\subseteq TM$ such that for every $\bar p\in M$ there exists smooth vector fields $X_1, \ldots, X_m$ defined on some neighborhood $U$ of $\bar p$ such that
\begin{equation}\label{eq_tangent_subbundle} \Delta_p := \Delta \cap T_p M= \Span \{ X_1(p), \ldots, X_m(p) \}, \qquad \forall p\in U.\end{equation}
Distributions of tangent subspaces are also simply referred to as {\em distributions}.
Furthermore, if there exists $r\in \N$ such that $r=\dim \Delta_p$,
for all $p\in M$, then we say that $\Delta$ has {\em constant rank} with {\em rank} equal to $r$.\index{rank} 
Distributions of rank $r$ are also called {\em distributions of $r$-planes} or {\em $r$-plane fields}.
Constant rank distributions are also called {\em polarizations} or 
{\em tangent subbundles}.
The pair $(M, \Delta)$ of a manifold $M$ and a polarization $\Delta$ is called {\em polarized manifold}, and $\Delta$ is referred to as the {\em horizontal subbundle} of the polarized manifold. 
 \index{field of distributions} \index{tangent! -- subbundle} \index{distribution}\index{horizontal! -- subbundle}\index{polarization|see {distribution}}
\end{definition}
Notice that each tangent subbundle is indeed a subbundle of the tangent bundle: 
A \emph{subbundle} $E$\index{subbundle} of a vector bundle $F$ (see Section~\ref{Vector bundles}) over a manifold $M$ is a collection of linear subspaces $E_p$ of the fibers $F_p$ of $F$ at each point $p$ in $M$ that forms a vector bundle in its own right. 
 In particular, a tangent subbundle of rank $r$ on an $n$-manifold is a manifold of dimension $n+r$.

Here is a simple example of a polarization on the 3-dimensional manifold $\R^3$, with coordinates $x,y,z$. 
Let $f, g: \R^3\to \R$ be smooth functions.
Then, the two smooth vector fields
\begin{eqnarray*}
X_1(x,y,z) &:=& \partial_x + f(x,y,z) \partial_z, \\
X_2(x,y,z) &:=& \partial_y + g(x,y,z) \partial_z 
\end{eqnarray*} 
are linearly independent at every point $(x,y,z)$ and define a rank-2 tangent subbundle $\Delta$ on $\R^3$ as
\begin{eqnarray*}
\Delta_{(x,y,z)} &:=& \{ a X_1 (x,y,z) + b X_2 (x,y,z) : a,b\in \R\}\\
&=& \{ \left(a, b, a f(x,y,z) +b g(x,y,z)\right) : a,b\in \R\}.
\end{eqnarray*} 

\begin{definition}
Here is some notation and terminology that is commonly used for distributions and families of vector fields:\index{$\Gamma(\Delta)$}
\begin{itemize}
\item The set of smooth vector fields on a manifold $M$ is denoted by $\Vec(M)$ or $\Gamma(TM)$. In fact, an element of $\Gamma(TM)$ is a smooth section $X: M\to TM$ of the bundle $TM\to M$. \index{tangent! vector field -- to a distribution}
\item A vector field $X:M\to TM$ is said to be {\em tangent} to a distribution $\Delta\subseteq TM$ at a point $p\in M$ 
 if $X(p)\in \Delta$.
\item 	Given a distribution $\Delta\subset TM$, we denote by $\Gamma(\Delta)$ 
the set of smooth
vector fields of $M$ tangent to $\Delta $ at every point of $M$.
	\item
	\index{$\scr F_p$}
		Given a family $\scr F\subset\Gamma(TM)$ of vector fields on $M$ and $p\in M$, we set $\scr F_p:=\large(F\large)_p:=\{X_p:X\in\scr F\}$.
\item 
\index{$\Lie(\scr F)$}
	Given a family $\scr F\subset\Gamma(TM)$ of vector fields on $M$, we denote by $\Lie(\scr F)$ the Lie algebra generated by $\scr F$ with respect to the Lie bracket of vector fields within $\Gamma(TM)$; see Section~\ref{Sec:vect_fields_brackets}. 
\end{itemize}
\end{definition}	

We specify that the set $\Lie(\scr F)$ is the smallest subset of $\Gamma(TM)$ containing $\scr F $ and satisfying the property 
 $$ X, Y\in \Lie(\scr F), a,b\in \R \implies [X,Y], aX+bY \in \Lie(\scr F).$$

We are now prepared to introduce a criterion on a polarization $\Delta$ that allows us to connect points with curves tangent to $\Delta$.
The following condition \eqref{Bracket: generating} goes by many names, including {\em H\"ormander's condition} or {\em Chow's condition}.\index{Chow's condition}\index{H\"ormander's condition}

\begin{definition}[Bracket generating]
 	A distribution $\Delta$ on a manifold $M$ is \emph{bracket generating} \index{bracket-generating} if
	\begin{equation}
\label{Bracket: generating}
	(\Lie(\Gamma(\Delta)))_p = T_pM, \qquad \forall p\in M.
	\end{equation}
\end{definition}

Next, we clarify the meaning of a curve tangent to a distribution:
\begin{definition}[Horizontal curve]\label{def horizontal curve}
 	Given a polarized manifold $(M, \Delta)$, 
	a curve $\gamma: [a,b]\to M$ is said to be 
	{\em $\Delta$-horizontal}\index{horizontal! -- curve}\index{curve! horizontal --}
	 if $\gamma$ is absolutely continuous (see Definition~\ref{def:AC})
	 and $\dot\gamma(t)\in\Delta_{\gamma(t)}$ for almost every $t\in[a,b]$.
	 Curves that are $\Delta$-horizontal are also said to be {\em horizontal with respect to} $\Delta$, or, simply, {\em horizontal} or {\em Legendrian}. 
	The terms {\em admissible curve} and {\em controlled path}\index{admissible! -- path}\index{controlled path}
	are also used to refer to such curves.
\end{definition}

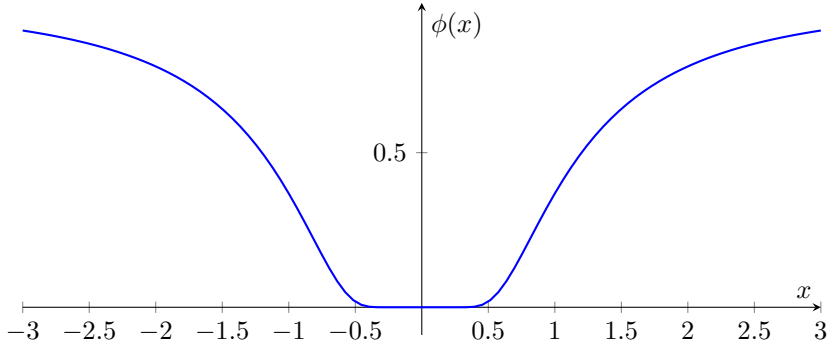
\begin{figure}[ht]
\centering
\begin{tikzpicture}
\begin{axis}[
 xlabel=$x$,
 ylabel=$\phi(x)$,
 domain=-3:3,
 samples=100,
 width=0.8\textwidth,
 height=6cm,
 axis lines=middle,
 enlarge y limits=true,
 clip=false
]
\addplot[blue, thick] {exp(-1/(x^2))};
\end{axis}
\end{tikzpicture}
\caption{
Plot of the function $\phi(x) = \exp\left(-\frac{1}{x^2}\right)$, with $\phi(0)=0$.
This function is smooth everywhere and has all derivatives equal to zero at 
$x=0$, but it is not identically zero.
It is commonly used to define the `Gaussian bump function'. 
}\label{fig:bump}
\end{figure}

\begin{remark} 
Special attention should be paid when verifying the condition \eqref{Bracket: generating} using a frame of vector fields.
Indeed, let $X_1, \ldots, X_m$ be vector fields spanning a distribution $\Delta$ on a manifold $M$, in the sense that \eqref{eq_tangent_subbundle} holds for all $p\in M$.
On the one hand, if 
	\begin{equation}
\label{Bracket: generating: frame}
	(\Lie(\{X_1, \ldots, X_m \}))_p = T_pM, \qquad \forall p\in M,
	\end{equation}
	 then $\Delta$ is bracket generating.
On the other hand, the converse implication may not hold: For instance, 
consider $C^\infty$ function $\phi:\R\to \R$ such that
$\phi(0)=0$ if and only if $x=0$, and
 $\frac{\dd^k}{\dd x^k}\phi(0)=0$, for all $k\in \N$, as shown in Figure~\ref{fig:bump} for an example. 
Consider on $\R^2$ with coordinates $(x,y)$ the vector fields 
$$X:= \partial_x \qquad \text{ and } \qquad Y:= \phi(x)\partial_y.$$
 Although \(X,Y\) do not satisfy \eqref{Bracket: generating: frame}, as demonstrated in Exercise~\ref{ex_Grushin_bad_frame}, they span the same distribution as the bracket-generating frame \(\partial_x, x\partial_y\), as shown in Exercise~\ref{ex_Grushin}.
\end{remark}

\subsection{Sub-Finsler structures of constant rank}
\begin{definition}[Sub-Finsler and sub-Riemannian manifolds of constant rank]\label{def_constant_rank_sub-Finsler}
 		 \index{sub-Finsler! -- manifold} \index{manifold! sub-Finsler --|see {sub-Finsler manifold}}
	A \emph{sub-Finsler manifold}
 is a triple $(M, \Delta, \|\cdot\|)$ where ${M}$ is a connected manifold,
$ \|\cdot\|$ is a continuously varying norm (recall Definition~\ref{def: continuously varying norm}), and
 $\Delta$ is a bracket-generating polarization on $M$, with the rank of $\Delta$ assumed to be constant. 	The pair $(\Delta, \|\cdot\|)$ is said to be a \emph{sub-Finsler structure} on $M$.
	 \index{sub-Finsler! -- structure}	
	If the norm $ \norm{\cdot}$ is given by a Riemannian scalar product $\langle \cdot, \cdot \rangle$, then $(M, \Delta, \langle \cdot, \cdot \rangle )$ is called a {\em sub-Riemannian manifold}.
		 \index{sub-Riemannian! -- manifold} \index{manifold! sub-Riemannian --|see {sub-Riemannian manifold}}
\end{definition}
		 \index{Riemannian! -- manifold} \index{manifold! Riemannian --|see {Riemannian manifold}}
		 		 \index{Finsler! -- manifold} \index{manifold! Finsler --|see {Finsler manifold}}
	 We consider Riemannian and Finsler manifolds as particular cases of sub-Riemannian and sub-Finsler manifolds, respectively, which occur when $\Delta$ is the entire tangent bundle.

Since only the values of the restriction $\left.|\cdot|\right|\Delta$ of $|\cdot|$ to $\Delta$ are important in what follows, we sometimes state that $(M, \Delta,|\cdot||_\Delta)$ is a sub-Finsler manifold with the sub-Finsler structure $(\Delta,|\cdot||_\Delta)$. Specifically, the length with respect to $|\cdot|$, as defined in \eqref{def_Finsler_length}, is considered only for curves that are horizontal with respect to $\Delta$.\index{length}


\begin{definition}[CC-distance]\label{def_cc_dist_mfd}
	\index{CC-distance|see {Carnot-Carath\'eodory distance}}
	\index{Carnot-Carath\'eodory! -- distance}
	\index{distance! Carnot-Carath\'eodory --}
	Given a sub-Finsler manifold $(M, \Delta, \|\cdot\|)$, the {\em Carnot-Carath\'eodory distance}
	between two points $p,q\in M$ is defined as:
	\begin{equation}\label{dist_CC}
	\dcc(p,q) := \inf\left\{\Length _{\|\cdot\|}(\gamma):\gamma\text{ is $\Delta$-horizontal curve from $p$ to $q$}\right\}.
		\end{equation}
	If the infimum is achieved by a curve $\gamma$, that is, if $\dcc(p,q) = \Length _{\|\cdot\|}(\gamma)$, then 
$\gamma$ is {\em length minimizing} among the horizontal curves joining the two points $p$ and $q$. 
The distance function $\dcc$ is also called {\em Carnot-Carath\'eodory metric}.\index{Carnot-Carath\'eodory! -- metric}	
\end{definition}

For a sub-Finsler manifold $(M, \Delta,|\cdot|)$, we also consider the associated Finsler distance.
If $d_F:= d_{\|\cdot\|}$ denotes the Finsler distance associated with $(M, \|\cdot\|)$ as defined in \eqref{def_Finsler_distance}, then it is evident that:\begin{equation}\label{fin<subfin}	
 d_F(p,q) \leq \dcc(p,q), \qquad \forall p,q\in M,
\end{equation}
because in the definition of $\dcc$, we infimize over a subset of the set used for $ d_F$. 
It is important to note that the same length functional is used to define both distances.

We anticipate that the aforementioned function $\dcc$ is indeed a finite distance. In fact, because $\Delta$ is assumed bracket generating and $M$ is assumed connected, we shall show the following result.
\begin{theorem}[Chow; see Section~\ref{sec:proof_chow}]\label{Chow2}
\index{Chow Theorem}
 	If $(M, \Delta, \|\cdot\|)$ is a sub-Finsler manifold, then $\dcc$ is finite and induces the manifold topology on $M$.
\end{theorem}
\begin{Rem}[Terminology]
\index{Carnot-Carath\'eodory! -- space}
			\index{distance! sub-Finsler --}	
		\index{sub-Finsler! -- distance}	
 	The Carnot-Carath\'eodory distance is sometimes referred to as {\em CC-distance} or {\em sub-Finsler distance}.
	A sub-Finsler manifold equipped with its Carnot-Carath\'eodory distance is called {\em Carnot-Carath\'eodory space}.
			If $\|\cdot\|$ is the norm coming from a Riemannian metric, and hence $(M, \Delta, \|\cdot\|)$ is a {sub-Riemannian manifold}, then $(\Delta, \|\cdot\|)$ is called a \emph{sub-Riemannian structure} and $\dcc$ is called \emph{sub-Riemannian distance}.
	 \index{sub-Riemannian! -- manifold}
	\index{sub-Riemannian! -- structure}
	\index{distance! sub-Riemannian --}	
		\index{sub-Riemannian! -- distance}	
		\end{Rem}

%
%
%


	\index{Finsler-Carnot-Carath\'eodory distance}\index{FCC distance! see {sub-Finsler distance}}
Some authors, ourselves included, refer to $\dcc$ as a {\em Finsler-Carnot-Carath\'eodory distance}, or {\em FCC distance}, for short, 	to highlight that in the context $\dcc$ might not necessarily be sub-Riemannian.
Sub-Riemannian metrics have been discussed in the literature under a variety of names, such as `singular Riemannian metric' or `non-holonomic Riemannian metric'. 
 They were also considered in the theory of hypoelliptic PDEs, though without a specific designation.
\index{singular Riemannian metric} \index{non-holonomic Riemannian metric}

\subsection{Control theory viewpoint}
In control theory, the focus lies on systems of differential equations of the form:
\begin{equation}
 \label{control-sys}\dot \gamma=\sum_{j=1}^m c_j(t) X_j{(\gamma)},
\end{equation}
where $X_1, \ldots, X_m$ are predetermined vector fields on a manifold $M$, and $c_1, \ldots,c_m$ are variable $L^1$ functions defined on a bounded interval. 
These functions are called {\em control functions} or {\em controls}. 
Paths obtained by integrating (\ref{control-sys}) are termed {\em controlled paths}.

When the rank of the system of vector fields $X_1, \ldots, X_m$ is constant, controlled paths coincide with the absolutely continuous paths that are tangent to the distribution $\Delta$ generated by $X_1, \ldots, X_m$ as
\begin{equation}\label{frame_Delta}
\Delta_p:= \Span_\R\left\{ X_1(p), \ldots, X_m(p) \right\}, \qquad \text{ for } p\in M.\end{equation}
 
Conversely, every rank-$m$ distribution $\Delta$ can, locally, be expressed as in \eqref{frame_Delta}. 
It is important to note that the adverb `locally' is necessary due to global topological constraints, for instance, for the tangent bundle $\Delta=T(\mathbb S^2)$ of the 2D sphere $\mathbb S^2$.

However, in many systems of interest in control theory, the vector fields $X_1, \ldots, X_m$ are not linearly independent at every point, and the distribution that they define does not have a constant rank. Nevertheless, a related distance can still be defined: for $p\in M$ and $v\in T_p M$, set
$$g_p(v):=\inf\left\{u^2_1+\cdots+u^2_m\;|\; u_1, \ldots,u_m\in \R, u_1 X_1(p)+\cdots+u_m X_m(p)=v\right\}.$$
We are using the convention that $\inf \emptyset =+\infty$.
We then have that $g_p$ is a positive-definite quadratic form on the subspace $$\Delta_p:=\Span_\R\left\{ X_1(p), \ldots, X_m(p) \right\}.$$
The {\em control distance associated with the system} $X_1, \ldots, X_m$ is defined for every $p$ and $q$ in $M$ as:\index{control distance}
\begin{equation} \label{d_control}
d(p,q):=\inf\left\{ \int_0^1 g_p(\dot\gamma(t))^{1/2}\dd t\;\Big{|}\; \gamma \text{ absolutely continuous path with } \gamma(0)=p, \gamma(1)=q\right\}.
\end{equation}

\subsection{The general definition with varying rank}\label{def_CC_varying_rank}
 The Carnot-Carath\'eodory distance (\ref{dist_CC}) and the control distance (\ref{d_control}) fit into a broader context. Specifically, using the language of vector bundles, we can provide a more general definition. The point here is that the distribution on the manifold is obtained as the image of a bundle, and not only the rank of the bundle can be strictly larger than the distribution, but also the rank of the distribution can be different at different points of the manifold.
 
\begin{definition}\label{rank-varying_structure}\index{possibly rank-varying sub-Finsler structure}\index{CC-bundle structure}\index{sub-Finsler! possibly rank-varying -- structure}
 A {\em CC-bundle structure}, also called {\em possibly rank-varying sub-Finsler structure}, on a manifold $M$ is
 a pair $(\sigma, N)$ of functions $\sigma:E\to TM$ and $N:E\to \R$
 where $E$ is a vector bundle over $M$, 
 $N$ is a continuous function such that for all $p\in M$, the restriction of $N$ to the fiber $E_p$ is a symmetric norm (refer to Definition~\ref{def: continuously varying norm}), and $\sigma$ is a smooth map that is a morphism of vector bundles lifting the identity, i.e., the following diagram commutes:\index{morphism! -- of vector bundles}
\[
\xymatrix{
E \ar[dr]_{\pi} \ar[rr]^\sigma & & \ar[dl]^{\pi} TM\\
&M 
		}
\]
and $\sigma|_{E_p}$ is a linear map from $E_p$ to $T_pM$.

For every such a CC-bundle structure $(\sigma, N)$, 
we set
$$\norm{v}:=\inf\left\{N(u)\,:\, u\in E_p, \sigma(u)=v\right\}, \qquad \forall p\in M, \forall v\in T_p M.$$
Analogously as before, the {\em sub-Finsler distance associated with the CC-bundle structure}, for every $p$ and $q$ in $M$, is defined as:
\begin{equation}\label{def_d_25may}
d(p,q):=\inf\left\{ \int_0^1 \norm{\dot\gamma(t)}\dd t\;\Big{|}\; \gamma \text{ absolutely continuous path with } \gamma(0)=p, \gamma(1)=q\right\}.
\end{equation}
\end{definition}

One can verify that, in the case of the inclusion $\sigma:\Delta \hookrightarrow T M$ of a subbundle of the tangent bundle, one recovers the Carnot-Carath\'eodory distance (\ref{dist_CC}). Similarly, for $E:=M\times\R^m$ and $\sigma(p,u):= u_1 X_1+\cdots+u_m X_m$, one recovers the control distance (\ref{d_control}); see also Exercise~\ref{ex_20may1006}.

\subsection{Equiregular distributions}\label{sec:equiregular}

Let $\Delta\subset TM$ be a subbundle.
For every $p\in M$ we define 
\begin{align*}
 	\Delta^{[0]}(p) &:= \{0\}\subset T_pM, \\
	\Delta^{[1]}(p) &:= \Delta_p , \\
	\Delta^{[2]}(p) &:= \Delta^{[1]}(p) + \Span\left\{[X,Y]_p:X,Y\in\Gamma(\Delta)\right\} .
\end{align*}
Then $\Delta^{[2]} := \bigcup_{p\in M}\Delta^{[2]}(p)$ is a subset of $TM$.
In general, the subset $\Delta^{[2]}$ may not be a subbundle since its rank may vary, i.e., the function $p\mapsto\dim\Delta^{[2]}(p)$ may not be constant.

\begin{example}[Non-equiregular distribution]\index{Martinet distribution}
	In $\R^3$ the \emph{Martinet distribution} is the subbundle $\Delta\subset T\R^3$ spanned by the vector fields:
$$X_1 = \partial_x+\frac {y^2}2 \partial_z \qquad \text{and}\qquad
		X_2 = \partial_y .$$
	By computing the Lie brackets, we find:
	\[
	X_3:= [X_2, X_1] = y\partial_z
	\qquad\text{and}\qquad
	X_4:= [X_2, X_3] = \partial_z .
	\]
	In this case, the subset $\Delta^{[2]}$ has fibers of varying dimension, since we have
	\[
	\Delta^{[2]}(p) = 
	\begin{cases}
	 	T_p\R^3 &\text{if } p_2\neq0, \\
		\Delta^{[1]}(p) &\text{if }p_2=0 .
	\end{cases}
	\]
\end{example}

\begin{Rem}
 	If $X_1, \dots, X_r$ is a frame for $\Delta$, then the collection of vector fields
	\[
	\left\{X_1, \dots, X_r\}\cup\{[X_i, X_j]:i,j\in\{1, \dots,r\}\right\}
	\]
	span $\Delta^{[2]}$ at every point.
	Indeed, if $X,Y\in\Gamma(\Delta)$, then $X=\sum_ia^iX_i$ and $Y=\sum_jb^jX_j$ for some smooth functions $a^i,b^j$.
	Thus, we have what we claimed:
$$
	 	[X,Y] = [a^iX_i,b^jX_j] = a^ib^j[X_i, X_j] + a^i (X_ib^j) X_j - b^j (X_ja^i)X_i .
	$$
\end{Rem}
\begin{definition}[$\Delta^{[k ]}$]\label{def:Delta_k}
\index{$\Delta^{[k]}$}
 Given a distribution $\Delta\subseteq TM$ on $M$, for each $k\in\N$ we shall define the subset $
\Delta^{[k ]}\subseteq TM$ by describing each of its
 fiber $\Delta^{[k ]}(p):= \Delta^{[k ]}\cap T_pM$ as $p$ varies in $M $. The fiber $\Delta^{[k ]}(p)$ is given by
\begin{equation}\label{Deltaj}
\Delta^{[k ]}(p) := 
\Span\left\{[X_1,[X_2, \dots,[X_{j-1}, X_{j}]\dots]](p)\,:\, j\in\{1, \ldots, k\}, \, X_1, \dots, X_{j}\in\Gamma(\Delta)\right\}.
\end{equation}
\end{definition}

The sets $\Delta^{[k ]}(p)$ can also be defined inductively by $\Delta^{[1]}=\Delta$ and, for all $k\geq 2$,
\begin{equation}
\Delta^{[k+1]}(p) = \Delta^{[k]}(p) + \Span\left\{[X_1,[X_2, \dots,[X_{k}, X_{k+1}]\dots]](p)\,:\, X_1, \dots, X_{k+1}\in\Gamma(\Delta)\right\} .
\end{equation}
\begin{definition}[Regular point for $\Delta$]\index{regular point for distribution}
 	If $\Delta$ is a distribution on $M$ and $p\in M$, we say that $p$ is \emph{regular} for $\Delta$ if for all $k\in\N$ the function 
	\begin{equation}\label{eq1243}
	q\longmapsto\dim\Delta^{[k]}(q)
	\end{equation} 
	is constant in a neighborhood of $p$.
\end{definition}
Notice that the functions \eqref{eq1243} is $\N$-valued. 
Hence, if it is locally constant, then it is constant on connected components.
\begin{definition}[Equiregular distributions]\label{def:equiregular}
 	Let $M$ be a manifold. A distribution $\Delta\subset TM$ is said to be \emph{equiregular} if every $p\in M$ is regular for $\Delta$.
	In this case, we call $(\Delta^{[k]})_{k\in \N}$, as in Definition~\ref{def:Delta_k}, the {\em flag of subbundles} for $\Delta$. \index{equiregular! -- distribution}\index{flag of subbundles}
\end{definition}
\begin{Rem} A distribution 
 	$\Delta\subset TM$ is equiregular if and only if, for all $k\in\N$, the set $\Delta^{[k]}$ is a subbundle.
\end{Rem}
Notice that if $\Delta$ is bracket generating and equiregular, then there is $s\in\N$ such that $\Delta^{[s]}=TM$.
The minimal such an $s$ is called \emph{step} of $\Delta$.
\begin{definition}[Equiregular sub-Finsler manifolds]
	A sub-Finsler manifold $(M, \Delta, \|\cdot\|)$ is called \emph{equiregular} if $\Delta$ is equiregular. 	
\end{definition}


\section{Chow's theorem and existence of geodesics}\label{sed_Chow_geodesics}
In this section, we explain how bracket-generating distributions allow the existence of horizontal curves connecting arbitrary points. Consequently, Carnot-Carath\'eodory distances on sub-Finsler manifolds do not alter the topology. Moreover, locally, sub-Finsler manifolds are geodesic spaces.
\subsection{Local transitivity and Sussmann's orbit theorem}
In this section, we highlight the fact that, because in every sub-Finsler manifold, the distribution is assumed bracket generating, then the Carnot-Carath\'eodory distance is finite.
The bracket-generating condition can be considered an infinitesimal form of transitivity. 
Chow's theorem states that this condition implies local transitivity:
\begin{theorem}[Chow]\label{Chow}
 If a subbundle $\Delta$ of the tangent bundle of a manifold is bracket generating at some point $p$ (i.e., \eqref{Bracket: generating} holds at $p$), then every point $q$ that is sufficiently close to $p$ can be joined to $p$ by an absolutely continuous curve almost everywhere tangent to $\Delta$.
\end{theorem} 
In fact, nearby points in a sub-Finsler manifold can be joined by horizontal curves with small Finsler length. 
This is precisely what Theorem~\ref{Chow2} asserts.

We first explain the validity of Theorem~\ref{Chow} taking for granted a theorem by Sussmann.
 We are omitting the proof of Sussmann's theorem, which is, in fact, the core of our first proof of Theorem~\ref{Chow}, but it is well presented in \cite{bellaiche}. 
The reader can write a second complete proof of the above Theorem~\ref{Chow} by following the hints in Exercise~\ref{proof: Chow}.
Later in the text, we will present a detailed proof of Theorem~\ref{Chow2}, a result of higher interest for us.
Also, the simpler case of Carnot groups, discussed in Section~\ref{effective_Chow_Carnot}, offers an elementary demonstration of Theorem~\ref{Chow2}.

\begin{theorem}[Sussmann \cite{Sussmann,Stefan,bellaiche}]
\label{thm_Sussmann}
Let $M$ be a manifold, $\Delta\subseteq TM$ a subbundle, and $p\in M$.
Let $\Sigma\subset M$ be the set of points that can be joined to $p$ by an absolutely continuous curve almost everywhere tangent to $\Delta$.
Then, the set $\Sigma$ is an immersed submanifold of $M$. 
\end{theorem}


\begin{proof}[A first proof of Theorem~\ref{Chow}, modulo Theorem~\ref{thm_Sussmann}]
In the assumptions of Theorem~\ref{Chow}, we employ Theorem~\ref{thm_Sussmann}.
Given a vector field $X\in\Gamma(\Delta)$ and a point $q\in \Sigma$, the flow line $t\mapsto \Phi^t_X(q)$ is tangent to $\Delta$ and lies in $\Sigma$. 
Thus, the vector $X_q$ is tangent to the submanifold $\Sigma$.  
Consequently, we have: 
$$\Gamma(\Delta)\subseteq \mathcal F :=\{ X\in \Vec(M) : X_q\in T\Sigma, \forall q\in \Sigma\}.$$
Since $\Sigma$ is a submanifold, the family
$\mathcal F$ is involutive on $\Sigma$, meaning 
$\Lie(\mathcal F) \subseteq \mathcal F$; see Exercise~\ref{Seba 3 oct 2024}.
Therefore, we deduce that $\Lie(\Gamma(\Delta))\subseteq\mathcal F$.
By the bracket-generating condition at $p$, we have
$$T_pM = \Lie(\Gamma(\Delta))_p \subseteq \mathcal F_p \subseteq T_p \Sigma.$$
This implies $\dim M=\dim \Sigma$, and thus $\Sigma$ is a neighborhood of $p$.
\end{proof}



\subsection{Reachable sets of bracket-generating distributions}
\label{sec:reachable}
\index{reachable set}
Let $\scr F\subset\Vec(M)$ be a family of smooth vector fields on a manifold $M$.
Define the \emph{reachable set for $\scr F$ from $p$ at time less than $T$} as
\[
\Phi^{<T}_\scr F(p) :=
\left\{
\Phi_{X_k}^{t_k}\circ\dots\circ\Phi_{X_1}^{t_1}(p):k\in\N,t_j>0, \sum_{j=1}^k t_j<T, X_j\in\scr F
\right\}.
\]

\begin{theorem}\label{thm: reachable}
Let $M$ be a manifold of positive dimension, and $\scr F$ be a family of vector fields on $M$. 
If $-\scr F=\scr F$ and $(\Lie(\scr F))_p=T_pM$ for all $p\in M$, then for all $T>0$ and for all $p\in M$, the set $\Phi^{<T}_{\scr F}(p)$ contains $p$ in its interior.
\end{theorem}
\begin{proof}
	Since $\dim(M)>0$, there is $X_1\in\scr F$ with $X_1(p)\neq0$. 
	Hence, there is $\epsilon_1\in(0,T)$ such that 
	\[
	M_1:= \{ \Phi^t_{X_1}(p):t\in(0, \epsilon_1)\}
	\]
	is a 1-dimensional submanifold of $M$.
	
	If $M$ is $1$-dimensional, the proof is concluded.
	If $\dim M>1$, then there is $X_2\in\scr F$ that is not tangent to $M_1$ (Otherwise, the family $\Lie(\scr F)$ of vector fields would be tangent to $M_1$ and not bracket-generating on points of $M_1$).
	Let $\hat t_1\in(0, \epsilon_1)$ such that 
	\[
	X_2(\Phi_{X_1}^{\hat t_1}(p)) \notin TM_1 .
	\]
	The map $(t_1,t_2)\mapsto \Phi_{X_2}^{t_2}\circ \Phi_{X_1}^{t_1}(p)$ has maximal rank (i.e., rank 2) at every point of the form $(\hat t_1,t_2)$ with $t_2$ sufficiently small, say $t_2\in(0, \epsilon_2)$ with $ t_1<t_1+\epsilon_2<T$. We have obtained an embedded parametrized surface, as shown in Figure~\ref{fig_flows}.

	\begin{figure}[ht]
\centering
\begin{tikzpicture}
\draw[thick] (0,3.1) to[out=-10, in=190] (6,1.7);
\foreach \x/\y/\offset/\curvecolor/\curvature/\label in {
0.5/3/0/blue!50!cyan/0.2/, 
1.5/2.67/0.2/blue!70!green/0.4/, 
2.5/2.3/0.4/blue!90!yellow/0.6/, 
3.5/1.9/0.6/blue/0.8/, 
4.55/1.65/0.8/blue!70!red/1/} {
 \draw[\curvecolor, thick] (\x, \y) to[out=30, in=180-\curvature*90] ++(2+\offset,0.5);
 \node[\curvecolor, anchor=west] at (\x+2.2+\offset, \y+0.3) {\textbf{\label}};
}
\node[circle, fill=black, inner sep=1pt, label=above:$p$] at (0.2,3.05) {};
\node[circle, fill=black, inner sep=1pt, label=below:$\Phi_{X_1}^{\hat t_1}(p)\, \, \,$] at (2.5,2.3) {};
\node[circle, fill=black, inner sep=1pt, label=right:$\, \Phi_{X_2}^{\hat t}(\Phi_{X_1}^{\hat t_1}(p))$] at (4.5,3.1) {};
\draw[thick, ->] (0.2,3.05) -- ++(1.4,{-0.25}) node[pos=0.4,below]{${X_1}$};
\draw[thick, ->] (2.5,2.3) -- ++(1.5,{-0.6}) node[pos=0.7,below]{${X_1}$};
\draw[thick, ->] (5.3,1.62) -- ++(1.5,{0.1}) node[pos=0.6,below]{${X_1}$};
\draw[blue!70!green, thick, ->] (1.5,2.67) -- ++(1.4,{0.9}) node[pos=0.3,above]{${X_2}\,$};
\draw[blue!90!yellow, thick, ->] (2.5,2.3) -- ++(1.3,{0.9}) node[pos=0.34,above]{${X_2}\,$};
\draw[blue, thick, ->] (3.5,1.9) -- ++(1.3,{0.8}) node[pos=0.4,above]{${X_2}\,$};
	\end{tikzpicture}
\caption{Composition of two flows to construct a surface within the reachable set of a point $p$: first, 
flowing along the vector field $X_1$, followed by flowing along $X_2$, which is transverse to $X_1$.}\label{fig_flows}
\end{figure}
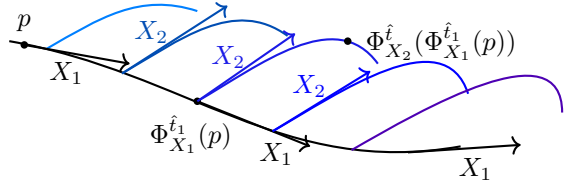
	
	Proceeding in this manner, for all $k$ with $k\leq\dim (M)$, we obtain vector fields $X_1, \dots, X_k\in\scr F$ such that the map
	\[
	F_k:(t_1, \dots,t_k)\longmapsto \Phi_{X_k}^{t_k}\circ\dots\circ\Phi_{X_1}^{t_1}(p)
	\]
has maximal rank $k$ at some point $(\hat t_1, \dots, \hat t_k)$ with $\hat t_j>0$ and $\sum_j\hat t_j<T$.
By the Constant-Rank Theorem, there exists a neighborhood $U_k$ of $(\hat t_1, \dots, \hat t_k)$ such that $M_k:=F_k(U_k)$ is an embedded submanifold.
	
	This procedure stops precisely when each element of $\scr F$ is tangent to $M_k$, i.e., when $M_k$ becomes an open subset of $M$.
Take $X_1, \dots, X_k\in\scr F$ such that the previously defined points $F_k(t_1, \dots,t_k)$ cover a neighborhood of a point $q\in M$ when varying $t_j>0$ with $\sum_jt_j<T$.
Notice that if $q$ is of the form $F_k(\bar t_1, \dots, \bar t_k)$, with $\bar t_j>0$ and $\sum_j \bar t_j<T$, then the map 
	\[
	q'\longmapsto \Phi^{\bar t_1}_{-X_1}\circ\dots\circ\Phi^{\bar t_k}_{-X_k}(q')
	\]
	is a diffeomorphism between some neighborhood of $q$ and its image, which is a neighborhood of $p$.
Notice that $-X_j\in-\scr F=\scr F$ by assumption.
Therefore
	\[
	(t_1, \dots,t_k)\longmapsto\Phi_{-X_1}^{\bar t_1}\circ\dots\circ\Phi_{-X_k}^{\bar t_k}\circ\Phi_{X_k}^{t_k}\circ\dots\circ\Phi_{X_1}^{t_1}(p)
	\]
	covers a neighborhood of $p$ when $t_j>0$ and $\sum_jt_j<T$.
	Thus $\Phi_{\scr F}^{<2T}(p)$ is a neighborhood of $p$.
\end{proof}

With the use of Theorem~\ref{thm: reachable}, one can provide an alternative proof of Theorem~\ref{Chow}; see Exercise~\ref{proof: Chow}.

\subsection{The metric version of Chow's theorem}\label{sec:proof_chow}\index{Chow Theorem}\index{Theorem! Chow --}
We are now prepared to prove Theorem~\ref{Chow2}. 
Namely, we show that Carnot-Carath\'eodory distances induce the manifold topology.

\begin{proof}[Proof of Theorem~\ref{Chow2}]
Let $\tau_M$ denote the manifold topology and $\tau_{CC}$ the topology induced by $\dcc$. We shall prove $\tau_{CC}=\tau_M$ by establishing the two containment containments.

Regarding the containment {$\tau_{CC}\subset\tau_M$}, let $U\in\tau_{CC}$ and $p\in U$.
By definition of $\tau_{CC}$, there exists $T>0$ such that $B_{\dcc}(p,T)\subset U$.
Set 
\[
\scr F := \{X\in \Gamma(\Delta):\|X(p)\|\le 1\ \forall p\in M\} \subset\Vec(M).
\]
With the notation of Section~\ref{sec:reachable},
observe that
\[
 \Phi^{<T}_{\scr F}(p)\subset B_{\dcc}(p,T).
\]
By Theorem~\ref{thm: reachable}, the point $p$ is in the $\tau_M$-interior of $ \Phi^{<T}_{\scr F}(p)$.
We deduce that $p$ is also in the $\tau_M$-interior of $U$.

Regarding the containment {$\tau_M\subset \tau_{CC}$}, let $U\in\tau_M$.
Together with the distance $\dcc$, we have a Finsler distance $d_F$ that satisfies \eqref{fin<subfin}.
Let $p\in U$. Then there is $r$ such that $B_{d_F}(p,r)\subset U$ because Finsler distances induce the manifold topology.
Since $d_F\le\dcc$, then $B_{\dcc}(p,r)\subset B_{{d_F}}(p,r)$.
Therefore $p$ is also in the $\tau_{CC}$-interior of $U$.
\end{proof}


\subsection{Comparison of length structures}\index{length! -- structure}
In certain situations, dealing with the case where $\Delta$ is a distribution with varying rank becomes more challenging. For instance, the following proposition remains valid when 
$(M, \Delta, \|\cdot\|)$ is a sub-Finsler manifold in the sense of Definition~\ref{rank-varying_structure}, as discussed in the work \cite{Antonelli_LeDonne_MR4645068}.
However, it is crucial to note that altering Definition~\ref{def_constant_rank_sub-Finsler} by considering rank-varying distributions (instead of polarizations) with norms defined on the entire tangent bundle, as in Definition~\ref{def_constant_rank_sub-Finsler}, may lead to the falsification of the conclusion of following proposition. In fact, in that setting, there are examples of smooth curves parametrized by arc length that are nowhere tangent to the distribution, as illustrated in Example~\ref{Grushin_bad_metric}.
Due to this issue, we confine our attention to sub-Finsler structures of constant rank, and we shall prove the following proposition exclusively in accordance with Definition~\ref{def_constant_rank_sub-Finsler}.
 
\begin{proposition}\label{prop_compare_length_structures}
 Let $(M, \Delta, \|\cdot\|)$ be a (constant-rank) sub-Finsler manifold equipped with its Carnot-Carath\'eodory distance $\dcc$.
Let $\gamma: [a,b]\to M$ be a curve.
	\begin{description}
\item[\ref{prop_compare_length_structures}.i.] 	
If $\Length_{\dcc}(\gamma)<\infty$, then the reparametrization by arc length of $\gamma$ is $\Delta$-horizontal.
\item[\ref{prop_compare_length_structures}.ii.] 	
If $\gamma$ is $\Delta$-horizontal, then $\Length_{\dcc}(\gamma)=\Length_{\|\cdot\|}(\gamma)$; and $\gamma$ is parametrized by arc length if and only if $\norm{\dot\gamma }=1$ almost everywhere.
	\end{description}
\end{proposition}
\proof[Proof of \ref{prop_compare_length_structures}.i.]
 Recall that in every metric space, every curve of finite length can be reparametrized by arc length; see Exercise~\ref{rect: PBAL}. 
Hence, we consider $\gamma$ to be parametrized by arc length with respect to the distance $\dcc$.	
 
Let $d_F$ be the Finsler distance for which we have \eqref{fin<subfin}.
Recall that $d_F$ is locally biLipschitz equivalent to every Riemannian distance, as stated in Proposition~\ref{prop_Finsler_Riem_top}.
Since $d_F\le \dcc$, we have: 
\begin{equation}\label{vaiubvieaurfho38}
	d_F(\gamma(s), \gamma(t)) \le \dcc(\gamma(s), \gamma(t)) \le \Length_{\dcc}(\gamma|_{[s,t]}) = |t-s|.
\end{equation}
This implies that $\gamma: [a,b]\to(M,d_F)$ is 1-Lipschitz. 
Consequently, in local coordinates, the curve $\gamma$ is (Euclidean) Lipschitz.
By the Rademacher Theorem, the curve $\gamma$ is absolutely continuous and hence differentiable almost everywhere.
	Let $t_0\in I$ be a point of differentiability for $\gamma$.
We shall prove that $\dot\gamma(t_0)\in\Delta_{\gamma(t_0)}$.
	
Assume, by contradiction, that $\dot\gamma(t_0)\notin\Delta_{\gamma(t_0)}$.
For simplicity, we work in coordinates and assume $t_0=0$, $\gamma(t_0)=0\in\R^n$, $\Delta_0=\R^k\times\{0\}^{n-k}$, and $\dot\gamma(t_0)=e_n=(0, \dots,0,1)$.
We then have: 
	\begin{equation}\label{ciao_2021}
	\gamma_n(t)> t/2, \qquad \text{ for small enough } t, \end{equation} where $\gamma_n(t)$ is the $n$-th component of $\gamma$.
	
	We claim that for all $\epsilon>0$ there exists $r_\epsilon>0$ such that: 
	\begin{equation}\label{marc_2021}
	p\in B_{d_F}(0,2 r_\epsilon), \, X
	\in\Delta_p, \,
	\|X\|\le 1 \implies| \langle \partial_n, X \rangle|<\epsilon, \end{equation}
	where we use the Euclidean scalar product making $\partial_1, \ldots, \partial_n$ orthonormal. Indeed, by contradiction, assume there exists $\epsilon>0$ and sequences $(p_j)_j $ in $M$ and $(X_j)_j$ in $TM$ with $X_j\in \Delta_{p_j}$ such that $p_j\to 0$, $\|X_j\|\le 1$,
	and 
	$ |\langle \partial_n, X_j\rangle| \geq\epsilon$. 
Let $c>0$ be a constant for which we have \eqref{equiv_to_Euclid} in some neighborhood of 0. Hence, eventually we have $\|X_j\|_{\mathbb E}\le c$. 
Therefore, being the sequence $X_j$ in a compact set, it converges to some $Y$, up to a subsequence.
Since $\Delta$ is assumed to be a polarization (hence a subbundle), it is a closed subset of $TM$; see Exercise~\ref{ex_subbundle_closed}. And since	$p_j\to 0$, we have that $Y\in \Delta_0$. Thus we get 	$$0=|\langle \partial_n, Y\rangle |= \lim_{j\to\infty}|\langle \partial_n, X_j\rangle | \geq\eps>0.$$
We inferred a contradiction, which gives the claim \eqref{marc_2021}.
		
Let $\epsilon>0$ and $r_\epsilon$ be chosen with the above property \eqref{marc_2021}. 
By definition of $\dcc$, we shall take a horizontal curve that almost realizes $\dcc(0, \gamma(r_\eps))$, which is not zero because of \eqref{ciao_2021}. 
In fact, there is a horizontal curve $\sigma: [0,b_\eps]\to M$ from $0$ to $\gamma(r_\eps)$ such that $\|\dot\sigma \|=1$ almost everywhere and $b_\eps=\Length_{\|\cdot\|}(\sigma )\le 2\dcc(0, \gamma(r_\eps))\le 2 r_\eps$, where in the last inequality we used \eqref{vaiubvieaurfho38}.
 Hence, first we have \begin{equation}\label{eq12580}
	\frac{b_{ \epsilon}}{r_\epsilon} 
	\le2,
 	\end{equation}
		 second, we have that 
	 the image of $\sigma $ is in $ B_{d_F}(0,2r_\eps)$.
	 Consequently, because $\sigma $ is horizontal and $\|\dot\sigma \|=1$ almost everywhere, from \eqref{marc_2021} we have that $|\dot \sigma _n|<\epsilon$, where $ \sigma_n$ is the $n$-th component of $ \sigma $, so $\dot \sigma _n = \langle \partial_n, \dot \sigma \rangle$.
	 We then infer that: 
	\[0< \frac{r_\eps} {2}\stackrel{\eqref{ciao_2021}}{<} 
	 \gamma_n(r_\eps)=
	 \sigma_n(b_{ \eps}) 
	= \int_0^{b_{ \eps}} \dot\sigma_n(s) \d s 
	\le \int_0^{b_{\eps}} |\dot\sigma_n (s)| \d s
	\le \epsilon b_\eps .
	\]

	Thus we just obtained a bound that contradicts \eqref{eq12580} for small enough $\epsilon$, since
	\begin{equation*}\label{eq1258}
	\frac{b_{ \epsilon}}{r_\epsilon} 
	\ge \frac1{2\epsilon} \to \infty
	\quad\text{ as }
	\epsilon\to0.
	\end{equation*}
	We deduce that $\gamma$ is horizontal. And the statement \ref{prop_compare_length_structures}.i is proven.
	
\proof[Proof of \ref{prop_compare_length_structures}.ii.] Let $\gamma$ be a horizontal curve. On the one hand, since $d_F\le\dcc$ and since
	$\Length_{\|\cdot\|}= 
	\Length_{d_F}$ by Theorem~\ref{thm: lengthscoincide}, then $\Length_{\|\cdot\|}\le \Length_{\dcc}$.
	On the other hand, since $\gamma$ is horizontal,
	\begin{align*}
	 	\Length_{\dcc}(\gamma) 
		&= \sup_{ (t_1, \dots,t_k)} \sum_{i=1}^{k-1} \dcc(\gamma(t_{i+1}), \gamma(t_i)) \\
		&\le \sup_{ (t_1, \dots,t_k)} \sum_{i=1}^{k-1} \Length_{\|\cdot\|}(\gamma|_{[t_i,t_{i+1}]}) \\
		&= \Length_{\|\cdot\|}(\gamma), 
	\end{align*}
	where the suprema are over all the partitions $(t_1, \dots,t_k)$ of the domain of $\gamma$.
\qed

With the previous Proposition~\ref{prop_compare_length_structures}, we deduce that in CC spaces, the two a priori different length structures coincide. Hence, it is unambiguous what we mean when we refer to them as length spaces, a concept introduced in Section~\ref{sec_length_sp}. 

\begin{corollary}\label{cor:CC:lenght}\index{length! -- space}
 Carnot-Carath\'eodory spaces are length spaces.
\end{corollary}
\begin{proof}
Let $(M, \Delta, \|\cdot\|)$ be a sub-Finsler manifold with its Carnot-Carath\'eodory distance $\dcc$.
We need to show that the distance $\dcc$ is finite and equal to the infimum of the length of curves joining points, as in \eqref{def_length_sp}. Here, the length is calculated with respect to the distance $\dcc$ itself.
First, the distance is finite by Chow's Theorem~\ref{Chow2}.
Second, Proposition~\ref{prop_compare_length_structures} tells us that the curves with finite length with respect to $\dcc$ are, up to reparametrizations, the $\Delta$-horizontal curves. Moreover, their length coincides with the integral of their velocity, as in \eqref{def_Finsler_length}. Therefore, we obtain \eqref{def_length_sp} from the definition \eqref{dist_CC} for $\dcc$.
\end{proof}

\subsection{Existence of geodesics in CC spaces}
The Riemannian analog of the following theorem is due to Heinz Hopf (1894--1971) and Willi Rinow (1907--1979), who independently proved versions of it.
The theorem guarantees the existence of geodesics between every two nearby points. Moreover, the two points can be chosen arbitrarily if the space is boundedly compact in the sense of Section~\ref{sec_length_sp}.\index{existence! -- of geodesics}

\begin{theorem}[Hopf-Rinow Theorem for CC spaces]\label{HRCV_CC}\index{Theorem! Hopf-Rinow --}\index{Hopf-Rinow Theorem}\index{geodesic}\index{locally compact}\index{boundedly compact}
Let $M$ be a CC space.
	\begin{enumerate}
	\item 	Every point in $M$ has a neighborhood in which every two points 
	can be joined by a curve that minimizes length with respect to the CC distance.
	\item 	If $M$ is boundedly compact, then it is a geodesic space.
	\end{enumerate}
\end{theorem}
\begin{proof}
 	By Chow's Theorem~\ref{Chow2}, since $M$ is connected and $\Delta$ is bracket generating, the distance function $\dcc$ is finite and the topology of 
	 $(M, \dcc)$ is locally compact. Moreover, the metric space $(M, \dcc)$ is a length space, as stated in Corollary~\ref{cor:CC:lenght}. 
	 
To find shortest paths, we shall use Proposition~\ref{exists: short}. 
This proposition ensures the existence of a shortest path between every two points that are sufficiently close. In fact, consider a point $p\in M$ and select $r>0$ small enough so that the closed ball $\bar B(p,r)$ is compact. We claim that every two points $p_1, p_2\in B(p,r/2)$ can be joined by a length-minimizing curve. 	 
	Indeed, by Proposition~\ref{exists: short}, there exists a curve $\sigma$ from $p_1$ to $ p_2$ that is one of the shortest among the curves contained in $\bar B(p,r)$. On the one hand, notice that the length of $\sigma$ is at most $r$, the reason being that each of $p_1$ and $ p_2$ can be connected to $p$ via a curve of length strictly less than $r/2$, which therefore remains inside $\bar B(p,r)$. 
	On the other hand, every competitor curve from $p_1$ to $ p_2$ that leaves $\bar B(p,r)$ has a length of at least $r$ because it starts in $ B(p,r/2)$, exit $\bar B(p,r)$, and returns into $ B(p,r/2)$. Therefore, the curve $\sigma$ is a length-minimizing curve. 	 
	  
	If, in addition, the metric space $(M, \dcc)$ is boundedly compact, we can conclude by applying Proposition~\ref{prop:length:geod}
\end{proof} 

\section{Ball-Box Theorem and Hausdorff dimension}

Discussing the Hausdorff dimension is one of the most direct ways to demonstrate that CC spaces are metrically equivalent to Riemannian spaces only when the polarization is the entire tangent bundle. 
 This metric dimension can be calculated by observing that the metric balls, with respect to CC distances, do not necessarily behave like cubes where all the edges have comparable lengths; instead, they behave like boxes with edges of differing magnitudes. Such a statement is made precise by the so-called Ball-Box Theorem.

\subsection{Ball-Box Theorem}
Let $(M, \Delta, \|\cdot\|)$ be an equiregular sub-Finsler manifold of topological dimension $n$ and step $s$. Consider the flag of subbundles:
\[
\Delta=\Delta^{[1]} \subset\Delta^{[2]}\subset\dots\subset\Delta^{[s]}=TM.
\]
Since the upcoming considerations will be local in nature, we assume the existence of a frame $X_1, \dots, X_n$ for $TM$ and numbers $m_1, \dots,m_s$ such that $X_1, \dots, X_{m_k}$ form a frame for $\Delta^{[k]}$.
In this case, we say that $X_1, \dots, X_n$ is an \emph{equiregular frame}. 
Equiregular frames are also known as \emph{adapted frames}.\index{equiregular! -- frame}\index{adapted! -- frame}

Notice that, for all $p\in M$, \begin{equation}\label{eq:mj} m_j=\dim \Delta^{[j]}(p).\end{equation}

We also say that $X_j$ has \emph{degree} $d_j$ if  \index{degree! -- of a vector field}
\begin{equation}\label{eq_def_degree}
X_j(p)\in\Delta^{[d_j]}\setminus\Delta^{[d_j-1]}, \qquad \forall p\in M,
\end{equation}
 i.e., $j\in\{m_{d_{j-1}}+1, \dots,m_{d_j}\}$. We might denote $d_j$ by $\deg(X_j)$.

The plan is to parametrize the manifold $M$ using the flow of linear sums of $X_1, \dots, X_n$. 
To such vector fields, we associate an \emph{exponential coordinate map} from a point $p\in M$ as
\index{exponential! -- coordinates}
\begin{equation}\label{def_phi_n_flows}
\Phi_p:\R^n\to M, \qquad 
(t_1, \dots,t_n)\longmapsto \Phi^1_{t_1X_1+\dots+t_nX_n}(p),
\end{equation}
where $\Phi^1_X(p)$ is the flow of $X$ at time $1$ starting from $p$.
Such a map might be defined only on a neighborhood of $0\in\R^n$. However, for the sake of simplicity and due to the fact that this is the case for Lie groups, we assume that $\Phi_p$ is globally defined.

We define the \emph{box} with respect to the numbers $d_1, \dots,d_n$ and radius $r>0$ as
\begin{equation}\label{def: box} \index{box}
\Bx(r) := \left\{(t_1, \dots,t_n)\in\R^n:|t_j|\le r^{d_j}\right\}.
\end{equation}

The following comparison theorem is due to several mathematicians, including Mitchell, Gershkovich, and Nagel-Stein-Wainger, as mentioned in \cite{Gromov}. It is known as the Ball-Box Theorem because it establishes a comparison between the box ${\rm Box}(r)$ in $\R^n$ and the ball $B(p,r')$ with respect to the CC distance, providing a biLipschitz relation between $r$ and $r'$.

\begin{theorem}[Ball-Box Theorem]\label{Ball-Box}\index{Theorem! Ball-Box --}
 	Let $(M, \Delta, \|\cdot\|)$ be a sub-Finsler manifold.
	Assume $\Delta$ is equiregular.
	Fix $\bar p\in M$ and an equiregular frame $X_1, \dots, X_n$ in a neighborhood of $\bar p$ with degrees $d_1, \dots,d_n$ and the corresponding boxes $\Bx(\cdot)$. 
	Then there exist a neighborhood $U$ of $\bar p$ in $M$ and constants $C>1$ and $\rho>0$ such that
	\[
	B_{\dcc}(p, r/C) \subset \Phi_p(\Bx(r)) \subset B_{\dcc}(p,Cr), \qquad \forall p\in U, \forall r\in(0, \rho).
	\]
\end{theorem}

The proof of the Ball-Box Theorem in this general context of manifolds will not be presented here. Instead, it will be demonstrated later in the more straightforward scenario of Carnot groups and then, more generally, for sub-Finsler groups; see Theorems~\ref{Ball-Box4Carnot} and~\ref{thm66b22d85}.

\begin{remark}
The Ball-Box Theorem~\ref{Ball-Box} provides a quantitative version of Chow's theorems~\ref{Chow} and~\ref{Chow2}.
\end{remark}

As of our current knowledge, there is no conclusive answer to the following natural question, except, perhaps, in the case of step-2 Carnot groups and contact $3$-manifolds – further investigation is required \cite{Baryshnikov_MR1771429, Barilari_MR4118144}.
 \begin{question}[Open!] Are all sufficiently small sub-Finsler balls and spheres homeomorphic to the usual Euclidean balls and spheres?
 \end{question}

Here is a first very useful consequence of the Ball-Box Theorem~\ref{Ball-Box}.
\begin{corollary}
[H\"older equivalence between CC and Euclidean metrics]
Locally, each equiregular sub-Finsler manifold is H\"older equivalent to a Riemannian manifold. 
Namely, if $s$ is the step, then locally around every point, there exists $C>1$ such that 
\begin{equation}\label{eq_BB_thm}
\frac{1}{C}(\dcc )^s \leq d_{\rm Riem}\leq C \dcc .
\end{equation}
\end{corollary}
\proof Let $(M, \Delta, \norm{\cdot})$ be the sub-Finsler manifold. Let $g$ be a Riemannian tensor with a norm smaller than $\norm{\cdot}$ and denote by $d_{\rm Riem}$ the induced Riemannian distance. Recall that every other Riemannian distance is locally Lipschitz equivalent to $d_{\rm Riem}$.
Consider the identity map $\id: M\to M$.
Obviously, the map
$$\id :(M, \dcc )\to (M,d_{\rm Riem})$$
is $1$-Lipschitz (and, therefore, H\"older); and hence we obtain the inequality on the right in \eqref{eq_BB_thm}.

For the other bound, we observe that the step $s$ of $\Delta$ is equal to the maximum of the degree $d_j$ of the vector fields of some equiregular frame $X_1, \ldots, X_n$, which we locally fix.
For $r\in(0,1)$, it is evident that
$$B_{E} (0,r^s)\subset\prod_{j=1}^n[-r^s, r^s]\subset {\rm Box}(r),$$
where $B_{E} $ denotes the Euclidean ball in $\R^n$.
Therefore, using the second inclusion of the Ball-Box Theorem~\ref{Ball-Box} and the fact that the exponential maps $\Phi_p$ are locally biLipschitz maps (locally uniformly in $p$), as shown in Exercise~\ref{ex_phi_n_flows}, we obtain:
\begin{eqnarray*}
B_{\dcc} (p,Cr) \supseteq \Phi_p({\rm Box}(r)) \supseteq \Phi_p(B_{E} (0,r^s))
 \supseteq B_{d_{\rm Riem} }(p,C ' r^s).
\end{eqnarray*}
Hence, the map
$$\id :(M,d_{\rm Riem})\to (M, \dcc )$$
is $1/s$-H\"older on compact sets, establishing the inequality on the left of \eqref{eq_BB_thm}.
\qed

\subsection{Dimensions of CC spaces}
With the aid of the Ball-Box Theorem, we can finally calculate the Hausdorff dimensions of sub-Finsler manifolds.
\begin{definition}[Homogeneous dimension]\label{def Homogeneous dimension manifolds}
 	If a distribution $\Delta$ on a manifold $M$ is equiregular, we define its \emph{homogeneous dimension} as the natural number
	\begin{equation}\label{DefOfQ}
	Q := Q_{\Delta} := \sum_{j=1}^n j \left( \dim\Delta^{[j]}(p) - \dim\Delta^{[j-1]}(p) \right), 
	\end{equation}
	which is independent of $p$ as it varies in $M$.\index{homogeneous! -- dimension} 
\end{definition}
In other words, in terms of the numbers $m_1, \ldots, m_s$ from \eqref{eq:mj}, we express $Q$ as:
\begin{equation}\label{DefOfQ:1}
Q = m_1 + 2 (m_2-m_1) + 3 (m_3-m_2) + \dots + s (m_s-m_{s-1}).
	\end{equation}
It is important to note that the box defined in \eqref{def: box} satisfies
\[
\LL^n(\Bx(r)) = r^Q, 
\]
where $\LL^n$ is the Lebesgue measure in $\R^n$.
In terms of the degrees of the vector fields of an adapted frame, as in \eqref{eq_def_degree}, we also have 
\begin{equation}\label{DefOfQ:2}Q=\sum_{j=1}^n d_j.\end{equation}

\begin{corollary}\label{HDim}
If a sub-Finsler manifold $(M, \Delta, \|\cdot\|)$ has an equiregular distribution, then the Hausdorff dimension of $(M, \dcc )$ equals the homogeneous dimension $Q$. 
Moreover, the $Q$-dimensional Hausdorff measure of $(M, \dcc )$ is locally biLipschitz equivalent to each volume form. 
\\
	In particular, if $TM\neq\Delta$, the Hausdorff dimension is strictly greater than the topological dimension.
\end{corollary}
\begin{proof}
We arbitrarily choose a Riemannian structure. Since all volume forms are locally biLipschitz equivalent, we assume that the volume form corresponds to the Riemannian volume form \( \vol \).

 Using the notation of the Ball-Box Theorem~\ref{Ball-Box}, let $k$ be the (locally uniform) biLipschitz constant	of the exponential map $\Phi_p$ with respect to the Riemannian distance on the $n$-manifold $M$ and the Euclidean distance on $\R^n$; see Exercise~\ref{ex_phi_n_flows}. 
Since $\vol$ (resp., the Lebesgue measure $ \LL^n$) is the $n$-dimensional Hausdorff measure of the Riemannian manifold $M$ (resp., of the Euclidean space $\R^n$), we have, for small enough $r$,
	\[
	\frac1{k^n} \LL^n(\Bx(r)) \le \vol(\Phi_p(\Bx(r))) \le k^n \LL^n(\Bx(r)).
	\]
	If $Q$ is the homogeneous dimension, according to the Ball-Box theorem,
	for small $r$ then have
	\[
	\frac1{k^nC^Q}r^Q 
	\le \vol(B_{\dcc}(p,r))
	\le k^nC^Q r^Q
	.
	\]
	By Theorem~\ref{thm: Ahlfors: regular: Hausdorf: dim} and Remark~\ref{rmk: Ahlfors: regular: Hausdorf: dim}, obtained in the section on Ahlfors regular measures, we can conclude.
\end{proof}

\subsubsection{Dimensions of submanifolds in CC spaces}
The question of computing the Hausdorff dimension and Hausdorff measure of submanifolds in sub-Finsler manifolds with respect to the Carnot-Carath\'eodory distance is a natural one. Gromov provided a general formula for the Hausdorff dimension of smooth submanifolds in equiregular Carnot-Carath\'eodory spaces in \cite[0.6~B]{Gromov}. This formula was later demonstrated to coincide with the degree of the submanifold, introduced in \cite{Magnani-Vittone}, in \cite{Magnani_2010}.

\begin{theorem}[{\cite[page104]{Gromov}}]\label{Gromov-HDim}
Let $(M, \Delta, \norm{\cdot})$ be a sub-Finsler manifold with an equiregular distribution $\Delta$ and Carnot-Carath\'eodory distance $\dcc $. 
Consider a smooth submanifold $\Sigma\subset M$.
Then, the Hausdorff dimension of $(\Sigma, \dcc )$ is 
$$\dim_{H}(\Sigma, \dcc )=\max\left\{
 \sum_{j=1}^n j \cdot\dim \left( T_p \Sigma \cap \Delta^{[j]}(p) ) {\Big{/}} (T_p \Sigma \cap\Delta^{[j-1]}(p) \right) \; : \; p\in \Sigma\right\}.$$
\end{theorem}

Nevertheless, questions concerning Hausdorff {\em measures} of smooth submanifolds remain unanswered.
In \cite{Magnani-Vittone}, Magnani and Vittone derived an integral formula for the spherical Hausdorff measure of submanifolds in Carnot groups under a suitable `negligibility condition'.
This negligibility condition has been established in step-two groups, \cite{Magnani_2010} using covering arguments, and in the Engel group, using blow-up arguments \cite{LeDonne-Magnani}.
However, the situation remains unclear in higher step groups and, more generally, in sub-Riemannian manifolds.
For further details on this problem and its connections with existing literature, we direct the reader to Magnani's works \cite{Magnani-Vittone, Magnani08, Magnani_2010, MR3947860}.


\section{Exercises}
\begin{exercise}[Grushin distribution]\label{ex_Grushin}\index{Grushin plane}
On $\R^2$ coordinates $(x,y)$ the vector fields
\begin{equation}\label{vec_Grushin}
X= \partial_x \qquad \text{ and } \qquad Y= x\partial_y,
\end{equation}
satisfy the generating condition \eqref{Bracket: generating: frame} and define a bracket-generating distribution whose rank is not constant.
\end{exercise}

\begin{exercise}[Grushin subRiemannian plane]\index{Grushin plane}
On $\R^2$ with coordinates $(x,y)$, consider the vector fields $X$ and $Y$ as in \eqref{vec_Grushin}.
Define the Carnot-Carathéodory distance $\dcc$ on $\R^2$ by setting $X$ and $Y$ orthonormal, as in Definition \eqref{def_d_25may}.
For every $v\in\R$, $\lambda>0$ and $p,q\in\R^2$, we have
\[
\dcc\large(p+(0,v),q+(0,v)\large) = \dcc(p,q)
\qquad\text{and}\qquad
\dcc(\delta_\lambda p,\delta_\lambda q) = \dcc(p,q) ,
\]
where $\delta_\lambda := \begin{pmatrix}\lambda&0\\0&\lambda^2\end{pmatrix}$, and
\[
\dcc\large((0,y_1),(0,y_2)\large) = C \sqrt{|y_1-y_2|}, \qquad \forall y_1,y_2\in\R,
\]
where $C := \dcc((0,0),(0,1))$.
In particular, the line $\{(0,y):y\in\R\}$ is isometric to the $1/2$-snowflake of the Euclidean line. 
\end{exercise}

\begin{exercise}\label{Grushin_bad_metric}
On $\R^2$ with coordinates $(x,y)$, consider the distribution $\Delta$ generated by the vector fields $X$ and $Y$ as in \eqref{vec_Grushin}.
Consider the continuously varying norm $\|\cdot\|$ given by the Euclidean norm for every tangent vector.
Then, we have
\\
(i) the CC distance $\dcc $ induced by $(\Delta, \|\cdot\|)$ is the Euclidean distance. 
\\
(ii) The curve $t\in \R \longmapsto (0,t)$ is parametrized by arc length, but at no point it is tangent to the distribution $\Delta$.
\end{exercise}

\begin{exercise}\label{ex_Grushin_bad_frame}
Let $\phi\in C^\infty(\R)$ be such that $\frac{\dd^k\phi}{\dd x^k}(0)=0$, for all $k\in \{0,1,2, \ldots\}$, and $ \phi(x) \neq 0 $, for all $x\neq0$, as in Figure~\ref{fig:bump}.
On $\R^2$ with coordinates $(x,y)$, for the vector fields 
$$X= \partial_x \qquad \text{ and } \qquad Y= \phi(x)\partial_y,$$
we have
$$\Lie(\{X, Y \})= \Span_\R\left\{\partial_x, \frac{\dd^k\phi}{\dd x^k}(x)\partial_y: k \in\{0,1,2, \ldots\}\right\}.$$
Consequently, the pair $X$ and $Y$ does not satisfy the generating condition \eqref{Bracket: generating: frame} at $p=(0,0)$. 
Still, the vector fields span the same distribution of Exercise~\ref{ex_Grushin}
\end{exercise}

\begin{exercise}\label{ex_subbundle_closed}
Every subbundle of a vector bundle is a closed subset. 
\end{exercise}

\begin{exercise}
Every Finsler distance induces the manifold topology.
\end{exercise}

\begin{exercise}
Every two Finsler distances on a compact set are biLipschitz equivalent.
\end{exercise}

\begin{exercise}
Carnot-Carath\'eodory distances, and in particular Riemannian and Finsler distances, are length distances.
\end{exercise}

\begin{exercise}
The Hausdorff dimension of a Riemannian $n$-manifold is $n$.
\end{exercise}

\begin{exercise}
	If $\gamma: I\to(M, \dcc)$ is a curve in a CC space that is parametrized by arc length, then $\|\dot\gamma\|=1$ a.e.
\end{exercise}

\begin{exercise}\label{ex Lipschitz gives bounded speed}
	If $\gamma:I\to(M, \dcc)$ is a $L$-Lipschitz curve in a CC space, then $\|\dot\gamma\|\leq L$ a.e.
\end{exercise}

\begin{exercise}
 Let $(M, \Delta, \norm{\cdot})$ be a sub-Finsler manifold. We denote by $\Length_{\dcc }$ and $\Length_{\norm{\cdot}}$ respectively the length with respect to the metric $\dcc $ and the length with respect to the Finsler norm $\norm{\cdot}$. Let $\gamma$ be a horizontal curve. We have
$$\Length_ {\norm{\cdot}}(\gamma)=\Length_{\dcc }(\gamma).$$
\end{exercise}
\begin{exercise}
 For every absolutely continuous curve $\gamma$ in a sub-Finsler manifold,
we have
$$\gamma\text{ is horizontal }\Longleftrightarrow \Length_{\dcc }(\gamma)<+\infty.$$
\end{exercise}
\begin{exercise}
Denote by $\Phi_{X_i}^{t_i}$ the flow at time $i$ with respect to a vector field $X_i$.
	Calculate the differential of 
	\[
	(t_1, \dots,t_k)\longmapsto \Phi_{X_k}^{t_k}\circ\dots\circ\Phi_{X_1}^{t_1}(p) .
	\]
\end{exercise}

\begin{exercise}\label{proof: Chow}
Use Theorem~\ref{thm: reachable} and the fact that the points where \eqref{Bracket: generating} holds is open to give a (second) proof of Theorem~\ref{Chow}.
\end{exercise}

\begin{exercise} 
Recall that $\Gamma(\Delta)$ denotes the smooth sections of a subbundle $\Delta$ of a tangent bundle of a manifold $M$. Define $\Span(\Delta):={\rm Lie}\text{-}\Span\{\Gamma(\Delta) \}.$
The H\"ormander's condition is equivalent to $\Span(\Delta)=TM$.
{\rm (What is not immediately obvious is that elements of the form $[[X_1, X_2],[X_3, X_4]]_p$, with $X_1, X_2, X_3, X_4\in \Gamma(\Delta)$, are contained in some $ \Delta^{[j]}(p)$.)}
\end{exercise}

\begin{exercise} \label{Ex-Delta} 
Let $ \Delta^{[j]}(p)$ be the vector space defined in \eqref{Deltaj}.
The set $\Delta^{[j]}(p)$ can be equivalently defined as the subspace of $T_p M$ spanned by all commutators of the $X_i$'s of order $\leq j$ (including, of course, the $X_i$'s). 
Namely, $X_i(p)$ has order $1$; $[X_i, X_j](p)$ has order $2$; $[X_i,[X_j, X_k]](p)$ has order $3$; but those of order $4$ are those in one of the two forms:
$$ [X_i,[X_j,[X_k, X_l]]](p)\qquad \text{ or }\qquad[[X_i, X_j],[X_k, X_l]](p) .$$
\end{exercise}

\begin{exercise} For the set $ \Delta^{[j]} $ as  defined in \eqref{Deltaj} we have the following properties:
\\(i). The set $ \Delta^{[j]}$ might not be a subbundle of $T M$.
\\{\it Hint:} Try the distribution given by the frame $X_1=\partial_1$, $X_2=\partial_2 +x_1^2 \partial_3$.
 \\(ii). If $ \Delta^{[j]}$ is a subbundle, then it makes sense to consider smooth sections $\Gamma(\Delta^{[j]})$ 
 and
\\$ \Delta^{[j+1]}(p)=\Delta^{[j]}(p)+\R\text{-}\Span\left\{[X,Y](p) \;:\;X\in\Gamma(\Delta),Y\in \Gamma(\Delta^{[j]})\right\}.$
\end{exercise}

\begin{exercise} If $(M, \Delta, \|\cdot\|)$ is a sub-Finsler manifold with induced distance $\dcc$, then the metric space $(M, \dcc)$ is homeomorphic to the manifold $M$ via the identity map.
\end{exercise}

\begin{exercise}\label{ex_phi_n_flows} The maps $\Phi_p:\R^n \to M$ from \eqref{def_phi_n_flows} are locally biLipschitz maps locally uniformly in $p$: 
Namely, fix a compact subset $K$ of $M$ and a Riemannian distance $d_{\rm Riem}$ on $M$, then there exists $C>1$ and exists a neighborhood $U$ of $0$ in $\R^n$ such that of all $p\in K$ the map $\Phi_p|_{U} $ is a $C$-biLipschitz homeomorphism between $U$ equipped with the Euclidean distance and its image equipped with $ d_{\rm Riem}$. 
\end{exercise}

\begin{exercise} 
Each smooth surface in the Heisenberg group has Hausdorff dimension $3$.
\\{\it Hint.} One may use Theorem~\ref{Gromov-HDim} or give a direct proof considering vertical planes as preliminary cases.
\end{exercise}

\begin{exercise} 
Give a proof of Theorem~\ref{Gromov-HDim}.
\end{exercise}

 \begin{exercise}\label{ex:Lip_for_CCmfds}
 Let $F: M_1 \to M_2$ be a smooth map between sub-Finsler manifolds $(M_1, \Delta^{M_1}, \norm{\cdot})$ and $({M_2}, \Delta^{M_2}, \norm{\cdot})$. Let $L>0$.
Assume, for all $ p\in M_1$, that $\dd \pi (\Delta^{M_1}_p) \subseteq \Delta^{M_2}_{F(p)}$ and that 
 $$ (\dd \pi)_p :(\Delta_p, \norm{\cdot} ) \to (\Delta_{F(p)}, \norm{\cdot} ) $$
 is $L$-Lipschitz.
 Then, the map $F: M_1 \to M_2$ is $L$-Lipschitz with respect to the respective sub-Finsler metrics.
 \end{exercise}

\begin{exercise}
For every CC-bundle structure $(\sigma, N)$ on a manifold $M$, as in Definition~\ref{rank-varying_structure},
given an absolutely continuous curve $\gamma: [0,1] \to M$ with $\norm{\dot\gamma(t)}<\infty$ for almost every $t\in [0,1]$, then there exists a measurable map $u: [0,1] \to E$ such that $\sigma(u(t)) = \dot\gamma(t)$ and $\norm{\dot\gamma(t)} = N(u(t)) $, for almost every $t\in [0,1]$.
\end{exercise}

\begin{exercise}\label{ex_20may1006}
As in Definition~\ref{rank-varying_structure}, let $M$ be a smooth manifold and 
$ f:M\times \mathbb R^m\to TM $
 a smooth $M$-bundle morphism. 
Let 
$ N:M\times \mathbb R^m\to [0,+\infty) $
be a continuous function such that $N(p, \cdot)$ is a norm for every $p\in M$.
Then, the associated sub-Finsler distance \eqref{def_d_25may}
between $p$ and $q$ in $M$ is
\begin{equation}\label{def_dist_bundles}
d_{(f,N)}(p,q):=\inf\left\{\int_0^1 N(\gamma(s),u(s))\dd s \, \left|\,
u\in L^\infty([0,1];\mathbb R^m),
\begin{array}{rl}
\gamma(0)&=p \\
\gamma(1)&=q \\
\dot\gamma(t) &= f(\gamma(t),u(t))
\end{array}\right.\right\}.
\end{equation}
\end{exercise}



\chapter{A review of Lie groups}\label{ch_LieGroups}
%
%

In the following chapter, we will review the theory of Lie groups. This revision serves two purposes:
First, Lie groups equipped with special sub-Finsler structures appear as tangent spaces of Carnot-Carathéodory spaces. 
Such Lie groups are the infinitesimal models for sub-Riemannian manifolds, and therefore, they play the same role as the Euclidean vector spaces in Riemannian geometry. 
Second, sub-Finsler structures on Lie groups are highly interesting and arise in various contexts, including geometric group theory, harmonic analysis, hyperbolic geometry, and furthermore in stochastic processes and mechanics. They are, in a sense, easier to study than general Carnot-Carathéodory spaces.

The prerequisites for understanding Lie groups and Lie algebras primarily lie in the realm of differential geometry. The results presented in this chapter are classical and are based on the references \cite{Warner, Corwin-Greenleaf, Hilgert_Neeb:book}.

\section{Lie groups, Lie algebras, and their morphisms}\label{Sec: Lie_objects}
\index{group}
\index{group! -- product}
\index{product! group --}
\index{associativity law}
\index{identity! -- element}
\index{inversion}
In this section, we will review the following concepts: Lie group, Lie algebra, Lie algebra associated with a Lie group, Lie subgroup, Lie subalgebra, Lie group homomorphism, Lie algebra homomorphism, and Lie algebra homomorphism induced by a Lie group homomorphism. 
We will also state certain results regarding these objects, but the proofs will be deferred to later sections.

For clarity, we provide a reminder that a {\em group} is a set $G$ equipped with a binary operation, referred to as its {\em product} or {\em group product} or {\em product law}, usually denoted by the symbol $\cdot$. The product is a function $(a,b)\in G\times G\longmapsto a\cdot b\in G$ that satisfies associativity, the existence of an identity element, and of an inversion map. The inversion map is denoted as $a\mapsto a^{-1}$. 
The identity element of a group $G$ is denoted by $1$. If there is a need to emphasize that $1$ is specifically the identity element of the group $G$, it may be denoted as $1_G$. Other texts or references may use alternative symbols such as $e$ or $e_G$.

 	Let $G$ be a group and $g\in G$. The \emph{left translation} by $g$ is the bijection
	\index{left! -- translation}
	\index{translation! left --}
	\[
	\begin{aligned}
	 	L_g :& 	&G& &\longrightarrow &&G\, \\
	 		&	&h& &\longmapsto &&gh.
	\end{aligned}
	\]
	The \emph{right translation} by $g$ is the bijection 
	\index{right translation}
	\index{translation! right --}
	\[
	\begin{aligned}
	 	R_g :& 	&G& &\longrightarrow &&G\, \\
	 		&	&h& &\longmapsto &&\, h g.
	\end{aligned}
	\]
		The \emph{conjugation} by $g$ is the bijection 
	\index{conjugation} 
	\[
	\begin{aligned}
	 	C_g :& 	&G& &\longrightarrow &&G\qquad \\
	 		&	&h& &\longmapsto &&\, g h g^{-1}.
	\end{aligned}
	\]


We shall focus on Lie groups, which are differentiable manifolds with a smooth group operation. However, some of the remarks will hold in the general setting of topological groups:\index{topological group}
A {\em topological group} is a group together with a Hausdorff topology for which the group product and the inversion map are continuous.
Lie groups are special topological groups:

\begin{definition}[Lie group]\index{Lie group}\index{Lie! -- group|see {Lie group}}\index{group! Lie --|see {Lie group}}  
 A Lie group is a differentiable manifold (second countable, but not necessarily connected) together with a group structure such that both 
 \begin{equation}\begin{array}{cccccccc}
\text{the product }& G\times G&\to& G			&\qquad\text{ and the inverse } \qquad&G&\to& G\\
 & (x,y)&\mapsto &x\cdot y 		&&	g&\mapsto& g^{-1}
 \end{array} \end{equation}
are $C^\infty$ maps.
\end{definition}


As in every manifold, the set $\Vec(G)$ of vector fields on $G$ forms a Lie algebra.
The general notion of Lie algebra is the following:
\begin{definition}[Lie algebra]\label{def: Lie_algebra}\index{Lie algebra}\index{Lie! -- algebra|see {Lie algebra}} \index{algebra! Lie --|see {Lie algebra}} 
 	A \emph{Lie algebra $\g$} ({\em over} $\R$) is a vector space (over $\R$) together with a bilinear operation $[\cdot, \cdot]: \g\times\g\to\g$ called \emph{Lie bracket},\index{bracket|see {Lie bracket}}\index{Lie! -- bracket} 
	such that for all $X, Y, Z\in\g$, the following two properties hold:
	
		$[X, Y] = -[Y, X]\qquad$ ({\em anti-commutativity}), \index{anti-commutativity}

	$[[X, Y],Z] + [[Y,Z], X] + [[Z, X],Y] = 0\qquad$ ({\em Jacobi identity}).\index{Jacobi identity}\index{identity! Jacobi --}
\end{definition}
Lie algebras are usually denoted by gothic letters.
The gothic letters for $g,h,n,o,l,p,s$ are $\mathfrak{g,h,n,o,l,p,s}$. 
Lie algebras can also be considered on other fields than $\R$. However, in this text, we shall only consider those over the real numbers.
The structure of a Lie algebra can be represented by expressing the Lie bracket using a basis. 
Namely, if $\g$ is a Lie algebra with bracket $[\cdot, \cdot]$ and $X_1, \dots, X_n$ is an ordered basis of $\g$ as vector space, then the \emph{structural constants}\index{structural constants} of $\g$ with respect to $X_1, \dots, X_n$ are the real numbers $c^k_{ij}\in\R$ with $i,j,k\in\{1, \dots,n\}$ such that
	\begin{equation}\label{eq1655}
	[X_i, X_j] = \sum_{k=1}^n c^k_{ij}X_k,
	\qquad\forall i,j\in\{1, \dots,n\}.
	\end{equation}
 The data $X_1, \dots, X_n$ and $(c^k_{ij})_{i,j,k\in \{1, \ldots, n\}}$ record the whole info about the Lie bracket; see Exercises~\ref{str_const_1}, \ref{str_const_2}, and~\ref{str_const_3}.

The importance of the concept of Lie algebra is that there is a special finite-dimensional Lie algebra intimately associated with each Lie group, and that properties of the Lie group are reflected in properties of its Lie algebra. 
We shall recall, for example, that simply connected Lie groups are completely determined (up to isomorphism) by their Lie algebras; see Corollary~\ref{simplyconn-isom}.

The Lie algebra associated with a group is isomorphic, as a vector space, to the tangent space $T_{1_G}G$ at the identity element $1_G$. In order to define a Lie bracket structure, one identifies $T_{1_G}G$ as a subset of the space $\Vec(G) 
$ of smooth vector fields on $G$ by suitably extending each vector to a vector field.
Forced to make a choice\footnote{Actually, we prefer to consider left-actions by a group (on itself) because we think of groups as transformations, and we are nowadays used to put symbols of functions on the left of variables, like $f(x)$.}, we follow the majority of the literature
focusing on the \emph{left} invariant vector fields,\index{left-invariant! -- vector field}
i.e., the vector fields $X\in \Vec(G)$ such that 
$(\dd L_g)_h X_h=X_{L_g(h)}$ for all $g,h\in G$. 
Thanks to \eqref{bir} with $F=L_g$, the class of left-invariant vector fields is easily seen to be closed under the Lie bracket; see Exercise~\ref{LIVFs is closed under Lie bracket}.
In other words, the set of left-invariant vector fields form a Lie algebra. 

Note that, after fixing a vector $v\in T_{1_G}G$, we can construct a left-invariant vector field $X$ defining $X_g:=(\dd L_g)_{1_G} (v)$ for $g\in G$. This construction is a linear isomorphism between the set of all left-invariant vector fields and $T_{1_G}G$, and proves that left-invariant vector fields form an $n$-dimensional subspace of $\Vec(G)$, where $n:=\dim G$.
We denote by $\g$ the vector space $T_{1_G}G$ equipped with the Lie bracket coming from the identification with the left-invariant vector fields.
Such a $\g$ is called the \emph{Lie algebra} of $G$, and it is occasionally denoted by ${\rm Lie}( G)$. We next summarize this definition:
\begin{definition}[Lie algebra of a Lie group]\label{def2610}\index{Lie algebra! -- of a Lie group, $\g$, $\Lie(G)$}
Let $G$ be a Lie group. The \emph{Lie algebra} of $G$, denoted by $\Lie(G)$, has two realizations:

Interpretation 1:	
	$\Lie(G)$ is the linear space ${\rm LIVF}(G)$ of left-invariant vector fields on $G$ endowed with the bracket of vector fields.\index{LIVF}

 Interpretation 2: 
 $\Lie(G)$ is the tangent space $T_{1_G}G$ equipped with the bracket
	\[
	[X, Y] := [\tilde X, \tilde Y]_{1_G},
	\quad\forall X, Y\in T_{1_G}G,
	\]
	where $\tilde X, \tilde Y$ are the left-invariant vector fields such that $\tilde X_{1_G}=X$ and $\tilde Y_{1_G}=Y$, respectively.
	We shall use both points of view.
\end{definition}

	\index{Lie! -- subgroup} 
	\index{subgroup! Lie --}
 	Let $G$ be a Lie group and $H < G$ a subgroup.
	We say that $H$ is a \emph{Lie subgroup} of $G$ if $H$ admits the structure of a Lie group such that the inclusion $H\into G$ is a smooth group homomorphism.
	It is a consequence that the inclusion is actually an immersion; see Exercise~\ref{Ex:homo:immersion}.
 	A Lie subgroup $H< G$ is said to be a {\em closed Lie subgroup}
\index{closed Lie subgroup}
if $H$ is topologically closed within $G$.
	It is a consequence that, in this case, the inclusion $H\into G$ is an embedding; see Theorem~\ref{teo1131}.
Closed Lie subgroups are also called {\em regular Lie subgroups}, 
	especially if one wants to stress that it is an embedded submanifold (not just immersed).
\index{regular Lie subgroup}
\index{Lie! regular -- subgroup}
\index{Lie! closed -- subgroup}
\index{closed Lie subgroup}

	A \emph{subalgebra} of a Lie algebra $\g$ is a vector subspace $\h\subset\g$ that is closed under the Lie bracket operation of $\g$.\index{Lie! -- subalgebra}\index{subalgebra|see{Lie subalgebra}}
Hence, if $H$ is a Lie subgroup of a Lie group $G$, then it is an exercise to show that $\Lie(H)$ is canonically isomorphic to a subalgebra of $\Lie(G)$. Vice versa, every subalgebra comes from a Lie subgroup:
	 \begin{theorem}[Existence of subgroups; see Theorem~\ref{teo1145}]\label{teo1145bis}
 	Let $G$ be a Lie group.
	For every subalgebra $\h\subset \Lie(G)$, there is a unique connected Lie subgroup $H$ with Lie algebra $\h$. 
\end{theorem}
A subspace $\mathfrak i \subseteq \g$ is an {\em ideal}, \index{ideal} if $[\mathfrak i, \g] \subseteq\mathfrak i$.
A subgroup $N$ of a group $G$ is {\em normal}, if $C_g(N)=N$, for all $g\in G$.
For ideals $\mathfrak i \subseteq \g$ and normal subgroups $N\subseteq G$ we write
 $\mathfrak i \vartriangleleft \g $ and $N \vartriangleleft G $, respectively. The two notions are connected: it is an exercise to show that the Lie algebra of a connected Lie subgroup is an ideal if and only if the subgroup is normal. 

Next, we discuss the maps of the categories in which the objects are the Lie groups and the Lie algebras, respectively.	A map $\varphi: G\to H$ between groups is a \emph{group homomorphism}, or, simply, a {\em homomorphism} or a {\em morphism}, if
	\index{group! -- homomorphism}\index{morphism}
	\index{homomorphism! group --}
	\[
	\varphi(g_1\cdot g_2) = \varphi(g_1)\cdot \varphi(g_2),
	\qquad\forall g_1, g_2\in G .
	\]	
If 
 	 $G,H$ are Lie groups, then 
	a homomorphism $\varphi: G\to H$ is called \emph{Lie group homomorphism} if it is smooth. 	
	If in addition $H=G$, then $\varphi$ is called \emph{Lie group endomorphism}.
	A bijective Lie group homomorphism is called \emph{Lie group isomorphism}.
	A bijective Lie group endomorphism is a \emph{Lie group automorphism}.
	\index{Lie group! -- homomorphism}
	\index{homomorphism! Lie group --}
	\index{Lie group! -- endomorphism}
	\index{endomorphism! Lie group --}
	\index{Lie group! -- isomorphism}
	\index{isomorphism! Lie group --}
	\index{Lie group! -- automorphism}
	\index{automorphism! Lie group --}

A map $\psi : \g\to\h$ 
 between Lie algebras is called {\em Lie algebra homomorphism}\index{Lie algebra! -- homomorphism}\index{homomorphism! Lie algebra --} if it is linear and preserves brackets: 
 $$ \psi([X, Y])=[ \psi(X), \psi(Y) ], \qquad \forall X, Y\in \g.$$
 	If in addition $\h=\g$, then $\psi$ is called \emph{Lie algebra endomorphism}.
	\index{Lie algebra! -- endomorphism}
	\index{endomorphism! Lie algebra --}
	A bijective Lie algebra homomorphism (resp.~endomorphism) is called \emph{Lie algebra isomorphism} (resp.~\emph{automorphism}).
	\index{Lie algebra! -- isomorphism}
	\index{isomorphism! Lie algebra --}
	\index{Lie algebra! -- automorphism}
	\index{automorphism! Lie algebra --}

 Each Lie group homomorphism induces a Lie algebra homomorphism: if $\varphi: G\to H$ is a Lie group homomorphism, note that $\varphi(1_G)=1_H$, and one can easily show that the differential at the identity element
\begin{equation}\label{eq:induced_algebra_homomorphism}
\varphi_*:=\dd \varphi_{1_G}: T_{1_g} G\to T_{1_H} H
\end{equation}
 commutes with the Lie bracket operation; see Exercise~\ref{phi_star_homo}.
 Namely, the map $\varphi_*: {\rm Lie}( G)\to {\rm Lie}( H)$ is a Lie algebra homomorphism, called the {\em Lie algebra homomorphism induced} by $ \varphi$. 	\index{induced! -- Lie algebra homomorphism}

Vice versa, in the case when $G$ is a Lie group that is simply connected as topological space, then each 
Lie algebra homomorphism comes from a Lie group homomorphism:
\begin{theorem}[Induced Lie group homomorphism; see Theorem~\ref{thm_induced_homo}]\label{thm_induced_homo0} 
\index{induced! -- Lie group homomorphism}
\index{simply connected}
 	Let $G$ and $H$ be Lie groups.
	Assume $G$ is simply connected.
	For every Lie algebra homomorphism $\psi: \Lie(G)\to\Lie(H)$, there exists a unique Lie group homomorphism $\varphi: G\to H$ with $\varphi_*=\psi$.
\end{theorem}

\begin{corollary}\label{simplyconn-isom}
If simply connected Lie groups $G$ and $H$ have isomorphic Lie algebras, then $G$ and $H$ are Lie group isomorphic.
\end{corollary}

As a consequence of a theorem due to Ado, see \cite[page 199]{Jacobson} and also Section~\ref{sec_Birkhoff}, for every Lie algebra $\g$ there exists a simply connected Lie group $G$ with Lie algebra $\g$. We then have the following correspondence.

 \begin{theorem}\label{Lie's Third Theorem}
 There is a one-to-one correspondence between isomorphism classes of Lie algebras and isomorphism classes of simply connected Lie groups.
\end{theorem}
We shall only prove the above theorem, together with Ado's result, for stratified Lie algebras since the proof is much easier, and it is what is needed for the Lie groups of our interest: the Carnot groups.
We refer to Section~\ref{sec Birkhoff for stratified}.

\section{Exponential map}
Let $M$ be a differentiable manifold. 
Consider a smooth vector field $X \in \Vec(M)$. Given a point $p\in M$, there exists a unique curve $t\mapsto \gamma(t)$ satisfying $\gamma(0) = p$ and having a tangent vector $\dot{\gamma}(t) = X_{\gamma(t)}$. 
We refer to this curve as the integral curve of $X$ passing through $p$. The {\em exponential} of $X$ is defined as $\Phi_X^1(p) = \gamma(1)$, which gives us the endpoint of the integral curve after a unit of time.
It should be noted that generally, the exponential of $X$ is defined 
only for $X$ is some small neighborhood of zero in $\Vec(M)$ and maps it to a neighborhood of $p$ in the manifold. This locality arises from the reliance on the theorem of existence and uniqueness of ordinary differential equations, which is itself local in nature.

In the theory of Lie groups, the {\em exponential map} is a map from the Lie algebra $\g$ to the group $G$: 
 $$\exp\colon \g \to G , \qquad X\mapsto \Phi_X^1(1_G).$$
 Here, elements of the Lie algebra $\mathfrak{g}$ are identified with left-invariant vector fields, and thus we have $\mathfrak{g} \subset \Vec(G)$. 
It can be shown that for every $X \in \mathfrak{g}$, the ordinary differential equation $\dot{\gamma}(t) = X_{\gamma(t)}$ has global solutions and that these integral curves $\gamma(t)$ correspond to Lie group homomorphisms from the additive group $\mathbb{R}$ to the group $G$. Such homomorphisms from $\mathbb{R}$ to $G$ are commonly referred to as {\em one-parameter subgroups}.

\subsection{One-parameter subgroups}
\begin{definition}[One-parameter subgroup]
 	Let $G$ be a Lie group. A Lie group homomorphism $\theta: \R\to G$ is called a \emph{one-parameter subgroup} (OPS, for short).
With abuse of terminology, sometimes we say that a one-parameter subgroup is the image $\theta(\R)\subseteq G$ of one such map. \index{one-parameter subgroup}
\index{subgroup! one-parameter --, OPS}
\index{OPS|see{one-parameter subgroup}}
\end{definition}

 	Equivalently, $\theta: \R\to G$ is a one-parameter subgroup if and only if
	\begin{enumerate}[label=\roman*).]
	\item 	$\theta$ is smooth,
	\item 	$\theta(0)=1_G$,
	\item 	$\theta(t+s)=\theta(t)\cdot\theta(s)$, for all $s,t\in\R$.
	\end{enumerate}
	
We will soon see that the one-parameter subgroups are exactly the integral curves from the identity element of the left-invariant vector fields (and also of the right-invariant vector fields).
\index{flow! -- of a vector field}
\index{flow! -- line}
Recall that we denote by $\Phi_X^t(p)$ the {\em flow} of a vector field $X$ at time $t$ starting from a point $p$.

\begin{proposition}\label{flow_vs_OPS}
 	Let $G$ be a Lie group and $X$ be a left-invariant vector field on $G$.
\\	\ref{flow_vs_OPS}.i.
	The flow line $t\longmapsto \Phi^t_X(1_G)$ of $X$ from $1_G$ is a one-parameter subgroup.
\\	\ref{flow_vs_OPS}.ii.	If $\theta: \R\to G$ is a one-parameter subgroup with $\dot\theta(0)=X_{1_G}$, then $\theta(t) = \Phi^t_X(1_G)$, for all $t\in\R$.
\end{proposition}
\begin{proof}[Proof of \ref{flow_vs_OPS}.i.]
 	Let $\sigma(t)=\Phi^t_X(1_G)$, which is defined for $t$ in some maximal interval $(-\epsilon, \epsilon)$. Fix $s\in(-\epsilon, \epsilon)$ and consider $\gamma(t):=\sigma(s)\cdot\sigma(t)$.
	We claim that $\gamma$ is the integral curve of $X$ from $\sigma(s)$.
	Indeed, we have
	\begin{eqnarray*}
	 	\dot\gamma(t) &=& \frac{\dd}{\dd t} (\sigma(s)\cdot\sigma(t)) \\
		&=& \frac{\dd}{\dd t} (L_{\sigma(s)}(\sigma(t)) ) \\
		&=& (\dd L_{\sigma(s)} )_{\sigma(t)} \sigma'(t) \\
		&=& (\dd L_{\sigma(s)} )_{\sigma(t)} X_{\sigma(t)}\\
		&=& X_{\sigma(s)\cdot\sigma(t)} 
		= X_{\gamma(t)} .
	\end{eqnarray*}
	By uniqueness of integral curves, we have $\gamma(t) = \sigma(s+t)$ and so $\sigma(s+t) = \sigma(s)\cdot\sigma(t)$.
	Moreover, since $\sigma$ can be prolonged by $L_{\sigma(s)}\circ\gamma$, then $\sigma$ is defined on the whole of $\R$.
\proof[Proof of \ref{flow_vs_OPS}.ii.] 	Being $\theta$ a one-parameter subgroup, we have $\theta(s+t)=\theta(s)\cdot\theta(t) = L_{\theta(s)}(\theta(t))$.
	Hence, since $\dot\theta(0)=X_1$, we have
	\begin{eqnarray*}
	 	\dot\theta(s)
		&=& \left. \frac{\dd}{\dd t} \theta(s+t) \right|_{t=0}\\
		&=& \left. \frac{\dd}{\dd t} L_{\theta(s)}(\theta(t)) \right|_{t=0} \\
		&=& \left( \dd L_{\theta(s)} \right)_{\theta(0)} \dot\theta(0)\\
		&=& \left( \dd L_{\theta(s)} \right)_1 X_1
		= X_{\theta(s)} .
	\end{eqnarray*}
	So $\theta$ is the integral curve of $X$ from $1_G$. \qedhere
\end{proof}
\begin{remark}If $\theta$ is a OPS, then its image $\theta(\R)$ is a Lie subgroup.
 	Indeed, if $\dot\theta(0)=0$, then $\theta$ is constantly equal to $1_G$, which is a Lie subgroup of dimension 0.
	If instead $\dot\theta(0)\neq0$, then $\dot\theta(t)\neq0$ for all $t\in\R$, and $\theta$ is an immersion. Hence, $\theta(\R)\subset G$ is a Lie subgroup of dimension 1.
\end{remark}

\subsection{Exponential map}
\begin{definition}[Exponential map]
 	Let $G$ be a Lie group and $\g$ its Lie algebra, seen as the space of left-invariant vector fields.	
	The \emph{exponential map} is defined as
	\[
	\exp: \g\to G
	, \quad
	X\in\g\longmapsto\exp(X) := \Phi^1_X(1_G),
	\]
	i.e., $\exp(X)$ is the flow of $X$ at time $1$ starting from $1_G$.
	\index{exponential! -- map of a Lie group}
\end{definition}
\begin{remark}

The exponential map may be different from the exponential map of Riemannian geometry. 
	In Exercise~\ref{ex:exp:no:sujective}, one can see that the exponential map of the Lie group $\GL^+(n, \R)$ is not a Riemannian exponential for any Riemannian metric.
However, if a Lie group is compact, then it has a Riemannian metric invariant under left and right translations, and the Lie group exponential map is the Riemannian exponential map of this Riemannian metric; see Section~\ref{sec:biinv_compact}.
\end{remark}

One first key property of the exponential map is the following -- see Exercise~\ref{exp: OPS} for other properties.
\begin{corollary}[of Proposition~\ref{flow_vs_OPS}]\label{Warner3.31}For every left-invariant vector field $X$
 the curve $t\mapsto\exp(tX)$ is a one-parameter subgroup and an integral curve of $X$. For every $g\in G$, 
 the curve $t\mapsto g\exp(tX)=L_g(\exp(tX))$ is the flow line of $X$ starting at $g$. 
\end{corollary}
Consequently, first, we infer that left-invariant vector fields are complete. 
Second, we proved that the flows of \emph{left} invariant vector fields are \emph{right} translations, as we next express. 
\begin{corollary}\label{prop:right_left}
 	Let $X$ be a left-invariant vector field on a Lie group $G$. Then
 \begin{equation}\label{flow_of_left}
		\Phi^t_X = R_{\exp(tX)}, \qquad \forall t\in\R. 
		\end{equation}
\end{corollary}

Likewise,
\index{right-invariant! -- vector field}
\index{vector field! right-invariant --}
if we let $X^\dag$ be the right-invariant vector field such that $(X^\dag)_1=X_1$,
	then we also have 
	\begin{equation}\label{flow_of_right}
	\exp(tX) = \Phi^t_{X^\dag}(1_G) \qquad \text{ and }\qquad
	\Phi_{X^\dag}^t = L_{\exp(tX)}, \quad \forall t\in\R.
	\end{equation}
	From the fact that we explicitly know the above flows \eqref{flow_of_left} and \eqref{flow_of_right}, we have many consequences; see Exercises~\ref{ex: LIVF: commutes: RIVF}, \ref{ex: commutative}, and~\ref{LieRIVF}.

	\index{left-invariant! -- vector field}
	\index{vector field! left-invariant --}
\index{flow! -- of a vector field}
\index{flow! -- line}
	\index{right translation}
	\index{translation! right --}

 We summarize the following three interpretations of the exponential map: 
 	\[
	\exp(X) = 
	\begin{cases}
	 	\text{flow at time 1 of the LIVF $X$}, \\
		\text{OPS at time 1 tangent to $X_{1_G}$ (or $X^\dag_{1_G}$) at time }0, \\
		\text{flow at time 1 of the RIVF $X^\dag$}.
	\end{cases}
	\]

The next result is a very important feature of the exponential map. It implies that $\exp$ gives a local parametrization of $G$ near $1_G$.

\begin{proposition}\label{exp:diffeo}
\index{local! -- diffeomorphism}
\index{differential! -- of $\exp$ at $0$}
 	Let $G$ be a Lie group with Lie algebra $\g $.
	Then $\exp: \g\to G$ is smooth and 
$(\dd\exp)_0$ is the identity map:
$$(\dd\exp)_0=\id_\g: \g=T_0\g \to T_{1_G}G=\g.$$
Consequently, the map $\exp$ is a diffeomorphism between some neighborhood of $0$ in $\g$ and some neighborhood of $1_G$ in $G$;

\end{proposition}
\begin{proof}
For the smoothness of $\exp$, we refer to Exercise~\ref{exp_smooth}.
	Regarding its differential, 	fix $X\in\g=T_{1_G}G$. Let $\sigma: \R\to\g$ be the curve $\sigma(t):=tX$ so that $\sigma'(0)=X$.
	Then 
	\begin{eqnarray*}
	(\dd\exp)_0(X) 
	&=& (\exp\circ\sigma)'(0)\\
	&=& \left.\frac{\dd}{\dd t} \exp(tX)\right|_{t=0} \\
	&=& \left.\frac{\dd}{\dd t} \Phi^1_{t\tilde X}(1)\right|_{t=0}\\
	&=& \left.\frac{\dd}{\dd t} \Phi^t_{\tilde X}(1)\right|_{t=0}\\
	&=& X,
	\end{eqnarray*}
	where $\tilde X$ is the left-invariant vector fields with $\tilde X_{1_G}=X$.
The last part of the statement of the proposition is a consequence of the Inverse Function Theorem.\qedhere
\end{proof}

The exponential map gives a first link between the Lie-group level and the Lie-algebra level: 
 \begin{proposition}\label{prop: expF: Fexp}
\label{Warner3.32}
\index{induced! -- Lie algebra homomorphism}
 	Let $\varphi: G\to H$ be a Lie group homomorphism.
	If $\varphi_*: \Lie(G)\to\Lie(H)$ is the induced Lie algebra homomorphism, as in \eqref{eq:induced_algebra_homomorphism}, then
	\[
	\exp\circ \varphi_*=\varphi\circ\exp,
	\]
	i.e., the following diagram commutes.
	\[
	\xymatrix{
	 \Lie(G)\ar[d]_{\exp}\ar[r]^{\varphi_*} & \Lie(H)\ar[d]^{\exp} \\
	 G\ar[r]_\varphi & H
	}
	\]
\end{proposition}
\begin{proof}
 	We need to show that for every left-invariant vector field $X$
	\[
	\varphi(\exp(X)) = \exp(\widetilde{(\dd \varphi)_1X_1}) .
	\]
	We plan to show that for every left-invariant vector field $X$ and for every $t\in\R$
	\[
	\sigma(t):= \varphi(\exp(tX)) = \exp(t\widetilde{(\dd \varphi)_1X_1} ) .
	\]
	Namely, we claim that the curve $t\mapsto\sigma(t)$ is the one-parameter subgroup in $H$ generated by $(\dd \varphi)_1X_1$.
	
	First, we check that $\sigma$ is a one-parameter subgroup:
	\begin{eqnarray*}
	 	\sigma(s)\sigma(t) 
		&=& \varphi(\exp(sX))\varphi(\exp(tX))\\
		&=& \varphi\left(\exp(sX)\exp(tX)\right) \\
		&=& \varphi\left( \exp((s+t)X) \right)\\
		&=& \sigma(s+t), 
	\end{eqnarray*}
	where we used that $\varphi$ is a homomorphism and that $t\mapsto\exp(tX)$ is a one-parameter subgroup.
	
	Second, the derivative at $0$ of $\sigma$ is
	\begin{eqnarray*}
	 	\left. \frac{\dd}{\dd t} \sigma(t) \right|_{t=0}
		&=& \left. \frac{\dd}{\dd t} \varphi(\exp(tX)) \right|_{t=0} \\
		&=& ( \dd \varphi)_{\exp(0\cdot X)} \left.\frac{\dd}{\dd t} \exp(tX) \right|_{t=0} \\
		&=& (\dd \varphi)_1 X_1 . 
	\end{eqnarray*}\qedhere
\end{proof}


\subsection{Exponential coordinates}\label{sec1120}
We draw a first consequence of Proposition~\ref{exp:diffeo}.
Let $X_1, \dots, X_n$ be a basis of the Lie algebra of a Lie group $G$.
The map $\alpha: \R^n\to G$, 
\[
\alpha(t_1, \dots,t_n) := \exp(t_1X_1+\dots+t_nX_n)
\]
is a diffeomorphism between some neighborhood of $0\in\R^n$ and some neighborhood of $1_G$ in $G$.
Such a map is called \emph{exponential local coordinate map} (or \emph{exponential local coordinates of the first kind}) with respect to $X_1, \dots, X_n$.
\index{exponential! -- coordinates! -- of the first kind}

The map $\beta: \R^n\to G$
\[
\beta(t_1, \dots,t_n) := \exp(t_1X_1)\cdots\exp(t_nX_n)
\]
is called \emph{exponential local coordinates of the second kind} with respect to $X_1, \dots, X_n$.
\index{exponential! -- coordinates! -- of the second kind}

One can consider intermediate examples. For example, given $k\in\{1, \ldots, n-1\}$, one can let $\beta_k: \R^n\to G$ be the map
\[
\beta_k(t_1, \dots,t_n) := \exp(t_1X_1+\cdots +t_k X_k)\exp(t_{k+1} X_{k+1}+\cdots + t_nX_n),
\]
which is called an \emph{exponential local coordinates of mixed kind} with respect to $X_1, \dots, X_n$.
\index{exponential! -- coordinates! -- of mixed kind}

Notice that $\beta$ and $\beta_k$ might depend on the ordering of the basis.
The maps $\beta$ and $\beta_k$ are indeed coordinate maps since, for them, the differential at $0$ is an isomorphism. Indeed, for $\beta$ we have
\begin{eqnarray*}
(\dd \beta)_0 (\partial_j|_0)&\stackrel{\text{def}}{=}& \left.\frac{\partial}{\partial t_j}\beta(t_1, \dots,t_n) \right|_{(t_1, \dots,t_n)=(0, \dots,0)}\\
&=& \left.\frac{\dd}{\dd t_j} \beta(0, \dots,0,t_j,0, \dots,0) \right|_{t_j=0} \\
&=& \left.\frac{\dd}{\dd t_j} \exp(t_jX_j) \right|_{t_j=0}\\
&=& X_j .
\end{eqnarray*}

Warning: 
There are examples of groups for which $\alpha$ and $\beta$ are not surjective.




\section{From continuity to smoothness}

\subsection{Smoothness of continuous homomorphisms}\label{sec smoothness of continuous homomorphisms}
\index{continuous homomorphisms}
\index{smoothness! -- of continuous homomorphisms}

This subsection aims to show that continuous homomorphisms between Lie groups are smooth. 
We begin with the 1-dimensional situation: the domain of the homomorphism is the additive group $\R$. 
\begin{theorem}\label{thm: 1Dcontinuous: smooth}
 	Let $G$ be a Lie group. Every 
	 continuous homomorphism $\theta: \R\to G$ is smooth.
\end{theorem}
\begin{proof}
 	Since $\theta(t+s)=\theta(t)\theta(s)=L_{\theta(t)}(\theta(s))$, it is enough to prove that $\theta$ is smooth at $0$.
	
	Let $U\subset \g$ be an open set that is star-shaped with respect to $0$ and such that $\exp|_U$ is a diffeomorphism between $U$ and its image.
	Let 
	\[
	U':=\tfrac12 U=\left\{\tfrac12X:X\in U\right\} .
	\]
	By continuity there is $t_1>0$ such that $\theta(t)\in\exp(U')$, for all $|t|\le t_1$.
	Let $X_1\in U'$ be such that $\exp(X_1)=\theta(t_1)$. For $n\in\N$, let $X_n\in U'$ be such that $\exp(X_n)=\theta({t_1}/{n})$.
	
	We claim that $nX_n\in U'$. Indeed, fix $n$ and assume by induction that $jX_n\in U'$ for some $j<n$.
	On the one hand, we have
	 $(j+1)X_n = \frac{j+1}{j}j X_n \in 2U' = U$, since $\frac{j+1}j<2$.
	On the other hand, recalling Exercise \ref{ex:exp:mx}, which we will use several times in this proof, we have
	\[
	\exp((j+1)X_n) = (\exp(X_n))^{j+1} =\left( \theta\left(\frac{t_1}n\right)\right)^{j+1} = \theta\left(\frac{j+1}n t_1\right),
	\]
	which is in $\exp(U')$ since $\frac{j+1}n\le 1$.
	Thus $(j+1)X_n \in U$ and $\exp((j+1)X_n) \in \exp(U')$.
	Since $\exp$ is injective on $U$, we get $(j+1)X_n\in U'$.
	By induction, we get the claim.
	
	Now since 
	\[
	\exp(nX_n) = (\exp(X_n))^n = \left(\theta\left(\frac{t_1}n\right)\right)^n = \theta(t_1) = \exp(X_1)
	\]
	and since $\exp$ is injective on $U'$, we get $X_n=\frac1n X_1$.
	
	Let $\frac mn$ be a rational number with $m\in \Z$ and $n>0$.
	Then
	\[
	\theta\left(\frac mn t_1\right) = \theta\left( \frac{t_1}n \right)^m
	= \exp(X_n)^m = \exp\left(\frac{X_1}n\right)^m 
	= \exp\left(\frac mn X_1 \right).
	\]
	By continuity of $\theta$ (and density of rational numbers), we get  
	\[
	\theta(st_1)=\exp(sX_1), \qquad\forall s\in \R. 
	\]
	Hence, the homomorphism $\theta$ is the exponential curve
	\[
	\theta(t) = \exp\left( \frac t{t_1}X_1 \right),
	\]
	 and thus, it is smooth.
\end{proof}

\begin{theorem}\label{thm: continuous: smooth}\index{continuous homomorphisms}
 	Let $G$ and $H$ be Lie groups. Every continuous homomorphism $\varphi: G\to H$ is smooth.
\end{theorem}
\begin{proof}
 	Consider $\beta: \R^n\to G$ the exponential local coordinates of the second kind with respect to a basis $X_1, \dots, X_n$, i.e.,
	\(
	\beta(t_1, \dots,t_n) := \exp(t_1X_1)\cdots\exp(t_nX_n) .
	\)
	For each $j\in\{1, \dots,n\}$, the map $t\in \R\mapsto \varphi(\exp(tX_j))$ is a continuous homomorphism, so by the previous theorem (Theorem \ref{thm: 1Dcontinuous: smooth}), it is smooth. 
	Hence, the map $\varphi\circ\beta$ is smooth, being the product of smooth maps:
	\begin{eqnarray*}
	\varphi\circ\beta(t_1, \dots,t_n) 
	&=& \varphi(\exp(t_1X_1)\cdot\dots\cdot\exp(t_nX_n)) \\
	&=& \varphi(\exp(t_1X_1))\cdot\dots\cdot\varphi(\exp(t_nX_n)) .
	\end{eqnarray*}
	Hence, in a neighborhood of $1_G$ the map $\varphi=(\varphi\circ\beta)\circ\beta^{-1}$ is smooth. 
	Let $g\in G$ be an arbitrary point. 
	Since $\varphi = L_{\varphi(g)}\circ\varphi\circ L_{g^{-1}}$ and $\varphi$ is smooth at $1_G$, then $\varphi$ is smooth at $g$.
\end{proof}
\begin{corollary}[Uniqueness of Lie structures]\label{Uniqueness of Lie structures}\index{differentiable structure!uniqueness of -- on Lie groups}
 	Each Lie group has only one differentiable structure of a Lie group with the same topology and group structure.
\end{corollary}
\begin{proof}
 	If the identity map between the group with two differentiable structures is a continuous homomorphism, then it is a diffeomorphism by Theorem~\ref{thm: continuous: smooth}.
\end{proof}

\subsection{Closed subgroups of Lie groups}\label{sec:Closed subgroups}
\index{closed subgroup}
\index{Lie! regular -- subgroup}
\index{regular Lie subgroup}
\begin{theorem}[Closed subgroups are regular]\label{teo1131}
 	Let $G$ be a Lie group and $H<G$ a subgroup.
	If $H$ is a closed subset of $G$, then $H$ is a regular Lie subgroup of $G$.
\end{theorem}
Before the proof of the above theorem, we present a preparatory lemma that will be used twice in the proof of the theorem.
\begin{lemma}\label{lem1131}
 	Let $H<G$ be closed.
	Let $X\in\Lie(G)$ and for each $j\in\N$ let $X_j\in\Lie(G)$ and $t_j>0$. 
	Assume $X_j\to X$ and $t_j\to 0$.
	If $$\exp(t_jX_j)\in H, \qquad \forall j\in\N,$$ then $$\exp(tX)\in H, \qquad \forall t\in\R.$$
\end{lemma}
\begin{proof}
 	Fix $t\in\R$. For all $j$ take $m_j\in\left[\frac t{t_j}-1, \frac t{t_j}\right]\cap\Z$.
	Hence, $m_jt_j\to t$ and $m_jt_jX_j\to tX$.
	Since $H$ is closed,
	\begin{equation*}
	 	\exp(tX) = \lim_{j\to\infty} \exp(m_jt_jX_j) 
		= \lim_{j\to\infty} (\exp(t_jX_j))^{m_j} \in H, 
	\end{equation*}
	where we used Exercise \ref{ex:exp:mx}.
\end{proof}

\begin{proof}[Proof of Theorem \ref{teo1131}]
 	We need to show that $H$ is an embedded sub\-ma\-ni\-fold. 
	Let $\g=T_{1_G}G$. Set 
	\[
	\h := \left\{X\in\g: \exists\sigma: \R\to H\text{ such that }\sigma'(0)=X \right\} .
	\]
	The subset $\h$ is a vector space:
	for all $X_1, X_2\in\h$ and all $\lambda_1, \lambda_2\in\R$ we have $\lambda_1X_1+\lambda_2X_2\in\h$, since if $\sigma_j: \R\to H$ with $\sigma'_j(0)=X_j$ for $j\in\{1,2\}$, then $\sigma_1(\lambda_1 t) \cdot \sigma_2(\lambda_2 t)\in H$ and, recalling formula \eqref{derivative:product: curves} of Exercise \ref{ex:derivative:product: curves},
	\[
	\left.\frac{\dd}{\dd t}(\sigma_1(\lambda_1 t) \cdot \sigma_2(\lambda_2 t))\right|_{t=0} = \lambda_1X_1+\lambda_2X_2 .
	\]
	
	We claim now that 
	\begin{equation}\label{eq1131}
	\h = \left\{X\in\g: \exp(tX)\in H\ \forall t\in\R\right\} .
	\end{equation}
	Indeed, it is obvious that $\h$ contains the right-hand side.
	For the opposite inclusion, we shall use Lemma \ref{lem1131}.
	If $X\in\h$, i.e., $X=\sigma'(0)$ for some $\sigma: \R\to H$, we set $\tau(t):=\exp^{-1}(\sigma(t))\in\g$, which is defined when $t$ is sufficiently small
	(recall Proposition~\ref{exp:diffeo}).
	Therefore, 
	\begin{eqnarray*}
	 	X&=&\sigma'(0)
		= \left.\frac{\dd}{\dd t}\sigma(t)\right|_{t=0} 
		= \left.\frac{\dd}{\dd t}\exp(\tau(t))\right|_{t=0} \\
		&=& (\dd\exp)_{\tau(0)} \frac{\dd}{\dd t}\tau(t)|_{t=0} 
		= \tau'(0)
		= \lim_{j\to\infty} j\tau(\frac1j) .
	\end{eqnarray*}
	Set $t_j:=\frac1j$ and $X_j:=j\tau(\frac1j)$, so that $\exp(t_jX_j)=\exp(\tau(\frac1j)) = \sigma(\frac1j)\in H$.
	By Lemma \ref{lem1131}, we conclude that $\exp(tX)\in H$ for all $t\in\R$, i.e., the missing inclusion of the claim \eqref{eq1131} is proved.
	
	The idea is now to use the exponential map to obtain a coordinate map in a neighborhood of $1_G$ in $H$ and then use left translations to get an atlas for $H$. 
	Let $V\subset\g$ be a vector subspace such that $\g=\h\oplus V$.
	We plan to show
	\begin{equation}\label{eq1155}
	 	\exists\Omega\subset V\text{ neighborhood of $0$ in $V$ such that }\exp(\Omega)\cap H=\{1_G\}.
	\end{equation}
	If this is not the case, there are $(Y_j)_{j\in\N}\subset V\setminus\{0\}$ with $Y_j\to 0$ and $\exp(Y_j)\in H$ for all $j\in\N$.
	Let $\|\cdot\|$ be any norm on $\g$.
	Set $t_j:=\|Y_j\|$ and $X_j:=\frac1{t_j}Y_j$. 
	So $X_j$ are unit vectors in $V$ and, up to passing to a subsequence, we have $X_j\to X$ for some $X\in V$.
	Notice that $t_j\to 0$ and $\exp(t_jX_j)=\exp(Y_j)\in H$ for all $j\in\N$.
	By Lemma~\ref{lem1131} we get $\exp(tX)\in H$ for all $t\in\R$. So $X\in\h$.
	Since $\|X\|=1$, we get a contradiction:
	\[
	0\neq X\in V\cap\h = \{0\}. 
	\]
	Hence, equation \eqref{eq1155} is proven.
	
	Let $\psi: \h\times V\to G$ be defined by 
	\[
	\psi(X, Y) := \exp(X)\cdot \exp(Y) , \qquad \forall X\in \h, \forall Y\in V.
	\]
	Then $\psi$ is a diffeomorphism between some neighborhood of $\Omega_1\times\Omega_2$ of $(0,0)$ in $\h\times V$ and some neighborhood $\Omega_3$ of $1_G$ in $G$ (since $\psi$ is an exponential local coordinate map of a mixed kind, c.f. Section \ref{sec1120}).
	We may assume that $\Omega=\Omega_2$, i.e., $\exp(\Omega_2)\cap H=\{1_G\}$.
	
	We plan to show that 
	\begin{equation}\label{eq1124}
	\Omega_3 \cap H = \exp(\Omega_1) .
	\end{equation}
	We have $\exp(\Omega_1)\subset\Omega_3$ by construction.
	Also $\exp(\Omega_1)\subset H$ follows from $\Omega_1\subset \h$ and the claim \eqref{eq1131}.
	Hence $\exp(\Omega_1)\subset\Omega_3\cap H$.

	Vice versa, let $h\in\Omega_3\cap H$, then there are unique $X\in\Omega_1$ and $Y\in\Omega_2$ such that $h=\exp(X)\cdot\exp(Y)$.
	Thus $\exp(Y)=\exp(-X)\cdot h \in\exp(\Omega_2)\cap H = \{1_G\}$.
	Therefore $Y=0$ and $h=\exp(X)\in\exp(\Omega_1)$.
	
	Then $\varphi:=(\psi|_{\Omega_1\times\Omega_2})^{-1}: \Omega_3\to \Omega_1\times\Omega_2$ 
	is a coordinate map for $G$ that is centered at $1_G$ and is adapted to $H$, i.e., $\varphi|_H$ is a coordinate map for $H$ into the vector space $ \h\subset\g$.
	
	To conclude, we consider the atlas $\{(L_h(\Omega_3), \varphi\circ L_{h^{-1}})\}_{h\in H}$.
\end{proof}

\section{General Linear Groups, its Lie algebra, and its exponential map}

The General Linear Group, denoted as $\GL(n, \mathbb K)$, consists of invertible $n\times n$ matrices over a given field $\mathbb K$. Its associated Lie algebra, denoted as $\gl(n, \mathbb K)$, consists of the set of all $n\times n$ matrices equipped with the commutator bracket operation.
The exponential map, defined on the Lie algebra, provides a way to exponentiate matrices and obtain elements in the General Linear Group. It plays a crucial role in Lie theory and connects the algebraic structure of the Lie algebra with the geometric properties of the Lie group. 

%
In our study, it is essential to work with finite-dimensional real vector spaces that are not explicitly identified with $\mathbb{R}^n$. 
Consequently, we consider general linear groups for vector spaces, i.e., sets of linear automorphisms. This abstraction enables us to consider structures like $\gl(n,\R)$ itself or, very importantly, the Lie algebra associated with a Lie group.

Throughout this chapter, all the vector spaces under consideration are defined over the field of real numbers. Similarly, the matrices we examine possess real coefficients. 

\subsection{$\GL(V)$ and $\gl(V)$}\label{sec_GL:g}
 	For $n\in \N$, we denote by
$
 \Mat_{n\times n}(\R) $ the space of $n\times n$ matrices with real entries. \index{$\Mat_{n\times n}(\R)$}
 
The $n$-th \emph{general linear group} is
	\index{general linear group, $\GL$}
	\index{$\GL(n, \R)$}\index{$\GL(n)$}
	\index{group! general linear --, $\GL$}
	\[
	\GL(n):=\GL(n, \R) := \{\text{$A\in \Mat_{n\times n}(\R)$ $\;| $ $\det A\neq0$}\}.
	\]
	This is a group when equipped with the row-column product of matrices.
	Slightly more generally, if $V$ is a (real) vector space, then 
		\[
	\GL(V) := \Aut(V):=\{ A:V\to V \;| \, A \text{ is an invertible linear transformation}\}.
	\]
		This is a group when equipped with the composition rule where the identity element is the identity transformation $\mathbb{I}: V\to V$. Noting that this product rule and the inversion rule are smooth, we infer that $\GL(n, \R) $ and $\GL(V)$ are Lie groups, assuming that $V$ is finite-dimensional. 
		Indeed, the Lie group $\GL(V)$ is Lie group isomorphic to $\GL(n, \R)$ for $n:=\dim(V).$

 	For $n\in \N$, we define
\index{$\gl(n, \R)$}\index{$\mathfrak {gl}(n)$}
\[
\gl(n):=\gl(n, \R) := \Mat_{n\times n}(\R) := \{\text{ $n\times n$ matrices with real entries}\}.
\]
If $V$ is a vector space, then 
\[
\gl(V) := \End(V) := \{\text{ linear transformations from $V$ to $V$ }\}.
\]
Clearly, we have $\GL(n, \R) = \GL(\R^n) $ and $\gl(n, \R) =\gl(\R^n) $.

For $A,B\in\gl(n, \R)$, with $n\in \N$, or, more generally, for $A,B\in\gl(V)$ for a vector space $V$, we consier the {\em commutator}:\index{commutator} 
\begin{equation}\label{del_bracket_gln}
[A, B]:= A B-B A.
\end{equation}
Such an operation is a Lie bracket that makes $\gl(V)$ into a Lie algebra.
And, as the choice of name suggests, this Lie algebra is the Lie algebra of $\GL(V)$; see Proposition~\ref{Prop_LieGL}.

\subsection{Matrix exponential}\label{sec Matrix exponential}

\index{matrix! exponential of a --}
\index{exponential! -- of a matrix}

In this subsection, we recall the matrix exponential: exponential of matrices. Since we shall consider linear endomorphisms of vector spaces, like, for example, the Lie algebra of a Lie group, we define the matrix exponential on the space $\gl(V)$.

\begin{definition}[Matrix exponential]
Let $V$ be a finite-dimensional vector space.
For each $A\in\gl(V) $, define {\em the matrix exponential of $A$} as
\begin{equation}
e^A := \mathbb{I} + A+ \frac12 A^2 + \frac1{3!}A^3+\dots = \sum_{k=0}^\infty \frac1{k!} A^k.
\end{equation}
\end{definition}
In fact, the series giving $e^A$ is absolutely converging; see Exercise~\ref{rmk:exp:converge}. 
Consequently, the function $A \mapsto e^A$ is smooth (in fact, analytic). Moreover, each $e^A$ is invertible with inverse $e^{-A}$; see Exercises~\ref{ex: abba} and~\ref{ex_eA_invertible}, so we have a map $A\in \gl(V) \mapsto
e^A \in \GL(V)$.

in the following discussion we use that $\GL(V)\subseteq \End(V)$ and $\End(V)$ is a vector space. Therefore, tangent vectors of curves into $\GL(V)$ are represented by elements in $\End(V)$.

It is easy to see (and a proof is in Proposition~\ref{prop: der: eta}) that for every linear map $A $, the curve $t\mapsto e^{tA}$ satisfies 
 $$\left.e^{tA} \right|_{t=0} =\mathbb{I}, \, \qquad \qquad \left. \frac{\dd }{\dd t} e^{tA}\right|_{t=0} =A.$$

Moreover, the map $\phi ^t$ that to every linear map $B$ on $V$ associates $\phi ^t(B):= Be^{tA}$ satisfies the following properties:
\begin{itemize}
\item $\phi ^t$ is a flow, i.e., $\phi ^t \circ \phi ^s = \phi ^{t+s}$ because for every $B$ we have that $(Be^{sA})e^{tA} = Be^{(t+s)A}$;
\item $\phi ^t$ is left invariant, i.e., $\phi ^t (MB)=M\phi ^t (B) $ because $(MB) e^{tA} = M(B e^{tA})$.
\end{itemize}
 Hence, this flow is the flow of its derivative at $0:$
 \begin{equation*}
\begin{aligned}
\left.\frac{\dd }{\dd t} \phi ^t (B)\right|_{t=0} = \left.\frac{\dd }{\dd t} Be^{tA}\right|_{t=0} = B \left.\frac{\dd }{\dd t} e^{tA}\right|_{t=0} =BA.
\end{aligned}
\end{equation*}

We summarize, see also Proposition~\ref{prop: der: eta}, the basic properties of the matrix exponential:

\begin{proposition}[Matrix exponential]\label{prop_matrix_exp}
Let $V$ be a finite-dimensional vector space.
\begin{description}
\item[\ref{prop_matrix_exp}.i.] The matrix exponential
\begin{equation*}
\begin{aligned}
\exp : \mathfrak {gl}(V) & \to \GL(V)\\
 A & \mapsto e^A,
\end{aligned}
\end{equation*}
 is an analytic map.
 \item[\ref{prop_matrix_exp}.ii.] For every $A \in \mathfrak {gl}(V),$ the curve $t \mapsto e^{tA}$ is a one-parameter subgroup.
 \item[\ref{prop_matrix_exp}.iii.] For every $A \in \mathfrak {gl}(V),$ the map 
 \begin{equation*}
\begin{aligned}
GL(V) & \to T\left(GL(V)\right)\\
 B & \mapsto BA,
\end{aligned}
\end{equation*}
 defines a left-invariant vector field on $\GL(V)$ whose flow $\R \times \GL(V) \to \GL(V)$ is defined by $(t,B ) \mapsto Be^{tA}.$ 
\end{description}
\end{proposition}

Rephrasing when $V=\R^n$, we have that for all $A\in \mathfrak {gl}(n, \R)\simeq T_\mathbb{I}GL(n, \R)$, the unique LIVF on $\GL(n, \R)$ that equals $A$ at $\mathbb{I}$ is 
$$B \in \GL(n, \R) \longmapsto BA\in \Mat_{n\times n}(\R).$$

In the next proposition, we spell out the argument that shows what is the derivative of the OPS $t\mapsto e^{tA}$. We shall refer to this proposition several times.
\begin{proposition}[Derivative of $e^{tA}$]\label{prop: der: eta}
\index{derivative! -- of $e^{tA}$}
For every finite-dimensional vector space $V$ and every $A\in \End(V)$,
	the curve $t\mapsto e^{tA}$ is the one-parameter subgroup of $\GL(V)$ such that
	\[
	\frac{\dd}{\dd t} e^{tA} = A e^{tA}
	\qquad\text{
	and
	}\qquad
	\left.\frac{\dd}{\dd t}e^{tA}\right|_{t=0} = A.
	\]
\end{proposition}
\begin{proof}
 	Recall that $A\mapsto e^A$ is smooth, that $e^{sA}\cdot e^{tA} = e^{(s+t)A}$, and that $e^0=\mathbb{I}$.
	Therefore $t\mapsto e^{tA}$ is a one-parameter subgroup of $\GL(V)$.
	For the last two claims, we have
	\begin{eqnarray*}
	 	\frac{\dd}{\dd t} (e^{tA}) 
		&=& \frac{\dd}{\dd t} \left(\sum_{k=0}^\infty \frac1{k!} (tA)^k \right) \\
		&=& \sum_{k=0}^\infty \frac1{k!} \frac{\dd}{\dd t}(t^kA^k) \\
		&=& \sum_{k=1}^\infty \frac1{k!} k t^{k-1} A^k \\
		&=& A\sum_{k=1}^\infty \frac1{(k-1)!} t^{k-1} A^{k-1} \\
		&=& A e^{tA} ,
	\end{eqnarray*}
where we could exchange derivative and summation because the series is absolutely convergent.
\end{proof}


\subsection{Lie algebras of general linear groups}\index{Lie algebra! -- of general linear group}
This section aims to show that $ \mathfrak {gl}(n, \R)$ when equipped with \eqref{del_bracket_gln} is the Lie algebra of the Lie group $ \GL(n, \R)$.

\begin{proposition}\label{Prop_LieGL}
The Lie algebra of $\GL(V)$ is isomorphic to the Lie algebra $\mathfrak {gl}(V).$ 

\end{proposition}

\begin{proof}
Recall from \ref{prop_matrix_exp}.iii that every element in $ \mathfrak {gl}(V)$ induces a left-invariant vector field on $\GL(V)$ for which we have a formula for the flow.
The key point of the proof is to show that for every $A, B \in \mathfrak {gl}(V)$ the vector field
$M \in \GL(V) \mapsto M(AB-BA)$
is the Lie bracket between 
\begin{equation*}
M \mapsto MA \qquad \mbox{ and } \qquad M \mapsto MB.
\end{equation*}
Thus, in terms of flows, by recalling \ref{def:Lie_bracket_vector_fields}.d, we need to show that
\begin{equation}\label{flow_uguali}
 \left.\frac{\dd}{\dd t} \phi ^t_{AB-BA} (M)\right| _{{t=0}} = \left.\frac{\dd}{\dd t} ( \phi _B^{-\sqrt t} \circ \phi _A^{-\sqrt t} \circ \phi _B^{\sqrt t} \circ \phi _A^{\sqrt t}) (M)\right|_{{t=0}}.
\end{equation}
On the one hand, the left-hand side is 
\begin{equation*}
 {\rm LHS\, of \,} \eqref{flow_uguali} = \left.\frac{\dd}{\dd t} Me^{t(AB-BA)} \right|_{{t=0}} = M(AB-BA).
\end{equation*}
On the other hand, recalling that $e^{\sqrt t A}=\mathbb{I} +\sqrt t A +\frac{tA^2} 2+ o(t),$ the right-hand side becomes
\begin{equation*}
\begin{aligned}
 {\rm RHS\, of \,} \eqref{flow_uguali} &=\left. \frac{\dd}{\dd t} M ( e^{\sqrt t A} e^{\sqrt t B} e^{-\sqrt t A} e^{-\sqrt t B})\right|_{{t=0}}\\
&= \left.\frac{\dd}{\dd t} M ( e^{\sqrt t A} e^{\sqrt t B} e^{-\sqrt t A} e^{-\sqrt t B})\right|_{{t=0}}\\
&= \frac{\dd}{\dd t} M \Bigg(\mathbb{I} +\sqrt t (A+B-A-B) \\
&\left. \qquad +t \left(\frac{A^2} 2 + \frac{B^2} 2+ \frac{A^2} 2 + \frac{B^2} 2+AB-A^2 -AB-BA-B^2+AB \right)+ o(t) \Bigg)\right|_{{t=0}}\\
&= M(AB-BA).
\end{aligned}
\end{equation*}
Hence, equation \eqref{flow_uguali} holds, as desired.
\end{proof}

Proposition~\ref{prop: der: eta}, together with Proposition~\ref{Prop_LieGL}, therefore clarified that 	the exponential of $\GL(n, \R)$ is the usual exponential of matrices $\exp:A\in \gl(n, \R) \mapsto e^A\in\GL(n, \R)$.
\begin{corollary}[of Proposition~\ref{prop: der: eta}]
\index{matrix! exponential of a --}
\index{exponential! -- of a matrix}
	For every finite-dimensional vector space $V$, the exponential map of the Lie group $\GL(V)$ is the matrix exponential $\exp: \gl(V)\to\GL(V)$, $A\mapsto e^A$.
\end{corollary}

%
%
%
%
%
%
%

\section{Adjoint representation}
\index{endomorphism! --s of a vector space, $\gl(V)$}
\index{$\gl(V)$}
The adjoint representation, also known as the adjoint action, of a Lie group $G$ provides a means of representing, possibly not injectively, the elements of the group as linear transformations of its Lie algebra, viewed as a vector space. Specifically, in the case of the general linear group $\GL(n, \mathbb{R})$, where the operations are linear, the adjoint representation corresponds to conjugation.
To obtain the adjoint representation for a Lie group, we linearize the group's action on itself through conjugation. Namely, we take the differentials. This natural representation captures the way the elements of the Lie group act on its Lie algebra, establishing a link between the group's abstract structure and the associated linear transformations.

\subsection{$\Ad$ and $\ad$}\label{sec Ad and ad}

 In this section, we shall consider a Lie algebra $\g$ as a vector space and then consider the spaces 
 $\gl(\g)$ and $\GL(\g)$, as in Section~\ref{sec_GL:g}.

\begin{definition}[Adjoint map]\label{def_ad}
 	Let $\g$ be a Lie algebra. The \emph{adjoint map} of $\g$ is the linear map
	\index{adjoint! -- map, $\ad$} 
	\[
	\ad: \g\to\gl(\g)
	\]
	given by
	\[
	\ad(X)(Y) := \ad_X(Y) := [X, Y], \qquad \forall X, Y\in\g.
	\]
\end{definition}
\begin{remark}\label{ad_homo}
 	The map $\ad_X: \g\to\g$ is indeed in $\gl(\g)$, i.e., it is linear, not necessarily invertible.
	Moreover, seeing $\gl(\g)$ as a Lie algebra, the map $\ad: \g \to\gl(\g) $ is a Lie algebra homomorphism: 
	for all $X, Y\in\g$ and all $s,t\in\R$ we have
	\begin{description} 
	\item[\ref{ad_homo}.i.] 	$\ad(sX+tY) = s\ \ad(X) + t\ \ad(Y)$,
	\item[\ref{ad_homo}.ii.] 	 	$\ad([X, Y]) = [\ad(X), \ad(Y)]$,
	\end{description}
see Exercise~\ref{ex_ad_homo} for the proof.
\end{remark}

\begin{definition}[Adjoint representation]
 	Let $G$ be a Lie group with Lie algebra $\g$.
	For $g\in G$ define
	\[
	\Ad(g) := \Ad_g := (\dd C_g)_{1_G},
	\]
	i.e., $\Ad(g): \g \to \g$ is the differential at the identity element of the conjugation $C_g:h\mapsto ghg^{-1}$.
	The map 
	\[
	\Ad: G\to \GL(\g)
	\]
	is called \emph{adjoint representation}.
	\index{adjoint! -- representation, $\Ad$}
	\index{conjugation}
\end{definition}
\begin{remark}\index{representation! -- of a group}\index{group! -- representation} 
The map $\Ad$ is indeed a representation, i.e., $\Ad$ is a group homomorphism into $\GL(\g)$:
\begin{equation}\label{ex_Ad_homo2}
	\Ad(gh) 
	= \Ad(g)\circ\Ad(h), \qquad \forall g,h\in G, 
	\end{equation}
see Exercise~\ref{ex_Ad_homo} for the proof.
\end{remark}

\subsection{Properties and formulas}
\begin{proposition}
 	Let $G$ be a Lie group with Lie algebra $\g$.
	The adjoint representation $\Ad: G\to\GL(\g)$ is a Lie group homomorphism, and the Lie algebra homomorphism associated with $\Ad$ is the adjoint map $\ad: \g\to\gl(\g)$, i.e., 
	\begin{equation}\label{eq Ad ad}
	(\Ad)_*=\ad,
	\end{equation}
i.e., we have $
	(\dd \Ad)_{1_G} (X) = \ad(X), $ for all $X\in\g .
	$
\end{proposition}
\begin{proof}
 	Since $t\mapsto\exp(tX)$ is a curve in $G$ that is tangent to $X$ at ${1_G}$, we have 
	\[
	(\dd\Ad)_{1_G}(X)(Y) 
	= \left.\frac{\dd}{\dd t} \Ad\left(\exp(tX)\right)(Y) \right|_{t=0}.
	\]
	Here, we consider $X, Y$ as elements in $T_{1_G}G$.
	We denote by $\tilde X, \tilde Y$ the left-invariant vector fields such that $\tilde X_{1_G}=X$ and $\tilde Y_{1_G}=Y$.
	We have
	\begin{eqnarray*}
	 	\Ad(\exp(tX))(Y) 
		&=& \left(\dd C_{\exp(tX)}\right)_{1_G} (Y) \\
		&=& \left(\dd R_{\exp(-tX)} \right)_{\exp(tX)}(\dd L_{\exp(tX)})_{1_G} (Y) \\
		&=& \left(\dd R_{\exp(-tX)} \right)_{\exp(tX)}(\tilde Y_{\Phi^t_{\tilde X}({1_G})})\\
		&=& \left(\dd \Phi^{-t}_{\tilde X}\right)_{\Phi^t_{\tilde X}({1_G})} \tilde Y_{\Phi^t_{\tilde X}({1_G})},
	\end{eqnarray*}
	where we used that the flow at time $t$ of the left-invariant vector field $\tilde X$ is the right translation by $\exp(tX)$; see Corollary~\ref{prop:right_left}.
	We get
	\begin{eqnarray*}
	 	(\dd\Ad)_{1_G} (X)(Y)
		&=& \left. \frac{\dd}{\dd t} \left( \dd\Phi^{-t}_{\tilde X} \right)_{\Phi_{\tilde X}^t({1_G})} \tilde Y_{\Phi^t_{\tilde X}({1_G})} \right|_{t=0} \\
		&\overset{\text{def.~of}}{\underset{\text{Lie deriv.}}=}& \left(\Lie_{\tilde X}(\tilde Y)\right)_{1_G}\\
		&=& [\tilde X, \tilde Y]_{1_G}\\
		&=& [X, Y]\\
		&=& \ad_X(Y) .
	\end{eqnarray*}\qedhere
\end{proof}
Recall that if $\varphi$ is a Lie group homomorphism and $\varphi_*$ is the Lie algebra homomorphism induced by $\varphi$, by Proposition~\ref{prop: expF: Fexp} we have the following first commutative diagram, and for $\varphi=C_g$ (resp.~$\varphi= \Ad$) we have the following second (resp.~third) commutative diagram.
%
\[
\xymatrix{
\g\ar[d]_{\exp}\ar[r]^{\varphi_*} & \h\ar[d]^{\exp} \\
G\ar[r]_{\varphi} & H
}
\qquad\qquad
\xymatrix{
\g\ar[d]_{\exp}\ar[r]^{(C_g)_* = \Ad_g} & \g\ar[d]^{\exp} \\
G\ar[r]_{C_g} & G
}
\qquad\qquad
\xymatrix{
\g\ar[d]_{\exp}\ar[r]^{(\Ad)_*=\ad} & \gl(\g)\ar[d]^{\exp = e^\cdot} \\
G\ar[r]_{\Ad} & \GL(\g)
}
\]
 \begin{formula}\label{Formula:Cg:Ad}
	Since $\Ad_g=(C_g)_*$ by definition, we have
	\[
	C_g(\exp(X)) = \exp(\Ad_gX), \qquad\forall X\in \g, \forall g\in G.
	\]
	Equivalently, 
		\[
	\exp(Y)\exp(X) \exp(-Y) = \exp(\Ad_{\exp(Y)}X), \qquad\forall X, Y\in \g.
	\]
	\end{formula}
 
 \begin{formula}\label{Formula: Ad: ad}
	Since $(\Ad)_* = \ad$ by the previous proposition, we have
	\[
	\Ad_{\exp(X)} = e^{\ad_X} .
	\]
	\end{formula}

 	In the above formula, $\ad_X$ is a linear transformation on $\g$, i.e., an element of $\gl(\g)$, which is the Lie algebra of $\GL(\g)$.
	We saw that $\exp: \gl(\g)\to\GL(\g)$ is given by the classical matrix exponential.
	Therefore
	\begin{align*}
	e^{\ad_X}(Y) &= \sum_{k=0}^\infty \frac{(\ad_X)^k}{k! }(Y) \\
	&= Y + [X, Y] + \frac12 [X,[X, Y]] + \frac1{3! } [X,[X,[X, Y]]] + \dots
	\end{align*}

 \begin{formula}[See Exercises~\ref{Ex:ad:BAB}, \ref{Formulae_GL2}, and~\ref{Formulae_GL3}]\label{Formulae_GL}
 Let $V$ be a vector space. For all $X, Y\in \mathfrak{gl}(V)$ and $B\in\GL(V)$ we have
\begin{description}
\item[\ref{Formulae_GL}.i.]	$
	\Ad_B(X) = B\cdot X\cdot B^{-1};
	$
\item[\ref{Formulae_GL}.ii.] $	 e^{\operatorname{ad} _X} Y = e^{X}Y e^{-X} $;
\item[\ref{Formulae_GL}.iii.]	
	$
	e^{BXB^{-1}} = Be^XB^{-1} .
	$
	\end{description}
	 \end{formula}

\begin{formula}[See Exercise~\ref{ex_fomula_flow_ad}]
	For every left-invariant vector fields $X, Y$ on a Lie group $G$, we have that $(\Phi^t_X)_*Y$ is a left-invariant vector field and 
\[
(\Phi^t_X)_*Y = e^{-\ad(tX)}Y, \qquad \forall t\in\R.
\]

 \end{formula}

\section{Semi-direct products}
In this section, our focus turns to the study of the semidirect product of Lie groups. 
We begin by considering group actions by group automorphisms, as well as Lie actions by derivations. These actions enable us to form semidirect products of Lie groups and, in a parallel manner, semidirect products of Lie algebras.
A key result will emerge: the Lie algebra of a semidirect product of Lie groups corresponds to a semidirect product of the Lie algebras associated with the constituent groups. 

In general, given a group \( G \) and a set \( X \), an \emph{action} of \( G \) on \( X \) is a map \( \Theta \colon G \times X \to X \) such that $\Theta(1_G,p)=p$, for all $p\in X$, and 
\begin{equation} \label{eq:def_left_action}
\Theta(g_1, \Theta(g_2,p)) = \Theta(g_1g_2,p), \qquad \forall g_1, g_2 \in G, \forall p\in X.
\end{equation}
 \index{group! -- actions}
\index{action!  group --}
\index{action!  left --}
\index{left! -- action}
\index{action}
Actions, as we have just defined, are also called {\em left actions} or {\em group actions}.
We write \( G \acts X \) in case of an action, and we write \( g.p \) for \( \Theta(g,p) \). The action defines a group homomorphism 
$G \to \mathfrak{S}(X) $, $ g \mapsto \Theta(g, \cdot) ,$
into the group $\mathfrak{S}(X)$ of permutations of $X$. 
Actions will be extensively considered again in Section~\ref{sec_actions}.





\subsection{Derivations and actions by automorphisms}

\begin{definition}[Group action by automorphisms]
 Let $H $ and $ G$ be groups.
An \emph{action} of $H$ {\em by automorphisms} of \( G \) 
is a group homomorphism $\theta:H\to \Aut(G)$.
Equivalently, it is a map $\theta:H\times G \to G$
such that, for $\theta_h:=\theta(h, \cdot)$, we have
$$\theta_h(g_1 g_2) = \theta_h(g_1)\theta_h(g_2), \qquad \forall g_1, g_2\in G, \forall h\in H$$
and
\begin{equation}\label{theta_homo}\theta_{h_1 h_2} = \theta_{h_1}\circ \theta_{h_2}, \qquad \forall h_1,h_2\in H.
\end{equation}
In case $G$ and $H$ are Lie groups, we say that an action $\theta: H\to \Aut(G)$ is {\em smooth} if it is smooth as a map $\theta: H\times G \to G$.
\end{definition}
Recall that every smooth element $\varphi\in\Aut(G)$ has an associated Lie algebra automorphism $\varphi_*$. 
Let us introduce a notation for the
 {\em space of Lie group automorphisms} of a Lie group $G$:\index{Lie group! space of -- automorphisms}
 $${{\Aut_{\Lie}}}(G) := \Aut(G)\cap C^\infty(G;G) = \{ \varphi : G\to G \text{ smooth automorphism}\} $$
 and the {\em space of Lie algebra automorphisms} of a Lie algebra $\g$:
$${{\Aut_{\Lie}}}(\g) := \{T \in\GL(\g) \;: \;T[u, v] = [Tu, Tv], \forall u, v \in \g\}.$$

We shall stress that ${{\Aut_{\Lie}}}(\g)$ is a closed Lie subgroup of $\GL(\g)$ with a Lie algebra that is a Lie subalgebra of $\gl(\g)$. The elements of this Lie subalgebra are the so-called derivations:

\begin{definition}[Derivation on a Lie algebra]\label{def:derivation}
 Let $\g$ be a Lie algebra.
 A {\em derivation} on $\g$ is a linear map $D: \g \to \g$ that satisfies Leibniz rule:\index{Leibniz rule}
\begin{equation}\label{eq_derivation}
 D([X, Y]) = [D(X),Y] + [X,D(Y)], \qquad \forall X, Y\in \g.
 \end{equation}
 Let ${\rm {Der}}(\g)$ be the {\em set of derivations} on $\g$. \index{derivation}\index{${\rm {Der}}(\g)$}
\end{definition}

We state the relative results for reference (see Exercises \ref{ex_deriv_0}--\ref{ex_deriv_5}):
\begin{proposition}\label{Prop_Lie_der} Let $G$ be a Lie group with Lie algebra $\g$.
\begin{description}
\item[\ref{Prop_Lie_der}.i.] The natural map 
\begin{eqnarray*}
{{\Aut_{\Lie}}}(G) &\longrightarrow & {{\Aut_{\Lie}}}(\g)\\
\varphi &\longmapsto& \varphi_*
\end{eqnarray*}
is an injective Lie group homomorphism, and if $G$ is simply connected, it is a Lie group isomorphism.
\item[\ref{Prop_Lie_der}.ii.] ${{\Aut_{\Lie}}}(\g)$ is a regular Lie subgroup of $\GL(\g)$ with
$$\Lie({{\Aut_{\Lie}}}(\g)) ={\rm {Der}}(\g).$$
\item[\ref{Prop_Lie_der}.iii.] The adjoint map $\ad$ and the adjoint representation $\Ad$ satisfy
\[\xymatrix{
\quad\g\ar[d]_{\exp}\ar[rr]^{ \ad} \quad && \quad {\rm {Der}}(\g)\ar[d]^{\exp} \subset\gl(\g)\\
\quad G\ar[rr]^{\Ad} \quad && \quad {{\Aut_{\Lie}}}(\g)\subset\GL(\g)
}
\]
\end{description}
\end{proposition} 

\subsection{Semi-direct products of Lie algebras and groups}\index{semi-direct product}
\begin{definition}[Semi-direct product of Lie algebras]\label{definition semi-direct Lie algebra}
Let $\g$ and $\h$ be Lie algebras, and let $\sigma: \h \to {\rm {Der}}(\g)$ be a Lie algebra homomorphism into the space of derivations of $\g$. On the direct sum 
$\g\oplus \h$ we consider the bracket that agrees with the brackets of $\g$ and $\h$, and additionally
$$[(0, Y), (X, 0)] := \sigma(Y)(X), \qquad \forall X\in \g, \forall Y\in \h.$$
More explicitly,
\begin{equation}\label{def_semi_direct_algebra}[(X, Y), (X', Y')] := ([X, X'] + \sigma(Y)(X') - \sigma(Y')(X), [Y, Y']), \qquad \forall X, X'\in \g, \forall Y,Y'\in \h.
\end{equation}
The resulting Lie algebra is the {\em semidirect product} of $\g$ and $\h$ with
respect to $\sigma$, and it is denoted by $\g \rtimes_\sigma \h$.
When $\sigma$ is understood, or there is no need to name it, we simply write $\g \rtimes \h$.
\end{definition}
\begin{remark}\label{rmk_semidir_algebra}
We have the following properties for the semi-direct product of Lie algebras.
\begin{description}
\item[\ref{rmk_semidir_algebra}.i.]
We have that $\g \rtimes \h$ is a Lie algebra and the maps
$X\in \g \mapsto (X,0)$ and $Y\in \h\mapsto (0,Y)$ give injective Lie algebra homomorphisms into $\g \rtimes \h$.
\item[\ref{rmk_semidir_algebra}.ii.] If $\sigma\equiv 0$, we call $\g \rtimes \h$ the {\em direct product} of $\g $ and $ \h$, and write it as $\g \times \h$.
\item[\ref{rmk_semidir_algebra}.iii.] In $\g \rtimes \h$, the Lie subalgebra $\g$ is an {\em ideal}, \index{ideal} i.e., $[\g \rtimes \h, \g] \subseteq \g$.
This is the reason for the choice of the symbol $\rtimes$ to resemble $\vartriangleleft$, and we write 
$\g \vartriangleleft \g \rtimes \h$.
If somewhere else you read $\g \ltimes \h$, then it means that in that setting, it is $\h$ that is an ideal and hence it is $\g$ that is acting on $\h$ by derivations.\index{$\rtimes$}
\item[\ref{rmk_semidir_algebra}.iv.] The map $\sigma$ represents the adjoint map in $\g \rtimes_\sigma \h$ of $\h$ on $\g$:
$$\ad_Y(X) = \sigma(Y)(X), \qquad \forall X\in \g, \forall Y\in \h;$$
recall that indeed every $\ad_Y$ is a derivation; see Exercise~\ref{ex_deriv_2}.
\end{description}
\end{remark}

\begin{definition}[Semi-direct product of groups]\label{def Semi-direct product of groups}
 	Let $G$ and $H$ be groups, and $\theta:H\to \Aut(G)$ an action of $H$ by 
automorphisms of $G$.
On the set $ \{(g,h):g\in G, \ h\in H\}$ we put the product
	\begin{equation}\label{prod:semi} 
	(g_1,h_1) \cdot (g_2,h_2) = (g_1\cdot\theta_{h_1}(g_2), h_1 h_2), \qquad \forall g_1, g_2\in G, \forall h_1,h_2\in H .
\end{equation}
The resulting group is the \emph{semi-direct product} of $G$ and $H$ with respect to $\theta$, and it is denoted by $	G\rtimes_\theta H$, or simply $	G\rtimes H$ if there is no need to write $\theta$ explicitly.
	 
\end{definition}

\begin{remark}\label{rmk_semidir_grp} Similarly to Remark~\ref{rmk_semidir_algebra}, we have the following properties for the semi-direct product of groups.
\begin{description}
\item[\ref{rmk_semidir_grp}.i.]
We have that $G \rtimes H$ is a group, and the maps $g\in G \mapsto (g,1_H)$ and $h\in \h\mapsto (1_G,h)$ give injective group homomorphisms into $G \rtimes H$.
\item[\ref{rmk_semidir_grp}.ii.] If $\theta\equiv \id_G$, we call $G \rtimes H$ the {\em direct product} of $G$ and $H$, and write it as $G \times H$.
\item[\ref{rmk_semidir_grp}.iii.] In $G \rtimes H$, the subset $G$ is a normal subgroup.
Thus, we have $G \vartriangleleft G \rtimes H$, which explains the symbol $\rtimes$.
 \item[\ref{rmk_semidir_grp}.iv.] To memorize the definition of the product law, one needs to understand its reason\footnote{A Finnish motto says: 
 {\it sitä, minkä ymmärtää, ei tarvitse muistaa}. 
 [{what you understand, you don't need to remember.}]}.
 Write the product of elements $g_1, g_2\in G $ and $ h_1,h_2\in H$ as
$$g_1h_1 g_2h_2 =g_1h_1 g_2h_1^{-1} h_1 h_2 = g_1 C_{h_1}(g_2) h_1 h_2.$$
In other words, the map $\theta $ represents the conjugation in $G \rtimes H$ of $H$ on $G$:
\begin{equation}\label{formula_conj_semidir}
C_h(g) = \theta_h(g), \qquad \forall g \in G, \forall h \in H;
\end{equation}
recall that indeed every $C_g$ is a group automorphism; see Exercise~\ref{ex_conj_auto}.
\item[\ref{rmk_semidir_grp}.v.] The element $(g,h)$ has inverse $(\theta_{h^{-1}}(g^{-1}),h^{-1})$.
\end{description}
\end{remark}

\subsection{Lie algebras of semi-direct products of Lie groups}
\begin{proposition}\label{Prop_Lie_semidirect}
Let $G$ and $H$ be Lie groups with Lie algebras $\g$ and $\h$, respectively, and $\theta:H\to \Aut(G)$ be a smooth action.
\begin{description}
\item[\ref{Prop_Lie_semidirect}.i.]
$G\rtimes_\theta H$ is a Lie group.
\item[\ref{Prop_Lie_semidirect}.ii.] The map $\tau:H\to {\Aut_{\Lie}}(\g) $ defined by $\tau_h:=(\theta_h)_*$, for $h\in H$, is a Lie group homomorphism.
\item[\ref{Prop_Lie_semidirect}.iii.]
For $\sigma:=\tau_*: \h \to {\rm {Der}}(\g)$, for the above $\tau$, we have
$$\Lie (G\rtimes_\theta H) = \g\rtimes_\sigma \h.$$
\end{description}
\end{proposition}

\begin{proof}
As a manifold, $G\rtimes_\theta H$ is the product of the manifolds $G$ and $H$. Moreover, the group structure is smooth by construction. Thus, $G\rtimes_\theta H$ is a Lie group.

By the chain rule applied to \eqref{theta_homo} we get $\tau_{h_1 h_2} = \tau_{h_1}\circ \tau_{h_2},$ for all $h_1,h_2\in H$. So also \ref{Prop_Lie_semidirect}.ii is proved.

For proving \ref{Prop_Lie_semidirect}.iii we need to calculate the Lie bracket $\ad_Y(X)=[Y, X]$ where $X\in \g $ and $ Y\in \h$. Hence, before calculating $\ad_Y(X)$ we calculate $\Ad_{\exp(Y)}(X)$ and before that $C_{\exp(Y)}(\exp(X))$.
For doing the calculation, we shall crucially use the relation between the exponential map and the induced morphisms; see Proposition~\ref{prop: expF: Fexp}. In fact, we have
\begin{equation}\label{teta_tau}\theta_h(\exp(X)) = \exp(\tau_h(X)), \qquad \forall X\in \g, \forall h\in H 
\end{equation}
and 
\begin{equation}\label{tau_sigma}\tau(\exp(Y)) = \exp(\sigma(Y)), \qquad \forall Y\in \h.\end{equation}
For all $X\in \g $, $ Y\in \h$, and $t\in \R$, we have
\begin{eqnarray*}
\exp(t \Ad_{\exp(Y)}(X)) &=&
\exp(\Ad_{\exp(Y)}(t X)) \\
&\stackrel{\rm F.~\ref{Formula:Cg:Ad}}{=}& C_{\exp(Y)}(\exp(tX)) \\
&\stackrel{\eqref{formula_conj_semidir}}{=}& \theta_{\exp(Y)}(\exp(tX)) \\
&\stackrel{\eqref{teta_tau}}{=}& \exp(\tau_{\exp(Y)}(t X))\\
&=& \exp(t\tau_{\exp(Y)}(X)),
\end{eqnarray*}
where we also used that both $\Ad_{\exp(Y)}$ and $\tau_{\exp(Y)}$ are linear.
We got an identity between OPSs. Therefore $ \Ad_{\exp(Y)}(X) 
= \tau_{\exp(Y)}(X))$, for all $X\in \g$, i.e., $ \Ad_{\exp(Y)} 
= \tau_{\exp(Y)}$. 
Consequently, 
$
e^{\ad_Y}	 \stackrel{\eqref{Formula: Ad: ad}}{=} \Ad_{\exp(Y)} = \tau_{\exp(Y)}
\stackrel{\eqref{tau_sigma}}{=} e^{\sigma(Y)}. $
Differentiating in $Y$ we get $\ad_Y= \sigma(Y)$.
\end{proof}

\begin{remark} Conversely, every semidirect product of Lie algebras is the Lie algebra of a semidirect product of Lie groups. Indeed, let $\g \rtimes_\sigma \h$ be a semidirect product of Lie algebras, and let $G$ and $H$ be simply connected Lie groups with Lie algebras $\g$ and $\h$, respectively, whose existence in ensured by Theorem~\ref{Lie's Third Theorem}. From Theorem~\ref{thm_induced_homo0}, since $H$ is simply connected, there is a Lie group homomorphism $\tau : H \to \Aut_{\Lie} (\g)$ such that $\tau_* = \sigma$. Then, again from Theorem~\ref{thm_induced_homo0}, since $G$ is simply connected, for every $h\in H$ there is a Lie group automorphism $\theta_h : G \to G$ such that $(\theta_h)_* =\tau_h$. Such a map induces a smooth action $\theta: H \to \Aut( G)$.
One can verify that $\Lie (G\rtimes_\theta H) = \g\rtimes_\sigma \h.$
\end{remark}

\section{From algebras to groups}

In this section, we revisit the discussion from Section~\ref{Sec: Lie_objects} regarding the relationship between objects at the level of the Lie algebra and their counterparts at the level of the Lie group.
Two key examples illustrate this relationship: the correspondence between Lie subalgebras and Lie subgroups, and the induction of Lie group homomorphisms from Lie algebra homomorphisms, provided that the Lie group in the source is simply connected.

\subsection{Existence of subgroups}
The next result shows the existence of Lie subgroups with given Lie subalgebras of Lie algebras of Lie groups. Together with Ado's theorem (see \cite[page 199]{Jacobson} and Theorem~\ref{Lie's Third Theorem}), we will deduce that for every abstract Lie algebra (real and finite-dimensional), there exists at least one Lie group with this Lie algebra.
\begin{theorem}[Existence of subgroups]\label{teo1145}\index{group! -- generated}\index{existence! -- of subgroups}
 	Let $G$ be a Lie group with Lie algebra $\g$.
	For every Lie subalgebra $\h\subset \g$, there is a unique connected Lie subgroup $H$ with Lie algebra $\h$; in fact, $H$ is the group generated by $\exp(\h)$. It may not be true that $H=\exp(\h)$.
\end{theorem}
\begin{proof}\index{Frobenius's theorem}
 	This is a consequence of Frobenius's theorem; see for example \cite[page~496]{Lee_MR2954043} or \cite[Sec.~3.7]{Abate_Tovena_geometria_differenziale}.
	We consider the subbundle $\Delta\subset TG$ defined by 
	\[
	\Delta_g := (\dd L_g)_{1}(\h) .
	\]
	Notice that $\Delta$ is left-invariant and involutive (since $\h$ is closed under the bracket).
	Frobenius's theorem implies that there exists a maximal connected sub\-ma\-ni\-fold $H$ of $G$ such that $1_G\in H$ and $TH =\Delta|_H$.
	By construction, $T_1H= \h$.
	We claim that since $\Delta$ is invariant under the maps $\{L_h\}_{h\in G}$, then $H$ is a subgroup. 
	Indeed, take $h_1, h_2\in H$ and observe that $L_{h_1^{-1}} H$ contains $1_G$ and is tangent to $\Delta$. By maximality, we have $h_1^{-1} h_2\in L_{h_1^{-1}} H\subseteq H$.
	
	Regarding uniqueness, if $\hat H$ is a connected subgroup with $\Lie(\hat H)=\h$, since $\exp(\h)$ is an open neighborhood of $1$ in $\hat H$, we have
	\[
	\hat H = (\hat H)^\circ = \langle \exp(\h) \rangle ,
	\]
where we used Exercise~\ref{Prop:generating}. The last assertion of the theorem comes from the fact that there are connected Lie groups with non-surjective exponential map; see Exercise~\ref{ex:exp:no:sujective}.
\end{proof}

\subsection{Existence of group homomorphisms}\index{existence! -- of group homomorphisms}
We shall show that every Lie algebra homomorphism between Lie algebras of Lie groups is induced by a Lie group homomorphism in case the source Lie group is simply connected. Moreover, this group homomorphism is unique. The existence fails in the case where the group is not simply connected. The uniqueness fails as long as the group is not connected.

We will use the following fact: when the Lie algebra homomorphism induced by a Lie group homomorphism is a bijection, then the group homomorphism is a covering map; see Exercise~\ref{prop1302}. For the basics of algebraic topology, such as the fact that a covering map onto a simply connected space is a homeomorphism, we direct to \cite{Munkres}.

\begin{theorem}[Induced Lie group homomorphism]\label{thm_induced_homo}
\index{induced! -- Lie group homomorphism}
\index{simply connected}
 	Let $G $ and $H$ be Lie groups.
	Assume $G$ simply connected.
	For each Lie algebra homomorphism $\psi: \Lie(G)\to\Lie(H)$, there exists a unique Lie group homomorphism $\varphi: G\to H$ with $\varphi_*=\psi$.
\end{theorem}
\begin{proof}
 	Let $\Lie(G)=\g$ and $\Lie(H)=\h$.
	Since $\psi$ is a homomorphism, its graph 
	\[
	\mathfrak k:=\{(X, \psi(X)):X\in\g\}\subset \g\times\h
	\]
	is a subalgebra of $\g\times\h=\Lie(G\times H)$: for $X, Y\in \g$ we have
	\begin{equation*}
	 	[(X, \psi(X)),(Y, \psi(Y))]
		= \left([X, Y],[\psi(X), \psi(Y)]\right)
		= \left([X, Y], \psi[X, Y]\right).
	\end{equation*}
	By Theorem~\ref{teo1145}, there is a unique connected Lie subgroup $K\subset G\times H$ with $\Lie(K)=\mathfrak k$.
	Let $\pi_1: G\times H\to G$ and $\pi_2: G\times H\to H$ be the projections, which are Lie group homomorphisms.
	Let 
	\[
	\phi:= \pi_1|_K:K\to G .
	\]
	We have that
	\begin{equation}\label{eq1139}
	 	(\dd \phi)_{(1_G,1_H)}(X, \psi(X)) = X , \quad\forall X\in\g .
	\end{equation}
	In particular, we get that $(\dd\phi)_{(1_G,1_H)}: \mathfrak k\to \g$ is injective, and therefore it is an isomorphism (since $\dim \mathfrak k = \dim\g$).
	By Exercise~\ref{prop1302}, we deduce that $\phi:K\to G$ is a covering map.
	Since $G$ is simply connected, the map $\phi$ is a bijection; see \cite[Theorem~54.4]{Munkres}. Hence, it is an isomorphism.
	Set $\varphi:=\pi_2|_K\circ\phi^{-1}: G\to H$, which is a Lie group homomorphism.
	From \eqref{eq1139} we also get that
	\[
	(\dd\varphi)_{1_G}(X) = (\dd\pi_2)_{(1_G,1_H)}\circ(\dd\phi^{-1})_{1_G}(X) = (\dd \pi_2)_{(1_G,1_H)}(X, \psi(X)) = \psi(X),
	\]
	that is $\varphi_*=\psi$.
	
	Regarding the uniqueness, if $\tilde\varphi$ is another homomorphism such that $(\tilde\varphi)_*=\psi$, we get that $\tilde\varphi\circ\exp=\exp\circ\psi=\varphi\circ\exp$. Since $\exp$ is invertible in a neighborhood $U$ of the identity element, we have $\varphi|_U=\tilde\varphi|_U$. 
	Since such a $U$ generates $G$ and since $\varphi$ and $\tilde\varphi$ are group homomorphisms, we get that $\tilde\varphi=\varphi$.
\end{proof}

\subsection{The Baker-Campbell-Dynkin-Hausdorff formula}\label{section BCH}

 The punchline of this subsection is that one can recover near the identity element the group product of every Lie group in terms of the operations on the Lie algebra.
Various formulas have been found over the years. 
In particular, we recall chronologically the work by Campbell (in the years 1897-8), by Baker (in 1901-5), by Hausdorff (in 1906), and, finally, by Dynkin (in 1953).
The Baker-Campbell-Hausdorff formula links Lie groups to Lie algebras by expressing $ \log(e^X e^Y)$ as an infinite sum in the iterated Lie brackets of elements $X$ and $Y$ in the Lie algebra. 
The logarithm is, by definition, the inverse of the exponential; in general, it is only locally defined in a neighborhood of the identity, thanks to Proposition~\ref{exp:diffeo}. However, for nilpotent simply connected Lie groups, the logarithm will be global by Theorem~\ref{CG-1.2.1}. 

For matrix Lie groups, a Baker-Campbell-Hausdorff formula can be obtained using formal series as follows.
Namely, for every $A,B \in \gl (n)$ 
if $A$ and $ B$ are enough near $0$ so that the series converge, we write
\begin{eqnarray}
\nonumber
\log (e^Ae^B) 
&=&\log (\mathbb{I}+(e^Ae^B-\mathbb{I})) \\
\nonumber
&=& \sum _{k=1}^{\infty} \frac{(-1)^{k+1}}{k} (e^Ae^B-\mathbb{I})^{k}, \\
\nonumber
 & =& \sum _{k=1}^{\infty} \frac{(-1)^{k+1}}{k} \left(\left(\sum_{i=0}^\infty \frac {A^i}{i!} \right) \left(\sum_{j=0}^\infty \frac {B^j}{j!} \right)-\mathbb{I} \right)^{k} \\
\nonumber
& = &\sum _{k=1}^{\infty} \frac{(-1)^{k+1}}{k} \left( \sum_{\stackrel{i,j\in \N\cup\{0\} }{(i,j) \ne (0,0)}} \frac {A^i B^j}{i!j!} \right)^k \\
\label{serie_GLn}
&= &	 \sum_{k=1}^\infty \dfrac{(-1)^{k+1}}{k}\sum_
{\stackrel{i_\ell,j_\ell\in \N\cup\{0\} }{(i_\ell,j_\ell) \ne (0,0)}}
\dfrac{A^{i_1}B^{j_1}\cdots A^{i_k}B^{j_k}}{i_1!j_2!\cdots i_k!j_k!}.
\end{eqnarray}
We seek for a formula that expresses the above series of products 
as a series of iterated adjoint maps, like
$$\left(\ad_X^{ r_1} \circ\ad_Y^{ s_1}\circ 
\ldots \circ\ad_X^{ r_n}\right)(Y) $$ 
$$= [ \underbrace{X,[X, \ldots[X}_{r_1} ,[ \underbrace{Y,[Y, \ldots[Y}_{s_1} , \, \ldots\, [ \underbrace{X,[X, \ldots[X}_{r_n} , Y ]\ldots]]. 
$$

We shall state the 
 {\em Baker-Campbell-Hausdorff Formula} ({\em BCH Formula}, for short) using what is called the Dynkin product:
 \begin{definition}[Dynkin product]\label{Dynkin_product}For elements $A$ and $B$ in a Lie algebra, we define their {\em Dynkin product} as\index{Dynkin product}
$$
A\star B := \sum _{k,m \geq 0, \, i_\ell +j_\ell \ne 0} \frac{(-1)^{k}}{(k+1) (j_1+\dots +j_k +1)} \frac{ \ad_A^{i_1} \ad_B^{j_1} \dots \ad_A^{i_k} \ad_B^{j_k} \ad_A^{m} (B)}{ i_1!\dots i_k! j_1!\dots j_k! m!},
$$
if the series converges.
\end{definition}

\begin{proposition}[BCH Formula]\label{Proposition BCH}
For every $A, B \in \gl (n)$ with $\|A\|, \|B\| < \frac 1 2 \log (2-\sqrt 2 /2)$ we have\index{Baker-Campbell-Hausdorff Formula}\index{BCH Formula}
\begin{equation}\label{Dynkin Formula}
A\star B = \log (e^Ae^B).
\end{equation}
\end{proposition}
A clear proof can be found in \cite[Proposition~3.4.5]{Hilgert_Neeb:book}. Moreover, there are other equivalent ways of writing the Dynkin product:
\begin{equation}\label{Dynkin Formula2}
X\star Y = \sum_{n>0}\frac {(-1)^{n-1}}{n} \sum_{ \begin{smallmatrix} {r_i + s_i > 0} \\ {1\le i \le n} \end{smallmatrix}} \frac{ 
\left(\ad_X^{ r_1} \circ\ad_Y^{ s_1}\circ \ad_X^{ r_2} \circ\ad_Y^{ s_2} \ldots \circ\ad_X^{ r_n}\circ \ad_Y^{ s_n-1}\right)(Y) 
}{r_1!s_1!\cdots r_n!s_n! \; \sum_{i=1}^n (r_i+s_i)} .
\end{equation}

The first terms of the series are
 \begin{eqnarray*}
\log(\exp X\exp Y)&=& X + Y + \frac{1}{2}[X, Y] + \frac{1}{12}\Big([X,[X, Y]] + [Y,[Y, X]]\Big) \\ 
&&\quad - \frac {1}{24}[Y,[X,[X, Y]]] \\ 
&&\quad - \frac{1}{720}\Big([[[[X, Y],Y],Y],Y] +[[[[Y, X], X], X], X] \Big) \\ 
&&\quad +\frac{1}{360}\Big([[[[X, Y],Y],Y], X]+[[[[Y, X], X], X],Y] \Big)\\ 
&&\quad + \frac{1}{120}\Big([[[[Y, X],Y], X],Y] +[[[[X, Y], X],Y], X] \Big) + \cdots \end{eqnarray*}

Some consequences of the BCH formula are the following: 
\begin{description}
\item[\ref{Proposition BCH}.i] {\em First expansion of BCH}: 
\begin{equation}\label{expansion of BCH}
e^ X e^ Y = \exp\left( X + Y + \frac{1}{2}[X, Y] + o (\norm{X}\cdot \norm{Y}) \right), \text{ as } \norm{X}, \norm{Y}\to 0 . 
\end{equation}
\item[\ref{Proposition BCH}.ii] {\em Trotter product formula}:\index{Trotter product formula}
\begin{equation*}
\lim_{k\to \infty} (e^{\frac A k} e^{\frac B k})^k =e^{A +B}, \qquad \forall A, B \in \gl (n).
\end{equation*}
\item[\ref{Proposition BCH}.iii] {\em Commutator formula}:\index{commutator! -- formula}
\begin{equation*}
\lim_{k\to \infty} (e^{\frac A k} e^{\frac B k} e^{-\frac A k} e^{-\frac B k} )^{k^2} =e^{A B-BA}, \qquad \forall A, B \in \gl (n).
\end{equation*}
 \end{description}

\section{Exercises}


\begin{exercise}\label{ex:L:R} For every elements $g,h$ in a group $ G$ we have

\begin{minipage}[c]{0.4\linewidth}
	 \ref{ex:L:R}.i. 	$\quad\, L_h\circ L_g = L_{hg}$,
		\\ \ref{ex:L:R}.iii. 	$\quad L_h\circ R_g = R_g\circ L_h$,
			\\ \ref{ex:L:R}.iv. 	$\quad (R_g)^{-1} = R_{g^{-1}}$,
\end{minipage}
\begin{minipage}[c]{0.4\linewidth}
	\ref{ex:L:R}.ii. 	$\quad R_h\circ R_g = R_{gh}$,

	\ref{ex:L:R}.iv. 	$\quad (L_g)^{-1} = L_{g^{-1}}$,

	\ref{ex:L:R}.vi. $\quad C_{g h} = C_g \circ C_h$,
\end{minipage}
\\where $L$, $R$, and $C$ denote the left translations, right translations, and conjugations.
\end{exercise}

For the next two exercises, for a subset $U$ of a group and an integer $n\in \N$, set
\[
	U^n := \{g_1\cdot\dots\cdot g_n:g_1, \dots, g_n\in U\}.
	\]

\begin{exercise}\label{ex1922}
	Let $G$ be a Lie group (or, more generally, a topological group).
	If $U\subset G$ is open, then $U^2$ is open.
\end{exercise}

 \begin{exercise}\label{Prop:generating}
\index{subgroup! smallest --}
\index{smallest subgroup}
\index{generating! -- subgroup}
Connected groups are generated by neighborhoods of the identity:
 	Let $G$ be a connected Lie group (or, more generally, a topological group) and $U\subset G$ an open subset with $1\in U$. 
	Then $G=\bigcup_{n=0}^\infty U^n$.
	In other words, $G$ is the smallest group containing $U$.
\\{\it Solution.}
 	Let $U^{-1}:=\{g^{-1}:g\in U\}$ and $V:=U\cap U^{-1}$.
	Then $V$ is open, $V^{-1}=V$, $e\in V$.
	Let $H:=\bigcup_{n=1}^\infty V^n\subset \bigcup_{n=1} U^n$.
	Observe that $H$ contains $V$ and is a union of the open sets $V^n$ (see Exercise~\ref{ex1922}).
	Moreover, $H$ is closed under multiplication and inversion, since $V^n\cdot V^m \subset V^{n+m}$ and $V^{-n}\subset V^n$.
	In other words, $H$ is an open subgroup of $G$.
		Note that $gH$ is open for all $g\in G$, so $\bigcup_{g\notin H}gH$ is an open set.
		Since $G$ is connected, $G=H\sqcup \bigcup_{g\notin H}gH$ and $H\neq\emptyset$, we conclude that $G=H$.
\end{exercise}
 

\begin{exercise}\label{ex_open_subgroup}
\index{open subgroup}
\index{closed subgroup}
\index{subgroup! open --}
\index{subgroup! closed --}
Let $G$ be a Lie group. 
\\ \ref{ex_open_subgroup}.i. If $ H$ is a subgroup of $G$ that, topologically, is open, then it is closed;
\\ \ref{ex_open_subgroup}.ii. 
If $ H$ is a subgroup of $G$ that has a nonempty interior, then it is open and closed.
\\ \ref{ex_open_subgroup}.iii. 
Let $G^\circ$ be the {\em identity component} of $G$, which is the connected component of $G$ containing $1_G$.\index{identity! -- componen}\index{$G^\circ$}
We have that $G^\circ$ is an open, closed, normal subgroup of $G$. Its Lie algebra is the same as $G$;
\\ \ref{ex_open_subgroup}.iv. 
Every neighborhood $U\subseteq G$ of $G^\circ$
generates $G^\circ$, i.e., every element in $G^\circ$ is the product of finitely many elements in $U$.
\end{exercise}

\begin{exercise}
On topological groups, right translations and left translations are homeomorphisms. While in Lie groups, they are smooth diffeomorphisms.
\end{exercise}


\begin{exercise}\label{str_const_1}
The anti-commutativity and Jacobi identity for a Lie algebra rephrase on the fact that the structural constants $c^k_{ij}$ as in \eqref{eq1655} satisfy:

$ 	\begin{array}{ccclcc}
	 \ref{str_const_1}.{\rm i}. &0&=&
	 	c^k_{ij} + c^k_{ji} & , &
		\forall i,j,k\in\{1, \dots,n\} ;\\
		 \ref{str_const_1}.{\rm ii}. &0&=&
		\sum_{r=1}^n \left(c_{ij}^r c_{rk}^s + c_{jk}^r c_{ri}^s + c_{ki}^rc_{rj}^s \right) &, &
		\forall i,j,k,s\in\{1, \dots,n\}.
	\end{array}$
\end{exercise}
\begin{exercise}\label{str_const_2}
	Let $c^k_{ij}\in \R$ satisfying \ref{str_const_1}.i-ii.
	Define $[\cdot, \cdot]$ by \eqref{eq1655}.
	Then $[\cdot, \cdot]$ uniquely extends into a Lie bracket turning $\Span\{X_1, \dots, X_n\}$ into a Lie algebra.
\end{exercise}
\begin{exercise}\label{str_const_3}
 	Let $\g$ and $\tilde{\g}$ be Lie algebras of dimension $n$.
	Let $c^k_{ij}$ be the structural constants of $\g$ with respect to a basis $X_1, \dots, X_n$, and let 
	$\tilde c^k_{ij}$ be the structural constants of $\tilde{\g}$ with respect to a basis $\tilde X_1, \dots, \tilde X_n$.
	If $c^k_{ij}=\tilde c^k_{ij}$,
	 then the map $\psi: \g\to\tilde{\g}$ defined by $\psi(X_i):=\tilde X_i$ is a Lie algebra isomorphism.
\end{exercise}

 \begin{exercise}\label{LIVFs is closed under Lie bracket}
	The Lie bracket of two left-invariant vector fields is left invariant.
\\{\it Solution.} For left-invariant vector fields $X, Y$ on a Lie group $G$ and $g\in G$,
we have	
$
	(L_g)_*[X, Y] = [(L_g)_*X,(L_g)_*Y] = [X, Y].
	$
\end{exercise}

\begin{exercise}[Right translation of LIVF]
\index{right translation! -- of LIVF}
	 Let $X$ be a left-invariant vector field on a Lie group $G$. Let $R_g$ be the right translation by an element $g\in G$.
Then, the vector field $(R_g)_*X $ is left-invariant.
\\{\it Solution.} 
Let $h\in G$. Then, using Exercise~\ref{ex:L:R}.iii and that $X$ is left-invariant, we have
\begin{eqnarray*}
 	\dd L_h\circ( (R_g)_*X )
	&=& \dd L_h \circ \dd R_g \circ X \circ R_g^{-1}\\
	&=& \dd(L_h\circ R_g) \circ X \circ R_g^{-1}\\
	&=& \dd( R_g\circ L_h) \circ X \circ R_g^{-1} \\
	&=& \dd R_g\circ \dd L_h \circ X \circ R_g^{-1}\\
	&=& \dd R_g\circ X \circ L_h \circ R_g^{-1}\\
	&=& \dd R_g\circ X \circ R_g^{-1}\circ L_h\\
	&=& (R_g)_*X \circ L_h .
\end{eqnarray*}
\end{exercise}

\begin{exercise}[Derivative of product of curves]\label{ex:derivative:product: curves}
\index{derivative! -- of product of curves}
\index{product! derivative of -- of curves}
 Let $G$ be a Lie group. 
 Let $\gamma, \sigma: \R\to G$ be two smooth curves into $G$.
 Consider the product of the two curves, i.e., the curve $$t\longmapsto \gamma(t)\sigma(t)$$
 and calculate the derivative of such a curve in terms of 
 $\gamma$, $\sigma$, and their derivatives.
In fact, a formula is
\begin{equation}\label{derivative:product: curves}
 \frac{\dd}{\dd t} \gamma(t)\sigma(t)=(\dd R_{\sigma(t)})_{\gamma(t)}\dot\gamma(t) + (\dd L_{\gamma(t)})_{\sigma(t)}\dot\sigma(t) .
\end{equation}
\\
{\it Solution.} 
Derivating one variable at a time, we get
\begin{eqnarray*}
 \left.\frac{\dd}{\dd t}\gamma(t)\sigma(t)\right|_{t=t_0}
 &=& \left.\frac{\dd}{\dd t}\gamma(t)\sigma(t_0)\right|_{t=t_0}+ \left.\frac{\dd}{\dd t}\gamma(t_0)\sigma(t)\right|_{t=t_0}\\
 &=& \left.\frac{\dd}{\dd t}(R_{\sigma(t_0)}\gamma(t))\right|_{t=t_0} + \left. \frac{\dd}{\dd t}(L_{\gamma(t_0)}\sigma(t))\right|_{t=t_0} 
\\&=&(\dd R_{\sigma(t_0)})_{\gamma(t_0)}\dot\gamma(t_0) + (\dd L_{\gamma(t_0)})_{\sigma(t_0)}\dot\sigma(t_0) .
\end{eqnarray*}
 \end{exercise}

\begin{exercise}\label{ex:derivative:inverse: curves}\index{derivative! -- of the inverse of a curve}
 Let $G$ be a Lie group. 
 Let $\gamma: \R\to G$ be a smooth curve into $G$.
 Consider the curve $$t\longmapsto \gamma(t)^{-1}$$
 and calculate the derivative at an arbitrary $t$ of such a curve in terms of 
 $\gamma$ and $\dot\gamma$. In fact, a formula is
 \begin{eqnarray}
\frac{\dd}{\dd t}(\gamma(t)^{-1}) &=&
- (\dd L_{\gamma(t)^{-1}})_{1_G} (\dd R_{\gamma(t)^{-1}})_{\gamma(t)}\dot\gamma(t) .
\end{eqnarray}
\\
{\it Solution.} 
From the fact that $1= \gamma(t) \gamma(t)^{-1}$, for all $t$, and 
formula
\eqref{derivative:product: curves},
we have
$$ 0=(\dd R_{\gamma(t)^{-1}})_{\gamma(t)}\dot\gamma(t) + (\dd L_{\gamma(t)})_{\gamma(t)^{-1}} \frac{\dd}{\dd t}(\gamma(t)^{-1}).$$
Thus
\begin{eqnarray*}
\frac{\dd}{\dd t}(\gamma(t)^{-1}) &=& \nonumber
- \left((\dd L_{\gamma(t)})_{\gamma(t)^{-1}}\right)^{-1} (\dd R_{\gamma(t)^{-1}})_{\gamma(t)}\dot\gamma(t) 
\\
&=&
- (\dd L_{\gamma(t)^{-1}})_{1_G} (\dd R_{\gamma(t)^{-1}})_{\gamma(t)}\dot\gamma(t) .
\end{eqnarray*}
 \end{exercise}

\begin{exercise}\label{ex:156}
	Let $\varphi: G\to H$ be a group homomorphism.
	\\ \ref{ex:156}.i. We have	$\varphi\circ L_g = L_{\varphi(g)} \circ \varphi$, for all $g\in G$;
	\\ \ref{ex:156}.ii. We have 	$\varphi\circ R_g = R_{\varphi(g)} \circ \varphi$, for all $g\in G$.
\end{exercise}

\begin{exercise}
\label{phi_star_homo}\index{induced! -- Lie algebra homomorphism}
 	Let $\varphi: G\to H$ be a Lie group homomorphism. 
	Given a left-invariant vector field $X$ on $G$, let 
	 $\varphi_*X$ be the left-invariant vector field on $H$ for which $(\varphi_*X)_{1_H}= (\dd \varphi)_{1_G}(X_{1_G})$.
	\begin{description} 
	\item[\ref{phi_star_homo}.i] 	The vector fields $X$ and $\varphi_*X$ are $\varphi$-{\em related}, in the sense that\index{related vector fields}
		$
		(\dd \varphi)_gX_g = (\varphi_*X)_{\varphi(g)} ,
		$ for all $g\in G$.
			\item[\ref{phi_star_homo}.ii] 	If $g, g'\in G$ are such that $\varphi(g)=\varphi(g')$,
			 then $(\dd\varphi)_gX_g = (\dd\varphi)_{g'}X_{g'}$.
			\item[\ref{phi_star_homo}.iii] For all $g\in G$, we have $(\dd\varphi)_g(X_g) = (\dd L_{\varphi(g)})_{1_H} (\dd\varphi)_{1_G} X_{1_G}$. Hence, $\varphi_*X$ is the left-invariant extension of the (a-priori-not-well-defined) vector field on $H$ given as the push forward of $X$ via $\varphi$.

				\item[\ref{phi_star_homo}.iv] 	$\varphi_*: \Lie(G)\to\Lie(H)$ is a Lie algebra homomorphism.

	\item[\ref{phi_star_homo}.v] 	$(\dd\varphi)_{1_G} : (T_{1_G},[\cdot, \cdot]) \to (T_{1_H}H,[\cdot, \cdot])$ is a Lie algebra homomorphism.
	\end{description}

{\it Hint:} From Exercise~\ref{ex:156}.(i), we have 
 	\begin{align*}
	(\varphi_*X)_{\varphi(g)} 
	&= (\dd L_{\varphi(g)})_{1_H} (\dd\varphi)_{1_G} X_{1_G} \\
	&= (\dd (L_{\varphi(g)}\circ\varphi))_{1_G}X_{1_G} \\
	&= (\dd(\varphi\circ L_g))_{1_G}X_{1_G} \\
	&= (\dd\varphi)_g(\dd L_g)_{1_G}X_{1_G} \\
	&= (\dd\varphi)_g X_g.
	\end{align*}
 For $X, Y\in\Lie(G)$, on the one hand $[X, Y]\in\Lie(G)$, on the other hand $[X, Y]$ and $[\varphi_*X, \varphi_*Y]$ are $\varphi$-related.
	Thus $\varphi_*[X, Y] = [\varphi_*X, \varphi_*Y]$. 
\end{exercise}


\begin{exercise}\label{ex2853}
	Let $X^\dag$ be a right-invariant vector field on $G$. A curve $\theta:\R\to G$ is a one-parameter subgroup with $\dot\theta(0)=X^\dag_{1_G}$ if and only if $\theta(t)=\Phi^t_{X^\dag}({1_G})$, for all $t\in\R$.
\end{exercise}

\begin{exercise}\label{ex4June1041}\index{one-parameter subgroup}
For a vector $X$ in the Lie algebra of a Lie group $G$, consider the map
$\psi: \Lie(\R)\to \Lie(G)$, 
$t \mapsto t X$.
\\ \ref{ex4June1041}.i. The map $\psi$ is a Lie algebra homomorphism.
\\ \ref{ex4June1041}.ii. There exists a one-parameter subgroup 
$\gamma: \R\to G$ with $\d\gamma=\psi$. 
\\ \ref{ex4June1041}.iii. We have $\dot\gamma(t) = X_{\gamma(t)}$.
\\{\it Solution.} 
Since $\R$ is simply connected, Theorem~\ref{thm_induced_homo} asserts that there exists such a 
$\gamma$. 
Regarding (iii), we have
\begin{eqnarray*}
 \dot\gamma(t) &=& \left.\dfrac{\dd}{\dd h}\gamma(t+h) \right|_{h=0}\\
&=& \left.\dfrac{\dd}{\dd h}\gamma(t)\gamma(h) \right|_{h=0}\\
&=& \left.\dfrac{\dd}{\dd h}L_{\gamma(t)}(\gamma(h)) \right|_{h=0}\\
&=& (L_{\gamma(t)})_*(\dd\gamma)_0(\partial_t)\\
&=& (L_{\gamma(t)})_*\psi(\partial_t)\\
&=& (L_{\gamma(t)})_*X\\
&=& X_{\gamma(t)}.
\end{eqnarray*}
\end{exercise}

\begin{exercise}\label{exp: OPS}\index{exponential! -- map of a Lie group}
\index{one-parameter subgroup}
 	Let $G$ be a Lie group and $X$ be a left-invariant vector field on $G$.
\\\ref{exp: OPS}.i.	$\Phi^1_{tX} = \Phi^t_X$, for all $t\in\R$;
\\\ref{exp: OPS}.ii.	$\exp(tX) = \Phi^t_X(1_G)$, for all $t\in\R$;
\\\ref{exp: OPS}.iii.	$\exp(sX+tX) = \exp(sX)\exp(tX)$, for all $s,t\in\R$;
\\\ref{exp: OPS}.iv.	$\exp(0)=1_G$;
\\\ref{exp: OPS}.v. $\exp(- X)=(\exp(X))^{-1};$
\\\ref{exp: OPS}.vi. 	 $t\mapsto\exp(t X)$ is a one-parameter subgroup and an integral curve of $X$;
\\\ref{exp: OPS}.vii. Equations \eqref{flow_of_left} and \eqref{flow_of_right} holds true.
\end{exercise}
%

\begin{exercise}\label{ex:exp:mx}
	 Let $G$ be a Lie group. For all $m\in \Z$ and $X\in \Lie(G)$ we have
 $$\exp(m X)=(\exp(X))^m.$$
 Pay attention that your proof is correct also when $m$ is negative.
\end{exercise}

\begin{exercise}
Let $X$ be a left-invariant vector field in a Lie group $G$. 
Then, we have
$$X f (\exp(tX)) = \dfrac{\rm d}{{\rm d} t} f(\exp(tX)) , \qquad \forall f\in C^\infty(G), \forall  t\in \R.$$
\end{exercise}

\begin{exercise}[Lie group structures on tangent bundles of Lie groups]\label{ex_TG_group}\index{Lie group! -- structures on tangent bundles of Lie groups}
Let $G$ be a Lie group. 
\\
\ref{ex_TG_group}.i.
On $TG$, consider the following operation:
$$
(g, X_g) * (h, Y_h) = d(R_h)_g(X_g) + d(L_g)_h(Y_h) \in T_{gh}G,
$$
for all $g, h \in G$, $ X_g \in T_gG $, and $ Y_h \in T_hG$.
Then, the pair $(TG, *)$ is a Lie group.
\\
\ref{ex_TG_group}.ii.
On $ \Lie(G)\times G$ consider the following operation:
$$(X, g) \bullet (Y, h) =
(X +\dd(C_g)_{1_G}Y, gh),$$
for all $g, h \in G$, $ X, Y \in \Lie(G)$.
Then, the pair $( \Lie(G)\times G, \bullet)$ is a Lie group.
\\
\ref{ex_TG_group}.ii. The two Lie groups $(TG, *)$ and $( \Lie(G)\times G , \bullet)$ are isomorphic.
\end{exercise}

\begin{exercise}\label{ex:gamma:k:k}
 Let $G$ be a Lie group. 
 Let $\gamma: \R\to G$ be a smooth curve into $G$ with $\gamma(0)=1_G$ and $\dot\gamma(0)=X$.
Then, we have
 $$ \lim_{k\to \infty} (\gamma(t/k))^k = \exp (tX), \qquad \forall t\in \R.$$
\\
{\it Solution.} 
For \( t \) small enough, one can consider \( \eta(t) := \exp^{-1}(\gamma(t)) \).
For fixed \( t \in \mathbb{R} \), we have
\begin{align*}  
t \dot\eta(0) 
&= t \lim_{h \to 0} \frac{\eta(h) - \eta(0)}{h}
=
t \lim_{k \to \infty} \frac{\eta(t/k) - \eta(0)}{t/k}
\\
&=
\lim_{k \to \infty} k \eta(t/k).
\end{align*}
Because of Proposition~\ref{exp:diffeo}, we write
\begin{equation*}  
\dot\eta(0) = \left.
\frac{\dd}{\dd t} \exp^{-1}(\gamma(t))\right|_{t=0}
= (\dd \exp^{-1})_0 \dot\gamma(0)
=(\dd \exp)_0^{-1} \dot\gamma(0) = \dot\gamma(0).
\end{equation*}
Thus, using the previous two formulas, we get
\begin{align*}  
\exp(tX) 
&= \exp(t \dot\gamma(0)) 
= \exp(t \dot\eta(0)) 
= \exp(\lim_{k \to \infty} k \eta(t/k)) 
\\
&= \lim_{k \to \infty} \exp(k \eta(t/k)) 
= \lim_{k \to \infty} (\exp(\eta(t/k))^k
= \lim_{k \to \infty} (\gamma(t/k))^k.
\end{align*} 
\end{exercise}

\begin{exercise}\label{ex: LIVF: commutes: RIVF}
	If $X$ is a LIVF and $Y$ is a RIVF, then $[X, Y]=0$.
\end{exercise}
\begin{exercise}\label{ex: commutative}
If $G$ is a commutative Lie group, then Lie$(G)$ is a commutative Lie algebra.
\\{\it Hint:} Use of Exercise~\ref{ex: LIVF: commutes: RIVF}.
\end{exercise}

\begin{exercise}[What happens if one uses right-invariant vector fields as Lie algebra]\index{Lie algebra! -- of RIVF}\label{LieRIVF}
 	For $X, Y\in T_{1_G}G$. Let $\tilde X, \tilde Y$ be the left-invariant vector fields such that $\tilde X_{1_G}=X$ and $\tilde Y_{1_G}=Y$.
	Let $X^\dag$ and $Y^\dag$ be the right-invariant vector fields with $(X^\dag)_{1_G}=X$ and $(Y^\dag)_{1_G}=Y$.
	\begin{enumerate}[label=(\roman*).]
	\item 	We have $[X^\dag,Y^\dag]_{1_G} = - [\tilde X, \tilde Y]_{1_G}$.
	\item 	Setting $[X, Y]_R := [X^\dag,Y^\dag]_{1_G}$, the two Lie algebras $\g=(T_{1_G}G,[\cdot, \cdot])$ and $(T_{1_G}G,[\cdot, \cdot]_R)$ are isomorphic Lie algebras via the map $X\mapsto -X$.
	\end{enumerate}
{\it Solution.}
 	Consider the map $J: G\to G^\dag$, $J(g)=g^{-1}$ from the group $G=(G, \cdot)$ to $G^\dag=(G,*)$, where 
	\[
	g*h := h\cdot g, \qquad \forall g,h\in G.
	\]
	Notice that $G^\dag$ is a Lie group.
	Observe that $J$ is a Lie group isomorphism: 
	\[
	J(g\cdot h) = (g\cdot h)^{-1} = h^{-1}\cdot g^{-1} = g^{-1}*h^{-1} = J(g)*J(h) .
	\]
	
	We claim that 
	\begin{equation}\label{eq0445}
	J_*\tilde X = -X^\dag, 
	\qquad\forall X\in T_{1_G}G .
	\end{equation}
	Indeed, using Corollary~\ref{prop:right_left}, for all $g\in G$ we have
	\begin{eqnarray*}
	 	(\dd J)_g\tilde X_g 
		&=&\left. \frac{\dd}{\dd t} J(g\exp(tX)) \right|_{t=0}\\
		&=& \left. \frac{\dd}{\dd t} \exp(-tX)\cdot g^{-1} \right|_{t=0} \\
		&=& \left. \frac{\dd}{\dd t} R_{g^{-1}}(\exp(-tX)) \right|_{t=0} \\
		&=& (\dd R_{g^{-1}})_{1_G}(-X)\\
		&=& -(X^\dag)_{g^{-1}}\\
		&=& -(X^\dag)_{J(g)} .
	\end{eqnarray*}
	This proves \eqref{eq0445}.
	
	The proof is thus complete, because $J_*:(T_{1_G}G,[\cdot, \cdot])\to(T_{1_G}G,[\cdot, \cdot]_R)$ is the Lie algebra isomorphism we were looking for.
\end{exercise}
 
\begin{exercise}\label{exp_smooth}\index{exponential! smoothness of the -- map}
Given a Lie group $G$ with Lie algebra $\g$, consider the vector field $Y$ on the manifold $G\times\g$ defined as follows.
For all $(g, X)\in\G\times\g$,
	\[
	Y_{(g, X)} := (X_g,0) \in T_gG\times T_X\g \simeq T_{(g, X)}(G\times\g).
	\]
		We have
	\[
	\Phi^t_Y((g, X)) = ( \Phi^t_X(g), X ) , \qquad \forall t\in \R , \forall g\in G, \forall X\in\g.
	\]
Thus, the map $X\mapsto \exp(X)$ is smooth because it is the projection of the flow at time $1$ of the smooth vector field $Y$.
\\
{\it Solution.}
	We have
	$
	 	\frac{\dd}{\dd t}(\Phi^t_X(g), X) = (X_{\Phi^t_X(g)},0 ) = Y_{(\Phi^t_X(g), X)}
	$
	and $\left.\left(\Phi^t_X(g), X\right)\right|_{t=0} = (g, X)$.
	We deduce that $\Phi^t_X(g)$, which is the first coordinate of the above flow, depends smoothly on the point $(g, X)$ and so does $\exp(X):=\Phi^1_X(1_G)$.
	\end{exercise}

\begin{exercise}
Let $\varphi_1, \varphi_2: G\to H$ be two Lie group homomorphisms such that the associated Lie algebra homomorphisms
$(\varphi_1)_*, (\varphi_2)_*$ coincide. Assume that $G$ is connected. We have that $\varphi_1=\varphi_2$. But, there are counterexamples when $G$ is not connected. 
\end{exercise}
\begin{exercise}\label{ex4June20241020}
	 Proposition~\ref{prop: expF: Fexp} implies that every Lie group homomorphism $F$ has the following properties: 
 \\ \ref{ex4June20241020}.i.	$F_*$ is injective (resp.~surjective) if and only if $F$ is locally injective (resp.~open) at $1_G$.
 \\ \ref{ex4June20241020}.ii.	Given $g\in G$, $F_*$ is injective (resp.~surjective) if and only if $F$ is locally injective (resp.~open) at $g$.	
 \\ \ref{ex4June20241020}.iii.	$F_*$ is injective (resp.~surjective) if and only if $F$ is locally injective (resp.~open).
 \\ \ref{ex4June20241020}.iv.	If $F$ is bijective, then $F^{-1}$ is smooth, hence $F$ is a diffeomorphism.
 \end{exercise}

\begin{exercise}\label{Ex:inv:homo:cont}
	Every bijective Lie group homomorphism has a continuous inverse. 
	\\{\it Hint:} Use Proposition~\ref{prop: expF: Fexp}.
\end{exercise}
\begin{exercise}\label{Ex:homo:immersion}
	Every injective Lie group homomorphism is an immersion (i.e., the differential is injective).\index{immersion}
	\\{\it Hint:} Use Proposition~\ref{prop: expF: Fexp}.
\end{exercise}

\begin{exercise}\label{ex: def: Lie: subgroup: cont}
In the definition of Lie subgroup, one can replace the requirement that the inclusion is an immersion by requiring that it is continuous.
Namely, a subgroup $H < G$ of a Lie group $G$
	 is a {Lie subgroup}
	of $G$ if $H$ admits the structure of a Lie group such that the inclusion $H\into G$ is a continuous homomorphism.
\end{exercise}

\begin{exercise}[Square root of a matrix]\label{no:square:root}\index{square root of a matrix}
 Let $G$ be a Lie group with Lie algebra $\mathfrak{g}$ and identity component $G^\circ$.

(i) We have $\exp(\mathfrak{g})\subset G^\circ$.

(ii) For each $A\in \exp(\mathfrak{g})$ there exists $B\in G$ such that 
$B^2=A$. (Every such a $B$ is called {\em a square root} of $A$).
\\
{\it Hint:} If $A=\exp(X) $ take $B:=\exp(\frac12X).$
\end{exercise}
\begin{exercise}
\label{prop1302}
 	Let $G, H$ be connected Lie groups, and $\varphi: G\to H$ a Lie group homomorphism. 
	The following are equivalent:
	\index{discrete kernel}
	\index{covering map}
	\index{local! -- diffeomorphism}
	\\ \ref{prop1302}.i. the map 	$\varphi$ is surjective and has a discrete kernel;
		\\ \ref{prop1302}.ii. 	the map $\varphi$ is a covering map;
		\\ \ref{prop1302}.iii. the map 	$\varphi_*$ is an isomorphism of Lie algebras;
		\\ \ref{prop1302}.iv. the map 	$\varphi$ is a local diffeomorphism.
\\{\it Hint.} The proof can be found in \cite[page~182, Proposizione~3.8.2]{Abate_Tovena_geometria_differenziale}.
\end{exercise}
 
 \begin{exercise}[Non-surjective exponential]\label{ex:exp:no:sujective}
 \index{non-surjective exponential}
 \index{exponential! non-surjective --}
 \index{non-Riemannian exponential}
 \index{exponential! non-Riemannian --}
 	\index{general linear group, $\GL$}
	\index{$\GL^+(n, \R)$}
	\index{group! general linear --, $\GL$}
Let $\GL^+(n, \R)$ be the subgroup of $\GL(n, \R)$ consisting of the matrices with positive determinant.

(i) The Lie group $\GL^+(n, \R)$ is open and connected, and it is the identity component of $\GL(n, \R)$.
%
%
%

(ii)
 Exercise~\ref{no:square:root} implies that for some $n\in \N$, the map $\exp : \mathfrak{gl}(n, \R) \to \GL^+(n, \R)$ is not surjective. 
\\
{\it Hint:} Try 
$
\begin{bmatrix}
 	-1 & 0 \\
	0 & -2 
\end{bmatrix}$ or 
$
\begin{bmatrix}
 	-1 & -1 \\
	0 & -1 
\end{bmatrix}\in\GL^+(2, \R) $.
 \end{exercise}
 


\begin{exercise}\label{ex_ad_homo}
The linearity of the Lie bracket and the Jacobi Identity give \eqref{ad_homo}.i and \eqref{ad_homo}.ii.
\\ {\it Solution of }\eqref{ad_homo}.ii. 
For $X, Y,Z\in\g$, we have
$\ad([X, Y])(Z) 
		= [[X, Y],Z]\\
		\overset{\rm Jacobi}= [X,[Y,Z]] - [Y,[X,Z]] 
		= \ad_X(\ad_Y(Z)) - \ad_Y(\ad_X(Z))
		= [\ad_X, \ad_Y](Z).$
\end{exercise}

\begin{exercise}\label{ex_Ad_homo}
 Exercise~\ref{ex:L:R}.vi implies \eqref{ex_Ad_homo2}. 
\\ {\it Solution}. 		Let $g,h\in G$ and differentiate at $1_G$ the identity $C_g\circ C_h = C_{gh}$, to get
	$
	\Ad(gh) 
	= (\dd C_{gh})_1
	= (\dd C_g)_1\circ(\dd C_h)_1
	= \Ad(g)\circ\Ad(h) .
	$
\end{exercise}

\begin{exercise}
For all $X, Y$ in the Lie algebra of a Lie group, we have
$$\exp({X})\exp(Y)\exp({-X} ) = \exp(e^{\operatorname{ad} _X} Y) .$$
{\it Solution.} 
Using Formula~\ref{Formula:Cg:Ad} first and then Formula~\ref{Formula: Ad: ad}, we have
\(
\exp(X)\exp(Y)\exp(-X) = \exp(\Ad_{\exp(X)}Y) = \exp(e^{\ad_X}Y).
\)
\end{exercise}

\begin{exercise}\label{Ex:ad:BAB}
	For every vector space $V$,   $A\in\gl(V)$, and   $B\in\GL(V)$, we have
	$
	\Ad_B(A) = B A B^{-1}.
	$
\\{\it Solution.} 
 	\begin{multline*}
	 	\Ad_B(A) 
		\stackrel{\text{def}}{=} (\dd C_B)_\mathbb{I}A
		= \left. \frac{\dd}{\dd t} C_B ( e^{tA} ) \right|_{t=0} \\
		= \left. \frac{\dd}{\dd t} B e^{tA}B^{-1} \right|_{t=0} 
		= \left. \frac{\dd}{\dd t} e^{t BAB^{-1}} \right|_{t=0}
		= BAB^{-1} .
	\end{multline*}
\end{exercise}

\begin{exercise}\label{Formulae_GL2}
For all $X, Y\in \mathfrak{gl}(V)$, we have $e^{\operatorname{ad} _X} Y = e^{X}Y e^{-X} .$
\\
{\it Solution.} 
Using Formula~\ref{Formula: Ad: ad} and then Exercise~\ref{Ex:ad:BAB}, we have $e^{\ad_X}Y = \Ad_{e^X}Y = e^X Y e^{-X} .$
\end{exercise}

\begin{exercise}\label{Formulae_GL3}
	For all $A\in\gl(V)$ and for all $B\in\GL(V)$, we have 
	$e^{BAB^{-1}} = Be^AB^{-1} $.
\\		{\it Hint:} Notice that $(BAB^{-1})^k = BA^k B^{-1}$, for $k\in\N$. Expand $e^{BAB^{-1}}$ in power series.
\end{exercise}

\begin{exercise}\label{ex_fomula_flow_ad}
	Let $X, Y$ left-invariant vector fields on a Lie group $G$.
For all $t\in\R$
\[
(\Phi^t_X)_*Y = e^{-\ad(tX)}Y.
\]
{\it Solution.} 
 	\begin{align*}
	 	(\Phi^t_X)_*(Y) 
		&= (R_{\exp(tX)})_*Y \\
		&= (R_{\exp(tX)})_*(L_{\exp(-tX)})_*Y \\
		&= (R_{\exp(tX)} \circ L_{\exp(-tX)})_*Y \\
		&= (C_{\exp(-tX)})_* Y \\
		&= \Ad_{\exp(-tX)} Y 
		= e^{\ad(-tX)}Y.
	\end{align*}
\end{exercise}

\begin{exercise} If $\gamma$ is a curve into a Lie group, then 
$$\frac{\dd}{\dd s} \Ad_{\gamma(s)} = \Ad_{\gamma(s)} \ad\big( (\dd L_{\gamma(s)}^{-1})_{\gamma(s)} ( {\frac{\dd}{\dd s} \gamma(s )})\big) $$
\\
{\it Solution.} 
Use twice that $\Ad_p \circ \Ad_q = \Ad_{pq}$ to obtain
\begin{eqnarray}
\nonumber
 \partial_s \Ad_{\gamma(s)}
&=& \partial_\eps \Ad_{\gamma(s+\eps)}|_{\eps=0} \\
\nonumber
&=& \partial_\eps\Ad_{\gamma(s)} \Ad_{\gamma(s)^{-1}} \Ad_{\gamma(s+\eps)}|_{\eps=0} \\
\nonumber
&=& \Ad_{\gamma(s)} \partial_\eps\Ad_{\gamma(s)^{-1} \gamma(s+\eps)}|_{\eps=0} \\
\nonumber
&=& \Ad_{\gamma(s)} \ad({\partial_\eps (\gamma(s)^{-1} \gamma(s+\eps))}|_{\eps=0}) \\
\nonumber
&=& \Ad_{\gamma(s)} \ad\big( (\dd L_{\gamma(s)}^{-1})_{\gamma(s)} ( {\partial_s \gamma(s )})\big).
\end{eqnarray}
 \end{exercise}

\begin{exercise}\label{rmk:exp:converge}
 	For all $A\in\Mat_{n\times n}(\R)$, entry by entry, the matrix exponential $e^A$ is an absolutely convergent series.
	\\
{\it Solution.} For each $M\in\Mat_{n\times n}(\R)$ set
	$
	\| M\| := \sup\{|Mv|:|v|\le 1\},
	$
	where $|\cdot|$ is the Euclidean norm in $\R^n$.
	Then 
	\[
	\left\| \sum_{k=N_1}^{N_2} \frac{1}{k!} A^k \right\|
	\le \sum_{k=N_1}^{N_2} \frac{1}{k!} \|A^k\|
	\le \sum_{k=N_1}^{N_2} \frac{1}{k!} \|A\|^k \ \overset{N_1,N_2\to\infty}\longrightarrow\ 0.
	\]
\end{exercise}

\begin{exercise}\label{ex: abba}
Let $A,B\in\gl(n\R)$.
If $AB=BA$, then $e^{A+B}=e^Ae^B = e^Be^A$.
\\
{\it Solution.} 
 	\begin{eqnarray*}
	 	 e^A \cdot e^B 
		&=& \left(\sum_{k=0}^\infty \frac1{k!}A^k \right) \cdot\left(\sum_{l=0}^\infty \frac1{l!}B^l \right)\\
		&=& \sum_{k,l} \frac1{k!}\frac1{l!} A^kB^l \\
		&=& \sum_{m=0}^\infty \sum_{j=0}^m \frac1{j!(m-j)!} A^jB^{m-j} \\
		&=& \sum_{m=0}^\infty \frac1{m!} \sum_{j=0}^m \binom mj A^j B^{m-j}
		\overset{(AB=BA)}= \sum_{m=0}^\infty \frac1{m!} (A+B)^m\\
		&=& e^{A+B}.
	\end{eqnarray*}
\end{exercise}

\begin{exercise}\label{ex_eA_invertible}
For every matrix $A$, the matrix $e^A$ is invertible. 
\\{\it Solution:} Use Exercise~\ref{ex: abba} and get $e^{A}e^{-A} = e^0 = I$.
\end{exercise}

\begin{exercise}
Calculate the exponential of the matrix
$
t \begin{bmatrix} 0 & a & c\\ 0 & 0 & b\\ 0 & 0 & 0\\ \end{bmatrix}.
$
\end{exercise}

\begin{exercise}Let $A,B\in\gl(n\R)$.
	Find $A,B$ such that $e^{A+B}\neq e^Ae^B \neq e^Be^A$.
Compare with Exercise~\ref{ex: abba}.
\end{exercise}

\begin{exercise}
	Given a matrix $A$ and an invertible matrix $B$, we have $e^{BAB^{-1}}=Be^AB^{-1}.$
\end{exercise}

\begin{exercise}\index{trace}
	The determinant function $\det: \GL(n, \R)\to (\R^*, \cdot)$ is a Lie group homomorphism,  the trace function ${\rm{tr}} : \mathfrak{gl}(n, \R)\to (\R,+)$ is a Lie algebra homomorphism, and  
	\[
	\det(e^A) = e^{\rm{tr}(A)}.
	\]
{\it Solution.} 
 	Given a matrix $A$, there is an invertible matrix $B$ such that $\tilde A = BAB^{-1}$ is upper triangular, i.e., of the form
	\[
	\tilde A=
	\begin{bmatrix}
	 	\alpha_1 & * & * & * \\
		0 & \alpha_2 & * & * \\
		0 & 0 & \ddots & \vdots \\
		0 & 0 &\dots & \alpha_n 
	\end{bmatrix} .
	\]
	For such matrices, we have
	\[
	\tilde A^k =
	\begin{bmatrix}
	 	\alpha_1^k & * & * & * \\
		0 & \alpha_2^k & * & * \\
		0 & 0 & \ddots & \vdots \\
		0 & 0 &\dots & \alpha_n^k 
	\end{bmatrix} \quad
	\text{ 
	and 
	}\quad
	e^{\tilde A} = 
	\begin{bmatrix}
	 	e^{\alpha_1} & * & * & * \\
		0 & e^{\alpha_2} & * & * \\
		0 & 0 & \ddots & \vdots \\
		0 & 0 &\dots & e^{\alpha_n} 
	\end{bmatrix} .
	\]
	Finally, using Formula~\ref{Formulae_GL}.iii we conclude
	\begin{eqnarray*}
	 	\det(e^A) &=& \det(Be^{A}B^{-1}) \\
		&=& \det(e^{BAB^{-1}}) \\
		&=& \det(e^{\tilde A}) \\ 
		&=&e^{\alpha_1}\cdots e^{\alpha_n} \\
		&=& e^{\sum_{i=1}^n\alpha_i} \\
		&=& e^{\rm{tr}(\tilde A)} \\
		&=& e^{\rm{tr}(BAB^{-1})} = e^{\rm{tr}(A)}.
	\end{eqnarray*}
\end{exercise}

\begin{exercise}\index{differential! -- of exponential map}
For all $X, Y\in \mathfrak{gl}(n, \R)$ the derivative of $e^X$ in the direction $Y$ has the formula:
$$\lim_{t\to 0 }\dfrac{e^{X+tY} - e^X}{t}
= \sum_{k=1}^\infty \dfrac{1}{k!} \sum_{i=1}^k X^{i-1} Y X^{k-i}.$$
\end{exercise}


\begin{exercise}\label{ex_deriv_0} Deduce
\ref{Prop_Lie_der}.i from Exercise~\ref{phi_star_homo}, the chain rule, and Theorem~\ref{thm_induced_homo0}.
\end{exercise}

\begin{exercise}\label{ex_deriv_1}
The space ${\rm {Der}}(\g)$ is a Lie subalgebra of $\gl(\g)$. 
\\
{\it Solution}. First ${\rm {Der}}(\g)$ is a linear subspace of $\gl(\g)$, because equation \eqref{eq_derivation} is linear in $D$. Second, it is easy to verify that if $D, D'$ are derivations, then $[D, D']:= D\circ D' - D'\circ D$ is a derivation: for all $X, Y\in \g$ we have
\begin{eqnarray*}
[D, D'] ([X, Y])&=&( D\circ D') ([X, Y]) - (D'\circ D)( [X, Y])\\
&=& D([D'X, Y] + [X,D'Y]) - D' ([DX, Y] + [X,DY])\\
&=& [DD'X, Y] +[D'(X),DY] + [DX,D'(Y)] + [X,DD'Y] \\
&&\hspace{2cm}- \left( [ D' DX, Y] +[ D(X),D'Y] + [D'X,DY]) + [X,D'DY]) \right)\\
&=& [(DD' -D' D )X, Y] + [X,(DD'-D'D)Y] \\
&=& [[D, D'] X, Y] + [X,[D, D'] Y] .
\end{eqnarray*}
\end{exercise}

\begin{exercise}\label{ex_deriv_2} Let $\g$ be a Lie algebra.
The map ${\ad}_X$ is a derivation on $\g$ (because of Jacobi identity) and 
 the map $X\mapsto {\ad}_X$ is a Lie algebra homomorphism of $\g$ into ${\rm {Der}}(\g)$.
\end{exercise}

 \begin{exercise}\label{ex_conj_auto}\index{conjugation}
 Let $G$ be a group, on which we denote by $C_g$ the conjugation by an element $g\in G$. The map $g\mapsto C_g$ is a group action of $G$ by automorphisms of $G$.
\end{exercise}

\begin{exercise}\label{ex_deriv_3} Let $G$ be a Lie group.
For all $g\in G$ we have ${\Ad}_g\in{\rm {Aut_{Lie}}}(\g)$ and the map $g\mapsto {\Ad}_g$ is a Lie group homomorphism of $G$ into ${\rm {Aut_{Lie}}}(\g)$.
\end{exercise}

\begin{exercise}\label{ex_deriv_4}
The space $
{\rm {Aut_{Lie}}}(\g) $ is a closed Lie subgroup of $\GL(\g)$ whose Lie algebra is ${\rm {Der}}(\g)$.
\end{exercise}

\begin{exercise}\label{ex_deriv_5} The previous exercises imply Proposition~\ref{Prop_Lie_der}.iii.
\end{exercise}



\chapter{Metric groups and homogeneous spaces}\label{ch_MetricGroups}

In differential geometry, the term {\em homogeneous space} refers to the quotient space of a Lie group modulo a closed subgroup, which results in a manifold with a smooth transitive action of the Lie group.
However, in metric geometry and, more generally, in analysis on metric spaces, the term homogeneous refers to functions that, when precomposed with dilations of the space, are multiplied by some constants;
see Exercise~\ref{ex:homog:dist}.


At this point, we are forced to distinguish between the terms.
If a metric space admits a transitive action by isometries, we refer to it as an isometrically homogeneous space; see Section~\ref{Isometrically homogeneous spaces}. If, instead, a space admits a self-map that non-trivially dilates the distance function, we refer to it as a self-similar space; see Section~\ref{sec: self-similar}.

In this chapter, we will see that every isometrically homogeneous space with mild topology has the structure of a Lie homogeneous space (in the standard sense of Section~\ref{sec: Lie: homog}). In Section~\ref{sec: self-similar2}, we will see that if, in addition, the space is self-similar, then it has the structure of a metric group in the sense of Section~\ref{sec: metric_group}.

\section{Lie homogeneous spaces}\label{sec: Lie: homog}
In this section, we review group actions and their quotients.
\subsection{The general viewpoint of group actions}\label{sec_actions}
The group of isometries of a metric space naturally acts on the space itself, and from this action, it inherits a topological structure. It is convenient to begin our discussions with the general notion of group action.
Group actions on sets were mentioned in \eqref{eq:def_left_action}. When, however, the sets and the groups have extra structures, it is natural to take these into account for the considered actions.

A {\em continuous action} of a topological group $G$ on a topological space $X$, generically denoted by $G\acts X$, is a continuous map $ G \times X \to X$, usually denoted by $(g,x)\mapsto g.x$, that obeys the {\em associativity law} and the {\em identity law}:\index{continuous action}\index{action! continuous --}\index{associativity law}\index{identity law}
$$ (gh).x = g.(h.x)\qquad \text {and}\qquad 1_G.x=x, \quad\forall g,h\in G, \forall x\in X.$$

\begin{definition}[Special types of actions]
Fix an action \( G \acts X \) of a group on a set.
The action is {\em transitive} if for every $x,y\in X$ there is $g\in G$ such that $g.x=y$.\index{transitive! -- action}\index{action! transitive --} It is {\em faithful} if for every distinct elements $g,h\in G$ there is $x\in X$ such that $g.x\neq h.x$. Faithful actions are also known as {\em effective} actions.\index{faithful action}\index{action! faithful --}\index{effective action|see {faithful action}}\index{action! effective --|see {faithful}}
The action is \emph{free} if \( g.x \neq x \), for all \( x \in X \) and all \( g \in G\setminus\{1_G\} \).\index{free! -- action}\index{action! free --}
Next, assume that the action
 \( G \times X \to X \) is a continuous action of a topological group on a topological space.
The action is \emph{proper} if the associated map\index{proper action}\index{action! proper --}
\[ G \times X \to X \times X 
\qquad
(g,p) \mapsto ( g.p, p)
\] 
is a proper map (i.e., inverse images of compact sets are compact); see also Exercise~\ref{Prop 21.5 from Lee} for alternative definitions. 
The action is \emph{properly discontinuous} if for all \( x,y \in X \) there exist a neighbourhood \( U_x \) of \( x \) and a neighbourhood \( U_y \) of \( y \) such that the set \( \{g \in G \mid (g.U_x) \cap U_y \neq \emptyset\} \) is finite.\index{properly discontinuous action}\index{action! properly discontinuous --}
\end{definition}

Given an action \( G \acts X \) and an element $x\in X$, the {\em orbit} of $x$ is defined as\index{orbit}
$$G.x: =\{g.x\,|\, g\in G\}.$$
Either two orbits $G.x$ and $G.x'$ coincide, or they are disjoint sets. Therefore, we have a well defined {\em quotient space} denoted by $\mfaktor{G}{X}$:\index{$\mfaktor{G}{X}$}\index{quotient! -- space} 
$$\mfaktor{G}{X}: = \{G.x \,|\,x\in X \}. $$
We may denote the map $x\in X \mapsto G.x\in\mfaktor{G}{X}$ by $\pi$ and call it the 
{\em quotient map}.\index{quotient! -- map}
While, fixed a point $x_0\in X$, we call the map 
$g\in G \mapsto g.x_0\in{X}$ the
{\em orbit map}.\index{orbit! -- map}

Given an action \( G \acts X \) and an element $x\in X$, the {\em stabilizer subgroup} of $G$ at $x$, also called the {\em isotropy subgroup} at $x$,
is the subgroup 
$G_x: =\{g \in G : g.x=x\}$.\index{isotropy subgroup}\index{stabilizer} 

\begin{example}\label{ex_subgroup_acting}
There is a natural action of each group onto itself by left translations; see Exercise~\ref{ex:action_left_translations}.
More generally, every subgroup $H$ of a group $G$ acts on the whole group. 
Actually, there are two possible natural actions, as we now review.
The map
\begin{eqnarray*}
H \times G &\to& G\\
(h,g)&\mapsto & hg
\end{eqnarray*}
is an action of $H$ on $G$, called the {\em action by left translations} of $H$ on $G$.\index{action! -- by left translations} Whereas, the map
\begin{eqnarray*}
H \times G &\to& G\\
(h,g)&\mapsto & g h^{-1}
\end{eqnarray*}
is an action of $H$ on $G$, called the {\em action by right translations} of $H$ on $G$.\index{action! -- by right translations} We observe that when $H\acts G$ by right translations, then 
$\mfaktor{H}{G} = G/H: =\left\lbrace gH \mid g\in G\right\rbrace.$
Moreover, the map $(g, g'H)\mapsto g g' H$ defines an action of $G$ on $G/H$.
\end{example}

If $G\acts X$ is a continuous action of a topological group on a topological space, then there is a unique topology, which we call the {\em quotient topology}, on 
$\mfaktor{G}{X}$ that makes the quotient map continuous and open; see Exercise~\ref{ex_orbit_topo}.

The following basic result gives a natural homeomorphism between orbits and groups modulo stabilizers. The general argument goes back to \cite{Arens}. See also \cite[Lemma~5.38]{Drutu-Kapovich} for another version.
In the following theorem, we shall begin to see the importance of dealing with spaces and groups that are locally compact (and are second countable; see Exercise~\ref{second countable VS sigma-compact}).

\begin{theorem}\label{Helgason_Theorem_3.2}\index{transitive! -- action}\index{stabilizer}\index{locally compact}
Let $X$ be a locally compact Hausdorff space. 
Let $G$ be a locally compact group with a countable basis. 
Let $G\acts X$ be a continuous action. 
Let $x\in X$. 
Then 
the stabilizer subgroup $G_x$ at $x$ is closed.
Moreover, if the action is transitive, then the map
$g G_x \mapsto g \cdot x$
is a homeomorphism between $G/ G_x$ and $X$.
\end{theorem}
\begin{proof}
Being $X$ a Hausdorff space, the set $\{x\}$ is closed.
 Since the orbit map $ g\mapsto g.x $ is continuous and $G_x$ is the preimage of $x$ under such a map, we deduce that 
$G_x $ is closed in $G$. 
Next, we additionally assume that the action is transitive, so the orbit map is surjective. 
We then consider the following commutative diagram of surjective maps:
\begin{center}
 \begin{tikzcd}
 G \ar[d,->>=.5] \ar[drr,->>=.5]{r}{g\mapsto g.x} && \\
G/ G_x \arrow[rr,->>=.5]{d}{gG_x\mapsto g.x\quad} &&X 
 \end{tikzcd}
\end{center}
Notice that the bottom map is a bijection.
Recalling that the quotient map $ G \to G/ G_x$ is open and continuous, we stress that it is enough to prove that the orbit map is open. 
Let $V$ be an open subset of $G$ and $g$ a point in $V$.
 Since $G$ is assumed locally compact, select a compact neighborhood $U$ of $1_G$ in $G$
such that $U =U^{-1}$ and $gU^2 \subseteq V$. 
Since $G$ is assumed to have a countable basis, 
there exists a sequence $(g_n)_{n\in\N} \subseteq G$ such that 
$G = \bigcup_{n\in\N} g_nU$. The group action being transitive implies 
$X = \bigcup_{n\in\N} g_nU x $.
 Each summand is compact; hence, it is a closed subset of $X$.
Recall that the Baire Category Theorem holds for locally compact Hausdorff spaces; see \cite[Theorem~11.7.3, p.394]{MR2723563} (or look at the general argument in \cite[p.161]{Folland_book}).\index{Baire Category Theorem} 
By such a theorem, some summand, and therefore $U .x$, contains an
inner point $u.x$ with $u\in U$. Then $x$ is an inner point of $u^{-1}U. x \subseteq U^2 . x$ and
consequently $g.x$ is an inner point of $V . x$. This shows that the orbit map is open.
\end{proof}

\label{Lie group actions}
\index{Lie group! -- action|see {Lie action}}
\index{action! Lie group --|see {Lie action}}
\index{group! Lie -- action|see {Lie action}}
For Lie groups acting on manifolds, actions that take into account the differential structures of the group and the manifold are considered.
\begin{definition}[Lie action]
Let \( G \) be a Lie group and \( M \) a smooth manifold. A \emph{Lie action} of \( G \) on \( M \) is an action \( G \times M \to M \) 
that is a smooth map.\index{Lie! -- action}\index{action! Lie --}

\end{definition}
 There is a natural Lie action of every Lie group onto itself by left translations; see Exercise~\ref{ex:action_left_translations}.
 Likewise, every Lie subgroup $H$ of a Lie group $G$ gives Lie actions \( H \acts G \), as in Example~\ref{ex_subgroup_acting}.

When a Lie action is proper and free, the quotient space is naturally a manifold; see Theorem~\ref{thm_quotient}. We shall not prove this theorem in this generality; we shall just refer to \cite[Theorem~21.10]{Lee_MR2954043}. We shall prove the respective result in the case of actions of closed subgroups of Lie groups; see Theorem~\ref{thm_diff_structure_quotient}.

\begin{theorem}[Quotient Manifold Theorem]\label{thm_quotient}\index{Theorem! Quotient Manifold --}
Let \( G \acts M \) be a proper and free Lie action of a Lie group $G$ on a manifold $M$. Then, there exists a unique differentiable structure on \( \mfaktor{G}{M} \) such that \( M \to \mfaktor{G}{M} \) is a smooth map with surjective differential.
\end{theorem}

 
\subsection{Lie coset spaces}\index{Lie! -- coset space}\index{quotient! -- space} 
Let $G$ be a Lie group and $H$ a closed subgroup of $G$.
We will consider manifold structures for the space $$G/H: =\left\lbrace gH \mid g\in G\right\rbrace.$$
When the set $G/H$ is equipped with the differentiable structure from the following theorem, it is called a \emph{homogeneous manifold}, or {\em Lie homogeneous manifold}, or {\em Lie coset space}.\index{Lie! -- homogeneous manifold} 

\begin{theorem}\label{thm_diff_structure_quotient}\index{manifold! -- structure}
 Let $G$ be a Lie group and $H<G$ a closed subgroup.
Then, the topological space $G/H$ admits the structure of a differentiable manifold.
Moreover, there is a unique manifold structure for which the action 
$G\times G/H \to G/H$ is smooth.
\end{theorem}

\begin{proof}[Sketch of the proof.]
Complete proofs are in \cite[Theorem 3.58]{Warner} and in \cite[Theorem~4.2]{Helgason}.
 We begin by stressing that $G/H$ is a second countable Hausdorff space; see Exercises~\ref{G/H Hausdorff} and \ref{G/H second countable}.
 We shall find an atlas for $G/H$ by constructing {\em local cross sections} in $G$.\index{local! -- cross section}


 Recalling that $H$ is a Lie group by Theorem~\ref{teo1131}, we write $\mathfrak g$ and $\mathfrak h$ for the Lie algebras of $G$ and $H$, respectively.
 Let $\mathfrak m$ be some vector subspace of $\mathfrak g$ that is in direct sum with $\mathfrak h$, that is, $\mathfrak g = \mathfrak m\oplus\mathfrak h$.
 Because of Proposition~\ref{exp:diffeo}, there are open neighbourhoods $U_{\mathfrak m}\subset \mathfrak m$ of $0_{\mathfrak m}$ and $U_{\mathfrak h}\subset\mathfrak h$ of $0_{\mathfrak h}$
 such that the map
 $\phi:\mathfrak m\oplus\mathfrak h\to G$, $\phi(x,y): =\exp(x)\exp(y)$,
 is a diffeomorphism between $U_{\mathfrak m}\times U_{\mathfrak h}$ and its image $U_G\subset G$,
 and such that $H\cap U_G = \exp(U_{\mathfrak h})$ 
 (the latter claim uses the hypothesis that $H$ is closed).

 For each $g\in G$, define the map
 $\psi_g:U_{\mathfrak m}\to G/H$, $\psi_g(x): = \pi(g\exp(x)) $, where $\pi(g): =gH$.
 One shows that, up to shrinking $U_{\mathfrak m}$, the map $\psi_e:U_{\mathfrak m}\to G/H$ is injective onto $\pi(U_G)$,
 and so $\psi_g$ is a homeomorphism $U_{\mathfrak m}\to\pi(gU_G)$ for all $g\in G$.
 We obtained the candidate atlas:
 \begin{equation}\label{eq66338472}
 \{(\pi(gU_G) , \psi_g^{-1})\}_{g\in G}
 \end{equation}
 To check that 
 each composition $\psi_{g_2}\circ \psi_{g_1}^{-1}$ is smooth on the appropriate domain,
 we locally write such a map as the composition 
 $\pi_{\mathfrak m}\circ\phi^{-1}\circ R_h\circ L_{g_2^{-1}g_1} \circ \phi$,
 where $\pi_{\mathfrak m}:\mathfrak m\oplus\mathfrak h\to\mathfrak m$ is the projection modulo $\mathfrak h$, and $h\in H$ is a suitable element.


 We stress that, with this differentiable structure~\eqref{eq66338472}, each map $\psi_g$ is a smooth section of the projection $\pi$.
%
 Consequently, the action $G\times G/H \to G/H$ is smooth.
 Also the uniqueness is an easy consequence.
 For these last details, we refer to \cite[page 122]{Warner}.
\end{proof}

\section{Isometrically homogeneous spaces}\label{Isometrically homogeneous spaces}
\begin{definition}\label{def isometrically homogeneous}
We say that a metric space $M$ is {\em isometrically homogeneous} if its group of isometries acts on the space transitively. Explicitly, this means that, for every $p,q\in M$, there exists a distance-preserving homeomorphism $f: M\to M$ such that $f(p)=q$.
An {\em isometrically homogeneous space} is a metric space that is isometrically homogeneous.\index{isometrically homogeneous}
\end{definition}

Examples of isometrically homogeneous spaces are groups equipped with left-invariant distances, which we call {\em metric groups}. This whole book is devoted to studying various types of metric groups.\index{metric! -- group}

From the general viewpoint, we shall study isometrically homogeneous spaces with mild topological assumptions. One of them is local compactness. Other topological assumptions are connectedness and local connectedness; see \cite{Munkres}.
 
 If \( M \) is a metric space, we consider its isometry group \( \Isom(M) \), that is, the set of self-isometries of \( M \) equipped with the composition rule and the pointwise-convergence topology.\index{isometry! -- group}
We stress that in the isometry group, the compact-open topology and the pointwise-convergence topology coincide, and the topological group \( \Isom(M) \) acts continuously on \( M \); see Remark~\ref{prop:3-topols}.

The main aim of this section is to show that isometry groups of isometrically homogeneous spaces (with mild topological assumptions) are Lie groups; see Theorem~\ref{Montgomery-Zippin}. 
Such a fact is a consequence of the solution of Hilbert's fifth
problem by Montgomery--Zippin, and Gleason, together with the observation that isometry groups are second-countable and locally compact. This latter property follows the Ascoli--Arzel\`a Theorem (Exercise~\ref{ex_AA}).
Regarding the axioms of countability, see Exercise~\ref{second countable VS sigma-compact}.

\subsection{Transitive actions by locally compact groups and Hilbert fifth problem}
\index{transitive! -- action}\index{locally compact}\index{Hilbert fifth problem}

We next state one of the results that are coming from Gleason-Montgomery-Zippin's theory, as it is stated in Montgomery-Zippin's book \cite[Corollary on page 243, Section 6.3]{Montgomery_Zippin52}:
{\it If a (separable) locally compact group $G$ satisfying Property A} (which we soon recall) {\it acts effectively and transitively on a locally compact, connected, (locally connected), and finite-dimensional space, then $G$ is a Lie group.} We added the words in parenthesis because we think the authors were implicitly assuming them. Moreover, in the language of Montgomery-Zippin \cite[Section 6.2]{Montgomery_Zippin52}, {\em Property A} means that for every neighborhood $V$ of the identity element in $G$, there exists a compact subgroup $K$ of $G$ such that $K\subset V$ and $G/K$, equipped with the quotient topology, is a Lie group.\index{Property A}\index{Montgomery-Zippin Theorem}\index{Theorem! Montgomery-Zippin --}

A clean argument for the above statement is not present in the literature. 
Since we want to present a proof of the Lie structure of locally compact groups transitively acting on `nice' spaces, we shall then use a stronger and more established version of the solution of the Hilbert fifth problem, attributed to Gleason and Yamabe.\index{Gleason-Yamabe Theorem}\index{Theorem! Gleason-Yamabe --}

\begin{theorem}[Gleason-Yamabe]\label{Thm_GY}
Let $G$ be a 
first countable, locally compact group. Then, there exists an open subgroup $G'<G$ that is the inverse limit of a sequence of Lie groups.
\end{theorem}

When using Theorem~\ref{Thm_GY} in the proof of Theorem~\ref{Thm_GMYZ}, we shall rephrase the concept of an inverse limit. We shall use the notion of inverse limit of topological groups and of topological spaces; as a reference, see \cite[Definition~4.2.5]{Tao_book}.\index{inverse limit}
We shall not prove Theorem~\ref{Thm_GY} but suggest to look at the book by Tao \cite[Chapter~6]{Tao_book}.\index{Tao}
We stress that in \cite[paragraph below Theorem 6.0.11]{Tao_book}, there is the requirement of $G$ being Hausdorff, but, for us, this is part of the definition of a topological group.
Before drawing consequences of interest to isometrically homogeneous spaces in Theorem~\ref{Montgomery-Zippin}, we first give proof of the following general result on continuous, faithful, and transitive actions. 

\begin{theorem}[after Gleason-Montgomery-Yamabe-Zippin]\label{Thm_GMYZ}\index{Theorem! Gleason-Montgomery-Yamabe-Zippin --}
Let $G$ be a topological group that is locally compact and second countable. 
Let $X$ be a Hausdorff topological space that is connected, locally connected, locally compact, and finite-dimensional.
Assume that $G$ acts on $X$ continuously, faithfully, and transitively. 
Then $G$ is a Lie group. 
\end{theorem}

 
\begin{proof}
A proof in the case $X$ is a manifold (and $G$ is $\sigma$-compact, and not assumed second countable) can be found in \cite[Section~6.4]{Tao_book}.
We present an argument along the same lines.
Pick a point $x_0\in X$ and let $H: =G_{x_0}$ be the stabilizer of $G$ at $x_0$.
Observe that $H<G$ is a closed subgroup, that $G/H$ is a Hausdorff topological space, and that the map 
$\phi: G/H \to X$, $\phi(gH): =gx_0$, is a homeomorphism (this is true because $G$ is locally compact and second countable and $X$ is locally compact and Hausdorff, and the action is transitive; see Theorem~\ref{Helgason_Theorem_3.2}.)


We apply the Gleason-Yamabe result, Theorem~\ref{Thm_GY}, recalling Exercise~\ref{second countable VS sigma-compact}.i.
Thus, there is a subgroup $G'<G$ and Lie groups $G_n$, for $n\in \N$, 
such that $G'$ is open in $G$ and $G'=\varprojlim G_n$.
The fact that $G'$ is the inverse limit of the sequence $G_n$ in the category of topological groups rephrases as saying that there exist compact normal subgroups $K_n\lhd G'$ such that $G_n=G'/K_n$ is a Lie group, for each $n\in \N$, and $K_n \searrow \{1\}$ as $n\to \infty$ (i.e., $K_n\supset K_{n+1}$ and $\bigcap_{n\in\N} K_n = \{1\}$).
By Exercise~\ref{open_component_trans}, since $X$ is connected and $G'$ open, then $G'$ acts on $X$ transitively.

Set $H'$ to be the stabilizer of $G'$ at the initial point $x_0$ 
and $H_n: =H'/K_n< G'/K_n=G_n$. 
We have two inverse limits: one in the category of topological groups and one in the category of topological spaces. They are expressed by the following diagrams.
\begin{center}
 \begin{tikzcd}
G' \ar[dr, two heads ] \ar[drr, two heads ]\ar[drrr, two heads] \ar[drrrrr, two heads] & & & & &\\
\cdots\arrow[r, two heads]& G_{n+1} \arrow[r, two heads] &G_n \arrow[r, two heads]&G_{n-1}\arrow[r, two heads]&\cdots\arrow[r, two heads]&G_1
 \end{tikzcd}
\end{center}
\begin{center}
 \begin{tikzcd}
G'/H'\simeq X \ar[dr, two heads ] \ar[drr, two heads ]\ar[drrr, two heads] \ar[drrrrr, two heads] & & & & &\\
\cdots\arrow[r, two heads]& G_{n+1} /H_{n+1}\arrow[r, two heads] &G_n /H_n\arrow[r, two heads]&G_{n-1}/H_{n-1}\arrow[r, two heads]&\cdots\arrow[r, two heads]&G_1/H_1
 \end{tikzcd}
\end{center}
 Notice that each $H_n$ is a closed subgroup of the Lie group $G_n$. Hence, the quotient $G_n/H_n$ is a smooth manifold.
 Moreover, the space $X$ is the inverse limit of $G_n/H_n$ in the category of topological spaces. In fact, the compact group $K_n$ acts continuously on $X$ with quotient map $$\pi_n:X \to X_n: = G_n/H_n. 
 $$
 For every $m\ge n$, the compact Lie group $K_n/K_m$ smoothly acts on $X_m$
with quotient map $\pi_{m,n}:X_m\to X_n$.

 Since, by assumption, the dimension of $X$ is finite, then the increasing sequence of dimensions $\dim(G_n/H_n)$ must stabilize.
 Note that the structure group of the quotient
 \begin{equation}\label{projection_GY}
 \begin{tikzcd}
 G_{n+1} /H_{n+1}\arrow[r, two heads] &G_n/ H_n
 \end{tikzcd}
\end{equation}
is $K_n/K_{n+1}$. 
By Exercise~\ref{ex_dim_quot}, since the dimensions stabilize, then eventually $K_n/K_{n+1}$ is zero-dimensional and compact, and hence it is finite. Hence, this projection map is a covering map between manifolds, whose covering multiplicity is the cardinality $\# K_n/K_{n+1}$.

We claim that the projections in \eqref{projection_GY} eventually stabilize as homeomorphisms, which is equivalent to saying that eventually $ K_n/K_{n+1}$ is trivial. To show this, we shall use that $G'/H'\simeq X$ is assumed locally connected.

%
%

Let $\bar n$ be such that, for every $m\ge n\ge \bar n$, the group $K_n/K_m$ is finite
and thus the projections $\pi_{m,n}:X_m\to X_n$ are covering maps. 
Fix a simply connected neighbourhood $U_{\bar n} $ of $\pi_{\bar n}(x_0)$ in $X_{\bar n}$.
For $m\ge\bar n$, define $U_m : = \pi_{m,\bar n}^{-1}(U_{\bar n})$.
 Notice that $U_m$ is a disjoint union of connected components, each homeomorphic to $U_{\bar n}$ via $\pi_{m,\bar n}$.
Since $X$ is locally connected,
there is a connected neighbourhood $V\subset \pi_{\bar n}^{-1}(U_{\bar n})$ of $x_0$ in $X$.
Choose $\hat n\ge \bar n$ such that
$K_{\hat n}\subset\phi^{-1}(V)\subset G'$,
so that $K_{\hat n}x_0\subset V$.

We then check that $K_{\hat n}=K_m$ for every $m\ge\hat n$.
To this aim, notice that the orbit $\Omega: =(K_{\hat n}/K_m)(\pi_m(x_0))$
is a finite set with cardinality $\#(K_{\hat n}/K_m)$.
On the one hand, every two points in $\Omega$ belong to different connected components of $U_m$.
On the other hand,
since $\Omega = \pi_m( K_{\hat n}x_0 ) \subset \pi_m(V) $
and since $\pi_m(V)$ is a connected subset of $U_m$,
the set $\Omega$ is contained in one connected component of $U_m$.
We conclude that $\#(K_{\hat n}/K_m) = 1$, i.e., $K_{\hat n}=K_m$.

Therefore, the quotient $ K_n/K_{n+1}$ is eventually trivial, which implies that the inverse limit $\varprojlim G_n$ stabilizes and $G'$ is its stabilization. Thus, we deduce that $G'$ is a Lie group, and $G$ is a Lie group since $G'$ is open in $G$.
\end{proof}


\subsection{Properties of isometrically homogeneous spaces}\index{isometrically homogeneous}\index{transitive! -- action}\index{locally compact}

\subsubsection{Two preliminary observations}
We begin with a first topological property of locally compact isometrically homogeneous spaces: completeness. 
\begin{lemma}\label{completeness_homog}\index{complete}
Every locally compact isometrically homogeneous space is complete.
\end{lemma}	
\begin{proof}
Pick a point $\bar x$ in a locally compact isometrically homogeneous space $M$. By local compactness, there is $\bar r>0$ such that the closed ball $\bar B(\bar x,\bar r)$ is compact.
Now, take a Cauchy sequence $(x_n)_{n\in \N}$ in $M$. For $N\in\N$ large enough we have $d(x_n, x_m)<\bar r$, for all $n,m\geq N$. Thus the sequence $(x_n)_{n >N}$ is in the set $\bar B( x_N,\bar r)$. This latter set is compact by isometric homogeneity. Thus $(x_n)_{n >N}$ has a convergent subsequence.
\end{proof}

The isometry group is equipped with the compact-open topology - review this notion in \cite[p.285]{Munkres}.\index{compact-open topology}
It is a closed subset of the space of homeomorphisms; see Exercise~\ref{limit of isometries}.
We shall see that it forms a locally compact group. However, this is not immediately clear from Ascoli--Arzel\'a theorem because it may happen that 
a
locally compact isometrically homogeneous space $M $ 
is not boundedly compact: even if small balls are compact, some large balls may not be compact. 
For example, this is the case for the distance $\min\{d_E, 1\}$ on $\R$, where $d_E$ denotes the Euclidean distance.
We shall solve this issue using the following result.

\begin{lemma}[Reduction to proper distances]\label{rmk:boundedy_cpt}\index{boundedly compact}\index{reduction to proper distance}
If \( (M,d) \) is a connected locally compact isometrically homogeneous space, then there exists a finite-valued distance \( \rho \) inducing the same topology such that \( (M,\rho) \) is an isometrically homogeneous space that is boundedly compact and \( \Isom(M,d) \) is a closed subgroup of \( \Isom(M,\rho) \).
\end{lemma}
\begin{proof}
Fix a point $o\in M$. Since the metric is locally compact, then there exists some $r_0>0$ such that the closed ball
$\bar B_d(o, r_0)$ is compact. Notice that, by isometric homogeneity, every other $r_0$-ball $\bar B_d(p, r_0)$, with $p\in M$, is compact.
Then, for each $p, q\in M$, we consider the value
$$\rho(p,q): =
\inf
\left\{ 
\sum_{i=1}^k d(p_{i-1}, p_i) \,:\, k\in \N,\, p_i\in M,\, p_0=p, \,p_k=q,\, d(p_{i-1}, p_i)\leq r_0/4
 \right\}.$$
First, we observe that this is a distance function that gives the same topology, since for all $r\in (0,r_0/4)$, we have that $\rho(p,q)<r$ if and only if $d(p,q)<r$.
 Consequently, this distance must be finitely valued since, otherwise, the connected components of points at finite distance would disconnect the space, which is assumed connected.
 Second, because the construction is done intrinsically, we have 
 $\Isom(M,d) \subseteq \Isom(M,\rho).$
Regarding its closeness, recall Exercise~\ref{limit of isometries}.
 Finally, we claim that for all $r>0$, the set $\bar B_\rho(o, r)$ is compact.
We prove by induction on $n\in \N$ that $\bar B_\rho(o, nr_0/8)$ is compact.
The base of the induction is $n=0$, and it is clear since the ball reduces to a point.
Assume that $\bar B_\rho(o, nr_0/8)$ is compact. Hence, it is a totally bounded set: there exists a finite set $Y_n\subseteq \bar B_\rho(o, nr_0/8)$ such that 
\begin{equation} \label{Y n almost covers}\bar B_\rho(o, nr_0/8)\subseteq \bigcup_{y\in Y_n} B_d(y, r_0/2).\end{equation} 
We claim that $$ \bar B_\rho(o, (n+1)r_0/8)\subseteq \bigcup_{y\in Y_n} \bar B_d(y, r_0),$$
from which we will deduce that $\bar B_\rho(o, (n+1)r_0/8)$ is compact.
To prove the claim, take $p\in \bar B_\rho(o, (n+1)r_0/8)$. By the definition of $\rho$, there exists $\bar p \in \bar B_\rho(o, nr_0/8)$ such that $d(\bar p, p)<r_0/2$ (in fact, one can take points $p_0=p,p_1,\ldots, p_k=o$ such that $\sum_{i=1}^k d(p_{i-1}, p_i) < (n+2)r_0/8$ and $d(p_{i-1}, p_i)\leq r_0/4$; then a possible $\bar p$ is the first point among $p_0 ,p_1,\ldots, p_k$ such that $d(\bar p, o )\leq nr_0/8$).
By \eqref{Y n almost covers}, there is some $y\in Y_n$ such that $d(\bar p, y)<r_0/2$.
Then, by triangle inequality $d(y, p)<r_0 $. Thus, the claim is proved.
\end{proof}

\subsubsection{Transitive isometry groups}\index{transitive! -- isometry group}\index{isometry! -- group}
We consider groups that act on metric spaces transitively and by isometries. 
In particular, we will clarify why the isometry group of a connected locally compact isometrically homogeneous space is locally compact.

\begin{remark}\label{prop:3-topols}
We stress that for every metric space $M$, the compact-open topology, the topology of uniform convergence on compact sets, and the topology of pointwise convergence agree on the isometry group \( \Isom(M) \) of $M$. Moreover, the group \( \Isom(M, d) \) becomes a topological group continuously acting on $M$. For a reference to the fact that these topologies agree on \( \Isom(M, d) \); see \cite[Lemmas~5.B.1 and 5.B.2]{Cornulier-Harpe} or \cite[p.232, Theorem 15]{Kelley_gen_top}.
The fact that this structure makes the isometry group a topological group is well known; van Dantzig and van der Waerden \cite{Dantzig-Waerden} showed this in the case where \( M \) is connected, locally compact, and separable, and a general proof for metric spaces can be found in \cite[Lemma~5.B.3]{Cornulier-Harpe}. The fact that the action is continuous is an immediate consequence of the definition of compact-open topology. 
\end{remark}

A well-known consequence of Ascoli--Arzel\'a argument, see \cite[Lemma~5.B.4]{Cornulier-Harpe}, is that if a metric space is boundedly compact, then the isometry group is locally compact; see also Exercise~\ref{Cornulier-Harpe5.B.4}. We record this property for further reference.\index{Ascoli--Arzel\'a Theorem}

\begin{proposition}\label{Iso_gp_loc_cpt}\index{stabilizer}
Let $M$ be a boundedly compact metric space. Then $\Isom(M)$ equipped with the equivalent topologies, 
compact-open or pointwise convergence, is a second-countable, $\sigma$-compact, locally compact group acting on $M$ continuously and properly. 
For every $o\in M$ the stabilizer ${\rm Stab}(o)$ of $\Isom(M)$ at $o$ 
is compact.
If, in addition, the space $M$ is isometrically homogeneous, 
 the orbit map $\phi:\Isom(M)\to M$, $\phi(f): =f(o)$, induces a homeomorphism between 
topological spaces
$$\faktor{\Isom(M)}{{\rm Stab}(o)}\, \longrightarrow \,\quad M.$$ 
\end{proposition}
\begin{proof}
Most of this standard fact is an exercise and can be mostly found in \cite[Proposition~5.B.5]{Cornulier-Harpe}; see Remark~\ref{prop:3-topols}. However, the result is just a consequence of the Ascoli-Arzel\`a argument (Exercise~\ref{ex_AA}). The homeomorphism is the one discussed in Theorem~\ref{Helgason_Theorem_3.2}.
\end{proof}

\begin{remark}\label{rmk:boundedy_cpt2}
In the case $M$ is a connected locally compact isometrically homogeneous space, then the conclusion of Proposition~\ref{Iso_gp_loc_cpt} is still valid. 
Indeed, if $d$ denotes the distance function of $M$, then by Lemma~\ref{rmk:boundedy_cpt} 
there is a proper distance $\rho$, i.e., $(M, \rho)$ is boundedly compact, for which \( \Isom(M,d) \) is a closed subgroup of \( \Isom(M,\rho) \). Clearly, also the stabilizer of \( \Isom(M,d) \) at a point $o\in M$ is a closed subgroup of the stabilizer at $o$ of \( \Isom(M,\rho) \).
Because we got closed subgroups, we get the same conclusions about \( \Isom(M,d) \) and its action from Proposition~\ref{Iso_gp_loc_cpt}.
 \end{remark}
 
%
 

As in the proof of Theorem~\ref{Thm_GMYZ}, it is essential to know when identity components of transitive groups still act transitively. This is the case for isometry groups of connected and locally compact spaces. The following result will be used a few times for the characterization of isometrically homogeneous spaces in Section~\ref{sec: self-similar}.

\begin{proposition}\label{identity_component_trans}
Suppose that $M$ is an isometrically homogenous space and that $M$ is connected and locally compact. Then, every open subgroup of $\Isom(M)$ acts transitively on $M$. If, in addition, the topological group $\Isom(M)$ has the structure of a Lie group, then the identity component $\Isom(M)^\circ$ of $\Isom(M)$ acts transitively on $M$.
\end{proposition}
\begin{proof}
 We shall apply Theorem~\ref{Helgason_Theorem_3.2}. Being a metric space, the space $M$ is Hausdorff. Regarding $\Isom(M)$, we showed that it is a locally compact group and has a countable basis; see Remark~\ref{prop:3-topols} and Exercise~\ref{Cornulier-Harpe5.B.4}. 

Let $H$ be an open subgroup of $\Isom(M)$.
First, we claim that for every $q\in M$, the orbit $H\cdot q$ of $q$ under $H$ is open. Indeed, this is because, called 
$G: = \Isom(M)$ and $G_p$ the stabilizer subgroup of $G$ at $p$,
 the projection $G\to G/G_p$ is open and the orbit action $G/G_p\to M$ is a homeomorphism by Theorem~\ref{Helgason_Theorem_3.2}. 
Next, fix a point $p\in M$, suppose by contradiction that $H\cdot p\neq M$. Hence,
$$ M=\left(H\cdot p\right)\bigsqcup\left( \bigcup_{q\notin H\cdot p}H\cdot q\right)$$
is a disjoint union of two non-empty open
sets of $M$. This contradicts the fact that $M$ is connected. Therefore, the subgroup $H$ acts transitively.
For the last part of the statement, recall that connected components of manifolds are open subsets. 
\end{proof}

We are ready to explain the Lie group structure of isometry groups of isometrically homogeneous spaces with mild topological assumptions.
\begin{theorem}[after Gleason-Montgomery-Yamabe-Zippin]\label{Montgomery-Zippin}\index{Theorem! Gleason-Montgomery-Yamabe-Zippin --}\index{manifold! -- structure}\index{Lie group}\index{isometry! -- group}
Let $M$ be a metric space that is connected, locally connected, locally compact, and has finite topological dimension.
Assume that the isometry group ${\rm Isom}(M)$ of $M$ acts transitively on $M$.
Then ${\rm Isom}(M)$ has the structure of a Lie group with finitely many
connected components, and $M$ has the structure of a smooth manifold for which the action ${\rm Isom}(M)\acts M$ is smooth.
\end{theorem}

\begin{proof}
We shall apply Theorem~\ref{Thm_GMYZ}. For another proof using the earlier works of Gleason, Montgomery, and Zippin,
\cite{Montgomery_Zippin52, Gleason}; see \cite[Chapter 16]{Drutu-Kapovich}. 
By Proposition~\ref{Iso_gp_loc_cpt}, together with Remark~\ref{rmk:boundedy_cpt2},
the topological group $\Isom(M)$ is locally compact and second countable, and it is acting continuously.
 Obviously, the group $\Isom(M)$ acts effectively. It acts transitively by assumption.
 Theorem~\ref{Thm_GMYZ} implies that $\Isom(M)$ is a Lie group, with finitely many components; see Exercise~\ref{ex finitely many components}. 
 
By Proposition~\ref{Iso_gp_loc_cpt}, and Remark~\ref{rmk:boundedy_cpt2}, each 
stabilizer ${\rm Stab}(o)$, with $o\in M$, is compact and hence a closed Lie subgroup.
Consequently, by Theorem~\ref{Helgason_Theorem_3.2} the metric space $M$ is homeomorphic to 
 ${\Isom(M)}/{{\rm Stab}(o)}$, which is a manifold on which ${\rm Isom}(M) $ acts smoothly, by Theorem~\ref{thm_diff_structure_quotient}.
\end{proof}



\section{Metric groups}\label{sec: metric_group}
In this section, we focus on metric groups. Most of the discussion extends to quotients: isometrically homogeneous spaces. For this more general viewpoint, we refer to \cite{LeDonne-Ottazzi, CKLNO}.

We generally refer to groups equipped with left-invariant distance functions as {\em metric groups}.\index{metric! -- group}
With the term \emph{metric Lie group}, we mean a Lie group equipped with a left-invariant distance function that induces the manifold topology.\index{metric! -- Lie group}\index{Lie group! metric --}
In general, when we have a topological space and we equip it with a distance function, we say that the distance function is {\em admissible} if it induces the topology of the space.\index{admissible! -- distance function}\label{def_admissible_distance}

\subsection{Smoothness of isometries between metric Lie groups}\label{sec:smooth}\index{smoothness! -- of isometry}
With the use of the results by Gleason-Montgomery-Yamabe-Zippin, Theorem~\ref{Montgomery-Zippin}, we deduce the differentiable regularity of isometries between metric Lie groups.
\begin{theorem}\label{smooth:thm}
Isometries between metric Lie groups are smooth maps.
\end{theorem}
Before giving the proof, we remark that in the Riemannian setting, the classical result of Myers and Steenrod gives smoothness of isometries; see \cite{Myers-Steenrod}, and more generally \cite{Capogna_LeDonne} for sub-Riemannian manifolds.
However, the following proof is different in spirit and, nonetheless, it will imply (see Theorem~\ref{prop:isometries_are_riemannian}) that such metric isometries are Riemannian isometries for some Riemannian structures. 

In the proof of Theorem~\ref{smooth:thm}, we will obtain a continuous isomorphism of Lie groups; hence, from what we saw in Section~\ref{sec smoothness of continuous homomorphisms}, this isomorphism is a smooth map. One can actually discuss analytic structures on Lie groups, and because of the BCH formula, every Lie group has a (unique) analytic structure, and in fact, we obtain that continuous group homomorphisms are analytic; see \cite[p.~117, Theorem~2.6]{Helgason}.
Hence, one can improve Theorem~\ref{smooth:thm} and prove the analyticity of isometries between metric Lie groups; see \cite[Theorem~1.1]{Kivioja_LeDonne_isom_nilpotent} and also \cite{LeDonne-Ottazzi}. We shall mainly focus on the $C^\infty$ regularity.

\begin{proof}[Proof of Theorem~\ref{smooth:thm}] 
Let \( F \colon M_1 \to M_2 \) be an isometry between metric Lie groups.
Without loss of generality, we may assume that 
\( F(1_{M_1})=1_{M_2} \) and that both \( M_1 \) and \( M_2 \) are connected, since left translations are smooth isometries and connected components of identity elements are open. By Theorem~\ref{Montgomery-Zippin}, for \( i \in \{ 1,2 \} \), the space \( G_i : = \Isom(M_i) \) is a Lie group smoothly acting on $M_i$; see Exercise~\ref{ex isometry group smoothly acting}. The conjugation map \( \mathrm{C}_F \colon G_1 \to G_2 \) defined as \( I \mapsto F \circ I \circ F^{-1} \) is a group isomorphism that is continuous with respect to the point-wise convergence. 
Hence, the map \( \mathrm{C}_F \) is smooth by Theorem~\ref{thm: continuous: smooth}.

Consider also the inclusion \( \iota \colon M_1 \to G_1 \), \( m \mapsto L_m \), which is smooth being a continuous homomorphism, and the orbit map \( \sigma \colon G_2 \to M_2 \), \( I \mapsto I(1_{M_2}) \), which is smooth since the action is smooth.
We deduce that \( \sigma \circ \mathrm{C}_F \circ \iota \) is smooth.
We claim that this map is \( F \). Indeed,
for every \( m \in M_1 \) it holds
\[
(\sigma \circ \mathrm{C}_F \circ \iota)(m) 
=\sigma(F \circ L_m \circ F^{-1})
=(F \circ L_m \circ F^{-1})(1_{M_2})
=F(m).
\qedhere
\]
\end{proof} 

\begin{remark}\label{analytic structure}
With the same techniques above, one can prove the following result, which is meant to summarize the smooth case and generalize it to the analytic category.
Let $M$ be a metric Lie group. Assume that $M$ is connected. Consider $M$ to be equipped with its unique analytic structure.
Then, the isometry group $\Isom(M)$ has the structure of a Lie group (finite-dimensional and with finitely
many connected components), which, equipped with its unique analytic structure, acts analytically on $M$.
Moreover, the stabilizers of the action $\Isom(M) \acts M$ are compact analytic Lie subgroups. 
\end{remark}


We get another consequence of the solution of Hilbert's fifth problem on Lie homogeneous spaces. Let $M=G/H$ be a Lie coset space of a Lie group $G$ modulo a compact subgroup $H$. Hence, the space $M$ has the structure of an (analytic-)smooth manifold, as from Theorem~\ref{thm_diff_structure_quotient}.
Assume that $M$ is equipped with a $G$-invariant admissible distance function. From Theorem~\ref{Montgomery-Zippin} with Exercise~\ref{ex isometry group smoothly acting}, it is immediate that the group of isometries ${\rm Iso}(M,d)$ of the manifold $M$ is a Lie group acting transitively on $M$.
Hence, the space $M$ admits an analytic structure for which ${\rm Iso}(M,d)$ acts by analytic maps.
One can prove that this analytic structure coincides with the initial one.

 \begin{theorem}[{\cite[Theorem 1.3]{LeDonne-Ottazzi}}]\label{regularity-global}
 Let $G/H$ be a homogeneous space of a Lie group $G$ modulo a compact subgroup $H$.
 Assume that $d$ is a $G$-invariant distance that induces the manifold topology.
 If $F:(G/H,d)\to (G/H,d)$ is an isometry, then $F$ is analytic.
\end{theorem}

\subsection{Submetries between metric groups}

Distance functions on metric groups do not necessarily descend to quotients. In fact, let us fix a metric group $(G,d)$ and a subgroup $H$ of $G$.
We stress that, by definition, the group $G$ acts on itself isometrically by left translations. 
 Naturally, the group $H$ acts on $G$ in two ways: by left translations or by right translations. Unless the group is normal, these actions have different orbits. The respective quotient spaces are\index{quotient! -- space} 
$$\mfaktor{H}{G} : = \{Hg\,:\, g\in G\}$$ 
and
$$G/H: = \{gH\,:\, g\in G\}.$$

On $G/H$ there is still a transitive $G$-action: $G\acts G/H$ as
$$\bar g .(g H) : = \bar g g H, \qquad \forall \bar g,g\in G.$$
However, there is no `adapted' metric geometry on $G/H$, in the sense that there may be no induced distance for which this action is by isometries; see Exercise~\ref{no_iso_act_Heis}.

Instead, regarding $\mfaktor{H}{G}$, the isometric action by left translations by elements of $G$ does not pass to this quotient. 
But, there is a good candidate for a distance function:
\begin{equation}\label{dist_HG}d(H g_1, H g_2): =\inf\left\{ d(h_1 g_1, h_2 g_2)\,:\, h_1, h_2 \in H\right\}, \qquad \forall h_1, h_2\in H.
\end{equation}
Note that this definition is in agreement with the notation that, for subsets $A$ and $ B $ of a metric space, we set $d(A, B): =\inf \{d(p,q);\; p\in A, \; q\in B\}$.

The function \eqref{dist_HG} is symmetric and satisfies the triangle inequality (see Proposition~\ref{prop:distance_quotient}); however, it may not be positively defined. In fact, in general situations, pathologies may be encountered. This is not the case if we assume that $H$ is a closed subgroup, so that $\mfaktor{H}{G}$ is a Hausdorff space, see Exercise~\ref{ex quotient Hausdorff}, and we assume that the right translations on $G$ are continuous with respect to the topology induced by the distance on $G$, as if $G$ is a topological group equipped with an admissible distance function. In this case, not only \eqref{dist_HG} metrizes the quotient, but in addition, the quotient map becomes a submetry in the sense of Definition~\ref{def submetry}.

There is an analog for metric spaces of this following result; see Exercise~\ref{prop: parallelfibers-submetry}.

\begin{proposition}\label{prop:distance_quotient}\index{submetry}
Let $G$ be a topological group equipped with an admissible left-invariant distance function $d$. Let $H$ be a closed subgroup of $G$. Then the function \eqref{dist_HG} satisfies
\begin{equation}\label{distance_quotient}
d (H g_1 , H g_2 ) 
=\inf\{ d( g_1 , h g_2 ) \,:\, h \in H\}, \qquad \forall g_1,g_2\in G,
\end{equation}
and is an admissible distance function on the quotient space $\mfaktor{H}{G}$ for which the projection $\pi: G\to \mfaktor{H}{G}$ becomes a submetry.
\end{proposition}
\begin{proof}
Because the distance function on $G$ is left-invariant, then for all $h_1\in H$ and all $g_1,g_2\in G$ we have 
\begin{equation*}
d (g_1 , H g_2 ) 
\stackrel{\rm def}{=} \inf\left\{ d( g_1 , h g_2 ) \,:\, h \in H\right\}=\inf\left\{ d( h_1g_1 , h g_2 ) \,:\, h \in H\right\}\stackrel{\rm def}{=} d(h_1 g_1, H g_2).
\end{equation*}
Taking the infimum over $h_1\in H$ we infer
\begin{equation}\label{distance_quotient_bis}
d (g_1 , H g_2 ) =d(H g_1 , H g_2 ) , \qquad \forall g_1,g_2\in G.
\end{equation}
We next check that it is a distance function. It is obviously symmetric and zero if $H g_1 = H g_2 $. Whereas, assume $H g_1 \neq H g_2 $, so $g_1\notin Hg_2$. Because we are assuming that $H$ is closed and that each right translation $R_{g_2}$ is continuous, we have that $Hg_2$ is closed. Therefore, there exists $r>0$ such that $B(g_1, r) \cap Hg_2 =\emptyset$. Hence, by \eqref{distance_quotient_bis} we infer
 $d (H g_1 , H g_2 )=d ( g_1 , H g_2 )\geq r> 0.$ So the distance function is positively defined.
 Regarding the triangle inequality, for each $g_1, g_2, g_3 \in G$ and for each $\eps>0$, let $h_\eps, h_\eps'\in H$ such that 
 $d(g_1, h_\eps g_2 ) \leq d (H g_1 , H g_2 ) + \eps$ and
 $d(h_\eps g_2, h_\eps h_\eps' g_3 )= d(g_2, h_\eps' g_3 ) \leq d (H g_2 , H g_3 ) + \eps$.
 Then 
 \begin{eqnarray*}
 d (H g_1 , H g_3 )&\leq& d(g_1, h_\eps h_\eps' g_3 )\\
 &\leq&d(g_1, h_\eps g_2 )+ d(h_\eps g_2, h_\eps h_\eps' g_3 )\\
 &\leq& d (H g_1 , H g_2 ) + d (H g_2 , H g_3 ) + 2\eps.
 \end{eqnarray*}
 By the arbitrariness of $\eps$, we get the triangle inequality.

Before proving that this distance is admissible, we check that $\pi$ satisfies \eqref{def:submetry} for all $p\in G$ and $r>0$.
For one inclusion, if $H q\in \pi(\bar B(p, r))=H \bar B(p,r)$ with $ q \in \bar B(p, r)$, then $d(H q,H p)\leq d( q, p)\leq r$ and so $H q\in \bar B(H p, r)=\bar B(\pi(p),r)$. 
Vice versa, if $H q\in \bar B(H p, r)$, i.e., $d( H q, p)\stackrel{\eqref{distance_quotient_bis}}{=}d( H q,H p)\leq r$,
 then there exists $h\in H$ such that for $\bar q: = hq $ we have $d(\bar q, p)\leq r$ and so $H q= H\bar q\in \bar B(H p, r)=\bar B(\pi(p),r)$.

Finally, we show that the distance defines the quotient topology on $\mfaktor{H}{G}$.
We stress that every submetry is an open map, is Lipschitz, and hence, is continuous. 
Since the quotient topology on $G/N$ is the only topology for which the projection $\pi: G\to G/N$ is continuous and open, the considered distance is admissible.
\end{proof}

Clearly, if we consider actions of normal subgroups, then we have quotient spaces that are homogenously metrized. Recall that $G/N=\mfaktor{N}{G}$ is a (Hausdorff) topological group when $N$ is a closed normal subgroup of a topological group $G$. Therefore, the group $G$ acts on the group $G/N$ by isometries with respect to the distance \eqref{dist_HG}. Thus, Proposition~\ref{prop:distance_quotient} gives the following consequence.

\begin{corollary}\label{prop:distance_quotient1}
Let $G$ be a topological group equipped with an admissible left-invariant distance function $d$. Let $N$ be a closed normal subgroup of $G$. Then the function
\begin{equation}\label{distance_quotient_normal}
d_{G/N}(g_1 N, g_2 N): =d(g_1 N, g_2 N) \stackrel{\rm def}{=}\inf\{ d(g_1 n_1, g_2 n_2) \,:\, n_1, n_2\in N\}
\end{equation}
is an admissible left-invariant distance function on the quotient group $G/N$ for which the projection $\pi: G\to G/N$ becomes a submetry.
\end{corollary}

 \subsection{Quasi-isometric equivalence of geodesic distances}
 
In this section, we establish the result that geodesic left-invariant metrics on the same group $G$ are quasi-isometric equivalent, according to Definition~\ref{def quasi-isometry}. Actually, we can relax the assumption that the metrics are geodesic and instead require them to be quasi-geodesic.

\begin{definition}[Quasi-geodesic space]\label{quasi-geodesic}\index{quasi-geodesic}
A metric space $M$ is said to be {\em quasi-geodesic} if there exist constants $C>0$ and $L>1$ such that every two points in $M$ can be join with a $(L,C)$-{\em quasi-arc} in the following sense:\index{quasi-arc} for all $x,x'\in M$, there exist $k\in \N$ and $x_0,x_1,\ldots, x_k\in M$ such that
$x_0=x$, $x_k=x'$,
$d(x_{i-1},x_i)\leq C$, for $i\in\{1,\ldots, k\}$, and
$\sum_{i=1}^k d(x_{i-1},x_i)\leq L d(x,x')+ C.$
Equivalently, there exists an $(L, C)$-{\em quasi-arc} joining $x$ and $x'$ if and only if there is an $(L, C)$-quasi-isometric embedding of an interval into the metric space, with $x$ and $x'$ in its image. See Definition~\ref{def quasi-isometry} for the notion of quasi-isometric embedding. 
\end{definition}

We stress that if two distances on a locally compact group are
 locally bounded (i.e., bounded on compact sets) and proper (i.e., the distance from a point is a proper map),
 then they have the same bounded subsets; see Exercise~\ref{ex compact = bounded}.

\begin{proposition}\label{quasiinvariant}\index{quasi-isometry}
Let $d$ and $d'$ be quasi-geodesic left-invariant distances on a group. Assume that $d$ and $d'$ have the same bounded subsets.
 Then there exist constants $c\geq 0$ and $L\geq 1$ such that
$L^{-1}d-c\leq d'\leq Ld+c.$
\end{proposition}

\begin{proof}
Since $d$ is quasi-geodesic, there are two constants $C_1\geq 0$ and $L_1\geq 1$ with the following property: 
Given a group element $g$, there are $g_1,\dots,g_n$ such that $d(1,g_i) \leq C_1$, for all $i$, and $g=g_1\cdot \ldots \cdot g_n$, while $\sum_1^n d( 1, g_i) \leq L_1d( 1, g)+C_1$. 
Grouping some $g_i$'s together if necessary, 
 we may assume that $C_1/2 \leq d( 1, g_i)$, for $i\in\{1,\ldots, n-1\}$, still having the weaker uniform condition $d( 1, g_i) \leq 2C_1$, for $i\in\{1,\ldots, n\}$.
We claim that \begin{equation}\label{3518435e468465}
n \leq \dfrac{2L_1}{C_1} d( 1, g)+3.
\end{equation} 
Indeed, using the lower bound on each $d( 1, g_i)$ we get
$$(n-1)C_1/2 \leq \sum_{i=1}^n d( 1, g_i) \leq L_1d( 1, g)+C_1,$$
from which the claim \eqref{3518435e468465} follows.

By our assumption on the distances, since the values $d( 1, g_i)$ are uniformly bounded by $2 C_1$, then
also 
the values
$d'( 1, g_i)$ are uniformly bounded, say by some constant $\tilde C>0$. Then 
$$d'( 1, g) \leq \sum_{i=1}^n d'( 1, g_i) \leq \tilde C n \stackrel{\eqref{3518435e468465}}{\leq} Ld( 1, g)+c,$$ for $L: =\dfrac{2\tilde CL_1}{C_1}$ and $c: =3\tilde C$. The proposition follows by 
left invariance and by
exchanging the roles of $d$ and $d'$.
\end{proof}

The above result will apply in particular to Finsler and sub-Finsler left-invariant metrics on Lie groups.
In fact, every pair of left-invariant subFisler distances on the same Lie group are quasi-isometric equivalent. For example, the Riemannian Heisenberg group and the sub-Riemannian Heisenberg group are quasi-isometric. In fact, quasi-isometries can completely change the local geometry. We conclude the section with a question.
\begin{problem} Are there translation-invariant quasi-geodesic metrics on $(\R^2,+)$ that are not quasi-isometric to a geodesic distance? Are there such metrics that, in addition, are boundedly compact?
 \end{problem}

\section{Haar measures and polynomial growth}\label{sec_Haar_poly_growth}

Every Lie group, as every locally compact group, has a natural class of measures: the {\it Haar measures}. A measure $\mu$ on a group $G$ is called {\em left-invariant}, if for every $g\in G$ and every set $E$ on which $\mu $ is defined, we have that $\mu $ is defined on $L_g^{-1}(E)$ and 
$$ \left((L_g)_\#\mu\right)(E): =\mu\left(L_g^{-1}(E)\right)=\mu(E).$$
If $G$ is equipped with a topology, then a Borel measure $\mu$ on $G$ is said to be a {\em Radon measure} if it is finite on compact sets, outer regular on Borel sets, and inner regular on open sets; see \cite[page~212]{Folland_book}.
If $G$ is a topological group, then a Borel measure $\mu$ on $G$ is called 
a {\em left-Haar measure}, or simply a \emph{Haar measure}, if it is left-invariant, Radon, and not the zero measure.
Similarly, a \emph{right-Haar measure} is a Radon non-zero measure that is right invariant.\index{Haar measure}

On every locally compact group $G$, there are Haar measures; see \cite[Theorem~11.8]{Folland_book}.
Moreover, Haar measures are unique in the following sense: 
	Every two left-Haar measures $\mu_1$ and $\mu_2$ on the same locally compact group differ by a constant multiplicative factor: $\mu_1= c \mu_2$, for some $c>0$. We refer to \cite[Theorem~11.9]{Folland_book} for a proof.
 
Using the Haar measure, we can consider the asymptotic growth of locally compact groups:

\begin{definition}[Polynomial growth \& exponential growth]\index{exponential! -- growth}\index{polynomial! -- growth}
A 
 Lie group $G$ is said to have {\em polynomial growth} if the Haar measures of powers of compact sets grow at most polynomially.
Namely, called $\mu$ a Haar measure of $G$, for every compact set $U\subseteq G$ there is 
$Q>0$
such that $\mu(U^n)\leq n^Q$, for all $n\in \N$. 
A 
 Lie group $G$ is said to have {\em exponential growth} if
 for every compact neighborhood $U\subseteq G$ of $1_G$ there is $k>0$ such that 
$\mu(U^n)\geq k^n$, for all $n\in \N$. 
\end{definition}

We stress that in the case $G$ is connected, to check whether it has polynomial growth, it is enough to check the property for some compact set $U$ that is a neighborhood of the identity element; see Exercise~\ref{ex_poly_growth_cpt}. 
Similarly,
a group has exponential growth if and only if for some compact neighborhood $U$ of the identity element there is $k>0$ such that $\mu(U^n)\geq k^n$, for all $n\in \N$.

\subsection{Type (R) and Guivarc’h-Jenkins Theorem}\label{sec type (R)}

In this section, we will present an algebraic characterization of groups with polynomial growth. We begin by introducing the algebraic property.

\begin{definition}[Type (R)]\label{def_typeR}\index{type (R)}
A Lie algebra $\g$ is said to have {\em type (R)} if the eigenvalues of $\ad_X$ are purely imaginary for each $X\in \g$. 
A Lie group $G$ is said to have {\em type (R)} if the eigenvalues of $\Ad_g$ are of absolute value 1 for each $g\in G$. 
\end{definition}

The two notions are related:
A connected Lie group has type (R) if and only if its Lie algebra has type (R).
A proof can be found in \cite[Proposition 1.3’]{MR0349895}.

Guivarc’h and Jenkins independently proved the following equivalence: 

\begin{theorem}[Guivarc’h-Jenkins, \cite{MR0369608, MR0349895}]\label{Guivarch_Jenkins_Theorem}\index{Guivarc’h-Jenkins Theorem}\index{Theorem! Guivarc’h-Jenkins --}
A connected Lie group has polynomial growth if and only if it is of type (R).
\end{theorem}

\begin{proof}[Sketch of a proof.]
If a group is not of type (R), then we will soon prove in Proposition~\ref{noR_exp_growth} that the group has exponential growth.
 
Vice versa, if the group has type (R), then it actually can be made isometric to a nilpotent Lie group, called its {\em nilshadow}. 
Nilpotent Lie groups have polynomial growth because they are asymptotic to Carnot groups. Proving this latter statement will be the main aim of Sections~\ref{sec_Pansu_asymp} and \ref{sec_largescale}.
\end{proof}

\begin{example}[Models of groups not of type (R)]\label{model_typeR}\index{$G_z$}
We present now some Lie groups, denoted by $G_z$, with $z\in \C\setminus i\R$, that are in some sense the reason for exponential growth. Namely, we will deduce that if a group has exponential growth, then it has one of such groups as a subgroup. Given a complex number $z$ that is not purely imaginary, we distinguish two cases: either $z$ is real or not. 

If $z\in \R\setminus\{0\}$, then we consider the 2D Lie algebra that has a basis $X,Y$ with operation $[X,Y,]= z Y$. Let $G_z$ be the associated simply connected Lie group. Such a Lie group does not have type (R) and, when metrized with a left-invariant Riemannian metric, it has constant negative curvature. When properly renormalized, it is the hyperbolic plane, which has exponential growth.

If instead $z = a+bi$ with $a,b\neq 0 $, then we consider the 3D Lie algebra that has a basis $X_1,X_2, X_3$ with operation 
\begin{equation*}
 \begin{array}{ccl}
 \quad [X_1,X_2 ] &=& a X_2 + b X_3,\\
 \quad [X_1,X_3 ] &=& a X_3 -b X_2,\\ 
 \quad [X_2,X_3]&=& 0.
 \end{array}
\end{equation*}
Let $G_z$ be the associated simply connected Lie group. Such a Lie group does not have type (R) and
 it has exponential growth; see Exercise~\ref{ex:Gdelta}. These are examples of Heintze groups, as we will discuss in Chapter~\ref{ch_Heintze}.
\end{example}

\begin{proposition}\label{noR_exp_growth}
If $G$ has not type (R), then there is $z\in \C$ and a subgroup $H$ of $G$ such that $H$ is isomorphic to $G_z$. Consequently, the group $G$ has exponential growth.
\end{proposition}

\begin{proof}
We begin by searching for a subalgebra of $\g$ that is isomorphic to the Lie algebra of some $G_z$.
Let $X_1\in \g$ for which $\ad_{X_1}$ had some eigenvalue in $\C\setminus i\R$. We distinguish two cases: either there is a real eigenvalue, or not.
If there is a real eigenvalue $z\in \R^*$ with eigenvector $X_2$, then the span of $X_1$ and $X_2$ is the Lie algebra of $G_z$. 

If, instead, there are no real eigenvalues for $\ad_{X_1}$, then take an eigenvalue $z = a+bi$ with $a,b\neq 0 $. Then there is $X_2$ and $X_3\in \g$ such that \begin{equation*}
 \begin{array}{ccl}
 \quad [X_1,X_2 ] &=& a X_2 + b X_3,\\
 \quad [X_1,X_3 ] &=& a X_3 -b X_2
 \end{array}.
\end{equation*}
 Using these relations and the Jacobi identity, we have
 $ \ad_{X_1}[X_2,X_3]= 2 a [X_2,X_3]$.
 Since we are assuming that there are no real eigenvalues and $a\neq 0$, we deduce 
 $ [X_2,X_3]= 0$. 
The span of $X_1$, $X_2$, and $X_3$ gives the Lie algebra of $G_z$. 

Once we have found a Lie algebra that isomorphic to the Lie algebra of $G_z$'s, we have a subgroup $H$ of $G$ with such a Lie algebra, recall Theorem~\ref{teo1145bis}. Being $G_z$ simply connected, then $H$ is a quotient of $G_z$. However, one can show that $G_z$ has no 
normal discrete subgroups (see Exercise~\ref{Normal_discrete_connected_central} and~\ref{ex:Gdelta})
Thus, the groups $G_z$ and $H$ are isomorphic.

We stress that $H$ is necessarily closed in $G$. In fact, the group $H$ is properly embedded in $G$ because, due to the possible $\ad_X$ in $G_z$, when elements $h_n\in H$ diverge in $H$ then $\Ad_{h_n}$ diverge as maps on the Lie algebra of $H$ and then as a map on the Lie algebra of $G$. Hence, they cannot converge in $G$.

We metrize $G$ with a geodesic left-invariant distance, e.g., Riemannian, and consider on $H$ the induced path distance. The intersection of the closed unit ball $\bar B_G(1,1)$ and $H$ is compact. Thus, for some $\delta>0$ we have $\bar B_H(1,\delta) \supseteq \bar B_G(1,1) \cap H$.
Since $G_z$ and $H$ have exponential growth, we have that maximal $\delta$-separated sets in the balls of radius $R$ in $H$ grow exponentially in $R$. Consequently, 
maximal $1$-separated sets in the balls of radius $R$ in $G$ grow at least exponentially in $R$. Therefore, recalling Exercise~\ref{ex_net_growth}, the group $G$ has exponential growth.
 \end{proof}

\section{Isometrically homogeneous spaces with dilations (first part)}\label{sec: self-similar}

Next, we present the notion of dilation on a metric space. 
Fixed a scalar number $\lambda$, we will denote by $\delta_\lambda$ one such a dilation by factor $\lambda$.
However, there is no such map for general metric spaces, and if it exists, it may not be unique. 
In some settings, like in vector spaces or, more generally, in Carnot groups, we have canonical maps given by the algebraic structure.

Let $M$ be a metric space with distance $d$, and let $\lambda\in (0,+\infty)$. A \emph{dilation of factor $\lambda$} on $M$ is a surjective map $\delta_\lambda:M\to M$ such that 
\[
d(\delta_\lambda(p),\delta_\lambda(q)) = \lambda d(p,q), \qquad \forall p,q\in M .
\]
We say that a metric space $X$ is {\em self-similar} if it admits a dilation of factor $\lambda\in(0,+\infty)\setminus\{1\}$. Such dilations are also called {\em homotheties}.\index{self-similar}\index{homothety}\index{dilation}

Clearly, the dilations of factor 1 are the isometries. Hence, isometrically homogeneous spaces have plenty of them.

In the next result, we show that, under some mild topological assumptions, if a metric space is isometrically homogeneous and self-similar, then it has the structure of a Lie coset space.
After the revision of the theory of nilpotent Lie groups done in Chapter~\ref{ch_Nilpotent}, we will reach a stronger conclusion: the space has the structure of a nilpotent metric Lie group for which the dilations fixing the identity element are group isomorphisms; see Theorem~\ref{teo05171840}.

\begin{theorem}[to be improved in Theorem~\ref{teo05171840}]\label{teo05171840baby}
	Let $M$ be a metric space that is
	
	 \begin{minipage}[c]{0.4\linewidth}
	 (1) locally compact,\\ (3) isometrically homogeneous, 
	 \end{minipage}
	 	 \begin{minipage}[c]{0.6\linewidth}
	(2) locally connected,\\ (4) self-similar. 
	 \end{minipage}
	 \\
	Then, $M$ has the structure of a smooth manifold. In fact, there is a connected Lie group $G$, a compact subgroup $S$ of $G$, and a $G$-invariant distance on the coset space $G/S$ such that $G/S$ is isometric to $M$.
\end{theorem}
 
Before proving Theorem~\ref{teo05171840baby}, we provide reasons why the considered space has finite topological dimension and is connected. 
The aim is to use Gleason-Montgomery-Yamabe-Zippin's result, Theorem~\ref{Montgomery-Zippin}. 
	
As a side remark, if a locally compact isometrically homogeneous space admits dilations of {\em every} factor, then it is connected; 
 see Exercise~\ref{connectivity_connected_balls_homotetic}. However, in Theorem~\ref{teo05171840baby}, we assume only the presence of at least {\em one} dilation (together with all its iterations, of course).

\begin{proposition}\label{prop:doubling_homog} 
Every locally compact, self-similar, isometrically homogeneous space is metrically doubling, and in particular, it has finite Hausdorff dimension and finite topological dimension.	\index{doubling! -- metric space} 
\end{proposition}
\begin{proof} We need to prove that for every such metric space $X$, there is $C>0$ such that every ball of radius $r>0$ in $X$ can be covered with $C$ balls of radius $r/2$. Since $X$ is locally compact, it has a closed ball $\bar B(x_0,r_0)$ that is compact.
Up to replacing the dilation by its inverse, we may assume that its factor 
$\lambda$ is greater than $1$. 
From Exercise~\ref{Ex_local_to_global_dilations}, the metric space is boundedly compact.
Hence, the closed ball $\bar B(x_0, \lambda r_0)$ can be covered with finitely many balls of radius $r_0/2$, say $C$ many of them. 
Therefore, for every $s\in [1,\lambda]$ each ball $B(x_0, s r_0)$ can be covered with $C$ balls of radius $s r_0/2$. By self-similarity and homogeneity, every other ball can be covered with less than $C$ balls of half radius. 
 Finally, recall that doubling metric spaces have finite Hausdorff dimension and, hence, finite topological dimension; see \cite[Theorem~8.14 and Exercise~10.16]{Heinonenbook}. 
\end{proof}

\begin{proof}[Proof of Theorem~\ref{teo05171840baby}]

Let $\delta_\lambda:M\to M$ be a dilation of factor $\lambda\in(0,+\infty)\setminus\{1\}$.
	Since $\lambda $ is surjective, it is a homeomorphism and its inverse $(\delta_\lambda)^{-1}:M\to M$ is a dilation of factor $\frac1\lambda$.
	Hence, up to replacing $\delta_\lambda$ with its inverse if necessary, we may assume $\lambda<1$.
	Observe that the space is complete (see Exercise~\ref{completeness_homog}) and apply the Banach Fixed Point Theorem: The contraction $\delta_\lambda$ has a unique fixed point, which we denote by $o$. We consider the orbit map with respect to $o$, i.e., the map $f\in \Isom(M,d)\mapsto \pi(f): =f(o)\in M$.

	

Recall Proposition~\ref{Iso_gp_loc_cpt} regarding the fact that
$\Isom(M,d)$ is a locally compact group acting continuously on $M$ with compact stabilizers. Moreover, by Theorem~\ref{Helgason_Theorem_3.2}, the quotient map $\pi$ 
induces a homeomorphism between $M$ and the manifold $ \Isom(M,d)/S$, where $S: ={\rm Stab}(o)=\{f\in \Isom(M,d):f(o)=o\}$.

	In order to apply Gleason-Montgomery-Yamabe-Zippin's Theorem~\ref{Montgomery-Zippin} to $(M,d)$, we need to observe that $M$ 
	has finite topological dimension and
	is connected; see
	 Proposition~\ref{prop:doubling_homog} 
		 and Exercise~\ref{Ex_local_to_global_dilations}, respectively.
	Hence, thanks to Theorem~\ref{Montgomery-Zippin}, $\Isom(M,d)$ and $M$ have some differential structures such that $\Isom(M,d)$ is a Lie group and the action of $\Isom(M,d)$ on $M$ is smooth.
	To conclude, we take $G: =\Isom(M,d)^\circ$ to be
	 the identity component of $\Isom(M,d)$ and the stabilizer $S$ as compact subgroup. 
\end{proof}

 
Suppose that $G$ is a Lie group, and let ${\rm {Der}}(\g)$ the space of derivations on its Lie algebra; see Definition~\ref{def:derivation}.
By Proposition~\ref{Prop_Lie_der}.ii, every multiplicative one-parameter group $\R_+\to\Aut(G)$, $\lambda\mapsto\delta_\lambda$, of Lie automorphisms is determined by some derivation $A\in{\rm {Der}}(\g)$ such that 
\begin{equation}\label{eq12061455}
(\delta_\lambda)_* = \lambda^A : = e^{(\log\lambda)A} .
\end{equation}
Such an $A$ is the called the {\em infinitesimal generator} of $\lambda\mapsto\delta_\lambda$, and such $(\delta_\lambda)_\lambda$ form the {\em (multiplicative) one-parameter subgroup of automorphisms determined} by the derivation $A$.\index{infinitesimal! -- generator}\index{one-parameter subgroup! -- of automorphisms determined by a derivation}

\begin{definition}[Homogeneous metric group]\label{homogeneous metric group}\index{homogeneous! -- group}\index{group! homogeneous --} 
Let $G$ be a Lie group and $A\in{\rm {Der}}(\g)$ a derivation on its Lie algebra.
Assume that $d$ is a left-invariant distance function on $G$ 
 for which
 \begin{equation}\label{eq12061454}
 d(\delta_\lambda x,\delta_\lambda y) = \lambda d(x,y)
 \qquad\forall x,y\in G,\ \forall \lambda>0 ,
 \end{equation}
 where $(\delta_\lambda)_\lambda$ is the one-parameter subgroup of automorphisms determined by $A$.
Then, we say that $d$ is a \emph{$A$-homogeneous} distance function on $G $ and $(G,d)$ is a \emph{$A$-homogeneous metric group}.\index{$A$-homogeneous distance}\index{homogeneous! -- metric group}
\end{definition}

\section{Exercises}

\begin{exercise}[Characterisations of proper actions]\label{Prop 21.5 from Lee}
 Let $M$ be a manifold and $G$ a Lie group acting continuously on $M$. The following are equivalent:
(a) the action is proper; (b) if $(p_i)_{i\in \N}$ is a sequence in $M$ and $(g_i)_{i\in \N}$ is a sequence in $G$ such that both $(p_i)_{i\in \N}$ and $(g_i.p_i)_{i\in \N}$ converge, then a subsequence of $(g_i)_{i\in \N}$ converges; (c) for every compact subset 
$ K \subseteq M$, the set $\{g\in G: (g.K)\cap K\neq \emptyset\}$ is compact.
\\{\it Hint.} A solution can be found in \cite[Proposition~21.5]{Lee_MR2954043}.
 \end{exercise}

\begin{exercise}\label{ex_orbit_topo}
For a continuous action $G\acts X$ of a topological group on a topological space, with quotient map $\pi$, declare a subset $U\subseteq \mfaktor{G}{X}$ to be open if and only if $\pi^{-1}(U)$ is open in $X$. 
This definition gives a topology on $X$, and it makes the map $\pi$ continuous and open.
Moreover, this is the only topology that makes $\pi$ continuous and open.
\end{exercise}

\begin{exercise}\label{ex quotient Hausdorff}
Let $G$ be a topological group and $H<G$ a topological subgroup.
The actions from Example~\ref{ex_subgroup_acting} are continuous.
If $H$ is closed, then for each of these actions, we have that $\mfaktor{H}{G} $ is a Hausdorff topological space when equipped with the quotient topology. 
\\{\it Hint.} Check Exercise~\ref{G/H Hausdorff}, using that topological groups are $T_3$-regular.
\end{exercise}

\begin{exercise}\label{ex:action_left_translations}
Given a group $G$, and denoting by $L_g$ its left translation by $g$, we have that 
 \( g \mapsto L_g \) is an action of \( G \) on \( G \).
 Moreover, if $G$ is a Lie group, then this is a Lie action.
\end{exercise}

\begin{exercise}\label{ex_dim_quot}
If \( G \acts M \) is a proper and free Lie action, then \( \dim \mfaktor{G}{M} = \dim M - \dim G \).
\end{exercise}


Exercises~\ref{ex:quotient:covering1}--\ref{extra_ex_actions} are extra exercises on Lie actions.

%

\begin{exercise}\label{ex:quotient:covering1}
Let \( \Gamma \) be a discrete Lie group, \( M \) a smooth manifold, and \( \Gamma \acts M \) a Lie action that is proper and free. Then, there exists a differentiable structure on \( \mfaktor{\Gamma}{M} \) such that \( M \to \mfaktor{\Gamma}{M} \) is a smooth covering map.
\end{exercise}

\begin{exercise}
Let \( G \) be a Lie group and \( \Gamma < G \) a closed subgroup. Then, the natural action \( \Gamma \acts G \) by left-translations is proper and free.
\\
{\it Solution.} 
The action by left-translations is free because \( L_g(p)=p \) implies \( g = 1_G \). To see that the action is proper, notice that the map 
\( \Psi \colon G \times G \to G \times G ,
(g_1,g_2) \mapsto (g_1 g_2, g_2) \),
has the inverse \( (g_1,g_2) \mapsto (g_1 g_2^{-1},g_2) \) and thus \( \Psi \in {\rm{Diffeo}}(G \times G) \). The map that is claimed to be proper is \( \tilde \Theta = \Psi|_{\Gamma \times G} \). Given a compact set \( K \subset G \times G \), we have 
\(\tilde \Theta^{-1}(K) = \Psi^{-1}(K) \cap (\Gamma \times G). \)
Since \( \Gamma \) is closed in \( G \), then \( \tilde \Theta^{-1}(K) \) is an intersection of a compact set with a closed set.
\end{exercise}


\begin{exercise}\label{ex:quotient:covering2}
 If \( G \) is a Lie group and \( \Gamma < G \) a discrete subgroup, then there exists a differentiable structure on \( \mfaktor{\Gamma}{G} \) such that the natural map \( G \to \mfaktor{\Gamma}{G} \) is a smooth covering map.
\end{exercise}


\begin{exercise}\label{extra_ex_actions}
Let \( G \) be a Lie group and \( M \) a simply connected smooth manifold such that \( \dim G \le \dim M \). If there exists a transitive Lie action \( G \acts M \), then \( M \) has a structure of a Lie group (i.e., \( M \) is diffeomorphic to a Lie group). 
\end{exercise}


\begin{exercise}\label{G/H Hausdorff} If $H$ is a closed subgroup of a Lie group $G$, then $G/H$ is Hausdorff.\\
{\it Solution.}
Consider the set 
 $$R: =\left\lbrace (h_1,h_2)\in G\times G\mid \textrm{ there exists } h\in H \textrm{ such that } h_1=h_2 h\right\rbrace,$$
which is a closed subset of $G\times G$. 
So, if $h_1H$ and $h_2H$ are two distinct cosets, then $(h_1,h_2)$ does not belong to $R$, thus one can find open sets $h_1\in V$ and $h_2\in W$ such that $(V\times W)\cap R=\emptyset$. Since the projection map $\pi$ is an open map, the sets $\pi(V)$ and $\pi(W)$ are disjoint and open in $G/H$, containing $h_1H$ and $h_2H$, respectively.
\end{exercise}

\begin{exercise}\label{G/H second countable} If $H$ is a subgroup of a Lie group $G$, then $G/H$ is second countable.
\\
{\it Hint.} Every countable basis for the topology of $G$ projects to a countable basis for $G/H$.
\end{exercise}


{\color{black}
\begin{exercise} 
Let $M$ be a sub-Finsler manifold.
Let $G\acts M$ be a Lie group action that is proper, free, and by isometries.
Then, there exists a sub-Finsler distance $d$ on $\mfaktor{G}{M}$ that makes $\pi$ a submetry. 
\\{\it Hint.} Regarding the existence of a submetric distance, check Exercise~\ref{ex_action_isometry_submetry}. Regarding the fact that it is sub-Finsler, check Proposition~\ref{prop_subFin_submetry_Lie}.
\end{exercise}
}


\begin{exercise}\label{ex:homog:dist}\index{homogeneous! -- distance}\index{dilation}
Let $(X,d)$ be a metric space. Within this exercise, we say that the distance function $d$ is a {\em homogeneous distance} if for every $\lambda>0$ there is a bijection $\delta_\lambda:X\to X$ such that 
$$d(\delta_\lambda (x),\delta_\lambda (x'))= \lambda d(x, x'), \qquad \forall x,x'\in X.$$
In this case, we also say that $\delta_\lambda$ is a {\em dilation} by $\lambda$. We have the following properties:
\\ (i) on $\R^n$, the Euclidean distance $d_E$ is homogeneous;
\\ (ii) on $(\R\times\{0\})\cup (\{0\} \times \R)$ as subset of $\R^2$, the distance function $d_E$ is homogeneous;
\\ (iii) on $\mathbb S^1$ as subset of $\R^2$, the distance function $d_E$ is not homogeneous.
\\Question: Which of these spaces are isometrically homogeneous, according to Definition~\ref{def isometrically homogeneous}?
\end{exercise}
\begin{exercise}\label{ex action compact open}
If $G\acts X$ continuously, then the induced map $G\to {\rm Homeo}(X)$ is continuous with respect to the compact-open topology.
 \end{exercise}

\begin{exercise}
If a non-necessarily Hausdorff topological group $G$ acts on a Hausdorff space continuously and faithfully, then $G$ 
 is Hausdorff. 
 \\{\it Hint.} Use Exercise~\ref{ex action compact open}.
 \end{exercise}

\begin{exercise}\label{second countable VS sigma-compact}
 Recall that a topological space $X$ is called {\em first countable} if its topology admits a countable local basis at every point. It is called {\em second countable} if its topology admits a countable basis.
 While $X$ is {\em $\sigma$-compact} if it is the countable union of compact sets.\index{second countable}\index{$\sigma$-compact}\index{first countable}
\begin{description}
\item[\ref{second countable VS sigma-compact}.i.]
Every second countable space is first countable. 
\item[\ref{second countable VS sigma-compact}.ii.]
Every second countable, locally compact space is $\sigma$-compact. 
\item[\ref{second countable VS sigma-compact}.iii.]
There is a topological space that is $\sigma$-compact and locally compact, but it is not second countable. {\it Hint.} Consider the excluded-point topology on the real numbers.
\end{description}
 \end{exercise}

\begin{exercise}\label{limit of isometries}
Pointwise limits of isometric embeddings are isometric embeddings.
\end{exercise}

\begin{exercise}\label{open_component_trans}
Let $G$ be a locally compact group with a countable basis acting continuously and transitively on a connected locally compact Hausdorff space $X$. 
 Then, every open subgroup of $G$ acts transitively on $X$.
 \\ {\it Hint}. The proof is very similar to the one we provided for Proposition~\ref{identity_component_trans}. 
\end{exercise}

\begin{exercise}\label{Cornulier-Harpe5.B.4}\index{boundedly compact}\index{isometry! -- group}
Let $M$ be a boundedly compact metric space. Equip the isometry group $\Isom(M)$ with the compact-open topology. 
\begin{description}
\item[\ref{Cornulier-Harpe5.B.4}.i.]
The action $\Isom(M)\acts M$ is continuous and proper. 
\item[\ref{Cornulier-Harpe5.B.4}.ii.]
The topology on $\Isom(M)$ is second-countable and locally compact, and hence, $\sigma$-compact. 
\item[\ref{Cornulier-Harpe5.B.4}.iii.]
Check that $\Isom(M)$ is a topological group.
\end{description}
{\it Hint.} One can find a proof in \cite[Lemma~5.B.4]{Cornulier-Harpe}. Recall also Exercise~\ref{second countable VS sigma-compact}.ii.
\end{exercise}

\begin{exercise}\label{ex finitely many components}
Let $G\acts M$ be a Lie group action that is transitive, faithful, and with compact stabilizers.
Then, the topological space $G$ has finitely many connected components.
\\{\it Hint.} On the one hand, the quotient $G/G^\circ$ is discrete since connected components of manifolds are open. On the other hand, this space $G/G^\circ$ is compact because, fixing $p\in M$, for every sequence $g_nG^\circ$ there exists $g_n'\in g_nG^\circ$ such that $g_n'.p=p$. (Recall that $G^\circ$ acts transitively by Proposition~\ref{identity_component_trans}). Thus, the sequence $g_n'$ subconverges.
\end{exercise}
\begin{exercise}
Let $M$ be a topological manifold. Let $G$ be a locally compact group with a countable basis.
Let $G\times M\to M$ be a continuous, effective, and transitive action of $G$ on $M$. 
Then $G$ is a Lie group, and $M$ can be equipped with a differentiable structure so that $G$ acts by diffeomorphisms.
\end{exercise}


\begin{exercise}\label{ex isometry group smoothly acting}
Let $M$ be a connected metric Lie group. Then, the topological group $\Isom(M)$ has the structure of a Lie group smoothly acting on $M$.
\\{\it Hint.} In addition to Theorem~\ref{Montgomery-Zippin} recall that Lie groups have unique smooth structures as in Corollary~\ref{Uniqueness of Lie structures}.
\end{exercise}

\begin{exercise}\label{ex_action_isometry_submetry}
Let $M$ be a metric space. Let $G\acts M$ be a continuous action by isometries with closed orbits. Then, the quotient map
$\pi: M\to \mfaktor{G}{M}$ becomes a submetry with respect to some unique distance on $\mfaktor{G}{M}$. 
\\{\it Hint.} Check Exercise~\ref{prop: parallelfibers-submetry}. 
\end{exercise}

\begin{exercise}\label{no_iso_act_Heis}\index{Heisenberg! -- group} 
Let $H$ be the Heisenberg group with Lie algebra spanned by $X, Y, [X, Y]$. Let $L: =\exp(X)$.
Then $H/L$ is a Lie homogenous space on which $H$ acts on the left. The stabilizer of this action is homeomorphic to $L$, which is not relatively compact inside the space of homeomorphisms of $H/L$. In particular, this action cannot be isometric. 
\end{exercise}


 \begin{exercise}\label{ex compact = bounded} 
 Let $X$ be a topological space equipped with a distance $d$ that does not necessarily induce the same topology. Assume that $d$ is locally bounded (in the sense that it is bounded on compact subsets of $X$) and proper (in the sense that the distance function from a point in $X$ is a proper map). Then, the sets that are bounded with respect to $d$ are exactly the precompact sets. 
 \end{exercise}
 

\begin{exercise}\label{ex_poly_growth_cpt}\index{polynomial! -- growth}
Let $G$ be a connected Lie group and let $\mu_G$ a Haar measure on $G$. Let $U\subseteq G$ be a compact
 neighborhood of $1_G\in G$. If there is $C>0$ such that $\mu_G(U^n) 
	\le C n^Q ,$ for all $n\in\N$, then $G$ has polynomial growth.
 \\{\it Hint.} By Exercise~\ref{Prop:generating}, for every other compact set $\tilde U\subseteq G$ there is $k\in \N$ such that $\tilde U\subseteq U^k$.
\end{exercise}

\begin{exercise}\label{ex_poly_growth}\index{polynomial! -- growth}\index{Haar measure}\index{doubling! -- measure}
 Let $M$ be a metric space as in Theorem~\ref{teo05171840baby}.
Let $G$ be the identity component of the isometry group $\Isom(M,d)$, equipped with the topology of the uniform convergence on compact sets. 
Let $\pi:G\to M $ be the orbit map $f\mapsto f(o)$. 
Let $\mu_G$ be a Haar measure on $G$, and consider the push forward $\mu_M: =\pi_*\mu_G$
\begin{description}
\item[\ref{ex_poly_growth}.i.] The measure $\mu_M$ is a $G$-invariant Radon measure on $M$ such that $\mu_M(E)=\mu_G(\pi^{-1}(E))$ for all open sets $E\subset M$.

\item[\ref{ex_poly_growth}ii.] 	The measure $\mu_M$ is a doubling measure (in the sense of Exercise~\ref{ex_doubling_measure}) and that there are $C>0$ and $Q>0$ such that for all $p\in M$ and all $r\ge 1$ 
one has
$	\mu_M(B_d(p,r)) \le C r^Q .$
\item[\ref{ex_poly_growth}.iii.] Deduce that the topological group $G$ has polynomial growth.
\end{description}
{\it Hint.} Recall Exercise~\ref{ex_doubling_measure2}. Consider the set $U: =\pi^{-1}(B_d(o,1))$. Then $\pi(U^n)\subset B_d(o,n)$ and 
	\(
	\mu_G(U^n) 
	\le \mu_G(\pi^{-1}(\pi(U^n)))
	= \mu_M(\pi(U^n))
	\le \mu_M(B_d(o,n))
	\le C n^Q .
	\)
\end{exercise}

\begin{exercise}\label{ex:Gdelta}\index{$G_z$}
For $z\in \C\setminus i\R$, let $G_z$ be the Lie group from Example~\ref{model_typeR}. 
Then, the center of $G_z$ is trivial
and $G_z$ has exponential growth.
\\{Hint.} See Exercise~\ref{ex:Gdelta2}
\end{exercise}

\begin{exercise}\label{Normal_discrete_connected_central}
Normal discrete subgroups of connected topological groups are central.
\end{exercise}

\begin{exercise}[Big separated net gives big growth]\label{ex_net_growth}
Let $G$ be a connected Lie group equipped with a left-invariant admissible geodesic distance.
Let $\rho(n)$ be the maximal cardinality for a 2-separated set in the ball $B(1_G,n)$, for $n\in \N$. 
Then $\rho(n)\mu(B(1_G,1))\leq\mu(B(1_G,1)^n) \leq \rho(n)\mu(B(1_G,2))$, for all $n\in \N$.
\\{\it Hint.} The distance being geodesic implies $B(1_G,1)^n =B(1_G,n)$. If two points have a distance of at least 2, their balls of radius $1$ are disjoint.
\end{exercise}


\begin{exercise}\label{Ex_local_to_global_dilations}\index{dilation}
	Let $M$ be a metric space that admits a dilation $\delta_\lambda: M\to M$ of factor $\lambda >1$ with a fixed point $o$.
	 If 
 	$M$ is locally connected (resp., locally simply connected; resp., locally compact; resp., locally contractible), then $M$ is connected
	 (resp., simply connected; resp., boundedly compact; resp., contractible).
	\\ {\it Hint:} Take a `good' neighborhood $U$ of $o$ and $r>0$ such that
	 $B(o,r)\subset U$. 
	Then, we have $M = \bigcup_{n=1}^\infty B(o,\lambda^{n}r) =\bigcup_{n=1}^\infty \delta_\lambda^{n}(B(o,r)) \subseteq \bigcup_{n=1}^\infty \delta_\lambda^{n}(U) $. Regarding contractibility, all homotopy groups are trivial, and one can use the Whitehead Theorem. 
\end{exercise}

 \begin{exercise}\label{connectivity_connected_balls_homotetic}
Let $(Y,d_Y)$ be a locally compact metric space and $y_0\in Y$. Assume that, for all $\lambda>0$, 
there is an isometry $f: (Y,d_Y) \to (Y,\lambda d_Y)$ with $f(y_0)=y_0$.
Then $Y$ is connected. In fact, every closed metric ball at $y_0$ is connected. 
\\
{\it Solution.}
We begin with a technical claim:
\begin{equation}\label{claim_connected_balls_homotetic}\forall \eps>0,\forall \bar y \in Y, \text{ with } \bar y\neq y_0 \qquad \exists y'\in Y \text{ such that } d(\bar y, y')<\eps \text{ and } d(y',y_0)<d(\bar y, y_0).
\end{equation}
For proving \eqref{claim_connected_balls_homotetic}, 
for each $\delta\in(0,1]$, we denote $\delta Y: =(Y,\delta d_Y)$ and consider the set
$$\Omega_\delta: =\left\{ f(\bar y)\;|\; f:\delta Y\to Y {\text{ isometry with }} f(y_0)=y_0\right\}.$$
We have $\Omega_\delta\subseteq\{y\in Y\;:\; d(y_0,y)=\delta d(y_0,\bar y)\}.$
The assumption that $(\delta Y,y_0)$ is isometric to $(Y,y_0)$ rephrases as $\Omega_\delta\neq \emptyset$.
The space is boundedly compact by Exercise~\ref{Ex_local_to_global_dilations}, so we use  Ascoli--Arzel\`a's Theorem, Exercise~\ref{ex_AA}.
For $\delta\in(1/2,1)$, pick $z_\delta\in \Omega_\delta$ and a respective $\delta$-homothety $f_\delta$. Using the Ascoli--Arzel\`a argument to the uniformly Lipschitz maps $f_\delta$, we have that there exists a sequence $\delta_n \nearrow 1$ for which the maps $f_{\delta_n}$ converge uniformly on compact sets to an isometry $f$. In particular,
$f(y_0)=y_0$ and $f_{\delta_n}(\bar y)\to f(\bar y)$, as $n \to \infty$. 
Set $y_n: =f^{-1}(f_{\delta_n}(\bar y))$. Observe that
$$d(\bar y, y_n)=d(f(\bar y), f_{\delta_n}(\bar y)) \to 0,\qquad \text{ as } n\to\infty,$$
and
$$
d(y_0,y_n)=d(f(y_0), f_{\delta_n}(\bar y))=d(y_0, f_{\delta_n}(\bar y))=d(f_{\delta_n}(y_0), f_{\delta_n}(\bar y))=\delta_n d(y_0, \bar y)<d(y_0, \bar y).$$
Thus, for $n$ large enough, one can take $y'$ as $y_n$ for obtaining \eqref{claim_connected_balls_homotetic}.

To deduce Exercise~\ref{connectivity_connected_balls_homotetic}, assume by contradiction that some closed ball $\bar B(y_0,r)$ is not connected.
Thus, there are  non-empty closed sets $K_1$, $K_2\subseteq \bar B(y_0,r)$ such that
$K_1\cap K_2=\emptyset$ and 
$K_1\cup K_2=\bar B(y_0,r)$.
Say $y_0\in K_1$.
Since $Y$ is boundedly compact (see Exercise~\ref{Ex_local_to_global_dilations}), both $K_1$ and $K_2$ are compact.
Consequently, first, the value $\eps: =d(K_1,K_2) $ is strictly positive. Second, there exists $\bar y\in K_2$ such that $d(\bar y, y_0)=d(K_2, y_0)$.
By \eqref{claim_connected_balls_homotetic}, there is some $y'\in Y$ such that $d(\bar y, y')<\eps$ and $d(y',y_0)<d(\bar y, y_0)$.
By the second inequality, we get that $y'\notin K_2$ and that $y'\in \bar B(y_0,r)$. 
By the first inequality, we get that $y'\notin K_1$.
We contradicted the fact that $K_1\cup K_2=\bar B(y_0, r)$. Hence, the set $\bar B(y_0,r)$ is connected. 
By Exercise~\ref{Ex_local_to_global_dilations}, the metric space $Y$ is connected. 
\end{exercise}





\chapter{Sub-Finsler Lie groups}\label{ch_SFLie}
This chapter is fundamental; it is the core of this book.
We consider sub-Finsler Lie groups and discuss their differential geometry and their metric geometry when equipped with Carnot-Carath\'eodory distances.
In Section~\ref{sec: 6.1}, we discuss the perspective of left-invariant sub-Finsler structures.
We see how, in sub-Finsler geometry, quotients can be seen as submetries, and hence, we can lift geodesics. 
We present Chow's theorem and prove a weakened version of the Ball-Box Theorem, which will be later refined in Theorem~\ref{thm66b22d85}.
In Section~\ref{sec: 6.2}, we discuss the endpoint map and its singular elements: the abnormal curves.
In Section~\ref{sec: 6.3}, in the setting of sub-Riemannian manifolds, we provide a first-order study of geodesics, leading to the so-called Pontryagin Maximum Principle.
In the final Section~\ref{sec: 6.4}, we show that geodesic homogeneous manifolds are Carnot-Carath\'eodory spaces.

\section{Left-invariant sub-Finsler structures on Lie groups}\label{sec: 6.1}\index{sub-Finsler! -- structure}\index{Lie group}
A standard assumption in the geometry of Lie groups is that the objects under consideration, such as distributions and sub-Finsler structures, are left-invariant. As a result, every considered distribution is a polarization, meaning that it has a constant rank. Similarly, as for Lie algebras of Lie groups, we will have two interpretations for polarizations and two for continuously varying norms because of the left-invariance.

\subsection{Left-invariant polarizations and horizontal curves}

In this section, we set-wise interpret the Lie algebra of each Lie group $G$ as the tangent space $T_{1}G$ at the identity element $1=1_G$.
	Let $G$ be a Lie group with Lie algebra $\g=\Lie(G)$. 
	Let $\Delta\subset TG$ be a distribution, as in Definition~\ref{def_distribution}. The distribution $\Delta$ is said to be {\em left-invariant} if\index{left-invariant! -- distribution}
	$$ \Delta_{gh}=(\dd L_g)_h \Delta_{h}, \qquad \forall g,h\in G,$$
	where, as in the previous chapter, we denoted by $L_g$ the left-translation by $g$.
	Each $(\dd L_g)_h$ is an isomorphism, and left-translations act transitively; thus, the rank of the subspaces $ \Delta_{g}$ is independent of $g\in G$. In other words, every {left-invariant} distribution is a polarization, following the terminology of Definition~\ref{def_distribution}. 
	
	Every {left-invariant} distribution $\Delta\subset TG$ determines a vector subspace $V: =\Delta_{1_G} \subseteq \g$. Vice versa, every vector subspace $V\subseteq \g$ of $\g$ determines a left-invariant distribution $\Delta$, by $ \Delta_{1_G}: =V$ and \begin{equation}\label{eq:induced:distribution}
	\Delta_g: =\left\{v \in T_gG: (\dd L_g)^{-1}v\in V\right\}, \qquad \forall g\in G.
	\end{equation}
		Observe, indeed, that the set $\Delta\subseteq TG$ is left-invariant and a polarization, whose rank equals $\dim(V)$, for which there is a global frame; see Exercise~\ref{ex_frame_left_inv_distribution}.\index{rank! of a polarization}
		In essence, there is a one-to-one correspondence between vector subspaces of $\Lie(G)$ and {left-invariant} distributions on $G$.

\begin{definition}[Polarized group]\index{polarized! -- Lie group}\index{polarized! -- group|see{polarized Lie group}}\index{polarization}\index{polarized! -- Lie algebra}
Given a Lie group $G$ and a vector subspace $V\subseteq \Lie(G)$, we say that the pair $(G, V)$ forms a
 {\em polarized Lie group}, or simply, a {\em polarized group}, and we refer to the distribution $\Delta$ defined in \eqref{eq:induced:distribution} as the {\em induced (left-invariant) distribution}. We can also refer to $V$ as a {\em polarization}, equating it with $\Delta$. We also say that the pair $( \Lie(G), V)$ is a {\em polarized Lie algebra}.
\end{definition}

Because every left-invariant distribution comes from a vector subspace of the Lie algebra, it is easy to verify whether the distribution is bracket generating. Indeed, we have two (equivalent) ways of calculating the iterated brackets and the flag of subbundles of Definition~\ref{def:Delta_k}. It should not be surprising that each left-invariant distribution is 
equiregular, according to Definition~\ref{def:equiregular}, because the flag of subbundles preserves the symmetry of being left-invariant.
	
	In the following, given two subsets $U, V$ of a Lie algebra $\g$, we use the notation
	\[
	[U, V]: =\mathrm{span}\left\{[u,v]:u\in U, v\in V\right\}\subseteq \g.
	\]
	Moreover, for every vector subspace $V\subseteq \g$, we iteratively define 
	\begin{equation}\label{LIdistr_flag}
	V^{[1]}: =V, \qquad V^{[k]}: =V^{[k-1]}+[V, V^{[k-1]}], \qquad \forall k=2,3, \ldots.
	\end{equation}
A first observation is that, for each $k\in \N$, the left-invariant distribution induced by $V^{[k]}$ is the $k$-th element $\Delta^{[k]}$ in the flag of subbundles associated with $\Delta$ as in Definition~\ref{def:Delta_k} (see Exercise~\ref{ex_LIdistr_flag}).
A second observation is that to determine the Lie algebra $\Lie(\Gamma(\Delta)) $ generated by the sections of $\Delta$, one can use left-invariant frames, and hence look at the Lie algebra generated by $V: =\Delta_{1_G}$ within $\g$. In other words, we have
$$(\Lie(\Gamma(\Delta)))_1 = \bigcup_{k\in\N} V^{[k]} \qquad \text{ and } \qquad
(\Lie(\Gamma(\Delta)))_g =(\dd L_g)_1(\Lie(\Gamma(\Delta)))_1 .$$
Moreover, we stress that the subspaces $V^{[k]}$ are nested and of integer dimension. 
Thus, the function $k\in\N \mapsto \dim(V^{[k]})$ is non-decreasing and it takes values in $\{1,2, \ldots,n\}$, where $n: =\dim (G)$.
Actually, unless $n\leq1$, we need to have $\dim(V)>1$ in order for $\Delta$ to be bracket generating.
 Thus, if $\Delta$ is not bracket generating then there exists $\bar k<n$ such that $V^{[\bar k]}=V^{[l]}$ for every $l\geq k$; see Exercise~\ref{ex dim stabilizes for V}.
We have proved the next result, which explains when a left-invariant distribution satisfies the bracket-generating condition, \eqref{Bracket: generating}.

\begin{proposition}[Criterion for bracket generation]\label{Prop:criterion_generating_LieGroup}\index{criterion for bracket generation}\index{bracket-generating}
If $(G, V)$ is a polarized Lie group of dimension $n\geq 2$, then we have the following dichotomy:
\begin{itemize}
	\item[(a)] either $V^{[n-1]}= \g$ and consequently the induced left-invariant distribution $\Delta$ is bracket generating with step less than $n$;
	\item[(b)] or $V^{[n-1]}\neq \g$, and in fact there exists a Lie subgroup $H<G$ with $\dim(H)<\dim( G)$ and the restriction $\Delta|_H$ is contained in $TH$ and is bracket generating $TH$. 
	Here $\Delta|_H: = \{v\in \Delta\;:\; \pi(v)\in H \}$, with $\pi:TG\to G$ the bundle projection. 
\end{itemize} 
Moreover, denoting by $\bar k$ the smallest integer for which 	$V^{[ \bar k]}=V^{[\bar k+1]}$, if 
$V^{[ \bar k]}\neq\g$, then $\bar k<n-1$.
\end{proposition}

We can further rephrase the notion of horizontal curve as from Definition~\ref{def horizontal curve}.
In a polarized group $(G, V)$ with induced distribution $\Delta$, an absolutely continuous curve $\gamma \colon I \rightarrow G$ defined on an interval $I$ is {$\Delta$-horizontal} if\index{horizontal! -- curve}\index{curve! horizontal --}
\begin{equation}\label{def_gamma_prime}\index{$\gamma'$}\index{$\dot\gamma$} 
 \gamma'(t): =(\dd L_{\gamma(t)})^{-1} \dot\gamma(t) \in V, \qquad \text{ for almost every } t\in I.
 \end{equation}
 Notice that we have just defined a $V$-valued curve $ \gamma': I\to V$ not to be confused with the derivative $ \dot\gamma: I\to TG$, which is $TG$-valued. However, the curve $ \gamma'$, together with the initial point, maintains the whole information about $\gamma$, as we next see.
 
 \begin{proposition}[Integration of the tangent vector]\label{prop_Integration_tangent_vector}\index{integration of tangent vector}\index{Carath\'eodory Theorem}\index{Theorem! Carath\'eodory --} 
 Let $G$ be a Lie group, let $[a,b]\subseteq\R$ be an interval, and let $u: [a,b]\to \g$ be integrable. Then, for every $p\in G$ there exists a unique absolutely continuous curve $\gamma \colon [a,b] \rightarrow G$ such that 
 $\gamma(a)=p$ and $u= \gamma'$, where the latter is defined in \eqref{def_gamma_prime}.
 \end{proposition}
 
 \begin{proof}
We consider the ODE
 \begin{equation}\label{lift_multiplicative_integral}
 \left\{ \begin{array}{ccl}
\dot\gamma(t)&=&(L_{\gamma(t)})_* u(t)\\
\gamma(a)&=&p.
\end{array} \right.\end{equation}
 The existence of an absolutely continuous solution to the ODE is a consequence of the general Carath\'eodory's theorem; see \cite[Theorem 3.4]{ORegan}. 
We stress that Lie groups equipped with Riemannian left-invariant metrics are complete; see Lemma~\ref{completeness_homog}. Consequently, short-time solutions stay in compact sets. Hence, we have global solutions; see \cite[page 43]{Coddington_Levinson}. 
The uniqueness can be shown proving that, if $\gamma_1(t)$ and $ \gamma_2(t)$ are two solutions, then 
$$\dfrac{\dd}{\dd t}\left(\gamma_1(t)\gamma_2(t)^{-1}\right)\equiv 0.$$
This last fact can be easily shown using Exercise~\ref{ex:derivative:product: curves} and Exercise~\ref{ex:derivative:inverse: curves}. 
\end{proof}

\begin{definition}[Development of a curve]\label{def development}
Let $T\geq0$ and $\gamma: [0, T]\to G $ be an absolutely continuous curve in a Lie group $G$.
The {\em development} of $\gamma$ is the curve $\sigma: [0, T]\to \g$ defined as
$$\sigma(t): = \int_0^t \gamma'(s)\dd s, \qquad \forall t\in [0, T],$$
where $\gamma' :[0,T] \to  \g$ is as in \eqref{def_gamma_prime}.\index{development}
\end{definition}
The above curve $\sigma$ is valued in the vector space $\g$, where we identify $\sigma'$ with $\dot \sigma$. Thus, we have the identity 
$$\sigma'=\gamma'.$$ 
We stress that if $\gamma$ is $\Delta$-horizontal, for some left-invariant polarization $\Delta $, then $\sigma$ is $\Delta_1$-valued.

We also recall that by Proposition~\ref{prop_Integration_tangent_vector}, every absolutely continuous curve $\sigma: [0, T]\to \g$ is the development of some 
absolutely continuous curve $\gamma: [0, T]\to G $. We shall refer to this curve as the multiplicative integral.
\begin{definition}[Multiplicative integral]\label{def multiplicative integral}
Let $G$ be a Lie group with Lie algebra $\g$, and $T>0$.
Let $\sigma: [0, T]\to \g$ be an absolutely continuous curve.
The absolutely continuous curve 
$\gamma: [0, T]\to G $ with $\gamma(0)=1_G$ and $\sigma'=\gamma'$, which exists and is unique by Proposition~\ref{prop_Integration_tangent_vector}, is called the {\em multiplicative integral} of $\sigma$.\index{multiplicative! -- integral}
\end{definition}
Up to only considering curves starting at $1_G$ or $0_\g$, respectively, the multiplicative integral is the inverse operation of the development of a curve, as of Definition~\ref{def development}.

\subsection{Left-invariant norms and distances on Lie groups}\label{sec sub-Finsler Lie groups}

The first aim of this subsection is to clarify that left-invariant continuously varying norms on $TG$ are in one-to-one correspondence with symmetric norms on $T_{1_G}G$. 
A function $N: TG \to \R$ on the tangent bundle of a Lie group $G$ is said to be {\em left-invariant} if $N\circ \dd L_g = N$ for all $g\in G$.\index{left-invariant! -- function on the tangent bundle}
In fact, via left translations, the tangent bundle $TG$ is trivializable as $G\times \g$ and, in the left-invariant case, every continuously varying norm $N: TG\simeq G\times \g\to \R $, $(g,v)\mapsto N(g,v)$, will be constant with respect to the point $g$ and will be a norm in the vector $v$. 


One can check that if a function $\|\cdot \|: TG \to \R$ is left-invariant and its restriction to $T_{1_G}G$ is a symmetric norm, then $\|\cdot \|$ is a continuously varying norm, in the sense of Definition~\ref{def: continuously varying norm}.
Moreover, every symmetric norm on $T_{1_G}G$ is the restriction to $T_{1_G}G$ of a unique
 left-invariant continuously varying norm $\|\cdot \|$.
 Indeed, if $\|\cdot \|_{{1_G}}$ is a symmetric norm on $\mathfrak{g}$, then 
\begin{equation} \label{cont_var_norm}
 \|v\|: = \|(\dd L_{g^{-1}})_gv \|_{{1_G}}, \quad \forall g \in G, \; \forall v \in T_gG,
\end{equation}
defines a left-invariant continuously varying norm (Exercise~\ref{ex_LI_norm}).

\begin{definition}[Sub-Finsler Lie group]\label{def_Sub-Finsler_Lie_group}\index{sub-Finsler! -- Lie group}\index{Lie group! sub-Finsler --}
A {\em sub-Finsler Lie group} is a triple $(G, V, \| \cdot \|)$ where $G$ is a Lie group, $V$ is a bracket-generating subspace of $T_{1_G}G$,
and $\|\cdot \|$ is a symmetric norm on $V$. 
Every sub-Finsler Lie group is naturally seen as 
a Carnot-Carathéodory space where the distribution $\Delta$ is the induced distribution of the polarized group $(G, V)$ as in \eqref{eq:induced:distribution} and $\|\cdot \|$ is extended on $TG$ by first arbitrarily extending $\|\cdot \|$ as a norm on $T_{1_G}G$ and then as a left-invariant continuously varying norm by \eqref{cont_var_norm}. 
 \end{definition}

In the definition of sub-Finsler Lie group, the left-invariant continuously varying norm restricts to $\Delta$ as a function $\Delta\to \R$. 
This function is the left-invariant extension of the norm on $V$.
We still call it {\em continuously varying norm}.\index{left-invariant! -- continuously varying norm}  
 
 Every sub-Finsler Lie group has an associated sub-Finsler metric, as in \eqref{dist_CC}, which can be formulated using the above double viewpoint for left-invariant structures. Moreover, we recall Proposition~\ref{prop: energy} about the energy of curves.
Thus, the Carnot-Carathéodory distance $\dcc $ between two points $p,q \in G$ is\index{Carnot-Carath\'eodory! -- distance}\index{distance! Carnot-Carath\'eodory --}
\begin{eqnarray*}\label{def:formula:dcc:groups}
 \dcc (p,q) &: =& \inf \big\{\Length_{\|\cdot \|}(\gamma) \, \colon \, \gamma \; \Delta\text{-horizontal curve from } p \text{ to } q \big\}\\
 &=& \inf \left\{\int \| {\gamma}'\|_{1_G} \, \colon \, \gamma \; \text{AC curve from } p \text{ to } q, \text{ with } \gamma' \in V \right\} \\
 &=& \inf \left\{\sqrt{ \int_0^1 \| {\gamma}'\|_{1_G}^2}\, \, \colon \, \gamma \; \text{AC curve on } [0, 1] \text{ from } p \text{ to } q, \text{ with } \gamma' \in V \right\}.\\
\end{eqnarray*}
In Definition~\ref{def_Sub-Finsler_Lie_group}, we made the choice of assuming that the polarization is bracket generating. 
Of course, one could also consider sub-Finsler metrics associated with non-bracket-generating polarizations. However, because of Proposition~\ref{Prop:criterion_generating_LieGroup}, if the polarization is not bracket generating, then one can just restrict to the Lie subgroup that it generates, as by Theorem~\ref{teo1145}.

The following are some basic metric properties of sub-Finsler Lie groups when they are seen as metric spaces with their Carnot-Carathéodory metrics.
\begin{theorem}\label{Properties_SFG}\index{complete}\index{geodesic! -- metric}\index{boundedly compact}\index{isometrically homogeneous}
Every sub-Finsler Lie group is a metric space that is
	\begin{description}
	\item[\ref{Properties_SFG}.i.] complete,
	\item[\ref{Properties_SFG}.ii.] geodesic, 
	\item[\ref{Properties_SFG}.iii.] boundedly compact,
	\item[\ref{Properties_SFG}.iv] isometrically homogeneous: the distance is left-invariant.
	\end{description}
\end{theorem}

\begin{proof} Let $(G, d)$ be a sub-Finsler Lie group equipped with its CC distance.
We begin by observing that $(G, d)$ is isometrically homogeneous since the group of left translations acts transitively and by maps that preserve the family of horizontal curves and their length.

Then, since the distance induces the manifold topology, there exists $r>0$ such that the closed metric ball $\bar B(1_G,r)$ is compact. By homogeneity, the metric space is locally compact and complete.

We conclude by invoking Proposition~\ref{prop:length:geod} or Theorem~\ref{HRCV}, or Theorem~\ref{HRCV_CC}.
\end{proof}

\begin{proposition}\label{biLipschitz_equivalence_CCnorms}\index{bi-Lipschitz! -- equivalent distances}
	If $(G, V)$ is a polarized Lie group,	then every two left-invariant CC distances induced by norms on $V$ are bi-Lipschitz equivalent, globally.
\end{proposition}
\begin{proof}
	The notion of length of a horizontal curve $\gamma$ (and hence the notion of the associated CC distance) depends on the norm ${\|\cdot \|}$ in the following way:
	$ \Length_{\|\cdot \|}(\gamma)= \int \| {\gamma}'\|_{1_G}$.
	Since $V$ is finite-dimensional, every choice of $\|\cdot\|_{1_G}$ is biLipschitz equivalent to every other. This produces a biLipschitz equivalence for CC distances.
\end{proof}
Because of Proposition~\ref{biLipschitz_equivalence_CCnorms}, whenever our interest is in metric spaces up to biLipschitz equivalence, we may assume that the norm $ \norm{\cdot} $
is coming from a scalar product.

\subsection{Quotients and submetries between sub-Finsler Lie groups}\index{quotient! -- of sub-Finsler Lie groups}\index{quotient! -- space}\index{submetry}


 \begin{proposition}\label{prop_subFin_submetry_Lie} 
Let $(G, \Delta, \norm{\cdot}) $ be a sub-Finsler Lie group and $H<G$ a closed subgroup.
Then, the distance function from \eqref{distance_quotient} on the quotient manifold $\mfaktor{H}{G}: =\{Hg: g\in G\}$, which makes the projection $\pi: G\to \mfaktor{H}{G}$ a submetry, is a sub-Finsler distance, where the distribution $\Delta$ on $ \mfaktor{H}{G}$  is the only one such that
 $$\Delta_p =\dd\pi(\Delta_{\tilde p} ), \qquad \forall p\in \mfaktor{H}{G},  \forall \tilde p \in \pi^{-1}(p),$$ with norm
\begin{eqnarray}
 \norm{v} &: =&\inf\{||w||: \tilde p\in \pi^{-1}(p), w\in T_{\tilde p}\Delta, \dd\pi (w)=v \}, \qquad \forall v\in 
 \Delta _p,
 \nonumber
 \\
 &=&\inf\{||w||: w\in T_{\tilde p}\Delta, \dd\pi (w)=v \}, \qquad \forall \tilde p \in \pi^{-1}(p), \forall v\in 
 \Delta _p .
 \label{d defnormadown2}
\end{eqnarray}
\end{proposition}
 \begin{proof} 
We begin by stressing that the distribution is well defined, and the definition of the norm has the property \eqref{d defnormadown2} because the group $H$ acts by isometries and transitively on fibers. Moreover, by construction we have that the map $ (\dd \pi)_{\tilde p}: (\Delta_{\tilde p}, \norm{\cdot} ) \to (\Delta_{\pi(\tilde p)}, \norm{\cdot} ), $ for $ \tilde p\in G$, is a submetry. 
 
We denote by $\dcc $ the sub-Finsler metric. 
Let $d$ be the distance function on $\mfaktor{H}{G}$ given by \eqref{distance_quotient}, so that $\pi: (G,\dcc) \to \left(\mfaktor{H}{G},d\right)$ is a submetry, by Proposition~\ref{prop:distance_quotient}.
The plan is to use Proposition~\ref{prop_lifting_gives_submetry} to show that also the map $\pi: (G,\dcc ) \to \left(\mfaktor{H}{G},\dcc \right)$ is a submetry. Since there is a unique distance function on the target for which a map is a submetry, we would deduce that the distance functions $d$ and $\dcc $ on $\mfaktor{H}{G}$ coincide.

For the argument, we will make use of a measurable selection theorem,  which we review at the end of this proof; see Theorem~\ref{thm_selection}. 
 Given $p\in \mfaktor{H}{G}$, if we take two preimages $\tilde p_1 $ and $\tilde p_2 \in \pi^{-1}(p)$, then the maps $\pi\circ L_{\tilde p_1}$ and $\pi\circ L_{\tilde p_2}$ coincide, and we denote such a map by $\hat \pi_p$.
 Consider the multifunction $\Psi$ taking each element $w \in \Delta_p$, with $p\in \mfaktor{H}{G}$, to the set
 \[\Psi(w) \coloneqq \{ v \in \Delta_{1_G} \colon (\hat \pi_p)_*(v)=w, \; \norma{v} = \norma{w} \}.\] Note that for every $w\in \Delta_p$, the set $\Psi(w)$ is a closed subset of $\Delta_{1_G}$,
 which is nonempty because of the definition of the norm on $ \mfaktor{H}{G}$. Let $U $ be an open subset of $\Delta_{1_G}$. We have
 \begin{align*}
 & \hspace{-1.5cm}\left\{ w\in \Delta\subseteq T({\mfaktor{H}{G}}) \, \colon\,   \Psi(w) \cap U \neq \emptyset \right\} 
 \\&= \left\{ w \, \colon\, \{ v \in \Delta_{1} \, \colon\,  (\hat \pi_p)_*(v)=w\in \Delta_p, p\in \mfaktor{H}{G},\; \norma{v} = \norma{w}\} \cap U \neq \emptyset \right\} \\
 &= \left\{ (\hat \pi_p)_*(v) \in \Delta_p \, \colon\, v \in U, p\in \mfaktor{H}{G}, \; \norma{v} = \norma{(\hat \pi_p)_*(v)} \right\} \\
 &= \pi_*\left(\left\{ X_{\tilde p} \in \Delta^G \, \colon\, (L_{\tilde p} )^* X_{\tilde p}\in U, \,
 \norma{X_{\tilde p}} = \norma{\dd \pi X_{\tilde p}}
 \right\}\right) .
 \end{align*}
 Notice that the set $\left\{ X_{\tilde p} \in \Delta \, \colon\, (L_{\tilde p} )^* X_{\tilde p}\in U 
 \right\}$ is open in $\Delta$ and $\left\{ X_{\tilde p} \in \Delta \, \colon\, \norma{X_{\tilde p}} = \norma{\dd \pi X_{\tilde p}}
 \right\}$ is closed. Hence, they and their intersection are a countable union of compact sets.
 The image under $\pi_*$, which is a continuous map, is again a countable union of compact sets. Hence, it is a Borel set. Thus we can apply Theorem~\ref{thm_selection}: there exists a Borel measurable selection $\psi: \Delta\subseteq T \left(\mfaktor{H}{G} \right) \rightarrow \Delta_{1_G}$ of $\Psi$.
 We observe that the constructed selection map $\psi$ has the following properties:
 \begin{equation}\label{16may2024}
 (\hat \pi_p)_*(\psi(w))=w \quad \text{ and } \quad \norma{\psi(w)} = \norma{w}, \qquad \forall p\in \mfaktor{H}{G}, \forall w \in \Delta_p.
 \end{equation}
 
 To use Proposition~\ref{prop_lifting_gives_submetry}, we consider  a horizontal curve $\gamma: [0, T] \to \mfaktor{H}{G}$ and $\tilde g \in \pi^{-1}(\gamma(0))$ and we want to construct a curve $\gamma$ with the properties \ref{prop_lifting_gives_submetry}.2.i-iii.
 Notice that $\psi \circ \dot \gamma: [0, T] \to \Delta_{1_G}$ is measurable.
 We consider the curve $\tilde\gamma \colon [0, T] \rightarrow G$ solution of 
 \begin{equation*}
 \begin{cases}
 \tilde\gamma(0) = \tilde g\\
 \dot{\tilde\gamma}(t) = (\dd L_{\tilde\gamma(t)})\psi ( \dot \gamma(t)),
 \end{cases}
\end{equation*}
which exists by Proposition~\ref{prop_Integration_tangent_vector}.
The curve $\gamma$ satisfies \ref{prop_lifting_gives_submetry}.2.i by construction. Regarding \ref{prop_lifting_gives_submetry}.2.ii, the two curves $\pi\circ\tilde\gamma$ and $ \gamma$ are equal because they have the same initial point and same derivatives:
 $$\frac{\dd}{\dd t} (\pi\circ\tilde\gamma)
 = \pi_*(\dd L_{\tilde\gamma(t)})\psi ( \dot \gamma)
 = (\pi_{\gamma(t)})_*\psi ( \dot \gamma)\stackrel{\eqref{16may2024}}{=}\dot\gamma.$$
 Regarding \ref{prop_lifting_gives_submetry}.2.iii, the two curves have the same length because they have the same speed: $$\norm{ \dot{\tilde\gamma}} = 
 \norm{ (\dd L_{\tilde\gamma(t)})\psi ( \dot \gamma)}
 =\norm{ \psi ( \dot \gamma)}
 \stackrel{\eqref{16may2024}}{=} \norm{\dot\gamma}.$$

Finally, Condition \ref{prop_lifting_gives_submetry}.1 is satisfied since $\pi$ is 1-Lipschitz; see Exercise~\ref{ex:Lip_for_CCmfds}.
	By Proposition~\ref{prop_lifting_gives_submetry}, we deduce that the two distances $d$ and $\dcc$ on $\mfaktor{H}{G}$ are the same. 
 \end{proof}

In the above proof, we used the following general result by Kuratowski and Ryll-Nardzewski: 
\begin{theorem}[{\cite{KRN}}]\label{thm_selection}\index{Kuratowski and Ryll-Nardzewski Theorem}\index{measurable selection theorem}\index{selection}
 Let $X$ be a separable completely metrizable topological space, let $\mathcal{B}(X)$ be the Borel $\sigma$-algebra of $X$, let $(\Omega, \mathcal{F})$ be a measurable space, and $\Psi$ a multifunction on $\Omega$ taking values in the set of nonempty closed subsets of $X$. If for every open subset $U$ of $X$ we have 
 \[\{\omega \colon \Psi(\omega) \cap U \neq \emptyset\} \in \mathcal{F}, \]
 then $\Psi$ has a selection that is $\mathcal{F}$-$\mathcal{B}(X)$-measurable.
\end{theorem}

\begin{remark}
When the distribution $\Delta$ of the sub-Finsler group $G$ is such that $\Delta_{1_G}\oplus\Lie(H)=\Lie(G)$, then the existence of the selection in Theorem~\ref{thm_selection} is a triviality. 
In fact, the map $ \dd \pi |_{\Delta_{p}}: \Delta_p \to \Delta_{\pi(p)} $ is a linear isomorphism, hence we take $\psi: = (\dd\pi |_{\Delta })^{-1}$.
Moreover, in this case, the lifted curves $\tilde\gamma$ are unique, given the initial point.
\end{remark}

\begin{remark}
 Let $H$ be a sub-Finsler Lie group and $G$ a Lie group.
 Suppose that $\pi: G\to H$ is a surjective Lie homomorphism.
 Then, there is a subFinsler structure on $G$ such that $\pi$ is a submetry.
 In fact, if we take a bracket-generating polarization $V\subseteq T_1 G$ such that $\pi_*V = \Delta^H_1$, e.g., $\pi_*^{-1}\Delta^H_1$, then there is a norm on $V$ such that $\pi$ is a submetry. The unit ball for such a norm can be taken as the intersection of the preimage of the unit ball on $\Delta_1^H$ and a strip that is transverse to the kernel of $\pi$.
\end{remark}



%
\begin{remark} 
Let $\pi: G\to H$ be a submetry and a Lie homomorphism between sub-Finsler Lie groups.
Thus, by Proposition~\ref{prop: same regularity2}, geodesics in $H$ can be lifted to geodesics in $G$. In particular, if all geodesics on $G$ are smooth, then so are the geodesics in $H$.
\end{remark}

\subsubsection{Open questions}\index{questions}
The following are some unsolved questions on sub-Finsler Lie groups.
\begin{question}
Let $(G, \Delta, \norm{\cdot})$ be a sub-Finsler Lie group.
If the unit ball of $(\Delta, \norm{\cdot}) $ is strictly convex, then does there exist, for every $p$ and $q\in G$, a smooth geodesic from $p$ to $q$?\end{question}
\begin{question}
Let $(G, \Delta, \norm{\cdot})$ be a sub-Finsler Lie group.
If the unit ball of $(\Delta, \norm{\cdot}) $ is polyhedral (i.e., it is the convex hull of finitely many points), then does there exist $K\in \N$ such that between every pair of points in $G$ there exists a geodesic whose control is piecewise constant with at most $K$ pieces?
\end{question}

\begin{question}
Let $G$ be a Lie group.
If $\rho$ is a left-invariant metric that is boundedly compact, locally bounded, and quasi-geodesic, is $\rho$ at a bounded distance from a sub-Finsler left-invariant metric?
\end{question}

\subsection{A direct, effective proof of Chow's theorem}\label{effective_Chow_Carnot}\index{Chow Theorem}\index{Theorem! Chow --}
In this section, we will give an explicit construction of a horizontal path connecting an arbitrary point $p$ in a bracket-generating polarized group $(G, V)$ to the identity element $1_G$. Moreover, when a norm is fixed on an $s$-step polarization, the path will have a length bounded by the $1/s$-power of a Riemannian distance between $p$ and $1_G$. This is the most important consequence of the Ball-Box Theorem: \eqref{eq_BB_thm}.

\subsubsection{Brackets as products of exponentials}
The philosophy behind the following discussion is that to go in a direction given as a bracket of two vector fields $X$ and $Y$ on a manifold $M$, one can go along a not-necessarily-closed quadrilateral constructed using the flows of the two vector fields. We will give a generalization of the following formula as in \ref{def:Lie_bracket_vector_fields}.d:
 $$
	[X, Y]_p=\left.\frac{1}{2}\frac{\dd^2}{\dd t^2}\left(\Phi^{-t}_Y\circ\Phi^{-t}_X\circ \Phi^t_Y \circ \Phi^t_X\right)(p)\right|_{t=0}, \qquad \forall p\in M
. $$ 
In the above formula, we denote by $\Phi^t_X$ the flow of the vector field $X$ at time $t$. 
By Corollary~\ref{Warner3.31}, for left-invariant vector fields in a Lie group $G$, we have
$$\Phi^t_X(p)=pe^{tX}, \qquad p\in G, \forall t\in \R, \forall X\in \g,$$
where, for the sake of clarity, we also use the notation $e ^ X: = \exp(X)$ for $X\in \g$.
Hence, for all $X$ and $Y$ in the Lie algebra $\g$ of the Lie group $G$, we have
\begin{equation}\label{eq_Bracket_Lie_groups}
	[X, Y]_{1_G} = \left.\frac{1}{2}\frac{\dd^2}{\dd t^2}
	e^{tX} e^{tY} e^{-tX}e^{-tY} 
	\right|_{t=0}, \qquad \forall X, Y\in \g.
\end{equation}

\begin{definition}[Maps $P_t$]\label{definition_pt_2024}\index{$P_t$}
Given a Lie group $G$, with Lie algebra $\g$, $X, Y\in \g$, and $t\in \R$, we define
$$P_t(X): =e^{tX}\quad \text{and }\quad P_t(X, Y): =e^{tX} e^{tY} e^{-tX}e^{-tY}.$$
By recurrence, 
 for $k\geq 2$ and $X_1, \ldots, X_{k+1} \in \g$, and $t\in \R$, we also define
$$P_t(X_1, \ldots, X_{k+1}): = P_t(X_1, \ldots, X_k) e^{tX_{k+1}}(P_t(X_1, \ldots, X_k))^{-1}e^{-tX_{k+1}}.$$
\end{definition}
\begin{remark}
As an immediate consequence of \eqref{eq_Bracket_Lie_groups} (or of BCH \eqref{expansion of BCH}), we have
$$P_t(X, Y)=\exp\left( {t^2[X, Y]+o(t^2)}\right), \qquad \text{ as } t\to0.$$
More generally, if we consider $X_1, \ldots, X_{d} \in \g$, with $d\in \N$, then we claim that 
%
\begin{equation}\label{expan:pt}
 P_t(X_1, \ldots, X_d)=\exp\left( {t^d [ X_1, \ldots, X_{d}] +o(t^d)}\right), \qquad \text{ as } t\to0,
\end{equation}
where we use the notation
$$[ X_1, X_2, X_3, \ldots, X_{d}]: = [\ldots[[X_1, X_2], X_3], \ldots, X_{d}] .$$
Indeed, we proceed by induction. The statement is true for $d\in \{1,2\}$. Assume it is true for an arbitrary $d$. Call $\omega(t)$ the $o(t^d)$ function such that
$P_t(X_1, \ldots, X_d)=e^{t^d [\cdots[[X_1, X_2], X_3], \ldots, X_{d}] +\omega(t)}.$
Then we have, by the BCH Formula \eqref{expansion of BCH},
\begin{eqnarray*}
 P_t(X_1, \ldots, X_{d+1})&= &P_t(X_1, \ldots, X_d) e^{tX_{d+1}}\left(P_t(X_1, \ldots, X_d)\right)^{-1} e^{-tX_{d+1}}\\
&= &e^{t^d [\ldots[X_1, X_2], \ldots, X_{d}] +\omega(t)} e^{tX_{d+1}} \left(e^{t^d [\ldots[X_2, X_1], \ldots, X_{d}] +\omega(t)}\right)^{-1} e^{-tX_{d+1}}\\
&= &e^{t^d [\ldots[X_1, X_2], \ldots, X_{d}] +\omega(t)} e^{tX_{d+1}} e^{-t^d [\ldots[X_2, X_1], \ldots, X_{d}] -\omega(t)} e^{-tX_{d+1}}\\
&\stackrel{\eqref{expansion of BCH}}{=} &e^{\left(tX_{d+1}+t^d [\ldots[X_1, X_2], \ldots, X_{d}] +\omega(t)+ \frac{1}{2} t^{d+1} [\ldots[X_1, X_2], \ldots, X_{d+1}]+o(t^{d+1})\right)} \\
&&\qquad\quad
\cdot\, e^{\left(-tX_{d+1}-t^d [\ldots[X_1, X_2], \ldots, X_{d}] -\omega(t)+ \frac{1}{2} t^{d+1} [\ldots[X_1, X_2], \ldots, X_{d+1}]+o(t^{d+1})\right)}\\
&\stackrel{\eqref{expansion of BCH}}{=} &e^{t^{d+1} [\ldots[[X_1, X_2], X_3], \ldots, X_{d+1}] +o(t^{d+1})}.
\end{eqnarray*}
\end{remark}
%

Each $P_t$ is, in fact, a product of elements of the form $ e^{\pm tX_i}$.
We now form a map that will help in constructing horizontal paths.

\begin{definition}\label{def_E_for_BB}\index{$E$}
Let $\{X_{j,k}\}_{j,k}\subset \g$, with $j\in \{1, \ldots, n\}$, $k\in \{1, \ldots, d_j\}$, and $ d_1, \ldots,  d_n\in \N$.
For each $j$, we consider the expression 
$$P^{(j)}(t): = P_t(X_{j, 1}, \ldots, X_{j,d_j}), \quad\forall t\in \R.$$
We finally define the map $E:\R^n \to G$ associated with $\{X_{j,k}\}_{j,k}$ as
$$E({\bf{t}}): =\prod_{j=1}^n 
P^{(j)}(t_j^{1/d_j})
: =P^{(1)}(\sqrt[d_1]{t_1})\cdots P^{(n)}(\sqrt[d_n]{t_n}),
\quad\forall \bf{t}\in \R^n
.$$
\end{definition}
We have used the notation $ t^\alpha ={\rm sgn}(t)|t|^\alpha$, for $t, \alpha\in \R$; so for example we have $\sqrt{-4}=-2$.

Such a map $E$ satisfies the following properties.
\begin{proposition}\label{prop66b0c659}
 Let $G$ be a Lie group of dimension $n\in \N$.
 Let $\{X_{j,k}| j\in\{1, \dots,n\}, \ k\in\{1, \dots,d_j\}\}$ with $d_j\in\N$ be a family of vectors in $\Lie(G)$.
 Assume that
 $([X_{j, 1}, X_{j,2}, \dots, X_{j,d_j}])_{j=1}^n$ span $\Lie(G)$.

 \begin{description}
 \item[\ref{prop66b0c659}.i.]
 The map $E:\R^n\to G$ associated with $(X_{j,k})_{j,k}$ as in Definition~\ref{def_E_for_BB} is a local $C^1$-diffeomorphism around $0$.
 \item[\ref{prop66b0c659}.ii.]
 If, in addition, the vectors $(X_{j,k})_{j,k}$ are of unit length with respect to a sub-Finsler structure on $ G$, then there exists $C$ such that, for all $\mathbf t\in\R^n$,
 \begin{equation}\label{eq66b0c509}
 \dcc (1_G,E(\mathbf t) ) \le C \sum_{j=1} |t_j|^{1/d_j} .
 \end{equation}
 \end{description}
\end{proposition}

\proof 
Regarding \ref{prop66b0c659}.i, we show that $(\dd E)_0$ is non-singular.
Let $X_j: =
[[\ldots[[X_{j, 1}, X_{j,2}], X_{j,3}], \ldots], X_{j,d_j}]$, for $j\in\{1, \ldots\,n\}$.
 From how $E$ has been defined and from \eqref{expan:pt}, we have
\begin{eqnarray*}
(\dd E)_0 \partial_j&=&\left.\dfrac{\dd}{\dd t_j}E({\bf{t}})\right|_{{\bf{t}}=0}\\
&\stackrel{\rm def}{=}&\left.\dfrac{\dd}{\dd t_j}P^{(j)}(\sqrt[d_j]{t_j})\right|_{t_j=0}\\
&\stackrel{\rm def}{=}&\left.\dfrac{\dd}{\dd t}P_{\sqrt[d_j]{t}}(X_{j, 1}, \ldots, X_{j,d_j})\right|_{t=0}\\
&\stackrel{\eqref{expan:pt}}{=}&\left.\dfrac{\dd}{\dd t}e^{t [\ldots[[X_{j, 1}, X_{j,2}], X_{j,3}], \ldots, X_{j,d_j}] +o(t )}\right|_{t=0}\\
&=&\left.\dfrac{\dd}{\dd t}e^{t X_j +o(t )}\right|_{t=0}\\
&=&X_j.
\end{eqnarray*}
In other words, the linear map $(\dd E)_0$ sends the basis $\partial_1, \ldots, \partial_n$ to the basis $X_1, \ldots, X_n$. 
We conclude by the Inverse Function Theorem.

Regarding \ref{prop66b0c659}.ii, recall, by Corollary~\ref{Warner3.31}, that the flow lines of a left-invariant vector field $X$ are the curves $ge^{tX}$, fixed $g\in G$ and varying $t\in\R$.
Now, since $P_t$ is a product of exponentials, then so is $E$. More explicitly, we have
$$E({\bf{t}})=\exp({\epsilon_1 t_{c_1}^{ a_1} X_{ b_1}})\cdots \exp({\epsilon_N t_{c_N}^{ a_N} X_{ b_N}}), \qquad \forall {\bf{t}}\in\R^n,
$$
for some $N\in\N$ and $\epsilon_i\in\{1,-1\}$, $ a_i^{-1}\in\N$, $ b_i\in
\{ (j,k): j\in \{1, \ldots\,n\}, k \in \{1, \ldots\,d_j\}\}
$, and $c_i\in\{1, \ldots\,n\}$, for each $i \in\{1, \ldots\,N\}$.
Notice that, if $c_i = j$, then $a_i = 1/d_j$.
Now it is enough to observe that, fixed $k$, the point
$$g: =\exp({\epsilon_1 t_{c_1}^{ a_1} X_{ b_1}})\cdots \exp({\epsilon_{k-1} t_{c_{k-1}}^{ a_{k-1}} X_{ b_{k-1}}})$$
can be connected to the point
$$\exp({\epsilon_1 t_{c_1}^{ a_1} X_{ b_1}})\cdots \exp({\epsilon_{k-1} t_{c_{k-1}}^{ a_{k-1}} X_{ b_{k-1}}})\exp({\epsilon_{k} t_{c_{k}}^{ a_{k}} X_{ b_{k}}})$$
by the path
$$g\exp({\epsilon_{k} s X_{ b_{k}}}), \quad \text{ for } s\in[0, |t_{c_{k}}^{ a_{k}}|],$$
which is tangent to $\pm X_{ b_{k}}$, thus horizontal
with length $|t_{c_{k}}|^{ a_{k}} = |t_{c_{k}}|^{1/ d_{c_{k}}}$.
Thus the curve obtained by joining these paths connects $1_G$ to $E(\mathbf t)$
with length at most $N\sum_{j=1}^n|t_j|^{1/d_j} $.
%
%
\qed

As a consequence of Proposition~\ref{prop66b0c659}, we obtain a quantitative bound for the Carnot-Carath\'eodory distance, which implies the H\"older equivalence between CC and Euclidean metrics; see \eqref{eq_BB_thm}. 
In Chapter \ref{ch_CCspaces}, we explained this weaker bound as a consequence of the more precise bound given by the Ball-Box Theorem~\ref{Ball-Box}. Here, we do not prove that theorem in the setting of manifolds. 
However, first with the previous Proposition~\ref{prop66b0c659}, we soon deduce 
 \eqref{eq_BB_thm} for Lie groups. Later, in Theorems~\ref{Ball-Box4Carnot} and~\ref{thm66b22d85}. we will prove the Ball-Box Theorem for Carnot groups and for sub-Finsler Lie groups, respectively. 
\begin{corollary}[Weak Ball-Box Theorem]\label{corol_weak_BB}\index{weak Ball-Box Theorem}\index{Theorem! weak Ball-Box}\index{H\"older equivalence}
 Let $ G$ be a sub-Finsler Lie group, equipped with the associated Carnot-Carath\'eodory metric $\dcc $ and a Riemannian distance $d_R$.
 Then there exists a neighborhood $U$ of $1_G$ and $C>0$ such that
 \begin{equation}\label{eq66b0b89a}
 \frac1C d_R \le \dcc \le C d_R^{1/s}
 \qquad\text{on $U$},
 \end{equation}
 where $s$ is the step of the sub-Finsler structure.
\end{corollary}
\begin{proof}
We may assume that $d_R $ is left invariant. 
 The first inequality in~\eqref{eq66b0b89a} is a direct consequence of the definition of CC distance; see~\eqref{fin<subfin}.
 For the second inequality,
since the subFinsler structure is bracket generating, then there are vectors 
 $([X_{j, 1}, X_{j,2}, \dots, X_{j,d_j}])_{j=1}^n$ that span $\Lie(G)$ and are such that each $X_{j,k}$ is horizontal and of unit length.
 Then, we use the associated map $E$ from Proposition~\ref{prop66b0c659}.ii.
 Firstly,
 since
by Proposition~\ref{prop66b0c659}.i,
 there are $R,C_1>0$ such that
 \[
 \|\mathbf t\| \le C_1 d_R(1_G,E(\mathbf t)), \qquad \forall 
 \mathbf t\in\R^n, \, \text{ with } \|\mathbf t\|<R.
 \]
 Secondly,
 by~\eqref{eq66b0c509}, there is $C_2>0$ such that, if $\|\mathbf t\| < R$, then
 \[
 \dcc (1_G,E(\mathbf t)) \le C_2 \|\mathbf t\|^{1/s} .
 \] 
 Let $U': = E(\{\mathbf t\in \R^n:\|\mathbf t\|<R \})\subset G$ which is a neighborhood of $1_G= E(\mathbf 0)$.
 If $p\in U'$, then there is $\mathbf t\in\R^n$ with $\|\mathbf t\|<R$ and $E(\mathbf t) = p$.
 So,
 \[
 \dcc (1_G,p) \le C_1 \|\mathbf t\|^{1/s}
 \le C_1 C_2^{1/s} d_R(1_G,p)^{1/s} .
 \]
 Finally, we take a neighborhood $U$ of $1_G$ with $U^2 \subset U'$,
 to obtain~\eqref{eq66b0b89a} by left invariance of both $\dcc $ and $d_R$.
\end{proof}

\section{Endpoint map on polarized groups}\label{sec: 6.2}
In this section, we begin with a parametrization of the space of those horizontal curves in a polarized group that start from the identity element.
The parameterization of this space of curves leads to a Hilbert-space structure, providing a function analytic framework for investigating Carnot-Carath\'eodory spaces.

\subsection{Endpoint map}
 Let $(G, V)$ be a polarized group. After fixing a basis $(e_1, \dots, e_r)$ for $V$ we can identify $V$ with $\R^r$, for $r: =\dim V$. We equip $\R^r$ with the standard Euclidean norm as an auxiliary tool to consider integrable functions. 
 Namely, we consider
 	 $\Omega: =L^2([0, 1];V)\cong L^2([0, 1];\mathbb R^r)$ and equip it with the $L^2$-norm
	\[
	\|u\|: =\left(\int_0^1 \sum_{i=1}^r u_i(t)^2 \dd t\right)^{\tfrac{1}{2}}.
	\]
	We refer to $\Omega$ as the {\em space of controls}.\index{space of controls}
 
 	For every element $u\in \Omega$, called {\em control}, the {\em trajectory associated} with $u$ is the curve $\gamma_u: [0, 1]\rightarrow G$ that is the solution of the ODE\index{control}\index{trajectory}
	\begin{equation}\label{ODE_horiz}
	\begin{cases}
	\gamma(0)=1_G, \\
	\dot\gamma(t)=\left(\dd L_{\gamma(t)}\right)u(t),& \text{ for a.e. }t\in [0, 1].
	\end{cases}
	\end{equation}
By Carath\'eodory Theorem on ODEs, see Proposition~\ref{prop_Integration_tangent_vector}, the equation is well posed, and in this way, each $u\in \Omega$ induces a $V$-horizontal curve $\gamma_u$ on $G$.\index{Carath\'eodory Theorem}\index{Theorem! Carath\'eodory --}  
Notice that such a $\gamma_u$ is absolutely continuous since $\dot \gamma$ is integrable by Cauchy--Schwarz inequality.
	Vice versa, every $V$-horizontal curve on $[0, 1]$ starting from $1_G$ is, up to reparametrization, of the form $\gamma_u$ for some $u\in \Omega$. 
	In fact, using the notation \eqref{def_gamma_prime}, if $\gamma$ is horizontal and reparametrized by constant speed, then $u: =\gamma'\in \Omega$ and $\gamma= \gamma_u$. 
	We call $u$ the {\em control} of $\gamma_u$.
	
The endpoint map is a key concept in sub-Finsler geometry, particularly from the perspective of control theory. It sends a control, and consequently the corresponding curve starting from a given base point, to the final point of the curve. This mapping allows for the analysis and optimization of trajectories. 

	The {\em endpoint map} is defined as\index{endpoint map}\index{map! endpoint --}
	\[
	\begin{aligned}
	\End: \Omega&\longrightarrow G\\
	u&\longmapsto \End(u): =\gamma_u(1),
	\end{aligned}
	\]
where $\gamma_u$ solves \eqref{ODE_horiz}.

\subsection{Differential of the endpoint map}
The differential of the endpoint map allows for sensitivity analysis, which examines how small changes in the control or initial conditions affect the reachable points. In sub-Finsler geometry, by analyzing this differential, one can derive necessary conditions for length-minimizing curves; see Section~\ref{sec_First-order necessary conditions}.

We shall not discuss why the endpoint map is smooth, for which we refer to \cite{Rifford:book}. We will directly calculate its (first) differential.
\begin{proposition}\label{prop:dif_end}\index{differential! -- of the endpoint map}\index{endpoint map! differential of the --}
Let $(G, V)$ be a polarized group with $\Omega$ as space of controls.
	For every $u\in \Omega$ the differential of $ {\End}$ at $u$ is given by
	\[
	\begin{aligned}
	\dd \End_u: \Omega&\longrightarrow T_{ \End(u)}(G ) \\
	v&\longmapsto \left(\dd R_{\gamma_u(1)}\right)_{1_G} \int_0^1 \Ad_{\gamma_u(t)}(v(t))\dd t, 
	\end{aligned}
	\]
	where $\Ad_g:\g\rightarrow \g$ is defined by $\Ad_g=(C_g)_*$ where $C_g(h)=ghg^{-1}$; see Section~\ref{sec Ad and ad}.
\end{proposition}
 \begin{proof} 
We first discuss the proof for matrix groups, i.e., we assume $G \subset \GL(n, \R)$, for some $n\in \N$, so we can interpret the Lie product as a matrix product and work in matrix coordinates. 
 Let $\gamma_{u+\eps v}$ be the curve with the control $u+\eps v$ and $\sigma(t)$ be the derivative of $\gamma_{u+\eps v}(t)$ 
 with respect to $\eps$ at $\eps=0$:
 $$
 \sigma(t): = \lim_{\eps\to0} \frac{1}{\eps} \left( \gamma_{u+\eps v}(t)-\gamma_{u }(t)\right)
 . $$
 Then one can verify (see Exercise~\ref{ex66b36bce}) that $\sigma$ satisfies the following ODE (which is the derivation with respect to $\eps$ of $(\ref{ODE_horiz})$ for $\gamma_{u+\eps v}$)
 \begin{eqnarray*}
 \frac{\dd\sigma}{\dd t} = \gamma(t) \cdot v(t) + \sigma(t) \cdot u(t).
 \end{eqnarray*}
Now it is easy to see that $t\mapsto \int_0^t \Ad_{\gamma(s)} ( v(s)) \dd s \cdot \gamma(t)$ satisfies the above equation with the same initial condition as $\sigma$; hence, it is equal to $\sigma$. 

To extend the proof to arbitrary Lie groups, one may use Ado's Theorem; see Exercise~\ref{ex7Aug2024}. References for the general case of manifolds are \cite[Section~20.3]{AS}, \cite[Section~5.2.2]{Montgomery}, or \cite[Section~8.1.1]{Agrachev_Barilari_Boscain:book}.
\end{proof}
\subsection{Singular curves}\label{sec: singular curves}\index{abnormal! -- curve}\index{singular curve}
In this subsection, we study those controls that are critical points for the endpoint map. The associated curves are called {\em singular curves} or {\em abnormal curves}. Studying these curves is one of the main difficulties in sub-Riemannian geometry. There are plenty of unsolved questions about them.

Fix a polarized group $(G, V)$, with $\Omega$ as space of controls.
If $u\in \Omega$ is a singular point for its endpoint map, then by definition, the linear map $\dd {\End}_u: \Omega\to T_{\End(u)}G$ is not surjective.
Namely, the subspace $\dd {\End}_u(\Omega)$ is a proper subspace of $T_{\End(u)}G$.
In this case, there is a nontrivial covector that annihilates its image, i.e., there exists an element 
$\xi$ in the dual space $(T_{\gamma_u(1)}G)^*$ such that 
$\xi \neq0$ and 
	\begin{equation}\label{eq1sep1139}
	\langle \xi | \dd {\End}_u(v)\rangle=0, \quad \forall v\in \Omega.
		\end{equation}
Here, we denote by $\langle \xi | w \rangle: = \xi(w)$ the evaluation of 
a dual element $\xi\in W^*$ at a vector $w\in W$, when $W$ is a vector space.

By Proposition~\ref{prop:dif_end}, Equation \eqref{eq1sep1139} is equivalent to 
	\begin{equation}\label{equations_abnormal1}
	0= \xi\left(\dd R_{\gamma_u(1)}\int_0^1 \Ad_{\gamma_u(t)} v(t) \dd t\right) =\lambda\left( \int_0^1 \Ad_{\gamma_u(t)} v(t) \dd t\right), \quad \forall v\in \Omega,
	\end{equation}
	where $\lambda \in \g^*\setminus\{0\}$ is defined as $\lambda: = \xi\dd R_{\gamma_u(1)}$.
	Fixing $t\in [0, 1]$ and $\bar v \in V$, let $v\in \Omega$ diverge to the distribution $\bar v \delta_t$ given by the Dirac mass at $t$ times $\bar v$, so to obtain $\lambda\left(\Ad_{\gamma_u(t)}\bar v\right)=0$. Hence, the equation
\begin{equation}\label{abnormal}
	\lambda\left(\Ad_{\gamma_u(t)}V\right)=\{0\}, \qquad \forall t\in [0, 1],
	\end{equation}
	is equivalent to \eqref{equations_abnormal1}.
	If $e_1, \ldots, e_r$ is a basis of $V$, then \eqref{abnormal} rephrases as a linear system of equations:
	A horizontal curve $\gamma $ is abnormal if and only if there exists $\lambda \in \g^*$ such that $\lambda\neq 0$ and
	\begin{equation}\label{abnormal1}\index{abnormal! -- equation}
	\lambda\left(\Ad_{\gamma(t)}(e_i) \right)=0, \qquad \forall i\in\{1, \ldots, r\}, \forall t\in [0, 1].
	\end{equation}
	We call \eqref{abnormal1} the {\em abnormal equation}. 
	We stress that if $\gamma$ satisfies \eqref{abnormal1} and $g\in G$, then the curve $t\mapsto g \gamma(t)$ still satisfies \eqref{abnormal1} with $\lambda$ replaced by $\lambda\Ad_g$. We still call all such curves {\em abnormal curves}, even if they do not start at $1_G$. Thus, left translations of abnormal curves are abnormal curves.
	Moreover, every restriction of an abnormal curve is an abnormal curve.
	
When, as in the case above, we have $\gamma(0)=1_G$, the last equation implies 
	\begin{equation}\label{abnormal1_at0}
	\lambda (e_i) =0, \qquad \forall i\in\{1, \ldots, r\}.
	\end{equation}
	Notice that, after we fix $i$ and $\lambda$, the function $g\mapsto\lambda\left(\Ad_{g}(e_i) \right)$ is smooth and \eqref{abnormal1} says that 
	$\gamma(t)$ lies in the zero-level set of such a function. 
	We shall notice that in nilpotent Lie groups (e.g., in Carnot groups), we have that, in exponential coordinates, the map $\Ad$ is polynomial; hence, these functions are polynomials.

	\begin{remark}
If $(G, V)$ is a polarized group with $V\neq \Lie(G)$, then constant curves are abnormal, and it may happen that the only abnormal curves are constant; see Exercise~\ref{ex_heis_abn}. Or it may happen that there are non-smooth examples; see Exercise~\ref{ex_nonsmooth_abn}.
		In Riemannian (and Finsler) geometry, i.e., if $V= \Lie(G)$, there are no abnormal curves. Indeed, there is no nonzero $\lambda$ as in \eqref{abnormal1_at0}. 
	\end{remark}

	\begin{definition}[Abnormal set and Sard Property]\index{abnormal! -- set}\index{Sard Property}
	The {\em abnormal set} of a polarized group $(G, V)$ is the subset ${\rm Abn}(G, V) \subset G$ of all singular values of the endpoint map. Equivalently, the set ${\rm Abn}(G, V)$ is the union of all abnormal curves passing through the identity element $1_G$.
If the abnormal set has measure $0$, then $(G, V)$ is said to satisfy the 
{\em Sard Property}.
		\end{definition}
	Proving the Sard Property in the general context of polarized manifolds is one of the major open problems in sub-Riemannian geometry; see the questions in \cite[Sec.~10.2]{Montgomery} and \cite[Problem~III]{Agrachev_problems}.

\section{Extrema in sub-Riemannian groups}\label{sec: 6.3}\index{extremal! -- curve}
 
\subsection{First-order necessary conditions for sub-Riemannian minimizers}\label{sec_First-order necessary conditions}\index{first-order necessary conditions}

	Let $G$ be a Lie group, and let $V\subseteq \g$ be a bracket-generating subspace.
	We shall consider sub-Riemannian structures on the polarized group $(G, V)$.
	A left-invariant sub-Riemannian structure is completely determined by the choice of a scalar product on $V$ and, hence, of an orthonormal basis $(e_1, \dots, e_r)$ for $V$. We next study conditions for length-minimizing curves for this sub-Riemannian structure.\index{length! --minimizer}\index{sub-Riemannian! -- structure}
	
	Recall from Proposition~\ref{prop: energy} that minimizing the length or minimizing the energy determines the same class of curves up to reparametrization with constant speed. And this is also why we can restrict to $L^2$ controls, i.e., elements in $\Omega: =L^2([0, 1]; V)$. Actually, because of Remark~\ref{anche_in_Lp} one can also take controls in $L^p$ with $p\in [1, \infty]$.
 
	
	
%
%
	We consider the {\em energy function}\index{energy! -- of a control}
	\[
	\begin{aligned}
	{\rm Energy}:\Omega&\longrightarrow \mathbb R\\
	u&\longmapsto {\rm Energy}(u): =\frac{1}{2}\|u\|^2.
	\end{aligned}
	\]
	This is the same functional that we saw in \eqref{eq:def_energy} for metric spaces, and now it is equal to 
$$	{\rm Energy}_{\dcc }(\gamma_u) \stackrel{\rm def}{=}\frac{1}{2} \int \| \dot{\gamma_u}\|^2=\frac{1}{2} \int \| {\gamma_u}'\|_{1_G}^2=\frac{1}{2}\|u\|^2, \qquad \forall u\in \Omega, $$
using notation \eqref{def_gamma_prime}.

	Together with the endpoint map, we form the {\em extended endpoint map}\index{extended endpoint map} 
	\[
	\begin{aligned}
	\widetilde{\End}:\Omega&\longrightarrow G\times \mathbb R	\\
	u&\longmapsto \big(\End(u), {\rm Energy}(u)\big).
	\end{aligned}
	\]
	
Given a point $p\in G$, minimizing the energy of curves between $1_G$ and $p$ rephrases as minimizing ${\rm Energy}(u)$ among all $u$ for which $\gamma_u(1)=p$. We shall say that $\gamma_u$ is a {\em minimizer for the energy}, or for short that $u$ is a {\em minimizer}, if for all $v\in \Omega$ we have\index{energy! -- minimizer}
$$\End(v)=\End(u)\implies {\rm Energy}(v)\geq {\rm Energy}(u).$$

\begin{remark}\label{End_is_singular}
	We claim that if $u_0$ is a minimizer for the energy, then $\widetilde{\End}$ cannot be open at any neighborhood of $u_0$ and therefore $u_0$ must be a singular point for
	 $\widetilde{\End}$. 
	Indeed, if there were a subset $U\subseteq \Omega$ for which	$ \widetilde{\End}(U)$ is a neighborhood of $\widetilde{\End}(u_0)$ within $G\times \mathbb R$, then we can find $\widetilde u \in U$ such that $\End(\widetilde u)=\End(u_0)$ and ${\rm Energy}(\widetilde u)<{\rm Energy}(u_0)$. This contradicts the minimality of $u_0$. 
	Moreover, if the differential of $(\dd \widetilde{\End})_{u_0}:\Omega \rightarrow T_{\widetilde \End(u)}(G\times \mathbb R)$ at $u_0$ were surjective, then we can take a vector subspace $W\subset\Omega$ for which  $(\dd \widetilde{\End})_{u_0}|_W:W \rightarrow T_{\widetilde \End(u)}(G\times \mathbb R)$ is an isomorphism. From the Implicit Function Theorem, we conclude that the map $\widetilde{\End}|_W: W \rightarrow G\times \mathbb R$ gives a diffeomorphism between some neighborhood of $ u_0$ within $W$ and some neighborhood of $\widetilde{\End}(u_0)$ within $G\times \mathbb R$. Such a fact contradicts the property that  $\widetilde{\End}$ cannot be open at $u_0$.
	\end{remark}
Because of this last remark, we need an expression for the differential of the extended endpoint map $\widetilde{\End}$.
After Proposition~\ref{prop:dif_end} and the standard calculation of the differential of the energy, see Exercise~\ref{ex differential of the energy}, we obtain the differential of the extended endpoint map:

\begin{corollary}\label{prop:dif_EXT_end}\index{differential! -- of the extended endpoint map}\index{extended endpoint map! differential of the --}\index{endpoint map! differential of the extended --}
Let $\widetilde{\End}$ be the extended endpoint map of a polarized group $(G, V)$ with $\Omega$ as space of controls.
	For every $u\in \Omega$ the differential of $\widetilde{\End}$ at $u$ is
	\[
	\begin{aligned}
	\dd\widetilde \End_u: \Omega&\longrightarrow T_{\widetilde \End(u)}(G\times \mathbb R)= T_{\End (u)}G\times \mathbb R \\
	v&\longmapsto \left(\left(\dd R_{\gamma_u(1)}\right)_{1_G} \int_0^1 \Ad_{\gamma_u(t)}(v(t))\dd t, \quad \langle u,v \rangle\right).
	\end{aligned}
	\]
\end{corollary}

	Assume now that, for some $u\in \Omega$, the curve $\gamma_u$ is length minimizing and parametrized with constant speed, so it is energy minimizing. By Remark~\ref{End_is_singular}, we deduce that $u$ is a critical point for $\widetilde{\End}$, i.e., the linear map
	$\dd\widetilde \End_u: \Omega \rightarrow T_{\End (u)}G\times \mathbb R $
	 is not surjective. Since then $\dd \widetilde \End_u( \Omega) $ is a strict subspace of $T_{\End (u)}G\times \mathbb R $,
 there exists $(\xi, \xi_0)\in\left(T_{\End(u)}G\right)^*\times \mathbb R = \left(T_{\End(u)}G\times \mathbb R\right)^*$ such that $(\xi, \xi_0)\neq (0,0)$ and 
	\[
	\langle(\xi, \xi_0) | \dd\widetilde{\End}_u(v)\rangle=0, \quad \forall v\in \Omega.
	\]
By Corollary~\ref{prop:dif_EXT_end}, this is equivalent to say that there exists $(\xi, \xi_0)\neq (0,0)$ such that 
	\begin{equation}\label{equation1}
	\xi\left(\dd R_{\gamma_u(1)}\int_0^1 \Ad_{\gamma_u(t)} v(t) \dd t\right)+\xi_0\langle u, v\rangle=0, \qquad \forall v\in \Omega.
	\end{equation}
	Since differentials of right translations give linear isomorphisms, Equation \eqref{equation1} is true if and only if there exist $\lambda \in \g^*$ and $\xi_0\in \mathbb R$ such that $(\lambda, \xi_0)\neq (0,0)$ and 
	\begin{equation}\label{equation2}
	\lambda\left( \int_0^1 \Ad_{\gamma_u(t)} v(t) \dd t\right)=\xi_0\langle u,v \rangle, \qquad \forall v\in \Omega.
	\end{equation}
	
	We now consider two cases: either $\xi_0\neq 0$ or $\xi_0= 0$. The first case is called {\em normal}, and the second one is called {\em abnormal}. We stress that if the codimension of $\dd \widetilde \End_u( \Omega) $ within $T_{\End (u)}G\times \mathbb R $ is strictly larger than 1, then there would be linearly independent choices for $(\lambda, \xi_0)$. Hence, some particular $u$ may have a normal pair $(\lambda, \xi_0)$ and a (different) abnormal pair $(\lambda', \xi_0')$.
		
	We separately consider the two cases:\index{normal! -- curve}\index{abnormal! -- curve}
	
	{\bf 	Normal case: $\xi_0\neq 0$.}
	Firstly, we suppose that $(\lambda, \xi_0)
	$ as in \eqref{equation2} is such that $\xi_0\neq 0$. Up to multiplying the equation by a constant, we can assume that $\xi_0=1$: 
		\begin{equation}\label{the first normal equation}
	\lambda\left( \int_0^1 \Ad_{\gamma_u(t)} v(t) \dd t\right)= \langle u,v \rangle, \quad \forall v\in \Omega.
	\end{equation}
	Fix a Lebesgue point $t$ of $u$, and let $v$ diverge to the distribution $e_i\delta_t $ given by the Dirac mass at $t$ times a basis vector $e_i$ of $V$. 
	Formally,
	we have
	\[
	\begin{aligned}
	\dot\gamma_u(t)&=\dd L_{\gamma_u(t)}u(t)=\dd L_{\gamma_u(t)}\sum_{i=1}^r\langle u, \delta_t e_i \rangle e_i\\
	&\stackrel{\eqref{the first normal equation}}{=}\dd L_{\gamma_u(t)} \sum_{i=1}^r\left(\lambda \int_0^1 \Ad_{\gamma_u(s)}(\delta_t e_i)\dd s\right) e_i\\
	&=\sum_{i=1}^r\lambda \left(\Ad_{\gamma_u(t)}(e_i)\right)X_i(\gamma_u(t)),
	\end{aligned}
	\]
where in the last equality we have used the identity $X_i(g)=\left(\dd L_g\right)e_i$.
	We therefore say that an absolutely continuous curve $\gamma$ satisfies the \emph{normal equation} (or the {\em sub-Riemannian geodesic equation}) with respect to the left-invariant orthonormal frame $X_1, \ldots, X_n$,
	if there exists $\lambda\in \g^*$ such that\index{normal! -- equation}\index{sub-Riemannian! -- geodesic equation}
	\begin{equation}\label{geodesic}
	\dot\gamma(t)=\sum_{i=1}^r\lambda \left(\Ad_{\gamma(t)}(e_i)\right)X_i(\gamma(t)), \qquad \text{ for a.e. } t\in [0, 1].
	\end{equation}
A solution to \eqref{geodesic} is called {\em normal curve}. By a bootstrap argument using \eqref{geodesic}, we deduce that the horizontal curve $\gamma$ and its control $u$ are $C^\infty$. In fact, they are analytic; see Exercise~\ref{ex analytic normal}.

Writing $ u = \sum_{i=1}^r u_i e_i $ and recalling \eqref{ODE_horiz}, we get another version of the normal equation:
	\begin{equation}\label{normal_final}
	u_i (t)= \lambda \left(\Ad_{\gamma_u (t)}(e_i)\right), \qquad \text{ for a.e. } t\in [0, 1] \text{ and }\forall i\in\{1, \ldots, r\}.
	\end{equation}
	In particular, since in our case $\gamma(0)=1_G$, the last equation implies 
	\begin{equation}\label{attenzione}
	u_i(0
	) =\lambda (e_i), \qquad \forall i\in\{1, \ldots, r\}.
	\end{equation}
In other words, the `horizontal' part of the covector $\lambda$ is the velocity at the initial time.

	{\bf 	Abnormal case: $\xi_0= 0$.}
	Secondly, we suppose that $(\lambda, \xi_0)
	$ as in \eqref{equation2} is such that $\xi_0=0$. 
	In other words, the curve $\gamma_u$ is singular for the endpoint map, as studied in Section~\ref{sec: singular curves}. These are the curves that we called abnormal. In \eqref{abnormal1}, we saw that	
	a horizontal curve is abnormal if and only if there exists $\lambda \in \g^*$ such that $\lambda\neq 0$ and
	\begin{equation}\label{abnormal1_unaltro}
	\lambda\left(\Ad_{\gamma(t)}(v) \right)=0, \qquad \forall t\in [0, 1], \forall v\in V.
	\end{equation}
	In other words, the set $\left\{ \Ad_{\gamma(t)} v\,:\,v\in V, \, t\in I\right\}$ of vectors in $\g$ is contained in a proper subspace of $\g$.
	
%

We summarise the obtained results with the following statement, which we attribute to Pontryagin because, in fact, he obtained a much more general result for control systems. 

\begin{theorem}[Pontryagin's maximum principle]\label{thm_Pontryagin}\index{Pontryagin Maximum Principle}
 Let $G$ be a sub-Riemannian group with polarization $V$. 
 Let
$\gamma: I \to G$ be an energy-minimizing curve with control $u$.
Then, at least one of the following happens:
\begin{description}
\item[(i)] there is $\lambda\in \g^*$ 
such that $ \langle u(t), X\rangle = 
\lambda \left(\Ad_{\gamma (t)}(X)\right)$
for all $X\in V$ and all $t\in I$, or 
\item[(ii)]
$ 
\Span \left\{ \Ad_{\gamma(t)} V\,:\, t\in I\right\} \neq \g$.
\end{description}
\end{theorem}

Each curve satisfying any of the above conditions is called an {\em extremal curve} or an {\em extremal}, and the equations are called {\em extremal equations}.\index{extremal! -- equation}\index{extremal! -- curve}

 	By a calibration argument, one can prove that every normal curve is length minimizing on short enough intervals; see \cite[Section~1.9.3]{Montgomery}. The converse is not true, and the first (surprising) example has been found by Montgomery \cite{Montgomery1994}.
	In other words, there are sub-Riemannian structures where it is possible to find
energy-minimizing curves that are not normal, and so are abnormal. 
They are called {\em strictly abnormal geodesics}.\index{strictly abnormal geodesics}	Be aware that there are geodesics that are normal and abnormal, and there are abnormal curves that are not geodesics; see Exercise~\ref{ex_nonsmooth_abn}.

Not much is known about strictly abnormal geodesics. We close this discussion with a second-order theorem, which we do not prove here. We refer to \cite[Chapter~12]{Agrachev_Barilari_Boscain:book}.
\begin{theorem}[Goh]\label{Goh_thm}\index{Goh! -- Theorem}\index{Theorem! Goh --}\index{Goh! -- condition}
 Let $G$ be a sub-Riemannian group with polarization $V$. 
 If $\gamma: I \to G$ is a strictly abnormal geodesic, then
$\Span \left\{ \Ad_{\gamma(t)} (V+[V, V])\,:\, t\in I\right\} \neq \g $.
\end{theorem}

\subsection{A distinguished class of functions}
 Both in the normal and the abnormal equations, \eqref{normal_final} and \eqref{abnormal1}, the functions 
$t\mapsto\lambda\left(\Ad_{\gamma(t)}(e_i) \right)$ are considered. We naturally see these functions as the composition of the curve with one of the following functions. 

For $\lambda\in \g^*$ and $Y\in \g$, define $P^\lambda_Y: G\to \R$ as\index{abnormal! -- polynomials}
\begin{equation}\label{abn_poly}
P^\lambda_Y(g): =\lambda\left(\Ad_{g}(Y) \right), \qquad \forall g\in G.
\end{equation}
These maps are linear both in $\lambda$ and in $Y$. Moreover, we will see that if $G$ is nilpotent, they are polynomials in $g$ when seen in exponential coordinates; see Section~\ref{sec Canonical coordinates}.
A useful formula that these functions satisfy is the following:
\begin{equation}\label{derfv_abn_poly}
X P^\lambda_Y = P^\lambda_{[X, Y]}, \qquad \forall X, Y\in \g, \forall \lambda\in \g^*.
\end{equation}
Indeed, fixing in addition $g\in G$, we first notice that 
\begin{eqnarray*}
 \left.\frac{\dd }{\dd t} \Ad_{g\exp(tX)}\right|_{t=0}
 &=& 
 \left.\frac{\dd }{\dd t} \Ad_g\Ad_{\exp(tX)}\right|_{t=0}\\
 &=& 
\Ad_g \left.\frac{\dd }{\dd t} \Ad_{\exp(tX)}\right|_{t=0} \stackrel{\eqref{eq Ad ad}}{=} \Ad_g \circ \ad_X, \\ 
 \end{eqnarray*}
 see also Exercise~\ref{ex_derivata_Lambda}.
Then, we deduce 
\begin{eqnarray*}
(X P^\lambda_Y )(g)& =& \left.\frac{\dd }{\dd t} \lambda \Ad_{g\exp(tX)}Y\right|_{t=0} \\
&=& (\lambda\circ \Ad_g \circ \ad_X)Y\\
&=& \lambda \Ad_g ( [X, Y])
= P^\lambda_{[X, Y]}(g).
 \end{eqnarray*}

If we have a normal curve $\gamma$ with covector $\lambda$, the {\bf normal equation} \eqref{geodesic} reads as\index{normal! -- equation}
\begin{eqnarray}	 \gamma'(t)&=&\sum_{i=1}^r\lambda \left(\Ad_{\gamma_u(t)}(e_i)\right)e_i
\nonumber
\\\label{normal_eq_20jun}
&=&\sum_{i=1}^r P^\lambda_{e_i}(\gamma (t) )e_i
, \qquad \text{ for a.e. } t\in [0, 1].\end{eqnarray}

We next deduce that normal curves are parametrized with constant speed: from \eqref{normal_eq_20jun}
we have that
the derivative of the squared speed is 
	\begin{eqnarray*}
	\frac{\dd }{\dd t}\norm{\dot\gamma(t)}^2&\stackrel{\eqref{normal_eq_20jun}}{=}&\frac{\dd }{\dd t} \sum_{i=1}^r \left( P^\lambda_{e_i}(\gamma (t))\right)^2 \\
	&=& \sum_{i=1}^r 2 P^\lambda_{e_i}(\gamma (t)) \frac{\dd }{\dd t} P^\lambda_{e_i}(\gamma (t)) \\
	&\stackrel{\eqref{derfv_abn_poly}}{=} & \sum_{i=1}^r 2 P^\lambda_{e_i}(\gamma (t)) 
	P^\lambda_{[\gamma' (t),e_i ]}(\gamma (t)) \\
		&\stackrel{\eqref{normal_eq_20jun}}{=} & \sum_{i=1}^r 2 P^\lambda_{e_i}(\gamma (t)) 
	P^\lambda_{[
	\sum_{j=1}^r P^\lambda_{e_j}(\gamma (t) )e_j
	,e_i ]}(\gamma (t)) \\
		&= & \sum_{i,j=1}^r 2 P^\lambda_{e_i}(\gamma (t)) 
	 P^\lambda_{e_j}(\gamma (t) )
	P^\lambda_{[
	e_j
	,e_i ]}(\gamma (t)) 
	=0
	,
	\end{eqnarray*}
	where at the end we used that the elements of this sum are antisymmetric in the indices due to the fact that 
	$[e_j, e_i ] =-[e_i, e_j ]$.
Hence, the speed of normal curves is constant.

In terms of the functions \eqref{abn_poly}, the {\bf abnormal equation} 
\eqref{abnormal1} for a curve $\gamma$ reads as
follows:
there exists $\lambda \in \g^*$ such that $\lambda\neq 0$ and
	\begin{equation}\label{abnormal1_poly}\index{abnormal! -- equation}
	P^\lambda_{e_i}\circ \gamma 
	\equiv 0, \qquad \forall i\in\{1, \ldots, r\}.
	\end{equation}


\subsection{Extremals in groups with rank-2 polarizations}\index{rank! ---two polarization}\index{extremal! -- equation}

We rephrase the extremal equations in the case of polarized groups where the polarization has rank 2.
We begin with the abnormal condition:

\begin{proposition}\label{prop_rank-2_abnormal}
Let $(G, V)$ be a sub-Riemannian group whose $V$ is spanned by two vectors $e_1, e_2$.
Let $\gamma: [0, 1] \to G$ be a horizontal curve. 
 Then $\gamma$ is 
 abnormal if and only if for some $\lambda\in \g^*$ with $\lambda\neq0$ and $\lambda (e_1) =\lambda (e_2) =0$ it
 satisfies \begin{equation}\label{abnormal12}\lambda(\Ad_{\gamma(t)} ([e_{1},e_2]))=0.\end{equation} 
 \end{proposition}

\begin{proof}
We use the notation $e_{12}: =[e_1,e_2]$.
We write $u(t): =\gamma'(t)= u_1 (t) e_1 + u_2(t) e_2 $ and we have
\begin{equation}\label{calcolo_commutatori_step2}[u(t),e_1] = - u_2 e_{12}\qquad \text{ and } \qquad
 [u(t),e_2]= u_1 e_{12} .	\end{equation}
 Taking the derivative of the abnormal equations, see \eqref{equazione_abnormale} in Exercise~\ref{ex:First derivative of the extremal equations}, we have that if
 $\gamma_u$ is an abnormal normal curve with covector $\lambda\in \g^*$, with $\lambda\neq0$, then 
\begin{equation}\label{abnormale1_2}u_2 \lambda(\Ad_{\gamma(t)} (e_{12}))= u_1 \lambda(\Ad_{\gamma(t)} (e_{12})) =0 .\end{equation}
In addition, notice that we may assume that $\gamma_u $ is parametrized by arc length, so $u$ has a constant nonzero norm almost everywhere. In particular, we have $(u_1,u_2)\neq(0,0)$ almost everywhere. 
Therefore, we can conclude that for such an abnormal curve, we have \eqref{abnormal12}.

Vice versa, assume $\gamma$ satisfies \eqref{abnormal12} for some $\lambda\in \g^*$ with $\lambda\neq0$. Then, it clearly satisfies \eqref{abnormale1_2} and, since we have $r=2$ and we have \eqref{calcolo_commutatori_step2}, we also have \eqref{equazione_abnormale}.
Then, look at each function $\lambda\Ad_{\gamma(t)}(e_i)$, for $i=1$ and $2$. 
On the one hand, because of \eqref{derivata_Lambda}, we have that its derivative is 0. On the other hand, if $\gamma(0)=1_G$ and if $\lambda$ satisfies \eqref{abnormal1_at0}, we have that the initial condition at time $t=0$ for \eqref{abnormal1} is satisfied. Hence, such a curve is abnormal. 
\end{proof}
 
 
We shall rephrase the normal equation in terms of a curvature condition for the development. See Definition~\ref{def development} for the notion of development, and Exercise~\ref{ex_oriented_curvature} for the oriented curvature.

 \begin{proposition}\index{development}\index{curvature! condition}
Let $(G, V)$ be a sub-Riemannian group whose $V$ is spanned by two orthonormal vectors $e_1, e_2$.
Let $\gamma: [0, 1] \to G$ be a horizontal curve with development $\sigma$. Assume $\sigma$ is $C^2$, with never vanishing speed, and let $\kappa$ be the oriented curvature of $\sigma$.
 Then $\gamma$ is 
 normal if and only if for some 
 $\lambda\in [\g, \g]^*$ we have that $\gamma$ satisfies 
 \begin{equation}\label{curvature_normal}
 \kappa =
 \frac{1}{\|\gamma'\|}\lambda(\Ad_{\gamma(t)} ([e_1,e_2]))
 .\end{equation}
 \end{proposition}

 \begin{proof}
Let $e_{12}: =[e_1,e_2]$. We take the derivative of the normal equations; see \eqref{equazione_normale}. By \eqref{calcolo_commutatori_step2}
 we have that if $\gamma=\gamma_u$ is a normal curve with covector $\lambda\in \g^*$ and control $u=\gamma'$, then 
\begin{equation}\label{uffa}\left\{\begin{array}{ccc}
\dot u_1 &=& -u_2 \lambda(\Ad_{\gamma(t)} (e_{12})), \\
\dot u_2 &=& u_1 \lambda(\Ad_{\gamma(t)} (e_{12})).
\end{array}\right. \end{equation}
Let $\sigma: [0, 1]\to \R^2$ be the planar curve such that $\dot\sigma= u$.
By Exercise~\ref{ex_oriented_curvature} its curvature satisfies
$$\kappa = \frac{ \sigma'_1 \sigma''_2- \sigma'_2 \sigma''_1 }{\|\sigma' \|^3} =
 \frac{ u_1 \dot u_2- u_2 \dot u_1 }{\|u\|^3} \stackrel{\eqref{uffa}}{=}
 \frac{ u_1 ^2 \lambda(\Ad_{\gamma} (e_{12}))+ u_2 ^2 \lambda(\Ad_{\gamma} (e_{12}) ) }{\|u\|^3} =
 \frac{1}{\|u\|}\lambda(\Ad_{\gamma(t)} (e_{12})).$$
 We observe that the element $ \Ad_{\gamma(t)}(e_{12}) $ is in $ [\g, \g]$, hence in the last equation we have lost the information of the value of $\lambda$ on $V$. Still, normal curves need to satisfy \eqref{abnormal1_at0}.
 
 Vice versa, we assume $\gamma$ is a horizontal curve in a rank-2 Carnot group and that for some $\lambda\in [\g, \g]^*$, we have that $\gamma$ satisfies \eqref{curvature_normal}.
 First, we observe that by bootstrapping \eqref{curvature_normal}, we have that $\gamma$ and its control $u$ are smooth. Then, we can extend $\lambda$ as an element of $\g^*$ so that we also have \eqref{attenzione}. Now that we have $\lambda \in \g^*$, we consider the normal curve (which is unique) associated with $\lambda$ with $\gamma(0)=1_G$, which we denote by $\gamma^\lambda$. We shall show that $\gamma=\gamma ^\lambda$. 
 
 The reason is that both curves satisfy the ODE \eqref{curvature_normal} with the same initial data. 
 We stress that such an ODE is not of the simple type $u'' (t) = F(t, u(t),u'(t))$, but it is of the type $F(t, u, u',u'')=0$. However, by the implicit function theorem, one can write the ODE of the second type in the form of the first type, as long as $\partial_{u''} F \neq 0$; see Exercise~\ref{ex ODE all'ODE}. 
 We then obtain an implicit function with a bound on the $C^1$ norm as long as the control $u$ is in a compact set. 
 \end{proof}
%
%
%


\subsection{Extrema in step-two nilpotent Lie groups}\index{step-two nilpotent Lie group}

In this section, we study a particular class of Lie groups. We denote each such group by $G_q$, where $q$ is a skew-symmetric bilinear form, which captures the Lie bracket.
We will see later that these groups $G_q$ are exactly the simply connected Lie groups that are nilpotent with nilpotency step 2; see Example~\ref{ex_group_step2_20jun}.
\begin{definition}[The Lie group $G_q$]\label{def_step_2_G}
Let $V$ and $W$ be finite-dimensional vector spaces and $q: V\times V\to W$ a skew-symmetric bilinear map.
Equip $V\times W$ with the group law 
\[
(v_1,w_1)\cdot(v_2,w_2): = \left( v_1+v_2, \quad w_1+w_2 + \frac12 q(v_1,v_2) \right),
\qquad \forall v_1,v_2\in V, \forall w_1,w_2\in W,
\]
whose identity element is $(0,0)$ and the inverse of $(v_1,w_1)$ is $(-v_1,-w_1)$.
We denote the group $(V\times W, \cdot)$ with $G_q$.
 We equip the vector space $V\times W$ with the Lie bracket
\[
[(v_1,w_1),(v_2,w_2)]: = (0, q(v_1,v_2)) .
\]
\end{definition}

Notice that $(V\times W,[\cdot, \cdot])$ is a Lie algebra and for every $x,y,z\in V\times W$, we have $[x,[y,z]] = 0$.
Moreover, when equipped with the differential structure of vector space, we have that $G_q$ is a Lie group, and its tangent at $(0,0)$ is naturally identified with $V\times W$.
The identity map $(V\times W,[\cdot, \cdot]) \to (V\times W, \cdot)$ is the exponential map; see Exercise~\ref{ex_5Jun2024}.

For $(v_1,w_1),(v_2,w_2)\in V\times W$ we have
\begin{eqnarray}
\Ad_{(v_1,w_1)}(v_2,w_2)\nonumber
&= &\Ad_{\exp(v_1,w_1)}(v_2,w_2) \\\nonumber
&=& e^{\ad_{(v_1,w_1)}}(v_2,w_2) \\\nonumber
&= &(v_2,w_2) + [(v_1,w_1),(v_2,w_2)] \\
&= &(v_2,w_2 + q(v_1,v_2)) .
\label{eq_Ad_gq}
\end{eqnarray}

\subsubsection{Abnormal curves in step-two nilpotent groups}
\begin{proposition}\label{existence: subgroup}
Let $G=G_q$ be a step-two nilpotent group equipped with a polarization. For each abnormal curve $\gamma$ in $G$ with $\gamma(0)=1_G$,
 there exists a proper subgroup $G'$ of $G$ containing $\gamma$, in which $\gamma$ is a non-abnormal horizontal curve, with respect to the induced polarization.
\end{proposition}

	\begin{proof}
Let $G'$ be the smallest subgroup of $G$ containing $\gamma$, which is still of the form $G_{q'}$; see Exercise~\ref{ex_sub_Gq}. 
Equip $G'$ with the induced polarization $\Delta$ from $G$, which by the minimality of $G'$ is bracket generating.
Assume by contradiction that $\gamma $ is abnormal in $G'$. 
By the abnormal equation \eqref{abnormal1}, 
 there exists $\lambda \in \Lie(G')^*$ such that $\lambda\neq 0$ and
	\begin{equation}\label{abnormal1_Gq}
	0= \lambda\left(\Ad_{\gamma(t)}( X ) \right) \stackrel{\eqref{eq_Ad_gq}}{=}
	\lambda(X + [\gamma(t), X] )
	\stackrel{\eqref{abnormal1_at0}}{=}
	\lambda( [\gamma(t), X] )
	, \qquad \forall X\in\Delta_1.
	\end{equation}
We consider the set 
$$Z: = \{ (v,w) \in V\times W \;:\; \lambda( [v, X] )=0
	, \forall X\in\Delta_1\}.$$
	Hence, the curve $\gamma$ is contained in $Z$. Moreover, the set $Z$ is a proper subgroup of $G_q$, because in the definition of $Z$ the dependence on $(v,w)$ is linear, $[\g,\g]\subseteq Z$, and $\lambda$ is not zero on $[\Delta_1, \Delta_1]$. We contradicted the minimality of $G'$. Hence, the curve $\gamma$ is not abnormal in $G'$.
	\end{proof}

	As a consequence, since every step-2 nilpotent group is a quotient of some $G_q$, we have that every geodesic curve in a group of nilpotency step 2 cannot be abnormal in the minimal group containing it. Since the curve is geodesic also within this subgroup, then, by PMP Theorem~\ref{thm_Pontryagin}, the curve is normal in the subgroup. We deduced, without using Goh Theorem~\ref{Goh_thm}, the following consequence:

\begin{corollary}
 Let $G $ be a step-two nilpotent group. 
 Equip $G$ with a left-invariant sub-Riemannian structure.
Then, every geodesic curve in $G$ is smooth.
\end{corollary}



\subsubsection{Normal curves in step-two nilpotent groups}\index{normal! -- curve! -- in step-two nilpotent group}

\begin{proposition}
 Let $G = G_q = V\times W$ be a step-two nilpotent group as in Definition~\ref{def_step_2_G}.
Equip $G$ with a left-invariant sub-Riemannian structure with polarization 
 $\Delta\subset V\times W$ such that 
 $V \subset\Delta$ and $V$ is the orthogonal of $W\cap \Delta$ in $\Delta$.
 Then,
 for every energy-minimizing curve \(\gamma = (v,w): [0, T] \to G\) with $\gamma(0)=(0,0)$,
 there are $M\in \mathfrak{so}(V)$, $b\in\ker(M)$, $c\in M(V)$, and $\zeta \in W \cap \Delta$ such that 
 \begin{equation}\label{ODE_step2}
\begin{cases}
	v(t) = c - e^{tM}c + tb, \\
	w(t) = \frac12\int_0^t [v(s), \dot v(s)] \dd s + t \zeta,
\end{cases} \qquad \forall t\in [0, T].
\end{equation}
\end{proposition} 
 \begin{proof}


By Pontryagin Theorem \ref{thm_Pontryagin}, the curve $\gamma$ is normal or abnormal. 
By Goh Theorem~\ref{Goh_thm}, the curve $\gamma$ is normal.
Let $e_1, \dots,e_r$ be an orthonormal basis of $V$
and $f_1, \dots,f_m$ an orthonormal basis of $W\cap\Delta$.
The normal equation \eqref{normal_eq_20jun} says that there is $\lambda \in \g^*$ such that 
\begin{equation}\label{eq6676c99f}
	\gamma'(t) = \sum_{i=1}^r \lambda(\Ad_{\gamma(t)}e_i) e_i 
		+ \sum_{i=1}^m \lambda(\Ad_{\gamma(t)}f_i) f_i .
\end{equation}

For $\gamma = (v,w)$ 
we decompose 
$\gamma' (t) =(\gamma'(t))_1+(\gamma'(t))_2$ with
 $(\gamma'(t))_1\in V$ and $(\gamma'(t))_2\in W$. Then, 
we have
\begin{align*}
(\dot v, \dot w)
&= \dot\gamma 
= \dd L_{\gamma} (\gamma' )
= \left( (\gamma')_1, \frac12[ v,(\gamma')_1] + (\gamma')_2 \right).
\end{align*}

Since $\Ad_{( v, w)} X
= X + [ v, X]$,
we obtain from~\eqref{eq6676c99f} that 
 the normal curve $\gamma = ( v, w)$ satisfies the ODE 
\[
\begin{cases}
\dot v&= \sum_{i=1}^r \lambda(e_i + [ v,e_i]) e_i \\
\dot w &= \frac12[ v, \dot v] + \sum_{i=1}^m \lambda(f_i) f_i .
\end{cases}
\]
We define $M:V\to V$ to be the linear map $Mx: = \sum_{i=1}^r \lambda([x,e_i]) e_i$
and $(\lambda_H,\zeta)\in V\times W$ to be the vectors $\lambda_H: = \sum_{i=1}^r \lambda(e_i)e_i$
and $ \zeta: = \sum_{i=1}^m \lambda(f_i) f_i$.
We can rewrite the ODE as
\begin{equation}\label{eq6676cf5b}
\begin{cases}
\dot v&= M v+ \lambda_H\\
\dot w &=\frac12 [ v, \dot v] + \zeta .
\end{cases}
\end{equation}
The linear transformation $M$ is skew-symmetric, because the matrix representation of $M$ in the orthonormal basis $e_1, \dots,e_n$ is $M_{ij} = \lambda([e_i,e_j])$.
The solution of~\eqref{eq6676cf5b} with $ v(0) = w(0) = 0 $ is \eqref{ODE_step2},
where $\lambda_H = b-M c$ is decomposed into $b\in\ker(M)$ and $c\in M(V )$.
 \end{proof}

%
%
%
%
%
%
 \begin{corollary}\label{cor infinite geodesics step 2}\index{one-parameter subgroup}\index{geodesic! -- ray}
Let \( G \) be a simply connected sub-Riemannian group of nilpotency step 2 with polarization $\Delta$.
Every isometric embedding \(\gamma:[0, \infty) \to G \) with \( \gamma(0) = 1_G \) 
 is the restriction of a OPS with $\dot\gamma(0)\perp (\Delta_1\cap[\g, \g])$.
\end{corollary}
\begin{proof}
Let   $G=G_q=V\times W$ be for some $q:V\times V \to W$ as in Definition~\ref{def_step_2_G}, with  $V \subset\Delta$ and $V\perp W\cap \Delta$.
Since the bracket given by $q$ is a bilinear map,  there exists $K>0$ such that 
\begin{equation}\label{Nalon linear}
    \lvert[v_1,v_2]\rvert \le K\lvert v_1 \rvert \lvert v_2 \rvert \quad \text{for every $v_1,v_2 \in \g$},
\end{equation}
where $|\cdot|$ is a fixed auxiliary norm on $\g$ that is equal to the sub-Riemannian norm on $\Delta_1$.

Let $\gamma \colon [0,T] \to G$ be isometric. Hence, it is of the form 
$\gamma(t)=v(t)+w(t)$, with $v(t) \in V$ and $w(t)\in W$ for every $t \in [0,T]$, satisfying   \eqref{ODE_step2}, for some
 $M\in\mathfrak{so}(V)$, $b \in \ker(M)$, $c \in M(V)$, and $\zeta \in W $.
Since $\gamma$ is parametrized by arc length, we have 
$\lvert Mc \rvert^2 + \lvert b \rvert^2 + \lvert \xi \rvert^2=1$.

Since $e^{T M}$ is an isometry and $c, e^{T M}c \in M(V) \subseteq  b^\perp$, from \eqref{ODE_step2} we get  
\begin{equation}\label{Eq Nalon 23 sept}
 d_\mathrm{sR}(0,v(T)) \stackrel{\rm Ex.\ref{ex OPS geodesic 2}}{=}\lvert v(T) \rvert \stackrel{\eqref{ODE_step2}}{=} \lvert c- e^{T M}c+T b \rvert \le \sqrt{|2c|^2 + |Tb|^2   }.
\end{equation} 
Next, we bound $w(T)$. 
By \eqref{ODE_step2}, we further have
\begin{eqnarray*}
    \lvert w(T) \rvert &\stackrel{\eqref{ODE_step2}}{\le}& \frac{1}{2}\bigg\lvert \int^T_0 [c-e^{tM}c+tb,-e^{tM}Mc + b]\dd t \bigg\rvert + T\lvert \zeta \rvert \\
    &\le& \frac{1}{2}\bigg\lvert \int^T_0 [c-e^{tM}c, e^{tM}Mc ]\dd t\bigg\rvert + \frac{1}{2}\bigg\lvert \int^T_0 [c-e^{tM}c,b]\dd t \bigg\rvert + \frac{1}{2}\bigg\lvert \int^T_0 [tb, e^{tM}Mc]\dd t \bigg\rvert + T\lvert \zeta \rvert
   \\
        & \stackrel{\eqref{Nalon linear}}{\le}&
         \frac{1}{2}  \int^T_0 K 2 |c| |Mc|  \dd t  +
          \frac{1}{2}  \int^T_0 2 | c | |b|  \dd t 
           + \frac{1}{2}\bigg\lvert \int^T_0 [ b, t e^{tM}Mc]\dd t \bigg\rvert + T\lvert \zeta \rvert
              \\
        &\stackrel{\eqref{Nalon linear}}{\le}&
            T  K   |c| |Mc|     +
           T   K | c | |b|   
           + \frac{1}{2}\bigg\lvert \left[ b , \int^T_0 te^{tM}Mc \dd t  \right]  \bigg\rvert + T\lvert \zeta \rvert
   ,
\end{eqnarray*}
where we used the linearity of the Lie bracket. 
Focusing on the third term, we 
integrate by part:
\begin{equation*}
\frac{1}{2}\bigg\lvert \left[ b , \int^T_0 te^{tM}Mc \dd t  \right]  \bigg\rvert 
=\frac{1}{2}\bigg\lvert \left[ b ,  Te^{T M}c- \int^T_0 e^{tM} c \dd t\right]\bigg\rvert 
\stackrel{\eqref{Nalon linear}}{\le}    T K\lvert c \rvert \lvert b \rvert 
.
\end{equation*}
All together,  we obtain a bound that is linear in $T$:
\begin{equation}\label{23 Sept 2024}
    \lvert w(T) \rvert \le   T (K \lvert c \rvert  \left( \lvert Mc \rvert + 2   \lvert b \rvert \right) +  \lvert \zeta \rvert).
\end{equation}
Lie groups of step 2 satisfy a global upper bound given by the box distance; see Exercise~\ref{Ex Nalon}. In particular, using \eqref{Ex Nalon eq}, we deduce that for some constants $C$ and $\tilde C= \tilde C(G, M,c,b)$ we have
\begin{eqnarray*}
T &=& d_\mathrm{sR}(0,v(T)+w(T))\\
&\le& d_\mathrm{sR}(0,v(T)) + d_\mathrm{sR}(0,w(T)) \\
&\stackrel{\eqref{Eq Nalon 23 sept} \& \eqref{Ex Nalon eq}}{\le}&  \sqrt{|Tb|^2 + |2c|^2} +  C \lvert w(T) \rvert^\frac{1}{2} \\
& \stackrel{\eqref{23 Sept 2024}}{\le}& \sqrt{|Tb|^2 + |2c|^2} + C \sqrt{T (K \lvert c \rvert  \left( \lvert Mc \rvert + 2   \lvert b \rvert \right) +  \lvert \zeta \rvert)} \\
&\le& T|b| + \tilde C (\sqrt{T}+T^{-1}) + \tilde C.
\end{eqnarray*}
If $|b|<1$, then the latter inequality cannot hold for arbitrarily large $T$.
We conclude that, if $\gamma: [0,\infty)\to G$ is an isometric embedding, then   $|b|=1$ and thus $c = \zeta= 0$, i.e., the curve $\gamma$ is a straight lines with directions orthogonal to $[\g, \g]$.
\end{proof}

\section{Characterization of geodesic left-invariant distances}\label{sec: 6.4}\index{characterization! -- of sub-Finsler metrics}
 Berestovskii's work \cite{b, b1, b2} clarified what are the possible isometrically homogeneous distances on manifolds that, in addition, are geodesic distances: They are sub-Finsler metrics.
\begin{theorem}[Berestovskii]\label{Berestovskii}\index{Lie! -- coset space}\index{Berestovskii! -- characterization}
Let $M: =G/H$ be the Lie coset space of a Lie group $G$ modulo a closed subgroup $H$. If $M$ is metrized by a geodesic distance that is $G$-invariant,
then the distance is a sub-Finsler metric, i.e., there is a $G$-invariant subbundle $\Delta$ on $M$ and a $G$-invariant norm on $\Delta$, such that the distance is given by the same formula \eqref{dist_CC}. 
\end{theorem}

We will only partially prove Theorem~\ref{Berestovskii}; one can check the original reference \cite[Theorem~2]{b}. We will discuss how to give a constructive characterization of the sub-Finsler structure.

At the beginning of the proof of Theorem~\ref{Berestovskii}, there is a crucial rectifiability result. Namely, one proves that the geodesics and, more generally, the curves of finite length for the geodesic distance on $G/H$ also have finite length for Riemannian distances. Consequently, such curves are absolutely continuous. Their derivatives will necessarily form the subbundle of the sub-Finsler structure; see \eqref{Berestovskii Delta}.

\begin{lemma}[Crucial lemma]\label{rectgeo}
Let $M: =G/H$ be the Lie coset space of a Lie group $G$ modulo a closed subgroup $H$.
Let $d$ be an admissible geodesic distance on $M$ that is $G$-invariant.
Then, for every Riemannian $G$-invariant distance $d_R$ on $M$ there is a constant $C>0$ such that
$$d_R \leq C d .$$
\end{lemma}
\begin{proof}
For simplicity of exposition, we consider only the case when $H$ is trivial and $M=G$.
As a preparation, via Proposition~\ref{exp:diffeo}, let $\delta>0$ such that there exists a neighborhood of $U$ of $ 0$ in $\g$ such that 
$\exp|_U: U \to B_{d_R}(1_G, \delta) $ is biLipschitz for some constant $k>1$, in particular,
\begin{equation}\label{eq_bilip4may} 
\frac1k\norm{\log(p)}\leq d_R( 1_G,p) \leq k\norm{\log(p)}, \qquad \forall p\in B_{d_R}(1_G, \delta).
\end{equation}
Also, since the topologies of the two metrics are the same, we can find $r>0$ such that 
\begin{equation}\label{4 May 1047}
\bar B_d(1_G,r) \subseteq B_{d_R}(1_G, \delta). 
\end{equation}

We prove the lemma by contradiction and assume that for each $n\in \N$, there exist $p_n, q_n\in M$ such that
\begin{equation}\label{1May2024_17:21}
 d_R(p_n,q_n)\geq n\,d(p_n,q_n).
\end{equation}

We shall show that we may suppose that $d_R(p_n,q_n)$ is bounded by the $\delta$ coming from \eqref{eq_bilip4may}. Indeed, for each $n\in N$, consider a curve $\gamma$ from $p_n$ to $q_n$ that is a geodesic with respect to $d$. 
By continuity, take points $a_0=p_n, a_1, a_2, \ldots, a_N=q_n$, for some $N\in \N$, along $\gamma$ such that
\begin{equation}\label{1May2024_17:31}
d_R(a_{j-1},a_j)<\delta, \qquad \forall j\in \{1, \ldots, N\}.
\end{equation}
 We claim that there exists an index $j\in \left\lbrace 1, \ldots, N \right\rbrace $ such that
\begin{equation}\label{4may1338}d_R(a_{j-1},a_j)\geq n\, d(a_{j-1},a_j).
\end{equation}
Indeed, if not, for every $j$, we would have 
$d_R(a_{j-1},a_j)< n\, d(a_{j-1},a_j)$; summing 
 and using that the points $a_j$'s are along a $d$-geodesic, we would have
$$d_R(p_n,q_n)\leq \sum_{j=1}^N d_R(a_{j-1},a_j)< n \sum_{j=1}^N d(a_{j-1},a_j)=n\,d(p_n,q_n).$$
We obtained a contradiction with \eqref{1May2024_17:21} and, thus, the claim is proved, which means that by replacing $p_n$ and $q_n$ with the just found $a_{j-1}$ and $a_{j}$, 
we have the bounds 
\begin{equation}\label{1May2024_17:33}
 d(p_n,q_n)\stackrel{\eqref{4may1338}}{\leq} \frac{1}{n}
 d_R(p_n,q_n)\stackrel{\eqref{1May2024_17:31}}{\leq} \frac{\delta}{n} 
, \qquad \forall n\in\N.
\end{equation}

By homogeneity of both distances, it is possible to assume that $q_n $ is always the point $1: =1_G$. 
Equation \eqref{1May2024_17:33} becomes
 \begin{equation}\label{2May2024_17:33}
 d(p_n, 1) \leq \frac{1}{n}
 d_R(p_n, 1) \leq \frac{\delta}{n} <\delta
, \qquad \forall n\in\N.
\end{equation}
Notice that, in particular, 
we have $ d(p_n, 1) \to 0$, i.e., the sequence $p_n$ converges to $1$.
Notice that, for the value $r$ from \eqref{4 May 1047}, the distance between 
$ B_d(1,r)$ and 
$M\setminus B_{d_R} (1, \delta)$ is positive. 
Therefore, for $n$ large enough, 
 then there exists $m_n\in \N$ such that 
\begin{equation}\label{eq_4may1150}(p_n)^{m_n}=\exp(m_n \log(p_m)) \in B_{d_R} (1, \delta)\setminus \bar B_d(1,r).\end{equation}
We stress that $\log( (p_n)^{m_n} )= m_n \log(p_n)$.
Now, using in addition the triangle inequality and the left invariance of the distance function, we bound
\begin{eqnarray*}
0<r&\stackrel{\eqref{eq_4may1150}}{<}&d(p_n^{m_n}, 1)
\\&\leq& \sum_{j=1}^{m_n} d(p_n^{j-1},p_n^{j})
\\&=& m_nd(p_n, 1)
\\&\stackrel{\eqref{2May2024_17:33}}{\leq} &\frac{1}{n}m_nd_R(p_n, 1)
\\&\stackrel{\eqref{eq_bilip4may}}{\leq} &\frac{k m_n}{n}\norm{\log(p_n)}
\\&= &\frac{k }{n}\norm{{m_n}\log(p_n)}
\\&= &\frac{k }{n}\norm{\log(p_n^{m_n})}
\\&\stackrel{\eqref{eq_bilip4may}}{\leq} &\frac{k^2 }{n}d_R (p_n^{m_n}, 1)
\stackrel{\eqref{eq_4may1150}}{\leq}\frac{k^2\delta}{n} \stackrel{n\to \infty}{\longrightarrow} 0.
\end{eqnarray*}
This contradiction proves the lemma.
\end{proof}

\subsection{Berestovskii's construction}
With the crucial help of Lemma \ref{rectgeo}, we will then explain the construction of the Berestovskii CC structure that appears in his Theorem \ref{Berestovskii}. 

\begin{proof}[Sketch of a proof of Theorem \ref{Berestovskii}]
From Lemma \ref{rectgeo}, we know that some dilation of the distance $d$ is greater than a Riemannian distance. 
So, rescaling the metric $d_R$, we assume to have $d_R\leq d $. 
We deduce that every curve $\gamma: I\to M$ that is $L$-Lipschitz with respect to $d$ is also $L$-Lipschitz with respect to $d_R$, since
\begin{equation}\label{10May2044}
d_R(\gamma (t_1), \gamma(t_2))\leq d(\gamma (t_1), \gamma(t_2)) \leq 
L |t_1-t_2|, \qquad \forall t_1, t_2\in I.
\end{equation}
%
Consequently, by Rademacher's Theorem we obtain that every $d$-rectifiable curve parametrized by a multiple of arc length is differentiable almost everywhere and absolutely continuous.
 


To prove Berestovskii's Theorem \ref{Berestovskii}, we need to find a sub-bundle $\Delta$ of the tangent bundle of $M$ and a norm on it. 
 As a result of \eqref{10May2044}, it makes sense to look at the set of velocities of $d$-rectifiable curves. 
 Denote by $\Lip([-1, 1];M)$ the collection of curves from $[-1, 1]$ to $M$ that are Lipschitz with respect to $d$.
We define the {\em Berestovskii bundle} as the subset of $TM$ given by\index{Berestovskii! -- bundle}
\begin{equation}\label{Berestovskii Delta}
\Delta: =\cup_{p\in M}\Delta_p\quad \text{ where } \quad \Delta_p: =\{\dot\gamma(0)\;:\; \gamma\in \Lip([-1, 1];M), \gamma(0)=p,
\exists \dot\gamma(0)
\}.
\end{equation}

One way to characterize the sub-Finsler norm that one needs to put on $\Delta$ is to describe the unit ball for the norm. One considers tangent vectors of curves that are Lipschitz with constant at most~1: \index{Berestovskii! -- norm}
$$\norm{v}\leq 1 \stackrel{\rm def}{\Longleftrightarrow} \exists \gamma: [-1, 1]\to M \text{ with } \gamma \, 1\text{-Lipschitz and } \dot\gamma(0)=v.$$

The fact that every fiber of Berestovskii's bundle is a vector subspace can be proved by constructing limits of zigzag curves obtained via translations of curves.
 Figure~\ref{fig_zigzag} pictures the infinitesimal construction: we begin with two Lipschitz curves $\gamma_1$, $\gamma_2$, giving tangent vectors $w_1$ and $w_2$, respectively. Given $\eps>0$ we follow $\gamma_1$ for time $\eps$, then we follow a translation of $\gamma_2$ for time $\eps$, then we follow a translation of $\gamma_1$ for time $\eps$, etc.
 As $\eps\to 0$, the limit curve will be rectifiable and will have a tangent equal to $\frac12(w_1+w_2)$ at $0$.
 Similarly, if $\gamma_1$ and $\gamma_2$ are 1-Lipschitz, then so are the zigzag curves and their limits. Hence, the defined norm is convex.
 \begin{figure}
 \centering
 \includegraphics[width=4.1in]{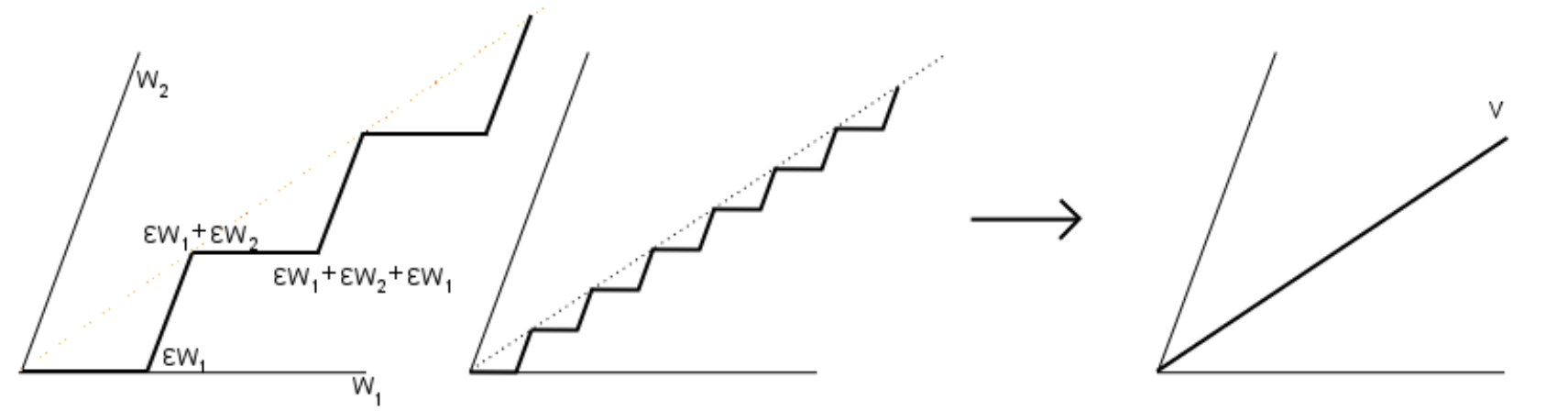}
	\caption{The sequence of $\eps$-zigzag curves converge (up to subsequences) to a curve whose tangent vector at the initial point is the average $v$ of the two tangent vectors $w_1$ and $w_2$ of the two initial curves.}
	\label{fig_zigzag}
\end{figure} 
 
 We stress that Berestovskii's distribution and norm are invariant by the (transitive) isometry group. In particular, the distribution is smooth and has a constant rank, and the norm is continuously varying. 
 We obtain a Carnot-Carath\'eodory space with length structure given by \eqref{def_Finsler_length} and distance $\dcc $ given by \eqref{dist_CC}. 

To obtain Berestovskii's result, one argues separately for the two inequalities: 
$\dcc \leq d$ and $d \leq \dcc .$
 The first inequality is straightforward: Given $p,q\in M$, let $\gamma$ be a $d$-geodesic from $p$ to $q$ that is parametrize by arc length with respect to $d$; so, for $\ell: =d(p,q)$, the curve $\gamma: [0, \ell ]\to M$ is 1-Lipschitz. By definition of the norm, we have $\norm{\dot \gamma(t)}\leq 1$.
Therefore, we have the bound:
$$\dcc (p,q) \stackrel{\rm def}{\leq} {\rm Length}_{\norm{\cdot}}(\gamma)\stackrel{\rm def}{=}\int_0^\ell \norm{\dot \gamma(t)}\dd t\leq \int_0^\ell1\dd t=\ell=d(p,q).$$

 The second inequality is more involved. Given a curve that is admissible for $\Delta$, one has to construct a $d$-rectifiable curve with almost the same length. This is done similarly with the zigzag method. We again refer to \cite{b}, but also suggest \cite[Secction~5.2]{LeDonne1}.
\end{proof}

The argument that we overviewed is quite flexible and can be used to study Lie coset spaces equipped with distances that are biLipschitz homogeneous.\index{bi-Lipschitz! -- homogeneous} 

\begin{theorem}[{\cite[Theorem 1.1]{LeDonne1}}]
 Let $M: =G/H$ be Lie coset space equipped with a geodesic admissible distance $d$.
 Suppose there exists a subgroup $G'$ of $G$ that acts transitively on $M$ and that acts by maps that are locally biLipschitz with respect to $d$. 
Then there exists a $G'$-invariant Carnot-Carath\'eodory structure on $M$ whose distance is locally biLipschitz equivalent to $d$.
\end{theorem}
 

\section{Exercises}
\begin{exercise}\label{ex_frame_left_inv_distribution} Let $\Delta $ be a left-invariant distribution on a Lie group. Then, there exists a global frame for $\Delta$ made of left-invariant vector fields. 
\end{exercise}

\begin{exercise}\label{ex_LIdistr_flag}
Let $\Delta$ be a left-invariant distribution on a Lie group $G$. 
For $\k\in \N$, let $\Delta^{[k]}$ be the $k$-th element in the flag of subbundles associated with $\Delta$ as in Definition~\ref{def:Delta_k}.
Then, the sequence of subspaces $V^{[k]}: =\Delta^{[k]}_{1_G}$, as $\k\in \N$, satisfies \eqref{LIdistr_flag}.
{\it Hint:} This is just the two interpretations of the Lie bracket for the Lie algebra.
\end{exercise}

\begin{exercise}\label{ex dim stabilizes for V}
Let $V$ be a subspace of a Lie algebra $\g$, with $n: =\dim \g$. Consider the spaces $V^{[k]}$ from \eqref{LIdistr_flag}

(i) If $\dim V \leq 1$, then $V^{[k]}=V$ for all $k\in \N$.

(ii) If $V$ is bracket generating and $V\neq \g$, then $\dim V \geq 2$, $n\geq3$, and $V^{[n-1]}=\g$.

(iii) If $V$ is not bracket generating and $n\geq3$, then $V^{[n-2]}=V^{[k]}$, for all $k\geq n-2$.
\end{exercise}

\begin{exercise}\label{ex_LI_norm}
Formula \eqref{cont_var_norm} defines a function that is left-invariant, continuous, and gives a continuously varying norm in the sense of Definition~\ref{def: continuously varying norm}.
\end{exercise}

\begin{exercise}\label{ex: surjective_for_CCgroups}
 Let $\varphi: G \to H$ be a Lie group homomorphism between sub-Finsler Lie groups $(G, \Delta^G, \norm{\cdot})$ and $(H, \Delta^H, \norm{\cdot})$.
Assume $\varphi_* (\Delta^G_1 ) \supseteq \Delta^H_1.$
 Then $\varphi: G \to H$ is surjective. 
\\{\it Hint.} The set
$\varphi_*(\g)$ is a Lie algebra containing the generating sub-space $\Delta^H_1$.
 \end{exercise}

\begin{exercise}\label{ex: Lip_for_CCgroups}
 Let $\varphi: G \to H$ be a Lie group homomorphism 
 between sub-Finsler Lie groups $(G, \Delta^G, \norm{\cdot})$ and $(H, \Delta^H, \norm{\cdot})$.
Assume $\varphi_* (\Delta^G_1) \subseteq \Delta^H_1.$
 Then $\varphi: G \to H$ is $L$-Lipschitz with respect to the respective sub-Finsler metrics, where $L$ is the Lipschitz constant of the linear map $$\varphi_*|_{\Delta^G_1}: (\Delta^G_1, \norm{\cdot}) \to (\Delta^H_1, \norm{\cdot}).$$
{\it Solution.} 
It is enough to observe that if 
$\gamma: [0, 1]\to G$ is a horizontal curve, 
then
\begin{eqnarray*}
\Length (\varphi\circ\gamma)
&=& \int_0^1 \norm{\dfrac{\dd}{\dd t}\left(\varphi(\gamma(t))\right)} \dd t\\
&=& \int_0^1 \norm{(L_{\varphi(\gamma(t))})^* \dfrac{\dd}{\dd t}\left(\varphi(\gamma(t))\right)}_{\Delta_1^H}\dd t\\
&=& \int_0^1 \norm{\varphi_*(\gamma'(t))}_{\Delta_1^H}\dd t\\
&\leq& L\int_0^1 \norm{\gamma'(t)}_{\Delta_1^G} \dd d\\
&=& L\int_0^1 \norm{\dot\gamma(t)}\dd t\\
 &=&L \Length (\gamma).
\end{eqnarray*}
\end{exercise}

\begin{exercise} Let $M_1$ and $M_2$ be sub-Riemannian manifolds with horizontal frames
$X_1, \ldots, X_m$ and $Y_1, \ldots, Y_m$, respectively.
Let $\pi: M_1\to M_2$ be smooth and surjective.
If $X_i$ is $\pi$-related to $Y_i$, for all $i\in \{1, \ldots, m\}$, then $\pi$ is a submetry. 
\end{exercise}

 \begin{exercise}\label{ex_braket_gen_quotients}
 Let $\pi:\g \to \mathfrak h$ be a surjective Lie algebra homomorphism. If $V\subseteq \g$ Lie generates $\g$, then $\pi(V)\subseteq \mathfrak h$ Lie generates $ \mathfrak h$.
 \end{exercise}

 \begin{exercise} 
Let $k\in \N$. Let $\pi: M_1\to M_2$ be a smooth map between sub-Finsler manifolds.
Assume that $M_1$ and $M_2$ are boundedly compact and their distributions have rank $k$.
 Assume $\pi$ surjective and such that 
 $$ (\dd \pi)_p: (\Delta_p, \norm{\cdot} ) \to (\Delta_{\pi(p)}, \norm{\cdot} ), \qquad \forall p\in M_1,$$
is a submetry, then $\pi$ is a submetry. 
$\skull$ What if the ranks are different?
\end{exercise}

\begin{exercise}\index{Heisenberg! -- group} 
 For the standard basis in the Heisenberg group, the map $E$ from Definition~\ref{def_E_for_BB} is 
$E({\bf{t}})=e^{t_1X} e^{t_2Y} e^{\sqrt{t_3}X} e^{\sqrt{t_3}Y} e^{-\sqrt{t_3}X}e^{-\sqrt{t_3}Y}.$
\end{exercise}
\begin{exercise}
For the standard basis in the Engel group (see \eqref{Engel_rel}), the map $E$ from Definition~\ref{def_E_for_BB} is 
{\small
\begin{eqnarray*}
E({\bf{t}})&=&
e^{t_2X} \; e^{t_2Y}\; e^{\sqrt{t_3}X} e^{\sqrt{t_3}Y} e^{-\sqrt{t_3}X}e^{-\sqrt{t_3}Y}\\&&\;\quad
e^{\sqrt[3]{t_4}X} \;e^{\sqrt[3]{t_4}Y} e^{-\sqrt[3]{t_4}X}e^{-\sqrt[3]{t_4}Y} e^{\sqrt[3]{t_4}X}
e^{\sqrt[3]{t_4}Y} e^{\sqrt[3]{t_4}X}e^{-\sqrt[3]{t_4}Y} e^{-\sqrt[3]{t_4}X}\;e^{-\sqrt[3]{t_4}X}.
\end{eqnarray*}
}
\end{exercise}

\begin{exercise}[Chow's theorem for sub-Finsler groups]\label{Chow4Carnot}\index{Chow Theorem! -- for sub-Finsler groups}
Every point $p\in G$ in a sub-Finsler group can be joined to the identity $1_G$ by a horizontal path.
Moreover, the CC distance induces the manifold topology. 
\\{\it Hint.} Use Corollary~\ref{corol_weak_BB}.
\end{exercise}
%


\begin{exercise}[Gr\"onwall Lemma (integral version)]\label{ex66b37cc8}\index{Gr\"onwall Lemma}
 Given $a<b$, let $\alpha, \beta,f: [a,b]\to\R$ be functions such that $\alpha$ is integrable, while $\beta$ and $f$ are continuous.
 Assume $\beta$ is non-negative and
 \[
 f(t) \leq \alpha(t) + \int_a^t \beta(s) f(s) \dd s .
 \]
 Then, for all $t\in[a,b]$,
 \[
 f(t) \leq \alpha(t) + \int_a^t \alpha(s) \beta(s) e^{\int_s^t\beta(r)\dd r} \dd s .
 \]
 If, moreover, the function $\alpha$ is non-decreasing, then
 \[
 f(t) \leq \alpha(t) e^{\int_a^t\beta(s)\dd s} .
 \]
\end{exercise}

\begin{exercise}\label{ex66b36bce}
 Let $u,v: [0, 1]\to \mathfrak{gl}(n)$ be measurable functions.
 For $\eps>0$, denote by $\gamma_{u+\eps v}$ the solution $[0, 1]\to\mathfrak{gl}(n)$
 of the Cauchy Problem
 \[
 \begin{cases}
 \dot\gamma_{u+\eps v}(t) &= \gamma_{u+\eps v}(t) \cdot(u(t)+\eps v(t))
 \qquad\forall t\in[0, 1], \\
 \gamma_{u+\eps v}(0) &= \Id_n .
 \end{cases}
 \]
 Let $\sigma: [0, 1]\to\mathfrak{gl}(n)$ be the solution to
 \[
 \begin{cases}
 \dot\sigma(t) &= \gamma_u(t)\cdot v(t) + \sigma(t)\cdot u(t) \\
 \sigma(0) &= 0 .
 \end{cases}
 \]
 Then
 \[
 \lim_{\eps\to0} \frac1\eps (\gamma_{u+\eps v}(t) - \gamma_u(t)) = \sigma(t),
 \qquad \forall t\in[0, 1] .
 \]
 {\it Hint:}
 Set $\eta_\eps: = \frac1\eps (\gamma_{u+\eps v}(t) - \gamma_u(t))$ and then
 use Gr\"ownwall Lemma from Exercise~\ref{ex66b37cc8} to show
 \[
 \left\| \sigma(t) - \eta_\eps(t) \right\|
 \leq 
 \eps 
 \alpha(t) e^{\int_a^t\beta(s)\dd s},
 \]
 where $\alpha(t): = \int_0^t\|\eta_\eps(s)\cdot v(s)\| \dd s$ and 
 $\beta(t): = \int_0^t \|u(s)\| \dd s$.
\end{exercise}

\begin{exercise}\label{ex66b38412}
 Given a Lie group $G$, there is a matrix group $H$ with the same universal covering group of $G$.
\\
 {\it Hint:} Use Ado's Theorem, or Birkhoff Theorem if $G$ is nilpotent.
\end{exercise}

\begin{exercise}\label{ex66b383da}
 Let $\pi: (\tilde G, \tilde V)\to (G, V)$ be a homomorphism of polarized Lie groups.
 Then
 \[
 \pi \circ \End_{(\tilde G, \tilde V)} = \End_{(G, V)}\circ \pi_*,
 \]
 where $\End_{(G, V)}$ denotes the end-point map for the polarized group $(G, V)$
 and $\pi_*:L^2([0, 1];\tilde V) \to L^2([0, 1]; V)$ is induced by the composition $\pi_*(u) := \pi_* \circ u$.
 \\ {\it Hint:} The curves $\pi(\gamma_u(t))$ and $\gamma_{\pi_*u}(t)$ satisfy the same ODE.
\end{exercise}

\begin{exercise}\label{ex7Aug2024}
 \begin{description}
 \item[i.]
 Let $\pi: (\tilde G, \tilde V)\to (G, V)$ be a homomorphism between polarized Lie groups
 such that $\pi_*:  \tilde{\g} \to  \g $ is a Lie algebra isomorphism and $ \pi_*  ( \tilde V) = V $.
Then, Proposition~\ref{prop:dif_end} holds for $\tilde G$ if and only if it does for $G$.
\\
 {\it Hint:} Use Exercise~\ref{ex66b383da}.
 \item[\eqref{ex7Aug2024}.ii.]
 We deduce a proof for Proposition~\ref{prop:dif_end}, which, in this text, has been proved only for matrix groups.
\\
 {\it Hint:} Use Exercise~\ref{ex66b38412}.
 \end{description}
\end{exercise}

	\begin{exercise}\label{ex_prod_polarized}
	If $(G_1, V_1)$ and $(G_2, V_2)$ are polarized groups, then 
	$(G_1 \times G_2, V_1 \times V_2) $ is a polarized group. Moreover, the polarization $V_1 \times V_2$ is bracket generating if and only if both $V_1$ and $ V_2$ are.
	\end{exercise}
	
\begin{exercise}\label{ex_prod_abn}
For 		polarized groups $G_1 $ and $G_2 $ consider the polarized group $G_1 \times G_2$ as in Exercise~\ref{ex_prod_polarized}.
 A curve $\gamma=(\gamma_1, \gamma_2): I\to G_1 \times G_2$ is abnormal if and only if at least one between $\gamma_1$ and $ \gamma_2$ is abnormal.
 \\{\it Hint.} The differential of the endpoint map splits into a block linear transformation that is block diagonal.
\end{exercise}

\begin{exercise}\label{ex_const_abn}
Let $(G, V)$ be a polarized group with $V\neq \Lie(G)$. Then, the constant curve $1_G$ is abnormal in $(G, V)$.	
\end{exercise}

\begin{exercise}\label{ex_heis_abn}
In the sub-Riemannian Heisenberg group, abnormal curves are constant.		
\end{exercise}

\begin{exercise}\label{ex_nonsmooth_abn}
Let $H$ be the sub-Riemannian Heisenberg group. Then, every absolutely continuous curve $t\mapsto (1_H, x(t), y(t) ) \in H\times \R^2$ is abnormal. Give an injective non-smooth example parametrized by arclength.
 \\{\it Hint.} Use Exercise~\ref{ex_prod_abn} and Exercise~\ref{ex_const_abn}.
\end{exercise}


\begin{exercise}\label{ex differential of the energy}
On the Hilbert space $L^2([0, 1];\R^n)$ with scalar product $\langle u,v \rangle: = \int_0^1u(t) \cdot v(t) \dd t$, consider the energy function given by 
$u\mapsto
\frac{1}{2}\|u\|^2: = \frac{1}{2} \langle u,u \rangle$.
This function is smooth and its differential at $u\in L^2(0, 1)$ is
$v\mapsto \langle u,v \rangle$.\index{energy! -- of a control}
\end{exercise}

\begin{exercise} 
In Riemannian Lie groups, length minimizers are reparametrizations of normal curves. 
\end{exercise}

\begin{exercise}\label{ex analytic normal}
Equation \eqref{geodesic} is an analytic ODE and its solutions are analytic.
 	\\{\it Hint.} Solutions of analytic ODEs are analytic; see \cite[pp. 121-128]{Birkhoff_Rota}.
\end{exercise}
\begin{exercise}
	In contact structures, as for example in every 3D sub-Riemannian Lie group, every abnormal curve is constant.
	\\{\it Hint.} See Proposition~\ref{prop_rank-2_abnormal}.
\end{exercise}
\begin{exercise}
	In every sub-Riemannian manifold of step $2$, every energy minimizer is normal. {\it Hint.} Use Goh result.
\end{exercise}

\begin{exercise}[Goh condition for rank-2 polarizations]\index{Goh! -- condition}
Let $G$ be a Lie group equipped with a rank-2 polarization $\Delta$. Let $\gamma: [0, 1]\to G$ be a $\Delta$-horizontal curve. Then, we have a stronger statement than Theorem~\ref{Goh_thm}: $ \left\{ \Ad_{\gamma(t)} ( [\Delta_1, \Delta_1]): t\in [0, 1] \right\} \subseteq \Span \left\{ \Ad_{\gamma(t)} (\Delta_1 ): t\in [0, 1]\right\}$.
\end{exercise}

\begin{exercise}\label{ex_derivata_Lambda}
Let $\gamma: I\to G$ be a curve in a Lie group. If $\lambda\in \g^*$ and $Y\in \g$, then
\begin{equation}\label{derivata_Lambda}
\frac{\dd }{\dd t}\left(\lambda \Ad_{\gamma(t)}(Y) \right)
=
 \lambda\Ad_{\gamma(t)} [\gamma'(t), Y], \qquad \forall t\in I.
\end{equation}
{\it Solution}.
\begin{eqnarray*}
\nonumber \frac{\dd }{\dd t}\left(\lambda \Ad_{\gamma(t)}(Y) \right)&=& 
\left.\frac{\dd }{\dd s}\left(\lambda \Ad_{\gamma(t+s)}(Y) \right)\right|_{s=0}\\
\nonumber&=& 
\left.\frac{\dd }{\dd s}\left(\lambda\Ad_{\gamma(t)} \Ad_{\gamma(t)^{-1}\gamma(t+s)}(Y) \right)\right|_{s=0}\\
\nonumber&=& 
 \lambda\Ad_{\gamma(t)} \ad_{ \left(\dd L_{\gamma(t)}\right)^{-1} \gamma(t)^{-1}\dot\gamma(t)}(Y) \\
 &=& 
 \lambda\Ad_{\gamma(t)} [\gamma'(t), Y], 
\end{eqnarray*}
where we used that $\Ad_{gh}=\Ad_g \Ad_h$, that $\lambda $ and $\Ad_g$ are linear. 
\end{exercise}

\begin{exercise}[First derivative of the extremal equations]\label{ex:First derivative of the extremal equations}\index{first derivative of the extremal equations}
Let $\gamma: I\to G$ be a horizontal curve in a polarized group $(G, V)$. Let $u: =\gamma'$.
If $\gamma$ is an abnormal curve with covector $\lambda\in \g^*\setminus\{0\}$, then 
\begin{equation}\label{equazione_abnormale}0= \lambda\Ad_{\gamma(t)} [u(t), X], \qquad \forall t\in I, \forall X\in V.
\end{equation} 
If $V$ is equipped with a scalar product, $e_1, \ldots, e_r$ are an orthonormal bases, $u = u_1 e_1+\ldots+ u_r e_r$,
and $\gamma$ is a normal curve with covector $\lambda\in \g^*$, then 
\begin{equation}\label{equazione_normale}
\dot u_i = \lambda\Ad_{\gamma(t)} [u(t),e_i], \qquad \forall t\in I, \forall i\in\{1, \ldots, r\}.\end{equation}
{\it Hint.} Take the derivative of the abnormal equation \eqref{abnormal1} and the normal equation \eqref{geodesic}, using Exercise~\ref{ex_derivata_Lambda}.
\end{exercise}

\begin{exercise}\label{ex_oriented_curvature}
 Let $\sigma: [0, 1]\to \R^2$ be a planar curve with $ u: =\sigma'$ of class $C^2$ and with never-vanishing speed.
The oriented curvature of $\sigma$ is 
$\kappa = \frac{ \sigma'_1 \sigma''_2- \sigma'_2 \sigma''_1 }{\|\sigma' \|^3} =
 \frac{ u_1 \dot u_2- u_2 \dot u_1 }{\|u\|^3} .$
\\ {\it Hint.} See \cite[Equation (1.11)]{Abate_Tovena_curves_surfaces}.
 \end{exercise}

\begin{exercise}\label{ex ODE all'ODE}
        The ODE~\eqref{curvature_normal} is of the form
        \[
        \begin{cases}
        F(t,\gamma(t),\gamma'(t),\gamma''(t)) = 0 ,\\
        \gamma(0) = 1_G ,\\
        \gamma'(0) = v_0 ,
        \end{cases}
        \]
        where $v_0 = \lambda(e_1) e_1 + \lambda(e_2) e_2 \neq 0$
        and $F:I\times G\times V\times V \to \R$ is
        \[
        F(t,p,v,w) = \lambda(\Ad_p e_{12}) - \frac{v_1w_2 - v_2w_1}{\|v\|^2} .
        \]
        The derivative of $F$ in $w$ is not zero at points $(t,p,v,w)$ with $v\neq0$.
        Therefore, for every $v_0\in V\setminus\{0\}$
        there are $\Omega\subset I\times G\times V$ neighborhood of $(0,1_G,v)$,
        and a smooth function $\tilde F: \Omega \to V$
        such that
         the ODE~\eqref{curvature_normal} can be written as
        \[
        \begin{cases}
        \gamma''(t) = \tilde F(t,\gamma(t),\gamma'(t)),\\
        \gamma(0) = 1_G ,\\
        \gamma'(0) = v_0 .
        \end{cases}
        \]
        We conclude that solutions to the ODE~\eqref{curvature_normal} are unique and smooth.
 \end{exercise}

\begin{exercise}\label{ex_5Jun2024}
On the group $G_q$ as in Definition~\ref{def_step_2_G}, for $(v_1,w_1),(v_2,w_2)\in V\times W$ we have that
$t\mapsto (t v_1, t w_1) $ is an OPS and 
$$
[(v_1,w_1),(v_2,w_2)]= 
\left.\frac{1}{2}\frac{\dd^2}{\dd t^2}
(t v_1, t w_1) (t v_2, t w_2) (-t v_1, -t w_1) (-t v_2, -t w_2) \right|_{t=0} = (0, q(v_1,v_2))
.$$
{\it Hint.} Recall Definition~\ref{def:Lie_bracket_vector_fields}.d.
\end{exercise}

\begin{exercise}\label{ex OPS geodesic 2}
On the group $G_q$ as in Definition~\ref{def_step_2_G},
 let $v\in V $ and consider a sub-Finsler norm on $G_q$ such that $|v|=1$.
Then, the line  $t\mapsto  t v   $ is a geodesic.
\\{\it Hint.} See Proposition~\ref{OPS geodesic}.
\end{exercise}

\begin{exercise}\label{ex_sub_Gq}
Every connected subgroup of a group $G_q$ as in Definition~\ref{def_step_2_G}, is of the form $G_{q'}$ for some skew-symmetric bilinear map $q':V'\times V'\to W'$.
\end{exercise}

\begin{exercise}
 Let $G $ be a sub-Riemannian nilpotent group of step two with
Lie algebra $\g$ and polarization $\Delta\subset\g$.
Define $V_2: = [\g, \g]$ and $V_1\subset\Delta$ as the orthogonal of $V_2\cap \Delta$ in $\Delta$.
Then, the polarized group $(G, \Delta_1)$ is isomorphic to some $G_q$ as in Definition~\ref{def_step_2_G}, with $q:V_1\times V_1\to V_2$ being the Lie bracket.
\end{exercise}

\begin{exercise}\label{Ex Nalon}
Let    $G=G_q=V\times W$ be a step-two nilpotent group as in Definition~\ref{def_step_2_G}. 
\\
{\ref{Ex Nalon}.i.} For $\lambda>0$, the map $\delta_\lambda: (v,w)\in V\times W \mapsto (\lambda v, \lambda^2 w)$ is a group automorphism of $V\times W$.
\\
{\ref{Ex Nalon}.ii.} Fixing a norm on $V$, the sub-Finsler distance satisfies 
$d_{\rm sF} ( \delta_\lambda(v+w), 0 ) = \lambda d_{\rm sF} ( v+w ,0  ),$  for all $v\in V$ and $ w\in W$.
\\
{\ref{Ex Nalon}.iii.} If $d_{\rm sF} ' $ is another sub-Finsler distance with polarization $\Delta$ such that $V\subseteq \Delta_1$, then for some $C$ we have
$d_{\rm sF} ' \leq C d_{\rm sF} $.
In particular, after fixing a norm on $W$, we have 
\begin{equation}\label{Ex Nalon eq}
d_{\rm sF} ' (0,w) \leq C \sqrt{w}  , \qquad \forall w\in W,
\end{equation}
for some $C>0$.
\end{exercise}


\begin{exercise}
Let $k>1$. Let $d$ be an admissible distance on $\R^2$ that is geodesic and has the property that the translations are $k$-biLipschitz maps. Then $d$ is biLipschitz equivalent to the Euclidean distance.
\\{\it Hint.} Start by following the argument in the proof of Lemma \ref{rectgeo}.
\end{exercise}
 
\begin{exercise}[Geodesic distances on $\R^2$]\index{geodesic! -- metric} 
The only isometrically homogeneous geodesic metric spaces topologically equivalent to a plane are the Euclidean 2-space and the hyperbolic plane, equipped with left-invariant Finsler metrics.
\\ {\it Hint.} Use Berestovskii Theorem~\ref{Berestovskii}, pass to a sub-Riemannian metric, and deduce that the curvature is constant.
\end{exercise}


\chapter{Riemannian Lie groups}\label{ch_RiemaLie}

In this chapter, we discuss a very classical family of sub-Finsler Lie groups: those where the polarization is the entire tangent bundle and the norm comes from a scalar product. We refer to them as {\em Riemannian Lie groups}, but a more precise name could be Lie groups equipped with left-invariant Riemannian structures.\index{Riemannian! -- Lie group}

The results that we present in Section~\ref{sec Left-invariant Riemannian metrics} are classical, and a reference is \cite{Milnor}. 
In Section~\ref{sec Isometries of metric groups}, we discuss the fact that isometry groups of metric Lie groups are subgroups of isometry groups of Riemannian Lie groups. Here, the main source is \cite{Kivioja_LeDonne_isom_nilpotent}.
Some other expository readings are \cite{Abate_Tovena_geometria_differenziale, Purho_master, Sopio_master}.

\section{Left-invariant Riemannian metrics}\label{sec Left-invariant Riemannian metrics}

In this section, we extensively use various notations from Lie group theory, as we recalled in Chapter~\ref{ch_LieGroups}.
For example, the maps $L_h$ and $R_h$ are the left translation and the right translation, respectively, by a group element $h$. We shall use $\Ad$ and $\ad$ from Section~\ref{sec Ad and ad}, and the notion of structural constants as in \eqref{eq1655}.
We will also discuss notions from Riemannian geometry.
In addition to what we presented in Section~\ref{sec_diff_geom}, we point out the books \cite{Lee_MR2954043, Lee_Riemannian}.
The notion of a Riemannian metric, also called Riemannian metric tensor, is introduced in Section~\ref{Riemannian and Finsler geometry}.

\begin{definition}[Left-invariant Riemannian metric]
 	A Riemannian metric $\langle \cdot, \cdot \rangle_{\cdot}$ on a Lie group $G$ is said to be \emph{left-invariant} if \index{left-invariant! -- Riemannian metric}
	\[
	\langle (L_h)_*v,(L_h)_*w \rangle_{hg} = \langle v,w \rangle_g, \qquad \forall g,h\in G, \forall v,w\in T_gG.
	\]
Similarly, we say that $\langle \cdot, \cdot \rangle_{\cdot}$ is \emph{right-invariant} if \index{right-invariant! -- Riemannian metric}
	\[
	\langle (R_h)_*v,(R_h)_*w \rangle_{gh} = \langle v,w \rangle_g, \qquad \forall g,h\in G, \forall v,w\in T_gG.
	\]
It is \emph{bi-invariant} if it is left-invariant and right-invariant.
\index{bi-invariant! -- Riemannian metric}
\end{definition}

We shall see various equivalent characterizations for those Riemannian metrics that are both left-invariant and right-invariant. One of them is the property that for all $Z\in\g$, the adjoint map
 $\ad_Z$ is a {\em skew-adjoint} transformation of $(T_{1_G}G, \langle \cdot, \cdot \rangle_{1_G})$,
 i.e., it is antisymmetric in the sense that 
		\index{skew-adjoint}\index{antisymmetric}\index{adjoint! -- map, $\ad$}
		\[
		\langle \ad_Z X,Y \rangle = -\langle X, \ad_ZY \rangle, \qquad \forall X,Y,Z\in\g.
		\]

\begin{theorem}\label{teo2136} 
\index{isometry}\index{adjoint! -- representation, $\Ad$}
\index{adjoint! -- map, $\ad$}
 	Let $G$ be a connected Lie group with Lie algebra $\g$ and with a left-invariant Riemannian metric $\langle \cdot, \cdot \rangle$.
	The following are equivalent:
	\begin{description}
	\item[\ref{teo2136}.i.] 	$\langle \cdot, \cdot \rangle$ is right invariant;
	\item[\ref{teo2136}.ii.] 	$\Ad_g$ is an isometry, for all $g\in G$;
	\item[\ref{teo2136}.iii.] 	$\ad_X$ is skew-adjoint, for all $X\in\g$.
	\end{description}
The connectedness of $G$ is required only for the implication $(iii)\THEN(ii)$.
\end{theorem}
\begin{proof}
	$(i)\IFF(ii)$: 
	Since $\langle \cdot, \cdot \rangle$ is left-invariant, item $(i)$ is equivalent to $C_g=R_{g^{-1}}\circ L_g$ being isometries, for all $g\in G$, which is further equivalent to $(\Ad)_g = (\dd C_g)_{1_G}$ being isometries for all $g$.
	
	$(ii)\THEN(iii)$:
	By Formula~\ref{Formula: Ad: ad} we have $\Ad_{\exp(X)} = e^{\ad_X}$, so $e^{\ad(tX)}=e^{t\ad X}$ is an isometry, for all $t\in\R$ and all $X\in\g$. 
	Recall that $\frac{\dd}{\dd t} e^{tA}|_{t=0}=A$
	from Proposition~\ref{prop: der: eta}
	 and take the derivative at $t=0$ of the identity
	\[
	\langle e^{\ad(tX)}Y,e^{\ad(tX)}Z \rangle_{1_G} = \langle Y, Z \rangle_{1_G}, \qquad \forall X,Y,Z\in\g. 
	\]
	We get
	\[
	\langle \ad_XY, Z \rangle_{1_G} + \langle Y, \ad_XZ \rangle_{1_G} = 0, \qquad \forall X,Y,Z\in\g.
	\]
	which is $(iii)$.
	
	$(iii)\THEN(ii)$:
	Recall that $\frac{\dd}{\dd t}e^{tA} = Ae^{tA}$, again from Proposition~\ref{prop: der: eta}.
	We calculate
	\begin{multline*}
	 	\frac{\dd}{\dd t}\langle e^{\ad(tX)}Y,e^{\ad(tX)}Z \rangle_{1_G}
		= \langle \ad(X)e^{\ad(tX)}Y, e^{\ad(tX)}Z \rangle_{1_G} +\\
			+ \langle e^{\ad(tX)}Y, \ad(X) e^{\ad(tX)}Z \rangle_{1_G}
		\overset{(iii)}= 0.
	\end{multline*}
	Hence, the function $t\mapsto \langle e^{\ad(tX)}Y,e^{\ad(tX)}Z \rangle_{1_G}$ is constant.
	Evaluating it at $t=0$ and $t=1$, we deduce that $e^{\ad_X}$ is an isometry.
	Hence, the map $\Ad_{\exp(X)}$ is an isometry, for all $X\in\g$.
	
	So $\Ad_g$ is an isometry, for all $g$ in a neighborhood $U$ of $1_G$ in $G$.
	Since, when $G$ is connected, every element in $G $ is a finite product $g=g_1\cdot\dots\cdot g_k$ of elements $g_1, \ldots, g_k\in U$ (see Exercise~\ref{Prop:generating}), then $\Ad_g=\Ad_{g_1}\circ\dots\circ\Ad_{g_k}$ is an isometry.
	
	We used the connectedness assumption of $G$ only to prove the implication $(iii)\THEN(ii)$.
	Without this assumption, there are counterexamples; see Exercise~\ref{Seba 4 oct 2024}.
\end{proof}

\subsection{Connections and geodesics on Lie groups}
\index{linear connection}
Recall from Section~\ref{Sec:vect_fields_brackets} that, given a manifold $M$, we denote by $\Vec(M)$ the space of smooth vector fields on $M$.

A \emph{linear connection} $\nabla$ on a manifold $M$ is a map 
\[
\begin{array}{rccc}
 	\nabla:&\Vec(M)\times\Vec(M) &\to &\Vec(M) \\
		&(X,Y) &\mapsto &\nabla_XY
\end{array}
\]
that is $C^\infty(M)$-linear in $X$, $\R$-linear in $Y$, and satisfies the {\em Leibniz rule}:
\[
\nabla_X(fY) = (Xf) Y + f\nabla_XY ,
\qquad \forall f\in C^\infty(M) .
\]
\begin{definition}[Left-invariant linear connection]
\index{linear connection! left-invariant --, $\nabla$}
\index{left-invariant! -- linear connection, $\nabla$}
	Let $G$ be a Lie group. 
	A linear connection $\nabla$ on $G$ is \emph{left invariant} if 
	\[
	(L_g)_*\nabla_XY = \nabla_{(L_g)_*X} (L_g)_*Y, \qquad \forall g\in G, \forall X,Y\in\Vec(G).
	\]
\end{definition}

On every Riemannian manifold, there is a unique linear connection that is compatible with the metric and is torsion-free; for its construction and properties, see \cite[Chapter~5]{Lee_Riemannian}. This connection is called {\em Levi-Civita connection}.
The Levi-Civita connection on a Riemannian manifold $M$ satisfies the {\em Koszul formula}: for all $X,Y,Z\in \Vec(M)$ 
	\begin{multline}\label{Koszul-formula}
	 	\langle \nabla_XY,Z \rangle
		= \frac12 \Big(
		X\langle Y,Z \rangle + Y\langle Z, X \rangle - Z\langle X,Y \rangle 
		+\langle [X,Y],Z \rangle	+ \langle [Z,Y], X \rangle + \langle [Z, X],Y \rangle
		\Big).
	\end{multline}

\begin{proposition}\label{prop0647}
 	Let $G$ be a Lie group with Lie algebra $\g$ of left-invariant vector fields.
	There is a one-to-one correspondence between the set of left-invariant linear connections $\nabla$ on $G$ and the set $\Mult(\g, \g;\g)$ of bilinear functions $\alpha:\g\times\g\to\g$ given by
	\[
	\alpha_{\nabla}(X,Y) = \nabla_XY, \qquad \forall X,Y\in\g.
	\]
\end{proposition}
\begin{proof} The result is clear once we notice that, fixing a frame $X_1, \dots, X_n$ of left-invariant vector fields on $G$, we can write every arbitrary pair of vector fields on $G$ as $\sum_ja^jX_j$ and $\sum_ib^iX_i $, for some $a^j, b^j\in C^\infty(G)$. Then, the connection is determined:
\[
 \nabla_{a^jX_j}(b^iX_i)
= a^jb^i \nabla_{X_j}X_i + a^j (X_jb^i) X_i
= a^jb^i \alpha_\nabla(X_j, X_i) + a^j (X_jb^i) X_i. 
\]
\end{proof}
\begin{proposition}\label{prop0640}
\index{one-parameter subgroup}
\index{geodesic}
 	Let $G$ be a Lie group with a left-invariant linear connection $\nabla$, and let $X$ be a left-invariant vector field.
	The following are equivalent:
	\begin{description}
	\item[\ref{prop0640}.i.] 	$\alpha_\nabla (X, X) = 0$ (i.e., $\nabla_XX=0$);
	\item[\ref{prop0640}.ii.] 	the one-parameter subgroup $t\mapsto\Phi^t_X({1_G})$ is a geodesic with respect to $\nabla$.
	\end{description}
\end{proposition}
\begin{proof}
 	The curve $t\mapsto \gamma(t):=\Phi^t_X({1_G})$ has derivative $\gamma'(t)=X_{\gamma(t)}$.
	By definition, the curve $\gamma$ is a geodesic with respect to $\nabla$ if and only if $\nabla_{\gamma'}\gamma'=0$, thus if and only if $(\nabla_XX)_\gamma \equiv0$.
	Since $\nabla_XX$ is a left-invariant vector field, then $\gamma$ is a geodesic if and only if $(\nabla_XX)_{1_G}=0$, if and only if $\nabla_XX=0$.
\end{proof}

For another characterization of when a one-parameter subgroup is a Riemannian geodesic, see also Corollary~\ref{Riemannian_geodesic_OPS}

\begin{example}
	Let $G$ be a Lie group and $c\in\R$.
	Then the map 
	\[
	(X,Y)\mapsto c[X,Y]
	\]
	is in $\Mult(\g, \g;\g)$.
	Hence, by Proposition~\ref{prop0647}, it induces a left-invariant linear connection on $G$. Notice that for this connection, the Christoffel symbols $\Gamma_{ij}^k$ with respect to a frame of left-invariant vector fields are precisely the structural constants $c_{ij}^k$ (see the definition in \eqref{eq1655}) with respect to the same frame, multiplied by $c$.\index{structural constants}
\end{example}

\begin{lemma}\label{lem:avrverafqwde}\index{Levi-Civita connection}
 	Let $\langle \cdot, \cdot \rangle$ be a left-invariant Riemannian metric on a Lie group $G$, and let $\nabla$ be the associated Levi-Civita connection. 
	\begin{description}
	\item[\ref{lem:avrverafqwde}.i] 	For all left-invariant vector fields $X,Y,Z$ on $G$, we have
	\begin{equation}\label{eq200819}
	\langle \nabla_XY,Z \rangle 
	= \frac12 \Big(
		\langle [X,Y],Z \rangle	+ \langle [Z,Y], X \rangle + \langle [Z, X],Y \rangle
		\Big).
	\end{equation}
	\item[\ref{lem:avrverafqwde}.ii] 	If $X_1, \dots, X_n$ are orthonormal left-invariant vector fields that form a basis of $\Lie(G)$ and $\alpha^k_{ij}$ are the corresponding structural constants, then
	\begin{align}
	 	\langle [X_i, X_j], X_k \rangle &= \alpha_{ij}^k , \\
		\langle \nabla_{X_i}X_j, X_k \rangle &= \frac12 \left(\alpha_{ij}^k - \alpha_{jk}^i +\alpha_{ki}^j\right), \\
		\nabla_{X_i}X_j &= \frac12 \sum_{k=1}^n \left(\alpha_{ij}^k - \alpha_{jk}^i+\alpha_{ki}^j\right) X_k.
	\end{align}
	\end{description}
\end{lemma}
\begin{proof}
Regarding 	 \ref{lem:avrverafqwde}.i, 	recall that the Levi-Civita connection satisfies Koszul Formula \eqref{Koszul-formula}.
	Note that if $X, Y$ are left-invariant vector fields, then their scalar product $\langle X, Y \rangle$ is constant along $G$.
	Hence, Koszul Formula simplifies to \eqref{eq200819}.
	%
	
Regarding \ref{lem:avrverafqwde}.ii, 	the structural constants are defined by the equation $[X_i, X_j]=\sum_k\alpha_{ij}^kX_k$. This implies $\langle [X_i, X_j], X_k \rangle = \sum_h\alpha_{ij}^h\langle X_h, X_k \rangle = \alpha_{ij}^k$, because the $X_i$'s are orthonormal.
	From 	 \ref{lem:avrverafqwde}.i, the rest follows. \qedhere
	 
\end{proof}

\begin{corollary}\label{Riemannian_geodesic_OPS}
Let $G$ be a Lie group endowed with a left-invariant Riemannian metric.
Then, the one-parameter subgroup in the direction of $X$ is geodesic if and only if $X$ is orthogonal to $[X, \g]$.
\end{corollary}
\begin{proof}
It follows from \eqref{eq200819} and Proposition~\ref{prop0640}.
\end{proof}

\begin{theorem}\label{thm2111}\index{one-parameter subgroup}
\index{geodesic}
\index{Levi-Civita connection}
 	Let $\langle \cdot, \cdot \rangle$ be a left-invariant Riemannian metric on a Lie group $G$ and $\nabla$ be the associated Levi-Civita connection.
	The following are equivalent:
	\begin{description}
			\item[\ref{thm2111}.i.] 
			The group exponential map $\exp$ coincides with the Riemannian exponential map 	$\exp_{1_G}$,	 
	i.e., the family of one-parameter subgroups is exactly the family of the geodesics from the identity element $1_G$;
	\item[\ref{thm2111}.ii.] 		If $X,Y\in\g$,	 then $\alpha_\nabla(X,Y) = \frac12[X,Y]$, i.e., 
		\[
		\nabla_XY=\frac12[X,Y], \qquad \forall X,Y\in\g;
		\] 
	\item[\ref{thm2111}.iii.] 	The map $\ad_Z$ is skew-adjoint, for all $Z\in \g$. 
		\index{skew-adjoint}\index{antisymmetric}\index{adjoint! -- map, $\ad$}
	\end{description}
\end{theorem}
\begin{proof}
	Note that the formula \eqref{eq200819} in Lemma~\ref{lem:avrverafqwde} can be written also as
	\begin{align*}
	\langle \nabla_XY,Z \rangle
		&= \frac12 \Big(
		\langle [X,Y],Z \rangle	+ \langle [Z,Y], X \rangle + \langle [Z, X],Y \rangle
		\Big) \\
		&= \frac12 \Big(
		\langle [X,Y],Z \rangle	+ \langle \ad_ZY, X \rangle + \langle \ad_ZX,Y \rangle
		\Big).
	\end{align*}
	The equivalence $(ii)\IFF(iii)$ easily follows from the last equality.
	Moreover, for $X=Y$ we get
	\[
	\langle \nabla_X X, Z \rangle 
	= \langle \ad_Z X, X \rangle ,
	\qquad\forall X, Z\in\g .
	\]
	Hence, the point $(i)$, which by Proposition~\ref{prop0640} is equivalent to $\nabla_XX=0$ for all $X\in\g$, is also equivalent to $\langle \ad_ZX, X \rangle=0$, for all $X,Z\in\g$, which (by an easy computation) is equivalent to $\ad_Z$ being skew-adjoint.
\end{proof}
\begin{remark}
 	Note that the
	equation $\tilde\nabla_XY:=\frac12[X,Y]$, for $X,Y\in \g$ always define a left-invariant connection $\tilde\nabla$ on $G$.
	Condition $(ii)$ of Theorem~\ref{thm2111} is that such a $\tilde\nabla$ is the Levi-Civita connection $ \nabla$.
\end{remark}


\subsection{Curvatures of left-invariant metrics}
\index{curvature}
\index{left-invariant! -- Riemannian metric}
\index{structural constants}
\index{Riemannian! -- curvature tensor}
\index{sectional curvature}
 	For the following discussion, be aware that our convention for the {\em Riemannian curvature tensor} is
	\begin{equation}\label{def_Riemannian_curvature_tensor}
	R(X,Y, \cdot) := \nabla_X\nabla_Y - \nabla_Y\nabla_X - \nabla_{[X,Y]} .
	\end{equation}
		Also, recall that 
		given two linearly independent tangent vectors $X,Y$ at the same point of a Riemannian manifold, the {\em sectional curvature} of their spanned plane is
	\begin{equation}\label{def_sectional_curvature}\index{curvature! sectional --}
	{\rm Sec}(X,Y) := \frac{ \langle R(X,Y,Y), X \rangle }{\|X\|^2\|Y\|^2 - \langle X,Y \rangle^2 } .
	\end{equation}
	
\begin{proposition}\label{Prop_sectional_curvature_structural_constants}
 	Let $G$ be a Lie group equipped with a left-invariant Riemannian metric.
	Let $X_1, \dots, X_n$ be orthonormal left-invariant vector fields that form a basis of $\Lie(G)$, and let $\alpha_{ij}^k$ be the corresponding structural constants.
	Then, the Riemannian curvature tensor satisfies
	\begin{align*}
	 	R(X_i, X_j, X_k) 
		= \sum_{\ell,h=1}^n \bigg[
		&\frac14 
			\left(\alpha_{jk}^\ell-\alpha_{k\ell}^j+\alpha_{\ell j}^k\right)
			\left(\alpha_{i\ell}^h-\alpha_{\ell h}^i+\alpha_{hi}^\ell \right) \\
		-&\frac14
			\left(\alpha_{ik}^\ell-\alpha_{k\ell}^i+\alpha_{\ell i}^k\right)
			\left(\alpha_{j\ell}^h-\alpha_{\ell h}^j+\alpha_{hj}^\ell \right) \\
		-&\frac12
			\alpha_{ij}^\ell\left(\alpha_{\ell k}^h-\alpha_{kh}^\ell+\alpha_{h\ell}^k\right)
		\bigg] X_h .
	\end{align*}
	The sectional curvature satisfies
	\begin{align}\label{eq2130}
	 	{\rm Sec}(X_1, X_2) =
		\sum_{\ell=1}^n \bigg[
		&- \frac12 \alpha_{12}^\ell\left(\alpha_{12}^\ell-\alpha_{2\ell}^1-\alpha_{\ell1}^2\right) \nonumber\\
		&- \frac14 \left(\alpha_{12}^\ell-\alpha_{2\ell}^1+\alpha_{\ell1}^2\right)
			\left(\alpha_{12}^\ell+\alpha_{2\ell}^1-\alpha_{\ell1}^2\right) \nonumber\\
		&- \alpha_{\ell1}^1\alpha_{\ell2}^2
			\bigg].
	\end{align}
\end{proposition}
\begin{proof} Recall that we defined $R$ by \eqref{def_Riemannian_curvature_tensor}.
	So
	$$
	 	R(X_i, X_j, X_k) 
		= \nabla_{X_i}\nabla_{X_j}X_k - \nabla_{X_j}\nabla_{X_i}X_k - \nabla_{\sum_{\ell}\alpha_{ij}^\ell X_\ell} X_k 
	$$
	\begin{multline*}	
		= \nabla_{X_i}\left( \frac12\sum_\ell(\alpha_{jk}^\ell-\alpha_{k\ell}^j+\alpha_{\ell j}^k)X_\ell \right) +\hfill \\
			- \nabla_{X_j}\left(\frac12\sum_\ell(\alpha_{jk}^\ell-\alpha_{k\ell}^i+\alpha_{\ell i}^k)X_\ell\right) + \\
\hfill		- \sum_\ell\alpha_{ij}^\ell\frac12 \sum_h (\alpha_{\ell k}^h - \alpha_{kh}^\ell + \alpha_{h\ell}^k) X_h \\
		= \frac12 \sum_\ell (\alpha_{jk}^\ell-\alpha_{k\ell}^j+\alpha_{\ell j}^k) \frac12\sum_h(\alpha_{i\ell}^h-\alpha_{\ell h}^i+\alpha_{hi}^\ell) X_h +\hfill \\
			- \frac12\sum_\ell(\alpha_{ik}^\ell-\alpha_{k\ell}^i+\alpha_{\ell i}^k) \frac12\sum_h(\alpha_{j\ell}^h-\alpha_{\ell h}^j + \alpha_{hj}^\ell) X_h + \\
\hfill		- \sum_\ell\frac12 \alpha_{ij}^\ell \sum_h(\alpha_{\ell k}^h-\alpha_{kh}^\ell+\alpha_{h\ell}^k)X_h \\
		= \sum_{\ell,h} \bigg[ \frac14(\alpha_{jk}^\ell-\alpha_{k\ell}^j+\alpha_{\ell j}^k)(\alpha_{i\ell}^h-\alpha_{\ell h}^i+\alpha_{hi}^\ell) + \hfill \\
			- \frac14(\alpha_{ik}^\ell-\alpha_{k\ell}^i+\alpha_{\ell i}^k) (\alpha_{j\ell}^h-\alpha_{\ell h}^j + \alpha_{hj}^\ell) + \\
			-\frac12 \alpha_{ij}^\ell(\alpha_{\ell k}^h-\alpha_{kh}^\ell+\alpha_{h\ell}^k) \bigg] X_h.
	\end{multline*}
	Regarding the sectional curvature \eqref{def_sectional_curvature}, since $X_1, X_2$ are orthonormal, we have
	\[
	{\rm Sec}(X_1, X_2) = \langle R(X_1, X_2, X_2), X_1 \rangle .
	\]
	So, using the above formula with $i=1$, $j=k=2$, and $h=1$, we have
	$$
	 	{\rm Sec}(X_1, X_2) 
		= \langle R(X_1, X_2, X_2), X_1 \rangle \\
	$$
	\begin{multline*}		
		= \sum_\ell\bigg[
			\frac14(-\alpha_{2\ell}^2+\alpha_{\ell2}^2)(\alpha_{1\ell}^1-\alpha_{\ell1}^1) + \hfill\\
			- \frac14 (\alpha_{12}^\ell-\alpha_{2\ell}^1+\alpha_{\ell1}^2)(\alpha_{2\ell}^1-\alpha_{\ell1}^2+\alpha_{12}^\ell) + \\
\hfill		-\frac12 \alpha_{12}^\ell(\alpha_{\ell2}^1-\alpha_{21}^\ell+\alpha_{1\ell}^2) \bigg] \\
\end{multline*}
	\begin{multline*}	 	
		= \sum_\ell \bigg[
			\frac14 \left(2\alpha_{\ell2}^2(-2\alpha_{\ell1}^1)\right) + \hfill\\
			- \frac14 (\alpha_{12}^\ell-\alpha_{2\ell}^1+\alpha_{\ell1}^2)(\alpha_{12}^\ell+\alpha_{2\ell}^1-\alpha_{\ell1}^2) + \\
\hfill		- \frac12\alpha_{12}^\ell(\alpha_{12}^\ell-\alpha_{2\ell}^1-\alpha_{\ell1}^2) \bigg] \\
		= \sum_\ell \bigg[ -\alpha_{\ell1}^1\alpha_{\ell2}^2 
			- \frac14(\alpha_{12}^\ell-\alpha_{2\ell}^1+\alpha_{\ell1}^2)(\alpha_{12}^\ell+\alpha_{2\ell}^1-\alpha_{\ell1}^2) 
-\frac12 \alpha_{12}^\ell(\alpha_{12}^\ell-\alpha_{2\ell}^1-\alpha_{\ell1}^2) \bigg] .
	\end{multline*}
\end{proof}
\begin{lemma}\label{lemma_skew_sec}
 	Let $X\in\Lie(G)$. 
	If $\ad_X$ is skew-adjoint, then ${\rm Sec}(X,Y)\ge0$, for each $Y\in\Lie(G)$ that is linearly independent from $X$.
\end{lemma}
\begin{proof}
 	Recall that $\ad_X$ being skew adjoint means
	\[
	\langle \ad_XY,Z \rangle = - \langle Y, \ad_XZ \rangle.
	\]
	Let $X\in\Lie(G)$ be such that $\ad_X$ is skew adjoint.
	Since $\ad_{\lambda X}=\lambda\ad_X$, we can assume that $\langle X, X \rangle=1$. 
	Let $Y\in\Lie(G)$ be such that $\langle X,Y \rangle=0$ and $\langle Y,Y \rangle=1$.
	Take an orthonormal basis $X_1, \dots, X_n$ of $\Lie(G)$ with $X_1=X$ and $X_2=Y$, and let $\alpha_{ij}^k$ be the corresponding structural constants, i.e., $\ad_{X_i}(X_j) = [X_i, X_j]=\sum_k\alpha_{ij}^kX_k$. Then
	\[
	\alpha_{1j}^k 
	= \langle \ad_{X_1}X_j, X_k \rangle 
	= -\langle X_j, \ad_{X_1}X_k \rangle 
	= - \alpha_{1k}^j .
	\]
	Thus, we have $\alpha_{\ell 1}^2=-\alpha_{1\ell }^2=\alpha_{12}^\ell$ and $\alpha_{\ell 1}^1=0$. 
	Therefore, formula \eqref{eq2130} simplifies to
	\begin{eqnarray*}
		{\rm Sec}(X_1, X_2) 
		&=&\sum_\ell 
		- \frac12 \alpha_{12}^\ell(-\alpha_{2\ell}^1) 
		 -\frac14 (2\alpha_{12}^\ell -\alpha_{2\ell}^1)(\alpha_{2\ell}^1) 
			\\
		&=& \sum_\ell 
		 \frac12\alpha_{12}^\ell\alpha_{2\ell}^1
			- \frac12 \alpha_{12}^\ell\alpha_{2\ell}^1 
			+ \frac14 (\alpha_{2\ell}^1)^2 
				\\
		&=&\sum_\ell \frac 14 (\alpha_{2\ell}^1)^2 \ge0 . 
	\end{eqnarray*}
\end{proof}

\subsection{Bi-invariant metrics}\label{sec:biinv_compact}
\index{compact group}
\index{bi-invariant! -- Riemannian metric}
Recall that Riemannian metrics of Lie groups that are left invariant and right invariant are said to be bi-invariant.
We obtained characterizations for bi-invariance of metrics in Theorems~\ref{teo2136} and \ref{thm2111}.

\begin{corollary}\label{teo2136+thm2111}
\index{right-invariant! -- Riemannian metric}
\index{isometry}\index{adjoint! -- representation, $\Ad$}
\index{adjoint! -- map, $\ad$}
 	Let $G$ be a connected Lie group with Lie algebra $\g$ and with a left-invariant Riemannian metric $\langle \cdot, \cdot \rangle$ with Levi-Civita connection $\nabla$.
	The following are equivalent:
	\begin{description}
	\item[\ref{teo2136+thm2111}.i.] 	$\langle \cdot, \cdot \rangle$ is bi-invariant;
		\item[\ref{teo2136+thm2111}.ii.] 	$\Ad_g$ is an isometry, for all $g\in G$;
		\item[\ref{teo2136+thm2111}.iii.] 	$\ad_X$ is skew-adjoint, for all $X\in\g$;
	 \index{one-parameter subgroup}
\index{geodesic}
\index{Levi-Civita connection}
 
		\item[\ref{teo2136+thm2111}.iv.] 	$\exp_{1_G}=\exp$; 
		\item[\ref{teo2136+thm2111}.v.] 		 
		$
		\nabla_XY=\frac12[X,Y], $ for all 
		$X,Y\in\g$;
		\index{skew-adjoint}\index{antisymmetric}\index{adjoint! -- map, $\ad$}
	\end{description}
\end{corollary}

Lemma~\ref{lemma_skew_sec} gives a non-trivial property of bi-invariant metrics:
\begin{corollary}
 	If a connected Lie group is equipped with a bi-invariant metric, then all its sectional curvatures are nonnegative.
\end{corollary}

\begin{theorem}\label{compact bi-invariant}
 	Every compact Lie group $G$ admits a bi-invariant Riemannian metric.
\end{theorem}
\begin{proof}
 	Let $\langle \cdot, \cdot \rangle$ be a scalar product on $T_{1_G}G$.
	Let $\vol$ be a left-invariant probability measure on $G$. Define a new product $\langle\!\langle \cdot, \cdot \rangle\!\rangle$ averaging $(\Ad_g)_*\langle \cdot, \cdot \rangle$ with $\vol$, i.e.,
	\[
	\langle\!\langle X,Y \rangle\!\rangle= \int_G \langle \Ad_g(X), \Ad_g(Y) \rangle \dd\vol(g), 
	\qquad\forall X,Y\in T_{1_G}G .
	\]
	Notice that $\langle\!\langle \cdot, \cdot \rangle\!\rangle$ is finite since $\vol(G)<\infty$ and since $g\mapsto  \langle \Ad_g(X), \Ad_g(Y) \rangle$ is a continuous function on the compact group $G$ and thus it is bounded.
	Moreover, for all $g\in G$ the product $\langle\!\langle \cdot, \cdot \rangle\!\rangle$ is $\Ad_g$ invariant:
	\begin{eqnarray*}
	 	\langle\!\langle \Ad_gX, \Ad_gY \rangle\!\rangle 
		&=& \int_G \langle \Ad_h\Ad_gX, \Ad_h\Ad_gY \rangle \dd\vol(h) \\
		&=& \int_G \langle \Ad_{g'}X, \Ad_{g'}Y \rangle \dd\vol(g')\\
		&=& \langle\!\langle X,Y \rangle\!\rangle .
	\end{eqnarray*}
	Extending $\langle\!\langle \cdot, \cdot \rangle\!\rangle$ by left translations, we obtain a left-invariant Riemannian metric for which $\Ad_g$ are isometries.
	Theorem~\ref{teo2136} implies that this metric is right invariant.
\end{proof}
Because for each bi-invariant Riemannian metric, the volume measure 
is bi-invariant, we obtain a first consequence:
\begin{corollary}\index{probability measure}
 	Every compact Lie group can be equipped with a probability measure that is bi-invariant.
\end{corollary}

Because for each bi-invariant Riemannian metric, the exponential map coincides with the Riemannian exponential (see \eqref{teo2136+thm2111}.iv), and the Riemannian exponential on complete connected manifolds is surjective (see Hopf-Rinow Theorem~\ref{HRCV}), we obtain a second consequence: 
\begin{corollary} 
 	On every compact connected Lie group $G$, the exponential map $\exp:\Lie(G)\to G$ is surjective.
\end{corollary}

Here is another characterization of groups that admit bi-invariant metrics that are admissible, in the sense defined at page~\pageref{def_admissible_distance}.  
\begin{theorem}\label{Summary Valto} Let $G$ be a connected Lie group. Then the following are equivalent:
\begin{description}
\item[\ref{Summary Valto}.i.] There exists an admissible bi-invariant distance on $G$;
\item[\ref{Summary Valto}.ii.] The set $\Ad_G$ is a compact subset of 
$\GL(\g)$;
\item[\ref{Summary Valto}.iii.] There exists a bi-invariant Riemannian metric on $G$;
\item[\ref{Summary Valto}.iv.] The Lie group $G$ is the direct product of a compact group and a vector group, that is, $\R^n$ for some $n\in\N$.
\end{description}
\end{theorem}

\begin{proof}[Sketch of the proof]
If there exists an admissible bi-invariant distance on $G$, then by Lemma~\ref{rmk:boundedy_cpt}, we may additionally assume that this distance is boundedly compact. Thus, by Ascoli--Arzel\'a theorem, the space of isometries fixing $1_G$ is compact; see Proposition~\ref{Iso_gp_loc_cpt}. Since the distance is bi-invariant, the conjugation maps are isometries.
Therefore, the set $\Ad_G\subseteq \GL(\g) $ of their differentials is compact. Alternatively, one can directly construct a bi-invariant Riemannian metric on $G$ using Lemma~\ref{lemma:riem_metric}.

If the set $\Ad_G$ is compact in the space of linear transformations of the Lie algebra $\g$, then, 
one constructs a bi-invariant Riemannian metric on $G$, as in the proof of Theorem~\ref{compact bi-invariant}.

Clearly, every bi-invariant Riemannian metric on $G$ is an admissible bi-invariant distance. 
Thus, the first three items are equivalent.

Regarding the last item, if the Lie group $G$ is the product of a compact group and a commutative group, then, obviously, the set $\Ad_G$ is compact. 
The converse implication is more difficult: 
It is a result of Milnor that	every connected Lie group that admits a bi-invariant metric is the Cartesian product of a compact group and a commutative group; see 
 \cite[Lemma~7.5]{Milnor}.
\end{proof}


From the topological viewpoint, every Lie group is homeomorphic to the product of a compact group and a vector space. This result is due to Iwasawa; see Exercise~\ref{ex Iwasawa}.

\subsection{More results on curvature}\index{curvature}

In Milnor's article \cite{Milnor}, one can also find the following results on the curvature of left-invariant Riemannian metrics on Lie groups.
\index{Milnor}

\begin{theorem}[Milnor, {\cite[Theorem~3.3]{Milnor}}]
 	If $\Lie(G)$ is not commutative, then $G$ admits a left-invariant metric of strictly negative \emph{scalar} curvature.
\end{theorem}

\begin{theorem}[Milnor, {\cite[Theorem~2.2]{Milnor}}]
 	A connected Lie group admits a left-invariant metric with all \emph{Ricci} curvatures strictly positive if and only if it is compact with finite fundamental group.
\end{theorem}
Recall: ${\rm Ric}(X) := \sum_i {\rm Sec}(X, X_i)$, where $X_1, \dots, X_n$ is an orthonormal frame.
\begin{theorem}[Milnor, {\cite[Theorem~2.5]{Milnor}}]
 	If there are non-zero $X,Y,Z\in\Lie(G)$ such that $[X,Y]=Z$, then there is a left-invariant metric on $G$ such that ${\rm Ric}(X)<0<{\rm Ric}(Z)$.
\end{theorem}
\begin{theorem}[Milnor, {\cite[Corollary~7.7]{Milnor}}]
 	Every Lie group whose universal covering space is compact admits a bi-invariant metric of constant Ricci curvature $+1$.
\end{theorem}


\section{Isometries of metric groups as Riemannian isometries}\label{sec Isometries of metric groups}\index{isometry}
Recall from Section~\ref{sec: metric_group} that a metric Lie group is a Lie group equipped with a left-invariant distance function that induces the manifold topology.\index{metric! -- Lie group}\index{Lie group! metric --}
In this section, we prove that isometries between metric Lie groups are Riemannian isometries for some left-invariant Riemannian metric.
If $M$ is a Lie group and $\rho$ is a left-invariant Riemannian metric tensor on $M$, then one has an induced Riemannian distance $d_\rho$ and, by the theorem of Myers and Steenrod 
\cite{Myers-Steenrod}, the group $\Isom(M,d_\rho) $ of distance-preserving bijections coincides with the group $\Isom(M,\rho)$ of tensor-preserving diffeomorphisms.
In the following, we shall write $(M,\rho)$ to denote the metric Lie group $(M,d_\rho) $.

\begin{theorem} \label{prop:isometries_are_riemannian}
If \( (M_1,d_1) \) and \( (M_2,d_2) \) are connected metric Lie groups, then there exists left-invariant Riemannian metrics \( \rho_1 \) and \( \rho_2 \) on \( M_1 \) and \( M_2 \), respectively, such that 
\( \Isom(M_i,d_i) \subseteq \Isom(M_i,\rho_i) \) for \( i \in \{1,2\} \) and for every isometry
\( F \colon (M_1,d_1) \ra (M_2,d_2) \) the map 
\( F \colon (M_1,\rho_1) \ra (M_2,\rho_2) \) is a Riemannian isometry.
\end{theorem}
Before proving Theorem~\ref{prop:isometries_are_riemannian}, we provide an auxiliary result.
We first consider the case \( (M_1,d_1) = (M_2,d_2) \).
\begin{lemma}\label{lemma:riem_metric}
If \( (M,d) \) is a connected metric Lie group, then there is a left-invariant Riemannian metric \( \rho \) such that 
\(\Isom(M,d) \subseteq \Isom(M,\rho) \).
\end{lemma}

\begin{proof}
Let \( S \) be the stabilizer subgroup of \( \Isom(M,d) \) at the identity element \( 1=1_M \) of $M$.
By Proposition~\ref{Iso_gp_loc_cpt}, together with Lemma~\ref{rmk:boundedy_cpt}, the topological group $S$ is compact. 
 Let \( \mu_S \) be the probability Haar measure on \( S \).
 Moreover, from the discussions in Section~\ref{sec:smooth}, we know that $S$ acts smoothly on \( M \).

Fix a scalar product \( \bbraket{\cdot, \cdot} \) on the tangent space \( T_1 M \) at 1. 
 Consider for \( v,w \in T_1 M \) the value
\[ \langle{v,w}\rangle := \int_S \bbraket{\dd F v, \dd F w} \dd \mu_S(F) .\]
Then \( \langle{\cdot, \cdot}\rangle \) defines an \( S \)-invariant scalar product on \( T_1 M \), and one can take \( \rho \) as the left-invariant Riemannian metric that coincides with \( \langle{\cdot, \cdot}\rangle \) at the identity.
\end{proof} 
\begin{proof}[Proof of Theorem~\ref{prop:isometries_are_riemannian}] 
By Lemma~\ref{lemma:riem_metric} let \( \rho_2 \) be a Riemannian metric on \( M_2 \) with \begin{equation}\label{eq:incl2}
 \Isom(M_2,d_2) \subseteq \Isom(M_2,\rho_2) .\end{equation} 
 Fix an isometry \( \tilde F \colon (M_1,d_1) \ra (M_2,d_2) \); if there is none, then there is nothing else to prove. By Theorem~\ref{smooth:thm}, the map \( \tilde F \) is smooth, and we may define a Riemannian metric on \( M_1 \) by \( \rho_1:= \tilde F^* \rho_2 \). There are two things to verify: a) \( \Isom(M_1,d_1) \subseteq \Isom(M_1,\rho_1) \), which in particular implies that \( \rho_1 \) is left-invariant and b) every isometry \( F \colon (M_1,d_1) \ra (M_2,d_2) \) is an isometry of Riemannian manifolds.

Since by construction \( \tilde F \) is also a Riemannian isometry, the map \( I \mapsto \tilde F \circ I \circ \tilde F^{-1} \) is a bijection between \( \Isom(M_1,d_1) \) and \( \Isom(M_2,d_2) \) and between \( \Isom(M_1,\rho_1) \) and \( \Isom(M_2,\rho_2) \).
Therefore the inclusion
\eqref{eq:incl2} implies the inclusion
\( \Isom(M_1,d_1) \subseteq \Isom(M_1,\rho_1) \).

From \( F \circ \tilde F^{-1} \in \Isom(M_2,d_2) \subseteq \Isom(M_2,\rho_2) \) we obtain \( (F \circ \tilde F^{-1})^* \rho_2 = \rho_2 \). Consequently, we conclude \( F^* \rho_2 = \tilde F^* (F \circ \tilde F^{-1})^* \rho_2 = \rho_1 \).
\end{proof} 




\section{Exercises}


\begin{exercise}
Let $G$ be a Lie group equipped with a Riemannian metric tensor that induces a distance function $d$.
Then, the Riemannian metric tensor is left invariant
 if and only if left translations are isometries with respect to $d$.
\end{exercise}
\begin{exercise}
For each left-invariant distance function $d$ on a group $G$, the following are equivalent:
(i) $d$ is right invariant;
(ii) $d$ is inversion invariant, that is, $d(x, y) = d(x^{-1}, y^{-1})$, for all $x,y$ in the group;
 (iii) conjugations are isometries.
 
 \end{exercise}

\begin{exercise}\label{Seba 4 oct 2024}
Consider the Lie group $G: =\R\rtimes_\theta \Z$ the semi-direct product of groups given by 
$\theta (m) (t) : = 2^m t $ for $m\in \Z$ and $t\in \R$, as in Definition~\ref{def Semi-direct product of groups}.
Equip $G$ with any Riemannian left-invariant metric.
Then, every map $\ad_X$ is skew-adjoint, for all $X\in\Lie(G)$, but for some $g\in G$ the map
$\Ad_g$ is not an isometry.  Compare this example with Theorem~\ref{teo2136}.
\end{exercise}

\begin{exercise}
	A linear connection $\nabla$ on $G$ is left invariant if and only if for all left-invariant vector fields $X, Y$ the vector field $\nabla_XY$ is left invariant.
\end{exercise}
\begin{exercise}
	The Levi-Civita connection of a left-invariant Riemannian metric on $G$ is left-invariant.
\end{exercise}
\begin{exercise}\index{Christoffel symbols}
	A linear connection $\nabla$ on $G$ is left invariant if and only if the Christoffel symbols $\Gamma_{ij}^k$ defined by $\nabla_{X_i}X_j = \sum_k \Gamma_{ij}^k X_k$, with respect to some/every frame of left-invariant vector fields $X_1, \ldots, X_n$, are constant functions on the group.
\end{exercise}

\begin{exercise}\label{ex:curvature formula}
 Let $G$ be a Riemannian Lie group,  $\nabla$  its Levi-Civita connection, and $R$ its Riemannian curvature tensor. 
For each $ X, Y\in \g$, define
$U(X,Y) := \frac12( \nabla_XY + \nabla_Y X)$, which is the  symmetric part of $\nabla$.
Then, for all $ X, Y, Z\in \g$ one has
\begin{description}
\item[\ref{ex:curvature formula}.i.] $\nabla_XY =U(X,Y)+\frac12[X,Y], $

\item[\ref{ex:curvature formula}.ii.] $\langle U ( X , Y ) , Z \rangle = \frac12 \langle X , [ Z , Y ] \rangle + \frac12 \langle Y , [ Z , X ] \rangle $

\item[\ref{ex:curvature formula}.iii.] $\langle R(X, Y )Y, X\rangle = \norm{U(X,Y)}^2 -\langle U(X, X),U(Y,Y)\rangle
-\frac34 \norm{[X,Y]}^2 - \frac12\langle[X,[X,Y]],Y\rangle - \frac12\langle[Y,[Y, X]], X\rangle.$
\end{description}
{Hint.} Use Koszul's formula and the fact that the Levi-Civita connection has zero torsion.
\end{exercise}

%


\begin{exercise}
	Every Lie group admits a left-invariant Riemannian metric and a left-invariant measure, as follows.
	Let $X_1, \dots, X_n$ be left-invariant vector fields forming a basis of $\Lie(G)$.
	Consider the Riemannian metric that makes $X_1, \dots, X_n$ orthonormal and the differential $n$-form $\vol$ for which $\vol(X_1, \dots, X_n)=1$.
\end{exercise}
\begin{exercise}
	If $G$ is a compact Lie group, then it admits a left-invariant Riemannian volume form that gives a probability measure.
\end{exercise}


\begin{exercise}\label{ex:Gdelta2}\index{$G_z$} 
Calculate the sectional curvatures of the Lie groups $G_z$ from Example~\ref{model_typeR} when equipped with some left-invariant Riemannian metric.
\end{exercise}

\begin{exercise}\label{ex sec curvature complex hyperbolic plane}
Let $\g$ be the 4D Lie algebra with basis $X_1, \ldots, X_4$ and only non-zero structural constants, as in \eqref{eq1655},
$$ c_{12}^3 =1, \quad c_{41}^1 = 1, \quad c_{42}^2 = 1, \quad c_{43}^3 = 2. $$
Consider a Lie group $G$ with $\g$ as Lie algebra, and consider the left-invariant Riemannian structure on $G$ for which $X_1, \ldots, X_4$ are orthonormal.
Verify the following values for the sectional curvature:
\begin{eqnarray*} 
{\rm Sec}(X_1, X_2)= -7/4 ,		&\qquad&	{\rm Sec}(X_1, X_4)= -1 ,
\\ 	{\rm Sec}(X_1, X_3)= -7/4 , 	&& 	{\rm Sec}(X_2, X_4)= - 1,
\\ 	{\rm Sec}(X_2, X_3)= -7/4 , 	&&	{\rm Sec}(X_3, X_4)= - 4.
\end{eqnarray*}
{\it Hint.} Apply \eqref{eq2130}. 
\end{exercise}

\begin{exercise}\index{Heisenberg! -- group} 
Let $\mathfrak h$ be the Heisenberg Lie algebra with basis $X_1, X_2, X_3:=[X_1, X_2]$.
Consider the left-invariant Riemannian structure on the Heisenberg group for which $X_1, X_2, X_3$ are orthonormal.
Verify the following values for the sectional curvature:
\\ $	{\rm Sec}(X_1, X_2)= -3/4 ,$
\\ $	{\rm Sec}(X_1, X_3)= 1/4 , $
\\ $	{\rm Sec}(X_2, X_3)= 1/4 . $
\\{\it Hint.} In this basis, the only nontrivial structural constant is $ c_{12}^3 =1$. Apply \eqref{eq2130}.
\end{exercise}

\begin{exercise}[Iwasawa's Theorem; see {\cite[page 327]{Milnor}}]\label{ex Iwasawa}
 	Let $G$ be connected Lie group. Then:
\\\ref{ex Iwasawa}.i.
		Every compact subgroup is contained in a maximal compact subgroup $H$.
	\\\ref{ex Iwasawa}.ii.	Every maximal compact subgroup 
		is a connected Lie subgroup.
\\\ref{ex Iwasawa}.iii.
		Every two maximal compact subgroups are conjugate.
\\\ref{ex Iwasawa}.iv.
		As a topological space, $G$ is homeomorphic to the product of $H$ and some Euclidean space $\R^m$.
\end{exercise}



\chapter{Nilpotent Lie groups}\label{ch_Nilpotent}

Nilpotent Lie groups play a fundamental role in the study of geometric structures and differential equations. These groups are characterized by an algebraic property known as nilpotency. The nilpotency assumption imposes a certain level of commutativity for the group's structure, leading to a wealth of remarkable properties. By restricting our attention to nilpotent Lie groups, we gain a deeper understanding of their topology and geometry.
The nilpotency condition, which also reflects a condition on the Lie algebra, serves as a guiding principle that allows us to explore the interplay between algebra and geometry. It paves the way for applications in geometric analysis, geometric group theory, harmonic analysis, and control theory, as well as number theory, dynamics, representation theory, and, ultimately, physics. In this chapter, we delve into the world of nilpotent Lie groups.

A recent reference that deserves strong recommendation is \cite{Hilgert_Neeb:book}. Other valuable reading materials on this topic include \cite{Raghunathan, Jacobson, Warner, Corwin-Greenleaf, Knapp}. 
 While our exposition may not be as comprehensive as those references, we will focus on the necessary concepts to understand Carnot groups, as well as other sub-Finsler Lie groups such as boundaries of Heintze groups, Malcev closures of finitely generated nilpotent groups, and their asymptotic cones. 
 
 Throughout our discussion, we will maintain a perspective rooted in differential geometry and linear algebra. It is worth noting that one of the compelling aspects of nilpotent Lie groups is their appearance as tangent spaces of sub-Riemannian manifolds, similar to how Euclidean vector spaces serve as tangents to Riemannian manifolds. We will discover that, akin to vector spaces, these tangents possess nilpotency, simply connectedness, and dilation structures.

\section{Nilpotent Lie algebras}

We begin this section by introducing the concept of a nilpotent Lie algebra.
Nilpotent Lie algebras are those for which iterated brackets $[x_1,[x_2,[x_3,[\ldots ]]]]$ of sufficiently large order vanish.
We anticipate that for connected Lie groups, a Lie algebra is nilpotent if and only if the group is nilpotent as a group, according to Definition~\ref{def_nilpotent_group}.
This correspondence between nilpotent Lie algebras and nilpotent Lie groups is a result that we will further explore in Section~\ref{sec:Lie_groups_nilpotent_algebras}. 
Typical examples of nilpotent Lie algebras are Lie algebras of strictly upper-triangular matrices, where the diagonal elements are all zero. The first fundamental result that we present about nilpotent Lie algebras is Engel's Theorem, which translates nilpotency into a pointwise condition; see Section~\ref{sec Engel}.
We close the section with the statement of Birkhoff-Embedding Theorem, which we will only prove for positively graded Lie algebras in Section~\ref{sec Birkhoff for stratified}. Birkhoff's theorem will tell us that every nilpotent Lie algebra is isomorphic to some subalgebra of some space of strictly upper-triangular matrices; see Corollary~\ref{cor_nilpotent_triangular}.

\begin{definition}[Nilpotent Lie algebra]\label{def lower central series}
Let $\mathfrak {g}$ be a Lie algebra. The elements of the {\em descending central series} of $\g$, also called {\em lower central series} of $\g$, are inductively defined by\index{descending central series! -- for a Lie algebra}\index{lower central series! -- for a Lie algebra}
\begin{eqnarray*}
C^1(\g) :=\g^{(1)}:= \g, &\quad& C^2(\g):=\g^{(2)} :=[\g, \g]= [\g, C^1(\g)], \\
& & C^{n}(\g) :=\g^{(n)}:= [\g, C^{n-1}(\g)], \qquad \forall n\in\N.
\end{eqnarray*}
Here, for $V, W \subseteq \g$, we use the notation $ [V,W]:=\Span\{ [v,w]\,:\, v\in V, w\in W \}.$
The space $C^2(\g)$ is called {\em commutator subalgebra}.\index{commutator! -- subalgebra} 
The Lie algebra $\g$ is said to be {\em nilpotent} if there is $d\in \N$ such that $C^{d+1}(\g)=\{0 \}.$ If $d$ is minimal with this property, then it is called {\em nilpotency degree} (or {\em nilpotency step} or, simply, {\em step}) of $\g$, and $\g$ is said $d$-{\em step nilpotent}.\index{nilpotent! -- Lie algebra}\index{step|see {nilpotency step}}\index{nilpotency step}
\index{Lie algebra! nilpotent --}
\end{definition}

One can rephrase the definition by saying that a Lie algebra $\g$ is $s$-step nilpotent if and only if all brackets of at least $s +1$ elements of $\g$ are $0$ but not every bracket of order $s$ is.

\begin{remark}
Each $C^n(\g)$ is an ideal and actually
\begin{equation*}
[C^n(\g), C^n(\g)] \subseteq [C^n(\g), \g] \stackrel{{\rm def}}{=} C^{n+1}(\g) \subseteq C^{n}(\g), 
\end{equation*}
where the last inclusion holds by induction, noting that $C^{2}(\g) \subseteq C^{1}(\g).$ 
\end{remark}

 A nilpotent Lie algebra $\g$ has always non-trivial center $Z(\g)$ by 
 Proposition~\ref{prop_basic_properties_nilpotent_Lie_albegra}.iii. 
 In fact, if $\g$ is $s$-step nilpotent, the subalgebra $\g^{(s)}$ is central.
 Be aware that the center might be strictly larger than $\g^{(s)}$; see Exercise~\ref{ex:lagercenter}.

\begin{proposition}\label{prop_basic_properties_nilpotent_Lie_albegra}
Let $\g$ be a Lie algebra. 
\begin{description}
\item[\ref{prop_basic_properties_nilpotent_Lie_albegra}.i.] If $\g$ is nilpotent, then all its subalgebras and all homomorphic images of $\g$ are nilpotent.
\item[\ref{prop_basic_properties_nilpotent_Lie_albegra}.ii.] If $\mathfrak {a} < Z(\g) $ and $\g / \mathfrak {a}$ is nilpotent, then $\g$ is nilpotent. 
\item[\ref{prop_basic_properties_nilpotent_Lie_albegra}.iii.] If $\g \ne \{0\}$ and $\g$ is nilpotent of step $s$, then $ \{0\}\ne C^s(\g) \subseteq Z(\g) .$
\item[\ref{prop_basic_properties_nilpotent_Lie_albegra}.iv.] If $\g$ is nilpotent of step $s,$ then $\ad_x^s\equiv 0$, for every $x\in \g $.
\item[\ref{prop_basic_properties_nilpotent_Lie_albegra}.v.] If $\mathfrak {i}$ is an ideal of $ \g,$ then $C^n(\mathfrak {i})$ is an ideal of $\g,$ for every $n\in \N.$ 
\end{description}
\end{proposition}

\begin{proof}
$(i).$ Let $\mathfrak {h}$ be a subalgebra of $\g$. Then $[\mathfrak {h}, \mathfrak {h}]\subset [\g, \g]$ and so, by induction, $C^n(\mathfrak {h}) \subset C^n(\mathfrak {g}).$ Hence, the Lie algebra $\mathfrak {h}$ is nilpotent, if so is $\g.$ Moreover, if we consider a Lie algebra homomorphism $\alpha : \mathfrak {g} \to \mathfrak {h}$, then $[\alpha (\g), \alpha (\g)]=\alpha ([\g, \g ])$. By induction we have that 
\begin{equation}\label{equation1303}
C^n(\alpha (\g)) = \alpha (C^n(\g)), \qquad \forall n\in \N.
\end{equation}
 Consequently, the image $\alpha (\g)$ is nilpotent, if so is $\g.$

$(ii).$ If $\g / \mathfrak {a}$ is nilpotent, by definition we know that there is $n\in \N$ such that $C^n(\g / \mathfrak {a})=\{0\}$ in $\g / \mathfrak {a}$, i.e., $C^n(\g / \mathfrak {a})= \mathfrak {a} / \mathfrak {a}$. Now we apply \eqref{equation1303} with $\alpha $ equal to the projection $\g \to \g/\mathfrak {a}$ and we deduce that $C^n(\g) + \mathfrak {a} = C^n (\mathfrak {g}/\mathfrak {a})= \mathfrak {a},$ i.e., $C^n(\g) \subset \mathfrak {a} \subset Z(\g)$, where the last inclusion holds by assumption. This implies that $\g$ is nilpotent:
\begin{equation*}
C^{n+1}(\g)\stackrel{\rm def}{= }[\g, C^n(\g)]\subset[\g, Z(\g)]= \{0\}.
\end{equation*}

$(iii).$ By hypothesis, the natural number $s$ is such that $C^{s+1}(\g)= \{0\}$ and $C^{s}(\g) \ne \{0\}.$ Then, $ \{0\} \ne C^{s}(\g) \subseteq Z(\g)$, since $[C^{s}(\g), \g]\stackrel{\rm def}{= } C^{s+1}(\g)= \{0\}.$

$(iv).$ Since $\g$ is $s$-step nilpotent we have that $C^{s+1}(\g) = \{0\}$. So for every $x\in \g$
\begin{equation*}
(\mbox{ad}(x))^s(\g)= \underbrace{ [x,[x, \dots [x }_{\mbox{$s$ times}}, \g]\dots ] \subseteq C^{s+1}(\g) = \{0\}.
\end{equation*}

$(v).$ It follows from the general easy fact that if 
$\mathfrak {i}, \mathfrak {j}$ are ideals of $ \g $, then so is $[\mathfrak {i}, \mathfrak {j}] $.
\end{proof}

Proposition \ref{prop_basic_properties_nilpotent_Lie_albegra}.iv states that 
if $\g$ is nilpotent, then each
 map 
 $y \in \g \mapsto [X, Y]\in \g$
 is a nilpotent transformation in the sense of Definition~\ref{def_nilpotent_unipotent}.i. 
We shall see in Theorem \ref{Engel's theorem} that the inverse implication holds true.

\subsection{Examples of nilpotent Lie algebras}

We present some examples of nilpotent Lie algebras.
\begin{example}[Abelian Lie algebras]
{\em Commutative Lie algebras} are those for which $[X, Y]=[Y, X]$, for all elements $X$ and $Y$.
By anticommutativity, this is equivalent to
 $[\cdot, \cdot ]\equiv 0$.
 We also refer to commutative Lie algebras as {\em abelian Lie algebras}. 
Consequently, a Lie algebra $\g$ is commutative if and only if it is nilpotent with nilpotency step equal to 1, because $C^2( \g) \stackrel{{\rm def}}{=}[\g, \g ]$.\index{abelian Lie algebra}
\index{commutative Lie algebra}
 \end{example}
\begin{example}[Heisenberg Lie Algebra]\index{Heisenberg! -- algebra} 
The Heisenberg Lie algebra
\begin{equation*}
\mathfrak {nil}_3:= \left\{ 
\begin{bmatrix}
0& x & z\\
0& 0 & y \\
0& 0& 0\\
\end{bmatrix}
\,:\, x,y,z \in \R \right\} \subseteq \gl(3)
\end{equation*}
is nilpotent of step 2 because for such a Lie algebra $\g:=\mathfrak {nil}_3$
$$C^2(\g) =\left\{ 
\begin{bmatrix}
0& 0 & z\\
0& 0 & 0 \\
0& 0& 0\\
\end{bmatrix}
\,:\, z \in \R \right\} \qquad \text{ and } \qquad C^3(\g)=\{{\bf 0}\}.$$
 \end{example}
\begin{example}[Nilpotent Lie algebras of step 2]\label{free_nilp_step2}\index{step-two nilpotent Lie group}
We revise Definition~\ref{def_step_2_G} of nilpotent Lie algebras of step 2.
 Let $V, W$ be vector spaces and $q: V\times V \to W$ a skew-symmetric bilinear map, then $[(v_1,w_1),(v_2,w_2)]:=(0, q(v_1,v_2))$ is a Lie bracket on $V\times W,$ and $V\times W$ becomes a step-2 Lie algebra unless $q$ is identically equal to 0 in which case we get a commutative Lie algebra. Namely, we have $[[X, Y],Z]=0$ for every $X,Y,Z \in V\times W.$
\item As a more specific example, for $n\in \N$, we consider
\begin{equation}\label{model_free}
\begin{aligned}
V := \Lambda ^1 (\R^n) & =\{ 1\mbox{-forms on } \R^n\}, \\
W := \Lambda ^2 (\R^n)& =\{ 2\mbox{-forms on } \R^n\}, \\
q(v_1,v_2) &:=v_1 \wedge v_2,
\end{aligned}
\end{equation} where $\wedge$ is the wedge of 1-forms.
Then $\Lambda ^1 (\R^n) \times \Lambda ^2 (\R^n)$ becomes a step-2 Lie algebra called the {\em free-nilpotent Lie algebra} of rank $n$ and step 2 (we will generalize this example in Example~\ref{ex_free_nilpotent}).\index{rank! -- of a free Lie algebra}\index{free! -- nilpotent}
 \end{example}

One common convention in describing Lie algebras - and one that we shall often use - is the following. Suppose that $\g = \R$-$\Span \{X_1, \ldots, X_n\}$.
To describe the Lie algebra structure of $\g$, it suffices to give $[X_i, X_j]$ for all $i<j$ in terms of $X_1, \ldots, X_n$. 
We can shorten this description considerably by giving only the non-zero brackets; all others are assumed to be zero.

\begin{example}[Filiform algebras of the first kind]\label{ex_fili1}\index{filiform! -- algebra of the first kind}
The {\em $(n+ 1)$-dimensional filiform algebra of the first kind}
is the Lie algebra with basis $X, Y_1, Y_2, \ldots, Y_n$ and only non-trivial relations 
$$[X,Y_j]=Y_{j+1}, \quad\text{ for } j\in \{1, \ldots, n-1\}.$$
It is an $n$-step nilpotent Lie algebra and can be realized as a matrix algebra considering the matrices of the form:
$$\begin{bmatrix}
 0 	& x 	&0	&\cdot	& 0	&y_n\\
 	& \cdot &\ddots	&	.& .	&\vdots\\
 	& 	&\cdot &\ddots	& 0	&\vdots\\
 	&	&	&\cdot	&x	&y_2\\
 	&	&	&&\cdot		&y_1\\
0 	&	 & 	& 	&	&0
 \end{bmatrix}, \qquad \text{ for } x, y_1, \ldots, y_n\in\R.$$

 \end{example}

\begin{example}[Free-nilpotent algebras]\label{ex_free_nilpotent}\index{rank! -- of a free Lie algebra}\index{free! -- nilpotent}
The {\em free nilpotent Lie algebra of step $k$ and rank $n$} is defined to be the quotient algebra $\mathfrak{f}_n/\mathfrak{f}^{(k+1)}$, where $\mathfrak{f}_n$ is the
free Lie algebra on $n$ generators. It is not hard to see that it is finite-dimensional, as shown by the following reasoning. 
One abstractly considers $n$-many letters $X_1, \ldots, X_n$; then performs abstract Lie brackets of at most $k$-many of them, in every order,
e.g., $[X_1,[X_7,[X_5, X_3]]]$ but also $[[X_5, X_1], [X_7, X_3]] $;
 then one considers the span, imposing anti-commutativity and Jacobi identity. 
 \index{database! -- for bases of free Lie algebras}
 There is a web app for calculating bases of free Lie algebras available here: \url{https://coropa.sourceforge.io/#cgi}

For example, the Lie algebra of rank $2$ and step $3$ is given by the diagram
\begin{center}
 \begin{tikzcd}[end anchor=north]
 & X_1 \ar[dr, no head] \ar[dd, no head]& &X_2\ar[dl, no head]\ar[dd, no head] & 
 \\
 & &X_{3} \ar[dl, no head]\ar[dr, no head] &\quad\;&
 \\
 &X_{4} & &X_{5} &
 \end{tikzcd},
\end{center}
which has to be read as $[X_1, X_2]=X_3$, $[X_1, X_3]=X_4$, and $[X_3, X_2]=X_5$.
Each bracket relation is expressed by a V-shape in the diagram, and it should be read as \textit{left} arm of the V bracket the \textit{right} arm gives the vector at the bottom of the V.\label{first description of diagram}
 If the diagram should be read differently, we shall use another notation; see page \pageref{description of diagram}.
 \index{diagram! description of the --}
 \end{example}

\begin{example}[Strictly upper-triangular matrix algebras]\label{example Strictly upper-triangular matrix algebra}
The algebra of strictly upper-triangular $n\times n$ matrices is an $(n-1)$-step nilpotent Lie algebra of dimension $n(n- 1)/2$, and its center is one-dimensional. 
	We denote this important Lie algebra by\index{$ {\mathfrak{nil}}_n$}
	\[
	 {\mathfrak{nil}}_n := \left\{
	\begin{bmatrix}
	 	0 & * & * \\
		 0 & \ddots & * \\
		 0 & 0 & 0
	\end{bmatrix}
	\right\}
	\subset\mathfrak{gl}(n).
	\]
	See the later Example~\ref{example Strictly upper-triangular matrix algebra and group} for the associated simply connected Lie group ${\rm {Nil }}_n$.
	By a result of Birkhoff, we shall see that, up to isomorphism, the subalgebras of $ {\mathfrak{nil}}_n$ are the most general examples of nilpotent Lie algebras; see Corollary~\ref{cor_nilpotent_triangular}. 
	Moreover, by a result of Engel, if $\g$ is a Lie subalgebra of $\mathfrak{gl}(n) $ made of nilpotent transformations (in the sense of Definition~\ref{def_nilpotent_unipotent}), then, up to a change of basis of $\R^n$,  we have that $\g$ is a subalgebra of $ {\mathfrak{nil}}_n$; see Theorem~\ref{Engel_linear_Lie}.
	 \end{example}
	 
\begin{example}[Generalized Flag-Shifting Lie Algebras, $ \g_{\rm nil} (\mathcal F)$]\label{ex g nil F} Here is a slight coordinate-free generalization of the previous example.
Let 	$V $ be a finite-dimensional vector space and let $\mathcal F= (V_0, \dots, V_m)$ be a flag for $V$, i.e.,
$$V_0=\{0\} \subsetneq V_1\subsetneq \ldots \subsetneq V_m=V.$$
The set of flag-preserving transformations, denoted by $\g(\mathcal F)$, is not nilpotent; see Exercise~\ref{ex g F}. Whereas the set of {\em flag-shifting transformations}\index{flag-shifting transformations}
 $$ \g_{\rm nil} (\mathcal F) :=\big\{A\in \gl(V) : A(V_k)\subseteq V_{k-1}, \forall k\in\{1, \ldots, m\}\big\}$$ 
 is nilpotent of step at most $m-1$; see Exercise~\ref{ex 11 1 2024}. Actually, the Lie algebra $ \g_{\rm nil} (\mathcal F) $ can be seen as a subalgebra of $	 {\mathfrak{nil}}_n$, for $n:=\dim V$.
 \end{example}

\begin{example}\label{N51}The following Lie algebra is denoted as $\mathfrak{n}_{5,1}$; see \cite[page 162]{LeDonne_Tripaldi}.
The non-trivial brackets in $\mathfrak{n}_{5,1}$ with respect to some basis $X_1, \ldots, X_5$ are the following:\index{$\mathfrak{n}_{5,1}$}
\begin{equation*}
 [X_1, X_2]=X_3\,, \quad [X_1, X_3]=X_4\,, \quad[X_1, X_4]=[X_2, X_3]=X_5\,.
\end{equation*}
This is a nilpotent Lie algebra of rank 2 and step 4. 
 The Lie brackets can be pictured with the diagram:
\begin{center}
 
	\begin{tikzcd}[end anchor=north]
		X_1\ar[ddr,no head] \ar[dddr,no head, end anchor={[xshift=-2.5ex]north east}]\ar[ddddr,no head, end anchor={[xshift=-2.5ex]north east}]& &\\
		&& X_2\dlar[no head]\ar[->-=.5,no head,dddl, end anchor={[xshift=-.9ex]north east},start anchor={[xshift=-3.ex]south east }, end anchor={[yshift=-.5ex]north east}]\\
		& X_{3}\dar[no head, end anchor={[xshift=-2.5ex]north east},start anchor={[xshift=-2.5ex]south east }]\ar[-<-=.3,no head,dd, end anchor={[xshift=-0.9ex]north east},start anchor={[xshift=-0.9ex]south east }, end anchor={[yshift=-0.8ex]north east},start anchor={[yshift=0.8ex]south east }]&\\
		&X_4\dar[no head, end anchor={[xshift=-2.5ex]north east},start anchor={[xshift=-2.5ex]south east }]&\\
		&X_5& \quad\;.
	\end{tikzcd}

\end{center}
\end{example}

\subsection{Nilpotent and unipotent transformations}\index{nilpotent! -- tranformation}\index{unipotent! -- tranformation}

\begin{definition}\label{def_nilpotent_unipotent}
Let $V$ be a vector space, on which we denote by $\mathbb{I}$ the identity map. 
\begin{description}
\item[\ref{def_nilpotent_unipotent}.i.] We say that $A\in \mathfrak {gl}(V)$ is a {\em nilpotent transformation} of $V$ if there is $d\in \N$ such that $A^d\equiv 0.$ 
\item[\ref{def_nilpotent_unipotent}.ii.] We say that $B\in \mathfrak {gl}(V)$ is a {\em unipotent transformation} of $V$ if $B-\mathbb{I}$ is nilpotent.
\end{description}
\end{definition}

For the next result, recall that if $X \in \mathfrak {gl}(V)$ then the adjoint map for $x$, from Definition~\ref{def_ad}, is the element 
$\ad_X \in \mathfrak {gl}( \mathfrak {gl}(V) )$ such that
$\ad_X (Y):= [X, Y]=X Y -Y X$, for $y\in \mathfrak {gl}(V)$.
\begin{proposition}\label{prop4.6}
Let $V$ be a finite-dimensional vector space. If $X \in \mathfrak {gl}(V)$ is a nilpotent transformation of $V$, then $
\ad_X $
 is a nilpotent transformation of $\mathfrak {gl}(V)$. 
\end{proposition}

\begin{proof}
 First of all, one shows by induction from the definition of $\ad $ that there are constants $\{c_{k,j}\}_{j,k\in\N}\subset\Z$ such that
 \begin{equation}\label{eq6703a67f}
 \ad_X^kY = \sum_{j=0}^k c_{j,k} X^j Y X^{k-j}
 \qquad\forall X,Y\in\gl(V),\forall k\in\N.
 \end{equation}
 Next, let $X\in\gl(V)$ and $d\in\N$ such that $X^d=0$.
 Notice that for all $j\in\{0,\dots,2d\}$ we have $j\ge d$ or $2d-j\ge d$.
 Therefore, for all $Y\in\gl(V)$, if we write $\ad^{2d}_XY$ as in~\eqref{eq6703a67f}, we see that
 $X^j Y X^{2d-j}=0$ for all $j$, and thus $\ad^{2d}_XY=0$.
\end{proof}

%
%

The following proposition strengthens the result that nilpotent matrices only have 
 $0$ as an eigenvalue.
\begin{proposition}\label{0 only eigenvalue}
Let $V$ be a finite-dimensional vector space. If $A \in \mathfrak {gl}(V)$ is a nilpotent transformation, then
 there is a basis of $V$ in which $A$
is strictly upper triangular.
 \end{proposition}

%
%

\begin{proof}
 Given $A\in\gl(V)$, there is a basis of $V$ such that the representation of $A$ in this basis is in real Jordan form; see Exercise~\ref{Jordan Theorem}.
 Notice that $A$ is nilpotent if and only if each block of the Jordan representation of $A$ 
 is nilpotent.
 Moreover, every such block 
 is nilpotent if and only if its eigenvalue is zero.
 Therefore, the real Jordan form of $A$ is upper triangular.
\end{proof}

\subsection{Engel's Theorem}\label{sec Engel}
In this section, we prove Engel's Theorem. We begin by recalling the notion of representation of a Lie algebra, then we prove the statement for linear Lie algebras and, finally, the desired theorem.

 \begin{definition}\index{representation! -- of a Lie algebra}\index{Lie algebra! -- representation}
A {\em representation} of a Lie algebra $\g$ on a vector space $V$ is a Lie algebra homomorphism from $\g$ to $\mathfrak {gl}(V).$ Equivalently, it is a $\g$-module structure on $V,$ i.e., the map $\g \times V \to V, (x,v) \mapsto xv$ is bilinear and it holds
\begin{equation*}
[x,y]v = x(yv)-y(xv), \quad \forall x,y \in \g, \forall v\in V.
\end{equation*}
\end{definition}

The adjoint map $\ad: \g \to \gl (\g)$, which is a Lie algebra representation because of Jacobi identity, may not be injective.
In fact, its
 kernel is exactly the center $Z(\g)$ of $\g$ and thus $\ad(\g):=\{ \ad_x \,:\, x\in \g\}\simeq \g/Z(\g).$
 
 \begin{remark}\label{rmk_rometta_dic_2023}
 Another example of Lie algebra representation is given by subalgebras:
 if $\mathfrak {h}$ is a subalgebra of $\g$, then we have a representation of $\mathfrak {h}$ on $\g/ \mathfrak {h}$ as
\begin{equation*}
\ad_{\g/ \mathfrak {h}}:\mathfrak {h} \to \mathfrak {gl} (\g/ \mathfrak {h}),
\end{equation*}
defined as
\begin{equation}\label{def_ad_gh}
\ad_{\g/ \mathfrak {h}} (h) (y+ \mathfrak {h}):= [h,y] + \mathfrak {h}, \quad \forall h\in \mathfrak {h}, \forall y\in \g. 
\end{equation}
Indeed, note that, if $h\in \mathfrak {h},$ then $[h, \mathfrak {h}] \subseteq \mathfrak {h},$ and so $[h, y+\mathfrak {h}] + \mathfrak {h} =[h,y] +\mathfrak {h}$, leading to a well-defined $\ad_{\g/ \mathfrak {h}} $.
Moreover, we stress that if for $x\in \mathfrak {h} $ the transformation $\ad(x): \g \to \g $ is nilpotent, then $\ad_{\g/ \mathfrak {h}}(x): \g/ \mathfrak {h} \to \g/ \mathfrak {h}$ is nilpotent. Indeed, we have $(\ad_{\g/ \mathfrak {h}}(x))^k (y+ \mathfrak {h}) = (\ad(x))^k(y) + \mathfrak {h}.$
\end{remark}

\subsubsection{Engel's Theorem on linear Lie algebras} 
In the next theorem, we see that every Lie algebra of nilpotent transformations is a subalgebra of a Lie algebra of the form $ \g_{\rm nil} (\mathcal F)$ as defined in 
Example~\ref{ex g nil F}.
 \begin{theorem}[Engel's Theorem on linear Lie algebra]\label{Engel_linear_Lie}\index{Engel Theorem}\index{Theorem! Engel --}
Let $V\ne \{0\}$ be a finite-dimensional vector space and $\g \subseteq \mathfrak {gl}(V)$ a subalgebra. Assume that every $x\in \g$ is a nilpotent transformation of $V$. Then
 \begin{description}
\item[\ref{Engel_linear_Lie}.i.] there is $v_0\in V \setminus \{0\}$ such that $\g(v_0) = \{0\}$;
\item[\ref{Engel_linear_Lie}.ii.] there is a flag $\mathcal F= (V_0, \dots, V_n)$ for $V$ with $\dim(V_k)=k$ and $\g \subseteq \g_{\rm nil} (\mathcal F).$ 
\end{description}
Consequently,
\begin{description}
\item[\ref{Engel_linear_Lie}.iii.] there is a basis for $V$ relative to which elements of $\g$ are strictly upper-triangular matrices;
\item[\ref{Engel_linear_Lie}.iv.] $\g$ is a nilpotent Lie algebra.
\end{description}
\end{theorem}

\begin{proof}
We begin by proving \ref{Engel_linear_Lie}.i. It is proved by induction on the dimension of $\g.$

If $\dim \g=0,$ then any $v_0\in V\setminus\{0\}$ works. If $\dim\g=1,$ then $\g =\R x$ for some $x\in \mathfrak {gl} (V).$ Since $x$ is nilpotent, it has $0$ as its (only) eigenvalue (as shown in Proposition~\ref{0 only eigenvalue} by Jordan decomposition). Thus there is $v\ne 0$ such that $xv=0$ and so $txv=0$, for all $t\in \R$.

Now, assume that $\dim\g>1.$ Pick a subalgebra $\mathfrak {h}<\g $ with $\mathfrak {h} \ne \g$ of maximal dimension. Note that $\dim \mathfrak {h} \geq 1$ because 1-dimensional subalgebras always exist. By assumption, every $x\in \g \subset \mathfrak {gl}(V)$ is nilpotent. Hence, as previously shown, the transformation $\ad(x) \in \mathfrak {gl} (\mathfrak {gl} (V))$ is nilpotent (see Proposition~\ref{prop4.6}) and, as previously noted in Remark~\ref{rmk_rometta_dic_2023}, the transformation $\ad_{\g/\mathfrak {h}}(x) \in \mathfrak {gl}(\g/\mathfrak {h})$, as defined in \eqref{def_ad_gh}, is nilpotent. 

Now, since $\dim(\ad_{\g/\mathfrak {h}} (\mathfrak {h})) \leq \dim(\mathfrak {h}) <\dim (\g), $ we apply the inductive hypothesis to the representation $\ad_{\g/\mathfrak {h}}: \mathfrak {h} \to \mathfrak {gl}(\g/\mathfrak {h}).$ Then there is a non-zero element in $\g/ \mathfrak {h}$ say $x_0+\mathfrak {h}$ with $x_0\in \g$ and $x_0 \notin \mathfrak {h}$ such that $\ad_{\g/\mathfrak {h}} (\mathfrak {h}) (x_0+\mathfrak {h})=\mathfrak {h}.$ In other words, we have $[\mathfrak {h}, x_0] \subset \mathfrak {h}.$ Hence, the vector space $\mathfrak {h} + \R x_0$ is a subalgebra of $\g.$ By maximality of $\mathfrak {h},$ we infer $\mathfrak {h} +\R x_0 =\g $ and so $[\g, \mathfrak {h}] \subset \mathfrak {h},$ i.e., $\mathfrak {h}$ is an ideal. 

Next, we consider the representation $\mathfrak {h} \to \mathfrak {gl}(V)$ given by the inclusion. Recall that since $\dim(\mathfrak {h}) <\dim(\g), $ by the inductive hypothesis there is $v\in V\setminus \{0\}$ such that $\mathfrak {h}(v) = \{0\}$ and so we can consider the nontrivial subspace 
\begin{equation*}
V_0:=\{ v\in V \, :\, \mathfrak {h} (v)= \{0\}\}.
\end{equation*}
Note that
\begin{equation*}
\g(V_0) \subset V_0,
\end{equation*}
because for every $x\in \g,v\in V_0$, and every $y\in \mathfrak {h}$, we have that
\begin{equation*}
yxv = xyv -[x, y]v \in x\mathfrak {h} v +\mathfrak {h} v =\{0\}.
\end{equation*}
Here, we used that $[x, y] \in [\g, \mathfrak{h}] \subset \mathfrak{h}$ and that $\mathfrak{h}v=\{0\}$. 
Therefore, we can apply the base case of induction in dimension 1 to the representation $\R x_0 \to \mathfrak{gl}(V_0)$ defined as $tx_0 \mapsto (tx_0)_{|_{V_0}}$:
There is $v_0 \in V_0 \setminus \{0\} $ such that $x_0v_0 =0$ and thus, putting all together, we have
\begin{equation*}
\g (v_0)=\mathfrak{h}v_0 + \R x_0 v_0 =\{0\}.
\end{equation*}
This proves \ref{Engel_linear_Lie}.i.

Next, we show \ref{Engel_linear_Lie}.ii. By \ref{Engel_linear_Lie}.i, we choose $v_1 \in V \setminus \{0\}$ such that $\g (v_1)=0$ and we define $V_1 := \R v_1.$ Then, the map $\alpha : \g \to \mathfrak {gl} (V/V_1)$ defined as $\alpha (x)(v+V_1):= x(v)+V_1$ is well defined and gives a representation of $\g$ on $V/V_1.$ Still, the set $\alpha (\g)$ consists of nilpotent transformations. By induction on the dimension of $V$, it follows that $V/V_1$ possesses a complete flag $\mathcal F_1:= (W_0, \dots, W_{n-1})$ with $\alpha (\g) \subseteq \g_{\rm nil} (\mathcal F_1).$ Then, the singleton $\{0\}$ together with the preimage of the flag $\mathcal F_1$ in $V$ yields a complete flag $\mathcal F$ in $V$ with $\g \subseteq \g_{\rm nil} (\mathcal F).$ The consequences at the end of the statement of the theorem are immediate.
\end{proof}

\subsubsection{Engel's characterization theorem for nilpotent Lie algebra} 

 \begin{theorem}[Engel's theorem]\label{Engel's theorem}\index{Engel Theorem}\index{Theorem! Engel --}
 Let $\g$ be a finite-dimensional Lie algebra. Then
 $$\g \text{ is nilpotent}
\qquad \Longleftrightarrow \qquad \ad_x \text{ is nilpotent, for every } x\in \g.$$
\end{theorem}

\begin{proof}
 $[\Rightarrow]$ The forward direction has already been proved in Proposition~\ref{prop_basic_properties_nilpotent_Lie_albegra}.iv.
 
 $[\Leftarrow ]$ Regarding the backward direction, considering $\ad(\g) \subseteq \mathfrak {gl} (\g),$ we have that $\ad(\g) \cong \g / Z(\g).$ Moreover, the set $\ad(\g)$ is a Lie algebra of nilpotent transformations of $\g$, by assumption. By Engel's Theorem on linear Lie algebras, Theorem~\ref{Engel_linear_Lie}, we have that $\ad(\g)$ is nilpotent. Finally, by Proposition~\ref{prop_basic_properties_nilpotent_Lie_albegra}.ii, we get that $\g$ is nilpotent, as desired.
\end{proof}

In the above result, when using Proposition~\ref{prop_basic_properties_nilpotent_Lie_albegra}.ii, it is important that we are quotienting by something central.
The reader should be aware of the following fact: if $\mathfrak{a}\subseteq \g$ is an ideal and both $\mathfrak{a}$ and $ \g/\mathfrak{a}$ are nilpotent, then $\g$ may not be nilpotent; see Exercise~\ref{ex solvable not nilpotent}.

\subsection{The general Birkhoff-Embedding Theorem}\label{sec_Birkhoff}

 \begin{theorem}[Birkhoff-Embedding Theorem]\label{thm Birkhoff general}\index{Birkhoff-Embedding Theorem}\index{Theorem! Birkhoff-Embedding --}
Let $\g$ be a nilpotent finite-dimensional Lie algebra. Then there are a finite-dimensional vector space $V$ and an injective homomorphism $\iota: \g \to \mathfrak {gl} (V)$ such that for every $x\in \g$ the transformation $\iota(x)$ is nilpotent.
 \end{theorem}

We will only prove this theorem in the particular case of positively graded Lie algebras, and, hence, for Carnot algebras; see Section~\ref{sec Birkhoff for stratified}. A proof of the general theorem can be found in \cite[Theorem 1.1.11]{Corwin-Greenleaf}, but it relies on the construction of the universal enveloping algebra and the Poincar\'e-Birkhoff-Witt Theorem; see \cite[Chapter III]{Knapp}.
There is a more general result due to Ado, which generalizes the previous Theorem~\ref{thm Birkhoff general}, stating that every finite-dimensional Lie algebras has an injective finite-dimensional representation whose restriction to the maximal nilpotent ideal is nilpotent valued; see \cite[Theorem~B.8, p.663]{Knapp}.
 
 We recall that we also have Engel's characterization of Lie algebras of nilpotent transformations: Theorem~\ref{Engel_linear_Lie}. Hence, Birkhoff's Theorem can be stated as follows:
 \begin{corollary}\label{cor_nilpotent_triangular}
For every nilpotent finite-dimensional Lie algebra $\g$ there are $n\in \N$ and a subalgebra $\tilde \g$ of the space $ {\mathfrak{nil}}_n$ of strictly upper-triangular $n\times n$-matrices such that $\g$ and $\tilde \g$ are isomorphic.
 \end{corollary}
 
\section{Gradings and stratifications}
Carnot groups, the tangents of Carnot-Carath\'eodory spaces, will have nilpotent Lie algebras that have very special structures, called stratifications. For a Lie algebra $\g$, an \emph{$s$-step stratification of $\g$} is a direct-sum decomposition
\[
\g = V_1\oplus V_2\oplus\dots\oplus V_s
\]
of $\g$ with the property that 
\begin{equation}\label{def_stratific_step_s}
V_s\neq\{0\}\qquad \text{ and } \qquad [V_1,V_j]=V_{j+1}, \quad \forall j\in \{1, \dots, s\},
\end{equation}
 where we set $V_{s+1}:=\{0\}$. In Definition~\ref{defiCARNOT}, we shall provide other alternative equivalent definitions. In fact, it is useful to see a stratification as a special type of grading.
Hence, we discuss this broader concept next.

\subsection{Graded vector spaces and graded Lie algebras}

 We begin with graded vector spaces.
 
 \begin{definition}[Grading for a vector space]\label{def_linear_grading}
Let $A$ be a set (or, in many situations, an abelian group such as $\Z$ or $\R$) and let $V$ be a vector space. A {\em linear grading} of $V$ over $A$ is a collection of vector subspaces $(V_a)_{a\in A}$ of $V$ such that 
\index{grading! -- of a vector space}
\index{grading! linear --}
\index{linear grading}
\begin{equation*}
V= \bigoplus_{a\in A} V_a.
\end{equation*}
This means that $V=\Span\{V_a \,:\, a\in A \}$ and for every $a,a' \in A$ with $a\ne a'$ we have that $V_a \cap V_{a'} = \{0\}.$
We refer to a linear grading of $V$ also as a {\em grading of $V$ as a vector space}, or as a {\em grading}, for short.
When a grading of $V$ over $A$ is fixed, we shall say that $V$ is an {\em $A$-graded vector space}. If the grading is such that $A \subseteq \R$ and \index{graded! -- vector space}\index{vector space! graded --}
\begin{equation*}
V=V_{>0}:= \bigoplus_{a>0} V _a,
\end{equation*}
then $V$ is said to be {\em positively graded}. Given a grading $(V_a)_{a\in A}$ and $a\in A$, elements in $V_a$ are said to have {\em degree} $a$. Each $V_a$ is called a {\em layer}, or the {\em $a$-th layer}.\index{graded! positively --}\index{positively graded! -- vector space}\index{degree}\index{layer}\end{definition}

Next, we examine Lie algebras, for which we consider two more restrictive notions of linear gradings. We stress that a Lie algebra $\g$ is a vector space with the additional structure of a Lie bracket. 

 \begin{definition}[Compatible linear grading]\label{def_compatible_linear_grading}\index{compatible linear grading}\index{grading! compatible linear --}\index{linear grading! compatible --}
 Let $\mathfrak{g}$ be a Lie algebra.
A \textit{compatible linear grading} on $\mathfrak{g}$ is a linear decomposition of $\mathfrak{g}$ into vector subspaces $V_1, V_2, \ldots$ such that 
\begin{align}\label{def_compatible_lin_grad}
 \mathfrak{g} = \bigoplus_{i=1}^\infty V_i \quad \text{ and } \quad \mathfrak{g}^{(i)} = \mathfrak{g}^{(i+1)} \oplus V_i, \quad \forall i\in\N,
\end{align}
where the $\mathfrak{g}^{(i)} $'s are the lower central series elements from Definition~\ref{def lower central series}.
\end{definition}
A compatible linear grading is a particular $\Z$-grading with a mild interaction with the Lie algebra structure.
Clearly, every nilpotent Lie algebra admits a compatible linear grading (see Exercise~\ref{ex Nilpotency and compatible linear gradings}). Carnot algebras will have a stronger property.

 \begin{definition}[Lie algebra grading]
 Given an abelian group $A$ (for instance $\Z$ or $\R$) and a Lie algebra $\g$, 
a {\em grading over $A$} of $\g$ {\em as a Lie algebra}, or a {\em Lie algebra $A$-grading},
 is a linear grading $(V_a)_{a\in A}$ of $\g$ as a vector space 
 with the extra property that\index{grading! -- of a Lie algebra}\index{Lie algebra! -- grading}
\begin{equation*}
[V_a, V_b] \subseteq V_{a+b}, \quad \forall a,b \in A.
\end{equation*}
\end{definition}

\subsection{Stratified Lie algebras}\label{sec_strat}
We shall focus on a particular type of Lie algebra grading: stratifications. There are various equivalent definitions for them; see Remark~\ref{other def stratification}.
 \begin{definition}[Stratification]\label{defiCARNOT}
If $\g$ is a Lie algebra, a Lie algebra $\Z$-grading $(V_a)_{a\in \Z}$ of $\g$ is called a {\em stratification} of $\g$ if the smallest Lie subalgebra of $\g$ containing $V_1$ is $\g.$ The maximal $a$ for which $V_a\neq\{0\}$ is called the {\em step} of the stratification. A Lie algebra is {\em stratifiable} if it admits a stratification.
When one fixes a stratification of a stratifiable Lie algebra $\g$, we say that $\g$ is \emph{stratified},
or a {\em Carnot algebra}. 
\index{stratification}
\index{stratified! -- Lie algebra}
\index{stratifiable! -- Lie algebra}
\index{Lie algebra! stratified --}
\index{Lie algebra! stratifiable --}
\index{Carnot! -- algebra}
\index{step! -- of a stratification}
\index{stratification! step of a --}
\end{definition}
\begin{remark}[Equivalent definitions]\label{other def stratification}
We rephrase the definition by saying that 
a stratification of a Lie algebra $\g$ is a $\Z$-grading for which $\g$ is Lie generated by the elements of degree 1. Equivalently, this means that 
there is a direct-sum decomposition $\g = V_1 \oplus \dots \oplus V_s$ for which
\begin{equation*}
[V_1, V_j]=V_{j+1}, \quad \forall j\in\{1, \dots, s\}, \text{ with } V_{s+1}:=\{0\}.
\end{equation*}
The latter is the most common version of the definition of stratification.
\end{remark}

\begin{example}
	Every commutative Lie algebra $\g$ admits a 1-step stratification with $V_1=\g$.
 	\end{example}
\begin{example}	
	Let $\g$ be the Heisenberg Lie algebra spanned by $X, Y, Z$ with relation $[X, Y]=Z$.
	Then the subspaces $V_1:=\Span\{X,Y\}$ and $V_2:=\Span\{Z\}$ form a step-2 stratification.
\end{example}

 \begin{remark} 
The following {\bf non invertible} implications hold for finite-dimensional Lie algebras:
\begin{equation*}
 \text{ Carnot } \quad\stackrel{\nLeftarrow}{\THEN} \quad\text{ positively\, graded} \quad\stackrel{\nLeftarrow}{\THEN} \quad\text{ nilpotent }.
\end{equation*}
For the forward implications, see Exercises~\ref{ex stratifications are gradings} and \ref{ex positive gradings give nilpotency}.
For the fact that reverse implications do not hold, see Exercises~\ref{ex N51} and \ref{ex 5D non-gradable}. In fact, more examples can be found in \cite{Goodman, MR4394318}.
\end{remark}

\subsubsection{Uniqueness of stratifications}


 In the following proposition, we prove that every two stratifications on the same stratifiable Lie algebra differ by an automorphism.
\begin{proposition}[Uniqueness of stratifications]\label{isomorphic:stratifications}\index{stratification! uniqueness of --}
 	Let $\g$ be a stratifiable Lie algebra with two stratifications,
	\[
	V_1\oplus\dots\oplus V_s = \g = W_1\oplus\dots\oplus W_t .
	\]
	Then, we have the following properties:
	\begin{description}
	\item[\ref{isomorphic:stratifications}.i.] The steps coincide, i.e.,	$s=t$,
and $\g^{(k)}=V_k\oplus\dots\oplus V_s = W_k\oplus\dots\oplus W_s$, for all $k\in\{1, \ldots,s\}$;
	\item[\ref{isomorphic:stratifications}.ii.] 	there is a Lie algebra automorphism $A:\g\to\g$ with $A(V_k)=W_k$, for all $k\in\{1, \ldots,s\}$.
	\end{description}
\end{proposition}
\begin{proof}
 	The first point is simple; see Exercise~\ref{lem051025}.
	Then the quotient mappings $\pi_k:\g^{(k)}\to \g^{(k)}/\g^{(k+1)}$ induce linear isomorphisms $\pi_k|_{V_k}:V_k\to \g^{(k)}/\g^{(k+1)}$ and $\pi_k|_{W_k}:W_k\to \g^{(k)}/\g^{(k+1)}$, by a dimension argument.
	For $v\in V_k$ define $A(v):= (\pi_k|_{W_k})^{-1}\circ \pi_k|_{V_k}(v)$. 
	Notice that for $v\in V_k$ and $w\in W_k$ we have 
	\[
	A(v)=w \quad\iff\quad v-w\in \g^{(k+1)}.
	\]
	Extend $A$ to a linear map $A:\g\to\g$. 
	This is clearly a linear isomorphism.
	Next, we need to show that $A$ is a Lie algebra homomorphism, i.e., $[Aa,Ab]=A([a,b])$ for all $a,b\in \g$.
	Let $a=\sum_{i=1}^sa_i$ and $b=\sum_{i=1}^sb_i$ with $a_i, b_i\in V_i$.
	Then
$$
	 	A([a,b]) = \sum_{i=1}^s\sum_{j=1}^s A([a_i,b_j]) 
	\quad \text{ and } \quad
		[Aa,Ab] = \sum_{i=1}^s\sum_{j=1}^s [Aa_i,Ab_j].
$$
	Therefore we can just prove $A([a_i,b_j])=[Aa_i,Ab_j]$ for $a_i\in V_i$ and $b_j\in V_j$.
	Notice that $A([a_i,b_j])$ and $[Aa_i,Ab_j]$ both belong to $W_{i+j}$. 
	Therefore, writing $[a_i,b_j]$ as
	$[Aa_i,Ab_j]+ ([a_i,b_j]- [Aa_i,Ab_j])$, 
	 we have $A([a_i,b_j])=[Aa_i,Ab_j]$ if and only if 
	$[a_i,b_j] - [Aa_i,Ab_j]\in\g^{(i+j+1)}$.
	And in fact
	\[
	[a_i,b_j] - [Aa_i,Ab_j] 
	= [a_i-Aa_i,b_j] - [ Aa_i,Ab_j-b_j] \in \g^{(i+j+1)},
	\]
	because, on the one hand, $a_i-Aa_i\in \g^{(i+1)}$ and $b_j\in W_j$, so $[a_i-Aa_i,b_j]\in\g^{(i+j+1)}$, on the other hand, $Aa_i\in W_i$ and $Ab_j-b_j\in \g^{(j+1)}$, so $[ Aa_i,Ab_j-b_j] \in \g^{(i+j+1)}$.
	We deduced that $A$ is a Lie algebra homomorphism.
\end{proof}


\subsubsection{Induced grading on $\mathfrak {gl} (V)$}\index{induced! -- grading}\index{$\mathfrak {gl} (V)$}
In the next definition, we consider subspaces of the linear transformations on a vector space equipped with a linear grading that will form a grading themselves. 
 \begin{definition}[$\mathfrak {gl} (V)_a $ and $\mathfrak {gl} (V)_{>0}$]\index{$\mathfrak {gl} (V)_a $}
 Let $(V_a)_{a\in A}$ be a linear grading of a vector space $V$ over an abelian group $A.$ For every $a\in A$ we define
\begin{equation*}
\mathfrak {gl} (V)_a := \{ M \in \mathfrak {gl} (V)\,:\, M(V_b) \subseteq V_{b+a}, \forall b\in A \},
\end{equation*}
and if $A$ has an ordering (for instance, if $A\subseteq \R$) we define
\begin{equation*}
\mathfrak {gl} (V)_{>0} := \bigoplus_{a>0} \, \mathfrak {gl}(V)_a.
\end{equation*}
\end{definition}
We shall show that the collection $(\mathfrak {gl}(V)_a) _{a\in A}$ form an $A$-grading of $\mathfrak {gl}(V)$ as a Lie algebra. 
Moreover, if $V$ is a Carnot algebra, then $\mathfrak {gl} (V)_{>0} $ is a Carnot algebra.
\begin{proposition}\label{glVa_properties}
Let $V$ be a finite-dimensional vector space with a linear grading $(V_a)_{a\in A}$ over an abelian group $A.$ Then, we have the following properties:
\begin{description}
\item[\ref{glVa_properties}.i.] $(\mathfrak {gl}(V)_a) _{a\in A}$ is a Lie algebra grading of $\mathfrak {gl}(V)$.
\item[\ref{glVa_properties}.ii.] If $A<( \R, +)$, then $\mathfrak {gl} (V)_{>0} $ is a Lie subalgebra of nilpotent transformations.
\item[\ref{glVa_properties}.iii.] If $A=\Z$ and there are $\bar a, \bar b \in A$ such that
\begin{equation*}
V_m \ne\{0\}\quad \Leftrightarrow \quad m\in \Z \cap [\bar a, \bar b],
\end{equation*}
then $\mathfrak {gl} (V)_{>0} $ is a Carnot algebra.
\end{description}
\end{proposition}

\begin{proof}
(i). Fix $X_1, \dots, X_n$ a basis of $V$ adapted to the direct-sum decomposition $V= \bigoplus_{a\in A} V_a.$ So for every $i\in\{1, \dots, n\}$ there is $a_i \in A$ such that $X_i \in V_{a_i}.$ For each $i,j \in \{ 1, \dots, n\},$ let $E^i_j \in \mathfrak {gl}(V)$ be such that 
 $$E^i_j(X_i):=X_j \qquad \text { and }\qquad
 E^i_j(X_k) :=0, \quad\forall k\ne i.$$
Consequently, the elements $(E^i_j)_{ i,j \in\{1, \dots, n\}}$ form a basis of $\mathfrak {gl}(V).$ Moreover, every $E^i_j$ is such that $E^i_j (V_{a_i}) \subseteq V_{a_j}$ and for every $a\ne a_i$ we have that $E^i_j (V_a)=\{0\} \subseteq V_{a+(a_j-a_i)},$ i.e., $E^i_j\in \mathfrak {gl}(V)_{a_j-a_i}.$ Therefore, 
\begin{equation*}
\bigoplus_{a\in A}\, \mathfrak {gl}(V)_a = \mathfrak {gl}(V).
\end{equation*}
Finally, notice that if $M_1 \in \mathfrak {gl}(V)_a$ and $ M_2 \in \mathfrak {gl}(V)_b$ then
\begin{equation*}
M_1M_2, \quad M_2M_1 \in \mathfrak {gl}(V)_{a+b},
\end{equation*}
and so $[ \mathfrak {gl}(V)_{a}, \mathfrak {gl}(V)_{b}]\subseteq \mathfrak {gl}(V)_{a+b}.$ This completes the proof of \ref{glVa_properties}.i.

(ii). Clearly, the subset $ \mathfrak {gl}(V)_{>0}$ is a Lie subalgebra,
see Exercise~\ref{V>0}.
 Moreover, assuming $a_i\leq a_j$ for every $i\leq j$ and defining 
\begin{equation*}
W_0:=\{ 0 \} \qquad \text{ and } \qquad
W_i := \Span \{ X_{n-i+1}, \dots, X_n\}, \quad \forall i\in\{1, \ldots,n\},
\end{equation*}
then $\mathcal F:=(W_0, \ldots, W_m)$ defines a flag for $V$ such that $ \mathfrak {gl}(V)_ {>0} \subseteq \mathfrak {g}_ {\rm nil}(\mathcal F)$; see the last notation from Example~\ref{ex g nil F}.

(iii). We want to prove that if $V=V_{\bar a}\oplus \dots \oplus V_{\bar b}$ with $V_j \ne \{0\}$ for every integer $j={\bar a}, \dots, \bar b$, then $ \mathfrak {gl}(V)_ {>0}$ is generated by $ \mathfrak {gl}(V)_1.$
By induction, we shall prove that
\begin{equation}\label{equation7aprile}
[\mathfrak {gl} (V)_1, \mathfrak {gl}(V)_k ] =\mathfrak {gl}(V)_ {k+1}, \qquad \forall k\in\N.
\end{equation}
 It is enough to prove that for every $i,j$ such that $a_j-a_i =k+1$ (i.e., $E^i_j \in \mathfrak {gl} (V)_{k+1}$) we get that
\begin{equation*}
E^i_j \in [\mathfrak {gl} (V)_1, \mathfrak {gl}(V)_k ]. 
\end{equation*}
Since $X_i \in V_{a_i}, X_j \in V_{a_j}$, and $a_i-a_j=k+1 \ne 0$ we have that ${\bar a}\leq a_i < a_j\leq \bar b.$ Moreover, by assumption, there is a basis element $X_\ell$ with $a_\ell =a_i +1.$ Then we claim that
\begin{equation*}
-E^i_j=[E^i_\ell, E^\ell_j] = E^i_\ell E^\ell_j - E^\ell_j E^i_\ell. 
\end{equation*}
Indeed, we have $E^i_\ell E^\ell_j=0$ and $E^\ell_j E^i_\ell =E^i_j,$ recalling that $a_i<a_\ell \leq a_j.$ Since $E^i_\ell \in \mathfrak {gl} (V)_1$ and $E^\ell_j \in \mathfrak {gl} (V)_k$, we proved \eqref{equation7aprile}, as desired.
 \end{proof}

\subsection{Dilation structures}
For vector spaces, gradings are in correspondence with dilation structures.

We stress that, in the presence of an $\R$-grading $(V_a)_{a\in \R }$ of a
finite-dimensional vector space $V$, there is a finite subset $I\subset \R$ such that $V_i\neq \{0\}$ if and only if $i\in I$. Hence, every vector $v\in V$ can be written uniquely as $v=\sum_{i\in I}v_i$, with $v_i\in V_i$. With abuse of notation, we shall still write $v=\sum_{i\in \R}v_i$.
\begin{definition}[Dilations on graded vector spaces]\label{Dilations:algebras}
Let $(V_a)_{a\in \R }$ be an $\R$-grading of a finite-dimensional vector space $V$.
For every $\lambda>0$, the \emph{inhomogeneous dilation on $V$ of factor $\lambda$ relative to the grading} (or, simply, the {\em dilation} $\delta_\lambda$, when the grading is understood) is the linear map $\delta_\lambda:V\to V$ such that 
\index{dilation! -- in stratified algebra}
\begin{equation}\label{def_dilation_relative_to_grading}
\delta_\lambda v = \lambda^a v, \qquad \forall a\in \R, \forall v\in V_a. \end{equation}
\end{definition}

In addition, we notice that if the grading is a $\Z$-grading, then the above equation defines dilations $\delta_\lambda$, also for $\lambda<0$. If $V_{\leq0} $ is trivial, then for $\lambda=0$ we set $\delta_\lambda\equiv0$. Therefore, if  $V_{\leq0} =\{0\} $ and the grading is a $\Z$-grading, the map $(\lambda,v)\in \R\times V\mapsto \delta_\lambda(v)\in V$ is continuous.

\begin{proposition}\label{lemma OPS delta lambda} Given a Lie algebra grading $(V_a)_{a\in \R }$ of a Lie algebra $\g$, consider the inhomogeneous dilations relative to the grading, as in
\eqref{def_dilation_relative_to_grading}.
Then, for $\lambda\in \R_{>0}$,	the dilation $\delta_\lambda:\g\to\g$ is a Lie algebra automorphism.
	Moreover, the map $(\R_{>0}, \cdot)\to \Aut_{\Lie}(\g)$, $\lambda\mapsto\delta_\lambda$, is a one-parameter subgroup:
\begin{equation}\label{eq OPS delta lambda}
\delta_\lambda \circ \delta_\mu = \delta_{\lambda\mu}, \qquad \forall\lambda, \mu\in \R.
\end{equation}
\end{proposition}
\begin{proof}
We need to show that the map is a linear bijection and 
	$$\delta_\lambda([X, Y])=[\delta_\lambda X, \delta_\lambda Y], \qquad \forall X, Y\in \g.$$
 	Take $X,Y\in\g$ and decompose them as $X=\sum_{i\in \R} X_i$ and $Y=\sum_{i\in\R}Y_i$, with $X_i,Y_i\in V_i$.
	Since $[X_i,Y_j]\in[V_i,V_j]\subset V_{i+j}$, we get
	\begin{align*}
	 	[\delta_\lambda X, \delta_\lambda Y] 
		= \sum_{i,j} [\lambda^i X_i, \lambda^jY_j]
		= \sum_{i,j} \lambda^{i+j}[X_i,Y_j]
		= \sum_{i,j} \delta_\lambda([X_i,Y_j])
		= \delta_\lambda\left( \sum_{i,j} [X_i,Y_j] \right)
		= \delta_\lambda ([X, Y]).
	\end{align*}
	Moreover, the map $\delta_\lambda$ is invertible with inverse $\delta_{1/\lambda}$. Equation \eqref{eq OPS delta lambda} is trivial.
\end{proof}

 In the presence of a $\Z$-grading, the map $(\R\setminus\{0\}, \cdot)\to \Aut_{\Lie}(\g)$, $\lambda\mapsto\delta_\lambda$ is a subgroup homomorphism, also defined for negative $\lambda$'s.
 
 {\color{black}
 Vice versa, if we have a direct-sum decomposition $\g=\bigoplus_{a\in\R}V_a$ and the map $\delta_\lambda$ 
 as in
\eqref{def_dilation_relative_to_grading} is a Lie algebra automorphism, 
 then the decomposition is necessarily a Lie algebra grading. 
 Notice that \eqref{def_dilation_relative_to_grading} can be rewritten as 
 $$\delta_\lambda = \exp(\log(\lambda) \alpha),$$
 where $\alpha$ is the diagonal transformation multiplying by $a\in \R$ the space $V_a$.
 The assumption that the $\delta_\lambda$'s are Lie algebra automorphisms, rephrases as $\alpha$ being a derivation, in the sense of Definition~\ref{def:derivation}.
More generally, every OPS of automorphisms gives a grading.\index{derivation}
 \begin{proposition}\label{prop12061739}
Given a Lie algebra $\g$, let $(\R_{>0}, \cdot)\to \Aut_{\Lie}(\g)$, $\lambda\mapsto\delta_\lambda$, be a one-parameter subgroup of Lie algebra automorphisms. Then there exists $\alpha\in {\rm {Der}}(\g)$
such that 
$\delta_\lambda:=\exp((\log(\lambda) \alpha)$, for all $\lambda\in \R_{>0}$, and 
the Lie algebra $\g$ admits a Lie algebra grading where the layers are
	\[
	V_t := \g \cap\bigoplus_{s\in\R} E^{\alpha}_{t+is}, \qquad \text{ for }t\in\R, 
	\]
	where $E^{\alpha}_{t+is}$ is a generalized eigenspace of $\alpha$ corresponding to the eigenvector $t+is$; see \eqref{def_generalized_eigenspace}.
\end{proposition}
We leave the proof as an exercise since we will actually prove something stronger in a few pages; see Proposition~\ref{prop08220937}.

 }

Later, when discussing measures on Carnot groups, we shall need the following result.

\begin{lemma}
If $\g=V_1\oplus\dots\oplus V_s$ is a stratified Lie algebra, then 
the dilation $\delta_\lambda:\g\to\g$ as in
\eqref{def_dilation_relative_to_grading} has determinant equal to $\lambda^Q$ with 
	\[
	Q=\sum_{j=1}^s j\cdot\dim(V_j) .
	\]
\end{lemma}
\begin{proof}
 	Fix a basis $X_1, \dots, X_n$ adapted to the stratification, i.e., for every $i$, there is a $j$ such that $X_i\in V_j$.
	Then, in this basis, the map $\delta_\lambda$ is represented by the diagonal matrix with diagonal 
	\[
	(\underbrace{\lambda, \dots, \lambda}_{\dim V_1}, \underbrace{\lambda^2, \dots, \lambda^2}_{\dim V_2}, \dots, \underbrace{\lambda^s, \dots, \lambda^s}_{\dim V_s}) .
	\]
	Hence, the determinant is $\lambda^{\dim V^1}\cdot(\lambda^2)^{\dim V_2}\dots\cdot(\lambda^s)^{\dim V_s} = \lambda^Q$.
\end{proof}

\subsubsection{Associated Carnot algebra}\label{sec_associated_Carnot}

To every nilpotent Lie algebra $\g$, there is a canonical way to associate a stratified Lie algebra. Later in the book, in Section~\ref{sec:Pansu asymptotic structure}, we shall see that this associated Lie algebra also has a geometric meaning. We first define this associated Carnot algebra as an abstract direct sum of quotients. Later, we shall equivalently describe it via limiting dilations.

\begin{definition}[Associated Carnot algebra]\label{Graded:algebra}\index{Carnot! associated -- algebra}\index{associated Carnot algebra}
Let $\g$ be a Lie algebra that is nilpotent of step $s$.
Let $\g^{( i+1)} \stackrel{\rm def}{=} [\g, \g^{( i)}]$ be the descending central series of $\g$.
The {\em associated Carnot algebra} of $\g$ is the Lie algebra $\g_\infty$ given by the direct-sum decomposition\index{$\g_\infty$}
$$\g_\infty := \bigoplus_{i=1}^s\g^{( i)} /\g^{( i+1)}, $$
endowed with the unique Lie bracket $\llbracket\cdot, \cdot\rrbracket_\infty$ that has the property that, if $X\in\g^{( i)}$ and $Y\in\g^{(j)}$, the bracket is defined, modulo $\g^{( i+j+1)}$, as
$$ \left\llbracket X+ \g^{(i+1)},Y+ \g^{(j+1)}\right\rrbracket_\infty:= [X, Y] + \g^{(i+j+1)}.$$
\end{definition}

The associated Carnot algebra is stratified by $\left(\g^{( i)} /\g^{( i+1)}\right)_{i=1}^s$. Fixing a compatible linear grading, one can setwise identify this new Lie algebra with the original one, and its Lie bracket can be equivalently defined by the following result.

\begin{lemma}\label{limit_bracket}
Let $(\g, [\cdot, \cdot])$ be a nilpotent Lie algebra on a vector space $\g$. 
Consider the dilations $(\delta_\lambda)_{\lambda>0}$ relative to some compatible linear grading and
define the map
\begin{equation}\label{formula-g-inf}
\llbracket X,Y\rrbracket_\infty :=\lim_{\lambda\to +\infty }\delta_\lambda^{-1}[\delta_\lambda X, \delta_\lambda Y], \qquad \forall X,Y\in \g.
\end{equation}
Then $\llbracket\cdot, \cdot\rrbracket_\infty$ defines a Lie bracket on $\g$, such that 
 $$\llbracket\delta_\lambda X, \delta_\lambda Y\rrbracket_\infty = \delta_\lambda \llbracket X,Y\rrbracket_\infty, \qquad \forall X,Y \in \g,$$
 and $(\g, \llbracket\cdot, \cdot\rrbracket_\infty)$ is isomorphic to the associated Carnot algebra of $ (\g, [\cdot, \cdot])$.
\end{lemma}
\begin{proof}
Let $(V_j)_{j=1}^\infty$ be a compatible linear grading of $(\g, [\cdot, \cdot])$.
Since the $V_j$'s form a direct decomposition of $\g$, it suffices to consider \eqref{formula-g-inf} for $X\in V_i$ and $Y\in V_j$, for some $i,j$. 
Because we are in the presence of a compatible linear grading, the element $[X, Y]$ belongs to $ \bigoplus_{k=i+j}^\infty V_k$; see Exercise~\ref{ex Nilpotency and compatible linear gradings}. Thus,
 we have that
$$[X, Y]=Z_{i+j}+Z_{i+j+1}+\ldots+Z_{s},$$
for some vectors $Z_k\in V_k$.
Hence, we compute
\begin{eqnarray*}
 \delta_\lambda^{-1}[\delta_\lambda X, \delta_\lambda Y]&=& \delta_\lambda^{-1}[\lambda^i X, \lambda^j Y]\\
 &=&\lambda^{i+j} \delta_\lambda^{-1}(Z_{i+j}+Z_{i+j+1}+\ldots+Z_{s})\\
 &=& Z_{i+j}+\lambda^{-1}Z_{i+j+1}+\ldots+\lambda^{i+j-s}Z_{s},
\end{eqnarray*}
which goes to $Z_{i+j}$, as
$\lambda\to\infty$. The proof is concluded by observing that $Z_{i+j}$ is a vector that represent $[X, Y]$ modulo $\g^{(i+ j+1)}$.
\end{proof}

\subsubsection{Siebert theorem}

Lie groups whose Lie algebra admits a positive grading are precisely those Lie groups whose universal covering Lie group admits automorphisms that are topologically contractive, in a sense that we soon review. Such a characterization is due to Siebert, \cite{Siebert}. We present here his result only focusing on Lie groups.
More generally, Siebert showed that a connected locally compact group $G$ admits a contractible automorphism if and only if $G$ is a simply connected Lie group whose Lie algebra admits a positive grading, \cite[Corollary~2.4]{Siebert}. In addition, a locally compact group $G$ admits a contractible automorphism if and only if $G$ is topologically isomorphic to the direct product of two groups $G_C$ and $G_D$ admitting contractible automorphisms where
$G_C$ is connected and $G_D$ is totally disconnected, \cite[Proposition~4.2 and Corollary~4.3]{Siebert}. We shall not discuss these more comprehensive results here.

We use the term contractive in the following topological sense. A map $F: X\to X$ from a topological space into itself is said to be {\em contractive}\label{page_contractive}\index{contractive} if there exists $x_0\in X$ such that for every $x\in X$ one has
$\lim_{n\to \infty} F^n(x)=x_0$, where $ F^n$ denotes the $n$-times composition of the map $F$.
We stress that if $F$ is a Lie group automorphism $G\to G$ or a Lie algebra automorphism $\g\to \g$, then the (only possible choice for) $x_0$ is $1_G$ or $0\in\g$, respectively. 

We saw that if a Lie algebra is graded, then the relative dilations from Definition~\ref{Dilations:algebras} define maps $\delta_\lambda$ that are Lie algebra automorphisms, for all $\lambda>0$; see Proposition~\ref{lemma OPS delta lambda}. We stress that if $\lambda\in (0,1)$ and the grading is a positive grading, then the map $\delta_\lambda$ is contractive; see Exercise~\ref{ex positively graded contraction}.
Siebert's theorem exactly states the inverse implication.

\begin{theorem}[Siebert, \cite{Siebert}]\label{thm_Siebert}\index{Siebert Theorem}\index{Theorem! Siebert --}
For every simply connected Lie group $G$,
 the following are equivalent.
\begin{description}
\item[\ref{thm_Siebert}.i.] $\Lie(G)$ admits a positive grading;
\item[\ref{thm_Siebert}.ii.] $G$ admits a contractive Lie group automorphism;
\item[\ref{thm_Siebert}.iii.] $\Lie(G)$ admits a contractive Lie algebra automorphism.
\end{description}
\end{theorem}

%

Before the proof of Theorem~\ref{thm_Siebert}, we discuss how automorphisms of Lie algebras induce gradings; see the following proposition. For defining the gradings, it is more convenient to pass to complexifications as we next review.

We denote by $V_\C$ the \emph{complexification} of a real vector space $V$, 
so 
$V_\C:=V\oplus V$ with complex scalar multiplication given by $i\cdot(X,Y) := 
(-Y, X)$, for $X,Y\in V$. 
It is a complex vector space with conjugation $(X, Y)^*:=(X,-Y)$.
If $\phi:V\to V$ is an $\R$-linear map, then its \emph{complexification} is the $\C$-linear map $\phi_\C:V_\C\to V_\C$, $(X,Y)\mapsto \phi_\C(X,Y) := (\phi (X), \phi (Y))$.
The \emph{spectrum} of $\phi$ is defined by\index{spectrum}
\[
	{\rm Spec}(\phi): = \{\alpha\in\C: \det(\phi_\C-\alpha\mathbb{I}) = 0 \},
\]
where $\mathbb{I}$ is the identity map on $ V_\C$,
and the \emph{generalized eigenspace} of $\phi$ corresponding to $\alpha\in\C$ by\index{generalized eigenspace}\index{$E^\phi_\alpha$}
\begin{equation}\label{def_generalized_eigenspace}
	E^\phi_\alpha: = \left\{v\in V_\C:\exists n\in\N \quad (\phi_\C-\alpha\mathbb{I})^nv = 0\right\} .
\end{equation}
By Jordan Theorem (see Exercise~\ref{Jordan Theorem}), we have $V_\C =  \bigoplus_{\alpha\in{\rm Spec}(\phi)} E^\phi_\alpha $.
\begin{proposition}\label{prop08220937}
	Let $\phi$ be an automorphism of a Lie algebra $\g $.
	For all $\lambda\in(0,+\infty)\setminus\{1\}$ and $t\in\R$, define 
	\[
	V_t 
	:= V_t(\lambda, \phi):= \g \cap \left( \bigoplus \left\{E^\phi_\alpha :\alpha\in \C, |\alpha|=\lambda^t\right\} \right), 
	\]
	where the $E^\phi_\alpha $ are defined in \eqref{def_generalized_eigenspace} for $V:=\g$.
	Then $\{V_t\}_{t\in\R}$ is a Lie algebra grading of $\g $, with degrees in $ \log_\lambda(|{\rm Spec}(\phi)|)$, i.e.,
	\begin{equation}
\label{eq19Aug20240853}
 \g = \bigoplus_{t\in D} V_t, \qquad \text{ for }D:=\{ \log_\lambda(|\alpha|) : \alpha \in {\rm Spec}(\phi)\}.
 	\end{equation}
	Moreover,
	\[
	|\det(\phi)| 
	= \lambda^{\sum_{t\in\R} t\cdot\dim(V_t)}.
	\]
\end{proposition}
\begin{proof}
We stress that if $\alpha \notin {\rm Spec}(\phi)=:\sigma(\phi)$, then $E^\phi_\alpha=\{0\}$.
Since $\sigma(\phi)$ is a finite set, then only finitely many $V_t$'s are not trivial.
	For $\alpha\in\sigma(\phi)\subset\C$, define $U^\phi_\alpha := ( E^{\phi}_\alpha+ E^{\phi}_{\bar\alpha} )\cap\g $.
	We claim that
	\begin{equation}\label{eq12061733}
	\g = \bigoplus_{\alpha\in\sigma(\phi)} U^\phi_\alpha, 
	\end{equation}
	where the sum is direct up to the identification $U^\phi_\alpha=U^\phi_{\bar \alpha}$.
	Indeed, take $v\in \G$ and write $v=\sum_{\alpha}v_\alpha\in\g $ with $v_\alpha\in E^\phi_\alpha$ for all $\alpha$, which is possible by Jordan theorem.
Since $v= v^*$, we have 	 $v=\frac12(v+v^*)=\frac12 \sum_{\alpha}(v_\alpha+ v_\alpha^*)$, where $v_\alpha+ v_\alpha^*\in U_\alpha$.
		So, we have $\g = \sum_{\alpha\in\sigma(\phi)} U^\phi_\alpha$.
	Since $U^\phi_\alpha\cap U^\phi_\beta=\{0\}$ if $\alpha\notin\{\beta, \bar\beta\}$, the sum is direct.
	This proves claim~\eqref{eq12061733}.
	
	Since $\phi$ is injective, 
	then $U^\phi_0=\{0\}$.
	Therefore, by~\eqref{eq12061733}, we have $\g = \bigoplus_{t\in\R} V_t$.
	
	Using Exercise~\ref{lemBurba}, we have
	\begin{equation}\label{eq12061734}
	[U^\phi_\alpha,U^\phi_\beta] \subset U^\phi_{\alpha\beta}\oplus U^\phi_{\bar\alpha\beta}, \qquad \forall \alpha, \beta\in\C.
	\end{equation}
	
	If $X\in U^\phi_\alpha$ and $Y\in U^\phi_\beta$ with $|\alpha|=\lambda^t$ and $|\beta|=\lambda^s$, then $[X, Y]\in U^\phi_{\alpha\beta}\oplus U^\phi_{\bar\alpha\beta} \subset V_{t+s}$, because of~\eqref{eq12061734} and $|\alpha\beta| = |\bar\alpha\beta| = \lambda^{s+t}$.
	Therefore, $[V_s,V_t] \subset V_{s+t}$ and $\{V_t\}_{t\in\R}$ is a real grading of $\g $.
	Finally, if we set $\varepsilon_\alpha=1$ if $\alpha\in\R$ and $\varepsilon_\alpha=1/2$ if $\alpha\in\C\setminus\R$,
	\[
	|\det(\phi)|
	= \left| \prod_{\alpha\in\sigma(\phi)} \alpha^{\dim_{\C}(E_\alpha)} \right| 
	= \prod_{\alpha\in\sigma(\phi)} |\alpha|^{\varepsilon_\alpha\dim_\R (U_\alpha)} 
	= \prod_{t\in\R} \lambda^{t \cdot\dim(V_t)} . \qedhere
	\]
\end{proof}

\begin{proof}[Proof of Theorem~\ref{thm_Siebert}]
The equivalence between 
 \ref{thm_Siebert}.ii and 
 \ref{thm_Siebert}.iii is trivial; see Exercise~\ref{ex contractive contractive}.
 The easy implication \ref{thm_Siebert}.i$\implies$\ref{thm_Siebert}.iii is in Exercise~\ref{ex positively graded contraction}.
 
 Regarding \ref{thm_Siebert}.iii$\implies$\ref{thm_Siebert}.i, if $\phi$ is an automorphism of $\g$ that is contractive, then also $\phi_\C:\g_\C\to \g_\C$ is contractive.
 Therefore, the transformation $\phi_\C$ cannot have eigenvalues with norm $\geq 1$.
 For $\lambda:=1/2$, we consider the grading $V_t 
	:= V_t(\lambda, \phi)$ of Proposition~\ref{prop08220937}.
	We check that it is a positive grading. Indeed, for all $t\leq0$ we have that for all $\alpha\in \C$ such that $|\alpha|=1/2^t$
the set $E^\phi_\alpha $ is trivial brcause $1/2^t\geq 1$.
Thus $V_t=\{0\}$.
\end{proof} 
 
We conclude the subsection with the following observation that tells us when an automorphism coming from a derivation is contractive.
\begin{remark}\label{spec_deriv_auto}
Let $A$ be a derivation of a Lie algebra $\g$. For $\psi:=\exp(A)\in \Aut_{\Lie}(\g)$ the following are equivalent:
\\ \ref{spec_deriv_auto}.i. $\psi$ is contractive;
\\ \ref{spec_deriv_auto}.ii. ${\rm Spec}(\psi)\subseteq \{z\in \C : |z|<1\}$;
\\ \ref{spec_deriv_auto}.iii. ${\rm Spec}(A)\subseteq \{z\in \C : \Re(z)<0\}$.
\end{remark} 
 
\subsection{Birkhoff Theorem for stratified Lie algebras}\label{sec Birkhoff for stratified}


We will present a proof of Birkhoff's Theorem for Carnot algebras, as Y. Cornulier explained it to the author.
We begin with a Carnot algebra $\g$. Then, we perform a semidirect product $\g \rtimes \R$, on which we naturally put a grading. Consequently, the Lie algebra $\mathfrak {gl} (\g \rtimes \R )_{>0}$ will be a Carnot algebra, containing a copy of $\g$.

\subsubsection{Induced grading on $\g \rtimes \R$}\index{induced! -- grading}
We begin very general: Let $\g$ be a $\Z$-graded Lie algebra with grading $(V_m)_{m\in \Z}.$ 
 We consider the semidirect product $\g \rtimes \R$ where $1\in \R$ acts on $\g$ as the derivation that multiplies by $m$ the vectors in $V_m$.  
Namely, recalling the construction from \eqref{def_semi_direct_algebra}, the Lie bracket on the semidirect product $\g \rtimes \R$ is
 \begin{equation*}
[(X,s), (Y,t)] \stackrel{\rm def}{=}\left([X, Y] + \sum_{m\in \Z} sm Y_m - \sum_{m\in \Z} tmX_m,0\right), \qquad \forall X,Y, \in \g, \forall s, t\in \R,
\end{equation*}
if $X= \sum_{m\in \Z} X_m$ and $Y=
 \sum_{m\in \Z} Y_m$ with $X_m, Y_m\in V_m$, for $m\in \Z$.
 
The Lie algebra $\g \rtimes \R$ is $\Z$-graded by $(V'_m)_{m\in \Z}$ defined as
\begin{eqnarray}\label{grading_g_R}
&&V'_0:=V_0 \times \{0\}\oplus \{0\}\times \R, 
\nonumber
\\
 &&V'_m :=V_m \times \{0\}, \qquad\forall m \ne 0; 
\end{eqnarray} 
 see Exercise~\ref{ex_grading_g_R}.i. 
Moreover, if $\g$ is a non-trivial Carnot algebra, then $\g \rtimes \R$ has a trivial center; see Exercise~\ref{ex_grading_g_R}.ii.

\subsubsection{Proof of Birkhoff Theorem for Carnot algebras}

When $\g$ is a $\Z$-graded Lie algebra, then $\mathfrak {gl} (\g \rtimes \R )_{>0}$ is a Lie algebra of nilpotent transformations and $\g \rtimes \R$ is graded by \eqref{grading_g_R}; see Exercise~\ref{ex_grading_g_R}. 
Note that
$\g \simeq \g \times \{0\} \subset \g \rtimes \R,$ so $\ad(\g)$ can be seen as a subset of $\mathfrak {gl} (\g \rtimes \R )$. We next prove that such a map $\ad$ gives an injective representation of $\g$ on the vector space $\g \rtimes \R $.

 \begin{theorem}[Birkhoff-Embedding Theorem for Carnot algebras]\label{Birkhoff Theorem for Carnot}\index{Birkhoff-Embedding Theorem}\index{Theorem! Birkhoff-Embedding --}
Let $\g$ be a Carnot algebra. Then $\ad:\g \to \mathfrak {gl} (\g \rtimes \R )_{>0} \subset \mathfrak {gl} (\g \rtimes \R )$ is an injective Lie algebra homomorphism into the Carnot algebra $\mathfrak {gl} (\g \rtimes \R )_{>0}$ of nilpotent transformations.
 \end{theorem}

 \begin{proof}
We can assume that $\g \ne \{0\}.$ Since $\g$ is a Carnot algebra we have that $V_0=\{0\} $ and so $\g \rtimes \R$ has trivial center; see Exercise~\ref{ex_grading_g_R}.ii. Consequently, the map $\ad:\g \to \mathfrak {gl} (\g \rtimes \R )$ 
 is injective. Moreover, to see that the map is valued into $\mathfrak {gl} (\g \rtimes \R )_{>0}$, take $n,m\in \N$, $X \in V_n(\g)$, and $(Y,s) \in V_m(\g \rtimes \R ).$ 
On the one hand, if $m=0$ so $Y=0,$ then $\ad_X(Y,s)=[(X,0), (0,s)]= (-snX,0) \in V_n.$ Hence, $\ad_X$ increased the degree by $n.$ 
On the other hand, if $m\ne 0$ so $s=0$, then $\ad_X(Y,s) =[(X,0),(Y,0))] =([X, Y],0) \in V_{m+n}$ and consequently $\ad_X$ increased the degree by $n.$
Thus $\ad_X \in \mathfrak {gl} (\g \rtimes \R )_{>0}$ and the proof is completed recalling Proposition~\ref{glVa_properties}.
 \end{proof}

\section{Nilpotent Lie groups}

Nilpotency for groups can be equivalently defined using either the lower central series or
 the upper central series, like for Lie algebras in Definition~\ref{def lower central series} and Exercise~\ref{upper central series algebra}, respectively.
 For the central series in groups, recall the {notation}: 
 For $a,b\in G$ we write $[a,b]:=aba^{-1}b^{-1}$ and, instead, for 
 subsets $A,B\subseteq G$ we define $[A,B]$ to be the subgroup of $G$ generated by all those elements $[a,b] $ as $a$ varies in $ A$ and $b$ varies in $ B$. 

\begin{definition}[Lower central series]\label{def lower central series groups}
Let $G$ be a group. We iteratively define the elements of the {\em lower central series} $(C^i(G))_{i\in\N}$ of $G$, also called the {\em descending central series} of $G$, by\index{descending central series! -- for a group}\index{lower central series! -- for a group}
$$C^1(G):=G, \quad
 \text{ and } \quad C^{i+1}(G):= [G,C^i(G) ],\qquad \forall i\in \N .$$
\end{definition}

\begin{definition}[Upper central series]
Let $G$ be a group. We iteratively define the elements of the {\em upper central series} $(\zeta_i(G))_{i\in\N}$ of $G$ by
\index{upper central series! of a group}
$$\zeta_0(G):=\{1_G\} \quad 
\text{ and } \quad \zeta_{i+1}(G):=\left\{ g\in G: [g,G]\subseteq \zeta_i(G) \right\} ,\qquad \forall i\in \N.$$
\end{definition}

We make some simple observation about the elements $C^i(G)$ and $\zeta_i(G) $:
\begin{enumerate}
\item $\zeta_1(G)= Z(G)$ is the center of the group.\index{center! -- of a group} 
\item $C^1(G)=[G,G]$ is the commutator subgroup.\index{commutator! -- subgroup} 
 \item $\{1\}=\zeta_0(G)<\zeta_1(G)<\ldots<\zeta_{i-1}(G)<\zeta_i(G)<\ldots $.
 \item $G=C^1(G)>C^2(G)>\ldots>C^{i-1}(G)>C^i(G)>\ldots $.
\end{enumerate}

\begin{definition}[Nilpotent group]\label{def_nilpotent_group}\index{nilpotent! -- group}
 A group $G$ is {\em nilpotent} if there exists $d$ such that $C^{d+1}(G)=\{1\}$. The minimal $d$ is called the {\em nilpotency step} of $G$. 
 Equivalently, a group $G$ is {\em nilpotent} if there exists $d$ such that
 $\zeta_d=\{1\}$ and the minimal $d$ is called the {\em nilpotency step} of $G$. 
\end{definition}

There will be an easy way to construct nilpotent Lie groups as subgroups of ${\rm {Nil }}_n$, where the latter, and its Lie algebra ${\mathfrak{nil}}_n $, are defined in Example~\ref{example Strictly upper-triangular matrix algebra and group}.
Indeed, given any nilpotent Lie algebra $\mathfrak{n}$, by Birkhoff theorem, we can see it as a subalgebra of ${\mathfrak{nil}}_n $, for some $n\in \N$.
By Theorem~\ref{teo1145bis} there exists a unique connected Lie subgroup $N\subseteq {\rm {Nil }}_n $ (a priori, not closed) with $\Lie(N)=\mathfrak{n}$.
We shall see that, actually, every such an $N$ is closed.

%
%


The main aim of this section is to show the following results:
\begin{enumerate}
\item[a)] A connected Lie group is nilpotent if and only if its Lie algebra is nilpotent; see Section~\ref{sec:Lie_groups_nilpotent_algebras}.
\item[b)] Every nilpotent simply connected Lie group is isomorphic to a closed subgroup of ${\rm {Nil }}_n$ for some $n\in \N$; see Theorem~\ref{CG-1.2.1}.iii.
 In particular, every nilpotent simply connected Lie group is a linear group.
 We will also prove that every connected subgroup of ${\rm {Nil }}_n$ is closed and simply connected; see Proposition~\ref{subgroups_Nil}.
 From these facts, we will get plenty of consequences.
\end{enumerate}

\subsection{Examples of nilpotent Lie groups}

\begin{example}[Upper-triangular unipotent matrices]\label{example Strictly upper-triangular matrix algebra and group}\index{unipotent! upper-triangular -- matrix}\index{upper-triangular matrix}
For each $n\in \N$, we consider the matrix group 	\[
	 {\rm {Nil }}_n := 
	\left\{
	\begin{bmatrix}
	 	1 & * & * \\
		 0 & \ddots & * \\
		 0 & 0 & 1
	\end{bmatrix}
	\right\}
	\subset\GL(n, \R)
	\]
	formed by all the matrices that have 1s along the diagonal and zero entries below the diagonal.
	The Lie group $ {\rm {Nil }}_n $
is nilpotent of step $(n-1)$ and its Lie algebra is ${\mathfrak{nil}}_n$ as defined in Example~\ref{example Strictly upper-triangular matrix algebra}.
	 \end{example}
	 
	\begin{example} 
Every closed subgroup of ${\rm {Nil }}_n$ is a nilpotent Lie group. We shall see that every nilpotent simply connected Lie group is of this type; see Proposition~\ref{prop:exists_group_algebra_nilpotent1}.
	 \end{example}

\begin{example}[$ {\rm {Nil }} (\mathcal F)$]\label{ex Nil F}\index{$ {\rm {Nil }} (\mathcal F)$} Here is a slight generalization of the previous example.
Let 	$V $ be a finite-dimensional vector space and let $\mathcal F= (V_0, \dots, V_n)$ be a flag for $V$.
The Lie group $$ {\rm {Nil }} (\mathcal F) :=\big\{A\in \gl(V) : (A-\mathbb I)(V_k)\subseteq V_{k-1}, \forall k\in\{1, \ldots, n\}\big\}$$ 
 is nilpotent of step $(n-1)$ with Lie algebra $ \g_{\rm nil} (\mathcal F)$ as in Example~\ref{ex g nil F}.
 \end{example}

\begin{example}[Heisenberg groups]\label{example Heisenberg algebras}\index{Heisenberg! -- group} 
The {\em $(2n+1)$-dimensional Heisenberg Lie algebra} is the Lie algebra with basis $\{X_1, \ldots, X_n, Y_1, \ldots, Y_n, Z\}$, whose pairwise brackets are equal to zero except for 
$$[X_j,Y_j]=Z, \quad\text{ for } j\in\{1, \ldots, n\}.$$
It is a two-step nilpotent Lie algebra. One
way to realize it as a matrix algebra is to consider $(n+2)\times(n+2)$ upper-triangular matrices of the form
$$\begin{bmatrix}
 0 	& x_1	&\ldots	& x_n	&z\\
 \cdot 	& 0 &\cdot	&0	&y_1\\
\cdot	&	&\cdot	&\cdot	&\vdots\\
\cdot	&	&	&0	&y_n\\
0 	&	\cdot& \cdot	& \cdot	&0
 \end{bmatrix}, \qquad \text{ for }x_1, \ldots, x_n, y_1, \ldots, y_n, z\in\R.$$
The simply connected Lie group associated with this Lie algebra is called the {\em $n$-th Heisenberg group}, and as a matrix group, it is 
$$G=\left\{ \begin{bmatrix}
 1 	& x_1	&\ldots	& x_n	&z\\
 \cdot 	& 1 &	\cdot&	0&y_1\\
\cdot	&	&\ddots	&	\cdot&\vdots\\
\cdot	&	&	&1	&y_n\\
0 	&	\cdot& \cdot	& \cdot	&1
 \end{bmatrix} : x_1, \ldots, x_n, y_1, \ldots, y_n, z\in\R\right\}\subset \GL(n+2, \R).$$
Every Heisenberg group is nilpotent of step 2.
 \end{example}
 
 \begin{example}\label{ex_group_step2_20jun} If $\g$ is a step-2 Lie algebra, then $$ x\cdot y = x+y+\frac 1 2 [X, Y],$$ defines a group structure on $\g.$ Such a Lie group is nilpotent of step 2. Hence, every nilpotent simply connected Lie group is isomorphic to a group $G_q$ as in Definition~\ref{def_step_2_G}.
\end{example}
 
\subsection{The exponential function on nilpotent matrices}\index{exponential! -- map on nilpotent matrices} 

 Nilpotent simply connected Lie groups have the feature that their exponential maps are global diffeomorphisms. Recall from Section~\ref{sec Matrix exponential} that for every finite-dimensional vector space $V$ the exponential map $\exp : \gl (V) \to\GL(V) $ is defined as $$A \mapsto e^A := \sum_{ n=0}^{\infty} \frac{A^n}{n!}.$$ 
 We shall first present a local inverse around the identity transformation $\mathbb{I}=e^0$ in $\GL(V)$.
\begin{definition}[Logarithm function]\index{logarithm}
Fix $n\in \N.$ Using the operator norm on $\GL(\R^n)$, we define the set
\begin{equation*}
B_\mathbb{I} (1):= \{ M \in \GL(\R^n)\,:\, \| M-\mathbb{I}\| < 1 \} \subseteq \GL(\R^n)
\end{equation*}
and the map $\log :B_\mathbb{I} (1) \subseteq \GL(\R^n)\to \gl(\R^n)$ as $$M \mapsto \log (M) := \sum_{ k=1}^{\infty} \frac{(-1)^k}{k} (M-\mathbb{I})^k.$$
\end{definition}

\begin{remark}
Since we have uniform convergence for all $x\in (-1,1)$ of $\log (1+x)= \sum _{ k=1}^{\infty} \frac{(-1)^k}{k} x^k,$ the map $\log$ is smooth on the open set $B_\mathbb{I} (1).$ Moreover, as a consequence of a formal series inversion, we have that the map $\log$ is the inverse of $\exp $ in small enough neighborhoods of $\mathbb{I}$ and $0$, respectively; see Exercise~\ref{ex log inverse of exp}.
\end{remark}

In the next proposition, we shall consider the set $\mathcal{U}$ of unipotent matrices and the set $\mathcal{N}$ of nilpotent matrices.
 The set $\mathcal{U}$ is not a subgroup, nor $\mathcal{N}$ is a subalgebra; see Exercise~\ref{ex span of N}. 
 We shall then study how to equivalently express $\exp$ and $\log$ on these sets.

\begin{proposition}\label{exp log on U N}
Fix $n\in \N$ and consider the following sets:\index{$\mathcal{U}$}\index{$\mathcal{N}$}
\begin{equation}\label{def U N}
\begin{aligned}
\mathcal{U} & := \{ M \in\GL(\R^n)\,:\, (M-\mathbb{I})^n =0 \}, \\
\mathcal{N}& := \{ A \in \gl(\R^n) \,:\, A^n =0 \}.\\
\end{aligned}
\end{equation}
The map $\exp_{|\mathcal{N}} :\mathcal{N} \to \mathcal{U}$ is equivalently defined as $$A \mapsto e^A\stackrel{\rm def}{=} \sum_{ k=0}^{\infty} \frac{A^k}{k!} =\sum_{ k=0}^{{ {\color{black} n} }} \frac{A^k}{k!},$$ and is an homeomorphism with inverse $\log _{|\mathcal{U}} :\mathcal{U} \to \mathcal{N}$, which is equivalently defined as $$M \mapsto \log M\stackrel{\rm def}{=} \sum_{ k=1}^{\infty} \frac{(-1)^{k+1}}{k} (M-\mathbb{I})^k = \sum_{ k=1}^{ {\color{black} n}} \frac{(-1)^{k+1}}{k} (M-\mathbb{I})^k.$$
\end{proposition}
We stress that in the second sum above, the sum over $k$ stops at $n$.
\begin{proof}
Observe that when restricting to $\mathcal{U}$, respectively to $\mathcal{N}$, the finite sums give the same maps.
Moreover, since the second sums are finite, they define maps that are polynomial and globally defined on $\gl(\R^n)$.

We check that $\exp_{|\mathcal{N}}$ is $\mathcal{U}$-valued: Take $A\in \mathcal{N},$ so $A^n=0.$ Observe that $A$ and $\sum_{ k=0}^{{ {\color{black} n} }} \frac{A^{k-1}}{k!}$ commute. Thus
$$(e^A-I)^n=\left(\sum_{ k=1}^{{ {\color{black} n} }} \frac{A^k}{k!}\right)^n =A^n\left(\sum_{ k=1}^{{ {\color{black} n} }} \frac{A^{k-1}}{k!}\right)^n=0.$$

In a similar way we check that $\log _{|\mathcal{U}}$ is $\mathcal{N}$-valued: Take $M \in \mathcal{U}$ so $(M-\mathbb{I})^n=0$ and so, because $(M-\mathbb{I})$ commutes with its powers, we get $(\log M)^n = (M-\mathbb{I})^n ( \sum_{ k=1}^{ {\color{black} n}} \frac{(-1)^{k+1}}{k} (M-\mathbb{I})^{k-1})^n =0.$ 

Finally, we check that $\exp \circ \log $ is the identity map. Fix $A \in \mathcal{N}.$ Consider $t\in \R \mapsto \log _{|\mathcal{U}} (\exp _{|\mathcal{N}} (tA)).$ Notice that $tA \in \mathcal{N}$ for every $t\in \R$ and, for $t$ small, the element $tA$ is in a neighborhood of $0$ where $\exp$ has $\log$ as inverse; see Exercise~\ref{ex log inverse of exp}. Hence, for $t$ small we get
\begin{equation*}
\log (\exp (tA)) =tA.
\end{equation*}
 Since both functions are polynomial in $t,$ then they coincides for every $t$ and in particular for $t=1.$
 
Similarly, $\exp_{|\mathcal{N}} (\log_{|\mathcal{U}} (\mathbb{I}+tA) )$ and $\mathbb{I}+tA$ are polynomials in $t$ that coincide for $t$ small, hence they do at $t=1.$
\end{proof}

 
 We now consider   ${\mathfrak{nil}}_n $ and ${\rm {Nil }}_n $, from Examples~\ref{example Strictly upper-triangular matrix algebra} and~\ref{example Strictly upper-triangular matrix algebra and group}, respectively.
We see  ${\mathfrak{nil}}_n $  as a subset of $ \mathcal{N}$ and ${\rm {Nil }}_n $ as a subset of $ \mathcal{U}$.
 The exponential map on ${\mathfrak{nil}}_n $ is given by the finite sum shown in Proposition~\ref{exp log on U N}.
\begin{corollary}\label{exp_Nil_diffeo}
The map
$\exp :{\mathfrak{nil}}_n \to {\rm {Nil }}_n$ is a polynomial diffeomorphism with polynomial inverse $\log_{|{\rm {Nil }}_n} : {\rm {Nil }}_n \to {\mathfrak{nil}}_n.$ 
\end{corollary}

\section{Connected nilpotent Lie groups}

\subsection{Connected Lie subgroups of ${\rm {Nil }}_n$} 

We shall prove that every nilpotent simply connected Lie group is, in fact, a subgroup of some ${\rm {Nil }}_n$; see Proposition~\ref{prop:exists_group_algebra_nilpotent1}. Before that, we study connected subgroups of ${\rm {Nil }}_n$.
\begin{proposition}\label{subgroups_Nil}\index{${\rm {Nil }}_n$}
If $G$ is a connected Lie subgroup of ${\rm {Nil }}_n$, then $\exp (\Lie(G))=G$ and in particular $G$ is closed and simply connected.
\end{proposition}

\begin{proof}
Recall from Corollary~\ref{exp_Nil_diffeo} that on ${\rm {Nil }}_n $ and on its Lie algebra ${\mathfrak{nil}}_n $, the respective maps $\exp_{|\mathcal{N}}$ and $\log _{|\mathcal{U}}$ are polynomial and are inverse of each other.
Let $\g := \Lie (G) \subset {\mathfrak{nil}}_n.$ Clearly, $\exp_{|\mathcal{N}} (\g)= \exp (\g)\subset G.$ We prove the other inclusion. Let $U$ be a connected neighborhood of $0$ in $\g$ such that $\exp _{| {U}}: U \to \exp (U)$ is a diffeomorphism and $V:= \exp (U)$ is a neighborhood of $1_G$ in $G$. 
From a general argument valid for topological groups, since $G$ is connected, we have that $G = \bigcup _{m=1}^\infty V^m$; see Exercise~\ref{Prop:generating}. We shall prove by induction on $m\in \N$ that \begin{equation}\label{claim712024}V^m \subset \exp (\g ).\end{equation} The base of induction is that $V^1=V\stackrel{\rm def}{=}\exp (U) \subset \exp (\g).$
Now assume that the claim \eqref{claim712024} holds for $m\in \N. $ Take $w\in V^{m+1} =V^m\cdot V.$ Hence, there are $x\in \g $ and $ y \in U$ such that $w=\exp (x)\exp (y)$ and so using BCH formula (see Proposition~\ref{Proposition BCH}) we have that $w= \exp (x+y + \frac 1 2 [X, Y] +\dots ) $, where the sum is actually finite and it is by terms in $\g$, then $w\in \exp (\g).$ Consequently, we infer that $\exp (\g)=G.$ The last assertion of the proposition is a consequence of the fact that $\Lie(G)$ is a vector subspace (hence closed and simply connected) and $\exp$ is a diffeomorphism on ${\mathfrak{nil}}_n$; see Corollary~\ref{exp_Nil_diffeo}. 
\end{proof}

\subsection{Lie groups with nilpotent Lie algebras}\label{sec:Lie_groups_nilpotent_algebras}
In this section, we explain why a connected Lie group is nilpotent if and only if its Lie algebra is nilpotent.
What the reader should expect is that the lower central series are linked by $ C^m(G) = \exp( C^m(\g)), $ for $m\in \N$, at least for (nilpotent) simply connected Lie groups.
\subsubsection{Commutator subgroups}\index{commutator! -- subgroup} 

We point out a general fact about Lie groups:
 For every connected Lie group with Lie algebra $\g$, the commutator subgroup $[G,G]$ is a Lie subgroup whose Lie algebra is $[\g, \g].$ 
The hardest part of the mentioned fact is explaining why $[G, G]$ is a Lie subgroup. 
A warning is that $[G, G]$ may not be topologically closed. It is if $G$ is simply connected.
The reader can read about these facts in \cite[page 138]{Hochschild}. 
We will not further discuss this topic in this general context because when $G$ is nilpotent and simply connected, we will easily have that $[G, G]$ is closed; see Exercise~\ref{ex connected nilpotent gives closed and simply connected}.
\subsubsection{Lie algebras of nilpotent Lie groups}
It is easy to show that nilpotent Lie groups have nilpotent Lie algebras. Doing the increasingly difficult Exercises \ref{ex845} and \ref{ex846}, the reader should be able to come up with an argument for the next proposition. 

\begin{proposition}\label{nilp_gp_nilp_alg}
If $G$ is a nilpotent Lie group of step $s$, then $\Lie(G)$ is a nilpotent Lie algebra of step $s$.
\end{proposition}

\subsubsection{Lie groups with nilpotent Lie algebras}

Clearly, unless the group is connected, we cannot expect that some information on the Lie algebra can give a global information on the group. For example,
every finite group is a $0$-dimensional Lie group whose Lie algebra is nilpotent (and abelian, trivially). Hence, there are Lie groups that are not nilpotent but have a nilpotent Lie algebra.

We shall prove that if $G$ is a connected Lie group and its Lie algebra $\g$ is nilpotent, then $G$ is nilpotent. There are various ways to show the latter fact. A common point is to construct a nilpotent Lie group $N$ with Lie algebra isomorphic to $\g.$

One way is to use Birkhoff Embedding Theorem, Theorem~\ref{thm Birkhoff general}:

\begin{proposition}[after Birkhoff's theorem]\label{prop:exists_group_algebra_nilpotent1}
For every nilpotent Lie algebra $\g$, there exists $n\in \N$ and a nilpotent closed simply connected Lie subgroup $N\subseteq {\rm {Nil }}_n$ with Lie algebra isomorphic to $\g$.
\end{proposition}

Another way is to use the Baker-Campbell-Hausdorff-(Dynkin) formula of Proposition~\ref{Proposition BCH}:
\begin{proposition}[after BCH Formula]\label{prop:exists_group_algebra_nilpotent2}\index{Dynkin product}
For every nilpotent Lie algebra $\g$,
if we equip $\g$ with the Dynkin product $\star$, then $(\g, \star)$ is a nilpotent simply connected Lie group with Lie algebra isomorphic to $\g$.
 Moreover, the exponential map $\exp: \g\to(\g, \star)$ is the identity map.
\end{proposition}
We shall prove the above two propositions after the next corollary.

\begin{corollary}[Consequence of either of the propositions]\label{cor_groups_w_nilpotent_algebras}
Every connected Lie group with nilpotent
 Lie algebra is nilpotent.
\end{corollary}
\begin{proof}[Proof of Corollary~\ref{cor_groups_w_nilpotent_algebras}]
Let $G$ be a connected Lie group with nilpotent Lie algebra $\g$. 
Either from Proposition~\ref{prop:exists_group_algebra_nilpotent1} or from Proposition~\ref{prop:exists_group_algebra_nilpotent2}, 
let $N$ be a nilpotent simply connected Lie group with Lie algebra isomorphic to $\g$.
Thus, the Lie group $N$ and the universal cover $\tilde G$ of $G$ are simply connected with isomorphic Lie algebras. Thus, by Corollary~\ref{simplyconn-isom}, they are isomorphic. In particular, since $N$ is nilpotent, then so is $\tilde G$. Consequently, being the quotient of a nilpotent group, also the group $G$ is nilpotent, as desired.
 \end{proof}

\begin{proof}[Proof of Proposition~\ref{prop:exists_group_algebra_nilpotent1}]
By Birkhoff Theorem~\ref{thm Birkhoff general}, there are $n\in \N$ and a subalgebra $\mathfrak n \subseteq {\mathfrak{nil}}_n$ that is isomorphic to $\g$. By Theorem~\ref{teo1145}, let $N \subseteq {\rm {Nil }}_n$ be the connected Lie subgroup with Lie algebra $ \mathfrak n $. Then, by Proposition~\ref{subgroups_Nil}, the Lie group $N$ is closed and simply connected. Moreover, since ${\rm {Nil }}_n$ is nilpotent then so is $N$. 
\end{proof}

\begin{proof}[Proof of Proposition~\ref{prop:exists_group_algebra_nilpotent2}]
Being $\g$ a nilpotent Lie algebra, the Dynkin series is finite and, for every $x, y \in \g$, the product $x\star y = \log _{|U} (\exp_{|N} (x)\exp_{|N} (y))$ is polynomial in $x$ and $y$. 
Obviously, the manifold $(\g, \star)$ is simply connected because it is isomorphic to a vector space.
We shall check that $(\g, \star)$ is a Lie group, it is nilpotent, its Lie algebra is $\g$, and its exponential map is the identity. 
By Ado's theorem (see Theorem~\ref{Lie's Third Theorem}), there is\footnote{However, notice that for Corollay~\ref{cor_groups_w_nilpotent_algebras} we already know the existence of $G$.}
a group $G$ with $\g$ as Lie algebra. Moreover, by the BCH Formula, in a small enough neighborhood of $1_G$, the group product of $G$ is exactly the Dynkin product, up to the identification via the exponential map. In particular, since $G$ is a group, then the product $\star$ satisfies the axioms of group products near $1_{(\g, \star )}=0$ and also $x\star (-x) = 1_{(\g, \star )}$ is true near $0$.
These identities are polynomial maps, which are verified in a neighborhood of $0.$ By analytic continuation, they hold everywhere on $\g$.

By the BCH formula, the exponential map $\exp: \g \to G$ is a Lie group isomorphism between the Lie group $(\g,\star)$ and the Lie group $G$. Thus, the Lie algebra of $(\g,\star)$ is isomorphic to $\g$.

We next explain why $(\g, \star)$ is nilpotent. We stress that we do not know (yet) that $G$ is nilpotent. The reader might like to know that there is a result by Lazard saying that every polynomial group product in $\R^n$ gives a nilpotent group; see \cite{Lazard} or \cite{DeKimpe_2003}.\index{Lazard}
 Still, for the Dynkin product, it is clearer since we have that
 $$x\star y-(x+ y) \in [x, \g]+[y, \g]\subseteq [\g, \g], \qquad \forall x,y\in \g.$$
 Hence, we get $$ x\star y\star x^{-1}\star y^{-1}\in [\g, \g], \qquad \forall x,y\in \g.$$
 Since $\g$ is a nilpotent Lie algebra, iterating this group product, one obtains nilpotency of the group $(\g, \star )$.

The one-parameter subgroups for the group product $\star$ are the curves $t\in \R \mapsto tx$ for each 
 $x\in \g $, because $(sx)\star (tx)= (s+t)x.$ So the identity $\g\to\g$ is the exponential map $\exp: \g \simeq \Lie(\g, \star) \to(\g, \star)$.
\end{proof}

In the next section, we shall see that another consequence of Proposition~\ref{prop:exists_group_algebra_nilpotent2} is that 
every nilpotent simply connected Lie group is isomorphic to a closed subgroup of ${\rm {Nil }}_n$ for some $n\in \N.$ 

\subsection{Simply connected nilpotent Lie groups}

Simply connected Lie groups are uniquely determined by their Lie algebras.
Indeed, recall from Corollary~\ref{simplyconn-isom}
that if two simply connected Lie groups have isomorphic Lie algebras, then they are isomorphic.
For nilpotent groups, either Birkhoff's embedding theorem or the exponential map and the BCH formula provide a concrete identification.
We will see how one can completely work on the Lie algebra using such coordinates.
 
\begin{theorem}
\label{CG-1.2.1}
Every nilpotent simply connected Lie group $G$ has the following properties: 
\begin{description}
 \item[\ref{CG-1.2.1}.i.] The exponential map $\exp:\Lie (G)\to G$ is a diffeomorphism.
 \item[\ref{CG-1.2.1}.ii.] The Baker-Campbell-Hausdorff Formula \eqref{Dynkin Formula} holds globally. 
 \item[\ref{CG-1.2.1}.iii.] There exists $n\in \N$ and a closed nilpotent simply connected Lie subgroup of $ {\rm {Nil }}_n$ isomorphic to $G$.
 \end{description}
\end{theorem}

\begin{proof}
We begin by observing that \ref{CG-1.2.1}.iii is a direct consequence of Theorem~\ref{thm Birkhoff general}, 
 Theorem~\ref{teo1145}, and
 Proposition~\ref{subgroups_Nil}, 
 as we saw in the proof of 
Proposition~\ref{prop:exists_group_algebra_nilpotent1}.

In addition, on $ {\rm {Nil }}_n$ we also have that the two functions 
$A\star B $ and $ \log (e^Ae^B)$ are analytic maps that coincide in a neighborhood of 0, recalling Proposition~\ref{Proposition BCH} and Corollary~\ref{exp_Nil_diffeo}, and also Exercise~\ref{ex Dynkin polynomial}. Then \ref{CG-1.2.1}.ii holds for each $ {\rm {Nil }}_n$.

Once we know that the statements are valid for subgroups of $ {\rm {Nil }}_n$, they hold for arbitrary nilpotent simply connected Lie groups by Corollary~\ref{exp_Nil_diffeo}. 
\end{proof}

\index{Lie group! stratified --}
\index{Lie group! stratifiable --}
\index{Lie group! positively graded --}
\index{Lie group! graded --}
\index{stratified! -- Lie algebra}
\index{positively gradable! -- Lie algebra}
\index{stratifiable! -- Lie algebra}
\index{graded! -- Lie algebra}
\index{positively graded! -- Lie algebra}
\index{stratified! -- Lie group}
\index{positively gradable! -- Lie group}
\index{stratifiable! -- Lie group}
\index{graded! -- Lie group}
\index{positively graded! -- Lie group}
A Lie algebra is called {\em positively gradable} if it admits a positive grading, and {\em stratifiable} if it admits a stratification. A Lie algebra is referred to as $\R$-{\em graded} if it is equipped with an $\R$-grading, {\em positively graded} if it is equipped with a positive grading, and {\em stratified} if it is equipped with a stratification. Similarly, a Lie group is called {\em positively gradable}, {\em stratifiable}, $\R$-{\em graded}, {\em positively graded}, or {\em stratified} if it is simply connected and its Lie algebra is positively gradable, stratifiable, $\R$-graded, positively graded, or stratified, respectively.

If a Lie group $G$ is $\R$-graded, then its Lie algebra $\g$ has the dilations $\delta_\lambda:\g\to \g $ as in \eqref{def_dilation_relative_to_grading}, and, because $G$ is simply connected, by Theorem~\ref{thm_induced_homo0} we have that 
each Lie algebra endomorphism $\delta_\lambda$ induces a 
 Lie group endomorphism of $G$ whose induced Lie algebra homomorphism is $\delta_\lambda$. Such a 
 Lie group endomorphism is still 
 denoted by $\delta_\lambda:G\to G $, so $(\delta_\lambda)_*=\delta_\lambda$.\label{discuss_dilations20jun}
Moreover, from
Theorem~\ref{Warner3.32}
 	we have 
	\begin{equation}\label{delta_exp=exp_delta0}
	\delta_\lambda\circ\exp = \exp\circ\, \delta_\lambda .
	\end{equation}
When $G$, in addition, is simply connected and nilpotent, the map $\exp: \g \to G$ is a diffeomorphism by Theorem~\ref{CG-1.2.1}.i. Therefore, every element $g\in G$ can be represented as $\exp(X)$ for some unique $X\in\g$, and therefore uniquely written in the form
\begin{equation*}
g=\exp\left(\sum_{a\in \R} v_a\right), \qquad v_a\in V_a, \, \,a\in \R.
\end{equation*}
This representation allows to have the formula:
\begin{equation*}
\delta_\lambda\left(\exp\left(\sum_{a\in \R} v_a\right)\right)=
\exp\left(\sum_{a\in \R} \lambda^a v_a\right )
.	\end{equation*}

\subsection{Canonical coordinates}\label{sec Canonical coordinates}

One important application of Theorem~\ref{CG-1.2.1} involves the existence of good coordinates on nilpotent simply connected Lie groups. For these groups, since the exponential map $\exp:\g \to G$ is a diffeomorphism, we can use it to transfer coordinates from $\g$ to the nilpotent simply connected Lie group $G$. Some authors use $\exp$ to identify $\g$ with $G$. Then, the group multiplication can be calculated using the Baker-Campbell-Hausdorff formula.

We shall use different types of coordinate systems: exponential coordinates and Malcev coordinates. They will make us read the group operations polynomially, and the Lebesgue measure will be invariant both by left translations and by right translations. The following theorem summarises the content of this subsection.

\begin{theorem}\label{Malcev global coordinates}
In every nilpotent simply connected Lie group, exponential coordinates and Malcev coordinates are global coordinates. 
The group product in such coordinates is polynomial. Left translations and right translations have Jacobian one. Consequently, the Lebesgue measure
in these coordinates is a bi-invariant Haar measure. 
\end{theorem}
We call {\em bi-invariant Haar measure} every left-Haar measure that is also a right-Haar measure,
in the sense of Section~\ref{sec_Haar_poly_growth}.
For the proof of the theorem, we refer to the end of this subsection.\index{bi-invariant! -- Haar measure}

\subsubsection{Exponential coordinates}
\begin{definition}[Exponential coordinates: canonical coordinates of the $1^{st}$ kind]
Let $( X_1, \ldots, X_n )$ be an ordered basis for a nilpotent Lie algebra of a nilpotent simply connected group $G$. The coordinates given by the map
 $$\Phi:\R^n\longrightarrow G$$
$$\Phi(t_1, \ldots,t_n):=\exp(t_1X_1+\ldots+t_n X_n)$$
are called 
{\em exponential coordinates} with respect to $X_1, \ldots, X_n$. 
The map $\Phi$ is called {\em exponential coordinate system}.
Exponential coordinates are also known as {\em canonical coordinates of the first kind}.\index{exponential! -- coordinates}
\end{definition}
With the choice of the basis and the use of $\exp$, we are identifying $\R^n$ with $\Lie(G)$ and $G$. Moreover, the group product can be obtained through the BCH formula from Section~\ref{section BCH}:
\[
(s_1, \dots,s_n)\star(t_1, \dots,t_n) 
= 
\log\left(\exp\left(\sum_{j=1}^n s_jX_j\right)\exp\left(\sum_{j=1}^n t_jX_j\right)\right).
\]

\begin{example}
\index{Heisenberg! -- group} 
We consider the Heisenberg group, whose Lie algebra has a basis $X_1, X_2, X_3 $ with only nontrivial relation $[X_1, X_2]=X_3$.
The exponential coordinate system with respect to this basis is the map
\begin{equation}\label{Heisenberg exponential coordinates}
(x_1,x_2,x_3) \stackrel{\Phi}{\longmapsto} \exp(x_1X_1+x_2X_2+x_3X_3) .
\end{equation}
In these coordinates, the product can be expressed by the BCH formula:
\begin{multline*}
\exp(x_1X_1+x_2X_2+x_3X_3)\exp(x_1'X_1+x_2'X_2+x_3'X_3) \\
=\exp\Big( x_1X_1+x_2X_2+x_3X_3 + x_1'X_1+x_2'X_2+x_3'X_3 \hspace{2cm} \\
		\hspace{2cm} + \frac12 \big[x_1X_1+x_2X_2+x_3X_3, \, x_1'X_1+x_2'X_2+x_3'X_3 \big] \Big) \\
	= \exp\left((x_1+x_1')X_1+(x_2+x_2')X_2+(x_3+x_3')X_3 + \frac12(x_1x_2' - x_2x_1') X_3 \right),
\end{multline*}
where we used the Lie bracket relations of the Heisenberg Lie algebra.
We conclude that the product, when read in these coordinates, is
\[
(x_1,x_2,x_3)\star(x_1',x_2',x_3') = \left( x_1+x_1', x_2+x_2', x_3+x_3'+\frac12(x_1x_2' - x_2x_1') \right).
\]
This map has a Jacobian determinant equal to $1$.
\end{example}
 
\subsubsection{Malcev bases}
It is useful to choose bases that are better adapted to the Lie brackets. We shall then have other canonical coordinates in Definition~\ref{def canonical coordinates 2}.
\begin{definition}[Malcev basis for a Lie algebra]\label{def Malcev basis}\index{Malcev! -- basis}
Let $\mathfrak g$ be a Lie algebra. An ordered basis $(X_1, \dots, X_n)$ for $\mathfrak{g}$ is called a {\em strong Malcev basis}, or, simply, a {\em Malcev basis} if for every $k\in\{1, \dots,n\}$ the space
\[
\mathfrak{g}_k:=\mathrm{span}\{X_1, \dots, X_k\}
\]
is an ideal of $\mathfrak g$, i.e.,
	\(
	[\g, \g_k]\subset\g_k 
	\).
\end{definition}
 We remark that there exist non-nilpotent Lie algebras (e.g., in 2D) that admit Malcev bases.
 We shall not discuss the notion of weak Malcev basis, for which we refer to \cite[page 10]{Corwin-Greenleaf}. 
 
The spaces $\mathfrak{g}_k$ coming from a Malcev basis satisfy the following stronger property.
\begin{proposition}
If $\mathfrak g$ is a nilpotent Lie algebra and $(X_1, \dots, X_n)$ is a Malcev basis for $\g$, then \begin{equation}\label{bracket:gk}
[\mathfrak g, \mathfrak g_k]\subseteq \mathfrak g_{k-1}\subseteq \mathfrak{g}_k, \qquad \text{ for } k\in\{1, \dots,n\},
\end{equation}
where $\mathfrak{g}_k:=\mathrm{span}\{X_1, \dots, X_k\}
$ and $\mathfrak{g}_0:= \{0\}$.
\end{proposition}

\begin{proof} Fix $k\in\{1, \dots,n\}$.
By definition of Malcev basis, we have $ [\g, \g_k]\subseteq \g_k$ and also $ [\g, \g_{k-1}]\subseteq \g_{k-1}$. If the conclusion of the proposition were not true, then there would be some $Y\in \g$ and $a_1, \ldots, a_k\in \R$ with $a_k\neq 0$ such that
$$[Y, X_k]=a_k X_k +\sum_{i=1}^{k-1} a_i X_i.$$
Now we iterate bracketing by $Y$, i.e., 
we iterate the map
$
	\ad_{Y} = [Y, \cdot]
$. Thus, 
	 we get, for some $a_1^{(l)}, \ldots, a_{k-1}^{(l)}\in \R$, the value
$$(\ad_{Y}^l)(X_k)=a_k^l X_k +\sum_{i=1}^{k-1} a_i^{(l)} X_i\in a_k^l X_k +\g_{k-1}, \qquad\forall l\in \N,$$
which is never zero and so contradicts the nilpotency of $\g$.
\end{proof}

In the special class of Carnot groups, as considered in Chapter~\ref{ch_Carnot}, the existence of Malcev bases will be a triviality; see Exercise~\ref{ex Malcev basis in Carnot}. However, every nilpotent algebra has Malcev bases. 
For additional information, we refer to \cite[Theorem~1.1.13]{Corwin-Greenleaf} and the notes therein. 

\begin{proposition}\label{Prop exists Malcev bases}
In nilpotent Lie algebras, Malcev bases exist.
\end{proposition}
\begin{proof}
For ${i\in\N}$, let $\zeta_i:=\zeta_i(\g) $ be the $i$-th element of the upper central series of $\g$, as in Exercise~\ref{upper central series algebra}.
For some $s\in \N$ we have 
$$\{ 0 \} \neq \zeta_1=Z(\g) \subsetneq \ldots \subsetneq\zeta_{s-1}\subsetneq\zeta_s=\g.$$
To build a Malcev basis for $\g$, we begin with a basis $X_1, \ldots, X_{n_1}$ of $\zeta_1$. Next, we complete it to a basis $X_1, \ldots, X_{n_2}$ of $\zeta_2$. Iterating, given a basis of $\zeta_{j-1}$, we complete it to a basis $X_1, \ldots, X_{n_j}$ of $\zeta_j$. Because of nilpotency, we obtain a basis of $\g$ after $s$ steps.
In fact, the vectors $X_1, \ldots, X_{n_s}$ form a basis of $\g$ such that for every $k\in \{1, \ldots, {n_s}\}$ there is $j_k \in \{1, \ldots, s\}$ such that $X_k \in \zeta_{j_k}\setminus\zeta_{j_k-1}$. For such a $j_k$,
 we have
$$\zeta_{j_k-1}\subseteq \g_k:=\mathrm{span}\{X_1, \dots, X_k\} \subseteq \zeta_{j_k}.$$
Hence, we deduce
	\(
	[\g, \g_k]\subset[\g, \zeta_{j_k}]\stackrel{\rm {def}}{=} \zeta_{j_k-1}\subset\g_k 
	\).
\end{proof}

 \subsubsection{Malcev coordinates}
\begin{definition}[Malcev coordinates: canonical coordinates of the $2^{nd}$ kind]\label{def canonical coordinates 2}\index{Malcev! -- coordinates}\index{exponential! -- coordinates! -- of the second kind}
Let $G$ be a Lie group with Lie algebra $\mathfrak g$. Fix an ordered basis $(X_1, \dots, X_n)$ of $\mathfrak{g}$. The coordinates given by the map 
\[
\Psi:(s_1, \dots,s_n)\in \mathbb R^n \mapsto \exp(s_1X_1)\cdots\exp(s_nX_n)\in G,
\]
are called {\em exponential coordinates of second kind}, or {\em canonical coordinates of the second kind} with respect to $(X_1, \dots, X_n)$.
 If $(X_1, \ldots, X_n)$ is a Malcev basis,
then the defined coordinate system $\Psi$ is called {\em Malcev coordinate system}.
\end{definition}

We remark that by observing the injectivity of the differential at $\bf 0$ of the map $\Psi$ as in Definition~\ref{def canonical coordinates 2}, one can see that $\Psi$ gives local coordinates around $1_G$. We will soon see that it gives global coordinates when the basis is a Malcev basis.

\begin{example}[Malcev and Non-Malcev coordinates] In the Heisenberg group, we have the exponential coordinates $\Phi$, as in 
\eqref{Heisenberg exponential coordinates}, with respect to the standard basis $(X_1, X_2, X_3)$. 
Such a basis, in this order, is not a Malcev basis: we stress that $X_3$ is in the derived subalgebra, while $X_1$ is not.
Instead, the ordered basis $(X_3, X_2, X_1)$ is a Malcev basis.
We may not have good properties if we consider canonical coordinate systems of the second kind with respect to bases that are not Malchev bases. For example, the triple $(X_1, X_2, X_1+ X_3)$ is not a Malcev basis, and if $\Psi$ denotes the canonical coordinate system of the second kind with respect to this basis, then $\Phi^{-1}\circ\Psi$ is the map
$$(s_1,s_2,s_3)\in \mathbb R^3 \mapsto \left(s_1 + s_3, \,s_2, \, 
s_3+\frac{s_1 s_2- s_2 s_3}{2} \right)\in \R^3.$$
This map has a Jacobian determinant equal to $1-s_2$.
We should avoid canonical coordinates of the second kind with respect to non-Malcev bases. 
\end{example} 

 If $(X_1, \ldots, X_n)$ is a Malcev basis for the Lie algebra of a nilpotent Lie group, we can consider both canonical coordinate systems; we have that the Malcev coordinates are related to the exponential coordinates by a polynomial diffeomorphism. Such a polynomial map has an upper triangular form, and, in particular, its Jacobian determinant is constantly equal to $1$. The proof is an adaptation of \cite[Proposition 1.2.7]{Corwin-Greenleaf}.
\begin{proposition}[Change of Malcev coordinates]\label{CG1.2.7}
Let $G$ be a nilpotent simply connected Lie group with Lie algebra $\g$.
 Let $(X_1, \ldots, X_n)$ be a Malcev basis for $\g$. Let $\Psi:\R^n\to G$ be the Malcev coordinate system and $\Phi:\R^n\longrightarrow G$ be the exponential coordinate system associated with the basis.
Then
\begin{description}
 \item[\ref{CG1.2.7}.i.] 
 the map $P:= \Phi^{-1}\circ \Psi$ from $\R^n $ to $ \R^n$ is a polynomial diffeomorphism with polynomial inverse;
 \item[\ref{CG1.2.7}.ii.] writing the components $P=(P_1, \ldots, P_n)$, then $P_j(s)=s_j+\hat P_j(s_{j+1}, \ldots,s_n),$ where $\hat P_j$ is a polynomial depending only on the last $n- j$ variables, for $j\in\{1, \ldots, n\}$.
 \end{description}
\end{proposition}

In other words, each Malcev coordinate system is a global diffeomorphism and, for some polynomials $P_1, \ldots, P_n$, we have the relation:
$$ \exp(s_1X_1)\cdots\exp(s_nX_n) = \exp(P_1(s)X_1+\ldots+P_n(s) X_n).$$

\begin{proof}
Recall that we are considering
\[
\Psi(s_1, \ldots,s_n):=\exp(s_1X_1)\cdots \exp(s_nX_n), \qquad \Phi(s_1, \ldots,s_n):=\exp\left(\sum_{i=1}^n s_iX_i\right).
\]
Because of BCH Formula (Proposition~\ref{Proposition BCH}), the map 
\[
P:=\Phi^{-1}\circ\Psi:\mathbb R^n\to\mathbb R^n
\]
is a polynomial map
\[
s:=(s_1, \ldots,s_n)\mapsto (P_1(s), \ldots,P_n(s)).
\]
By induction on $n$, we shall prove 
\[
P_j(s)=s_j+Q_j(s_{j+1}, \ldots,s_n), \qquad\forall j\in\{1, \ldots,n\}.
\]
When $n=1$, it is trivial. For arbitrary $n$, we take the quotient 
$\pi:\mathfrak g\to \mathfrak g/\mathfrak g_1$
modulo $\mathfrak{g}_1:=\mathbb RX_1\triangleleft \mathfrak g$. 
Then, for every $j=2, \dots,n$, we set $\overline{X}_j:=\pi(X_j)$, and we notice that $(\overline{X}_2, \ldots, \overline{X}_n)$ is a Malcev basis for $\mathfrak g/\mathfrak g_1$. Now, in $G$ we have the equation
\[
\exp(s_1X_1)\cdots \exp(s_nX_n)= \exp\left(\sum_{i=1}^n P_i(s) X_i\right),
\]
which, after passing to the quotient $G/\exp(\R X_1)$, becomes
\[
\exp(s_2 \overline{X}_2)\cdots \exp(s_n\overline{X}_n) = \exp\left(\sum_{i=2}^n P_i(s) \overline{X}_i\right).
\]
Thus, by induction, we have 
\[
P_j(s)=s_j+Q_j(s_{j+1}, \ldots,s_n), \qquad\forall j=2, \ldots,n.
\]
Lifting to $G$ we get
\[
\exp(s_2X_2)\cdots \exp(s_nX_n) =\exp\left(Q_1(s_2, \ldots,s_n)X_1+\sum_{i=2}^n P_i(s)X_i\right),
\] for some polynomial $Q_1$.
Notice that $X_1$ is central, hence $\exp(X_1+Y)=\exp(X_1)\exp(Y)$ for every $Y\in \mathfrak g$. Thus, from the previous equation, we get
\[
\exp(s_1X_1)\cdots \exp(s_nX_n)=\exp\left((s_1+Q_1(s_2, \ldots,s_n))X_1+\sum_{i=2}^n P_i(s)X_i\right).
\]
Now, in order to finish, notice that given $t\in \R^n$, the equations 
\[
t_j=s_j+P_j(s_{j+1}, \ldots,s_n), \qquad\forall j\in\{1, \ldots,n\},
\]
can be solved in $s$ with a formula 
\[
s_k=t_k+\widetilde Q_k(t_{k+1}, \ldots, t_n),
\]
where $\widetilde Q_k$ are polynomials.
\end{proof}

\begin{remark} If $(X_1, \dots, X_n)$ is a Malcev basis of a nilpotent simply connected Lie group, then left and right cosets of $\mathrm{exp}(\mathfrak{g}_k)$ are affine planes (set-wise) in exponential coordinates of the second kind. This may not be true in exponential coordinates of the first kind. 
 \end{remark}

\subsubsection{Jacobian determinants in exponential coordinates}
The following result describes properties of a nilpotent group law expressed in exponential coordinates with respect to some Malcev basis. The proof is an adaptation of \cite[Proposition 1.2.9]{Corwin-Greenleaf}.

\begin{proposition}[Group law in Malcev coordinates]\label{CG-1.2.9}
Let $G$ be a nilpotent simply connected Lie group equipped with a Malcev basis.
On $G$, consider exponential coordinates of the first or second kind associated with the basis.
In this coordinate system, the group law has an upper triangular form:
$$(s_1, \dots,s_n)\cdot (t_1, \dots,t_n)= s+t +\sum_{j=1}^n Q_j(s,t) e_j, \qquad \forall s,t \in \R^n,$$
where each $Q_j$ is a polynomial that is not depending on $s_1, \dots,s_{j-1}$ nor $t_1, \dots,t_{j-1}$.
In particular, for each $g\in G$, the left translation $L_g$ and the right translation $R_g$ are maps whose Jacobian determinants are identically equal to $1$.
\end{proposition}

\proof We prove the statement for exponential coordinates of the first kind and left translations.
The case of right translations is similar. For Malcev coordinates, the result will also be true because, by Proposition~\ref{CG1.2.7}, they differ from exponential coordinates by a polynomial diffeomorphism, which has an upper triangular form.

The proof is based on the BCH formula and \eqref{bracket:gk}.
Let $\Phi$ be the exponential coordinate system, and $L_g$ be the left translation by $g\in G$. We need to calculate the differential of 
$\Phi^{-1}\circ L_g\circ \Phi$. Thus, we consider the diagram
\begin{center}
 \begin{tikzcd}
 \R^n\ar[rrrrr]{u}{\Phi} &&& && G\ar[ddd]{r}{L_g} \\
(t_1, \ldots,t_n)\ar[rrrr,|->]{u}{\Phi} 	&&&&\exp(\sum_{j=1}^n t_j X_j) \ar[ddd,|->, crossing over]{r}{L_g}& \\
&& &&&\\
 \R^n\ar[rrrrr]{u}{\Phi} &&&&& G\\
(s_1, \ldots,s_n)	\ar[rrr,|->]{u}{\Phi} &\,&\,&\,&
\hspace{-1.3cm}
\exp(\sum_{j=1}^n s_j X_j )=g\exp(\sum_{j=1}^n t_j X_j)
\hspace{-1.3cm}
& 
 \end{tikzcd}
\end{center}
and we solve the dependence of the $s_i$'s from the $t_j$'s.
Since by Proposition~\ref{CG1.2.7} Malcev coordinate systems are surjective, we can find $u_1, \ldots,u_n\in \R$ and write
$$g=\exp(u_1 X_1)\cdots\exp(u_n X_n).$$
It is enough to consider the case $g=\exp(u_k X_k)$ for an arbitrary $k\in \{1, \ldots,n\}$ and then conclude considering compositions.
Thus, we need to consider the system
$$\exp(\sum_js_j X_j )=\exp(u_k X_k)\exp(\sum_jt_j X_j).$$
By the BCH Formula (Proposition~\ref{Proposition BCH}), 
$$\sum_js_j X_j = u_k X_k + \sum_jt_j X_j +\dfrac{1}{2} [ u_k X_k, \sum_jt_j X_j] +\ldots.$$
Since we have chosen a Malcev basis, we have the property \eqref{bracket:gk}. Thus a bracket as $[ X_k, X_j]$ is only a combination of $\{X_1, \ldots, X_{j-1}\}$.
In other words, 
the function $s_j$ is of the form $t_j$ plus a polynomial that does not depend on the variables $t_1, \ldots, t_j$.
Thus, the differential is of the form
$$\dd(\Phi^{-1}\circ L_g\circ \Phi)=\begin{bmatrix}
 1 	& * 	&\ldots	&\ldots	&\ldots	&\ldots & *	\\
 0 	&\ddots &\ddots&*	&\ldots	&*	&\vdots	\\
\cdot 	& \cdot	&1 &*	&\ddots	& \vdots&\vdots\\
\cdot	&\cdot	&0	&1&*	&*	& \vdots\\
 \cdot	&\cdot	&\cdot	&0	&1&\ddots	& \vdots\\
\cdot	&\cdot	 &\cdot	& \cdot	&\cdot	&\ddots& *	\\
0	&\cdot	 &\cdot	&\cdot 	&\cdot	& 0 	&1
 \end{bmatrix} .$$
Thus, the Jacobian of a left translation in exponential coordinates with respect to a Malcev basis is $1$ at every point.
\qed

\begin{definition}[Polynomial coordinates]\label{def polynomial coordinates}
Let $G$ be a Lie group and $X_1, \ldots, X_n$ a basis of its Lie algebra.
 	If $P:\R^n\to\R^n$ is a diffeomorphism such that $P$ and $P^{-1}$ have polynomial components, then 
	\[
	(s_1, \dots,s_n) \longmapsto \exp(P_1(s)X_1+\dots+P_n(s)X_n)
	\]
	is called \emph{polynomial coordinate system} for $G$.\index{polynomial! -- coordinates}
\end{definition}
 Examples of polynomial coordinate maps are, obviously, exponential and, by Proposition~\ref{CG1.2.7}, Malcev coordinate maps.


The key observation is that the Jacobian of every polynomial diffeomorphism with polynomial inverse is a polynomial that is invertible inside the polynomial ring, so it is a constant. Thus, changing coordinates by a polynomial diffeomorphism with a polynomial inverse preserves those maps that preserve the Lebesgue measure.

\begin{corollary}\label{polynomial coordinates Jacobian 1}
 	In nilpotent simply connected Lie groups, left translations and right translations have Jacobian determinant 1 in
	polynomial coordinates.
\end{corollary}
\begin{proof}
 	If $P$ is a polynomial map, then $\Jac(P)$ is a polynomial.
	If $P$ and $P^{-1}$ are polynomial diffeomorphisms, then $1= \Jac(\Id) = \Jac(P^{-1}\circ P) = (\Jac(P^{-1})\circ P )\cdot \Jac(P )$.
	
	Hence, the two maps $\Jac(P)$ and $\Jac(P^{-1})\circ P$ are polynomials whose product is constant. 
	Thus they are constant, and $\Jac(P^{-1})$ is constantly equal to $1/\Jac(P)$.
	
	If $\Phi$ is an exponential coordinate map, then 
	\begin{multline*}
	 	\Jac(P^{-1}\circ\Phi^{-1}\circ L_g\circ\Phi\circ P )_x = \\
		= \Jac(P^{-1})_{(\Phi^{-1}\circ L_g\circ\Phi\circ P)(x)} \cdot \Jac(\Phi^{-1}\circ L_g\circ\Phi)_{P(x)}\cdot \Jac(P)_x\\
		= \Jac(\Phi^{-1}\circ L_g\circ\Phi)_{P(x)} 
		= 1,
	\end{multline*}
where we finally used that by Proposition~\ref{CG-1.2.9} left-translations in exponential coordinates have Jacobian determinant equal to 1.	Similarly, it holds for the right translations. 
\end{proof}
\begin{Rem}
 	If a map $F:\R^n\to\R^n$ has Jacobian $1$, then it preserves the Lebesgue $n$-measure
	(because of the change-of-variables formula).
\end{Rem}

\subsubsection{Measures in nilpotent Lie groups}

On nilpotent groups, we have an interpretation in coordinates for the Haar measure; we discussed such a notion in Section~\ref{sec_Haar_poly_growth}.
An immediate consequence of Corollary~\ref{polynomial coordinates Jacobian 1} is the following result.\label{bi-invariant Haar measure}
\begin{corollary}
\label{CG-1.2.10}
Let $G$ be an $n$-dimensional connected, simply connected, and nilpotent Lie group. 
Every polynomial coordinate system pushes forward the Lebesgue measure on $\R^n$ to a bi-invariant Haar measure on $G$. 
\end{corollary}
It is not always true that left-Haar measures are also right-Haar measures; groups with such property are called {\em unimodular}. However, in every nilpotent Lie group, Haar measures are both left and right-invariant. Corollary~\ref{CG-1.2.10} shows such property for nilpotent simply connected Lie groups.

To summarize, we provide the proof of Theorem~\ref{Malcev global coordinates}.

\begin{proof}[Proof of Theorem~\ref{Malcev global coordinates}] 
In every nilpotent simply connected Lie group, exponential coordinates are global coordinates by Theorem~\ref{CG-1.2.1}.i, and by BCH formula, the group product in such coordinates is polynomial. In Proposition~\ref{CG-1.2.9}, we saw that, in exponential coordinates, left and right translations have Jacobian determinants equal to $1$, and so they preserve the Lebesgue measure.

Regarding Malcev coordinates, in Proposition~\ref{CG1.2.7}, we saw that they are polynomial coordinates in the sense of Definition~\ref{def polynomial coordinates}; so in particular, they also are global coordinates, and still, the group product in such coordinates is polynomial., 

In Corollary~\ref{polynomial coordinates Jacobian 1}, we saw that, in polynomial coordinates, translations have Jacobian determinant equal to 1, and hence they preserve the Lebesgue measure, which is, therefore, a left-Haar measure and a right-Haar measure. 
\end{proof}

\subsection{Structure of connected nilpotent Lie groups}
In this section, we discuss connected nilpotent groups that are not necessarily simply connected.
Recall that nilpotent Lie algebras can be equipped with the Dynkin product, as in Definition~\ref{Dynkin_product}, which is a polynomial group structure by Proposition~\ref{prop:exists_group_algebra_nilpotent2}.
\begin{theorem}\label{thm:Structure}\index{Dynkin product}
Let $G$ be a connected nilpotent Lie group with Lie algebra $\mathfrak{g}$. Equip $\mathfrak{g}$ with the Dynkin product $\star$. Then, we have the following properties:
\begin{description}
 \item[\ref{thm:Structure}.i.] The universal cover $\Tilde{G}$ of $G$ is isomorphic to $(\mathfrak{g}, \star)$, and
 we have Lie algebra isomorphisms $\mathfrak{g} \cong \Lie(\Tilde{G})\cong \Lie(\mathfrak{g}, \star) $.
 In fact, the exponential map $\mathrm{exp}:\mathfrak{g}\longrightarrow G$ is a covering (the universal covering) and $\mathrm{exp}:(\mathfrak{g}, \star)\longrightarrow G$ is a Lie group homomorphism. In particular, the map $\mathrm{exp}$ is surjective.
 \item[\ref{thm:Structure}.ii.] The center $Z(G)$  equals $\mathrm{exp}(Z(\mathfrak{g}))$, and it is connected.
 \item[\ref{thm:Structure}.iii.] There exists a discrete subgroup $\Gamma < Z(\mathfrak{g})$ such that $$G\cong (\mathfrak{g}, \star)/ \Gamma .$$ In fact, the group $\Gamma$ is the kernel of $\mathrm{exp}:(\mathfrak{g}, \star)\longrightarrow G$. In particular, the manifold $G$ is diffeomorphic to the abelian Lie group $(\mathfrak{g},+)/\Gamma$, as manifolds.
 \item[\ref{thm:Structure}.iv.] The set $\mathfrak t:=\mathrm{span} (\Gamma)\subseteq Z(\mathfrak{g})$ is a central Lie subalgebra for which $T:=\mathrm{exp}(\mathfrak t)$ is a torus. Finally, the manifold $G$ is diffeomorphic to the product manifold $(G/T)\times T$.
\end{description}
\end{theorem}
\begin{proof}
Let us prove the items.
\begin{enumerate}
 \item[i)] We saw in Proposition~\ref{prop:exists_group_algebra_nilpotent2} that $(\mathfrak{g}, \star)$ is a simply connected Lie group with Lie algebra $\mathfrak{g}$. Using Theorem~\ref{thm_induced_homo}, let $\phi:(\mathfrak{g}, \star)\longrightarrow G$ be the unique Lie group homomorphism such that
 $$
 \begin{tikzcd}
 \mathfrak{g}=\mathrm{Lie}(\mathfrak{g}, \star)\ar[r,"\id"]\ar[d,swap,"\mathrm{exp}_{(\mathfrak{g}, \star)}"] &\mathfrak{g}=\mathrm{Lie} (G)\ar[d,"\exp_G"]\\
 (\mathfrak{g}, \star)\ar[r,"\phi"] &\qquad G\qquad.
 \end{tikzcd}
 $$
 
 Since we have that $\mathrm{exp}_{(\mathfrak{g}, \star)}=\mathrm{id}$, then $\phi=\mathrm{exp}_G$.
 \item[ii)] The fact that $\mathrm{exp}(Z(\mathfrak{g})) \subseteq Z(G) $ is a general fact for connected Lie groups (Exercise~\ref{ex equivalence commutativity}).
Vice versa, if  $ g \in  Z(G)$ since $\exp$ is surjective there is $X\in \g$ such that 
$g=\exp(X)$. By the central proprety we have $\Id = \Ad_g = e^{\ad_X}$, by Formula~\ref{Formula: Ad: ad}.
By Proposition~\ref{prop_basic_properties_nilpotent_Lie_albegra}.iv, the transformation $\ad_X$ is nilpotent. 
By Proposition~\ref{exp log on U N}, we conclude that $\ad_X=0$, so $X\in  Z(\mathfrak{g})$.
  Being $Z(\mathfrak{g})$ a vector subspace, it is connected, and so is $\mathrm{exp}(Z(\mathfrak{g}))$.
  
 \item[iii)] Since $\mathrm{exp}:(\mathfrak{g}, \star)\longrightarrow G$ is a covering of Lie groups, then there exists $\Gamma < (\mathfrak{g}, \star)$ with $\Gamma\cong \pi_1(G)$, such that $\Gamma$ is discrete and central, and $$G\cong(\mathfrak{g}, \star)/\Gamma .$$ Indeed, we take $\Gamma:=\mathrm{exp}^{-1}(1_G)$. We stress that normal discrete subgroups of connected groups are central (see Exercise~\ref{Normal_discrete_connected_central}). 
 The quotient $(\mathfrak{g}, \star)/\Gamma$ is a Lie group isomorphic to $G$. However, for $x\in \mathfrak{g}$ and $\gamma\in\Gamma\subseteq Z(\mathfrak{g})$ we have $$ x\star\gamma=x+\gamma, $$ where the extra brackets in Dynkin product are zero since $\gamma$ is central. 
 Thus, we may identify the manifold $G\cong (\mathfrak{g}, \star)/\Gamma$ with the abelian Lie group $(\mathfrak{g},+)/\Gamma,$ only the group structures may be different.
 \item[iv)] We have that $\Gamma$ is a discrete set that generates $\mathfrak{t}:=\mathrm{span}( \Gamma)$. Thus $T:=\mathfrak{t}/\Gamma\cong \mathrm{exp}(\mathfrak{t})$ is a torus. Moreover, the subgroup $T$ is central (hence normal), so that $G/T\cong (\mathfrak{g}, \star)/\mathfrak{t}$ is also a Lie group, which is isomorphic to $ (\mathfrak{g}/\mathfrak{t}, \star)$.
 Let $\mathfrak{m}\subseteq \mathfrak{g}$ be a subspace such that $\mathfrak{g}=\mathfrak{m}\oplus\mathfrak{t}$, then for $x\in\mathfrak{m}$ and $z\in\mathfrak{t}\subseteq Z(\mathfrak{g})$ we have $x\star z =x+z$. Thus the map $\mathfrak{m}\times\mathfrak{t}\longrightarrow \mathfrak{g}$ given by $(x,z)\mapsto x\star z =x+z$ is a diffeomorphism. Hence $$\begin{aligned}
 &\mathfrak{m}\times \mathfrak{t}/\Gamma\longrightarrow G\\
 &(x,z+\Gamma)\longmapsto x\star(z+\Gamma)=x+z+\Gamma
 \end{aligned} $$ is a diffeomorphism. Since the quotient map $\mathfrak{g}\longrightarrow \mathfrak{g}/\mathfrak{t}$ induces a linear isomorphism $\mathfrak{m}\longrightarrow\mathfrak{g}/\mathfrak{t}$, we obtain diffeomorphisms $$G\cong \mathfrak{m}\times \mathfrak{t}/\Gamma \cong G/T\times T.$$ In fact, we conclude that $G\cong (\mathfrak{g}, \star)/\Gamma\cong \mathfrak{m}\times \mathfrak{t}/\Gamma\cong \mathfrak{g}/\mathfrak{t}\times T\cong G/T\times T.$\qedhere
\end{enumerate}
\end{proof}
\begin{corollary}
Every compact connected nilpotent Lie group is abelian
\end{corollary}
\begin{proof}
From the previous result, Theorem~\ref{thm:Structure}, each connected nilpotent Lie group $G$ is diffeomorphic to $\mathfrak{g}/\Gamma$ for some $\Gamma\subseteq Z(\mathfrak{g})$. If in addition $\mathfrak{g}/\Gamma$ is compact, then $\Gamma$ must span $\mathfrak{g}$ and then $\mathfrak{g} \subseteq Z(\mathfrak{g})$. Hence, $\mathfrak{g}$ is abelian, and so is $G$, recall Exercise~\ref{ex equivalence commutativity}.
\end{proof}
\begin{remark}
The only abelian connected compact Lie groups are isomorphic as Lie groups to $\R^n/\Z^n=:\mathbb{T}^n$, for $n\in \N$.
\end{remark}

\subsection{Torsion in nilpotent Lie groups}\label{sec:torsion}\index{free! torsion --}\index{torsion}
\begin{definition}[Torsion]\index{torsion}Given a group $G$ and an element $g\in G$, we say that $g$ is a {\em torsion element} if it has finite order, say $k\in\mathbb N$, i.e., $g^k=1_{G}$. (Notice that $1_G$ is a torsion element).
The subset $\mathrm{Tor}(G)$ of $G$ consisting of the torsion elements of $G$ is called {\em torsion} of $G$. If $\mathrm{Tor}(G)=\{1_G\}$ we say that $G$ is {\em torsion free} or that it has {\em no torsion}.
\end{definition}

We warn that the set $\mathrm{Tor}(G)$ may not be a subgroup for general groups. The set $\mathrm{Tor}(G)$ is a subgroup when $G$ is nilpotent. 
We close this section with a characterization: a connected nilpotent Lie group is simply connected if and only if it has no torsion.

\begin{proposition}\label{simply_con_Torsion} $ $
 \begin{description}
 \item[\ref{simply_con_Torsion}.i] 
 Nilpotent simply connected Lie groups, and their subgroups, are torsion-free.
 \item[\ref{simply_con_Torsion}.ii] Every connected nilpotent Lie group that has no torsion is simply connected.
\end{description}
\end{proposition}
\begin{proof}
Regarding {\ref{simply_con_Torsion}.i}, recall, from Theorem~\ref{CG-1.2.1}.i, that in every nilpotent simply connected Lie group, the map $\exp$ is a diffeomorphism. Assume we have an element expressed in the form $\exp(X)$ and for some $m\in \N$ we have
\begin{equation}\label{eqn:Torsion}
\mathrm{exp}(0)=\mathrm{exp}(X)^m=\mathrm{exp}(mX). 
\end{equation}
Then $mX=0$ and so $X=0$.
 Regarding part \ref{simply_con_Torsion}.ii, by Theorem~\ref{thm:Structure}.iii we have that $G\cong (\mathfrak g, \star)/\Gamma$. 
 Assume by contradiction that $\Gamma$ is not trivial, and take $\gamma\in\Gamma$ with $\gamma\neq 1$. Thus there exists $X\in\mathfrak{g}$ such that $\gamma=\mathrm{exp}(X)$ and $X\neq 0$. By the fact that $\Gamma$ is discrete, there exists $k\in\mathbb N$ such that $g_k:=\mathrm{exp}(\frac{1}{k}X)\notin\Gamma$. Hence, the element $g_k$ has finite order $ k$, which is a contradiction.
\end{proof}

\section{Exercises}
\begin{exercise}Let $\g^{(i)}$ be the $i$-th element in the lower central series of a Lie algebra $\g$.
\\
i. 	We have $\g^{(i+1)}\subset\g^{(i)}$ for all $i\in \N$.
\\
ii.	If $\g^{(i)}=\g^{(i+1)}$ for some $i$, then for all $j>i$ we have $\g^{(j)}=\g^{(i)}$.
\end{exercise}

\begin{exercise} 
For finite-dimensional Lie algebras, the nilpotency of $\g$ is equivalent to the vanishing of the ideal
$\g^{(\infty)} (\g) := \bigcap_{i\in \N} \g^{(i)}(\g)$.
\end{exercise}

 \begin{exercise}[Upper central series of Lie algebras]\label{upper central series algebra}
Let $\g$ be a Lie algebra. One iteratively defines the elements of the {\em upper central series} $(\zeta_i(\g))_{i\in\N}$ of $\g$ by
\index{upper central series! of a Lie algebra}
$$\zeta_0(\g):=\{0_\g\} \quad 
\text{ and } \quad \zeta_{i+1}(\g):=\left\{ X\in \g: [X, \g]\subseteq \zeta_i(v) \right\} .$$
The subset $\zeta_1(\g)$ is the center $ Z(\g)$. The sequence $(\zeta_i(\g))_{i\in\N}$ is weakly increasing, and $[\zeta_i, \g]\subseteq \zeta_{i-1} $.
\end{exercise}


\begin{exercise}\label{ex g F}
 Given a flag $\mathcal F= (V_0, \dots, V_m)$ for a vector space $V$, the {\em set of flag-preserving transformations}\index{flag-preserving transformations}
 $$ \g(\mathcal F) :=\big\{A\in \gl(V) : A(V_k)\subseteq V_{k}, \forall k\in\{1, \ldots, m\}\big\}$$ 
is a Lie algebra that is not nilpotent as soon as $\dim V>1$. Moreover, the subset $ \g_{\rm nil} (\mathcal F)$, as defined in Example~\ref{ex g nil F}, is an ideal of $ \g(\mathcal F)$.
\end{exercise} 

\begin{exercise}\label{ex 11 1 2024} The set $ \g_{\rm nil} (\mathcal F)$ of flag-shifting transformations for a flag $\mathcal F= (V_0, \dots, V_m)$, as in Example~\ref{ex g nil F}, is a Lie algebra and 
$$C^\ell( \g_{\rm nil} (\mathcal F)) V_k \subseteq V_{k-\ell}, \qquad \forall k, \ell\in\N.$$
Deduce that 
$C^m( \g_{\rm nil} (\mathcal F)) =\{0\}$, 
and that $ \g_{\rm nil} (\mathcal F)$ is nilpotent, with step at most $m-1$.
\end{exercise} 

\begin{exercise} For $n=2,$ the Lie algebra $\Lambda ^1 (\R^2) \times \Lambda ^2 (\R^2)$ from \eqref{model_free} gives the first Heisenberg Lie algebra of Example~\ref{example Heisenberg algebras}.
\end{exercise}

\begin{exercise}\index{free! -- nilpotent}
	Every function defined on the set of generators of a free Lie algebra to a nilpotent Lie algebra of no larger step extends uniquely to a Lie algebra homomorphism.
\end{exercise}

\begin{exercise}
Let $A, B\in \mathfrak {gl}(V).$ Assume that $A$ and $B$ commute. 
If $A$ and $ B $ are nilpotent transformations, then so is $A+B$. 
If $A$ and $ B $ are unipotent transformations, then so is $AB$. 
\end{exercise}

\begin{exercise}[Real Jordan Theorem]\label{Jordan Theorem}\index{Jordan form}
 If $A\in\gl(V)$, there is a basis of $V$ such that the matrix representation of $A$ in this basis is in \emph{real Jordan form}, i.e.,
 it is written in blocks
 \begin{equation}\label{eq6703c1e5}
 \begin{pmatrix}
 A_1 & 0 & \dots \\
 0 & A_2 & \dots \\
 \vdots & & \ddots
 \end{pmatrix}
 \end{equation}
 where each $A_j$ is a matrix of one of the following four forms:

 \begin{equation}\label{eq6703c0c7}
 \begin{pmatrix}
 \alpha & 0 & \dots \\
 0 & \alpha & \dots \\
 & \vdots & \ddots
 \end{pmatrix} ,
\qquad
 \begin{pmatrix}
 \alpha & 1 & 0 & \dots \\
 0 & \alpha & 1 & \dots \\
 & \vdots & & \ddots \\
 \end{pmatrix} ,
 \end{equation}

 \begin{equation}\label{eq6703c0d4}
 \begin{pmatrix}
 \alpha & -\beta & 0 & 0 & \dots \\
 \beta & \alpha & 0 & 0 & \dots \\
 0 & 0 & \alpha & -\beta & \dots \\
 0 & 0 & \beta & \alpha & \dots \\
 & \vdots & & & \ddots
 \end{pmatrix} ,
\qquad
 \begin{pmatrix}
 \alpha & -\beta & 1 & 0 & 0 & 0 & \dots \\
 \beta & \alpha & 0 & 1 & 0 & 0 & \dots \\
 0 & 0 & \alpha & -\beta & 1 & 0 & \dots \\
 0 & 0 & \beta & \alpha & 0 & 1 & \dots \\
 & \vdots & & \vdots & & & \ddots
 \end{pmatrix} .
 \end{equation}
 The blocks as in~\eqref{eq6703c0c7} correspond to a real eigenvalue $\alpha$, while the blocks as in~\eqref{eq6703c0d4} correspond to a complex eigenvalue $\alpha+i\beta$. Consequently, the matrix $A$ is nilpotent if and only if all blocks are upper triangular.
\end{exercise}

\begin{exercise}\label{ex solvable not nilpotent}
Use any of the following points to deduce that if $\mathfrak{a}\subseteq \g$ is an ideal of a Lie algebra and both $\mathfrak{a}$ and $ \g/\mathfrak{a}$ are nilpotent, then $\g$ may not be nilpotent.
\begin{description}
\item[\ref{ex solvable not nilpotent}.i.] 
 Let $\g$ be the Lie algebra spanned by two vectors $X$ and $Y$ with relation $[X, Y]=X$. Then, the subset $\mathfrak{h} :=\R X$ is a commutative ideal, $\g/\mathfrak{h}$ is commutative, but $\g$ is not commutative, nor nilpotent.
\item[\ref{ex solvable not nilpotent}.ii.] Let $\g$ be a solvable Lie algebra (see Exercise~\ref{ex_def_solvable} for the definition).
Then, both $[\g, \g]$ and $\g/ [\g, \g]$ are nilpotent.\index{solvable}\index{commutator! -- subgroup} 
\end{description}
\end{exercise}


\begin{exercise}
	Every stratification is completely determined by $V_1$.
\end{exercise}

\begin{exercise}
Every compatible linear grading that is also a Lie algebra grading is a stratification. 
\end{exercise}

\begin{exercise}\label{ex stratification 2023} If a Lie algebra $\g$ has a $s$-step stratification $\g=V_1\oplus \cdots\oplus V_s $,
then
 \\\ref{ex stratification 2023}.i. $V_s$ is contained in the center of $\g$;
 \\\ref{ex stratification 2023}.ii. $V_k\oplus \cdots\oplus V_s $ is an ideal in $\g$;
 \\\ref{ex stratification 2023}.iii. $(V_k\oplus \cdots\oplus V_s )/(V_{k+1}\oplus \cdots\oplus V_s)$ is contained in the center of $(V_1\oplus \cdots\oplus V_s)/(V_{k+1}\oplus \cdots\oplus V_s)$.
 \\\ref{ex stratification 2023}.iv. $\g/V_s$ is a stratifiable Lie algebra of step $s-1$.
\end{exercise}

\begin{exercise}
Let $\g$ be a Lie algebra that admits a grading.
Assume that the elements of degree $1$, namely $V_1$, generate $ \g$, as a Lie algebra, then $\g$ is stratified by $V_1, \ldots, V_s$.
\\{\it Hint:} If the bracket of $V_1$ with itself were smaller than $V_2$, then $V_1$ would not generate because the Lie subalgebra it generates will not contain all of $V_2$. Proceed by induction.
\end{exercise}

\begin{exercise}\index{Heisenberg! -- algebra} 
Let $\h$ be the Heisenberg Lie algebra generated by the vectors $X$, $Y$, and $Z$ with only non-trivial relation $[X, Y]=Z$.
Then, the Lie algebra $\h$ can be stratified as
$\h= \Span\{X,Y\}\oplus \Span\{Z\},$
and this stratification has step 2.
\end{exercise}
\begin{exercise}\label{ex:lagercenter}
Let $\g:=\R\times \h$ be the direct product of $\R$ with the (above) Heisenberg Lie algebra $\h$.
Then, the Lie algebra $\g$ can be stratified as
$$\g= (\R\times\Span\{X,Y\})\oplus (\{0\}\times\Span\{Z\});$$
still, its center is $\R\times\Span\{Z\}$, which is strictly bigger than $V_2$.
\end{exercise}

\begin{exercise}
	Every step-two nilpotent Lie algebra is stratifiable.
\end{exercise}

\begin{exercise}[Positively graded algebras are nilpotent]\label{ex positive gradings give nilpotency}
Every finite-dimensional Lie algebra that admits a positive grading is nilpotent.\\
{\it Hint.} Let $\g= \bigoplus_{t>0} V _t =\bigoplus_{a\in [a,b]} V _t $, with $0<a<b$. Then if $ b/a<m\in \N$, then $C^m(\g)=\{0\}$.
\end{exercise}

\begin{exercise}\label{ex:v2v2}For every Lie-algebra stratification $V_1\oplus\dots\oplus V_s$, we have $[V_2,V_2] \subseteq V_4$.
\\{\it Solution}. We have 
	\begin{multline*}
	 	[V_2,V_2] = [[V_1,V_1],[V_1,V_1]] 
		= \Span\left\{[[X_1, X_2],[X_3, X_4]]: X_i\in V_1\right\} \subset \\
		\overset{\text{(Jacobi)}}\subset \Span\left\{[X_1,[X_2,[X_3, X_4]]]:X_i\in V_1\right\}
		= [V_1,[V_1,[V_1,V_1]]] = V_4,
	\end{multline*}
where we used the properties of the stratification and Jacobi identity from Definition~\ref{def: Lie_algebra}.
\end{exercise}

\begin{exercise}[Stratifications are positive gradings]\label{ex stratifications are gradings}
 	Let $\g=V_1\oplus\dots\oplus V_s$ be a stratified Lie algebra.
	Then, setting $V_k:=\{0\}$ for $k>s$, we have
	\[
	[V_i,V_j]\subset V_{i+j}, \qquad\forall i,j\in\{1, \dots,s\}.
	\]
{\it Solution.} 	The proof is by induction on $i$. 
	If $i=1$ we know that $[V_1,V_j]\subset V_{j+1}$ for all $j$.
	Suppose that $[V_i,V_j]\subset V_{i+j}$ for all $j$ and a fixed $i$.
	We shall check that $[V_{i+1},V_j]\subset V_{i+1+j}$ for all $j$.
	Indeed, the space $V_{i+1}$ is generated by the elements $[v_1,v_i]$ where $v_1\in V_1$ and $v_i\in V_i$, and for these elements we have for all $v_j\in V_j$ the Jacobi identity:
	$
	[[v_1,v_i],v_j] = -[[v_i,v_j],v_1] - [[v_j,v_1],v_i], 
	$
	where $[v_i,v_j]\in V_{i+j}$ by the inductive hypothesis and so $-[[v_i,v_j],v_1] = [v_1,[v_i,v_j]]\in[V_1,V_{i+j}]= V_{i+1+j}$, and $ - [[v_j,v_1],v_i] = [v_i,[v_j,v_1]] \in [V_i,V_{j+1}] \subset V_{i+1+j}$ by the inductive hypothesis again.
	All in all, $[[v_1,v_i],v_j]\in V_{i+1+j}$ and therefore $[V_{i+1},V_j]\subset V_{i+1+j}$.
	\end{exercise}

\begin{exercise}[Elements of lower central series in terms of stratifications]\label{lem051025} Let $\g$ be a Lie algebra with a $s$-step stratification $\g=V_1\oplus \cdots\oplus V_s $. Then, we have
 $$\g^{(k)}=V_k\oplus \cdots\oplus V_s .$$
{\it Solution.}
 	The proof is by induction.
	For $k=1$, it is trivial.
	Suppose it is true for $k$, then 
	\begin{multline*}
	\g^{(k+1)} 
	= [\g, \g^{(k)}] 
	= [V_1\oplus\dots\oplus V_s, V_k\oplus\dots\oplus V_s ] = \\
	= \sum_{i=1}^s \sum_{j=k}^s [V_i,V_j] 
	= \sum_{j=k}^s [V_1,V_j] + \sum_{i=2}^s \sum_{j=k}^s [V_i,V_j] = \\
	= V_{k+1}\oplus\dots\oplus V_s + \sum_{i=2}^s \sum_{j=k}^s [V_i,V_j]
	= V_{k+1}\oplus\dots\oplus V_s,
	\end{multline*}
	where $\sum_{i=2}^s \sum_{j=k}^s [V_i,V_j] \subset \sum_{i=2}^s \sum_{j=k}^s V_{i+j}\subset V_{k+1}\oplus\dots\oplus V_s$, by Exercise~\ref{ex stratifications are gradings}.
\end{exercise}

\begin{exercise}
 If a Lie algebra $\g$ admits an $s$-step stratification, then $\g$ is $s$-step nilpotent. 
\\{\it Hint.} Check Exercise~\ref{lem051025}.
 \end{exercise}

\begin{exercise}[Positively gradable, non-stratifiable Lie algebra]\label{ex N51}
The Lie algebra $\mathfrak{n}_{5,1}$ from Example~\ref{N51} is positively gradable, yet not stratifiable.
\end{exercise}

\begin{exercise}[Another nilpotent nonstratifiable Lie algebra]\label{nonCarnot-example}
Consider the $7$-dimensional Lie algebra $\mathfrak{h}$ generated by $X_1, \ldots, X_7$ with only nontrivial brackets
$$ 
\begin{array} {rcl c rcl}
{[X_1, X_2]} &= &X_{ 3} &\quad&
{[X_1, X_3]} &=&2 X_{ 4} \\
{[X_1, X_4]} &= &3 X_{ 5} &&
{[X_2, X_3]} &= &X_{ 5} \\
{[X_1, X_5]} &= &4 X_{ 6} &&
{[X_2, X_4]} &= &2 X_{ 6} \\
{[X_1, X_6]} &= &5 X_{ 7} &&
{[X_2, X_5]} &= &3 X_{ 7} \\
{[X_3, X_4]} &=&X_{ 7} 
\end{array}
$$
The following facts hold true:
\\\ref{nonCarnot-example}.i. $\g$ is a Lie algebra
\\\ref{nonCarnot-example}.ii. $\g$ is nilpotent 
\\\ref{nonCarnot-example}.iii. $\g$ admits a grading, with $V_i=\R X_i$.
\\\ref{nonCarnot-example}.iv. For every grading of $\g$, the elements of degree $1$, $V_1$, do not generate $ \g$.
\\\ref{nonCarnot-example}.v. $\g$ does not admit any stratification.
\end{exercise}

\begin{exercise}\label{ex 5D non-gradable}
[Nilpotent Lie algebras that are not positively gradable]
Consider the following $7$-dimensional Lie algebras denoted by 
 $\mathfrak{h}_{12457G}$ and $\mathfrak{h}_{12357B}$. They are generated by $X_1, \ldots, X_7$ with only nontrivial brackets
$$ 
\begin{array} {l}
\text{ {\it Relations for } 12457G}\\
{[Y_1, Y_2]}	=	Y_5\\
{[Y_1, Y_3]}	=	Y_6\\
{[Y_1, Y_4]}	=	Y_6\\
{[Y_1, Y_7]}	=	-Y_2\\
{[Y_2, Y_3]}	=	-Y_6\\
{[Y_2, Y_7]}	=	-Y_3\\
{[Y_3, Y_7]}	=	-Y_4\\
{[Y_5, Y_7]}	=	-Y_6
\end{array}
\quad, \qquad 
\begin{array} {l}
\text{ {\it Relations for } 12357B}\\
 {[Y_1, Y_6]}	=	-Y_2\\
{[Y_1, Y_7]}	=	Y_3 + Y_5\\
{[Y_2, Y_6]}	=	-Y_3\\
{[Y_2, Y_7]}	=	Y_4\\
{[Y_3, Y_6]}	=	-Y_4\\
{[Y_3, Y_7]}	=	Y_5\\
{[Y_4, Y_6]}	=	-Y_5\\
\quad
\end{array}.
$$
The following facts hold true:
\\\ref{ex 5D non-gradable}.i. both $\mathfrak{h}_{12457G}$ and $\mathfrak{h}_{12357B}$ are Lie algebras,
\\\ref{ex 5D non-gradable}.ii. they are nilpotent,
\\\ref{ex 5D non-gradable}.iii. $\skull$ none of them admit a positive grading. 
\\{\it Hint.} Check \cite{MR4394318}.
\end{exercise}

\begin{exercise}\label{ex family non-stratifiable}
Fix a positive
integer $n\geq 7$, and consider the $n$-dimensional Lie algebra $\mathfrak{h}$ generated by $X_1, \ldots, X_n$ with
$$[X_i, X_j] =\left\{
\begin{array} {cc}
( j - i ) X_{ i + j}, &\text{if}\, i+j\leq n, \\
0, &\text{otherwise.}
\end{array} \right.
$$
The following facts hold true:
\\\ref{ex family non-stratifiable}.i. it is a Lie algebra
\\\ref{ex family non-stratifiable}.ii. it is nilpotent
\\\ref{ex family non-stratifiable}.iii. it does not admit any stratification.
\end{exercise}


\begin{exercise}[A nontrivial filiform algebra]\label{ex nontrivial filiform}\index{filiform! -- algebra of the second kind}
Consider the $6$-dimensional Lie algebra $\mathfrak{g}$ 
 given by $\Span \{ y_0,y_1,y_2, y_3,y_4, y_5\}$ with only 
non-zero brackets
$$\begin{array}{rcl c rcl}
[y_0,y_1]&=&y_2,&\quad&
{[y_0,y_2]}&=&y_3, \\
 {[y_0,y_3]}&=&y_4,&&
 { [y_0,y_4]}&=&y_5, \\
{[y_1,y_4]}&=&-y_5, &&
{[y_2,y_3]}&=&y_5.
\end{array}$$
The following facts hold true:
\\\ref{ex nontrivial filiform}.i. It is a Lie algebra, i.e., Jacobi identity is satisfied.
\\\ref{ex nontrivial filiform}.ii. It admits a stratification.
\\\ref{ex nontrivial filiform}.iii. It is a filiform algebra (i.e., the dimensions of the subspaces of the stratification are the smallest possible, namely $2,1, \dots, 1$).
\end{exercise}

\begin{exercise}[Suggested by E. Breuillard]\label{ex Breuillard strange V 1}
Let $\g$ be the quotient of the free nilpotent   Lie algebra of $3$-step and generated by $e_1,e_2,e_3$ modulo the ideal generated by $[e_2,e_3] $, so we have the additional relation $[e_2,e_3]=0$. 
The Lie algebra $\g$ has dimension $10$. Moreover, the following is a stratification of $\g$:
$$\begin{array}{rcl}
V_1&:=& \Span\{ e_1, e_2, e_3\}, \\
V_2&:=&\Span\{ [e_1,e_2], \quad [e_1,e_3]\}, \\
V_3&:=&\Span \{ [e_1,[e_1,e_2]], \quad [e_2,[e_1,e_2]], \quad[e_3,[e_1,e_2]], 
\quad [e_3,[e_3,e_1]], \quad [e_1,[e_1,e_3]] \}.
\end{array}$$
The subspace $V'_1:= \Span\{ e_1,e_2 + [e_1,e_2],e_3\}$ is not the first layer of a stratification. Still, it is in direct sum with $[\g, \g]$. \\ {\it Hint.} The space $[V'_1,V'_1]$ is not in direct sum with $[\g,[\g, \g]]$, because it contains $[e_3, e_2 + [e_1,e_2]]=[e_3,[e_1,e_2]]$, and in fact $V'_2:=[V'_1,V'_1]$ has dimension $3$, not $2$.
\end{exercise}

\begin{exercise}\label{V>0}
 Let $(V_a)_{a\in \R}$ be a Lie algebra grading of a Lie algebra $\g$. 
 Then, the subspace $V_{>0}:= \bigoplus_{a>0} V _a$ is a subalgebra.
 Moreover, if the grading is {\em nonnegative}, \index{nonnegative grading}
 in the sense that $\g = \bigoplus_{a\geq0} V _a$, then $V_{>0}$ is an ideal.
 \\{\it Solution.}	If $a\geq 0$ and $b > 0$, then obviously $a+b>0$, and so $[V_a,V_b]\subseteq V_{a+b} \subseteq V_{>0}.$
\end{exercise}

 \begin{exercise} Every one-parameter subgroup of Lie algebra automorphisms leads to a Lie algebra grading as in Proposition~\ref{prop12061739}.
 \\{\it Hint.} Use Propositions \ref{Prop_Lie_der}.ii and \ref{prop08220937}.
 \end{exercise}
 
 \begin{exercise}
Let $\g= V_1\oplus \cdots\oplus V_s$ be a stratified algebra.
For all $\lambda\geq 0$, let $\delta_\lambda$ be the dilation of factor $\lambda$ as defined in Definition~\ref{Dilations:algebras}.
Then, we have  
$$\delta_\lambda\left(\sum_{i=1}^s v_i\right) =\sum_{i=1}^s\lambda^i v_i, $$
where $X=\sum\limits_{i=1}^s v_i$ with $v_i\in V_i$, $1\leq i\leq s$.
\end{exercise}


\begin{exercise}[Nilpotency and compatible linear gradings]\label{ex Nilpotency and compatible linear gradings}
Let $\g$ be a finite-dimensional Lie algebra.
We have that $\g$ is nilpotent if and only if it admits a compatible linear grading. 
Moreover, if $\mathfrak{g} = \bigoplus_{i=1}^\infty V_i$ is a compatible linear grading, 
then 
$$[V_i, V_j]\subseteq [\g^{(i)}, \g^{(j)}]\subseteq \mathfrak{g}^{(i+j)}= \bigoplus_{k=i+j}^\infty V_k.$$ 
\end{exercise}

\begin{exercise}\label{stratifiable:isomorphic:graded}
Every stratifiable Lie algebra is isomorphic to its associated Carnot algebra
(as defined in Definition~\ref{Graded:algebra}).
\end{exercise}

\begin{exercise}\label{ex positively graded contraction}\index{contractive}\index{positively graded! -- Lie algebra} If $\g$ is a finite-dimensional positively graded Lie algebra, then for each 
 $\lambda\in (0,1)$ the dilation $\delta_\lambda$ is contractive.
 If $\g$ is $\N$-graded, then 
 the dilation $\delta_\lambda$ is contractive for 
 $\lambda\in (-1,1)$. 
 \\ 
 {\it Hint.} We have $\delta_\lambda^n(v) = \lambda^{jn} v \to 0$, as $n\to \infty$, for every $v$ of degree $j$.
\end{exercise}

\begin{exercise}
Let $\delta: G\to G$ be an automorphism of a Lie group.
Then, the sequence of maps $(\delta^n)_{n\in \N}$ pointwise converges to the constant map $1_G$ if and only if it converges uniformly on compact sets to $1_G$. \\{\it Hint.} Sequences of linear maps that converge on the elements of a basis converge uniformly on compact sets.
\end{exercise}

\begin{exercise}\label{ex contractive contractive}
Let $G$ be a 
connected Lie group with Lie algebra $\g$.
Let $\delta: G\to G$ be a Lie group automorphism with induced Lie algebra automorphism $\delta_*:\g \to \g$.
Then, the map $\delta$ is contractive if and only if so is $\delta_*$. (Recall Proposition~\ref{prop: expF: Fexp})
\end{exercise}

\begin{exercise} If a Lie group $G$ admits a contractive automorphism, then it is simply connected.
\\{\it Hint.} Manifolds are locally simply connected.
\end{exercise}

\begin{exercise}\label{lemBurba}
	If $\phi\in\Aut_\C(\g _\C)$ and $\alpha, \beta\in\C$, then for the $E^{\phi}_\alpha$'s defined in \eqref{def_generalized_eigenspace} we have
	\begin{equation*}
	[E^{\phi}_\alpha,E^{\phi}_\beta] \subset E^{\phi}_{\alpha\beta}, \qquad \forall \alpha, \beta\in\C .
	\end{equation*}
{\it Solution.}	The proof is elementary after one proves by induction on $n\in \N$ that
	\[
	(\phi-\alpha\beta\mathbb{I})^n[v,w] 
	= \sum_{\substack{j+k=n\\j,k\ge0 }} \binom{n}{j} [\alpha^k(\phi-\alpha\mathbb{I})^jv, \phi^j (\phi-\beta\mathbb{I})^k w] 
	\]
	holds for all $v,w\in\g _\C$, all $\alpha, \beta\in\C$ and all $n\in\N$.
	See also \cite[Ch.7, sec 1.4, Prop.12, p. 11]{Bourbaki_Lie_7_9}. 
\end{exercise}


\begin{exercise}\label{ex_grading_g_R}
 Let $\g$ be a $\Z$-graded Lie algebra.
 \begin{description}
 \item[\ref{ex_grading_g_R}.i.]
The Lie algebra $\g \rtimes \R$ is $\Z$-graded by \eqref{grading_g_R}.
 \item[\ref{ex_grading_g_R}.ii.] If $V_0= \{0\}$ and $\g \ne \{0\},$ then $\g \rtimes \R$ has trivial center.
 \item[\ref{ex_grading_g_R}.iii.]
 If $\g=\{0\}$, then $\g \rtimes \R=\R$, which is abelian.
 \item[\ref{ex_grading_g_R}.iv.] 
If $\g$ is Carnot of step $s$, then the grading of $\g \rtimes \R$ satisfies
\begin{equation*}
V'_m \ne\{0\}\quad \Leftrightarrow \quad m\in \Z \cap [0,s] .
\end{equation*}
\end{description}
{\it Solution.} 
(i). It is enough to check that
 \begin{equation*}
[\{0\} \times \R, \{0\} \times \R ] \subset V'_0, \quad \mbox{and} \quad [V_m \times \{0\}, \{0\} \times \R ] \subset V'_m.
\end{equation*}
In fact, one has $[(0,s), (0,t)] =(0,0) \in V'_0$ and $[(X,0), (0,t)] =(-tmX,0) \in V'_m$ if $X\in V_m.$

(ii). Let $(X,s) \in Z(\g \rtimes \R).$ Write $X$ as $\sum _{m\in \Z} X_m$, with $X_m\in V_m$. Then $(0,0)= [(X,s), (0,1)] = (-\sum _{m\in \Z} mX_m,0).$ Hence, for every $m\ne 0$ we have that $X_m=0$ and so $X=0$ since $V_0= \{0\}.$

Moreover, take $Y \in V_m \setminus \{0\}$ for some $m \in \Z \setminus \{0\},$ which exists since $\g$ is not trivial. Then $(0,0) =[(X,0), (Y,0)] =[(0,s), (Y,0)] = (smY,0)$ and consequently $s=0.$ 

(iii). It trivially follows from the basic Definition~\ref{definition semi-direct Lie algebra}.

(iv). From the definition of $V'_a$, if the non-trivial layers of the gradings of $\g$ are $V_1, \ldots, V_s$, then the
non-trivial layers of the gradings of $\g \rtimes \R$ are $V'_0, V'_1, \ldots, V'_s$.
 \end{exercise}
 
 
 \begin{exercise}\label{ex:UCS_bracket}
Each element $\zeta_i:=\zeta_i(G)$ of the upper central series is normal in $G$ and satisfies $[\zeta_i,G]\subseteq \zeta_{i-1} .$ 
 \end{exercise}
 \begin{exercise} The elements $\zeta_i:=\zeta_i(G)$ of the upper central series satisfy
 $\zeta_{i+1}/\zeta_i=Z\left(G/\zeta_i\right)$.
\end{exercise}

\begin{exercise} 
The Lie group ${\rm {Nil }}_n$ is nilpotent and simply connected with Lie algebra ${\mathfrak{nil}}_n$.
\\ {\it Hint.} Write the obvious diffeomorphism with $\R^{n(n- 1)/2}$. 
 \end{exercise}

\begin{exercise}
	The group ${\rm {Nil }}_3$, as from Example~\ref{example Strictly upper-triangular matrix algebra and group}, is the Heisenberg group, as defined in Section~\ref{sec:Heisenberg}.
\end{exercise}


\begin{exercise}[The maps $\log$ and $\exp $ are the inverse of each others, locally]\label{ex log inverse of exp} 
For all $n\in \N$, we have the following properties:
\begin{description}
\item[\ref{ex log inverse of exp}.i.] For every $A \in \gl (n) $ with $\|A\| < \log 2$ we have that $\log (e^A)=A.$
\item[\ref{ex log inverse of exp}.ii.] For every $M \in\GL(\R^n)$ with $\|M-\mathbb{I}\| < 1$ we have that $e^{\log M}=M.$
\end{description}
{\it Hint.} Check \cite[Proposition 3.3.2]{Hilgert_Neeb:book}.
\end{exercise}

 \begin{exercise}\label{ex span of N} For $\mathcal{N}$ as defined in \eqref{def U N}, the set
 $\Span(\mathcal{N})$ consist of the matrices in $ \gl(n)$ that are 0 in the diagonal. In particular, for all $n>1$, some of them are not in $\mathcal{N}$.
 \end{exercise}
 
 \begin{exercise}
Let $V$ be a finite-dimensional vector space, $v\in V$, and $A \in \gl (V)$ nilpotent transformation. We have that $Av=0$
if and only if $e^Av =v$.
\\{\it Solution.} One direction is true even when $A$ is not nilpotent. For the other one, if $(e^A-\mathbb{I})v=0$ then
$
Av = \log e^A v = \sum_{k=1}^\infty \frac{(-1)^{k+1}}{k} (e^A-\mathbb{I})^{k}v = 0.
$
\end{exercise}

\begin{exercise}\label{ex Dynkin polynomial}\index{Dynkin product} On every nilpotent Lie algebra the Dynkin product 
$(A,B)\in \g\mapsto A\star B $ is a polynomial map.
\end{exercise}

\begin{exercise} In ${\rm {Nil }}_n$, the group product in exponential coordinates is polynomial: $(A,B) \mapsto \log (e^Ae^B).$
\end{exercise}

\begin{exercise} 
Let $A, B \in \gl (V)$ be nilpotent transformations of a vector space $V$. We have that $[A,B]=0$ if and only if
 $e^Ae^B =e^Be^A$.
 \end{exercise}
 

\begin{exercise}\label{ex equivalence commutativity} Let $G$ be a Lie group with Lie algebra $\g$.
Then, the center $Z(G)$ is a closed Lie subgroup with $\Lie(Z(G)) = Z(\g).$ Consequently, if $G$ is connected then $G$ is commutative if and only if $\g$ is commutative. 
\end{exercise}

\begin{exercise}\label{ex845}
If $ [G,G]=\{1_G\},$ then $[\g, \g]=\{0\}.$
\end{exercise} 

\begin{exercise}\label{ex846}
If $[G, [G,G]]=\{1_G\},$ then $[\g, [\g, \g]]=\{0\}.$
\\ {\it Hint}. Use that $[e^{ tx+o(t)}, e^{ ty+o(t)}] = e^{ t^2 [X, Y] + o(t^2)}$, as $t\to0$, for every $x, y \in \g $.
\end{exercise}
%


\begin{exercise}
Let $G$ be a nilpotent simply connected Lie group with Lie algebra $\g$.
The two lower central series, from Definition~\ref{def lower central series} and Definition~\ref{def lower central series groups}, satisfy 
$$ C^m(G) = \exp( C^m(\g)), \qquad \forall m\in \N.$$
Deduce that ${\rm step}(G)={\rm step}(\g)$.
\end{exercise}

\begin{exercise}\label{ex connected nilpotent gives closed and simply connected}
 Every connected Lie subgroup of a nilpotent simply connected Lie group is closed and simply connected.\\ {\it Hint.}
 Check Proposition~\ref{subgroups_Nil} and Theorem~\ref{CG-1.2.1}.iii.
\end{exercise}

\begin{exercise}\label{CG90_1.2.3}
 Every nilpotent simply connected Lie group has an embedding as a closed subgroup of ${\rm {Nil }}_n$, for some $n\in \N$. 
 \\{\it Hint.} Check Theorem~\ref{CG-1.2.1}.iii, Proposition~\ref{subgroups_Nil}, and Theorem~\ref{teo1131}.
\end{exercise}



\begin{exercise}\label{ex Malcev basis in Carnot}
 	Every positively graded Lie algebra admits a Malcev basis. In fact, every ordered basis such that the reversed-ordered basis is adapted to the grading is a  Malcev basis.
\\{\it Solution.}
The solution is a variation of the one from Proposition~\ref{Prop exists Malcev bases}. 	
Let $\g=\bigoplus_{a>0} V _a,$ be a positively graded Lie algebra.
Let $(X_1, \dots, X_n)$ be an ordered basis of $\g$ such that the reversed-ordered basis $(X_n, \dots, X_1)$ is adapted to the direct sum.
	We claim that $(X_1, \dots, X_n)$ is a Malcev basis.
	Indeed, set $\g_k:=\Span\{X_1, \dots, X_k\}$.
	Thus $X_k\in V_b$ for some $b>0$, then 
	\[
	V_{> b} \subset\g_k\subset V_{\geq b} .
	\]
Because we have a positive grading, we conclude
$	 	[\g, \g_k] 
		\subset [\bigoplus_{a>0} V _a, V_{\geq b} ] 
		\subset V_{> b}  
		\subset \g_k$.
\end{exercise}

 
 \begin{exercise} For each prime number $p$, the group
 \[\left\{
 \begin{bmatrix}
1 & x & z\\
0 & 1 & y\\
0 & 0 & 1
\end{bmatrix}: x,y,z\in\mathbb Z/p\mathbb Z\right\}
 \]
is not a subgroup of a connected nilpotent Lie group.
\end{exercise}

\begin{exercise}\label{ex10May1833}
Nilpotent simply connected Lie groups do not have nontrivial compact subgroups. 
 \end{exercise}
 
\begin{exercise}
If $N$ is a nilpotent connected Lie group and $K$ is a compact subgroup of $N$, then $K$ is central. 
\\
{\it Hint.} Take the universal cover, use Theorem~\ref{thm:Structure}, and Exercise~\ref{ex10May1833}.
 \end{exercise}
 



 \chapter{Metrics on nilpotent groups}\label{ch_MetricNilpotent}
 

 In this chapter, we consider metrics on nilpotent Lie groups.
 In Section~\ref{sec:sub-Finsler nilpotent groups}, we consider Carnot-Carath\'eodory metrics.
 In Section\ref{sec: self-similar2}, we complete the description of isometrically homogeneous spaces with dilations, which we started in Section~\ref{sec: self-similar}.
 In Section~\ref{sec:affinity of isometries}, we briefly mention results on the isometry group of nilpotent Lie groups. 
 In Section~\ref{sec:seminorms}, we discuss seminorms on nilpotent Lie groups and show some preliminary results useful for then obtaining Pansu Asymptotic Theorem 
 in Section~\ref{sec_Pansu_asymp}.
 \section{Sub-Finsler nilpotent groups}\label{sec:sub-Finsler nilpotent groups}

Let us now turn to nilpotent simply connected Lie groups equipped with sub-Finsler metrics. 
The notion of sub-Finsler Lie group was introduced in Chapter~\ref{ch_SFLie}.
In Chapter~\ref{ch_Nilpotent}, we discussed the theory of nilpotent Lie groups.

Let $G$ be a Lie group, which we will assume to be simply connected and nilpotent in this chapter.
 Let $V\subseteq T_1 G$ be a vector subspace of the tangent space at the identity element $1_G$ of $G$. Let $\Delta$ be the left-invariant distribution on $G$ with $\Delta_1=V$.
 We shall consider bracket-generating distributions. 
 Fixing a norm $\norm{\cdot}$ on $V$, we extend it to a left-invariant norm on $\Delta$. Then, the triple $(G, \Delta, \norm{\cdot})$ is a sub-Finsler manifold, and we denote by $d_{sF}$ its sub-Finsler metric, which is left-invariant by construction.

 
 The main aim of this section is to show that nilpotent simply connected sub-Finsler groups have a special submetry map: the projection on the abelianization. Through this map, we will have a new viewpoint for developments and multiplicative integrals.

\subsection{A submetry onto the abelianization} 
\label{sec:sFLie_proj}
\subsubsection{The abelianization of nilpotent simply connected Lie groups}

Assume that $G$ is a Lie group that is nilpotent and simply connected, whose Lie algebra is denoted by $\g$. 
In this case, the subgroup $[G, G]$ is connected, normal, and closed; see Exercise~\ref{ex connected nilpotent gives closed and simply connected}. Thus, the quotient $G/[G, G]$ is a Lie group that is abelian and simply connected; see Exercise~\ref{ex G/H simply connected}. Hence, it is a vector space, and it can be identified with its Lie algebra, which is naturally identified with $\g/[\g, \g]$.

The group ${\rm Ab}(G):=G/[G, G] $ is called the {\em abelianization} of the group $G$, while $\g/[\g, \g]$ is called the {\em abelianization} of the Lie algebra $\g$.\index{abelianization} When $G$ is a nilpotent simply connected Lie group, the Lie algebra of ${\rm Ab}(G)$ is isomorphic to $\g/[\g, \g]$ via the diagram
\begin{center}
 \begin{tikzcd}
 \g\ar[r]{u}{\exp}\ar[d,->>] &G\ar[d,->>] &\\
 \hspace{-1cm}{\rm Ab}(G)\simeq \g/[\g, \g] \ar[r]{u}{\exp}&G/[G, G] & \hspace{-1cm} =: {\rm Ab}(G),
 \end{tikzcd}
\end{center}
which is commutative; see Exercise~\ref{ex same ab reason}. In addition, recall from Theorem~\ref{CG-1.2.1}.i that the exponential map is a diffeomorphism. For this reason, we shall identify the two viewpoints for the abelianization and denote by $\pi_{\rm ab}$ either of the projections:
$$\pi_{\rm ab}: \g \twoheadrightarrow{\rm Ab}(G) \quad \text{ or } \quad \pi_{\rm ab}: G \twoheadrightarrow{\rm Ab}(G) ,$$
which we call {\em abelianization maps}.\index{abelianization! -- map}

\subsubsection{A sub-Finsler submetry on the abelianization} 
In addition to $G$ being a nilpotent simply connected Lie group, we assume that it is equipped with a left-invariant 
sub-Finsler structure $(\Delta, \norm{\cdot})$. Thus, the triple
 $(G, \Delta, \norm{\cdot})$ is a sub-Finsler Lie group, as in Section~\ref{sec sub-Finsler Lie groups}.
 We stress that $V:=\Delta_{1_G}$ is assumed to be a bracket-generating subspace of $\g$ and, if we project it via the abelianization map 
 $\pi_{\rm ab}: \g \to{\rm Ab}(G) $, we claim that 
 \begin{equation*} \pi_{\rm ab}(V)={\rm Ab}(G).
 \end{equation*}
 Indeed, the set $\pi_{\rm ab}(V) \subseteq \g/[\g, \g]$ is a bracket-generating subspace of an abelian Lie algebra, so it is necessarily equal to the whole Lie algebra, recall Exercise~\ref{ex_braket_gen_quotients}.
 

\begin{definition}[Abelianization norm]\label{def abelianization norm}
Let $(G, \Delta, \norm{\cdot})$ be a sub-Finsler nilpotent simply connected Lie group.
Let $\norm{\cdot}_{\rm ab}$ be the norm on the abelianization ${\rm Ab}(G) $ whose unit ball at the origin is the image in ${\rm Ab}(G) $ under $\pi_{\rm ab} $ of the unit ball at the origin for $ \norm{\cdot}$ in $\Delta_{1_G}$, that is,
\begin{equation}\label{def_norm_abelianization}
\pi_{\rm ab}\left(B_{(\Delta_{1_G}, \norm{\cdot})}(0,1)\right) = B_{({\rm Ab}(G), \norm{\cdot}_{\rm ab} )}(0,1).
\end{equation}
\end{definition}

For another expression of the norm $\norm{\cdot}_{\rm ab} $ see Exercise~\ref{ex20030242050}.
Notice that $\norm{\cdot}_{\rm ab} $ is actually the only norm that makes the projection 
$\pi_{\rm ab}|_{\Delta_{1_G}}: \Delta_{1_G}\subseteq \g \to {\rm Ab}(G) $
a submetry, in the sense of Definition~\ref{def submetry}, where on $\Delta_{1_G}$ we consider the distance induced by $\norm{\cdot}$ and on ${\rm Ab}(G) $ the distance induced by $\norm{\cdot}_{\rm ab} $.
We shall actually show that we obtained a submetry between sub-Finsler Lie groups.

\begin{proposition}[$\pi_{\rm ab}$ is a submetry]\label{prop submetry nilpotent}
 Let $G$ be a nilpotent simply connected Lie group metrized with a left-invariant sub-Finsler metric.
 Equip the abelianization ${\rm Ab}(G)$ with the norm $\norm{\cdot}_{\rm ab}$ from Definition~\ref{def abelianization norm}. 
 Then $({\rm Ab}(G) , \norm{\cdot}_{\rm ab})$ is a 
 normed vector space and 
 the abelianization map $\pi_{\rm ab}: G \rightarrow {\rm Ab}(G)$ becomes a group homomorphism
 and a submetry.
\end{proposition}

\begin{proof} 
Being the abelianization map the projection modulo the normal subgroup $[G, G]$, the map $\pi:=\pi_{\rm ab}$ is a group homomorphism; see also
Exercise~\ref{ex pi ab homo}. As we explained at the beginning of the section, the group ${\rm Ab}(G)$ is a vector space equipped with the norm from Definition~\ref{def abelianization norm}.

Regarding the fact that $\pi $ is a submetry, this is a consequence of the more general Proposition~\ref{prop_subFin_submetry_Lie}.
Since, actually, in this case, the proof is more elementary, we directly check that $\pi$ satisfies the condition in Definition~\ref{def submetry}. Fix $r>0$ and $g \in G$. 
Let $\bar B_{sF}(g,r)$ be the closed ball with respect to the sub-Finsler distance $d_{sF}$ on $G$ with center $g$ and radius $r$. 
We need to show that $\pi\left(\bar B_{sF}(g,r)\right)=\bar B_{\norm{\cdot}_{\rm ab}}(\pi(g),r)$. From the left invariance of $d_{sF}$, we may assume that $g=1_G$. By definition of $\norm{\cdot}_{\rm ab}$, we have $\norm{\pi(X)}_{\rm ab} \leq \norm{X}$ for every $X \in \Delta_{1_G}$; see Exercise~\ref{ex20030242050}. 
Integrating this inequality for the velocities along each horizontal path $\gamma$, it follows that
$$ {\norm{(\pi\circ \gamma)'}_{\rm ab}} = {\norm{\pi( \dot\gamma)}_{\rm ab}}=
{\norm{\pi( \gamma')}_{\rm ab}} \leq {\norm{ \gamma'} }
={\norm{ \dot\gamma} } ,$$
where we use that $\pi$ is a homomorphism and that the norms are left-invariant (and we also used the notation \eqref{def_gamma_prime}).
Therefore, the map $\pi$ 
does not increases
 distances and we have one inclusion: $\pi(\bar B_{sF}(1_G ,r)) \subset \bar B_{\norm{\cdot}_{\rm ab}}(0,r)$.

Regarding the other inclusion, take $Y\in {\rm Ab}(G) $ satisfying $\norm{Y}_{\rm ab}\leq r$, then, by definition of the norm $\norm{\cdot}_{\rm ab}$,
there exists $X\in \Delta_{1_G}$, such that $ {\pi(X)} = {Y} $ and $\norm{X}\leq r$.
The curve $t\in [0,1]\mapsto\exp(tX)$ is a horizontal path connecting $1_G $ and $\exp(X)$ with length at most $r$.
Hence, we have $d_{sF}(1_G , \exp(X)) \leq r$.
 Finally $\pi(\exp(X))=\pi(X)=Y$, so we have proved the opposite inclusion: $\bar B_{\norm{\cdot}_{\rm ab}}(0,r) \subset \pi(\bar B_{sF}(1_G ,r))$.
\end{proof}

\subsubsection{Lifts of curves}

\begin{proposition}\label{prop lift}
Let $G$ be a nilpotent simply connected Lie group metrized with a left-invariant sub-Finsler metric.
Equip the abelianization ${\rm Ab}(G)$ with the norm $\norm{\cdot}_{\rm ab}$ from Definition~\ref{def abelianization norm}. 
Then for every $T>0$ and every $w :[0,T] \to {\rm Ab}(G)$ measurable with
$\norm{w}_{\rm ab}= 1$ almost everywhere, 
 there exists a horizontal curve $\gamma: [0,T] \to G$ with 
 $\gamma(0)=1_G$, $ \pi_{\rm ab}(\dot \gamma) = w$, and $\norma{\dot \gamma} = 1$, for almost every in $ [0,T].$
\end{proposition}

\begin{proof} 
We get this existence result as a consequence of the fact that $\pi$ is a submetry. 
Let $\sigma: [0,T] \to {\rm Ab}(G)$ be the absolutely continuous curve with $\sigma(0)=0$ and $\sigma'=w$. Notice that since ${\rm Ab}(G)$ is a vector space, we just have $\sigma(t) := \int_0^t w(s)\dd s$, and moreover the speed is
$\norma{\dot\sigma}_{\rm ab} = \norm{w}_{\rm ab}= 1$, so the curve $\sigma$ is 1-Lipschitz, recalling \eqref{norm_less_than_length}.
The curve $\gamma$ for which we are looking will be a lift of $\sigma$ under the abelianization map $ \pi:=\pi_{\rm ab}$.

In Proposition~\ref{prop submetry nilpotent}, we saw that $\pi: G \rightarrow {\rm Ab}(G)$ is a submetry, in addition to a group homomorphism.
 Hence, via Proposition~\ref{lem:pathlifting} we can lift $\sigma$. Namely, there is a 
 $1$-Lipschitz curve $ {\gamma}:[0,T]\rightarrow G$ such that $ {\gamma}(0)=1_G$ and $\pi \circ {\gamma} = \sigma$. 
 
 First, we observe that $\pi ( {\dot\gamma} )= \dot\sigma=\sigma' =w$, where we recall that we are identifying $\pi$ with its differential, and using the notation \eqref{def_gamma_prime}. 
 Second, we notice that, on the one hand, since $\gamma$ is 1-Lipschitz, then $\norma{\dot \gamma} \leq 1$; see Exercise~\ref{ex Lipschitz gives bounded speed}. On the other hand, since $\pi$ is 1-Lipschitz (being a submetry), and $\norma{\pi ( \dot\gamma)}_{\rm ab} = \norma{\dot\sigma}_{\rm ab}=1$,
then $\norma{\dot \gamma} \geq 1$. Thus, we conclude that $\norma{\dot \gamma} = 1$.
\end{proof}

The proof of Proposition~\ref{prop lift} actually tells us that we can lift rectifiable curves from the abelianization to the group. We leave it to the reader to write the variation of the proof of the following consequence; see Exercise~\ref{ex_second_proof for prop lift}.
\begin{corollary}\label{cor prop lift}
Let $G$ be a nilpotent simply connected Lie group metrized with a left-invariant sub-Finsler metric.
Equip the abelianization ${\rm Ab}(G)$ with the norm $\norm{\cdot}_{\rm ab}$ from Definition~\ref{def abelianization norm}. 
Then for every $T>0$ and every AC curve $\sigma :[0,T] \to {\rm Ab}(G)$ with
bounded speed, 
 there exists a horizontal curve $\gamma: [0,T] \to G$ with 
 $\pi_{\rm ab}(\gamma(0))=\sigma(0)$, $ \pi_{\rm ab}\circ \gamma = \sigma$, and $\norma{\dot \gamma} = \norma{\dot \sigma}_{\rm ab}$, almost everywhere in $ [0,T].$
\end{corollary}
We picture Corollary~\ref{cor prop lift} by the following diagram:
\begin{center}
 \begin{tikzcd}
 &G\ar[d,->>]{u}{\pi_{\rm ab}} \\
 {[0,T]} \ar[r]{u}{\sigma}\ar[ru]{u}{\gamma}& {\rm Ab}(G).
 \end{tikzcd}
\end{center}

 
The abelianization map behaves naturally under homomorphisms:
 \begin{lemma}\label{induced_abelianization}
 Let $G$ and $H$ be nilpotent simply connected Lie groups. Let $\varphi: G \to H$ be a Lie group homomorphism.
 There exists a unique linear homomorphism $\varphi_{\rm ab}: {\rm Ab}(G) \to {\rm Ab}(H)$ for which the following diagram commutes:
\begin{center}
 \begin{tikzcd}
G\ar[r]{u}{\varphi}\ar[d,->>]{l}{\pi_{\rm ab}} &H\ar[d,->>]{u}{\pi_{\rm ab}} \\
 {\rm Ab}(G) \ar[r]{u}{\varphi_{\rm ab} } & {\rm Ab}(H).
 \end{tikzcd}
\end{center}
Moreover, if $\varphi_{\rm ab}$ is surjective, then so is $\varphi$.
 \end{lemma}
\begin{proof} The first statement is completely algebraic. For the diagram to commute, we need to define
$$\varphi_{\rm ab}(g[G, G]):= \varphi(g)[H, H], \qquad \forall g\in G.$$
Since $\varphi([G, G])\subseteq [H, H]$, then the map is well defined and is a homomorphism since 
$[G, G]$ and $ [H, H]$ are normal subgroups.

If in addition $\varphi_{\rm ab}$ is surjective, then we have $\h = \varphi_*(\g) + [\h, \h].$
Hence, the Lie algebra $ \varphi_*(\g)$ is bracket generating in $\h$, by Exercise~\ref{V1generates2}.
Hence $\varphi_*(\g) = \h$, and hence $\varphi(G) = H$.
 \end{proof}

 \subsection{A special sub-Finsler geometry on nilpotent groups}

In this subsection, we consider particular Carnot-Carath\'eodory metrics on nilpotent Lie groups: 
Given a nilpotent simply connected Lie group $G$ with Lie algebra $\g$, we consider left-invariant distributions $\Delta$ with the property that 
 \begin{equation}\label{special: nilpotent}\g= \Delta_1\oplus[\g, \g]. \end{equation}
This will be the case in the metric groups of our main interest: Carnot groups.

Such a $\Delta_1$ is not unique, not even up to isomorphism; see Exercise~\ref{Delta 1 not unique}.
However, every $\Delta_1$ satisfying \eqref{special: nilpotent} Lie generates the Lie algebra $\g$, recall Exercise~\ref{V1generates}. Moreover, no smaller subspace would Lie generate; see Exercise~\ref{smaller_no_generates}. Thus, these distributions may be called {\em minimal bracket-generating polarizations}.\index{minimal bracket-generating polarization}\index{bracket-generating! minimal -- polarization}


Next, we shall see how one can check via abelianizations whether a group homomorphism is Lipschitz. 
In connection to the next part \ref{lemma4homo2}.i, we also point out that there is a more general reason why the following map is Lipschitz; see Exercise~\ref{ex: Lip_for_CCgroups}.

\begin{proposition}[Criterion for Lipschitz or submetry property]\label{lemma4homo2}
Let $G$ and $H$ be nilpotent simply connected Lie groups equipped with left-invariant sub-Finsler structures: $(G, \Delta^G, \norm{\cdot})$ and $(H, \Delta^H, \norm{\cdot})$. Assume that $\mathfrak{h}= \Delta^H_1\oplus[\mathfrak{h}, \mathfrak{h}]$.
 Let $\varphi: G \to H$ be a Lie group homomorphism.
\begin{description}
 \item[\ref{lemma4homo2}.i.] If 
$\varphi_*(\Delta_1^{G}) \subseteq \Delta_1^{H}$,
then $\varphi$ is Lipschitz with respect to the sub-Finsler distances. In fact, we have $\Lip(\varphi)=\Lip( \varphi_{\rm ab})$. 
 \item[\ref{lemma4homo2}.ii.] If $\varphi_*(\Delta_1^{G}) \subseteq \Delta_1^{H}$ and $\varphi_{\rm ab}: {\rm Ab}(G) \to {\rm Ab}(H)$ is a submetry,
then $\varphi$ is a submetry.
 \end{description}
\end{proposition}
\begin{proof}
Regarding \ref{lemma4homo2}.i, the fact that $\varphi$ is Lipschitz is a more general fact; see Exercise~\ref{ex: Lip_for_CCgroups}.
We prove that, in this case, the Lipschitz constant equals $\Lip( \varphi_{\rm ab})$.
Let $\gamma$ be a $\Delta^{G}$-horizontal curve, parametrized with speed 1.
Because of the assumption, the curve $\varphi\circ \gamma$ is $\Delta^{H}$-horizontal.
Since $\pi:=\pi_{\rm ab}: G\to {\rm Ab}(G)$ is a submetry (and hence 1-Lipschitz), the curve $\pi\circ \gamma$ is 1-Lipschitz.
By Lemma~\ref{induced_abelianization}, the map $\varphi_{\rm ab}$ exists, and, being linear, it is Lipschitz.
Let $L$ be the Lipschitz constant of $\varphi_{\rm ab}$.
Then $\varphi_{\rm ab}\circ \pi\circ \gamma$ is $L$-Lipschitz.
The curve $\varphi \circ \gamma$ is a lift of 
 $\varphi_{\rm ab}\circ \pi\circ \gamma = \pi \circ\varphi \circ \gamma$, and, since $\Delta^H$ is in direct sum with $[\h, \h]$, this is the only horizontal lift; check the first commutative diagram \eqref{two diagrams for the proofs}.
 Hence, by Corollary~\ref{cor prop lift}, the curves $\varphi \circ \gamma$ and $\varphi_{\rm ab}\circ \pi\circ \gamma$ have the same speeds so, in particular, we deduce that $\varphi \circ \gamma$ is $L$-Lipschitz.
 Since the distances are length distances, we infer that $\varphi$ is $L$-Lipschitz.
 \begin{equation}\label{two diagrams for the proofs}
 \begin{tikzcd}
I\ar[r]{u}{\gamma}&G\ar[r]{u}{\varphi}\ar[d,->>]{l}{\pi_{\rm ab}} &H\ar[d,->>]{u}{\pi_{\rm ab}} \\
 &{\rm Ab}(G) \ar[r]{u}{\varphi_{\rm ab} } & {\rm Ab}(H)
 \end{tikzcd}
 \qquad 
 \begin{tikzcd}
I\ar[rd]{u}{\gamma}\ar[rrd, bend left=10]{u}{\sigma}\ar[rdd, bend right=10]{l}{\eta}\\
&G\ar[r]{u}{\varphi}\ar[d,->>]{l}{\pi_{\rm ab}} &H\ar[d,->>]{u}{\pi_{\rm ab}} \\
 &{\rm Ab}(G) \ar[r]{u}{\varphi_{\rm ab} } & {\rm Ab}(H).
 \end{tikzcd}
\end{equation}
 
 Regarding \ref{lemma4homo2}.ii, the map $\varphi$ is $1$-Lipschitz from \ref{lemma4homo2}.i.
 To deduce that $\varphi$ is a submetry, we can just show that for every absolutely continuous curve $\sigma: I\to H$ with $\norm{\dot \sigma}\equiv 1$ there is an absolutely continuous curve $\gamma: I\to G$ with $\norm{\dot \gamma}\equiv 1$ and $\varphi\circ \gamma=\sigma$.
By Corollary~\ref{cor prop lift}, since $ \Delta^H_1$ is in direct sum with $[\mathfrak{h}, \mathfrak{h}]$, the curve $\pi_{\rm ab}\circ \sigma :I\to {\rm Ab}(H)$ is such that $\norm{(\pi_{\rm ab}\circ \sigma)'} \equiv1$.
Similarly, since $\varphi_{\rm ab}$ is a submetry by assumption, there is $\eta:I\to {\rm Ab}(G)$ such that 
$\varphi_{\rm ab}\circ \eta = \pi_{\rm ab}\circ \sigma$ and $\norm{\dot\eta} \equiv1$.
Since $\pi_{\rm ab}: G\to {\rm Ab}(G) $ is a submetry, there exists $\gamma:I\to G$ with $\norm{\dot \gamma}\equiv 1$ and 
$\pi_{\rm ab} \circ \gamma=\eta$.
We conclude by checking that 
$\varphi\circ \gamma=\sigma$.
In fact, the curves $\varphi\circ \gamma$ and $\sigma$ are equal because they are the unique horizontal lift of 
$\pi_{\rm ab} \circ \sigma$, because
$$ \pi_{\rm ab} \circ \varphi\circ \gamma=
\varphi_{\rm ab}\circ \pi_{\rm ab}\circ \gamma=
\varphi_{\rm ab}\circ \eta
=\pi_{\rm ab} \circ \sigma
. $$
In other words, we showed the commutativity of the diagrams in \eqref{two diagrams for the proofs}.
\end{proof}
\subsubsection{Projection on the good first layer}
\begin{definition}[The projection $\pi_{\Delta_1}$]\label{def_proj_pi} 
Under the assumption that $\Delta_1$ is in direct sum with $[\g, \g]$, i.e., \eqref{special: nilpotent}, 
let ${\rm proj}_{\Delta_1}:\g\to \Delta_1$ be the projection onto $\Delta_1$ with kernel $[\g, \g] $ and define
\begin{eqnarray}
\pi_{\Delta_1}: &G\to&\Delta_1\nonumber\\
 &g\mapsto&\pi_{\Delta_1}(p):={\rm proj}_{\Delta_1}(\exp^{-1}(g)).\label{projection_pi1}
\end{eqnarray}
\end{definition}

\begin{remark}\label{remark same maps} We stress that, under the assumption \eqref{special: nilpotent}, the quotient space $\g/[\g, \g]$ can be identified with $\Delta_1$, as vector spaces, and the projection $\g \twoheadrightarrow\g/[\g, \g]$ can be identified with $\g \twoheadrightarrow\Delta_1$, modulo $[\g, \g]$. Hence, the map $\pi_{\Delta_1}: G \twoheadrightarrow\Delta_1$ is nothing else than the abelianization map $\pi_{\rm ab}: G \twoheadrightarrow{\rm Ab}(G)$, in the case \eqref{special: nilpotent}.
\end{remark}


\begin{lemma}\label{lemma:pi1} Let $G$ be a nilpotent simply connected Lie group and let $\Delta_1$ be a polarization on $G$ such that 
\eqref{special: nilpotent} holds. For the map $\pi_{\Delta_1}$ of Definition~\ref{def_proj_pi}, we have:
 \begin{description}
\item[\ref{lemma:pi1}.i.] The map $\pi_{\Delta_1}: G\to (\Delta_1,+)$ is a Lie group homomorphism.
\item[\ref{lemma:pi1}.ii.] The differential of $\pi_{\Delta_1}$ is the identity when restricted to $\Delta_1$:
$$ (\pi_{\Delta_1} )_* |_{\Delta_1}=\id_{\Delta_1}.$$
\end{description}
\end{lemma}
\proof[Proof of (i).] 
Because of Remark~\ref{remark same maps}, this is just a restatement of Proposition~\ref{prop submetry nilpotent}. We spell out another proof that avoids identifications.
We write $\pi$ for $\pi_{\Delta_1}$ and ${\rm proj}$ for ${\rm proj}_{\Delta_1}$.
By Theorem~\ref{CG-1.2.1}, since $G$ is simply connected and nilpotent, for all $p$ and $q\in G$, exist $X$ and $Y\in \g$ such that $\exp(X)=p$ and $\exp(Y)=q$. On the one hand, by BCH formula and assumption \eqref{special: nilpotent}, we have
$$\pi (p\cdot q)\stackrel{\eqref{projection_pi1}}{=}{\rm proj}(\exp^{-1}(pq))
={\rm proj}\left(\exp^{-1}(\exp(X)\exp(Y))\right)$$
$$\stackrel{\eqref{expansion of BCH}}{=}{\rm proj}\left(X+Y+\dfrac{1}{2}[X,Y]+\ldots\right)
\stackrel{\eqref{special: nilpotent}}{=}{\rm proj}(X+Y).$$
On the other hand, we have
$$\pi (p)+\pi ( q)={\rm proj}(\exp^{-1}(p))+{\rm proj}(\exp^{-1}(q))=
{\rm proj}(\exp^{-1}(p)+ \exp^{-1}(q))={\rm proj}(X+Y).$$
\proof[Proof of (ii).] By Proposition~\ref{exp:diffeo}, we have
\begin{eqnarray*}
 \dd\pi|_{\Delta_1}&\stackrel{\rm def}{=}&(\dd ({\rm proj}\circ \exp^{-1}))_{1_G} |_{\Delta_1}\\
 &=&(\dd \,{\rm proj})_0 |_{\Delta_1}\\
 &=&(\dd \id)_0 |_{\Delta_1} =\id_{\Delta_1}.
 \end{eqnarray*}
\qed

Here is another restatement of Proposition~\ref{prop submetry nilpotent}: 
 \begin{corollary}\label{pi1_submetry}
 Let $G$ be a nilpotent simply connected Lie group polarized by $\Delta$ so that $\Delta_1\oplus[\g, \g] = \g$. Fix a norm $\norm{\cdot}$ on $\Delta_1$. The map 
 $\pi_{\Delta_1}:(G, \dcc )\to (\Delta_1, \norm{\cdot})$ is a submetry, where $\dcc $ is the metric on the sub-Finsler Lie group $(G, \Delta, \norm{\cdot})$.
 \end{corollary}

The map $\pi_{\Delta_1}$ from \eqref{projection_pi1} is useful since it gives a second link between the tangents of horizontal curves and vectors in $\Delta_1$. Recall the notions of development and multiplicative integral from Definitions~\ref{def development} and \ref{def multiplicative integral}, respectively.

\begin{proposition}[Development as projection on abelianization]\label{prop development as projection}
Let $G$ be a nilpotent simply connected Lie group polarized by $\Delta$ so that $\Delta_1\oplus[\g, \g] = \g$. Fix $T>0$.
If $\gamma:[0,T]\to G $ is a $\Delta$-horizontal curve, then the development of $\gamma$ is $ \pi_{\Delta_1}\circ \gamma.$
Vice versa, if $\sigma:[0, T]\to \Delta_1$ is an absolutely continuous curve, then the multiplicative integral of $\sigma$ is the only $\Delta$-horizontal curve $\gamma$ with $\sigma:= \pi_{\Delta_1}\circ \gamma$.
\end{proposition}


\begin{proof} 
We need to prove the formula 
\begin{eqnarray}\label{gamma'formula}
 \gamma'(t)\stackrel{\rm def}{=}(L_{\gamma(t)})_*^{-1}\dot\gamma(t)=\dfrac{\dd}{\dd t}\left(\pi_{\Delta_1}\circ \gamma \right)(t), \qquad \forall t\in [0,T],
\end{eqnarray}
as elements of $\Delta_1$.
Using Lemma~\ref{lemma:pi1}, and that $\pi_{\Delta_1}(1_G)=0$, we get
\begin{eqnarray*}
 \dfrac{\dd}{\dd t}\left(\pi_{\Delta_1}\circ \gamma \right)(t) &=&
\lim_{h\to 0}\dfrac{\pi_{\Delta_1}(\gamma(t+h))-\pi_{\Delta_1}(\gamma(t)) }{h}\\
&\stackrel{\rm\ref{lemma:pi1}.i}{=}&
\lim_{h\to 0}\dfrac{\pi_{\Delta_1}(\gamma(t)^{-1} \gamma(t+h)) }{h}\\
&=&
\lim_{h\to 0}\dfrac{\pi_{\Delta_1}(L_{\gamma(t)}^{-1} \gamma(t+h)) }{h}\\
&=&
\lim_{h\to 0}\dfrac{\pi_{\Delta_1}(L_{\gamma(t)}^{-1} \gamma(t+h)) -\pi_{\Delta_1}(L_{\gamma(t)}^{-1} \gamma(t)) }{h}\\
&=& \left.\dfrac{\dd}{\dd h}\left((\pi_{\Delta_1}\circ L_{\gamma(t)}^{-1}\circ \gamma \right)(t+h)\right|_{h=0}\\
 &=&(\pi_{\Delta_1})_*\circ( L_{\gamma(t)}^{-1} )_* \dot\gamma(t)\\
& \stackrel{\rm \ref{lemma:pi1}.ii}{=}&\id ( \gamma'(t))= \gamma'(t).
\end{eqnarray*}
The vice versa is obvious since multiplicative integrals are the inverse operation of developments.
\end{proof}

\begin{corollary}\label{length_for_projection}
Let $G$ be a nilpotent simply connected Lie group polarized by $\Delta$ so that $\Delta_1\oplus[\g, \g] = \g$. Fix a norm $\norm{\cdot}$ on $\Delta_1$.
The length of the horizontal curves equals the length of their projections:
$$\Length (\gamma)=\Length (\pi_{\Delta_1}\circ\gamma), \qquad \forall \gamma \text{ horizontal curve in } (G, \Delta) ,$$
where the first length is with respect to the sub-Finsler metric and the second one is in the normed space $(\Delta_1, \norm{\cdot})$.
\end{corollary}
 \begin{proof}
 By Proposition~\ref{prop development as projection} (see \eqref{gamma'formula}), one has
\begin{eqnarray*}
 \Length (\pi \circ\gamma)&= &\int_0^1 \norm{\dfrac{\dd}{\dd t}\left(\pi\circ \gamma \right)(t)} \dd t\\
&= &\int_0^1 \norm{ \gamma'(t) } \dd t \\
&=&\int_0^1 \norm{(L_{\gamma(t)})_*^{-1}\dot\gamma(t)} \dd t \\
&=&\int_0^1 \norm{\dot\gamma(t)} \dd t \\
&=& \Length (\gamma).
\end{eqnarray*}
\end{proof}



\subsubsection{Horizontal lines as geodesics}
 \begin{definition} 
Let $(G, \Delta, \| \cdot \|)$ be a sub-Finsler Lie group. 
 Let $X\in \Delta_1$. The curve $ t \mapsto\exp(t X)$ is the one-parameter subgroup tangent to the vector $X$, and it is called the {\em horizontal line} in the direction of $X$.
 \end{definition}
 
 On each sub-Finsler Lie group $(G, \Delta, \| \cdot \|)$,
 the curve $ t \mapsto\gamma(t):=\exp(t X)$ is horizontal with respect to $\Delta$, because, recalling Corollary~\ref{Warner3.31},
we have $\dot\gamma(t)=X_{\gamma(t)}
 = \dd L_{\gamma(t)} X
 \in\Delta_{\gamma(t)}.$ 
We claim that the length of $t \mapsto \gamma(t)$, for $t\in[0,T]$, with respect to the sub-Finsler structure of $(G, V, \norm{\cdot})$ is $T\norm{X}$. Indeed, we have
 \begin{eqnarray*}
\Length (\gamma)&=& \int_0^T\norm{\dot\gamma(t)} \dd t\\
 &=& \int_0^T\norm{X_{\gamma(t)}} \dd t\\
 &=& \int_0^T\norm{\left(L_{\gamma(t)}\right)_*X_{1}} \dd t\\
 &=& \int_0^T\norm{X} \dd t\\
 &=&T \norm{X},
 \end{eqnarray*}
 where we used that both $X$ and the norm are left-invariant.
 In fact, we get the formula
 \begin{equation}\label{length:horiz:line}
 \Length \left(\exp(tX)|_{t\in[a,b]}\right )=|b-a|\cdot \norm{X} \qquad \forall a,b\in \R.
 \end{equation}

 In Lie groups endowed with left-invariant Riemannian metrics, 
 one-parameter subgroups may not be geodesics.
 For instance, in the semidirect product $\R\rtimes \R_+$, which is the hyperbolic half-plane, then the distance between the identity element
 $(0,1)$ and $(t,1)$ equals $ 2 {\rm arcsinh}(t/2)$, for $t\in \R$; see Exercise~\ref{ex distance on half space}.
 

Regarding nilpotent groups, in the Riemannian Heisenberg group with orthonormal frame $X, Y, Z=[X, Y]$, the one-parameter subgroup in the direction $Z$ is a Riemannian geodesic, but not globally length-minimizing. Whereas, the one-parameter subgroup in the direction $X+Z$ is not even a Riemannian geodesic, recall Proposition~\ref{Riemannian_geodesic_OPS}.

\begin{proposition}\label{OPS geodesic}
Consider a nilpotent simply connected Lie group $G$ endowed with a left-invariant sub-Finsler distance with respect to some distribution $\Delta$ such that 
 $\Delta_1\oplus[\g, \g]=\g$. 
Then, one-parameter subgroups of horizontal vectors are length-minimizing.
\end{proposition} 
\proof 
Take $X\in \Delta_1$. We shall prove that the curve $t\in \R \mapsto \exp(tX)$ is an homothetic embedding:
$$\dcc (\exp(aX), \exp(bX))= |b-a|\cdot \norm{X}, \qquad \forall a,b\in \R.$$
Recall from Corollary~\ref{pi1_submetry} that the map $\pi_{\Delta_1}$ is a submetry. 
We have the following bounds: for all real numbers $a<b$ 
\begin{eqnarray*}\dcc (\exp(aX), \exp(bX))&\leq&\Length \left(\exp(tX)|_{t\in[a,b]}\right )
\\&\stackrel{ \eqref{length:horiz:line}}{=}&(b-a) \norm{X}\\&=&\norm{(b-a) X}
\\&=&\norm{\pi_{\Delta_1}((b-a) X)}
\\&\leq& \dcc (\exp(aX), \exp(bX)),
\end{eqnarray*}
where in the last inequality we used the 1-Lipschitz property of $\pi_{\Delta_1}$.
\qed

\subsection{Lift of regularity} 

With the notion of submetry, we can transport several regularity results about geodesics from a sub-Finsler Lie group to its quotients. Said with a different perspective, if a nilpotent sub-Finsler Lie group may admit pathological geodesics, then the 
free-nilpotent Lie group with the same rank and step admits the same pathologies.
Next, we will discuss the property of joining points by smooth geodesics. For other regularity results, we point out to Exercises~\ref{prop: same regularity0} and~\ref{ex:same regularity2}, and Proposition~\ref{prop: same regularity2}.

\begin{proposition}
Let $G$ and $H$ be nilpotent simply connected Lie groups equipped with left-invariant sub-Finsler structures: $(G, \Delta^G, \norm{\cdot})$ and $(H, \Delta^H, \norm{\cdot})$ 
such that $\mathfrak{h}= \Delta^H_1\oplus[\mathfrak{h}, \mathfrak{h}]$.
 Let $\varphi: G \to H$ be a Lie group homomorphism.
 Assume that each pair of points in $G$ can be joined by a smooth geodesic.
If $$\varphi_*|_{V_1^{(G)}} :\Delta_1^G \to \Delta_1^H$$
is a submetry of normed spaces, then each pair of points in $H$ can be joined by a smooth geodesic.
\end{proposition}

\proof
The proof is an immediate consequence of Proposition~\ref{lemma4homo2}.ii. and Proposition~\ref{prop: same regularity2}.
\qed

\section{Isometrically homogeneous spaces with dilations (second part)}\label{sec: self-similar2}\index{self-similar}
This section is the continuation of Section~\ref{sec: self-similar}, where we already considered isometrically homogeneous spaces with dilations and proved that they, at least, have a structure of Lie coset spaces. 
We will now see that, in fact, they have the structure of nilpotent positively graduable Lie groups.
Such a group structure comes from the nilradical of the isometry group. See Exercise~\ref{ex: nilradical} for the definition of nilradical.
	
	\begin{theorem}\label{teo05171840}
	Let $(M,d)$ be a metric space such that
	\begin{enumerate}
	\item 	$(M,d)$ is locally compact;
	\item 	$(M,d)$ is locally connected;
	\item 	the isometry group $\Isom(M,d)$ of $(M,d)$ acts transitively on $M$;
	\item 	there is a dilation on $M$ of factor in $(0,+\infty)\setminus\{1\}$.
	\end{enumerate}
	Then there is a unique nilpotent Lie group $N$ with a left-invariant distance $d_N$ on $N$ such that $(N,d_N)$ is isometric to $(M,d)$
and every dilation of $(N,d_N)$ is an affine map of $N$.
	Moreover, this Lie group $N$ is 
	positively graded
	and
	isomorphic to the nilradical of $\Isom(M,d)$. 
\end{theorem}

\begin{proof}
Recall what we already saw in the proof of Theorem~\ref{teo05171840baby}:
the identity component $G:=\Isom(M,d)^\circ$ of the isometry group $\Isom(M,d)$ is a connected Lie group.
We fixed a point $o\in M$ and a dilation $\lambda: M\to M$ of factor $\lambda\in(0,1)$ fixing $o$.
We defined $S:={\rm Stab}(o)$ to be the stabilizer of $G$ at $o$, which is a compact Lie subgroup of $G$. 
Via the orbit map $\pi$, the Lie coset space $G/S$ is homeomorphic to $M$.


	{\color{black}
	Define for $f\in \Isom(M,d)$
	\[
	Tf := \lambda\circ f\circ\lambda^{-1}.
	\]
We claim that the map $T$ is valued into $\Isom(M,d)$, that $T:\Isom(M,d)\to \Isom(M,d)$ is a Lie group automorphism, and $T(S)=S$.
Indeed, one easily sees that if $f$ is an isometry of $M$, then $Tf: M\to M$ is an isometry as well, and that $T(f\circ g)=Tf\circ Tg$ for all $f,g\in \Isom(M,d)$.
Moreover, $T$ is continuous: if $\{f_k\}_{k\in\N}$ is a sequence converging in $\Isom(M,d)$ to $f$, i.e., $f_k\to f$ uniformly on compact sets, then $Tf_k$ converges to $Tf$ uniformly on compact sets as $k\to\infty$, because for every compact $K\subseteq M$
	\[
	\sup\left\{ d(Tf_k(p),Tf(p)) : p\in K
	 \right\}
	= \lambda \sup\left\{ d(f_k(p),f(p)) : p\in \lambda^{-1}(K)
	\right\}.
	\]
Since $T$ is a continuous automorphism, then $T$ is a Lie group automorphism; see Theorem~\ref{thm: continuous: smooth}.
	Finally, since $\lambda(o)=o$, then $T(f)\in S$, for all $f\in S$.
	
	}
	
	
	
	Let $\g :=\Lie(\Isom(M,d))$ be the Lie algebra of $\Isom(M,d)$ and $\mathfrak s:=\Lie (S)$ the Lie algebra of the stabilizer $S$.
	We consider the bilinear form 
	$$B:\g \times\g \to \R, \qquad B(v,w) := {\rm Trace}(\ad(v)\ad(w)), $$
	 which is the bilinear symmetric form called the Killing form; see Exercise~\ref{ex_def_Killing}.
	
We claim that $B$ restricted to $\mathfrak s$ is negative definite.
For doing this, fix $v\in \mathfrak s$, for which we want to show that $B(v,v)<0$ unless $v=0$.
	Since $\Ad_S$ is compact in $\GL(\g)$, then
 the eigenvalues of $\ad_v$ are purely imaginary, see Exercise~\ref{ex8may1429}. 
So the eigenvalues of $\ad_v^2$ are real and nonpositive. 
Consequently, because the trace is the sum of the eigenvalues, we have that $B(v,v)= {\rm Trace}(\ad_v^2) \leq0$.
	It is left to prove that,
	\begin{equation}\label{B_positive}
	 \text{if } B(v,v)=0, \quad \text{ then }v=0. 
	 \end{equation}
	 Assume that ${\rm Trace}(\ad_v^2) =B(v,v)=0$, then since again the eigenvalues of $\ad_v^2$ are nonpositive, we deduce that then all the eigenvalues of $\ad_v^2$, and hence the ones of $\ad_v$, are zero.
	But the only diagonalizable transformation with all zero eigenvalues is the zero transformation (recall Exercise~\ref{ex8may1429}). We deduced that $\ad_v=0$.
	Therefore, from Formula~\ref{Formula: Ad: ad} we have
	$\Ad_{\exp(tv)} = e^{t\ad_v}=\id_{\g }$ and $C_{\exp(tv)} =\id_G$, for all $t\in \R$, recalling that $G $ is the identity component of $\Isom(M,d)$. Therefore, since $G$ is connected, we have that $\exp(tv)$ is both in the center $Z(G)$ of $G$ and also in $S$ for all $t$.
	We claim that $Z(G)\cap S$ is trivial. Indeed, take $f\in Z(G)\cap S$ and $p\in M$. Since the action of $G$ on the connected space $M$ is transitive by Proposition~\ref{identity_component_trans}, there is $g\in G$ such that $g(o)=p$. We deduce that $f(p)=f(g(o))=g(f(o))=g(o)=p$. Since this holds for every $p\in M$, we deduce that $f=\Id_M$.
	The claim about the triviality of $Z(G)\cap S$ is proven. Consequently, $\exp(tv)=\Id_M$ for all $t\in \R$. Hence $v=0$. We proved \eqref{B_positive}, and, therefore, that $B$ is negative definite on $\mathfrak s$.

We next consider the map 
$$\phi:\g \to \mathfrak s^*, \qquad\phi(v):=B(v, \cdot),$$ and define the orthogonal of $\mathfrak s$ with respect to the Killing form:
	\begin{equation}\label{def_n_proof}
\mathfrak n :=\Ker\phi \stackrel{\rm def}{=} \{v\in\g : B(v,w)=0, \forall w\in\mathfrak s \}.
	\end{equation}
We claim that $\g =\mathfrak n\oplus\mathfrak s$, as linear spaces. 
We begin by noting that $\mathfrak n\cap\mathfrak s=\{0\}$ by \eqref{B_positive}. 
To deduce that $\g =\mathfrak n\oplus\mathfrak s$, we stress that $\mathfrak n $ is the kernel of the map $\phi$ and that $\phi$ restricted to $\mathfrak s$ is an isomorphism because $B$ is negatively definite on $\mathfrak s$.

	 {\color{black} We next claim that $T_*(\mathfrak n)=\mathfrak n$. Indeed, since $T_*:\g \to\g $ is a Lie algebra automorphism, we have $B(T_*v,T_*w)=B(v,w)$ for all $v,w\in\g $.
	Moreover, $T_*(\mathfrak s)=\mathfrak s$.
	Therefore, $T_*(\mathfrak n)=\mathfrak n$.
	}

 	{\color{black}
For $t\in \R$, let $V_t: = V_t(e^{-1}, T_*)$ be the Lie algebra grading of $\g$ induced by $T_*$ as in Proposition~\ref{prop08220937}.
Since both $\mathfrak s$ and $\mathfrak n$ are invariant subspaces for $T_*$, we have $V_t = (V_t\cap\mathfrak s) \oplus (V_t\cap \mathfrak n)$ for each $t\in\R$.

We claim that $\mathfrak s\subset V_0$.
Indeed, we have seen that $T_*$ is preserving $B$ and that $-B$ is a scalar product on $\mathfrak s$.
Therefore, the linear transformation $T_*|_{\mathfrak s}$ is orthogonal with respect to $-B$ and thus every eigenvalue $\alpha$ of $T_*|_{\mathfrak s}$ has norm 1. Thus, every element in ${\mathfrak s}$ has degree $t$ with $t\stackrel{\rm def}{=}-\log(|\alpha|)=- \log 1=0$; see \eqref{eq19Aug20240853}.

Next, we claim that $\mathfrak n\subset V_{>0}:= \bigoplus_{t>0} V_t$.
Arguing by contradiction, suppose that there is an eigenvalue $\alpha\in\C$ such that $V_t\cap\mathfrak n \neq\{0\}$ with $t:=-\log(|\alpha|)<0$, i.e., with $|\alpha|\ge 1$.
Let $\mathscr B$ be a basis of $\mathfrak g$ that represents $T_*$ in Jordan form and let $|\cdot|$ the Euclidean norm on $\mathfrak g$ so that $\mathscr B$ is an orthonormal basis.
Let $U$ be an open ball in $\mathfrak g$ centered at $0$ such that
the restriction of $\exp:\mathfrak g\to G$ to $2U$ is a diffeomorphism between $2U$ and its image and such that $\exp(\mathfrak s\cap U) = S \cap \exp(U)$.
In a first case, if $\alpha\in\R$, then there is a nonzero $a\in \mathfrak n\cap U$ such that $T_*a=\alpha a$, while, in the second case, if $\alpha\notin\R$ there there are nonzero $a,b\in \mathfrak n\cap U$ such that $|a|=|b|$,
 $T_*a = \Re(\alpha) a + \imaginary(\alpha) b$ and $T_*b = -\imaginary(\alpha) a + \Re(\alpha) b$.
Notice that, in both cases, we have that $|T_*^k a| = |\alpha|^k |a|$ for all $k\in\N$.

In a first case, if $|\alpha|=1$, then there is a diverging sequence $k_j $ of natural numberrs with $\lim_{j\to\infty}T_*^{k_j}a = a \in U$.
Hence
\[
\exp(a)(o)
= \lim_{j\to\infty} \exp(T_*^{k_j}a)(o)
= \lim_{j\to\infty} \delta^{k_j}\circ \exp(a) \circ \delta^{-k_j} ( o)
= \lim_{j\to\infty} \delta^{k_j}(\exp(a)(o))
= o,
\]
that is $\exp(a)\in S\cap \exp(U)$.
Thanks to our conditions on $U$, we obtain $a\in\mathfrak s\cap\mathfrak n = \{0\}$, in a contradiction since $a\neq 0$.

In the second case, if $|\alpha|>1$, then $|T_*^{-k}a| = |\alpha|^{-k}|a|$ is a vanishing sequence as $k\to\infty$.
Therefore,
\[
o = \Id(o)
= \lim_{k\to\infty} \exp(T_*^{-k}a)(o)
= \lim_{k\to\infty} \delta^{-k}\exp(a)\delta^k(o)
= \lim_{k\to\infty} \delta^{-k}\exp(a)(o),
\]
where $ \lim_{k\to\infty} d(o, \delta^{-k}\exp(a)(o)) = +\infty$, since $\delta^{-1}$ is an expanding dilation.
In both cases, we have got a contradiction.
We have thus proven our claim that $\mathfrak n\subset V_{>0}$.

We conclude that the grading is nonnegative, $\mathfrak s = V_0$, and $\mathfrak n = V_{>0}$. 
Consequently, the linear space $\mathfrak n$ is a nilpotent ideal in $\mathfrak g$; recall Exercises~\ref{ex positive gradings give nilpotency} and \ref{V>0}.

Therefore, by the definition of nilradical, we have $\mathfrak n\subset\mathfrak{nil}(\mathfrak g)$.
We stress that the Killing form is negative definite on $\mathfrak s$, while it is zero on $\mathfrak{nil}(\mathfrak g)$; see Exercise~\ref{ex_Killing_on_nilrad}.
We infer $\mathfrak s\cap \mathfrak{nil}(\mathfrak g) = \{0\}$ and hence $\mathfrak n = \mathfrak{nil}(\mathfrak g)$. 
}

Now that we know that $\mathfrak n$ is a Lie algebra, we can consider $N<\Isom(M,d)$ to be the connected Lie subgroup with Lie algebra $\mathfrak n$, by Theorem~\ref{teo1145bis}.
We claim that the orbit map $\pi|_N: N\to M$ restricted to $N$ is a homeomorphism. 
Indeed, we firstly recall that the identity component $G$ of the isometry group acts transitively on the connected space $M$; see Proposition~\ref{identity_component_trans}.
Secondly, 	because $\g =\mathfrak n\oplus\mathfrak s$, the group $N$ acts transitively on $M$, i.e., $\pi(N)=M$.
Thirdly, $N\cap S$ is discrete, because $\mathfrak n\cap\mathfrak s$ is trivial.
So, $\pi|_N: N\to M$ is a covering map. 
Since the manifold $M$ is simply connected by Exercise~\ref{Ex_local_to_global_dilations}, the map $\pi|_N$ is actually a homeomorphism.
Therefore, $N$ is simply connected, and it is a graded nilpotent Lie group.
	
Finally, we make $N$ into a metric Lie group by pulling back the distance $d$ from $M$ to $N$ via $\pi$:
	\[
	d_N(f,g) := d(\pi(f), \pi(g)) = d(f(o),g(o)), \qquad\forall f,g\in N.
	\]
	 {\color{black}
We stress that for every $f\in N$ by definitions we have
	 $\pi(Tf) =Tf(o)= (\lambda f \lambda^{-1})(o)
	 =(\lambda f)(o) = \lambda \pi(f)$.
We conclude that $\pi|_N$ is an isometry between $(N,d_N)$ and $(M, d)$ that relates the dilation $\lambda$ on $M$ with the automorphism $T$ on $N$. }


{\color{black}

For the uniqueness, suppose that $M$ has another Lie group structure so that $d$ is a left-invariant admissible distance, and $\delta$ is a Lie automorphism.
Then the set $M^L$ of left translations is a subgroup of $G$ whose Lie algebra $\mathfrak m$ is complementary to $\mathfrak s$, i.e., $\mathfrak g = \mathfrak s \oplus \mathfrak m$,
and invariant under $T_*$, i.e., $T_* \mathfrak m = \mathfrak m$.
Take $m\in \mathfrak m$ and decompose it as $m= m_{\mathfrak s} +m_{\mathfrak n} $ with $m_{\mathfrak s}\in \mathfrak s $ and $m_{\mathfrak n}\in \mathfrak n$. Then, since $T_*$ is contractive on $\mathfrak n$, up to subsequence, we get
$$\mathfrak m\ni \lim_{k\to\infty} T_*^k m = \lim_{k\to\infty} T_*^k( m_{\mathfrak s} +m_{\mathfrak n}) = \lim_{k\to\infty} T_*^k m_{\mathfrak s} \in \mathfrak s.$$
Thus $\lim_{k\to\infty} T_*^k m_{\mathfrak s}=0$ and, since $T_* $ restricted to $\mathfrak s$ preserves the norm given by $-B$, we deduce that $m_{\mathfrak s}=0$.
We conclude that $\mathfrak m = \mathfrak n={\rm nil}(\mathfrak g)$, and thus the original group structure is equal to the one given by $N$.}
\end{proof}

\section{Affinity of isometries of nilpotent Lie groups}\label{sec:affinity of isometries}
We say that a map between groups is \emph{affine} if it is the composition of a left translation and a group homomorphism.
Equivalently, it is the composition of a group homomorphism and a left translation.
For nilpotent groups, we have the following result on their isometries, which we will not prove but will discuss only partially.

\begin{theorem}[{\cite[Theorem 1.2]{Kivioja_LeDonne_isom_nilpotent}}, after Wolf \cite{Wolf62, Wilson}]\label{thm_Ville} 
Isometries between nilpotent connected metric Lie groups are affine.
\end{theorem}

Here are some consequences:
\begin{description}
\item[{[\ref{thm_Ville}.i]}] Two isometric nilpotent connected metric Lie groups are isomorphic;
\item[{[\ref{thm_Ville}.ii]}] 
Given a connected metric Lie group \( N \), its isometry group $\Isom(N)$, which always is a Lie group, is a semidirect product if \( N \) is nilpotent.
Namely, 
\begin{equation*}\label{semidirect}
\Isom (N) = N \rtimes \operatorname{AutIsom}(N),
\end{equation*}
where $N$ is seen inside $ \Isom (N)$ as left translations and \( \operatorname{AutIsom}(N)\) denotes the group of automorphisms of \( N \) that are isometries. 
\end{description}
Moreover, with the above notation, we have 
\begin{description}
\item[{[\ref{thm_Ville}.iii]}] 
The group $N$ is a maximal connected nilpotent subgroup of $ \Isom (N)$, and the Lie algebra of $N$ is the nilradical of the Lie algebra of $ \Isom (N)$.
\end{description}

We will not prove Theorem~\ref{thm_Ville} in this text. Historically, this has been proved in \cite{Kivioja_LeDonne_isom_nilpotent} using Theorem~\ref{prop:isometries_are_riemannian} and reducing it to the case of Riemannian manifolds where Wolf has proved the result in \cite{Wolf62}.
In Wolf's proof, there is an unclear step, which has been later clarified by Wilson in \cite{Wilson}. Nowadays, we have more direct proofs of Theorem~\ref{thm_Ville}. See, for example, \cite{CKLNO}.

Given a metric group $(M,d)$, we denote by $M^L$ the group of the left translations inside the isometry group $\Isom(M,d)$ and by 
$\mathrm{Stab}_1 (\Isom(M,d)) $ the stabilizer of the identity element $1=1_M$.\index{stabilizer}
We denote by $\mathrm{Aff}(M)$ the group of affine maps from $M$ to $M$ and by $\mathrm{Aut}(M) $ the group of automorphisms of $M$.
Then the following properties are equivalent:
\begin{enumerate}
\renewcommand{\theenumi}{\alph{enumi}}
\item \( M^L \lhd \Isom(M,d) \); 
\item
$\Isom(M,d)< \mathrm{Aff}(M)$; 
\item
 \( \mathrm{Stab}_1 (\Isom(M,d)) < \mathrm{Aut}(M) \); 
\item
 \( \Isom(M,d) = M^L \rtimes \mathrm{Stab}_1 (\Isom(M,d)) \); 
 \item 
 \( \Isom(M,d) = M^L \rtimes (\Isom(M,d)\cap \mathrm{Aut}(M))\). 
\end{enumerate}

All of the properties hold when $(M,d)$ is nilpotent and connected.

\section{Guivarc'h seminorms on nilpotent Lie groups}\label{sec:seminorms}

On nilpotent simply connected Lie groups, there are special (coarse) distance functions that are more linked to the Lie algebra nilpotency properties. 
They are called homogeneous quasi-norms since they behave well with respect to dilations.
The most important examples of homogeneous quasi-norm are the following:
\begin{enumerate}
 \item Guivarc'h seminorms, which are present in every nilpotent simply connected Lie group; see Lemma~\ref{lem_Guivarch}.
 \item The distance $d(1, \cdot)$ from the identity element in every Carnot group; see Proposition~\ref{dilations_on_Carnot2024}. 
 \item The asymptotic distance $d_\infty(1_G, \cdot)$ from the identity element in every simply connected nilpotent sub-Finsler Lie group $G$, after the identification of $G$ with its associated Carnot group $G_\infty$,
 see Corollary~\ref{thm_guiv_bis}.
\end{enumerate}

\begin{definition}\label{def_homogeneous_quasi_norm}\index{homogeneous! -- quasi-norm}\index{quasi-norm! homogeneous --}
Let $G$ be a nilpotent simply connected Lie group. Fix on its Lie algebra $\g$ a linear grading $(V_a)_{a\in \R }$, as in Definition~\ref{def_linear_grading} and let $\delta_\lambda$ be the inhomogeneous dilation on $\g$ of factor $\lambda$ relative to the grading, as in Definition~\ref{Dilations:algebras}.
 An {\em homogeneous quasi-norm} is a function $|\cdot|: G\to \R_{\geq0}$ such that 
 \begin{description}
\item[\ref{def_homogeneous_quasi_norm}.i.] it a continuous function 
\item[\ref{def_homogeneous_quasi_norm}.ii.] 
 $|x| = 0$ if and only if $x=1_G$, and 
 \item[\ref{def_homogeneous_quasi_norm}.iii.] $|\delta_\lambda(x)| = \lambda|x|$ for every $x\in G$ and every $\lambda \in \mathbb{R}_+$. 
\end{description}
\end{definition}
In particular, for each homogeneous quasi-norm, in exponential coordinates, we have that 
$$ | \lambda x| = \lambda^{1/a} |x|, \qquad \forall a\in \R,\forall x\in V_a.$$

As previously mentioned, important examples were constructed by Guivarc'h in the presence of compatible linear gradings as in Definition~\ref{def_compatible_linear_grading}. 
We will present the proof taken from \cite[Lemma 2.5]{Breuillard_MR3267520}, which, however, is an alternative exposition of the original one \cite[Lemma II.1]{MR0369608}. 
\begin{lemma}[Guivarc'h]\label{lem_Guivarch}\index{Guivarc'h! -- quas-inorm}\index{quasi-norm! Guivarc'h --}
 Let $G$ be a nilpotent simply connected Lie group. On its Lie algebra $\g$ fix a compatible linear grading $\g= V_1 \oplus \ldots \oplus V_s$, 
 and denote by $((x)_1, \ldots,(x)_s)$ the decomposition of each $x\in \g$ with respect to this grading.
 Then for every $\varepsilon > 0$ there exists a norm $\norma{\cdot}$ on $\g$ such that the function 
 \begin{align}\label{eq_Guivarch_norm}
 g\in G\longmapsto |g| \coloneqq \max_{j\in \{1, \ldots,s\}} \left\| (\log g)_j \right\|^{1/j}
 \end{align}
 is a homogeneous quasi-norm that satisfies
 \begin{align}\label{eq_lem_guiv}
 |g\cdot h| \leq |g|+|h|+\varepsilon, \quad \forall g,h \in G.
 \end{align}
\end{lemma}

\begin{proof}
 We fix some scalar product on $\g$ that makes the compatible linear grading orthogonal. We shall not change the induced norm $\norm{\cdot}$ on $V_1$, while on every other $V_j$ we will replace the norm with the norm $\lambda _{j}\norm{\cdot}$ for some $\lambda _{j}\in \R_+$ so that \eqref{eq_lem_guiv} would hold. 
The $\lambda _{j}$'s will be taken to be smaller and smaller as $j$ increases. 

We work in exponential coordinates.
We set $|x|_{\lambda}:=\max_{j}\left\| \lambda _{j}(x)_j\right\|^{1/j}$ for $\lambda\in \R^s_+$. 
We want that for every index $j\in \{1, \ldots, s\},$
\begin{equation}
\lambda _{j}\left\| (x\star y)_j\right\| \leq \left( |x|_{\lambda
}+|y|_{\lambda }+\varepsilon \right) ^{j} , \qquad \forall x,y \in \g. \label{zooo}
\end{equation}
By BCH Formula, we have $ (x\star y)_j= (x)_j+ (y)_j+P_{j}(x,y)$, where $P_{j}$ is a polynomial map into $V_j$ depending only on $(x)_i$ and $(y)_i$ with $i\in \{1, \ldots, j-1\}$ such that
\begin{equation*}
\left\| P_{j}(x,y)\right\| \leq C_{j}\cdot \sum_{l,m\geq 1,l+m\leq
j}M_{j-1}(x)^{l}M_{j-1}(y)^{m} ,
\end{equation*}
where $M_{k}(x):=\max_{i\leq k}\left\| (x)_i \right\| ^{1/i}$ and $C_{j}>0$ is a constant depending on $P_{j}$ and on the norms $\left\| \cdot\right\| _{i}$'s. 
Since $\varepsilon >0$, when expanding the right-hand side of (\ref{zooo}) all terms of the form $|x|_{\lambda }^{l}|y|_{\lambda }^{m}$ with $l+m\leq j$ appear with some positive coefficient, say $\varepsilon_{l,m}$. 
The terms $|x|_{\lambda }^{j}$ and $|y|_{\lambda }^{j}$ appear with coefficient $1$ and cause no trouble since we always have $\lambda_{j}\left\| (x)_j\right\| \leq |x|_{\lambda }^{j}$ and $\lambda_{j}\left\| (y)_j\right\| \leq |y|_{\lambda }^{j}$. 
Therefore, for \eqref{zooo} to hold, it is sufficient that
\begin{equation*}
\lambda _{j}C_{j}M_{j-1}(x)^{l}M_{j-1}(y)^{m}\leq \varepsilon
_{l,m}|x|_{\lambda }^{l}|y|_{\lambda }^{m},
\end{equation*}
for all remaining $l$ and $m.$ 
However, clearly $M_{k}(x)\leq \Lambda_{k}\cdot |x|_{\lambda }$ where $\Lambda _{k}:=\max_{i\leq k}\{1/\lambda_{i}^{1/i}\}\geq 1$.  
Hence, a sufficient condition for \eqref{zooo} to hold is
\begin{equation*}
\lambda _{j}\leq 
\left({C_{j}\Lambda _{j-1}^{j}}\right)^{-1}{\min_{l+m\leq j} \varepsilon _{l,m}}.
\end{equation*}
Since $\Lambda_{j-1}$ depends only on $\lambda _1, \ldots, \lambda _{j-1}$, 
such a set of conditions can be fulfilled by a suitable $s$-tuple $\lambda$. 
\end{proof}
 
Every two homogeneous quasi-norms are biLipschitz equivalent; see Exercise~\ref{ex_homog_seminorm_biLip}. 
Moreover, for every left-invariant geodesic distance, the distance from the identity element is coarsely biLipschitz equivalent to each homogeneous quasi-norm, in the following sense.
 
\begin{theorem}[Guivarc'h]\label{thm_guiv}\index{Guivarc'h! -- Theorem}\index{Theorem! Guivarc'h --}
 Let $G$ be a nilpotent simply connected Lie group equipped with a left-invariant geodesic distance $d$. 
 Let $|\cdot|$ be a homogeneous quasi-norm with respect to a compatible linear grading.
 Then there exists a constant $C >1$ such that
 \begin{align}\label{eq_thm_guiv}
 \frac{1}{C}|g| - C \leq d(1,g) \leq C |g| + C, \quad \forall g \in G.
 \end{align}
\end{theorem}

\begin{proof}
Recall that Guivarc'h seminorms, as defined in Lemma \ref{lem_Guivarch}, are homogeneous quasi-norms.
Then, since homogeneous quasi-norms (with respect to the same dilations) are biLipschitz equivalent, see Exercise~\ref{ex_homog_seminorm_biLip}, it is enough to show \eqref{eq_thm_guiv} when $|\cdot| $ is some Guivarc'h seminorm, for some $\eps>0$ as to satisfy \eqref{eq_lem_guiv}. 
 
 We begin with the first inequality: $\frac{1}{C}|x| - C \leq d(1, x)$. 
 Let \[C \coloneqq \max \{|x| \colon d(1, x) \leq 1\} + \varepsilon.\]
 Take $x \in G$ with $L \coloneqq d(1, x)$. Take a geodesic from $1$ to $x$ and subdivide it into $k \coloneqq \lfloor L \rfloor + 1$ pieces $x_0 = 1, x_1, \ldots , x_k = x$ with $d(x_{i-1}, x_i) \leq 1$, for every $i = 1,2, \ldots ,k$. Then we have 
 \begin{align*}
 d(1, x_{i-1}^{-1}x_i)= d(x_{i-1}, x_i)
 \leq 1, \quad \forall i = 1, 2, \ldots , k, 
 \end{align*}
 and thus, by definition of $C$, we have
\begin{align}\label{guiv_norm_segment}
 |x_{i-1}^{-1}x_i| \leq C - \varepsilon, \quad \forall i = 1, 2, \ldots , k.
 \end{align}
 Hence, using that $|\cdot|$ comes from Lemma~\ref{lem_Guivarch}, we can bound 
 \begin{eqnarray*}
 |x| \, \, &=& \, \, |x_1 \cdot x_1^{-1}x_2 \cdot \cdots \cdot x_{k-1}^{-1}x_k| \\
 &\stackrel{\text{\eqref{eq_lem_guiv}}}{\leq} &\!|x_1| + |x_1^{-1}x_2| + \cdots + |x_{k-1}^{-1}x_k| + k\varepsilon \\
 &\stackrel{\text{\eqref{guiv_norm_segment}}}{\leq} &k (C-\varepsilon)+k\varepsilon \\
 &=& kC \leq C(L+1) = Cd(1, x) + C.
 \end{eqnarray*}
We then prove the second inequality: $d(1, x) \leq C|x| + C$. 
The proof is by induction on the nilpotency step $s$ of $G $. 
In the abelian case, it is clear. 
We assume the result is proved for groups up to step $s-1$, and we want to prove it for an $s$-step group $G$. 
We consider the quotient modulo the normal subgroup given by the last non-trivial element $C^s(G)$ of the lower central series.
Note that $\mfaktor{C^s(G)}{G}=\faktor{G}{C^s(G)}$ is a Lie group of step $s-1$ that is equipped with the distance \eqref{dist_HG}, which in this case is left-invariant and geodesic since $d$ is geodesic; see Proposition~\ref{proj_geod_geod}. 
Hence, the result is valid in $\faktor{G}{C^s(G)}$. Namely, there exists $C>0 $ such that for every $x \in G$ we have
\[d(1, xC^s(G)) \leq C|xC^s(G)| + C.\]
It follows that, there exists $z \in C^s(G)$ such that 
\begin{align}\label{ineq_guiv_xz}
 d(1, xz) \leq C|xC^s(G)| + C\leq C|x|+C,
\end{align}
where we used the fact that projections reduce Guivarc'h seminorms. 
On the one hand, by triangle inequality and left invariance, we have 
\begin{eqnarray*}
 d(1, x) &\leq& d(1, xz) + d(xz, x) \\
&\stackrel{\text{\eqref{ineq_guiv_xz}}}{\leq} &\!C|x|+C + d(1,z).
\end{eqnarray*}
On the other hand, by Lemma~\ref{lem_Guivarch}, we also have
\begin{eqnarray*}
 |z| \, \,&=& \, \,|x^{-1} \cdot xz| \\
&\stackrel{\text{\eqref{eq_lem_guiv}}}{\leq}& \!|x^{-1}| + |xz| + \varepsilon \\
 &\leq & \, \,|x| + Cd(1, xz) + C + \varepsilon \\
 &\stackrel{\text{\eqref{ineq_guiv_xz}}}{\leq}& \!|x| + C(C|x|+C)+C+ \varepsilon \\
 &= &(C^2 + 1)|x|+C^2+C+\varepsilon,
\end{eqnarray*}
where in the second inequality, we used that in exponential coordinates, the inverse of $x$ is $-x$, and the inequality proved in the first part of this proof.
Therefore, we just need to bound $d(1,z)$ affinely in $|z|$, for $z \in C^s(G)$. 

Let $c \coloneqq \dim C^s(G)$ and $\Omega \coloneqq B_d(1_G,1)$. Then the set 
\[\Omega' \coloneqq \left\{ e_1 \cdots e_c \colon \forall i = 1,2, \ldots,c \; \exists \; u_{i1}, \ldots,u_{is} \in \Omega \colon e_i = [u_{i1}, \ldots,u_{is}] \right\}\]
is a neighborhood of $1$ in $C^s(G)$; see Exercise \ref{lem_exercise_1}. 
Then there exists $m \in \mathbb{N}$ such that $\log (\Omega')$ contains the $\frac{1}{m^s}$-ball with respect to the norm $\norma{\cdot}$ on $V_s$ giving the Guivarc'h seminorm. 
For every $z \in C^s(G)$, let $k \in \mathbb{N}$ be such that
\begin{equation}\label{eq_Guiv_kz}
 k-1 \leq |z| \stackrel{\text{{\rm def}}}{=} \norma{z}^{1/s} \leq k.
\end{equation}
Consequently, without distinguishing the group with the Lie algebra via the exponential map, we have 
\[\frac{1}{m^sk^s}z \in \Omega'.\]
Hence, we can write 
\[\frac{1}{m^sk^s}z = e_1 \cdot \ldots \cdot e_c, 
\]
with $e_i = [u_{i1}, \ldots,u_{is}]$ and $u_{ij} \in B_d(1_G,1)$, for every $i,j$. Then, using that each $[u_{i1}, \ldots ,u_{is}]$ is a central element, we get
\begin{align*}
 z &= \left(\frac{1}{m^sk^s}z\right)^{m^sk^s} \\
 &= \left(\prod_{i=1}^c[u_{i1}, \ldots ,u_{is}] \right)^{m^sk^s} \\
 &= \prod_{i=1}^c[u_{i1}, \ldots ,u_{is}]^{m^sk^s} \\
 &= \prod_{i=1}^c[u_{i1}^{mk}, \ldots ,u_{is}^{mk}],
\end{align*}
where in the last equation, we use this general property of nilpotent groups; see Exercise~\ref{lem_exercise_2}.
Then, iteratively using the triangle inequality and the fact that $u_{ij}$ are in the unit ball, we infer that for some constant $C_{c,s}$, we can bound
\begin{eqnarray*}
 d(1,z) &=& d\Big(1, \prod_{i=1}^c[u_{i1}^{mk}, \ldots ,u_{is}^{mk}]\Big) \\
 &\leq &C_{c,s} \, \max_{i,j}d(1, u_{ij}^{mk}) \\
 &\leq &C_{c,s} \,mk \, \max_{i,j}d(1,u_{ij}) \\
 &\leq &C_{c,s,m} \, k \stackrel{\eqref{eq_Guiv_kz}}{\leq} C_{c,s,m} (|z| + 1). 
\end{eqnarray*}
\end{proof}

\section{Exercises}


\begin{exercise}\label{ex G/H simply connected}
Let $H$ be a closed subgroup of a topological group $G$. Assume that $G$ is simply connected and $H$ is connected. Then, the topological space $G/H$ is simply connected.
\\
{\it Hint.} Every lift of a loop has extremal points joining elements in $H$.
\end{exercise} 

\begin{exercise}\label{ex same ab reason}
For every nilpotent simply connected Lie group $G$, one has
$$ \exp(X+[\g, \g]) =
\exp(X)\exp([\g, \g]) = \exp(X)[G, G], \qquad\forall X\in \g.$$
{\it Hint.} Use Baker-Campbell-Hausdorff formula.
\end{exercise} 

\begin{exercise}\label{ex20030242050}
For every sub-Finsler nilpotent simply connected Lie group $(G, \Delta, \norm{\cdot})$, as in Definition~\ref{def abelianization norm}, 
the norm given in \eqref{def_norm_abelianization} is
$$\norm{v}_{\rm ab} = \min\{ \norm{w} : w\in \Delta_{1_G}, \pi_{\rm ab}(w) = v\}, \qquad \forall v\in {\rm Ab}(G).$$
\end{exercise}

\begin{exercise}\label{ex pi ab homo} 
Let $G$ be a simply connected and nilpotent Lie group. 
See ${\rm Ab}(G)$ as $\g/[\g, \g]$ with abelianization map $\pi_{\rm ab}: \g \rightarrow {\rm Ab}(G)$.
Consider the map $\tilde \pi : G \rightarrow {\rm Ab}(G)$ given by $\tilde \pi := \pi_{\rm ab}\circ\log$.
Then, we have $$\tilde\pi(xy)=\tilde\pi(x)+\tilde\pi(y), \qquad \forall x,y\in G.$$ 
{\it Hint.} Baker-Campbell-Hausdorff formula gives $\exp^{-1}(xy)\equiv \exp^{-1}(x)+\exp^{-1}(y)$ modulo $[\g, \g]$.
\end{exercise}

\begin{exercise}\label{ex_second_proof for prop lift}
Using selection arguments, like Theorem~\ref{thm_selection} as used in Proposition~\ref{prop_subFin_submetry_Lie}, give a proof of Corollary~\ref{cor prop lift}.
\end{exercise}

\begin{exercise} \label{V1generates}
 Let $\g$ be a nilpotent Lie algebra and let $\Delta_1$ be a subspace of $\g$ such that
\begin{equation}\label{special Delta 1}
\g= \Delta_1+[\g, \g].
\end{equation}
For the lower central series $(\g^{( i)})_{i\in \N}$ of $\g$,
we have $$\g^{( 2)}= [\Delta_1, \Delta_1]+\g^{( 3)} $$
and, more generally, we have
$$\g^{( i)}= [\Delta_1,[\Delta_1,[\ldots, [\Delta_1, \Delta_1]\ldots]]+\g^{( i+1)}, $$
where in the above brackets $\Delta_1$ appears $i$-many times.
Finally, deduce that every $\Delta_1$ with property~\eqref{special Delta 1} is Lie generating $\g$.
\end{exercise}

\begin{exercise} \label{V1generates2}
 Let $\g$ be a nilpotent Lie algebra and let $\Delta_1$ be a subspace of $\g$.
Then, the set $\Delta_1$ Lie generates $\g$ if and only if 
$\g= \Delta_1+[\g, \g].$
\end{exercise}

\begin{exercise}\label{smaller_no_generates}
Let $\g$ be a nilpotent Lie algebra. Let $V\subset\g$ be a subspace such that $V+ [\g, \g]\neq \g$. 
Then $V$ is not bracket-generating $\g$. 
\\ {\it Hint.} Consider the vector subspace $\pi_{\rm ab}(V)$.
\end{exercise}

\begin{exercise}\label{Delta 1 not unique}
Find a nilpotent simply connected Lie group $G$ and sub-spaces $\Delta_1$ and $\tilde \Delta_1$ such that
$$\g= \Delta_1\oplus[\g, \g]=\tilde \Delta_1\oplus[\g, \g],$$
such that there is no Lie algebra isomorphism $\phi:\g \to\g$ such that $\phi(\Delta_1)=\tilde \Delta_1$. Compare it with Proposition~\ref{isomorphic:stratifications}. 
\\{\it Hint.} Check Exercise~\ref{ex Breuillard strange V 1}.
\end{exercise}

\begin{exercise}\label{lemma4homo}
Let $\g$ and $\mathfrak h$ be nilpotent Lie algebras.
Let $\Delta^\g_1$ (resp., $\Delta^H{\mathfrak h}_1$) be a linear subspace of $\g$ (resp., $\mathfrak h$) such that $\g= \Delta^\g_1\oplus[\g, \g]$ (resp.,
$\mathfrak h= \Delta^{\mathfrak h}_1\oplus[\mathfrak h, \mathfrak h]$).
 Let $\phi:\g\to\h$ a Lie algebra homomorphism.
 If $\phi$ has the property that 
$\phi(\Delta^\g_1) \subseteq \Delta^{\mathfrak h}_1 $,
then 
$$ {\rm proj}_{\Delta^{\mathfrak h}_1}\circ \phi=\phi\circ{\rm proj}_{\Delta^\g_1} ,
\qquad \text{ i.e., }\quad
 \begin{tikzcd}
\g\ar[r]{u}{\phi}
\ar[d,->>]{d}{{\rm proj}} 
&\h
\ar[d,->>]{u}{{\rm proj}}
 \\
\Delta^\g_1 \ar[r]{u}{\phi } & \Delta^{\mathfrak h}_1
 \end{tikzcd},
$$
where ${\rm proj}^\g:\g\to\g$ and ${\rm proj}^{\mathfrak h}:\h\to\h$ are the projections onto ${\Delta^\g_1}$ and ${\Delta^{\mathfrak h}_1}$ respectively with kernels $[\g, \g]$ and $[\h, \h]$ respectively. 
\\{\it Solution.}
If $X\in \Delta^\g_1$, then $(\phi\circ{\rm proj})(X)=\phi(X)$. Since by assumption we also have $\phi(X)\in \Delta^{\mathfrak h}_1$,
then $({\rm proj}\circ \phi)(X)=\phi(X)$.
So $ {\rm proj}\circ \phi$ and $\phi\circ{\rm proj}$ are two homomorphisms that coincide on $ \Delta^\g_1$.
Since $\Delta^\g_1$ generates the algebra $\g$, then the two homomorphisms are equal.
\end{exercise}

\begin{exercise}
Let $G$ be a nilpotent simply connected Lie group polarized by $\Delta$ so that $\Delta_1\oplus[\g, \g] = \g$.
 The projection map $\pi:=\pi_{\Delta_1} $ from \eqref{projection_pi1} has the following properties:
 For every Lipschitz curve $\sigma$ in $\Delta_1$ with $\sigma(0)=0$, there exists a unique Lipschitz horizontal curve $\gamma$ with $\pi(\gamma)=\sigma$ and $\gamma(0)=1_G$, and such a curve is the solution of the ODE
 \begin{equation}\label{lift}
 \left\{ \begin{array}{ccl}
\dot\gamma(t)&=&(L_{\gamma(t)})_* \dot\sigma(t)\\
\gamma(0)&=&1_G.
\end{array} \right.\end{equation}
{\it Solution.} 
We discussed the existence and the uniqueness of the ODE in the proof of Proposition~\ref{prop_Integration_tangent_vector}.
Let $\gamma(t)$ be the solution.
Then $$\gamma'(t)=(L_{\gamma(t)})^*\dot\gamma(t)=\dfrac{\dd}{\dd t}(\sigma(t)).$$
Because of Formula \eqref{gamma'formula}, we have that $\pi\circ\gamma$ and $\sigma$ are two curves in $\Delta_1$ with same starting point $\pi(\gamma(0))=0=\sigma(0)$ and same derivative:
$ \dfrac{\dd}{\dd t}( \pi\circ\gamma)=\dfrac{\dd}{\dd t}\sigma$.
Therefore $\pi\circ\gamma=\sigma$.
\end{exercise}

\begin{exercise}\label{prop: same regularity0}
Let $G$ and $H$ be sub-Finsler Lie groups and $\varphi: G \to H$ a Lie group homomorphism that is also a submetry.
If geodesics in $G$ are analytic, then so are those in $H$. 
\\{\it Hint.} See Proposition~\ref{prop: same regularity2}.
\end{exercise}

\begin{exercise}\label{ex:same regularity2}
 Let $G$ be a nilpotent simply connected Lie group and let $\Delta_1$ be a polarization on $G$ such that \eqref{special: nilpotent} holds. 
If $\sigma: I\to \Delta_1$ is a curve of regularity $C^k$, for some $k\in \N$, then the multiplicative integral of $\sigma$ is $C^k$.
\\{\it Hint.} Consider the free-nilpotent Lie group $\tilde G$ with same rank and step of $G$, so that there is a quotient map $\tilde G\to G$. In this situation, it is clear that the multiplicative integral of $\sigma$ in $\tilde G$ is $C^k$ because the system \eqref{lift_multiplicative_integral} can be integrated one stratum at a time. 
The multiplicative integral of $\sigma$ in $ G$ is the projection of the one in $\tilde G$.
\end{exercise}


\begin{exercise}\label{ex8may1429}
Given a Lie algebra $\g$ and $v\in \g$, the transformation $\ad_v$ is diagonalizable over $\C$ and its eigenvalues are purely imaginary if and only if $\{e^{t\ad_v}:t\in\R\}$ is precompact in $\GL (\g)$. 	
\end{exercise}



\begin{exercise}[Killing form]\label{ex_def_Killing}\index{Killing form}
\index{form! Killin --}
Let $\g$ be a Lie algebra. The {\em Killing form} on $\g$ is defined for all $X,Y\in \g$ as 
$$B(X, Y) := \mathrm{Trace}(\mathrm{ad}_X\mathrm{ad}_Y).$$
The Killing form $B$ is bilinear, symmetric, and
$B(\ad_X Y, Z) - B(Y, \ad_X Y) =0,$ for all $X,Y, Z\in \g$
and invariant under automorphisms of the algebra g, that is,
$B(\psi(X), \psi(Y)) = B(X,Y) $ for $\psi \in \Aut_{\rm Lie}(\g)$, and $X,Y\in \g$.
 \end{exercise}

 \begin{exercise}\label{nilpotent: Killing}
\index{nilpotent! -- Lie algebra}
\index{Lie algebra! nilpotent --}
The Killing form of a nilpotent Lie algebra is identically zero -- be aware that the inverse is not true already for some 3D Lie algebra.
\end{exercise}

\begin{exercise}[Nilradical]\label{ex: nilradical}\index{nilradical}\index{${\rm nil}(\mathfrak g)$}
The {\em nilradical} of a Lie algebra $\g$,
denoted by ${\rm nil}(\mathfrak g)$, is defined as the largest nilpotent ideal of $\g$.
Then, the set ${\rm nil}(\mathfrak g)$ can also be defined as the sum of all
nilpotent ideals of $g$; see \cite[Definition~5.2.10]{Hilgert_Neeb:book}.
\end{exercise}

\begin{exercise}\label{ex_Killing_on_nilrad0}
Let $X$ be an element in ${\rm nil}(\mathfrak g)$ of a Lie algebra $\g$.
Then $\ad_X $ is a nilpotent transformation of $\g$.
\end{exercise}

\begin{exercise}\label{ex_Killing_on_nilrad}
Given a Lie algebra $\g$, the Killing form restricted to the nilradical of $\g$ is identically 0.
\\{\it Hint.} 
From Exercise~\ref{ex_Killing_on_nilrad0}, for $X\in {\rm nil}(\mathfrak g)$, we have $B(X, X)=0$.
\end{exercise}

\begin{exercise}[Solvable Lie algebra]\label{ex_def_solvable}
Let $\g$ be a Lie algebra. We have the following equivalent definitions:

(i) $\g$ is {\em solvable}, in the sense that, defined the {\em derived series} of $\g$ as
$\g^{(0)} := \g$ and
$$\g^{(n)} := [\g^{(n-1)}, \g^{(n-1)}] \quad n \in \N,$$
then the derived series terminates in the zero subalgebra.

(ii) $[\g, \g]$ is nilpotent.


%
%
%

(iii) the Killing form $B$ of $\g$ satisfies $B(X, Y) = 0$ for all $X\in \g$ and $Y \in [\g, \g]$.
 This is {\em Cartan's criterion for solvability}.\index{Cartan! -- criterion}\index{criterion! Cartan --}
\\ 
 {\it Hint.} Check \cite[page 31]{Knapp} and \cite[Proposition~1.39]{Knapp}.
\end{exercise}


%
%
%


\begin{exercise}\label{Lemma316}
Let $M_1$ and $M_2$ be groups. 
Suppose $ F \colon M_1 \to M_2$ is a map such that $ F\circ M_1^L \circ F^{-1} = M^L_2 $, where $M^L$ denote the group of left translations on a group $M$. 
Then $F$ is affine.
\end{exercise}


\begin{exercise}\label{property quasinorms}
Let $|\cdot|: G\to \R$ be a homogeneous quasi-norm on a positively graded Lie group $G$.
Then, the function $|\cdot|: G\to \R$ is proper, i.e., the pre-image of compact sets is compact.
 It induces the same topology, in the sense that $g\to \hat g$ in $G$ if and only if $|\hat g g^{-1}|\to 0$.
Moreover, there exists a constant $M\geq 1 $ such that 
$|g h| \leq M (|g| +|h|)$, for $g,h\in G $.
\end{exercise}

\begin{exercise}\label{ex_homog_seminorm_biLip}
 Let $G$ be a nilpotent simply connected Lie group whose Lie algebra is equipped with a linear grading.
 Let $|\cdot|$ and $\norm{\cdot}$ be homogeneous quasi-norms with respect to the grading.
Then there exists a constant $C >1$ such that
 \begin{align}
 \frac{1}{C} |g| \leq \norm{g} \leq C |g| , \qquad \forall g \in G.
 \end{align}
 {\it Hint.} For each inequality, take as constant the maximum of a quasi-norm on the unit ball of the other quasi-norm. Then, use homogeneity.
\end{exercise}

\begin{exercise}\label{lem_exercise_1}
 For every neighborhood $\Omega$ of $1$ in a nilpotent simply connected Lie group $G$, the set 
\[\Omega' \coloneqq \left\{ e_1 \cdots e_c \, \colon\, \forall i = 1,2, \ldots,c \; \exists \; u_{i1}, \ldots,u_{is} \in \Omega \colon e_i = [u_{i1}, \ldots,u_{is}] \right\}\]
is a neighborhood of $1$ in $C^s(G)$, where $c \coloneqq \dim C^s(G)$ and $s$ is the nilpotency step of $G$.
\\{\it Solution.} 
We work in exponential coordinates, identifying the Lie algebra $\g$ with $G$, and $C^s(\g)$ with $C^s(G)$.
Being $C^s(\g)$ a $c$-dimensional vector space, for $i \in \{ 1, \ldots,c\}$, there are $g_ {i1}, \ldots,g_{is} \in \mathfrak{g}$ such that $[g_{11}, \ldots,g_{1s}], \ldots,[g_{c 1}, \ldots,g_{cs}]$ form a basis of $\mathfrak{g}^s$. 
 Up to shrinking $\Omega $, we assume that $\Omega $ is a convex and symmetric neighborhood of $0$ in $\g$.
 We choose $\varepsilon > 0$ sufficiently small so that $t g_{ij}\in \Omega$, for all $ i\in \{ 1, \ldots,c\},j \in \{ 1, \ldots,s\} $ and all $t\in (-\varepsilon , \varepsilon)$.
 Define $u_{ij} \coloneqq \varepsilon g_{ij} \in \Omega$ and the map 
 \begin{align*}
 \psi \colon [-1,1]^c &\rightarrow \mathfrak{g}^s \\
 (t_1, \ldots,t_c) &\mapsto \sum_{i=1}^c t_i[u_{i1}, \ldots,u_{is}].
 \end{align*}
The image of $\psi$ is contained in $\Omega'$, since, assuming $t_1, \ldots, t_c>0$ and leaving the general case to the reader, we have
 \begin{align*}
 \sum_{i=1}^c t_i[u_{i1}, \ldots,u_{is}] &= 
 \sum_{i=1}^c [t_i^{1/s}u_{i1}, \ldots,t_i^{1/s}u_{is}] \\
 &= [t_1^{1/s}u_{11}, \ldots,t_1^{1/s}u_{1s}] \star \cdots \star [t_c^{1/s}u_{c1}, \ldots,t_c^{1/s}u_{cs}] \in \Omega',
 \end{align*}
 where in the final inclusion we use the fact that $t_i^{1/s} u_{ij} \in \Omega$ for every $i$ and $j$. In addition, the image of the differential of $\psi$ at the identity is the span of $[u_{11}, \ldots,u_{1s}], \ldots , [u_{c1}, \ldots,u_{cs}]$, hence $(\dd\psi)_0$ is surjective since these vectors form a basis of $\mathfrak{g}^s$. It follows that there exists an open neighborhood $U \subset[0,1]^c$ of $0$ such that $\psi(U) \subset \Omega'$ is an embedded $c$-dimensional submanifold of $\mathfrak{g}^s$, hence it is open.
\end{exercise}

\begin{exercise}\label{lem_exercise_2}
 Let $G$ be a nilpotent group of step $s$ and $m_1, \ldots,m_s \in \mathbb{N}$. Then 
\[[x_1^{m_1}, \ldots, x_s^{m_s}] = [x_1, \ldots, x_s]^{m_1 \cdots m_s}.\]
{\it Solution.}
 We begin by showing that if $[x_1, x_2]$ commutes with $x_1$ and $x_2$, then for every integers $m_1,m_2$ we have 
 \begin{align}\label{fact_commutator}
 [x_1^{m_1}, x_2^{m_2}]=[x_1, x_2]^{m_1m_2}.
 \end{align}
 Indeed, we have 
 \begin{align*}
 [x_1, x_2]^m &= x_1x_2x_1^{-1}x_2^{-1}[x_1, x_2]^{m-1}\\
 &=x_1[x_1, x_2]x_2x_1^{-1}x_2^{-1} [x_1, x_2]^{m-2}\\
 &=x_1x_1x_2x_1^{-1}x_2^{-1}x_2x_1^{-1}x_2^{-1} [x_1, x_2]^{m-2}. \\
 &=x_1^2x_2x_1^{-2}x_2^{-1} [x_1, x_2]^{m-2} \\
 &=x_1^2[x_1, x_2]x_2x_1^{-2}x_2^{-1} [x_1, x_2]^{m-3} \\
 &=x_1^3x_2x_1^{-3}x_2^{-1} [x_1, x_2]^{m-3} \\
 &=\ldots \\
 &=x_1^mx_2x_1^{-m}x_2^{-1} = [x_1^m, x_2].
 \end{align*}
As a consequence, we obtain 
\begin{align*}
 [x_1^{m_1}, x_2^{m_2}] = [x_1, x_2^{m_2}]^{m_1} = [x_2^{m_2}, x_1]^{-m_1} = [x_2, x_1]^{-m_1m_2} = [x_1, x_2]^{m_1m_2}.
\end{align*}

We solve the exercise by induction on the nilpotency step of $G$. In the abelian case, the result is trivial. Now assume that the result holds for groups up to step $s-1$. Hence, the result is valid in $\faktor{G}{C^s(G)}$ which has step $s-1$. Thus, for every $x_2, \ldots, x_{s} \in G$ and every $m_2, \ldots,m_{s} \in \mathbb{N}$ we have 
\begin{align*}
 [x_2^{m_2}C^s(G), \ldots, x_{s}^{m_{s}}C^s(G)] = [x_2C^s(G), \ldots, x_{s}C^s(G)]^{m_2\cdots m_{s}}.
\end{align*}
It follows that 
\[[x_2^{m_2}, \ldots, x_{s}^{m_{s}}]C^s(G) = [x_2, \ldots, x_{s}]^{m_2\cdots m_{s}}C^s(G), \]
and hence there exists a central element $z\in C^s(G)$ such that 
\begin{align}\label{eq_commutator}
 [x_2^{m_2}, \ldots, x_{s}^{m_{s}}] = [x_2, \ldots, x_{s}]^{m_2 \cdots m_{s}}z.
\end{align}
Let $x_1 \in G$ and $m_1 \in \mathbb{N}$, then 
\begin{eqnarray*}
 [x_1^{m_1}, x_2^{m_2}, \ldots, x_s^{m_s}] \;&\stackrel{\text{\eqref{eq_commutator}}}{=} &\![x_1^{m_1},[x_2, \ldots, x_s]^{m_2 \cdots m_s}z] \\
 &=& [x_1^{m_1},[x_2, \ldots, x_s]^{m_2 \cdots m_s}] \\
 &\stackrel{\text{\eqref{fact_commutator}}}{=} &\![x_1, x_2, \ldots, x_s]^{m_1 \cdots m_s},
\end{eqnarray*}
where in the second equality we use that for $x,y \in G$ and $z\in C^s(G)$ we have $[x,yz] = [x,y]$.
\end{exercise}

{\color{black}
Extra exercises:
\begin{exercise}[Radical]\index{radical}\index{solvable radical|see {radical}}\index{${\rm rad}(\mathfrak g)$}\label{def_radical}
Every Lie algebra $\g$ has a maximal solvable ideal, and it is unique.
This ideal is called the {\em radical} (also called {\em solvable radical} and denoted by ${\rm rad}(\mathfrak g) $) of $\g$.
\end{exercise}

\begin{exercise}[Simple Lie algebra]\label{def: simple: algebra}
By definition, a {\em simple} Lie algebra is a non-abelian Lie algebra whose only ideals are 0 and itself.
We have that the word `non-abelian' can be replaced with `dimension at least 2'.
\end{exercise}

\begin{exercise}[Semisimple]\label{def_semisiple}
A Lie algebra $\g$ is {\em semisimple} if $[\g, \g]=\g$.
We have that a Lie algebra is semisimple if and only if it is a direct sum of simple Lie algebras.
\end{exercise}

\begin{exercise}[Criteria for semi-simplicity]\label{Cartan criterion}\index{Cartan! -- criterion}\index{criterion! Cartan --}
Let $\g$ be a Lie algebra. We have the following equivalent definitions:

(i) $\mathfrak g$ is semisimple, as in Exercise~\ref{def_semisiple};

(ii) $\mathfrak g$ has no non-zero abelian ideals;


(iii) the solvable radical ${\rm rad}(\mathfrak g) $ is zero;

(iv)
$\skull$
the Killing form $B$ of $\g$ is non-degenerate.
 This is {\em Cartan's criterion for semi-simplicity}.
\end{exercise}

\begin{exercise}[Levi decomposition]\label{ex_Levi_decomp}
\index{Levi decomposition}
\index{decomposition! Levi --}$\skull$
Every Lie algebra $\g$ is the semidirect product of a solvable ideal and a semisimple subalgebra. 
In fact, 
then there exists a subalgebra ${\mathfrak l}$, called a {\em Levi factor}, that is semisimple and for which
$\g={\rm rad}(\g) \rtimes{\mathfrak l}.$
\\{\it Hint.} See \cite[Appendix B, Section 1]{Knapp}.
\end{exercise}

\begin{exercise}\label{ex: Killing_radicals}
$\skull$
Given a Lie algebra $\g$, consider its nilradical ${\rm nil}(\mathfrak g)$, its radical ${\rm rad}(\mathfrak g)$, and a Levi factor $\mathfrak l$ in $\mathfrak g$. Then the Killing form $B$ satisfies:

\ref{ex: Killing_radicals}.i. $B({\rm rad}(\mathfrak g),[\g, \g])=\{0\}$; (See \cite[I.5.5,p.48, Prop.5b]{Bourbaki_Lie_1_3})

 \ref{ex: Killing_radicals}.ii. $B({\rm nil}(\mathfrak g), \mathfrak g)=\{0\}$. (See \cite[I.4.4,p.42, Prop.6b]{Bourbaki_Lie_1_3})
 \end{exercise} 
 
 \begin{exercise}\label{contenimenti_radicali}
Let $\mathfrak g$ be a Lie algebra with radical $ {\rm rad}(\mathfrak g)$ and nilradical ${\rm nil}(\mathfrak g)$. 
Then, we have that $[\mathfrak g,{\rm rad}(\mathfrak g)]\subset {\rm nil}(\mathfrak g)$, and consequently, every vector space $V$ with ${\rm nil}(\mathfrak g)\subset V\subset {\rm rad}(\mathfrak g)$ is an ideal of $\mathfrak g$.
 \end{exercise}

\begin{exercise}\label{semisimple_poly}$\skull$
Let $G$ be a connected semisimple Lie group.
Then $G$ is of type (R) if and only if it is compact. 
\\{\it Hint.} Consult \cite[Proposition II.4.8.III]{MR2000440}, or \cite{MR0349895}, \cite{MR0316625}, \cite{Breuillard_MR3267520}.
\end{exercise}

\begin{exercise}
In the proof of Theorem~\ref{teo05171840baby}, we have that the nilradical of the connected component of $\Isom(M,d)$ acts almost simply transitive on $M$, i.e., orbits are open immersions.
\end{exercise}

 }



\chapter{Carnot groups}\label{ch_Carnot}
Carnot groups are specific examples of Carnot-Carathéodory spaces.
 They are simply connected nilpotent Lie groups whose Lie algebras admit particular gradings: stratifications. 
 The polarizations correspond to first layers of stratifications with norms on them. 
 Every such sub-Finsler Lie group is self-similar with respect to its Carnot-Carathéodory metric. Specifically, there is a natural family of metric dilations.
 
This chapter starts with Section~\ref{sec Definition of Carnot groups}, where we discuss the definition of Carnot groups, Carnot bases, and examples.
In Section~\ref{sec consequences of dilations}, we present some fundamental consequences of dilations, such as the presence of good coordinates, an easy proof of the Ball-Box Theorem, and properties of Haar measures.

Section~\ref{sec Pansu-Rademacher Theorem} is dedicated to Pansu-Rademacher Differentiation Theorem for Lipschitz maps between Carnot groups.
In Section~\ref{sec characterization of Carnot}, we provide a complete characterization of the metric spaces that are isometric to Carnot groups, entirely in terms of metric geometry.

In Section~\ref{sec Extremal curves in Carnot}, we explore extremal curves in Carnot groups. We outline some properties of abnormal curves and of normal geodesics. We also present a sublinear isometric property of projections of geodesics, 
which shows that geodesics cannot form corners. Additionally, we discuss some open problems.

Lastly, we include supplementary material:
Section~\ref{sec self-similar sub-Finsler spaces} on the Lie coset structure of self-similar sub-Finsler spaces, which are submetry images of Carnot groups; and 
Section~\ref{sec Local isometries of Carnot groups} on isometries between open sets in Carnot groups.



\section{Definition of Carnot groups}\label{sec Definition of Carnot groups}

{\em Carnot groups} are sub-Finsler Lie groups (as in Definition~\ref{def_Sub-Finsler_Lie_group}) such that the polarization is the first layer of a stratification (as in Definition~\ref{defiCARNOT}).
We spell out this definition to have a slightly more self-contained presentation.

Let $G$ be a simply connected Lie group.
Assume its Lie algebra $\g:=\Lie(G) $ admits a stratification $ \g= V_1\oplus\dots\oplus V_s$, i.e., 
$V_1,\ldots, V_s$ are vector subspaces in direct sum such that 
$V_{j+1}= [ V_j,V_1]$ with $V_{s+1}:=\{0\}$.
We recall from Section~\ref{sec_strat} that each stratification gives a Lie algebra grading 
$$[V_i, V_j]\subseteq V_{i+j},$$
see Exercises~\ref{ex stratifications are gradings}, and the layer $V_1$ of degree 1 is bracket generating; see Remark~\ref{other def stratification}.

The vector space $V_1$ is called the {\em horizontal stratum} and is seen as a subset of the tangent space $T_{1_G}G$ of $G$ at the identity element $1_G$ of $G$.
As in \eqref{eq:induced:distribution}, it
induces a left-invariant subbundle $\Delta$, called the {\em horizontal bundle}, of the tangent bundle $TG$:
\begin{equation}\label{distribution:Carnot}
\Delta_g := (L_g)_*V_1,\qquad\forall g\in G .
\end{equation}
Fix a norm $\|\cdot\|$ on the vector space $V_1$.
As in \eqref{cont_var_norm}, the norm on $V_1$ induces a norm on every $\Delta_g$ as
\begin{equation}\label{norm:Carnot}
\|v\| := \|(L_g)^* v\|,
\qquad\forall v\in\Delta_g,
\quad\forall g\in G.
\end{equation}
 	The triple $(G,\Delta,\|\cdot\|)$ is a Carnot-Carathéodory space, and indeed, a sub-Finsler Lie group, which has an induced distance function as in Definition~\ref{def_cc_dist_mfd}:
	\begin{equation}\label{dista_in_Carnot} \dcc (p,q) :=d_{V_1,\norm{\cdot}}(p,q) := 
 \inf \left\{\int \| {\dot\gamma}\| \, \colon \, \gamma \; \text{AC curve from } p \text{ to } q, \text{ with } \dot\gamma \in \Delta \right\},\qquad \forall p,q\in G.\end{equation}
 

\begin{definition}[Carnot group]
\index{Carnot! -- group}\index{Carnot! -- algebra}
Let $G$ be a simply connected Lie group whose Lie algebra admits a stratification.
Given a first stratum $V_1$ of the stratification of $\Lie(G)$ and a norm on it,
let $\Delta$ and $\|\cdot\|$ be defined by \eqref{distribution:Carnot} and \eqref{norm:Carnot}, respectively.
Let $\dcc$ be the Carnot-Carath\'eodory distance associated with $\Delta$ and $\|\cdot\|$ as in \eqref{dista_in_Carnot}.
Both the sub-Finsler manifold $(G,\Delta,\|\cdot\|)$ and the metric space $(G,\dcc)$ are called \emph{Carnot groups}.
In accordance with Definition~\ref{defiCARNOT}, we call {\em Carnot algebras} the Lie algebras of Carnot groups.
\end{definition}

Given a stratification $\Lie(G) = V_1\oplus\dots\oplus V_s$ of the Lie algebra of a Carnot group $G$, with $V_s\neq\{0\}$,
the number $s$ is called the {\em step} of $G$ and the number $\dim V_1$ is called {\em rank} of $G$.\index{rank! -- of a Carnot group}\index{step! -- of a Carnot group}
The topological dimension of $G$
is $n:=\sum_i{\rm dim\,}V_i$
and the {\em homogeneous dimension} 
is\index{homogeneous! -- dimension}
\begin{equation}\label{DefOfQ:3}
Q:=\sum_{i=1}^s i\,{\rm dim\,}V_i.
\end{equation}

Each Carnot group $(G,\Delta,\|\cdot\|)$ is indeed a sub-Finsler Lie group and an equiregular Carnot-Carath\'eodory space of step $s$. Indeed,
one has that, for each $k\in\{1,\ldots, s\}$, the subset $
\Delta^{[k]}$ in the flag of subbundles for $\Delta$ as in Definition~\ref{def:equiregular} is the left-invariant subbundle for which 
$$
\Delta^{[k]}({1_G})=V_1\oplus \cdots\oplus V_k.$$

 Because of Proposition~\ref{biLipschitz_equivalence_CCnorms}, other choices of norms would not change the biLipschitz equivalence class of the CC metric. 

 	
 \subsection{Dilations on Carnot groups}

As in every $\R$-graded Lie group, in Carnot groups we have a canonical one-parameter family of dilations on the Lie algebras and on the groups; recall the discussion on page \pageref{discuss_dilations20jun}:
	
\begin{definition}[Dilations on stratified groups]\label{Dilations:groups}
\index{dilation! -- in Carnot group}
 	Let $G$ be a Carnot group.
	Let $\delta_\lambda:\Lie(G)\to\Lie(G)$ be the dilation of factor $\lambda$ associated with the stratification as in \eqref{def_dilation_relative_to_grading}.
	Then the \emph{dilation} $\delta_\lambda:G\to G$ of the group of factor $\lambda$ is the only group automorphism such that $(\delta_\lambda)_*=\delta_\lambda$.
	Such maps are also called the {\em intrinsic dilations} of the Carnot group or {\em Carnot dilations}. \index{Carnot! -- dilation}\index{intrinsic! -- dilations} 
\end{definition}
We have kept the same notation $\delta_\lambda$ for both dilations
(in $\Lie(G)$ and in $G$) because no ambiguity will arise since the two maps have different domains.

On Carnot groups, the intrinsic dilations satisfy the following formulas:
	\begin{eqnarray}\label{delta_exp=exp_delta1}
&	\delta_\lambda\circ\exp = \exp\circ\delta_\lambda,\qquad &\forall \lambda\in \R;
\\
\label{fomula414392}
&\delta_{\lambda}\circ\delta_{\eta}=\delta_{\lambda\eta},\qquad &\forall \lambda,\eta\in \R;
\\
%
& \label{dilat-prod}
\delta_\lambda(xy)=\delta_\lambda (x)\delta_\lambda (y)
\qquad&\forall x,\,y\in G, \forall \lambda\in \R.
\end{eqnarray}

\subsubsection{Relations between dilations and CC distances}
The Carnot-Carath\'eodory distance is well-behaved under the intrinsic
dilations, in the sense that such dilations multiply distances of a constant factor.

\begin{proposition}\label{dilations_on_Carnot2024}
 	If $(G,\dcc)$ is a Carnot group with dilations $(\delta_\lambda)_{\lambda\in \R}$,	then
	\begin{equation}\label{dilat-dist}
\dcc(\delta_\lambda p,\delta_\lambda q) = \lambda\dcc(p,q), \qquad\forall p,\,q\in G, \forall \lambda\in \R.
\end{equation}
\end{proposition}
\begin{proof}
 	Since $\delta_\lambda|_{V_1}$ is the multiplication by $\lambda$, we have that $\|\delta_\lambda v\| = \lambda \|v\|$, for all $v\in\Delta$.
	If $\gamma$ in a horizontal curve from $x$ to $y$, then $\delta_\lambda \circ \gamma$ is a curve going from $\delta_\lambda x$ to $\delta_\lambda y$ whose tangent vectors are, for almost all $t$,
\begin{equation}\label{dilated-tang-vec}
 (\delta_\lambda)_*\dot\gamma(t)=\delta_\lambda (\dot\gamma(t))=\lambda \dot\gamma(t),
\end{equation}
which are horizontal since $\dot\gamma(t)$ is horizontal. Moreover, from \eqref{dilated-tang-vec}, the length of $\delta_\lambda \circ \gamma$ is $\lambda$ times the length of $\gamma$, i.e., 
 for every horizontal curve $\gamma$,
	\[
	\Length_{\|\cdot\|}(\delta_\lambda\circ \gamma) = \lambda \Length_{\|\cdot\|}(\gamma) .
	\]
	 Thus, by \eqref{dista_in_Carnot} we get \eqref{dilat-dist}.
\end{proof}

\subsection{Good bases for Carnot groups}\label{sec Good bases for Carnot}
Let $G$ be a Carnot group with stratification $\g= V_1\oplus \cdots\oplus V_s. $
We want to construct a basis for $\g$ that is structured with respect to the stratification, is a Malcev basis, and each element of the basis that is not in $V_1$, is the bracket of two vectors of such a basis.

Start by picking a basis $X_1,\ldots,X_m$ of $V_1$. Then consider all brackets $[X_i,X_j]$, for $i,j\in\{1,\ldots,m\}$. Since $[V_1,V_1]=V_2$, we can find among such brackets a basis for $V_2$; see Exercise~\ref{basisV2}. Pick one such a basis and call its elements $X_{m+1},\ldots,X_{m_2}$.
Iterate the method: extract a basis $X_{m_2+1},\ldots,X_{m_3}$ of $V_3$ from the set $[X_i,X_j]$, for $i\in\{1,\ldots,m\}, j\in\{m+1,\ldots,m_2\}$. And so on. 
We have constructed a basis $X_1,\ldots,X_n$ of $\g$ and natural numbers $m_1\ldots, m_s$ such that
\begin{enumerate}
 \item $X_{m_{j-1}+1},\ldots,X_{m_j}$ is a basis of $V_j$,
\item For every $i\in\{m+1,\ldots,n\}$, there exist $d_i$, $l_i$, and $k_i$ such that $X_i\in V_{d_i}$, $X_{l_i}\in V_{1}$, $X_{k_i}\in V_{d_i-1}$, and
\begin{equation}
X_i=[ X_{l_i}, X_{k_i}].\label{Carnotbasis}
\end{equation}
\item The order-reversed basis $(X_n,\ldots,X_1)$ is a Malcev basis as in Definition~\ref{def Malcev basis}; in other words,
$$[\g,\Span\{ X_k,\ldots,X_n\} ] \subseteq \Span\{ X_{k+1},\ldots,X_n\}, \qquad \forall k\in \{1,\ldots, n\} .$$
\end{enumerate}
We suggest the terminology `{\em Carnot basis}' for a basis satisfying the above three conditions.
The reader should notice that the above property 1 implies the property 3. See Exercise~\ref{CarnotMalcev}.

To describe a Carnot algebra, we prefer to give a Carnot basis with a hierarchical diagram as follows:

\begin{center}\index{Heisenberg! -- group}
 \begin{tikzcd}[end anchor=north]
 V_1:&& X \ar[dr, no head]& &Y\ar[dl, no head] & \text{ for the 3D Heisenberg algebra;}\\
 V_2:&& &Z &\quad\;&
 \end{tikzcd}
\end{center}

\begin{center}\index{Engel! -- group}
 \begin{tikzcd}[end anchor=north]
 V_1:&& X \ar[dr, no head] \ar[dd, no head]& &Y\ar[dl, no head] & \text{ for the Engel algebra}\\
 V_2:&& &Z \ar[dl, no head] &\quad\;&\text{(the step-3 filiform algebra);}\\
 V_3:&&W & &\quad\;&
 \end{tikzcd}
\end{center}

\begin{center}\index{Cartan! -- group}
 \begin{tikzcd}[end anchor=north]
 V_1:&& X_1 \ar[dr, no head] \ar[dd, no head,-<-=.5]& &X_2\ar[dl, no head]\ar[dd, no head] & \text{ for the Cartan Lie algebra, in Hall's basis}\\
 V_2:&& &X_{21} \ar[dl, no head,->-=.5]\ar[dr, no head] &\quad\;&\text{(the rank 2 and step 3 algebra).}\\
 V_3:&&X_{211} & &X_{212} &
 \end{tikzcd}
\end{center}
As explained at page \pageref{first description of diagram}, in the diagram, we obtain the bracket relations by reading the arms \textit{from left to right} unless there is an arrow in which case they are read \textit{from right to left}.\index{diagram! description of the --}\label{description of diagram}
The $j$-th line in the diagram lists the vectors that span the stratum $V_j$. The black lines express the non-trivial brackets.
In the Lie algebra structure, there might be more relations than just those in \eqref{Carnotbasis}, as, for example, in the first quaternionic Heisenberg group; see page \pageref{quaternionic_Heisenberg}. In \cite{LeDonne_Tripaldi}, there are several other uses of these diagrams to represent Carnot groups in low dimensions. 

\subsection{Examples of Carnot groups and Carnot algebras}\label{examples_Carnot}
Let $G$ be a Carnot group with Lie algebra $\g$.
Given a basis $X_1,\ldots, X_n$ of $\mathfrak{g}$,
 we will use the identification
\begin{eqnarray}\label{identification}
\mathbb{R}^n&\longleftrightarrow& G \nonumber\\
 (x_1,\ldots,x_n)&\longmapsto&\exp\Big(\sum_{i=1}^nx_iX_i\Big)\,.
\end{eqnarray}
This identification
allows us to write the group product using Dynkin product \eqref{Dynkin Formula}. In fact, for all $\mathbf{x},\mathbf{y}\in\mathbb{R}^n$, there exists a unique $\mathbf{z}\in\mathbb{R}^n$ such that\index{Dynkin product} 
\begin{equation}\label{group law in G}
 \exp\Big(\sum_{i=1}^nx_iX_i\Big)\exp\Big(\sum_{i=1}^ny_iX_i\Big)=\exp\Big(\sum_{i=1}^nz_iX_i\Big)\,.
\end{equation}
Via the identification (\ref{identification}), one can write the group law in (\ref{group law in G}) as
\begin{equation}
 \mathbf{x}\star\mathbf{y}=\mathbf{z}\,.
\end{equation}
Hence, we have a group law $\star$ on $\mathbb{R}^n$ that makes $(\mathbb{R}^n,\star)$ a simply connected Lie group with Lie algebra $\mathfrak{g}$, whose identity element is $\mathbf{0}$.
\subsubsection{Carnot groups in dimension 4}

In dimension 4 there are only the following nilpotent simply connected Lie groups, which are all stratifiable:
$\R^4$, $\R\times N_3$, and $N_{4,2}$, where $N_3$ is the Heisenberg group and $N_{4,2}$ is the {\em Engel group} as we now recall.
 The {\em Engel Lie algebra} is spanned by 4 vectors $X_1, \ldots, X_4$ with 
 only non-trivial brackets 
\begin{eqnarray}\label{Engel_rel}
 [X_1,X_2]=X_3\,,\,[X_1,X_3]=X_4\,.
\end{eqnarray}
This is a nilpotent Lie algebra of rank 2 and step 3 that is stratifiable. It is also known as the filiform Lie algebra of dimension 4; see Exercise~\ref{ex_fili1}. The Lie brackets can be pictured with the diagram:
\begin{center}
 \begin{tikzcd}[end anchor=north]
 X_1\ar[dr, no head]\ar[ddr, no head] & &X_2\ar[dl, no head]\\
 &X_3\ar[d, no head] & \\
 & X_4 & \quad\;.
 \end{tikzcd}
\end{center}



The Lie group $N_{4,2}$ is the only simply connected Lie group with such a Lie algebra. The group law (\ref{group law in G}) of the Engel group in exponential coordinates is given by:
\begin{equation}\label{prod_Engel_2024}
 \mathbf{x}\star\mathbf{y}=\mathbf{z} \quad\Longleftrightarrow\quad 
\begin{cases}
 z_1=x_1+y_1 &\\
 z_2=x_2+y_2 &\\
 z_3=x_3+y_3+\frac{1}{2}(x_1y_2-x_2y_1) &\\
 z_4=x_4+y_4+\frac{1}{2}(x_1y_3-x_3y_1)+\frac{1}{12}(x_1-y_1)(x_1y_2-x_2y_1) &
\end{cases}.
\end{equation}


\subsubsection{Cartan group $\mathbb F_{2,3}$}\label{free 5-dim}\index{Cartan! -- group} 
 {\em Cartan Lie algebra} is another name for the free-nilpotent Lie algebra of step 3 and 2 generators. The simply connected Lie group with this Lie algebra is known as the {\em Cartan group}, and it is sometimes denoted by $\mathbb F_{2,3}$ or by $N_{5,2,3}$ as in \cite{Gong_Thesis}. 
 With respect to some basis $X_1, \ldots, X_5$, the non-trivial brackets are the following:
\begin{eqnarray*}\label{N523}
 [X_1,X_2]=X_3\,,\,[X_1,X_3]=X_4\,,\,[X_2,X_3]=X_5\,.
\end{eqnarray*}
The Lie brackets can be pictured with the diagram:
 
 \begin{center}
 
	\begin{tikzcd}[end anchor=north]
		X_1 \ar[dr, no head]\ar[no head, dd] & & X_2\;\;\ar[dl, no head]\ar[no head,dd, ->-=.5, end anchor={[xshift=-3.9ex]north east},start anchor={[xshift=-3.9ex]south east }] \\
		& X_3\ar[no head, dl]\ar[no head, dr, -<-=.5, end anchor={[xshift=-3.9ex]north east}] & &\\
		 X_4 & & X_5\;.
	\end{tikzcd}
	\end{center}
	


 The group law (\ref{group law in G}) of the Cartan group in exponential coordinates is given by:
 \begin{equation}\label{prod_Cartan_2024}
 \mathbf{x}\star\mathbf{y}=\mathbf{z} \quad\Longleftrightarrow\quad 
\begin{cases}
 z_1=x_1+y_1&\\
 z_2=x_2+y_2&\\
 z_3=x_3+y_3+\frac{1}{2}(x_1y_2-x_2y_1)&\\
 z_4=x_4+y_4+\frac{1}{2}(x_1y_3-x_3y_1)+\frac{1}{12}(x_1-y_1)(x_1y_2-x_2y_1)&\\
 z_5=x_5+y_5+\frac{1}{2}(x_2y_3-x_3y_2)+\frac{1}{12}(x_2-y_2)(x_1y_2-x_2y_1)&\\
\end{cases}.
\end{equation}


\subsubsection{Filiform groups}

Filiform groups are those simply connected Lie groups whose Lie algebra is filiform, in the sense that the Lie algebra is dimensionwise the smallest among those Lie algebras with the same nilpotency step. 
We saw in Example~\ref{ex_fili1} those of the first kind. 
 With respect to some basis $X_1,\ldots,X_{s+1}$, with $s\in \N$, the non-trivial brackets of the $(s+1)$-dimensional example of the first kind are the following:
 \begin{equation*}
 [X_1,X_i]=X_{i+1}\,,\qquad \text{ for } i\in \{2,\ldots, s\}.
 \end{equation*}
This is a nilpotent Lie algebra of rank 2 and step $s$ that is stratifiable.
The Lie brackets can be pictured with the diagram:
\begin{center}
 
	\begin{tikzcd}[end anchor=north]
		 X_1\ar[dr, no head]\ar[ddr, no head]\ar[dddr, no head]\ar[ddddr, no head]\ar[dddddr, no head] & & X_2\ar[dl, no head]\\
		 & X_3\ar[d, no head] & \\
		 & X_4\ar[d, no head] &\\
		 & \vdots\ar[d, no head]&\\
		 & \quad X_{s}\ar[d, no head] & \\
		 & X_{s+1} &\quad\;.
		 \end{tikzcd}
	\end{center}
	
	In \cite{Vergne} M. Vergne classified all filiform Lie algebras. In addition to the ones of the first type, 
	in each even dimension starting from dimension 6, there is exactly one more filiform Lie algebra, called {\em filiform Lie algebra of the second kind}.\index{filiform! -- algebra of the second kind} 
	In dimension $6$, we have the example discussed in Exercise~\ref{ex nontrivial filiform}.
	With respect to that basis $ Y_0,Y_1,Y_2, Y_3,Y_4, Y_5 $, the bracket diagram is the following:
	\begin{center}
	\begin{tikzcd}[end anchor=north]
		Y_0\drar[no head] \ar[ddr,no head, end anchor={[xshift=-2.5ex]north east}]\ar[dddr,no head, end anchor={[xshift=-2.5ex]north east}]
		\ar[ddddr,no head, end anchor={[xshift=-2.5ex]north east}]
		& & Y_1\;\,
		\dlar[no head]
		\ar[ddddl,bend right=-20,swap, 
		 no head, end anchor={[xshift=.5ex]north east}]
		 \\
		& Y_{2}\dar[no head, end anchor={[xshift=-2.5ex]north east},start anchor={[xshift=-2.5ex]south east }]
		\ar[ddd, bend right=-30,swap,->-=.5, no head, end anchor={[xshift=-1ex]north east}]& \\
		&Y_3
		\dar[no head, end anchor={[xshift=-2.5ex]north east},start anchor={[xshift=-2.5ex]south east }]
		\ar[dd,bend right=-30,swap, no head, -<-=.5,end anchor={[xshift=-1ex]north east}]& \\
		&Y_4\dar[no head, end anchor={[xshift=-2.5ex]north east},start anchor={[xshift=-2.5ex]south east }]
		\dar[
		bend right=-20,swap,no head, end anchor={[xshift=.5ex]north east}]& \\
		& Y_5 & \; .
	\end{tikzcd}
\end{center}	
Another presentation (see $N_{6,2,2}$ in \cite{Gong_Thesis}), in a basis $X_1, \ldots, X_6$ is given by with the diagram:
\begin{center}
	\begin{tikzcd}[end anchor=north]
		X_1\drar[no head] \ar[ddr,no head, end anchor={[xshift=-2.5ex]north east}]\ar[dddr,no head, end anchor={[xshift=-2.5ex]north east}]& & X_2\;\,\dlar[no head]\ar[dddd, ->-=.5, no head, end anchor={[xshift=-4.5ex]north east}] \\
		& X_{3}\dar[no head, end anchor={[xshift=-2.5ex]north east},start anchor={[xshift=-2.5ex]south east }]\ar[dddr, bend right=-10,swap,
		 no head, end anchor={[xshift=-3ex]north east}]& \\
		&X_4\dar[no head, end anchor={[xshift=-2.5ex]north east},start anchor={[xshift=-2.5ex]south east }]\ar[ddr, no head,bend right=-10,swap, 
		end anchor={[xshift=-3ex]north east}]& \\
		&X_5\drar[-<-=.5,no head,end anchor={[xshift=-4.5ex]north east}]& \\
		& & X_6\; .
	\end{tikzcd}
\end{center}


\subsubsection{The second Heisenberg group} 
\index{Heisenberg! -- group} 
The Lie brackets of the second Heisenberg algebra can be pictured with the diagram:
 
 \begin{center}
 
	\begin{tikzcd}[end anchor=north]
		X_1 \ar[no head, drr,end anchor={[xshift=-3.ex]north east}] & X_2 \ar[no head, dr,end anchor={[xshift=-3.ex]north east}] & X_3\ar[no head, d, end anchor={[xshift=-2.ex]north east},start anchor={[xshift=-2.ex]south east}] & X_4\ar[no head, dl, end anchor={[xshift=-2.ex]north east}]\\
 & & X_5 &\quad\;.
	\end{tikzcd}
	\end{center}
	We will encounter this group again in Chapter~\ref{ch_ROSS} discussing complex hyperbolic spaces; see 
	Exercise~\ref{n_Heisenberg_A}. It will be called the 
{\em $2$-nd $\mathbb{C}$-Heisenberg group} and will be denoted by $\mathcal{N}^\mathbb{C}_{2}$. 
	
	\subsubsection{The first quaternionic Heisenberg group}\label{quaternionic_Heisenberg} 
	The Lie algebra of the first quaternionic Heisenberg group 
 $\mathcal{N}^\mathbb{H}_{1}$ 
	can be characterized as the only 7D Carnot algebra of rank 4 and step 2 where every element in the first stratum has maximal rank, i.e., for every nonzero $X$ in the first stratum, the map $\ad_X$ is surjective. 
The Lie brackets can be pictured with the diagram:
\begin{center}
 
	\begin{tikzcd}[end anchor=north]
			X_1\ar[ddrr, no head,end anchor={[xshift=-3.5ex]north east},start anchor={[xshift=-1.3ex]south east},start anchor={[yshift=.5ex]south east}]\ar[ddr, no head,end anchor={[xshift=-3.3ex]north east},start anchor={[xshift=-3.3ex]south east}]\ar[ddrrr, no head,end anchor={[xshift=-1.6ex]north east},end anchor={[yshift=-.1ex]north east}] & X_2\ar[dd,no head,start anchor={[xshift=-3.8ex]south east},start anchor={[yshift=.7ex]south east},end anchor={[xshift=-3.3ex]north east}]\ar[ddr, no head, end anchor={[xshift=-1.5ex]north east},start anchor={[yshift=.7ex]south east},start anchor={[xshift=-2.ex]south east}] \ar[ddrr,-<-=.3, no head, end anchor={[yshift=-.5ex]north east},end anchor={[xshift=-2.9ex]north east}, start anchor={[xshift=-0.7ex]south east},start anchor={[yshift=1.9ex]south east}]& & X_3\ar[dd, no head,end anchor={[yshift=-.5ex]north east},end anchor={[xshift=-2.9ex]north east}, start anchor={[yshift=.5ex]south east},start anchor={[xshift=-2.9ex]south east},->-=.7]\ar[ddl, no head,end anchor={[xshift=-3.5ex]north east},start anchor={[xshift=-3.5ex]south east},start anchor={[yshift=.8ex]south east}]\ar[ddll,,-<-=.2,no head, end anchor={[xshift=-1.9ex]north east}] & X_4\ar[ddl, no head,end anchor={[xshift=-1.6ex]north east},start anchor={[xshift=-2.2ex]south east},start anchor={[yshift=.7ex]south east},end anchor={[yshift=-.1ex]north east}]\ar[ddll,no head,start anchor={[xshift=-3.5ex]south east},start anchor={[yshift=0.9ex]south east},end anchor={[xshift=-1.5ex]north east}]\ar[ddlll,->-=.4, no head,end anchor={[xshift=-1.9ex]north east}]\\
		& & & &\\
		 &X_5 & X_6 & X_7 &\quad\;.
	\end{tikzcd}
	\end{center}

	\subsubsection{Some free-Carnot groups: $\mathbb F_{n2}$, $\mathbb F_{24}$, $\mathbb F_{25}$, and $\mathbb F_{33}$}
	Free-nilpotent Lie algebras are stratifiable. Hence, the associated simply connected Lie groups 
	are the free objects in the category of 
	Carnot groups. Here are some examples:
	
	The free-Carnot group $\mathbb F_{n2}$ of rank $n$ and step 2 has been discussed in Example~\ref{free_nilp_step2}, recall also Example~\ref{ex_group_step2_20jun}.
	Setwise $\mathbb F_{n2} := \Lambda ^1 (\R^n) \oplus \Lambda ^2 (\R^n) $, while the group law is
$ x\cdot y = x+y+ \frac12 x\wedge y$.

The Lie group	
	$\mathbb F_{24}$
 is the Carnot group whose Lie algebra is free-nilpotent with 2 generators and nilpotency step 4. It has dimension 8.
The non-trivial brackets in some basis $X_1,\ldots, X_8$ are:
\begin{equation*}
 \begin{aligned}
 &\qquad[X_1,X_2]=X_{3}\,,\,[X_1,X_3]=X_4\,,\,[X_2,X_3]=X_{5}\,,\\ &[X_1,X_4]=X_6\,,\,[X_1,X_5]= [X_2,X_4]=X_7\,,\,
 [X_2,X_5]=X_8\,.
 \end{aligned}
 \end{equation*}

 The Lie group $\mathbb F_{25}$
 is the Carnot group whose Lie algebra is free-nilpotent with 2 generators and nilpotency step 5. It has dimension 14.
The non-trivial brackets in some basis $X_1,\ldots, X_{14}$ are:
\begin{equation*}
 \begin{aligned}
 &\;[X_1,X_2]=X_{3}\,,\,[X_1,X_3]=X_4\,,\,[X_2,X_3]=X_{5}\,,\,[X_1,X_4]=X_6\,,\\
 &
 [X_2,X_5]=X_8\,,\,[X_1,X_5]= [X_2,X_4]=X_7\,,\,[X_1,X_7]=X_{10}+X_{13}\,,\\
 &\qquad[X_1,X_6]=X_9\,,\,[X_1,X_8]=X_{11}+X_{14}\,,\,[X_2,X_6]=X_{10}\,,\\&[X_2,X_7]=X_{11}\,,\,[X_2,X_8]=X_{12}\,,\,[X_3,X_4]=X_{13}\,,\,[X_3,X_5]=X_{14}\,.
 \end{aligned}
 \end{equation*}

 The Lie group
 $\mathbb F_{33}$
 is the Carnot group whose Lie algebra is free-nilpotent with 3 generators and nilpotency step 3. It has dimension 14.
The non-trivial brackets in some basis $X_1,\ldots, X_{14}$ are:
\begin{equation*}
 \begin{aligned}
 &\;\quad[X_1,X_2]=X_{4}\,,\,[X_1,X_3]=X_5\,,\,[X_2,X_3]=X_{6}\,,\,[X_1,X_4]=X_7\,,\\ &\quad[X_1,X_5]=X_8\,,\,[X_1,X_6]=X_9\,,\, [X_2,X_4]=X_{10}\,,\,
 [X_2,X_6]=X_{11}\,,\\
 &[X_3,X_4]=X_{12}\,,\,[X_3,X_5]=X_{13}\,,\,[X_3,X_6]=X_{14}\,,\,[X_2,X_5]=X_9+X_{12}\,.
 \end{aligned}
 \end{equation*}

Carnot groups are completely classified up to dimension 7. However, in dimension 7, there are infinitely many non-isomorphic classes. 
For a cornucopia of examples, we refer to \cite{LeDonne_Tripaldi}.	

\section{Simple consequences of dilations}\label{sec consequences of dilations}
\subsection{Canonical coordinates}

 The next proposition holds in the general setting of nilpotent simply connected Lie groups; see Theorem~\ref{Malcev global coordinates}.
We will give a simplified proof for Carnot groups.
 
\begin{proposition}\label{prop66b1e96e}
 Let $G$ be a Carnot group with dilations $(\delta_\lambda)_{\lambda\in\R}$.
 Let $W_1\oplus\cdots\oplus W_m = \Lie(G)$ be a direct-sum decomposition by dilation invariant subspaces.
 Let
 \[
 \begin{aligned}
 \Psi:\g\simeq& W_1\times\cdots\times W_m &&\to &&G, \\
 &(X_1, \dots, X_m) &&\mapsto &&\prod_{j=1}^m \exp(X_j) .
 \end{aligned}
 \]
 Then, the map $\Psi:\g\to G$ gives a global coordinate system.
 In particular, on every Carnot group, Malcev coordinates and exponential coordinates exist, globally.
\end{proposition}
Above, like in other places in the book, we used the notation
\[
\prod_{j=1}^m g_j := g_1\cdots g_m .
\]

\begin{proof}[Proof of Proposition~\ref{prop66b1e96e}]
 The linear map $(\dd\Psi)_0$ is an isomorphism, because
 \[
 (\dd\Psi)_0 X = X,
 \qquad \forall j\in\{1, \dots, m\},\,\forall X\in W_j .
 \]
 Hence, the map $\Psi:\g\to G$ is a diffeomorphism between some neighborhood of $0$ in $\g$ and some neighborhood of $1_G$ in $G$.
 By assumption, we have
 \begin{equation}\label{eq66b1ea98}
 \delta_\lambda W_j = W_j,
 \qquad\forall \lambda>0,\forall j\in\{1, \dots, m\}.
 \end{equation}
 Consequently, we claim
 \begin{equation}\label{eq66b1eae1}
 \Psi\circ\delta_\lambda = \delta_\lambda\circ\Psi,
 \qquad\forall\lambda\in\R .
 \end{equation}
 Indeed, for all $(X_1,\cdots,X_m)\in W_1\times\cdots\times W_m$, we have
 \begin{eqnarray*}
 \Psi(\delta_\lambda(X_1+\cdots+X_m))
 &=& \Psi( \delta_\lambda X_1+\cdots+ \delta_\lambda X_m ) \\
 &\stackrel{\eqref{eq66b1ea98}}{=}& \prod_{j=1}^m \exp(\delta_\lambda X_j) \\
 &=& \prod_{j=1}^m \delta_\lambda\exp( X_j) \\
 &=& \delta_\lambda \prod_{j=1}^m \exp( X_j) \\
 &=& \delta_\lambda\Psi(X_1+\cdots+X_m) .
 \end{eqnarray*}
 By~\eqref{eq66b1eae1}, we deduce that $\Psi$ is a global diffeomorphism.
 Regarding Malcev coordinates, recall that we proved the existence of Malcev basis in Section~\ref{sec Good bases for Carnot} (and more generally in Proposition~\ref{Prop exists Malcev bases}).
\end{proof}

Group laws in exponential coordinates with respect to Carnot bases have triangular forms.
This result follows from a more general result about nilpotent simply connected Lie groups.
Indeed, since Carnot bases are order-reversed Malchev bases, we have the following consequence of Proposition~\ref{CG-1.2.9}. We stress that Proposition~\ref{CG-1.2.9} was just a consequence of the BCH formula; thus, the reader should try to prove Proposition~\ref{Group law in Carnot coordinates} directly, as an exercise. 

\begin{proposition}\label{Group law in Carnot coordinates}
Let $G$ be a Carnot group equipped with a Carnot basis $(X_1, \dots, X_n)$.
On $G$, consider exponential coordinates of the first or second kind associated with the basis.
Then, in these coordinate systems, the product law has a lower triangular form:
$$(s_1, \dots, s_n)\cdot (t_1, \dots, t_n)= s+t +\sum_{j=1}^n Q_j(s,t) e_j, \qquad \forall s,t \in \R^n,$$
where each $Q_j$ is a polynomial that is not depending on $s_{j+1}, \dots, s_n$ nor on $t_{j+1}, \dots, t_n$.
\end{proposition}
\begin{proof}
We have that the reversed-ordered basis $(X_n, \dots, X_1)$ is a Malcev basis.
Hence, we invert the order of the coordinates: $(s_n, \dots, s_1)$.
Then Proposition~\ref{CG-1.2.9} gives the result, replacing upper triangular with lower triangular.
\end{proof}

\subsection{A proof of Ball-Box Theorem for Carnot groups}

Let $G$ be a Carnot group with stratification $V_1, \dots, V_s$.
Let $X_1, \dots, X_n$ be a basis of $\Lie(G)$ adapted to the stratification, hence, for all $j$ let $d_j$ such that $X_j\in V_{d_j}$.
The {\em boxes} with respect to the numbers $(d_1, \ldots, d_n)$ are 
\begin{equation}\label{def_box_2024}
\Bx(r) := \left\{(t_1, \dots, t_n)\in\R^n\;:\;|t_j|<r^{d_j}\right\}, \qquad \forall r\geq 0.
\end{equation}
The {\em dilations} in $\R^n$ with respect to the numbers $(d_1, \ldots, d_n)$ are the maps $\delta_\lambda:\R^n\to\R^n$ defined as
\begin{equation}\label{def_delta_2024}
\delta_\lambda(t_1, \dots, t_n) := (\lambda t_1, \dots, \lambda^{d_j}t_j, \dots, \lambda^s t_n), \qquad \forall {(t_1, \dots, t_n)} \in \R^n, \forall \lambda\in \R.
\end{equation}

Let $\Phi:\R^n\to G$ be the exponential coordinate map with respect to the basis $X_1, \dots, X_n$, i.e., $\Phi({\bf{t}}):=\exp(\sum_jt_jX_j)$.
Given $p\in G$, the {\em exponential coordinate map from $p$ with respect to $X_1, \dots, X_n$} is $\Phi_p		:= L_p\circ \Phi$.

\begin{theorem}[Ball-Box for Carnot groups]\label{Ball-Box4Carnot}
 	Let $G$ be a Carnot group.
	Fix a basis adapted to the stratification.
	Then there is $C>1$ such that for all $p\in G$ and all $r>0$
\begin{equation}\label{B-Box}
	B(p,\frac\lambda C) \subset \Phi_p(\Bx(\lambda)) \subset B(p,\lambda C), 
\end{equation}
	where $\Phi_p$ is the exponential coordinate map from $p$ with respect to the basis, and balls are with respect to the Carnot-Carath\'eodory distance.
\end{theorem}
\begin{proof}
Each set $\Phi(\Bx(r))$, with $r>0$ is a bounded neighborhood of $1_G$ in $G$, recall Proposition~\ref{prop66b1e96e}.
Let $\dcc$ be the Carnot-Carath\'eodory distance function on the Carnot group $G$.
By Chow Theorem~\ref{Chow2}, or alternatively see Corollary~\ref{corol_weak_BB}, the distance $\dcc$ induces the manifold topology.
Hence, there is $r_0,r_1>0$ such that
\[
B(1_G,r_1) \subset \Phi(\Bx(r_0)) \subset B(1_G,1), 
\] 
Recalling that $\delta_\lambda(B(1_G,r)) = B(1_G,\lambda r)$, by Proposition~\ref{dilations_on_Carnot2024}, and applying $\delta_\lambda$, we get
\[
B(1_G, \lambda r_1) \subset \delta_\lambda\Phi(\Bx(r_0)) \subset B(1_G,\lambda ). 
\]
Moreover, we have
\begin{align*}
 	\delta_\lambda\left(\Phi(\Bx(r_0)) \right)
	&= \delta_\lambda\left( \Phi\left\{(t_1, \dots, t_n):|t_j|<r_0^{d_j}\right\} \right) \\
	&= \delta_\lambda\left\{\exp\left(\sum_jt_j X_j\right) : |t_j|<r_0^{d_j} \right\} \\
	&= \left\{\exp\left(\delta_\lambda\sum_jt_j X_j\right) : |t_j|<r_0^{d_j} \right\} \\
	&= \left\{\exp\left(\sum_j \lambda^{d_j}t_j X_j\right) : |t_j|<r_0^{d_j} \right\} \\
	&= \left\{\exp\left(\sum_j s_j X_j\right) : |s_j|<\lambda^{d_j}r_0^{d_j} \right\} \\
	&= \Phi(\Bx(\lambda r_0 )).
\end{align*}
Therefore, we deduce that
\begin{equation}\label{eq291244}
 B(1_G, \lambda r_1) \subset \Phi(\Bx(\lambda r_0)) \subset B(1_G,\lambda ), \qquad \forall \lambda>0.
\end{equation}
	Since $\dcc $ is left-invariant, applying $L_p$ to \eqref{eq291244}, we obtain \eqref{B-Box}
	for all $p\in G$ and all $\lambda>0$.
\end{proof}

\subsection{Canonical measures}

In every Carnot group, there are few natural choices of measures:
Haar, Hausdorff of top dimension, and Lebesgue measures in exponential coordinates. 
We introduced the Haar measures in Section~\ref{sec_Haar_poly_growth}, the Hausdorff measures in Section~\ref{subsec_Hausdorff}, and the exponential coordinates in Section~\ref{sec Canonical coordinates}.
In this section, we will prove that, up to a scalar factor, these measures are the same, and show some other properties. 

Let $G$ be a Carnot group. 
For every $k>0$, let ${\mathscr H}^k$ be the $k$-dimensional Hausdorff measure. We can also consider the $k$-dimensional spherical Hausdorff measure $\S^k$; see \cite[page~75]{Mattila}.
Because left-translations are isometries, the measures ${\mathscr H}^k$ and $\S^k$ are left-invariant.
We shall see that, when $k$ equals the homogeneous dimension $Q$, these measures are Radon measures, and therefore, they are Haar measures.

If we consider exponential coordinates $\exp: \R^n \to G$, linearly identifying the Lie algebra $\g$ of $G$ with $\R^n$ via the choice of a basis, then we can push the Lebesgue measure denoted by $\LL^n$ from $ \R^n$ to $G$.
We shall see that this measure, which obviously is a Radon measure, is left invariant and right invariant. Hence, it is a bi-invariant Haar measure.

When we fix a Haar measure on $G$, we denote it by $\vol$, or by $\vol_{G}$ if more than a group is considered.
As discussed in Section~\ref{sec_Haar_poly_growth}, every pair of Haar measures differs by a constant.
Hence, the measures ${\mathscr H}^Q$, $\vol$, and $\LL^n$ are a multiple of each other. 

Moreover, Carnot groups are unimodular, in the sense that left-Haar measures are right-Haar measures, and vice versa.
This latter fact holds for all nilpotent Lie groups; see Theorem~\ref{Malcev global coordinates}.

\begin{definition}[Homogeneous dimension for Carnot groups]\label{def Homogeneous dimension Carnot}
 	If $G$ is a Carnot group and $V_1, \dots, V_s$ is the stratification of its Lie algebra, we call {\em homogeneous dimension} of $G$ the integer number
	\[
	Q := \sum_{j=1}^s j\cdot\dim V_j .
	\]
\end{definition}

\begin{proposition}\label{prop_vol_Carnot}
 	Let $G$ be a Carnot group of homogeneous dimension $Q$.
	\begin{description} 
	\item[\ref{prop_vol_Carnot}.i.] 	
	If $\vol$ is a Haar measure of $G$, then 
		\[
		\vol(B(p,r)) = r^Q \vol(B({1_G},1)) .
		\]
	\item[\ref{prop_vol_Carnot}.ii.] 	
	Every Haar measure of $G$ is Ahlfors $Q$-regular, the Hausdorff dimension of $G$ is $Q$, and the Hausdorff $Q$-measure is a Haar measure.
	\item[\ref{prop_vol_Carnot}.iii.] 	
	In exponential coordinates, the Lebesgue measure is the Hausdorff $Q$-measure up to a multiplication by a constant.
	\end{description}
\end{proposition}
\begin{proof}
By Theorem~\ref{Malcev global coordinates} (see also Proposition~\ref{CG-1.2.9}), in exponential coordinates, the Lebesgue measure $\LL^n$ is both left and right-invariant. Moreover, every other Haar measure is a multiple of it; see \cite[Theorem~11.9]{Folland_book}.
	In exponential coordinates, the inhomogeneous dilations $\delta_\lambda$ have Jacobian $\lambda^Q$, i.e., 
	 $\LL^n(\delta_\lambda(E))=\lambda^Q \cdot\LL^n(E)$, for every measurable set $E\subseteq \R^n$.
	Hence 
	\[
	\LL^n(B(p,\lambda)) =\LL^n(B(1_G,\lambda)) = \LL^n(\delta_\lambda(B(1_G,1))) = \lambda^Q \LL^n(B(1_G,1))
	\]
	By Theorem~\ref{thm: Ahlfors: regular: Hausdorf: dim}, or Corollary~\ref{Corol_Ahlfors_Hausdorff}, the metric measure space $(G, \vol)$ is Ahlfors regular of dimension $Q$. In particular, the Hausdorff $Q$-measure ${\mathscr H}^Q$ is Radon, and hence, it is a Haar measure.
	The last part of the proposition follows since both $\LL^n$ and ${\mathscr H}^Q$ are Haar measures.
\end{proof}


\section{Pansu-Rademacher Theorem}\label{sec Pansu-Rademacher Theorem}
We would like to observe that the classical Rademacher Theorem states not only the almost-everywhere existence of a tangent map (called the differential) but also its realizability as a linear map, meaning as a group homomorphism that is compatible with the respective groups of dilations. 
Expressed in these terms, the theorem holds for general equiregular sub-Finsler manifolds as well; see \cite{Margulis-Mostow}.
The aim of this section is to explain the content of such a differentiability result and to give a complete proof of it in the case of Carnot groups.


\subsection{Pansu's theorem}
We shall prove Pansu's version of the Rademacher Theorem. 
\begin{definition}[Pansu differentiability]\index{Pansu! -- differential}
Let $G$ and $H$ be Carnot groups. We denote by $\delta_h$ the dilations of factor $h $ in both of the groups.
 If $f: G\to H$ is a map, then its {\em Pansu differential} at a point $x\in G$ is, if it exists, the limit
$$ Df_x:= \lim_{h\to0^+} \delta_{1/h}\circ L_{f(x)}^{-1} \circ f \circ L_{x}\circ \delta_{h},$$
where the limit is with respect to the convergence on compact sets.
 Moreover, we say that $f$ is {\em Pansu differentiable} if $ Df_x$ exists and is a homogeneous group homomorphism.
 \end{definition}
For a map $f:G\to H$ between Carnot groups, the value $Df(x;v):=\lim_{h\to0^+} \delta_{1/h}(f(x)^{-1} f (x \delta_{h}(v)))$ may be called {\em partial Pansu derivative} of $f$ at $x\in G$ along $v\in G$. 
Notice that if $Df(x;v)$ exists, then $Df(x;\delta_\lambda v)$ exists for all $\lambda >0$ and $Df(x;\delta_\lambda v) = \delta_\lambda Df(x;v)$.
The value $Df(x;v)$ may exist for all $v\in G$, but the limit may not be uniform on compact sets.
Moreover, even if the map $ Df_x$ exists, it may not be a group homomorphism from $G$ to $H$.

\begin{theorem} [Pansu's generalization of Rademacher Theorem]\label{Thm:Pansu:Rademacher}\index{Pansu! -- Rademacher Theorem}\index{Theorem! Pansu Rademacher --}
Let $f: G\to H$ be a Lipschitz map between sub-Finsler Carnot groups.
Then, for almost every $x\in G$, the map $f$ is Pansu differentiable at $x$.
 \end{theorem}

 \subsubsection{Preliminaries to the proof of Pansu's theorem}
 In the proof of Theorem~\ref{Thm:Pansu:Rademacher}, we will only take for granted a few classical results to which we give hints to the proofs and references in the exercise section.
 \begin{theorem}[Rademacher Theorem in 1D; see {\cite[Section 3.5]{Folland_book}}]\label{Rademacher1D}
If $\gamma :[0,1] \to \R^n$ is Lipschitz with respect to the Euclidean distance on $\R^n$, then the derivative $\dot \gamma(t)$ exists for almost every $t$ and
$$\gamma(t)=\gamma(0) + \int_0^t \dot \gamma(s)\d s, \qquad \text{ for all } t\in [0,1].$$
 \end{theorem}

\begin{theorem}[Egorov Theorem for metric spaces; see Exercise~\ref{Ex:Egorov}]\label{Egorov's}
Let $(X,\mu)$ be a measure space with $\mu (X)<\infty$ and let $Y$ be a separable metric space.
Let $(f_t)_{t>0}$ be a family of measurable functions from $X$ to $Y$ depending on $t\in \R_{>0}$.
Suppose that $(f_t)_t$ converges almost everywhere to some $f$, as $t\to 0$.
Then for every $\eta > 0$, there exists a measurable subset $K\subset X$ such that the $\mu(\Omega \setminus K)<\eta$ and $(f_t)_t$ converges to $f$ uniformly on $K$.
\end{theorem}

\begin{theorem}[Consequence of Lebesgue Differentiation Theorem for doubling metric spaces; see Exercise~\ref{Ex:LDT}]\label{LDT}
 If $(X,d,\mu)$ is a doubling metric measure space and $K$ is a measurable set in $X$ then $\mu$-almost every point of $K$ has density 1, that is,
$$\lim_{r\to 0^+} \frac{\mu(K\cap B(x,r) )}{\mu( B(x,r) )}=1,\qquad \text{ for a.e. } x\in K.$$
\end{theorem}
 
\subsubsection{A proof of Pansu's theorem} 

As in Pansu's original proof, we first deal with the case of curves. We shall prove that every Lipschitz curve into a Carnot group is Pansu differentiable almost everywhere.

\begin{proposition}[Case of curves]\label{prop:diff:curves}
Let $G$ be a Carnot group and $\gamma:[0,1] \to G$ a Lipschitz curve.
Then $\gamma$ is Pansu differentiable almost everywhere and for almost every $x\in [0,1]$ we have that for all $v\in \R$
$$D\gamma( x ; v) := \lim_{t\to 0} \delta_{1/t} \left(\gamma( x )^{-1} \gamma( x + t v)\right)=
\exp\left(v (L_{\gamma( x )})^* \dot \gamma( x )\right).$$
\end{proposition}
Here are a few remarks before the proof. First, we notice that the above curve $\gamma$ is, in particular, Euclidean Lipschitz, so the tangent vector $\dot \gamma( x )$ exists for almost every $ x $ by Theorem~\ref{Rademacher1D}.
We also stress that Pansu's differentiability for curves is stronger than Euclidean differentiability. Namely, if we consider the curve in a rank-$r$ Carnot group in exponential coordinates $\gamma(t)=(\gamma_1(t), \ldots, \gamma_n(t))$ with $\gamma( 0 )=0$ and $ 0 $ is a point of Euclidean differentiability for $\gamma$,
then 
$\dot \gamma(0) = \lim_{t\to 0} \gamma(t)/t =
\lim_{t\to 0} (\gamma_1(t)/t, \ldots, \gamma_n(t)/t)
= (h_1, \ldots, h_r, 0,\ldots, 0)$.
However, we have to consider 
$$\delta_{1/t} \gamma(t) = (\gamma_1(t)/t, \ldots, \gamma_n(t)/t^s)$$
and we need to prove that every coordinate $\gamma_j(t)$, with $j$ greater than the rank, in fact, vanishes not just faster than $t$ but faster than $t$ to the power of the degree of the coordinate.

\proof[Proof of Proposition~\ref{prop:diff:curves}]
For simplicity, we take $v=1$. We take a basis $X_1,\ldots, X_r$ of the first layer of the stratification of $\Lie(G)$.
Let $h_1, \ldots, h_r \in L^\infty([0,1];\R)$ be such that 
\begin{equation}\label{Boston20172}\dot \gamma(t)=\sum_{j=1}^r h_j(t) X_j(\gamma(t)), \qquad \text{ for almost all } t\in [0,1].\end{equation}
Since $\gamma$ is $L$-Lipschitz, we may take $|h_j(t)|\leq L$, for all $t$.
Let $ x \in [0,1]$ be both a point of Euclidean differentiability for $\gamma$ and a Lebesgue point for all $h_j$, i.e., 
$$\dfrac{1}{|t- x |}\int_ x ^t |h_j(s) - h_j(t)| \dd s \to 0, \quad \text{ as } t\to x .$$
Up to replacing $\gamma$ with the curve $t\mapsto \gamma( x )^{-1} \gamma(t+ x )$, we may assume that $ x =0$ and $\gamma( x )=0$.

We identify the group $G$ with its Lie algebra via the exponential map. Our aim is now to show that
$$\lim_{t\to0} \delta_{1/t} \gamma(t) = \dot \gamma(0),$$
where the latter equals $ \sum_{j=1}^r h_j(0) X_j(0)$ since 0 is a Lebesgue point for all $h_j$.

Set $\eta_t (s) := \delta_{1/t} \gamma(t\, s)$, so each $\eta_t:[0,1]\to G$ is a curve starting at 0 that is $L$-Lipschitz:
$$d( \eta_t(s), \eta_t(s') )
=
d( \delta_{1/t} \gamma(t\, s), \delta_{1/t} \gamma(t\, s') )
\leq
\dfrac{L}{t} |ts - ts'|= 
L|s-s'|.$$
Consequently, every sequence $(\eta_{t_k})_k$ has a uniformly converging subsequence.
Moreover, we claim we have the equality
\begin{equation}\label{Boston20171}
\dot\eta_t (s) = \sum_{j=1}^r h_j(ts) X_j(\eta_t(s)).
\end{equation}
Indeed,
\begin{eqnarray*}
        \dot\eta_t(s)
        &=& \frac{\dd}{\dd s} \delta_{1/r}\gamma(ts)
        = (\dd\delta_{1/t})_{\gamma(ts)} (t\dot\gamma(ts)) \\
        &\overset{\eqref{Boston20172}}{=}& t \sum_{j=1}^r h_j(ts) (\dd\delta_{1/t})_{\gamma(ts)} (X_j(\gamma(ts)) \\
        &\overset{\rm Ex.~\ref{ex X delta 2024}}= &t \sum_{j=1}^r h_j(ts) (\frac1t X_j(\delta_{1/t}\gamma(ts)))
        \overset{\eqref{Boston20171}}{=} \sum_{j=1}^r h_j(ts)  X_j(\eta_t(s)) ,
\end{eqnarray*}
which gives \eqref{Boston20171} from \eqref{Boston20172}.

We claim that $\eta_t$ uniformly converges to $\eta_0$, as $t\to 0$, where $\eta_0(t):= t \dot \gamma(0) $. This claim will complete the proof since in particular, $\eta_t(1) = \delta_{1/t} \gamma(t) \to \dot \gamma(0) $.
To prove the claim, we shall show that for all sequences $t_k \to 0$, there exists a subsequence $t_{k_i}$ such that $\eta_{t_{k_i}} \to \eta_0.$
Indeed, by Ascoli--Arzel\`a Theorem~\ref{AATHM}, there exists a subsequence $t_{k_i}$ and there exists $\xi:[0,1]\to G$ such that $\eta_{t_{k_i}} \to \xi$ uniformly.
We want to show  
$$\dot \xi(s) = \sum_{j=1}^r h_j(0) X_j(\xi(s))
, \qquad \text{ for almost every } s\in [0,1].$$
Let $\sigma$ be the curve such that $\sigma(0)=0$ and $\dot \sigma(s) = \sum_{j=1}^r h_j(0) X_j(\xi(s))$.
We integrate from $0$ to an arbitrary $v\in (0,1)$:
\begin{eqnarray*}
\left| \sigma(v) - \eta_{t_{k_i}}(v) \right|&=&
\left| 
\int_0^v  \dot\sigma(s) \dd s -\int_0^v \dot \eta_{t_{k_i}}(s) \dd s \right| \\
&=&
\left| 
\int_0^v \sum_{j=1}^r h_j(0) X_j(\xi(s)) \dd s - 
\int_0^v \sum_{j=1}^r h_j(t_{k_i}s) X_j(\eta_{t_{k_i}}(s)) \dd s 
 \right| \\&\leq&
\int_0^v \sum_{j=1}^r |h_j(0) -h_j(t_{k_i}s)| \left| X_j(\xi(s)) \right| \dd s +\\
&&\hspace{4.5cm}+
\int_0^v \sum_{j=1}^r |h_j(t_{k_i}s)| \left|X_j(\xi(s)) - X_j(\eta_{t_{k_i}}(s)) \right|\dd s, 
\end{eqnarray*}
where we used \eqref{Boston20171}.
As $i
\to 
\infty$, by continuity of $X_j$ and boundedness of $h_j$, we have that the last summand goes to $0$. Regarding the one before the last, we observe that
\begin{eqnarray*}
\int_0^v \sum_{j=1}^r |h_j(0) -h_j(t_{k_i}s)| \dd s
&\leq&
\int_0^1 \sum_{j=1}^r |h_j(0) -h_j(t_{k_i}s)| \dd s
\\&=&
\dfrac{1}{t}\int_0 ^t |h_j(0) - h_j(u)| \dd u \to 0,
\end{eqnarray*}
since $0$ is a Lebesgue point.
\qed

\proof[\bf Proof of Theorem~\ref{Thm:Pansu:Rademacher}]

Let $F: G\to H$ be a Lipschitz map.
	Define $$F_{p,\epsilon}(x) := \delta_{1/\epsilon} (F(p)^{-1}F(p\delta_\epsilon x)),\qquad \text{ for } p,x\in G \text{ and }\epsilon>0.$$
	Fix $X_1, \ldots X_m$ a basis of $V_1$. 
	For the entire proof, $j\in\{1, \dots, m\}$ and $R_j:=\exp(\R X_j)$. 
	
	Let $\tilde F^j_{p,\epsilon}$ be the restriction $F_{p,\epsilon}|_{R_j}:R_j\to H$. 
	By Proposition~\ref{prop:diff:curves}, for every $p\in G$ the maps 
	$F\circ L_p|_{ R_j}$ are almost everywhere differentable on $ R_j$.
	By Fubini's theorem,
	there is a subset $E\subset G$ of full measure such that, for all $p\in E$, the limit $\tilde F^j_{p,0}:=\lim_{\epsilon\to0^+} \tilde F^j_{p,\epsilon}$ exists and is a Lipschitz group homomorphism $R_j\to H$.
	The limit is uniform on compact subsets of $R_j$.
		
	Let $L$ be a Lipschitz constant of $F$. We shall consider the space ${\rm Lip}^L(R_j;H)$ of $L$-Lipschitz functions from $R_j$ to $H$, with a separable distance that metrizes the uniform convergence on compact sets; see Exercise~\ref{sec01081010}.
		
	We have $\tilde F^j_{p,\epsilon}\in {\rm Lip}^L(R_j;H)$ for every $p\in G$ and $\epsilon\ge0$.
	We can apply Egorov Theorem~\ref{Egorov's} to the functions $p\in G\mapsto 
	\tilde F^j_{p,\epsilon}\in {\rm Lip}^L(R_j;H)$. 
	Therefore, for every $\tau,r>0$ there exists a set $E_{\tau,r}\subset E\cap B(1_G,r)$ such that $|B(1_G,r)\setminus E_{\tau,r}|<\tau$ and
	\begin{equation}\label{domenica_mattina} 
			\begin{array}{l}
			 \{p_\epsilon\}_{\epsilon>0}\subset E_{\tau,r} 
			 \\
			 \lim_{\epsilon\to0}p_\epsilon = p\in E_{\tau,r}
		\end{array}
		\THEN\quad
			\begin{array}{l}
		 \tilde F^j_{p_\epsilon,\epsilon}\to \tilde F^j_{p,0}, \text{ as }\eps\to 0,
		 \\
		 \text{uniformly on compact sets of } R_j.	\end{array}
\end{equation}	
	Let $E_{\tau,r}^\circ\subset E_{\tau,r}$ be the set of density points of $E_{\tau,r}$. Since we are in a doubling metric space, $E_{\tau,r}^\circ$ has full measure within $E_{\tau,r}$.
	
	For the next few paragraphs we fix $p\in E_{\tau,r}^\circ$.
	 We claim that 
		 for all $v\in R_j$, and $\{p_\epsilon, q_\epsilon\}_{\eps>0} \subseteq E_{\tau,r}^\circ$
	\begin{equation}\label{eq602e0fea}
		\begin{array}{l}
			\lim_{\epsilon\to 0} \delta_{1/\epsilon}(p_\epsilon^{-1}q_\epsilon) = v\\
			\lim_{\epsilon\to 0} p_\epsilon = p
		\end{array}
		\THEN\quad
		\lim_{\epsilon\to0}\delta_{1/\epsilon}(F(p_\epsilon)^{-1}F(q_\epsilon)) = \tilde F^j_{p,0}(v) .
	\end{equation} 
%
%
	Indeed, \eqref{eq602e0fea} is a consequence of $p_\epsilon\to p$ in $E_{\tau,r}$:
	\begin{align*}
		d( \delta_{1/\epsilon}(F(p_\epsilon)^{-1}F(q_\epsilon)), \tilde F^j_{p,0}(v) )
		&= d( F_{p_\epsilon,\epsilon}(\delta_{1/\epsilon}(p_\epsilon^{-1}q_\epsilon)), \tilde F^j_{p,0}(v) ) \\
		&\le d( F_{p_\epsilon,\epsilon}(\delta_{1/\epsilon}(p_\epsilon^{-1}q_\epsilon)), F_{p_\epsilon,\epsilon} (v) )
			+ d( F_{p_\epsilon,\epsilon} (v), \tilde F^j_{p,0}(v) ) \\
		&\le L d( \delta_{1/\epsilon}(p_\epsilon^{-1}q_\epsilon), v )
			+ d( \tilde F^j_{p_\epsilon,\epsilon} (v), \tilde F^j_{p,0}(v) ) \to 0,
	\end{align*}
	where at the end we used the first assumption of \eqref{eq602e0fea} and \eqref{domenica_mattina}.


Define
	\[
	\scr D_p := \left\{
	v \in G : 
	\forall \{q_\epsilon\}_{\eps>0} \subseteq E_{\tau,r}^\circ
	\begin{array}{l}
	\text{ if } \lim_{\epsilon\to0}\delta_{1/\epsilon}(p^{-1}q_\epsilon) 
	v \\
	\text{ then }
	\delta_{1/\epsilon}(F(p)^{-1}F(q_\epsilon) ) \text{ converges }
	\end{array}
	\right\} .
	\]
	Therefore, for all $v\in \scr D_p$ there exists an element in $ H$, which we denote by $F_{p,0}(v) $, such that if 
	$q_\epsilon \in E_{\tau,r}^\circ$
	are such that 
$\lim_{\epsilon\to0}\delta_{1/\epsilon}(p^{-1}q_\epsilon) = v$,
	then
\[	F_{p,0}(v):= \lim_{\epsilon\to0}\delta_{1/\epsilon}\left(F(p)^{-1}F(q_\epsilon)\right) .\]


We claim that
	\begin{equation}\label{eq602e136d}
	g\in\scr D_p,\ v\in R_j
	\THEN
	gv\in \scr D_p, 
	\end{equation}
	with
		\begin{equation}\label{eq602e1a54}
F_{p,0}(gv) = F_{p,0}(g) \tilde F^j_{p,0}(v).
	\end{equation}
	Indeed, let $\{q_\epsilon\}_{\epsilon>0} \subset E_{\tau,r}^\circ$ be such that
	$\lim_{\epsilon\to0}\delta_{1/\epsilon}(p^{-1}q_\epsilon) = gv$; notice that since $p\in E_{\tau,r}^\circ $, then such family $q_\epsilon$ exists.
	Since $p\in E_{\tau,r}^\circ$
	then there is $\{p_\epsilon\}_\epsilon \subset E_{\tau,r}^\circ$ such that 
	$\lim_{\epsilon\to0}\delta_{1/\epsilon}(p_\epsilon^{-1}q_\epsilon) = v$.
	So,
	\begin{align*}
	\lim_{\epsilon\to0}\delta_{1/\epsilon}(F(p)^{-1}F(q_\epsilon))
	&= \lim_{\epsilon\to0}\delta_{1/\epsilon}(F(p)^{-1}F(p_\epsilon) )
		\delta_{1/\epsilon}(F(p_\epsilon)^{-1}F(q_\epsilon)) \\
	&\stackrel{\eqref{eq602e0fea}}{=} F_{p,0}(g) \tilde F^j_{p,0}(v) .
	\end{align*}
	

	Next we observe the easy fact $1_G\in \scr D_p$, and therefore from \eqref{eq602e136d} we infer
		\begin{equation}\label{eq602e0feaEnrico}R_1,\ldots, R_m\subset\scr D_p.
	\end{equation}

	From \eqref{eq602e0feaEnrico} and \eqref{eq602e136d}, together with the assumption that $R_1\cup \ldots \cup R_m$ finitely generates $G$, we get that 
	$\scr D_p = G$.
	From \eqref{eq602e1a54} the map $F_{p,0}:G\to H$ is a group homomorphism.
		Notice that if $v\in \scr D_p$,
	then for every
	 sequence $\epsilon_m \searrow 0$ such that 
$	F_{p,\epsilon_m}$ converges uniformly, as $m\to \infty$, we have
\begin{equation}\label{mappe_coincidono}
	F_{p,0}(v)= \lim_{m\to \infty} F_{p,\epsilon_m}(v). \end{equation}
	From \eqref{mappe_coincidono},
	we conclude that every blowup of $F$ at $p$, which exists by Ascoli--Arzel\'a theorem, coincides with the map
 $F_{p,0}:G\to H$. 
	We have proved that $F$ is Pansu differentiable at $p$.
	Since $\bigcup_{\tau,r>0} E_{\tau,r}^\circ$ has full measure in $G$, the map $F$ is differentiable almost everywhere on $G$.
	\qed

 \subsubsection{Original proof of Pansu's theorem}
 
In this subsection, we present the original proof by Pansu together with some extra explanation from Monti's thesis \cite{Monti_thesis}.
 For the proof, we introduce the {\em difference quotients}:
$$R(x; v, t):=  \delta_{{1}/{t}} \left(f(x)^{-1} f(x \delta_t v)\right),$$
so that $
Df(x;v):=\lim_{t\to 0^+} R(x;v,t).$ 

 We start with a preliminary result. It states that if we have partial derivatives in two directions, then we also have the partial derivative in the direction of the product of the two directions. Both the assumption and the conclusion are valid almost everywhere.
 
 \begin{proposition}\label{prop:aewrfqercrqrer}
 Let $f: G\to H$ be a Lipschitz map between sub-Finsler Carnot groups.
 Let $v,w\in G$.
 If $ Df(x; v)$ and $ Df(x; w)$ exist for almost every $x\in G$, then $Df(x; v w)$ exists and $Df(x; v w) = Df(x;v) Df(x; w)$ for almost every $x\in G$.
 \end{proposition}


 \proof
 Let $\Omega\subset \G$ be open and with finite measure.
 Let $\eta > 0$.
 By Egorov's theorem for metric spaces (see Theorem~\ref{Egorov's}), there exists a measurable subset $K=K(\eta)\subset \Omega$ such that the measure of $\Omega \setminus K$ is less than $\eta$ and $R(x;w,t) \to Df(x;w),$ as $t\to 0$, uniformly in $x$ on $K$.
Moreover, since the measure is regular, we may assume that $K$ is compact.

We claim that to conclude the proof, it is enough to show 
\begin{eqnarray}\label{oreuceoirubcwei}		
 R(x\delta_t v;w,t)\stackrel{t\to0}{\longrightarrow} Df(x;w)
, \qquad\text{
for almost every } x\in K.
\end{eqnarray}
 Indeed, in this case, for $x\in K$, we have
\begin{eqnarray*}
R(x;vw,t)&=&  \delta_{1/t} 
\left(f(x)^{-1} f(x \delta_t(vw))\right)\\
&=& \delta_{1/t} 
\left(f(x)^{-1} f(x \delta_tv) \right)  \delta_{1/t} \left( f(x \delta_tv)^{-1} f(x \delta_t v \delta_{t} w)\right) \\
&\stackrel{\eqref{oreuceoirubcwei}}{=}& R(x;v,t) R(x\delta_t v;w,t)\stackrel{t\to0}{\longrightarrow} Df(x;v) Df(x;w).
\end{eqnarray*}
Then one concludes by taking the union of the sets $K=K(\eta)$ when $1/\eta$ varies in $\N$, which forms a full measure set.

 For showing \eqref{oreuceoirubcwei}, take as $x$ a point of density for $K$, (recall that from Theorem~\ref{LDT} these points are of full measure in $K$). 
For $t>0$, let $x_t\in K$ be one projection of $x\delta_tv$ on $K$, i.e., such that $d(x\delta_t v, x_t)=d(x\delta_t v, K)=:r_t$.
Then $r_t\leq d(x\delta_t v, x)= t d( v, 0).$
We claim that $r_t/t\to 0$. Indeed,
$$\dfrac{ r_t^Q}{ (2 t d(v,0))^Q}= 
\dfrac{|B_d(x\delta_tv, r_t)|}{|B_d(x, 2 d(x,x\delta_t v))|} 
\leq
\dfrac{|B_d(x, 2 d(x,x\delta_t v) )\setminus K|}{|B_d(x, 2 d(x,x\delta_t v))|} \to 0,$$
because $x$ is a density point for $K$.

We now calculate the following quantity, defining the three points $A_t, B_t, C_t$ in $H$:
\begin{eqnarray*}
R(x\delta_t v;w,t)&=& \delta_{1/t} \left( f(x \delta_tv)^{-1} f(x \delta_t v \delta_{t} w)\right) \\
&=&\underbrace{
 \delta_{1/t}\left(
f(x \delta_t(v)) ^{-1} f(x_t)\right)}_{A_t}\,
\underbrace{
 \delta_{1/t} \left(f(x_t)^{-1} f(x_t \delta_t w) \right)}_{B_t}\,
\underbrace{
 \delta_{1/t} \left(
f(x_t \delta_t w)^{-1} f(x \delta_tv \delta_{t} w)\right)}_{C_t}
.
\end{eqnarray*}

First, we claim that the point $A_t$   tends to the identity element $ 1_H$, as $t\to 0$. Indeed, 
$$ d(1_H,A_t)=\dfrac{1}{t}  d ( f(x_t), f(x \delta_t v) )\leq \dfrac{L}{t} d(x_t, x \delta_t v) = Lr_t/t \to 0.$$
Second, we note that, since $x_t\in K$, $x_t\to x$, and on $K$ the convergence is uniform, we have that $B_t=R(x_t; w, t) \to Df(x;w)$, as $t\to 0$. 
Third, we claim that $C_t \to 1_H$, as $t\to 0$. Indeed, 
\begin{eqnarray*} d(1_H,C_t)&=&\dfrac{1}{t}  d ( f(x_t \delta_t w), f(x \delta_t v \delta_t w) )\\
&\leq &\dfrac{L}{t} d(x_t \delta_t w, x \delta_t v \delta_t w ) \\&=& L 
d(\delta_{1/t}(x_t) w, \delta_{1/t}(x \delta_t v) w)
\to 0,\end{eqnarray*} where we used that 
$d(\delta_{1/t}(x_t), \delta_{1/t}(x \delta_t v) )=\dfrac{d(x_t, x\delta_t v)}{t}\to 0$, and that right translations are continuous.
\qed

\proof[Another Proof of Theorem~\ref{Thm:Pansu:Rademacher}]

Let $X_1,\ldots, X_r$ be a basis of the first layer of the stratification of $\Lie(G)$.

We first claim that the set $$E:=\large\left\{ p\in \G : Df(p;\exp(X_i)) \text{ and } Df(p;\exp(-X_i)) \text{ exists for all } i\in \{1, \ldots, r\} \large\right\}$$ has full measure.
Indeed,
complete to a basis $X_1,\ldots, X_n$ of $\Lie(G)$.
For $j\in \{ 1, \ldots, r\}$, define $\phi_j :\R^n\to \G$ as
$\phi_j(x_1, \ldots, x_n) =
\exp(\sum_{i\neq j} x_i X_i)\exp(x_j X_j).$
Then $\phi_j $ is a diffeomorphism and for all $x\in \R$ the curve $t\mapsto \phi_j(x+te_j)$ is the flowline of $X_j$ starting at $\phi_j(x)$.
Set $$\tilde E_j := \{x\in \R^n : t\mapsto f(\phi_j(x_1, \ldots, x_{j-1},t,x_{j+1},\ldots, x_n)) \text{ is Pansu differentiable at } t=x_j\}.$$
By Fubini's theorem for the Lebesgue measure and by Proposition~\ref{prop:diff:curves},
the set $\tilde E_j $ has full measure.
Then, the set $E=\bigcap_{j=1}^r \phi_j ( \tilde E_j) $ has full measure.

Let $S:=\{v\in \G : d(0,v)=1\}$ be the unit sphere in $\G$.
For all $m\in \N$ there exist $v_1^m, \ldots, v_{j_m}^m$ such that
$S\subseteq \bigcup_{i=1}^{j_m} B_d(v_i^m, 1/m)$.
We then claim that each set 
$$E_m := \big\{ p\in E : Df(p; v_i^m ) \text{ exists for all } i\in\{1,\ldots,j_m\}\big\}$$
 has full measure.
Indeed, since $\mathcal G:= \big\{\exp(\lambda X_i) : \lambda\in \R, i\in\{1, \dots, r\}\big\}$ generates $\G$, then for all $i$ and all $m$ there exists $w_1,\ldots, w_k\in \mathcal G$ such that $v_i^m=
w_1\ldots w_k$.
Hence, from Proposition~\ref{prop:aewrfqercrqrer} for almost every $p\in \G$ we have that
$Df(p; v_i^m ) $ exists. Thus $E_m $ has full measure.

We finally claim that if $p\in \bigcap_{m\in \N}E_m$, then 
$R(p;v,t)$ converges uniformly in $v\in S$, as $t\to 0$.
Indeed, we want to show that for all $m\in \N$ there exists $\delta>0$ such that for all $s,t\in (0,\delta)$ and all $v\in S$
$$ d(R(p;v,t), R(p;v,s)) \leq \dfrac{1+2L}{m}.$$
Let $m\in \N$. Then there exists $\delta>0$ such that 
for all $i\in \{1,\ldots, i_m\}$ and all $s,t\in (0,\delta)$
$$ d(R(p;v_i^m,t), R(p;v_i^m,s)) \leq \dfrac{1}{m}.$$
Let $v\in S$. Then there exists $i$ such that 
$d(v, v_i^m)\leq \dfrac{1}{m}.$
Then for all $s,t\in (0,\delta)$
\begin{eqnarray*}
 d(R(p;v,t), R(p;v,s)) & \leq &
 d(R(p;v,s), R(p;v_i^m,s), )
+d(R(p;v_i^m,s), R(p;v_i^m,t))
+
 d( R(p;v_i^m,t), R(p;v,t)) \\
&\leq&
\dfrac{1}{s}  d( f(p \delta_sv_i^m), f(p \delta_sv))+\dfrac{1}{m} + \dfrac{1}{t}  d( f(p \delta_tv_i^m), f(p \delta_tv))
 \\
&\leq&
\dfrac{L}{s} s d(v_i^m,v) +\dfrac{1}{m} + \dfrac{L}{t} t d(v_i^m,v) 
\leq \dfrac{L+1+L}{m}.
\end{eqnarray*}
\qed

\subsection{Applications to non-embeddability}

 It was observed by Semmes, 
\cite[Theorem 7.1]{Semmes2}, that
Pansu's differentiation Theorem~\ref{Thm:Pansu:Rademacher} implies that
a Lipschitz embedding of the Heisenberg group with its CC distance into an Euclidean space, 
cannot be bi-Lipschitz. The same holds for every non-commutative Carnot group.

\begin{theorem} There is no bi-Lipschitz embedding from an open set in a non-commutative Carnot group into any Euclidean space. 
\label{non-embed}
\end{theorem}

\proof
Let $G$ be a non-commutative Carnot group.
Suppose that such an embedding $f:U\subseteq G \to \R^n$ exists, for some open set $U$ and $n\in \N$. 
 Pansu Rademacher Theorem~\ref{Thm:Pansu:Rademacher} implies that there exists at least one point in $U$ at which $f$ is Pansu differentiable, and whose tangent map $f_0$ is a group homomorphism. The blowing-up procedure
used to define the tangent map scales in the natural way, i.e.,
if $f$ is $L$-bi-Lipschitz, then each rescaled $f_\lambda$ is $L$-bi-Lipschitz and
 so the tangent map $f_0$ is bilipschitz too. In particular, the map $f_0$ is injective. We now get a contradiction,
because $f_0$ is a group homomorphism
between $G$
 and the abelian $\R^n$.
However, every homomorphism from a Lie group
 into $\R^n$ must have a kernel that contains at least 
the commutator subgroup.
Therefore, the subgroup $[G,G]$ is mapped to 0 via $f_0$,
and hence $f_0$ cannot be injective.
\qed

The same proof gives the following generalization.
\begin{corollary}
Let $G$ and $H$ be Carnot groups.
If no subgroup of $H$ is isomorphic to $G$, then there is no bi-Lipschitz embedding of $G$ in $H$. 
\end{corollary} 


%

A challenging task is to characterize the Banach spaces into which Carnot groups can be biLipschitz embedded.
On the one hand, we have that Kuratowski embedding from \cite{Kuratowski_1935} (see also \cite[page~99]{Heinonenbook}) implies that every separable metric space can be embedded isometrically into $\ell^\infty(\N)$.
But, on the other hand, we also have that every Carnot group that embeds biLipschitz into $\ell^1(\N)$ is commutative; see \cite{zbMATH07702900}.




\section{A metric characterization of Carnot groups}\label{sec characterization of Carnot}\index{Carnot! -- characterization}

 The purpose of this section is to give a more axiomatic presentation of Carnot groups from the viewpoint of Metric Geometry. In fact, we shall see that Carnot groups are the only locally compact and geodesic metric spaces that are isometrically homogeneous and self-similar.
 Such a result follows the spirit of Gromov's approach of `seeing {C}arnot-{C}arath\'eodory spaces from within', \cite{Gromov1}.

 \begin{samepage}
 \begin{theorem}\label{characterization_Carnot}\index{characterization! -- of Carnot groups}
The sub-Finsler Carnot groups are the only metric spaces that are 
\begin{enumerate}
\item locally compact, 
\item geodesic, 
 \item isometrically homogeneous, and
 \item self-similar (i.e., admitting a dilation).
\end{enumerate}
 \end{theorem}
 \end{samepage}
 
 
 Theorem~\ref{characterization_Carnot} provides a new equivalent definition of Carnot groups. 
 Obviously, (1) can be slightly strengthened assuming that the space is {\em boundedly compact} (the term {\em proper} is also used), i.e., closed balls are compact.

We point out that each of the four conditions in Theorem~\ref{characterization_Carnot} is necessary for the validity of the result. Indeed, 
let us mention examples of spaces that satisfy three out of the four conditions but are not Carnot groups:
every infinite-dimensional Banach space;
every snowflake of a Carnot group, e.g., $(\R, \sqrt{\norm{\cdot}})$;
many cones such as the usual Euclidean cone of cone angle in $(0,2\pi)$ or the union of two spaces such as
$\{(x,y)\in \R^2: xy\geq0\}$;
every compact homogeneous space such as $\mathbb S^1$.

Other papers focusing on metric characterizations of Carnot groups are
\cite{LeDonne6}, \cite{Buliga11}, \cite{Freeman12} (which is based on \cite{LeDonne1}), and \cite{Buyalo_Schroeder14}.

 \subsection{Proof of the characterization}\index{Berestovskii! -- characterization}\index{Theorem! Montgomery-Zippin --}\index{Theorem! Mitchell Tangent --}
 The proof of Theorem~\ref{characterization_Carnot} is an easy consequence of three hard theorems: 
 Montgomery-Zippin Theorem~\ref{Montgomery-Zippin}, 
 Berestovskii Theorem~\ref{Berestovskii}, and Mitchell Theorem~\ref{mitchellTheorem}, which will be discussed in the next chapter.

\proof[Proof of Theorem~\ref{characterization_Carnot}]
Let us verify that we can use Theorem~\ref{Montgomery-Zippin}. A geodesic metric space is obviously connected and locally connected.
Regarding finite dimensionality, recall from Proposition~\ref{prop:doubling_homog} that every locally compact, self-similar, isometrically homogeneous space $X$ is doubling and hence finite topological dimension. 
Therefore, by Theorem~\ref{Montgomery-Zippin} the isometry group $G$ is a Lie group.

Since the distance is geodesic, Theorem~\ref{Berestovskii} implies that our metric space is a sub-Finsler homogeneous manifold $G/H$. Since the sub-Finsler structure is $G$ invariant, in particular it is equiregular. Hence, on the one hand, because of 
Theorem~\ref{mitchellTheorem}
the tangents of our metric space are sub-Finsler Carnot groups.
On the other hand, the space admits a dilation, hence, iterating the dilation, we have that there exists a metric tangent of the metric space that is isometric to our original space.
Then the space is a sub-Finsler Carnot group.
\qed


\section{Extremal curves in Carnot groups}\label{sec Extremal curves in Carnot}

\subsection{Expected regularity of sub-Riemannian geodesics}

\subsubsection{Conjectures}
Given two points in a sub-Finsler group, 
the existence of a geodesic between them is ensured by Ascoli--Arzel\`a Theorem; see Proposition~\ref{exists: short}.
Being Lipschitz, we know that every such a curve is differentiable almost everywhere.

It is expected that when the metric is sub-Riemannian, then every geodesic is $C^1$, or, in fact, $C^\infty$. 
When instead the norm comes from a polytope, i.e., its unit ball is the convex hull of finitely many points, then we expect that there exists a constant $N\in\N$ such that each pair of points can be connected with a geodesic made of $N$ smooth pieces. 

\begin{conjecture}[Regularity conjecture for sub-Reimannian manifolds]
 If $(M,\Delta, \langle \cdot,\cdot\rangle)$ is a
sub-Riemannian manifold, then each geodesic for the $CC$-distance is $C^\infty$.
\end{conjecture}

\begin{conjecture}[Weak regularity conjecture for sub-Reimannian Carnot groups]
 If $\G$ is a Carnot group, then each pair of points can be connected by a $C^1$ geodesic.
\end{conjecture}

\begin{conjecture}[Regularity Conjecture for sub-Finsler Carnot groups]
 If $(\G,V_1, \norm{\cdot}_1 )$ is a Carnot group where $\norm{\cdot}_1$ is the $\ell^1$ norm, then there exists a constant $K$ such that each pair of points can be connected by a geodesic that is the concatenation of at most $K$ horizontal lines.
\end{conjecture}

\subsubsection{Various partial results}


The following is a collection of some partial results from \cite{Leonardi-Monti}, \cite{Tan_Yang_step_3, LMPOV}, \cite{LeDonne_Paddeu_Socionovo}, and \cite{Sussmann:regularity}, respectively:
\begin{theorem}\label{thm_parziali_geodetiche} 
Let $G$ be a sub-Riemannian Carnot group and $\gamma:[0,1]\to G$ be an energy minimizing curve.
\begin{description}
\item[\ref{thm_parziali_geodetiche}.i.] If $G$ has rank 2 and step $\leq 4$, then $\gamma$ is $C^\infty$.
\item[\ref{thm_parziali_geodetiche}.ii.] If $G$ has step $\leq 3$, then $\gamma$ is $C^\infty$.
\item[\ref{thm_parziali_geodetiche}.iii.] If $[G,G]$ is abelian, then $\gamma$ is $C^1$.
\item[\ref{thm_parziali_geodetiche}.iv.] The curve $\gamma$ is analytic on some open dense subset of its domain $[0,1]$.
\end{description}
 \end{theorem}
For sub-Riemannian manifolds, there are several other statements for which it is not clear if there is an analog for Carnot groups. 
For example, in \cite{Chitour-Jean-Trelat} 
the authors proved that for generic sub-Riemannian structures $(M,\Delta, \langle \cdot,\cdot\rangle)$ of rank at least $ 3$, 
energy minimizers are $C^\infty$.
Here, generic means on an open dense set of structures, with respect to the Whitney $C^\infty$ topology.
 
Regarding the case of polygonal norms very little is know. The following result is a summary from \cite{Breuillard_MR3267520, BBLDS, Ardentov_LeDonne_Sachkov1,Ardentov_LeDonne_Sachkov2} 
\begin{theorem}\index{Heisenberg! -- group}\index{Engel! -- group}\index{Cartan! -- group}  
Let $G$ be either the Heisenberg group, the Engel group, or the Cartan group, with first strata $V_1\simeq \R^2$. 
Let $\norm{\cdot}_1$ be the $\ell^1$ norm on $V_1 $. 
Then the geodesics with respect to the $CC$-distance of $(G,V_1, \norm{ \cdot}_1)$
are made of, at most, a concatenation of $14$ pieces of horizontal lines.
\end{theorem}

\subsection{Sublinear isometric property of projections}\index{sublinear isometric property}

In addition to the first-order analysis given by the PMP Theorem~\ref{thm_Pontryagin} (together with the second-order analysis of Goh's Theorem~\ref{Goh_thm} and its generalizations \cite{Boarotto_Monti_Palmurella, zbMATH07826620}), there is another method to deduce some type of regularity for geodesics.
The idea is to approximate geodesics in Carnot groups with lifts of geodesics in quotient groups. In brief, we have that for each geodesic, one of its quotients to a group of lower step is a geodesic up to a sublinear error. The following statement is the precise formulation; afterward, we will draw a list of consequences.

\begin{theorem}\label{thm:quantified tangent}
 
 Let $G$ be a Carnot group and $\gamma:[0,1]\to G$ be a geodesic. 
 Let $H$ be the smallest Carnot subgroup containing $H$ and $s$ be the step of $H$. 
 Then, there exists $C>0$ such that
	 \begin{equation}
 \label{eq:almost-geodesic} \abs{a-b}-C\abs{a-b}^{\frac{s}{s-1}}\leq d(\pi_s(\gamma(a)),\pi_s(\gamma(b)))\leq \abs{a-b},\qquad \forall a,b\in [0,1], \end{equation}
 where $\pi_s:H\to H/\exp(V_s(H))$ is the canonical submetry and $C^{s}(H)= \exp(V_s(H))$.
\end{theorem}

\begin{proof}
Since $\pi_s $ is a submetry, we have the
 second inequality of \eqref{eq:almost-geodesic}. 
To prove the first one, by compactness of $[0,1]^2$, 
it is enough to show that for all $\bar t\in[0,1]$ there exists $C,\delta>0$ such that \eqref{eq:almost-geodesic} holds for all $a,b\in [\bar t-\delta,\bar t+\delta]$. 
Fix $\bar t\in[0,1]$. Write $V_j:=V_j(H)$. Being $H$ the minimal Carnot subgroup containing $\gamma$, there exist increasing $t_1,\ldots,t_n\in [0,1]\setminus \{ \bar t\}$ and $Y_1,\ldots,Y_n\in V_1\oplus\cdots\oplus V_{s-1}$ such that $\{\Ad_{\gamma(t_i)}Y_i \}_{i\in\{1,\ldots,n\}}$ is a basis of $\mathfrak{h}:=\Lie(H)$.

Define $\phi:\mathbb{R}^n\to \mathfrak{h}$ by setting
\begin{equation}
 \label{eq:def_phi}
 \phi(y):=\log\left(\prod_{j=1}^n \mathrm{C}_{\gamma(t_j)}(\exp(y_j Y_j))\right), \quad \forall y\in\mathbb{R}^n.
\end{equation}
The point $0\in \R^n$ is a regular point for $\phi$, since 
 \begin{equation*}
 \mathrm{span}\{\Ad_{\gamma(t_i)}Y_i \}_{i\in\{1,\ldots,n\}}=\mathfrak{h}.
\end{equation*}
Hence, after fixing
 norms $\|\cdot \|$ on $\mathbb{R}^n$ and $\mathfrak{h}$, we apply 
the Inverse Mapping Theorem:
  there exist $\epsilon,C'>0$ such that 
\begin{equation}
\label{eq:open_map_thm}
 Z\in \mathfrak{h}, \quad \|Z\|\leq \epsilon\implies \begin{gathered}
 \exists y\in\mathbb{R}^n \text{ such that}\\
 \phi(y)=Z \text{ and }
 \|y\|\leq C'\|Z\|.
 \end{gathered} 
\end{equation}
 
Set $\delta:= \min\left\{|t_1-\bar t|, \dots, |t_n-\bar t|,  \frac{\epsilon^{ {1}/{s}}}{2 C_{\text{BB}}}  \right\}$, where $C_{\text{BB}}$ is the constant coming from the Ball-Box Theorem.
Take $a,b\in (\bar t-\delta,\bar t + \delta)$.
Consider a geodesic $\bar \sigma:[a,b]\to  H/\exp(V_s)$ between $\pi_s(\gamma(a))$ and $\pi_s(\gamma(b))$ and let $\sigma:[a,b]\to H$ be the lift of $\bar \sigma$ from $\gamma(a)$. Set 
\begin{equation}\label{eq Z centrale}
Z:=\log(\gamma(b)\sigma(b)^{-1})\in V_s . 
\end{equation} 
Notice, in particular, that $Z$ is in the center of $\h$.
We have,
\begin{equation}
\label{eq:stima-norma-Z}
 \|Z\|^{\frac{1}{s}}\leq C_{\text{BB}} d(1,\exp(Z))=C_{\text{BB}} d(\sigma(b),\gamma(b))\leq 2C_{\text{BB}} |a-b|.
\end{equation}

Since $|a-b|\leq \delta \leq \frac{\epsilon^{ {1}/{s}}}{2 C_{\text{BB}}} $, by \eqref{eq:stima-norma-Z} we have $\|Z\|\leq \epsilon$. Thus by \eqref{eq:open_map_thm} there exists $y\in\mathbb{R}^n$ with \begin{eqnarray}
\label{eq:phi=Z}
 \phi(y)=Z,\\
 \label{eq:norma_y}
 \|y\|\leq C'\|Z\|.
\end{eqnarray}

 \begin{figure}
 \centering
 \includegraphics[width=4in]{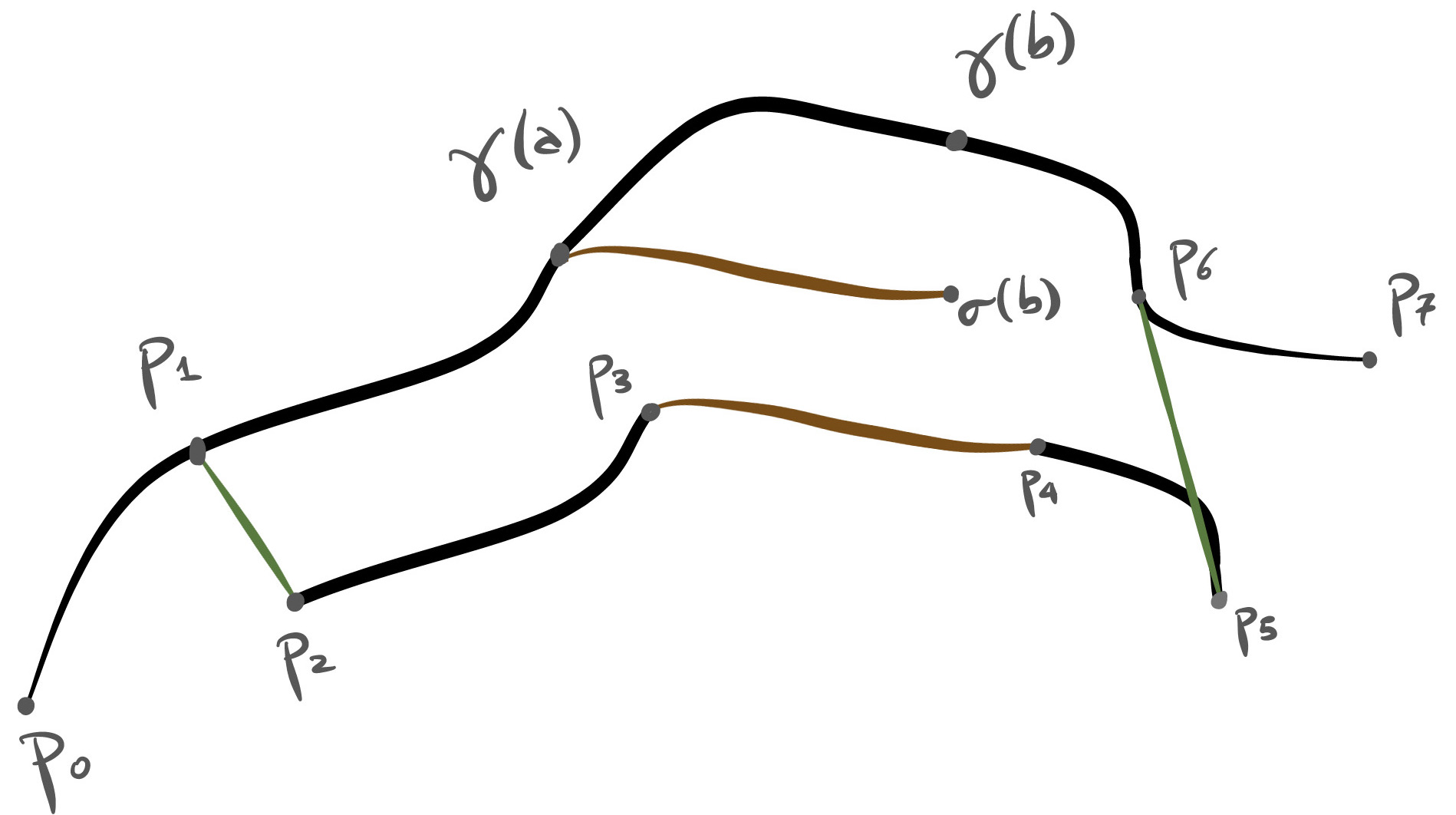}
	\caption{A sketch of the points used in the triangle inequality in the proof of Theorem~\ref{thm:quantified tangent}.
	For simplicity, we assume $\{t_1, \ldots, t_n\}=\{t_1, t_2^*\}$, and $t_1<a<b<t_2^*$. We consider the following points:
\\
 $		p_0 = \gamma(0) $; \\
$		p_1 = \gamma(t_1) $; \\
$			p_2 = p_1 \exp(y_1 Y_1) $; \\
$				p_3 = p_2 \gamma(t_1)^{-1} \gamma(a)$;\\
$					p_4 = p_3 \sigma(a)^{-1} \sigma(b)$; \\
$						p_5 = p_4 \gamma(b)^{-1} \gamma(t_2^*)$; \\
$							p_6 = p_5 \exp(y_2 Y_2)
=\gamma(t_1) \exp(y_1 Y_1) \gamma(t_1)^{-1} \gamma(a)\sigma(a)^{-1} \sigma(b)\gamma(b)^{-1} \gamma(t_2^*)\exp(y_2 Y_2)
$\\ $\mbox{\hspace{0.4cm}}=C_{\gamma(t_1) }(\exp(y_1 Y_1) ) \sigma(b)\gamma(b)^{-1} 
 C_{\gamma(t_2^*)}(\exp(y_2 Y_2))\gamma(t_2^*)
$\\ $\mbox{\hspace{0.4cm}}= \gamma(t_2^*) $; \\
$								p_7 = \gamma(1) $.
	}
	\label{fig:7_points}
\end{figure}

Let $k:= \max\{j\in \{1,\ldots,n\}\mid t_j<a\}$.
By \eqref{eq:phi=Z} we get
\begin{eqnarray*}
 \gamma(1) &\stackrel{\eqref{eq:phi=Z}}{=}& \exp(\phi(y))\exp(Z)^{-1}\gamma(1)\\
 &\stackrel{{\eqref{eq:def_phi}, \eqref{eq Z centrale}}}{=}&\left( \prod_{j=1}^n \mathrm{C}_{\gamma(t_j)}(\exp(y_jY_j)) \right)\sigma(b)\gamma(b)^{-1}\gamma(1)\\
 & 
 {=}&\left( \prod_{j=1}^k \mathrm{C}_{\gamma(t_j)}(\exp(y_jY_j)) \right) \gamma(a)\gamma(a)^{-1}\sigma(b)\gamma(b)^{-1} \left(\prod_{j=k}^n \mathrm{C}_{\gamma(t_j)}(\exp(y_jY_j)) \right)\gamma(1),
\end{eqnarray*}
 where in the last equality we used that $\sigma(b)\gamma(b)^{-1}$ is in the center of $H$. Evaluating the distance from $\gamma(0)$, applying a triangular inequality and using the fact that $\gamma$ is a geodesic, we get

 \begin{eqnarray*}
 1 &=& d(\gamma(0),\gamma(1))\\
 &=& d(1,\gamma(0)^{-1}\prod_{j=1}^k \mathrm{C}_{\gamma(t_j)}(\exp(y_jY_j)) \gamma(a)\gamma(a)^{-1}\sigma(b)\gamma(b)^{-1} \prod_{j=k}^n \mathrm{C}_{\gamma(t_j)}(\exp(y_jY_j))\gamma(1))\\
 &\leq&\sum_{j=1}^nd(1,\exp(y_jY_j))+d(\gamma(0),\gamma(t_1))+ \sum_{j=1}^{k-1} d(\gamma(t_j),\gamma(t_{j+1}))+d(\gamma(t_k),\gamma(a))\\
 &+& d(\gamma(a),\sigma(b))+ d(\gamma(b),\gamma(t_{k+1}))+\sum_{j=k+1}^{n-1} d(\gamma(t_j),\gamma(t_{j+1}))+d(\gamma(t_n),\gamma(1))\\
 &\stackrel{{\rm geodesic}}{=} & \sum_{j=1}^nd(1,\exp(y_jY_j)+1-|b-a|+d(\sigma(a),\sigma(b))\\
 &=& \sum_{j=1}^nd(1,\exp(y_jY_j)+1-|b-a|+d(\pi_s(\gamma(a)),\pi_s(\gamma(b)))\\
 &\stackrel{\text{BB}}{\leq} &n C_{\text{BB}} \norm{y}^{\frac{1}{s-1}}+1-|b-a|+d(\pi_s(\gamma(a)),\pi_s(\gamma(b)))\\
 &\stackrel{\eqref{eq:norma_y}}{\leq} & C'' \norm{Z}^{\frac{1}{s-1}}+1-|b-a|+d(\pi_s(\gamma(a)),\pi_s(\gamma(b)))\\
 &\stackrel{\eqref{eq:stima-norma-Z}}{\leq} & C|b-a|^{\frac{s}{s-1}}+1-|b-a|+d(\pi_s(\gamma(a)),\pi_s(\gamma(b)))
 .
 \end{eqnarray*}
\end{proof}

Here are some important consequences of Theorem~\ref{thm:quantified tangent}:
\begin{itemize}
\item 
The projection modulo the last layer $\exp(V_s)$ of every blowup of a geodesic in a step-s Carnot group is a geodesic. 
\item
In every sub-Riemannian manifold of step $s$, blowing up a geodesic $s$ times gives lines. 
\item
At every point, among the blow-ups of a geodesic, there is a line.

\item Corners cannot be geodesics.

\item
Every piecewise $C^1$ geodesic is $C^1$.

\end{itemize}
For more on the subject, we refer to \cite{Hakavuori_LeDonne_2016, Hakavuori_LeDonne_2018, Monti_Pigati_Vittone, zbMATH06944124}.




\section{Abnormal curves in Carnot groups}

Non-constant abnormal curves are present in step-2 Carnot groups. For example, there is one in $\R$ times the Heisenberg group; recall Exercise~\ref{ex_prod_abn}. There are none exactly when the polarization is {\em strongly bracket generating}, in the sense that for every $X\in V_1\setminus \{0\}$ we have $[X, V_1]=V_2$; see \cite{LeDonne_NicolussiGolo}.\index{bracket-generating! strongly --}
However, 
the abnormal curves starting from the identity element are confined in a subvariety, which has Haar measure 0.
This general result is proved in \cite{LMPOV}. However, next, in this section we discuss the proof in the case of free groups.

\subsubsection{Expected Sard property for the endpoint map}\index{Sard Property}
Recall that the abnormal curves (or, better, their final points) are precisely the critical values of the endpoint map.
Hence, in analogy with the classical Sard Theorem \cite[page~16]{zbMATH01950480}, it is expected that they form a zero-measure set. 

The query is still open, and there are several conjectures for possible results: 
 

\begin{conjecture}[Strong Sard Property]
Let $G$ be a nilpotent simply connected Lie group and $V\subseteq\Lie(G)$ a polarization. Then, the set of points in $G$ that are connected from $1_G$ by some abnormal curve is contained in a proper real algebraic subvariety of $G$.
\end{conjecture}

It is not clear whether the validity of the above conjecture is actually feasible. The next version is more likely to be true.

\begin{conjecture}[Weak Sard Property]
Let $(G,V)$ be a sub-Riemannian Carnot group. 
Then the set of points in $G$ that are connected from $1_G$ by some abnormal energy minimizing curve is contained in a zero-measure set
 of $G$.
\end{conjecture}

%
 \subsubsection{Extremals in step-two free-nilpotent Lie groups}
 
 In step-two Carnot groups, the Strong Sard Conjecture holds. In fact, the subvariety containing the abnormal curves is explicit; see \cite{zbMATH07643478}. We present it in the case of free groups, together with the fact that geodesics are normal curves, without using Goh result.
\begin{theorem}\label{thm_free_2step_extremals}
Let $\mathbb{F}_{n2}$
be the free Carnot group of rank $n$ and step 2, equipped with the standard sub-Riemannian Carnot metric. Then, we have the following two properties:
\begin{description}
\item[\eqref{thm_free_2step_extremals}.i.]
Every geodesic in $\mathbb{F}_{n2}$ is normal.
\item[\eqref{thm_free_2step_extremals}.ii.]
The Strong Sard Property on $\mathbb{F}_{n2}$ holds: Seeing $\mathbb{F}_{n2}$ as $\R^n \times \Lambda^2(\R^n)$, 
	abnormal curves starting from $(0,0)$ 
	are contained in the set
	\begin{align*}
	&\{(\theta,\omega)\in \R^n\times\bigwedge^2(\R^n) : \omega^{\frac{n}{2}} = 0 \}
	& \text{ if $n$ is even, or}\\
	&\{(\theta,\omega)\in \R^n\times\bigwedge^2(\R^n) : \theta\wedge\omega^{\frac{n-1}{2}} = 0 \}
	& \text{ if $n$ is odd.}
	\end{align*}
\end{description}
\end{theorem}

\begin{proof}
Let $\gamma:[0,T]\to \mathbb{F}_{n2}$ be a geodesic. We also assume $\gamma(0)=1_{\mathbb{F}_{n2}}$.
By PMP Theorem~\ref{thm_Pontryagin}, then it is a normal or an abnormal curve. Assume it is abnormal. 
Then, by Proposition~\ref{existence: subgroup} there exists a subgroup $H$ containing $\gamma $ and in which $\gamma$ is not abnormal.
We see $\mathbb{F}_{n2}$ as $\R^n \times \Lambda^2(\R^n)$.
There is a subspace $  V\subset \R^n$ such that  $ H=V\times \Lambda^2(V)$.
Notice that $H$ admits a complementary normal $N$ subgroup in $\mathbb{F}_{n2}$. 

Back to the curve $\gamma$, it is also a geodesic within $H$. Since it is not abnormal, it is normal in $H$, say with covector $\lambda\in \Lie(H)^*$. 
Hence, the covector $\lambda$ can be extended to be 0 in the Lie algebra of $N$.
With this extension, the curve $\gamma$ is normal in $\mathbb{F}_{n2}$.
Thus, we proved \eqref{thm_free_2step_extremals}.i.

Regarding \eqref{thm_free_2step_extremals}.ii, let $\gamma=(\theta,\omega)$ be an abnormal curve in $\R^n\times \Lambda^2(\R^n)$.
	By the discussion above, the curve $\theta$ is contained in a proper subspace $V$ of $\R^n$ and $\omega$ is contained in $\bigwedge^2V$.
	If $(\theta,\omega)\in V\times\bigwedge^2V$,
	then the degree of $\omega^{\frac{n}{2}}$ or $\theta\wedge\omega^{\frac{n-1}{2}}$ is $n>\dim(V) $.
	This implies that $\omega^{\frac{n}{2}}=0$ or $\theta\wedge\omega^{\frac{n-1}{2}}=0$.
\end{proof}


\section{\it Supplementary material}

\subsection{Self-similar sub-Finsler spaces}\label{sec self-similar sub-Finsler spaces}\index{self-similar}

In sub-Finsler geometry, self-similar spaces are well characterized. As differentiable manifolds they have a Lie homogeneous structure of a quotient space of a Carnot group modulo the action of a dilation-invariant subgroup via left-multiplication.
The Lie group metric quotient is the one we saw in Proposition~\ref{prop:distance_quotient} and then in Proposition~\ref{prop_subFin_submetry_Lie}. 
Hence, the well-defined action on the quotient is on the right, and it may not be by isometries. In fact, these quotients may not be isometrically homogeneous spaces. They are still called homogeneous because they admit dilations. To avoid this double use of the word homogeneity, since Section~\ref{sec: self-similar}
we adopted the term self-similar spaces.

\begin{definition}\label{def_selfsimilar}
A \emph{self-similar sub-Finsler space} (or, better, a \emph{self-similar constant-rank sub-Finsler space}) is a sub-Finsler manifold obtained as the quotient space of a (left-invariant) sub-Finsler Carnot group with respect to the left-action of a dilation-invariant subgroup, and it is equipped with the quotient distribution and metric, as in Proposition~\ref{prop_subFin_submetry_Lie}. Namely, 
assume $G$ is a Carnot group with distribution $\Delta$ and left-invariant norm $||\cdot||$ and $H<G$ is a closed dilation-invariant subgroup, for which $T_1H \cap \Delta_1=\{0\}$. 
On the quotient manifold $\mfaktor{H}{G}:=\{H g : g\in G\}$ we consider the sub-Finsler structure that makes the projection $\pi: G\to \mfaktor{H}{G}$ a submetry: we take $\Delta_{\mfaktor{H}{G}}:=\pi_*\Delta$ as (constant-rank) distribution and we define the continuously varying norm on $\mfaktor{H}{G}$ setting for all $p\in \mfaktor{H}{G}$ and for all $v\in 
(\pi_*\Delta)_p\subseteq 
T_p(\mfaktor{H}{G})$
\begin{equation}
 ||v||_{\mfaktor{H}{G}} :=\inf\{||w||_q : q\in \pi^{-1}(p), w\in T_q\Delta, \dd\pi_q(w)=v \}.
 \label{d defnormadown}
\end{equation}
\end{definition}

We stress that, if $G$ is a Carnot group with dilations $(\delta_\lambda)_{\lambda\in \R_+}$ then $H$ is dilation-invariant if $\delta_\lambda(H)=H$, for all $\lambda\in \R_+$, and in this case
the map
\begin{equation}\label{sentendo_Renato_Zero} \delta_\lambda ( Hg):= H\delta_\lambda (g), \qquad \forall g\in G, \forall \lambda\in \R_+
\end{equation}
is well defined. Metrically, each map $\delta_\lambda$ is a dilation of factor $\lambda$, in the sense of Section~\ref{sec: self-similar}.
We, therefore, clarified that metrically, every self-similar sub-Finsler space is self-similar.

Vice versa, we claim that every sub-Finsler space that admits a non-trivial dilation is of the above form. Indeed, for every dilation 
of factor different than 1, we can assume the factor is in the open interval $(0,1)$, up to replacing the map with its inverse.
Hence, being a contraction, the map has a fixed point. Consequently, the space (pointed at this point) is isometrically to a dilation of it. Therefore, the space is isometric to one of its metric tangent spaces. 
 By the work of Bella{\"\i}che (\cite{bellaiche}; see also \cite{Antonelli_LeDonne_MR4645068} for the sub-Finsler case), we know that the metric tangents of (constant-rank) sub-Finsler manifolds are (constant-rank) self-similar spaces.

\subsection{Local isometries of Carnot groups}\label{sec Local isometries of Carnot groups}\index{isometry}

From what we discussed in Theorem~\ref{teo05171840} and more generally in Theorem~\ref{thm_Ville}, isometries of sub-Finsler Carnot groups are affine maps:

\begin{corollary}
Let $G$ and $ H$ be sub-Finsler Carnot groups with normed first layers $(V_1(G), \norm{\cdot})$ and $(V_1(H), \norm{\cdot})$, respectively. 
A map $F: G\to H$ is an isometry with respect to the sub-Finsler metrics if and only if there is $g\in G$ and a Lie group automorphism $\phi: G\to H$ such that $\dd\phi|_{V_1}$ is an isometry between $(V_1(G), \norm{\cdot})$ and $(V_1(H), \norm{\cdot})$ and $F= \phi \circ L_g$.
\end{corollary}

 In sub-Finsler Carnot groups we additionally have that local isometries are restrictions of global isometries.
\begin{theorem}[\cite{LeDonne-Ottazzi}]\label{localisometries}
Let $G_1$ and $ G_2$ be sub-Finsler Carnot groups, and for $i=1,2$, consider $\Omega_i\subset G_i$ open sets. If $F:\Omega_1\to \Omega_2$ is an isometry, then such that $F$ is the restriction to $\Omega_1$ of an isometry $\tilde G_1\to G_2$.
\end{theorem}

Note that in the statement above, we require the domain $\Omega_1$ to be open.
 The hypothesis that the set $\Omega_1$ is open cannot be dropped. 
Unlike in the Euclidean space, Theorem \ref{localisometries} cannot be generalized to arbitrary subsets. We present a counterexample: take the sub-Riemannian Heisenberg group $(\mathbb H, d_{SR})$ in standard exponential coordinates $(x,y,z)$ with respect to the basis of its Lie algebra given by vectors $X, Y$ and $[X, Y]=: Z$. Consider the set given by the $xy$ plane together with the $z$ axis, 
namely,
$$E:=\exp (\R X\oplus\R Y) \cup \exp (\R Z).$$
Then the map 
 $(x,y,z)\mapsto (x,y,-z)$  
%
 is an isometry of $E$ into itself. However, this map is not the restriction of an isometry, and, actually, nor of a bi-Lipschitz map. 
 Indeed, notice that every homeomorphism that extends $F$ is reversing the topological orientation, while the isometries (and the bi-Lipschitz homeomorphisms) of the Heisenberg group preserve the topological orientation; see Exercise~\ref{ex_orientation_affine_Heis}.

\section{Exercises}


\begin{exercise} If $G$ is a Carnot group and $\Delta$ is the left-invariant distribution with $\Delta_{1_G}=V_1$, then $(\Delta^{[j]})_{1_G}=V_1\oplus \cdots\oplus V_j $.
\end{exercise}

\begin{exercise}
Let $\delta_\lambda$ be the dilation of factor $\lambda$ as defined either at the group level in Definition~\ref{Dilations:groups} or at the Lie-algebra level in Definition~\ref{Dilations:algebras}.
Then, we have that $(\delta_\lambda)^{-1} =\delta_{1/\lambda}$.
\end{exercise}

\begin{exercise}[Carnot morphism]
Let $\varphi: G \to H$ be a Lie group homomorphism between Carnot groups. 
Then, we have $\varphi_*(V_1(G))\subseteq V_1(H)$ if and only if $\varphi\circ \delta_\lambda = \delta_\lambda \circ\varphi$, for all $\lambda\in \R$. In this case, we say that $\varphi$ is a {\em homogeneous group homomorphism} or that it is a {\em Carnot morphism}.\index{homogeneous! -- group homomorphism}\index{Carnot! -- morphism}
\end{exercise}

\begin{exercise}\label{ex X delta 2024}
Let $G$ be a Carnot group. For all $ X\in \Lie(G), 
 u\in C^\infty(G),
 \lambda\geq 0, $ and $ 
 g\in G $ we have
$X(u\circ\delta_\lambda)(g)=(\delta_\lambda X)u(\delta_\lambda g)$.
\end{exercise}

\begin{exercise}
	Let $G$ be a Carnot group. For all $p\in G$, for all $r>0$
	\[
	B_{\dcc}(p,r) = L_p(\delta_r(B_{\dcc}(1_G,1))) .
	\]
\end{exercise}


\begin{exercise}
Consider a Carnot basis $X_1,\ldots, X_n$ of a Carnot algebra. Each element $X_j$ of the basis is such that
$$X_j=
[\ldots[[X_{j,1},X_{j,2}],X_{j,3}], \ldots,X_{j,d_j}],$$
where $X_{j,1},\ldots,X_{j,d_j}$ are basis elements in $V_1$, and $d_j$ is such that $X_j\in V_{d_j}$, in other words, it is the degree of $X_j$.
\\ {\it Hint.} Iterate \eqref{Carnotbasis}.
\end{exercise}

\begin{exercise}
In every Carnot group, there is a (strong) Malcev basis.
\end{exercise}

\begin{exercise}\label{basisV2}
Let $V$ and $W $ be vector subspaces of a Lie algebra $\g$ with 
$X_{1},\ldots,X_{l}$ and $Y_{1},\ldots,Y_{m}$ basis of $V$ and $W$, respectively. Then, the vectors
$[X_i,Y_j]$, for $i\in\{1,\ldots,l\}, j\in\{1,\ldots,m\}$ span the subspace $[V,W]$, thus one can extract a basis among such brackets.
\end{exercise}

\begin{exercise} \label{CarnotMalcev} Let $\g= V_1\oplus \cdots\oplus V_s $
be a stratification of a Lie algebra.
Assume that $X_{m_j+1},\ldots,X_{m_j}$ is a basis of $V_j$, for all $j\in \{1,\ldots, s\}$, then the order-reversed basis $(X_n,\ldots,X_1)$ is a Malcev basis.
\end{exercise} 


\begin{exercise}For the Engel group, with group product \eqref{prod_Engel_2024}, verify that the differential at $\bf 0$ of the left translation by $\bf x$ is
\\$
 \mathrm{d} (L_{\bf x})_{\bf 0}= \left[\begin{matrix} 
 1 & 0 & 0& 0 \\
 0 & 1 & 0& 0\\
 -\frac{x_2}{2} & \frac{x_1}{2} &1& 0\\
 -\frac{x_1x_2}{12} -\frac{x_3}{2} & \frac{x_1^2}{12} &\frac{x_1}{2}& 1 
 \end{matrix}\right]\,.
$

\end{exercise} 

\begin{exercise}For the Engel group,
consider the exponential coordinates of the second kind with respect to the ordered basis $(X_1, \ldots, X_4)$ with relations as in \eqref{Engel_rel}. In these coordinates, the left-invariant vector fields are
\\$\begin{cases}
\tilde X_1=\partial_{x_1}&\\
\tilde X_2=\partial_{x_2}+x_1\, \partial_{ x_3}+\frac{x_1^2}{2}\partial_{x_4}&\\
 \tilde X_3= \partial_{ x_3}+x_1\partial_{x_4}&\\
 \tilde X_4=\partial_{x_4}.
\end{cases}.$

\end{exercise}

\begin{exercise}\index{Cartan! -- group}  For the Cartan group, with group product \eqref{prod_Cartan_2024}, 
the induced left-invariant vector fields are:
\\$
\begin{cases}
X_1={\partial}_{x_1}-\frac{x_2}{2}{\partial}_{ x_3}-\big(\frac{x_3}{2}+\frac{x_1x_2}{12}\big){\partial}_{x_4}-\frac{x_2^2}{12}{\partial}_{x_5}&\\
 X_2={\partial}_{ x_2}+\frac{x_1}{2}{\partial}_{x_3}+\frac{x_1^2}{12}{\partial}_{x_4}-\big(\frac{x_3}{2}-\frac{x_1x_2}{12}\big){\partial}_{x_5}&\\
 X_3={\partial}_{x_3}+\frac{x_1}{2}{\partial}_{x_4}+\frac{x_2}{2}{\partial}_{ x_5}&\\
 X_4={\partial}_{x_4}&\\
 X_5={\partial}_{x_5}&
\end{cases}.$\\
While, the induced right-invariant vector fields are given by the columns of the following matrix:
\\$
 \mathrm{d} (R_{\bf x})_{\bf 0}= \left[\begin{matrix} 
 1 & 0 & 0& 0 & 0\\
 0 & 1 & 0& 0& 0\\
 \frac{x_2}{2} & -\frac{x_1}{2} &1& 0& 0\\
 \frac{x_3}{2}-\frac{x_1x_2}{12} & \frac{x_1^2}{12} &-\frac{x_1}{2}& 1 & 0\\
 -\frac{x_2^2}{12}& \frac{x_3}{2}+\frac{x_1x_2}{12}&-\frac{x_2}{2} & 0& 1
 \end{matrix}\right]\,.
$

\end{exercise}

\begin{exercise}\index{Cartan! -- group}  For the Cartan group, with group product \eqref{prod_Cartan_2024}, verify that the left-invariant 1-forms dual to the standard basis are:
\\$\begin{cases}
 \theta_1=\mathrm{d}x_1 &\\
 \theta_2=\mathrm{d}x_2&\\
 \theta_3=\mathrm{d}x_3-\frac{x_1}{2} \mathrm{d}x_2+\frac{x_2}{2}\mathrm{d}x_1&\\
 \theta_4=\mathrm{d}x_4-\frac{x_1}{2}\mathrm{d}x_3+\frac{x_1^2}{6}\mathrm{d}x_2+\big(\frac{x_3}{2}-\frac{x_1x_2}{6}\big)\mathrm{d}x_1&\\
 \theta_5=\mathrm{d}x_5-\frac{x_2}{2}\mathrm{d}x_3+\big(\frac{x_3}{2}+\frac{x_1x_2}{6}\big)\mathrm{d}x_2-\frac{x_2^2}{6}\mathrm{d}x_1&
\end{cases}.$

\end{exercise}


\begin{exercise}
For $s\in \N$, 	the following vector fields in $\R^{s+1}$ generate a Lie algebra isomorphic to the $s$-step filiform
 Lie algebra of the first kind 
\\$\begin{cases}
 X_1 =\partial_{x_1}&\\
 X_2 =\partial_{x_2}+ \sum_{j=1}^{s-1} x_1^j\, \partial_{ x_{j+2}}.&
\end{cases}.$

\end{exercise}

\begin{exercise}
 Extend the proof of Proposition~\ref{prop66b1e96e} to the case $G$ is a positively graded group.
\end{exercise}

\begin{exercise} 
	For the $\Bx$ as in \eqref{def_box_2024} and $ \delta_\lambda$ as in \eqref{def_delta_2024}, we have 
	$\delta_\lambda(\Bx(r)) = \Bx(\lambda r)$, for all $r,\lambda>0$ 
\end{exercise}


\begin{exercise}
	On each Carnot groups, the notion of homogeneous dimension from Definition~\ref{def Homogeneous dimension Carnot} agrees with the one on sub-Finsler manifolds in Definition~\ref{def Homogeneous dimension manifolds}.
\end{exercise}
%

\begin{exercise}
If $X_1,\ldots, X_n$ is a Carnot basis of a Carnot group and $\vol$ is a Haar measure, then for some constant $c$ we have
\\$
\vol\bigl(\{\exp(\sum_{i=1}^n x_iX_i):\ (x_1,\ldots,x_n)\in
A\}\bigr)=c \LL^n(A),\qquad\text{for all Borel sets $A\subseteq\R^n$}.
$
\end{exercise}

\begin{exercise}
For all Borel subsets $A$ of each Carnot group $ G$ equipped with a Haar measure $\vol$, we have
$\vol(\delta_\lambda(A))=\lambda^Q\vol(A)$.
\end{exercise}

%
%

\begin{exercise}[Homogeneous lines are snowflakes]\index{snowflake! -- of a metric space}
Let $G$ be a Carnot group with $V_j$ the $j$-th stratum of its stratification. Let $X\in V_j$. Then, restricting the metric of $G$ on  $\exp(\R X)$ gives a metric space that is isometric to the $1/j$-snowflake of the Euclidean line.
\end{exercise}


\begin{exercise}
Let $G$ be a Carnot group. For $f:\R\to G$ of class $C^1$, we have that $f$ is Pansu differentiable if and only if it is a horizontal curve.
\end{exercise}

\begin{exercise}
Let $G$ be a Carnot group. Every $f: G\to \R$ of class $C^1$ is Pansu differentiable.
\end{exercise}

\begin{exercise}
Let $G, H$ be Carnot groups. A $C^1$ map $f: G\to H $ is Pansu differentiable if and only if its (standard) differential preserves the horizontal bundle.
\end{exercise}

\begin{exercise}
Let $\mathbb H$ be the Heisenberg group equipped with its Carnot structure.
Find an example of a $C^1$ curve $f:\R\to \mathbb H$ that at $0$ is horizontal but it is not Pansu differentiable at $0$. 
\end{exercise}


\begin{exercise}\label{Ex:Egorov} Fill in the details in the following argument to prove Theorem~\ref{Egorov's}.
For $k\in \N$ and $t\in (0,\infty)$, let
$$
E_t(k) := \bigcup_{s\in (t,\infty)} \{x : |f_s(x) - f(x)| > k^{-1}\}.
$$
Then, for fixed $k$, the set $E_t(k)$ decreases as $t$ decreases, 
and $\mu\left(\bigcap_{t \in (0,\infty)} E_t(k) \right)= 0$, 
so since
$\mu(X) < \infty$ we conclude that 
$\mu(E_t(k))\to 0$ as 
$t\to 0$. 
Given $\eta > 0$ and $k\in \N$,
choose $t_k$ so large that $\mu(E_{t_k}(k)) < \eta 2^{- k}$ and let $E = \bigcap_{k\in\N} E_{t_k}(k)$. Then
$\mu(E) < \eta$, and we have $|f_t(x) - f(x)| < k^{-1}$ for $t\in (0, t_k)$ and $x \notin E$.
 Thus $(f_t)_t$ converges to $f$ uniformly on $X\setminus E$.
\end{exercise}
\begin{exercise}[Lebesgue Differentiation Theorem for doubling metric spaces]\label{Ex:LDT}
 If $(X,d,\mu)$ is a doubling metric measure space and $f\in L^1(X,\mu)$, then 
 for $\mu$-almost every $x\in X$ we have
 $$\dfrac{1}{\mu(B(x,r))}\int_{B(x,r)} |f(y)-f(x)| \d \mu (y) \to 0, \qquad \text{ as } r\to 0.$$
 In particular, if $K\subseteq X$ is measurable, then $\mu$-almost every point of $K$ has density 1.
 See \cite[Theorem~1.8]{Heinonenbook}.
\end{exercise}

\begin{exercise}Every sub-Riemannian group of nilpotency step 2 is at a bounded distance from its asymptotic cone. That is, the two metric spaces are \((1, C)\)-quasi-isometric for some \( C>1 \).
\end{exercise}

\begin{exercise}\label{ex: step 2 infinite geodesics}
	The infinite geodesics in sub-Riemannian Carnot groups of step 2 are precisely the left translations of the horizontal one-parameter subgroups.
	\\{\it Hint.} Recall Corollary~\ref{cor infinite geodesics step 2}
\end{exercise}

\begin{exercise}
	Isometric embeddings of Carnot groups into sub-Riemannian Carnot groups of step 2 are affine.
	\\{\it Hint.} Recall Exercise~\ref{ex: step 2 infinite geodesics}, or Corollary~\ref{cor infinite geodesics step 2} for the isometric embeddings of $\R$. The extension to higher dimensional domains can be found in \cite{Hakavuori_2019}. 
\end{exercise}


\begin{exercise}\label{ex12Aug1640}
Let $	\widetilde{\End}:\Omega \longrightarrow G\times \mathbb R$ be the extended endpoint map on a sub-Riemannian group $G$.
A control $u\in \Omega$ is strictly abnormal if and only if 
$\{0 \}\times \R \subseteq \dd 	\widetilde{\End}_u (\Omega) \neq T_{\End(u)} G \times \R.$
\end{exercise}

\begin{exercise}
Let $G$ be a free Carnot group and $H<G$ a Carnot subgroup.
Let $\gamma\subset H$ be a normal curve. Then $\gamma$ is normal as a curve in $G$.
\\
{\it Hint.} Find a normal subgroup $N$ such that $G=N \rtimes H$. Extend the normal covector on $\mathfrak h$ to be 0 on $\mathfrak n$.
\end{exercise}

\begin{exercise}
Let $G$ be a Carnot group of step at most 3. Then all energy minimizing curves in $G$ are smooth.
\\
{\it Hint.} By Goh Theorem~\ref{Goh_thm}, every strictly abnormal curve is contained in a proper subgroup. Hence, every energy minimizing curve is contained in a subgroup within which it is normal.
\end{exercise}

\begin{exercise}
Let $\pi: G\to H $ be a Carnot morphism between Carnot groups.
Let $\gamma:[0,1]\to G$ be a strictly abnormal geodesic. Assume that $\pi\circ \gamma$ is energy minimizing.
Then $\pi\circ \gamma$ is strictly abnormal.
\\
{\it Hint.} Use Exercise~\ref{ex12Aug1640}.
\end{exercise}

\begin{exercise}$\skull$ 
Let $\pi: G\to H $ be a Carnot morphism between Carnot groups.
Let $\gamma:[0,1]\to G$ be a normal geodesic. Assume that $\pi\circ \gamma$ is energy minimizing.
Is then $\pi\circ \gamma$ normal?
\end{exercise}


\begin{exercise}\label{ex_orientation_affine_Heis}
Every affine map of the first Heisenberg group has everywhere positive Jacobian determinant.
\end{exercise}


\chapter{Limits of CC spaces}\label{ch_limits}
\section{Limits of metric spaces}
Sub-Riemannian Carnot groups emerge as limit metric spaces, both as distinguished asymptotic spaces and as tangent metric spaces. In most cases, after some change of coordinates, the study can be reduced to distances that uniformly converge on compact sets. However, it can also be valuable to regard such convergence as a specific instance of Gromov-Hausdorff convergence.
 
\subsection{A topology on the space of metric spaces}

Recall the definition of quasi-isometric embedding from Definition~\ref{def quasi-isometry}.
In the following, by the term {\em pointed metric space}, we mean a pair of a metric space together with a point in it. Also, a map $\phi:(X,x)\to (Y,y)$ between pointed metric spaces is a map $\phi:X\to Y$ such that $\phi(x)=y$.\index{pointed metric space}

\begin{definition}[Hausdorff approximating sequence]
 Let $(X_j,x_j)$ and $(Y_j,y_j)$, for $j\in \N$, be two sequences of pointed metric spaces.
A sequence of maps $\phi_j:(X_j,x_j)\to (Y_j,y_j)$ is said to be \emph{Hausdorff approximating}\index{Hausdorff! -- approximating} if for all $R>0$ and all $\delta>0$ there exists a sequence $\epsilon_j>0$ such that
	\begin{enumerate}
	\item 	$\epsilon_j\to 0$ as $j\to \infty$;
	\item 	$\phi_j|_{B(x_j,R)}$ is a $(1, \epsilon_j)$-quasi-isometric embedding;
	\item 	$\phi_j(B(x_j,R))$ is an $\epsilon_j$-net for $B(y_j,R-\delta)$.
	\end{enumerate}
\end{definition}

\begin{definition}
We say that a sequence of pointed metric spaces $(X_j,x_j)$ {\em converges} to a pointed metric space $(Y,y)$, and that $(Y,y)$ is the {\em limit} of the sequence if there exists a Hausdorff approximating sequence $\phi_j:(X_j,x_j)\to (Y,y)$.\index{Gromov-Hausdorff convergence}\index{convergence! Gromov-Hausdorff --}
\end{definition}
Nowadays, this notion of convergence is called \emph{Gromov-Hausdorff convergence}. 
It was originally introduced by Edwards in 1975, \cite{Edwards_1975}, and then later independently by Gromov in 1981, \cite{Gromov-polygrowth}.
This is a generalization of the Hausdorff convergence from \cite{Hausdorff1914}.
The following criterion will be used to show the convergence of particular sequences of CC spaces.

\begin{proposition}\label{criterion_GH}
 	Let $d_j$ be a sequence of distances on a set $X$ that converge to a distance $d_\infty$ uniformly on bounded sets with respect to $d_\infty$.
	Let $x_0\in X$.
	If
	\begin{equation}\label{eq1639}
	\diam_{d_\infty}\big(\bigcup_{j\in\N}B_{d_j}(x_0,R)\big) < \infty,
	\quad\forall R>0,
	\end{equation}
	then $\id:(X, d_j,x_0)\to(X, d_\infty,x_0)$ is a Hausdorff approximating sequence and $(X, d_\infty,x_0)$ is the limit of $(X, d_j,x_0)$.
\end{proposition}
\begin{proof}
Fix $R, \delta>0$. By \eqref{eq1639} there is $R'>R$ such that $B_{d_j}(x_0, R )\subseteq B_{d_\infty}(x_0, R' )$, for all $j\in \N$.
Define $$\bar \eps_j: = \sup \{ | d_j(X, Y) - d_\infty(X, Y)| : X, Y \in B_{d_\infty}(x_0, R' )\}.$$
Take $\eps_j: = 
	\begin{cases}
	 	\max\{R, 	\bar \eps_j \} & \text{ if } \bar \eps_j\geq \delta \\
		\bar \eps_j & \text{ if } \bar \eps_j< \delta.
	\end{cases}
	$
	
	By assumption of uniform convergence, we have $ \eps_j \to 0$.
	To check that $\id: B_{d_j}(x_0,R)\to(X, d_\infty)$ is a $(1, \epsilon_j)$-quasi-isometric embedding, take $X, Y\in B_{d_j}(x_0,R)\subseteq B_{d_\infty}(x_0, R' )$, then $|d_j(X, Y) - d_\infty(X, Y)| <\bar \eps_j\leq \eps_j$, by the definition of $\eps_j$.

Regarding the fact that $ B_{d_j}(x_0,R)$ is an $\epsilon_j$-net for $B_{d_\infty}(x_0,R-\delta)$, we consider the two cases: Either $\bar \eps_j\geq \delta$ or not. 
In the first case, by definition, we have $\eps_j\geq R$, then, the $\eps_j$-neighborhood, in the metric $d_\infty$, of $x_0$ contains $B_{d_\infty}(x_0,R)$, which obviously contains $B_{d_\infty}(x_0,R-\delta)$.
If, instead, we have $ \bar \eps_j< \delta$, then for all $x\in B_{d_\infty}(x_0,R-\delta)$ we have
\begin{eqnarray*}
d_j(x,x_0) &\leq & d_\infty(x,x_0) +\eps_j\\
&<& R-\delta +\delta = R.
\end{eqnarray*}
So 
$B_{d_\infty}(x_0,R-\delta) \subseteq B_{d_j}(x_0,R)$.
\end{proof}

\begin{example}
	The following example shows that condition \eqref{eq1639} is necessary for the last proposition.
	For $n\in\N$ define $\gamma_n:\R\to\R^2$ by
	\[
	\gamma_n(t) : = 
	\begin{cases}
	 	(t,0) & t\le n \\
		(n,t-n) & n\le t\le n+1 \\
		(n-(t-n-1),1) & n+1 \le t
	\end{cases}
	\]
	\begin{tikzpicture}
 \coordinate (A) at (4,0);
 \coordinate (B) at (4,1);
 \coordinate (C) at (-9,0);
 \coordinate (D) at (-9,1);
 \coordinate (E) at (-1,0); 
 
 \draw[fill=black] (A) circle (2pt) node[below right] {(n,0)};
 \draw[fill=black] (B) circle (2pt) node[above right] {(n,1)};
 \draw[fill=black] (E) circle (2pt) node[below left] {(0,0)};
 
 \draw (A) -- (B);
 \draw (C) -- (A);
 \draw (D) -- (B);
\end{tikzpicture}
\\
	These mappings induce metrics $d_n$ on $\R$ by 
	\[
	d_n(X, Y): =d_{\R^2}(\gamma_n(x), \gamma_n(y)), 
	\qquad\forall X, Y\in\R .
	\]
	Here, as $n\to \infty$, the sequence $d_n(X, Y)$ converges to $d_\infty(X, Y): =|x-y|$, for $X, Y\in\R$. 
	The convergence is uniform on compact sets but not in the Gromov-Hausdorff sense.
	The sequence $(\R, d_n, 0)$ has a Gromov-Hausdorff limit, which however is isometric to $\left(\R\times \{0,1\}, d_{\R^2}, (0,0) \right)$.
\end{example}

\subsection{Asymptotic cones and tangent spaces}
If $X=(X, d)$ is a metric space and $\lambda>0$, we set $\lambda X : = (X, \lambda d)$.
\begin{definition}\index{asymptotic! -- cone}\index{tangent! -- metric space}
 	Let $X, Y$ be metric spaces, $x\in X$, and $y\in Y$.
	We say that $(Y,y)$ is the {\em asymptotic cone} of $X$ if for each infinitesimal sequence $\lambda_j\to0$ we have $(\lambda_j X,x)\to (Y,y)$, as $j\to \infty$.
	We say that $(Y,y)$ is the {\em tangent metric space} of $X$ at $x$ if for each diverging sequence $\lambda_j\to\infty$, $(\lambda_j X,x)\to (Y,y)$, as $j\to \infty$.
\end{definition}

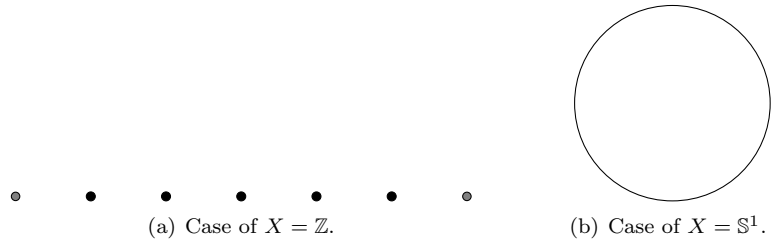
\begin{figure}[h]
\centering
\subfigure[Case of $X=\Z$.] 
{\begin{tikzpicture}

 \coordinate (A) at (-2,0);
 \coordinate (B) at (-1,0);
 \coordinate (C) at (0,0);
 \coordinate (D) at (1,0);
 \coordinate (E) at (2,0); 
 \coordinate (F) at (3,0);
 \coordinate (G) at (4,0); 
 
 \draw[fill=gray] (A) circle (1.6pt) ;
 \draw[fill=black] (B) circle (1.7pt) ;
 \draw[fill=black] (C) circle (1.7pt) ;
 \draw[fill=black] (D) circle (1.7pt) ;
 \draw[fill=black] (E) circle (1.7pt) ;
 \draw[fill=black] (F) circle (1.7pt) ;
 \draw[fill=gray] (G) circle (1.6pt) ; 
 
\end{tikzpicture}}
 \hspace{1.01cm}
\subfigure[Case of $X=\mathbb S^1$.] 
{
\begin{tikzpicture}
 \draw (0,0) circle (1.3cm);
 
\end{tikzpicture}} 
 
 \caption{The asymptotic cone of $\Z$ is $\R$, while at each point the tangent is a singleton.
 The asymptotic cone of $\mathbb S^1$ is a singleton, while at each point, the tangent is $\R$.}
\end{figure} 
\begin{Rem}
 	In general, asymptotic cones and tangent metric spaces may not exist.
 	Within the space of boundedly compact metric spaces, limits are unique up to isometries. 
 	The notion of asymptotic cone is independent from the base point $x$.
\end{Rem}

The following two theorems serve as the central focus of this chapter.
\begin{theorem}[Pansu; \cite{Pansu-croissance}; see Theorem~\ref{Pansu Asymptotic Theorem bis}]\label{Pansu Asymptotic Theorem}\index{Pansu! -- Asymptotic Theorem}\index{Theorem! Pansu Asymptotic --}
 	Let $G$ be a nilpotent simply connected Lie group equipped with a left-invariant sub-Finsler distance.
	Then, the asymptotic cone of $G$ exists and is a sub-Finsler Carnot group.
\end{theorem}
\begin{theorem}[Mitchell; \cite{Mitchell}; see Theorems~\ref{mitchellTheorem bis} and \ref{thm Mitchell mnf}]\label{mitchellTheorem}\index{Mitchell Tangent Theorem}\index{Theorem! Mitchell Tangent --}
 	Let $G$ be a sub-Finsler Lie group or, more generally, an equiregular sub-Finsler manifold and $p\in G$.
	Then, the tangent metric space of $G$ at $p$ exists and is a sub-Finsler Carnot group.
\end{theorem}

\section{Varying CC structures}\label{sec:Varying_bundle}

When taking limits of sub-Riemannian structures, it is important to note that the rank of the distribution may change. This can be observed in the example of the Riemannian Heisenberg group, whose distribution has rank 3, while its asymptotic cone, the sub-Riemannian Heisenberg group, has a distribution of rank 2. 
To study the limits of CC spaces effectively, it is advantageous to adopt the perspective of sub-Finsler structures with possibly varying rank, as in Section~\ref{def_CC_varying_rank}.


\subsection{Definition of structures with parameter}
As in Definition~\ref{rank-varying_structure}, we use the language of bundles to have sub-Finsler structures where the rank could change.
We consider families of CC-bundle structures depending on a parameter.

\begin{definition}[Varying CC-bundle structure]\index{varying CC-bundle structure}\index{CC-bundle structure! varying --}
 Let $\Lambda\subseteq \mathbb R$ be a set. Let $M$ be a smooth manifold. 
 Let $f:\Lambda\times M\times \mathbb R^m\to TM$ and $N:\Lambda\times M\times\mathbb R^m \to [0,+\infty)$ be maps such that for every $\lambda\in\Lambda$ we have that $(f_\lambda,N_\lambda)$ is a CC-bundle structure, where $f_\lambda: =f(\lambda, \cdot, \cdot)$ is the bundle morphism 
 and $N_\lambda: =N(\lambda, \cdot, \cdot)$ is the continuously varying norm, as in Definition~\ref{rank-varying_structure}. We say that the family $\{(f_\lambda,N_\lambda)\}_{\lambda\in\Lambda}$ is a {\em smoothly varying CC-bundle structure} if
 \begin{itemize}
 \item $f\in C^\infty(\Lambda\times M\times\mathbb R^m)$;
 \item $N\in C^0(\Lambda\times M\times\mathbb R^m)$.
 \end{itemize}
\end{definition}
The above definition can be generalized to `continuously varying Lipschitz-vector-fields structures', for which the results in this chapter have analogs; see Section~\ref{sec_general_limits} and \cite{Antonelli_LeDonne_MR4645068}.


\subsection{Divergence bound by Grönwall Lemma}
%

When comparing two sub-Finsler structures, as we will do when studying their convergence,
we need to estimate the displacement of the endpoints of two curves in terms of the difference between their derivatives.
These types of estimates are all consequences of the following general result, the so-called {\em Grönwall Lemma}.\index{Gr\"onwall Lemma}

\begin{lemma}[Grönwall Lemma]\label{lem1053}
 Let $\Omega\subset\R^n$ and $X, Y:\Omega\times[0,T]\to\R^n$. Let $\|\cdot\|$ be 
 a Euclidean norm on $\mathbb R^n$. 
 Suppose that there are $E,K>0$ such that for all $p,q\in\Omega$ and all $t\in[0,T]$
 \begin{equation}\label{eq664dd12e}
 \|X(p,t)-Y(q,t)\|\le E + K \|p-q\| .
 \end{equation}
 Let $\gamma, \eta:[0,T]\to \Omega$ be two absolutely continuous curves such that 
 \begin{equation*}
 \begin{cases}\gamma(0)=\eta(0), &\\
 \dot\gamma(t) = X(\gamma(t),t),& \qquad\text{ for almost every~}t\in[0,T] \text{ and} \\
 \dot\eta(t) = Y(\eta(t),t), &\qquad \text{ for almost every~}t\in[0,T].
 \end{cases}\end{equation*}
 Then
 \begin{equation}\label{eq664dd27f}
 \|\gamma(t)-\eta(t) \| \le E \frac{ e^{Kt}-1 }{K},
 \qquad\forall t\in[0,T] .
 \end{equation}
\end{lemma}
\begin{proof}
 Define $f(t) : = \|\gamma(t)-\eta(t)\|$.
 Notice that $f:[0,T]\to\R$ is absolutely continuous and $f(0)=0$.
 Moreover, for almost every $t\in[0,T]$ we have
 \begin{multline*}
 2f(t) f' (t) 
 = \frac{\dd}{\dd t}(f(t)^2) \\
 = 2 \langle \gamma(t)-\eta(t), \dot\gamma(t)-\dot\eta(t) \rangle 
 \le 2\cdot \|\gamma(t)-\eta(t)\|\cdot\|\dot\gamma(t)-\dot\eta(t)\| \\
 = 
 2f(t)\cdot\|X(\gamma(t),t)-Y(\eta(t),t)\| \\
 \le 2f(t)\cdot (E + K\|\gamma(t)-\eta(t)\|) 
 \le 2f(t)\cdot(E + K f(t)).
 \end{multline*}
 So, whenever $f(t)\neq0$ we have
 \[
 f'(t) \le E + K f(t) .
 \]

 Let $g(t): = e^{-Kt}f(t)$. Then, whenever $g(t)\neq0$, i.e., whenever $f(t)\neq0$, we have
 \[
 g'(t) = -Ke^{-Kt}f(t) + e^{-Kt}f'(t)
 \le -Ke^{-Kt}f(t) + e^{-Kt} (E+Kf(t))
 = e^{-Kt}E.
 \]

 We claim
 \begin{equation}\label{eq1143}
 g(t) \le \int_0^t e^{-Ks}E \dd s, \text{ for almost every }t\in[0,T].
 \end{equation}
 Indeed, if $g(t)=0$, then there is nothing to show because the right-hand side is nonnegative.
 If $g(t)>0$ instead, since $g$ is absolutely continuous, there is a maximal $\hat t< t$ such that $g(\hat t)=0$, and we have 
 \[
 g(t) = g(\hat t) + \int_{\hat t}^t g'(s) \dd s \le \int_{\hat t}^t e^{-Ks}E \dd s \le \int_0^t e^{-Ks}E \dd s, 
 \]
 and~\eqref{eq1143} is proved.

 Finally, we obtain
 \[
 e^{-Kt}f(t) = g(t) \le \int_0^t e^{-Ks}E \dd s = \frac EK (1-e^{-Kt}),
 \]
 which gives \eqref{eq664dd27f}.
\end{proof}


The following is a first application of the Grönwall Lemma that we will use in the next examples on Lie groups.
We use the notation \eqref{def_gamma_prime} for $\gamma'$.

\begin{lemma}[Consequence of Grönwall Lemma]\label{lem2051257}\index{Gr\"onwall Lemma}
 Let $G$ be a Lie group, $\|\cdot\|$ be a norm on $T_{1_G}G$, $d$ a left-invariant Riemannian distance on $G$, and $\nu>0$.
 Then there is $C>0$ such that for all $\epsilon>0$, for all absolutely continuous curves $\gamma, \sigma:[0,1]\to G$ such that $\gamma(0)=\sigma(0)$, $\|\gamma'\|, \|\sigma'\|\le \nu$ a.e., and $\|\gamma'-\sigma'\|<\epsilon$ a.e., then
 \[
 d(\gamma(1), \sigma(1)) \le C\epsilon .
 \]
\end{lemma}
\begin{proof}
 Since $d$ is left-invariant, we assume $\gamma(0)=\sigma(0)=1_G$.
 Let $U$ be a bounded neighborhood of $1_G$ that is a domain of exponential coordinates.
 We claim that there exists $\hat C>0$ such that for all $\epsilon>0$ and all absolutely continuous curves $\sigma, \gamma:[0,1]\to U$ valued into $U$ with $\|\gamma'(t)\|, \|\sigma'(t)\|\le \nu$ and $\|\gamma'(t)-\sigma'(t)\|<\epsilon$ for a.e.~$t\in[0,1]$, we have
 \begin{equation}\label{eq664de259}
 d(\gamma(1), \sigma(1)) \le \hat C\epsilon .
 \end{equation}
 We define the time-dependent vector fields $X$ and $Y$ on $U$ by
 \[
 X(g,t) : = (L_g)_*\gamma'(t),
 \qquad
 Y(g,t) : = (L_g)_*\sigma'(t) .
 \]
 Working in exponential coordinates, we will show an estimate as in~\eqref{eq664dd12e}.
 Since the map $(g,v)\mapsto (L_g)_*v$ is smooth, then it is Lipschitz on the bounded set $U\times B_{\norm{\cdot}}(0, \nu)$.
 Hence there is $K>0$ such that, for all $p,q\in U$ and for a.e.~$t\in[0,1]$,
 \begin{multline*}
 \|X(p,t)-Y(q,t)\|
 = \|(L_p)_*\gamma'(t) - (L_q)_*\sigma'(t)\| \\
 \le \| (L_p)_*\gamma'(t) - (L_p)_*\sigma'(t) \| + \| (L_p)_*\sigma'(t) - (L_q)_*\sigma'(t)\| \\
 \le K\|\gamma'(t) - \sigma'(t)\| + K \|p-q\| 
 \le K\epsilon + K \|p-q\| .
 \end{multline*}
 So we have~\eqref{eq664dd12e} with this $K$ and $E: =K\epsilon$.
 With these constants in~\eqref{eq664dd27f}, 
 we conclude that $\| \gamma(1) - \sigma(1)\| \le C\epsilon$ for some $C$.
 We finally obtain~\eqref{eq664de259} from the fact that $d$ and the norm $\|\cdot\|$ are biLipschitz equivalent on $U$.

 Next, since $U$ is a neighborhood of $1_G$, 
 there is $\hat\nu>0$ such that, if $\gamma:[0,1]\to G$ is an absolutely continuous curve with $\gamma(0)=1_G $ and $\|\gamma'(t)\|<\hat\nu$ for a.e.~$t\in[0,1]$, then $\gamma([0,1])\subset U$.

In order to prove the general statement of the lemma, for curves
$\gamma, \sigma:[0,1]\to G$
 that are not necessarily valued into $U$,
 let $n\in\N$ be such that $(n-1)\hat\nu <\nu \le n\hat\nu$.
 For $j\in\{1, \dots,n\}$ define $\gamma_j$ and $\sigma_j:[0,1]\to G$ by
 \[
 \gamma_j(t) : = \gamma\left(\frac{j-1}{n}\right)^{-1} \gamma\left(\frac{j-1}{n} + \frac{t}{n} \right),
 \qquad
 \sigma_j(t) : = \sigma\left(\frac{j-1}{n}\right)^{-1} \sigma\left(\frac{j-1}{n} + \frac{t}{n} \right).
 \]
 Notice that 
 \[
 \gamma_j'(t) = \frac1n \gamma'\left(\frac{j-1}{n} + \frac{t}{n} \right),
 \]
 and similarly for $\sigma_j'$.
 Therefore, $\|\gamma_j'\| < \nu/n \le \hat\nu$, and $\|\sigma_j'\| < \nu/n \le \hat\nu$, and $\|\gamma_j'-\sigma_j'\| < \epsilon/n$.
 It follows that the images of the curves $\gamma_j$ and $\sigma_j$ are in $U$.
 This shows, in particular, that the images of both $\gamma$ and $\sigma$ are in $U^n$ because both curves are concatenations of left translations of the arcs $\gamma_j$ and $\sigma_j$.
 Since $U^n$ is bounded, there is $L$ such that, for all $p\in U$ and all $X, Y\in U^n$,
 we have $d(xp,yp) \le L d(X, Y)$.

 From the above claim, we have 
 $d(\gamma(1/n), \sigma(1/n)) = d(\gamma_1(1), \sigma_1(1)) \le \hat C\epsilon/n$.
 Next, assume that $j\in\{2, \dots,n\}$ and that $C_{j-1}$ is such that
 $d(\gamma((j-1)/n), \sigma((j-1)/n)) \le C_{j-1} \epsilon$.
 Then, 
 \begin{align*}
 d(\gamma(j/n), \sigma(j/n)) 
 &\le d\left( \gamma\left(\tfrac{j}{n} \right), \gamma\left(\tfrac{j-1}{n}\right) \sigma\left(\tfrac{j-1}{n}\right)^{-1} \sigma\left(\tfrac{j}{n} \right) \right) \\
 &\qquad + d\left( \gamma\left(\tfrac{j-1}{n}\right) \sigma\left(\tfrac{j-1}{n}\right)^{-1} \sigma\left(\tfrac{j}{n} \right), \sigma\left(\tfrac{j}{n} \right) \right) \\
 &= d\left( \gamma_{j-1}(1), \sigma_{j-1}(1) \right) 
 + d\left( \gamma\left(\tfrac{j-1}{n}\right) \sigma_{j-1}(1), \sigma\left(\tfrac{j-1}{n}\right) \sigma_{j-1}(1) \right) \\
 &\le \hat C\epsilon + L d\left( \gamma\left(\tfrac{j-1}{n}\right), \sigma\left(\tfrac{j-1}{n}\right) \right) \\
 &\le (\hat C + L C_{j-1})\epsilon .
 \end{align*}
 By iterating, we obtain that there is a positive number $C$ that depends on $\nu$ but not on $\gamma$ or $\sigma$ such that $d(\gamma(1), \sigma(1)) \le C\epsilon$.
\end{proof}


\subsection{Analytically varying Lie structures}\label{sec Analytically varying Lie structures}

An example of varying CC structures is given by a normed vector space $(V, \norm{\cdot})$ and a one-parameter family of Lie group structures $\star_\lambda$ on $V$, or at least defined locally in a neighborhood of $0$ in $V$. Each product gives rise to a left-invariant extension of the norm $\norm{\cdot}$, seen initially as a norm on $T_{1_{\star_\lambda}}V\simeq V$.
Instead of parametrizing the group structures, we prefer to parametrize the Lie bracket structures of the Lie algebras and consider the group products via the Dynkin formula; see Definition~\ref{Dynkin_product}. 
We will assume that these structures depend analytically on the parameter and call them {\em analytically varying Lie-algebra structures}.
 \begin{lemma}\label{lem664f4e91}\index{Dynkin product}
 Let $V$ be a vector space, $\Lambda\subset\R$ open, and $([\cdot, \cdot]_\lambda)_{\lambda\in\Lambda}$ a family of Lie brackets on $V$.
 Assume that the map $\lambda\mapsto[\cdot, \cdot]_\lambda$ is analytic.
 Denote by $\star_\lambda$ the Dynkin product associated with $[\cdot, \cdot]_\lambda$.
Then, the maps
 \begin{equation}\label{eq6650a1ac}
 (\lambda, X, Y)\mapsto x\star_\lambda y
 \qquad\text{and}\qquad
 (\lambda, X, Y)\mapsto \left.\frac{\dd}{\dd s} x\star_\lambda (sy)\right|_{s=0}
 \end{equation}
 are analytic on some open subset of $\Lambda\times V\times V$ that contains $\Lambda\times\{0\}\times\{0\}$.

 Moreover, if the Lie algebras $(V,[\cdot, \cdot]_\lambda)$ are nilpotent with uniform step and if $\lambda\mapsto[\cdot, \cdot]_\lambda$ is polynomial,
 then both of the above maps are polynomials from $\Lambda\times V\times V $ to $ V$.
\end{lemma}

\begin{proof}
 Recall the alternative formula \eqref{Dynkin Formula2} for the Dynkin product:
 for $X$ and $Y$ in $V$,
 $$
 X\star_\lambda Y = \sum_{n>0}\frac {(-1)^{n-1}}{n} \sum_{ \begin{smallmatrix} {r_i + s_i > 0} \\ {1\le i \le n} \end{smallmatrix}} \frac{
 \left((\ad^\lambda_X)^{ r_1} \circ(\ad^\lambda_Y)^{ s_1}\circ (\ad^\lambda_X)^{ r_2} \circ(\ad^\lambda_Y)^{ s_2} \ldots \circ(\ad^\lambda_X)^{ r_n}\circ (\ad^\lambda_Y)^{ s_n-1}\right)(Y)
 }{r_1!s_1!\cdots r_n!s_n! \; \sum_{i=1}^n (r_i+s_i)}, 
 $$
 where we denote $\ad^\lambda_XY : = [X, Y]_\lambda$.
 We define for $\alpha, \beta\in\mathfrak{gl}(V)$,
 $$
 \mathscr{D}(\alpha, \beta) : =
 \sum_{n>0}\frac {(-1)^{n-1}}{n} \sum_{ \begin{smallmatrix} {r_i + s_i > 0} \\ {1\le i \le n} \end{smallmatrix}} \frac{
 \alpha^{ r_1} \circ\beta^{ s_1}\circ \alpha^{ r_2} \circ\beta^{ s_2} \ldots \circ\alpha^{ r_n}\circ \beta^{ s_n-1}
 }{r_1!s_1!\cdots r_n!s_n! \; \sum_{i=1}^n (r_i+s_i)} .
 $$
 Using the fact that the series in \eqref{serie_GLn} has a positive radius of convergence,
 we deduce that also the above series that defines $\mathscr{D}(\alpha, \beta)$ has a positive radius of convergence.
 Indeed, both series have the same monomials with different coefficients, and the coefficients for $\mathscr{D}(\alpha, \beta)$ are smaller. 
 Therefore, $\mathscr{D}$ is a $\mathfrak{gl}(V)$-valued analytic map defined an open subset $\Omega$ of $ \mathfrak{gl}(V)\times \mathfrak{gl}(V)$ that contains $\{0\}\times\{0\}$.


 Since the function $(\lambda, X)\mapsto \ad^\lambda_X\in\mathfrak{gl}(V)$ is analytic, 
 then the composition $(\lambda, X, Y)\mapsto \mathscr{D}(\ad^\lambda_X, \ad^\lambda_Y)$ and its derivatives are analytic on the open set
 $$
 \left\{(\lambda, X, Y)\in\Lambda\times V\times V : (\lambda, \ad^\lambda_X, \ad^\lambda_Y)\in\Omega\right\}, 
 $$
 which clearly contains $\Lambda\times\{0\}\times\{0\}$.

 If all Lie-algebra structures $[\cdot, \cdot]_\lambda$ are nilpotent of step at most $s$, then $(\ad^\lambda_X)^s = 0$ for all $X\in V$ and $\lambda\in\Lambda$.
 Then, in the previous reasoning, we can truncate the series defining $\mathscr D$ to a finite sum.
 If also $(\lambda, X)\mapsto \ad^\lambda_X$ is polynomial, then the maps in~\eqref{eq6650a1ac} are polynomial.
\end{proof}
 
\begin{remark}\label{rem29May1538}
Let $(\star_\lambda)_{\lambda\in\Lambda}$ be the Dynkin product associated with an analytically varying Lie-algebra structure (as in Lemma~\ref{lem664f4e91}) on a vector space $V$, equipped with a norm $\|\cdot\|$.
For every fixed $\lambda\in\Lambda$ and $y\in V$, the map $p\in V\mapsto \left.\frac{\dd}{\dd s} p\star_\lambda (s y)\right|_{s=0} $ gives a vector field that is left-invariant with respect to $\star_\lambda$. Moreover, 
 for $u:[0,1]\to V$ and $\lambda\in\Lambda$,
 we define the time-dependent vector field
 \begin{equation}\label{def_Xulambda_29May}
 X^{u, \lambda}(t,p) : = \left.\frac{\dd}{\dd s} p\star_\lambda (su(t))\right|_{s=0} .
 \end{equation}
 As a consequence of Lemma~\ref{lem664f4e91} and Grönwall Lemma~\ref{lem1053}, for every compact sets $K\subseteq V$ and $\Lambda'\subset\Lambda$ and every positive constants $C_1$, $ C_2$, and $\ell$ there exists $C$ with the following property:
 Let $\eps>0$, $u,v\in L^\infty([0,1];V)$ with $\|u-v\|_{\infty}< C_1 \eps$, $\|u\|_{\infty}, \|v\|_{\infty}\leq\ell$ and $a,b\in\Lambda'$ with $|a-b| < C_2 \eps$.
 Suppose $\alpha$ and $\beta:[0,1]\to K$ are
 absolutely continuous integral curves of $X^{u,a}$ and $Y^{v,b}$, respectively, with $\alpha(0)=\beta(0)$.
 Then 
 \begin{equation}
 \| \alpha(1) - \beta(1) \| \le C \eps .
 \end{equation}
\end{remark}



\section{Examples of limits}
In this section, we present some examples that show the ideas for the more general strategy for proving Mitchell's and Pansu's theorems. 
The techniques are based on Gr\"onwall estimates and on quantitative Chow's theorems, as in Section~\ref{effective_Chow_Carnot}.
\subsection{The sub-Riemannian roto-translation group}
From the neurogeometry point of view, the most important sub-Riemannian manifold that is not a Carnot group is the rototranslation group; see page~\pageref{roto_transl_neuro} or \cite{Sarti-Citti-Petitot}. For this reason, we present a special case of the proof of Mitchell Theorem~\ref{mitchellTheorem} for such a space.
\begin{theorem}[See Proposition~\ref{thm20may}]\label{thm:tangent:roto}\index{Heisenberg! -- group}\index{rototranslation group} 
 	The tangent metric space of the sub-Riemannian rototranslation group is the sub-Riemannian Heisenberg group.
\end{theorem}

The statement will be restated and proved in Proposition~\ref{thm20may}. 
For the argument, we will use the following result, which implies Chow's theorem and the Ball-Box theorem for the sub-Riemannian roto-translation group.
Moreover, it gives a uniform estimate for sequences of structures.
	 For $t<0$, we denote by $\sqrt{t}$ the value $-\sqrt{-t}$. We denote the Euclidean balls by $B_E$.
\index{quantitative Chow's theorem}
\begin{proposition}[Uniform Chow Theorem]\label{prop131301}
 	Let $X_\lambda,Y_\lambda$ be a pair of vector fields in $\R^3$ that depend smoothly on $\lambda\in[0,1]$.
	Assume $X_\lambda,Y_\lambda,[X_\lambda,Y_\lambda]$ is a frame of $\R^3$ for all $\lambda$.
	For each $p\in \R^3$ and $\lambda\in[0,1]$, consider the map (composition of flows)
	\[
	\Phi_\lambda^p(t_1,t_2,t_3)
	: = (\Phi^{-\sqrt{t_3}}_{Y_\lambda}\circ\Phi^{-\sqrt{t_3}}_{X_\lambda}\circ\Phi^{\sqrt{t_3}}_{Y_\lambda}\circ\Phi^{\sqrt{t_3}}_{X_\lambda}\circ\Phi^{t_2}_{Y_\lambda}\circ\Phi^{t_1}_{X_\lambda})(p), \qquad \forall(t_1,t_2,t_3)\in\R^3.
	\]
	\begin{description}
	\item[\ref{prop131301}.i.] 	The map $\Phi^p_\lambda$ is smooth and $(\dd\Phi^p_\lambda)_0$ has maximal rank and
	the biLipschitz constant of $(\dd\Phi^p_\lambda)_0$ is 
	uniformly bounded for $(\lambda,p)$ on compact sets of $[0,1]\times \R^3$.
\item[\ref{prop131301}.ii.] 	For every $R>0$ there exists $C>0$ such that 
		\[
		\Phi^p_\lambda(B_E(0,Cr)) \supset B_E(p,r), \qquad \forall \lambda\in[0,1], \forall r\in(0,R), \forall p\in B_E(0,R)
 .
		\]
	\item[\ref{prop131301}.iii.] 	If $d_\lambda$ is the sub-Riemannian distance for which $X_\lambda$, $Y_\lambda$ are orthonormal, then for every $R>0$ there exists $C>0$ such that 
		\[
		d_\lambda(p,q) \le C\sqrt{d_E(p,q)}, \qquad \forall p,q\in B_E(0,R), \forall\lambda\in[0,1].
		\]
	\end{description}
\end{proposition}
\begin{proof}[Proof of \ref{prop131301}.i]
We already saw the following strategy in Section~\ref{effective_Chow_Carnot}.
 	We have $(\partial_{t_1}\Phi^p_\lambda)(0) = X_\lambda(p)$, $(\partial_{t_2}\Phi^p_\lambda)(0) = Y_\lambda(p)$, and $(\partial_{t_3}\Phi^p_\lambda)(0) = [X_\lambda,Y_\lambda](p)$.
	Hence, the linear map $(\dd\Phi^p_\lambda)(0):\R^3 \to \R^3$ has rank 3.
	Moreover, for every compact set $U$ there is $C>0$ such that 
	\[
	\|(\dd\Phi^p_\lambda)(x)(v)\|\ge C\|v\|, \qquad \forall x\in U, \forall v\in \R^3, \forall \lambda\in[0,1].
	\]
	By continuity in $\lambda$, we can take $C$ uniform when $\lambda\in[0,1]$.
 
 	\proof[Proof of \ref{prop131301}.ii] It follows from \ref{prop131301}.i and the Inverse Mapping Theorem.
	
	\proof[Proof of \ref{prop131301}.iii] Notice that for some $K>0$
	\begin{align*}
	 	d_\lambda(p, \Phi^p_\lambda(t_1,t_2,t_3)) 
		&\le |t_1|+|t_2|+4\sqrt{|t_3|} \\
		&\le K \sqrt{\|(t_1,t_2,t_3)\|_E}, \qquad \forall t_1,t_2,t_3\in(0,1).
	\end{align*}
	Let $R$ as in \ref{prop131301}.ii, take $p,q\in B_E(0, \frac R2)$ so for $r: =d_E(p,q)$
	\[
	q\in \overline B_E(p,r) \subset \Phi^p_\lambda(\overline B_E(0,Cr)),
	\]
	i.e., there are $t_1,t_2,t_3$ with $\|(t_1,t_2,t_3)\|_E < Cr$ such that $q=\Phi^p_\lambda(t_1,t_2,t_3)$.
	Hence,
	$
	d_\lambda(p,q) \le K\sqrt{\|(t_1,t_2,t_3)\|_E}
	\le K\sqrt{Cr}
	= K\sqrt C \sqrt{d_E(p,q)}.
	$
\end{proof}

Because of Proposition~\ref{criterion_GH}, the next proposition implies Theorem~\ref{thm:tangent:roto}.
\begin{proposition}\label{thm20may}
 	In $\R^3$ with coordinates $X, Y, \theta$ let 
	\begin{equation*}
		\begin{array}{llccll}
	 	X &: = \cos\theta\partial_x+\sin\theta\partial_y &,&& Y &: =\partial_\theta, \\
		X_\infty &: = \partial_x+\theta\partial_y &,&& Y_\infty &: = \partial_\theta, \\
		X_n &: = \cos\frac\theta n \partial_x + n\sin\frac\theta n\partial_y &,&& Y_n &: = \partial_\theta, \qquad\forall n\in\N.
	\end{array}
	\end{equation*}
	Let $d$ (resp. $d_n$, resp $d_\infty$) be the sub-Riemannian distance for which $X, Y$ (resp. $X_n,Y_n$, resp. $X_\infty,Y_\infty$) are orthonormal.
	Then
	\begin{description}
	\item[\ref{thm20may}.i.] 	$(\R^3,n d)$ is isometric to $(\R^3, d_n)$, for each $n\in \N$.
	\item[\ref{thm20may}.ii] 	
		$d_n\to d_\infty$ uniformly on compact sets, as $n\to\infty$.
	\end{description}
\end{proposition}
\begin{proof}[Proof of \ref{thm20may}.i.]
 	The distance $nd$ is the sub-Riemannian distance associated with the orthonormal frame $\frac1nX, \frac1nY$.
	Let 
	$
	\delta_n:(X, Y, \theta)\longmapsto (nx, n^2y, n\theta) .
	$
	Then 
$
	 	\dd\delta_n(\frac1n X) = \cos\theta\partial_x + n\sin\theta \partial_y = X_n\circ\delta_n 
		$ and 
		$\dd\delta_n(\frac1n Y) = Y_n\circ\delta_n $.
	So $\delta_n$ is an isometry between $(\R^3,nd)$ and $(\R^3, d_n)$.
\proof[Proof of \ref{thm20may}.ii.]	
	Fix $R>0$. Take $p,q\in B_{d_\infty}(0,R)$.
	Let $\sigma:[0,1]\to\R^3$ be a $d_\infty$-geodesic from $p$ to $q$, parametrized by constant speed. Thus $\|\dot\sigma\|_\infty<2R$ and for some measurable functions $a$ and $b$ bounded by $2R$ we have
	$
	\dot\sigma = a X_\infty + b Y_\infty.
	$ 
	Let $\gamma$ be a solution of $\dot\gamma = aX_n+b Y_n$ with $\gamma(0)=p$.
	Then
	\begin{eqnarray*}
	|\dot\sigma-\dot\gamma| 
	&\le &|a| |X_\infty\circ\sigma-X_n\circ\gamma| + |b| |Y_\infty\circ\sigma-Y_n\circ\gamma| \\
	&\le &2R (K|\sigma-\gamma|+\|X_\infty-X_n\|_{L^\infty(B_{d_\infty}(0,R))}) \\
	&\le &2RK |\sigma-\gamma| + 2RK \bar\epsilon_n,
	\end{eqnarray*}
	where $\bar\epsilon_n: =\sup_{B_{d_\infty}(0,R)}|X_n-X_\infty|$.
	Notice that $\bar\epsilon_n\to 0$, because $X_n\to X_\infty$ uniformly on compact sets.
	From Grönwall Lemma~\ref{lem1053}, we get
	$
	|\gamma(1)-\sigma(1)| = o(1).
	$
	Then, by Proposition~\ref{prop131301}
	\begin{eqnarray*}
	 	d_n(p,q) 
		&\le& d_n(p, \gamma(1)) + d_n(\gamma(1), \sigma(1)) \\
		&\stackrel{\ref{prop131301}.iii}{\le}& L_{d_n}(\gamma) + C\, \sqrt{|\gamma(1)-\sigma(1)|} \\
		&\le& L_{d_\infty}(\sigma) + o(1) \\
		&=& d_\infty(p,q) + o(1) .
	\end{eqnarray*}
	In particular, we have $d_n(p,1)\le 3R$ for $n$ large enough. One bound has been proved.
	
	Regarding the other bound, we proceed similarly: Let $\gamma$ be a $d_n$-geodesic from $p$ to $q$, $\gamma:[0,1]\to\R^3$ with 
	$
	\dot\gamma = aX_n+bY_n, 
	$
	for some measurable functions $a$ and $b$ bounded by $3R$. 
	Let $\sigma$ be such that $\dot\sigma=aX_\infty+bY_\infty$ and $\sigma(0)=p$. Then as before $|\gamma(1)-\sigma(1)|=o(1)$. We then bound
	\begin{align*}
	 	d_\infty(p,q) 
		&\le d_\infty(p, \sigma(1)) + d_\infty(\sigma(1), \gamma(1)) \\
		&\le L_{d_\infty}(\sigma) + C\, \sqrt{|\sigma(1)-\gamma(1)|} \\
		&\le L_{d_n}(\gamma) + o(1) \\
		&= d_n(p,q) + o(1) .
	\end{align*}
\end{proof}

\subsection{Sub-Riemannian Carnot group as asymptotic spaces of Riemannian groups}

This section aims to present some instances in which sub-Riemannian manifolds appear as limiting objects of Riemannian manifolds. Examples are asymptotic cones, and a general result is Pansu Theorem~\ref{Pansu Asymptotic Theorem}. We will also see that every sub-Riemannian manifold is the limit of some sequence of Riemannian manifolds; see Theorem~\ref{thm_SR_monotone_limit}.

\subsubsection{The Riemannian Heisenberg group}
\index{Heisenberg! -- group} 
\begin{theorem}\label{Asympt_Riem_Heis}
The asymptotic cone of the Riemannian Heisenberg group is the sub-Riemannian Heisenberg group.
 	\end{theorem}
	We explain a stronger result that we will actually prove: Let $X, Y, Z$ be a basis of the Lie algebra of the Heisenberg group $H$ with only relation $[X, Y]=Z$.
	For $n\in\N$, let $d_n$ be the Riemannian distance for which $X, Y, \frac1n Z$ form an orthonormal frame.
	Let $\dcc$ be the sub-Riemannian distance for which $X, Y$ are orthonormal.
	Then we will prove that for all $R>0$ there is $C>0$ such that 
	\begin{equation}\label{eq_sympt_Riem_Heis}
	d_n(p,q)\le \dcc(p,q) \le d_n(p,q) +C\frac{1}{\sqrt{n}}, \qquad \forall p,q\in B_{CC}(1_H,R), \forall n\in \N .
	\end{equation}
	In particular, we have $d_n\to d_\infty$ uniformly on compact sets.
Consequently, by Proposition~\ref{criterion_GH} if $d$ is the Riemannian distance for which $X, Y,Z$ are orthonormal, then the sequence $(H, \frac1n d)$, which is isometric to $(H, d_n)$, converges to $(H, \dcc)$.

\begin{proof}[Proof of Theorem~\ref{Asympt_Riem_Heis}]
 	The fact that $d_n\le \dcc$ is clear since every horizontal curve for $\dcc$ has the same length with respect to $d_n$.
	For the other inequality, take $R>0 $ and $p,q\in B_{CC}(1_H,R)$. 
	Let $\gamma_n:[0,1]\to H$ be a curve from $p$ to $q$ that minimizes the length with respect to $d_n$.
	Decompose $\dot\gamma$ as
	\[
	\dot\gamma(t) = a_1(t)X + a_2(t)Y + a_3(t) Z,
	\]
	with $a_3(t)$ not necessarily $0$.
	Let $\sigma:[0,1]\to H$ be the curve such that $\sigma(0)=p$ and $\dot\sigma(t) = a_1(t)X+a_2(t)Y$.
	Let $\bar q: =\sigma(1)$.
	Let $\eta:[0,1]\to H$ be the curve such that $\eta(0)=\bar q$ and $\dot\eta(t)=a_3(t)Z$.
	
$$	\begin{tikzpicture}

\draw[thick] plot [smooth, tension=1] coordinates {(-1,0) (1,0.5)(2,-0.5) (3,-0.5)};
\draw[thick] plot [smooth, tension=1] coordinates {(-1,0) (1,1.5)(1.9,1) (3,1.5)};
\draw[thick] plot [smooth, tension=1] coordinates {(3,1.5) (3,-0.5)};

\fill (-1,0) circle (2pt) node[left] {$p = \gamma(0) = \sigma(0)$};
\fill (3,1.5) circle (2pt) node[above right] {$q = \gamma(1) = \eta(1)$};
\fill (3,-0.5) circle (2pt) node[below right] {$\bar{q} = \sigma(1)$};

\node at (-0.2,1) {$\gamma$};
\node at (3.4,0.7) {$\eta$};
\node at (1.4,-0.2) {$\sigma$};

\end{tikzpicture}$$

	We claim that
	\begin{equation}\label{3Maggio2021}
	\eta(t) = (L_{\bar q}\circ L_{\sigma(t)}^{-1} ) (\gamma(t)), \qquad \forall t\in[0,1] .
	\end{equation}
	Since 
$	(L_{\bar q}\circ L_{\sigma(0)}^{-1} ) (\gamma(0)) = L_{\bar q}\circ L^{-1}_p(p) = \bar q=\eta(0)
$,
	it is enough to show that 
	\[
	\frac{\dd}{\dd t}\left(L_{\bar q}\circ L_{\sigma(t)}^{-1} \circ\gamma(t))\right) = \dot \eta(t).
	\]
	For doing this, we consider exponential coordinate as in Section~\ref{sec:Heisenberg}:
	\[
	 	\dot\gamma = a_1X+a_2Y+a_3Z
		= \left(a_1,a_2,a_3-\frac{\gamma_2}2a_1 + \frac{\gamma_1}2a_2 \right)
	\]
	and
	\[
	\dot\sigma = \left(a_1,a_2,-\frac{\sigma_2}2a_1+\frac{\sigma_1}2a_2\right).
	\]
	Thus $\gamma_1 (t) =\sigma_1(t)=p_1+\int_0^ta_1$ and $\gamma_2(t)=\sigma_2(t)=p_2+\int_0^ta_2$.
	We have
	\begin{eqnarray*}
	 	\sigma(t)^{-1}\gamma(t) 
		&\stackrel{\eqref{Hprod}}{=}& \left(\gamma_1-\sigma_1, \gamma_2-\sigma_2, \gamma_3-\sigma_3-\frac12(\sigma_1\gamma_2-\sigma_2\gamma_1) \right)\\
		&= &(0,0, \gamma_3-\sigma_3).
	\end{eqnarray*}
	Thus
	\[
	\frac{\dd}{\dd t}\sigma(t)^{-1}\gamma(t) = (0,0, \dot\gamma_3-\dot\sigma_3) = a_3 Z.
	\]
	The claim \eqref{3Maggio2021} is proved, and, in particular, we have 
	$
	\eta(1) = \bar q\bar q^{-1} q = q .
	$
	
	We need to bound the length $L_{d_1}(\eta)$.
	Since $X, Y,Z$ are orthogonal and $\|\frac1n Z\|_n = 1$, we have
	\begin{multline*}
	 	\int_0^1 n\cdot|a_3| 
		= \int_0^1\|a_3Z\|_n 
		\le \int_0^1 \|a_1X+a_2Y+a_3Z\|_n
		= L_{d_n}(\gamma)
		= d_n(p,q)
		\le \dcc(p,q)
		\le 2R.
	\end{multline*}
	Then
	\[
	L_{d_1}(\eta) = \int_0^1\|a_3Z\|_1 = \int_0^1 |a_3| \le \frac{2R}n .
	\]
	Using the Ball-Box Theorem for the sub-Riemannian Heisenberg group, proved in Proposition~\ref{BB_for_Heis}, we bound
	\begin{multline*}
	 	\dcc(\bar q,q) 
		\stackrel{\eqref{eq_BB_thm}}{\le} K d_1(\bar q,q)^{1/2} 
		\le K (L_{d_1}(\eta))^{1/2}
		\le K \left(\frac{2R}n\right)^{1/2} = C\frac{1}{\sqrt{n}}, 
	\end{multline*}
	for the suitable constant $C$.
	Since $\dcc(p, \bar q)\le L_{CC}(\sigma) \le L_{d_n}(\gamma) = d_n(p,q)$, we conclude that
	\[
	\dcc(p,q) \le \dcc(p, \bar q) + \dcc(\bar q,q) \le d_n(p,q) + C\frac1{\sqrt{n}} .
	\]
	We proved \eqref{eq_sympt_Riem_Heis}, and we conclude the theorem recalling Proposition~\ref{criterion_GH}.
\end{proof}

\subsubsection{Asymptotic cones of some Riemannian groups}

\begin{example}\label{teo:limit_dist_gruppi}
 	Let $G$ be a Lie group and let $\Delta\subset TG$ be a bracket-generating left-invariant distribution. 
	Let $\Delta^\perp$ be a left-invariant distribution such that 
$
 	\Delta_{1_G}\oplus\Delta^\perp_{1_G} = T_{1_G}G .
$
	Let $(\langle \cdot, \cdot \rangle_n)_{n\in \N\cup\{0\}}$ be a sequence of left-invariant Riemannian metrics on $G$ such that
	\begin{description}
	\item[\ref{teo:limit_dist_gruppi}.i.] 	 $\Delta_{1_G}$ is orthogonal to $\Delta_{1_G}^\perp$ with respect to $\langle \cdot, \cdot \rangle_n$, for all $n\in \N\cup\{0\}$;
	\item[\ref{teo:limit_dist_gruppi}.ii.]	 $\|X\|_n=\|X\|_0$, for all $n\in \N$, for all $X\in\Delta$; 
		\item[\ref{teo:limit_dist_gruppi}.iii.]	for all $X\notin\Delta$, we have $\|X\|_n\to +\infty$, as $n\to \infty$.
	\end{description}
	Let $\dcc$ be the sub-Riemannian distance associated with $\Delta$ and $\langle \cdot, \cdot \rangle_0$, and $d_n$ the Riemannian distance associated with $\langle \cdot, \cdot \rangle_n$.
	Then $d_n$ converges uniformly on compact sets to $ \dcc $.
	In fact, for all $R>0$ there exists an infinitesimal sequence $\epsilon_n $ such that 
	\begin{equation}\label{distanze_strizzate_sequenza}
	d_n(p,q) \le \dcc(p,q) \le d_n(p,q) + \epsilon_n, \qquad \forall p,q\in B_{CC}(1_G,R).
	\end{equation} 
\end{example}
\begin{proof}
The left-hand side of \eqref{distanze_strizzate_sequenza} 	is obvious from \ref{teo:limit_dist_gruppi}.ii. 
For the right-hand side, begin by noticing that 
the unit tangent bundle of $ \Delta_{1_g}^\perp$
 is compact. Consequently from \ref{teo:limit_dist_gruppi}.iii, there exists a diverging sequence $K_n$ in $\R$ 
 such that 
 \begin{equation}\label{ci_riusciro}K_n\cdot\|Z\|_0\leq \|Z\|_n, \qquad \forall Z\in \Delta_{1_g}^\perp, \forall n\in \N.\end{equation} 
	Take $R>0$ and $p,q\in B_{CC}(1_G,R)$, so $p,q\in B_{d_n}(1_G,R)$.
	Let $\gamma=\gamma_n:[0,1]\to G$ be a constant-speed curve from $p$ to $q$ such that $L_{d_n}(\gamma)=d_n(p,q)$. Consequently, we have $\|\dot\gamma\|_n< 2R$.	
For all	$t\in[0,1]$ we decompose $\gamma': =(L_\gamma)^*\dot\gamma$ as $\gamma'(t)=X(t)+Z(t)$ with $ X(t)\in\Delta_{1_G}$ and $Z(t)\in\Delta^\perp_{1_G}$. 
From \ref{teo:limit_dist_gruppi}.i we know that $Z\perp X$ for every of the Riemannian metrics.
	Let $\sigma:[0,1]\to G$ be a solution of $\sigma(0)=p$ and $\dot\sigma(t)=(L_\sigma)_*X(t)$.
	Then
	\[
	K_n\cdot\|Z\|_0 \stackrel{\eqref{ci_riusciro}}{\le} \|Z\|_n \stackrel{Z\perp X}{\le} \|X+Z\|_n 
	= \|\dot\gamma\|_n <2R .
	\]
	Let $\xi(t): =\diam_{\dcc}(\overline B_{d_0}(1_G,r))$ as in Exercise~\ref{ex8051255}.
	Then we are going to use Grönwall Lemma~\ref{lem2051257} since $\|\gamma'\|_0, \|\sigma'\|_0\le \|\gamma'\|_n<2R$ and $\|\gamma'-\sigma'\|_0=\|Z\|_0 < \frac{2R}K_n$ and get 
	\begin{align*}
	\dcc(p,q) 
	&\le \dcc(p, \sigma(1)) + \dcc(\sigma(1), \gamma(1)) \\
	&\le L_{CC}(\sigma) + \xi(d_0(\sigma(1), \gamma(1))) \\
	&\le L_{d_n}(\gamma) + \xi\left(C\cdot\frac{2R}{K_n}\right) \\
	&= d_n(p,q) + o(1)\quad\text{, as $n\to \infty$,}
	\end{align*}
	where we used that $K_n\to \infty$ and that $\xi(r)=o(1)$ by Exercise~\ref{ex8051255}.
\end{proof}
\begin{corollary}
 	Let $G$ be a stratified group equipped with a Riemannian structure for which the stratification is orthogonal.
	Then, the asymptotic cone of $G$ is a Carnot group.
	In fact, if $d$ is the Riemannian distance, then there exist Riemannian distances $d_\lambda$ on $G$ such that $d_\lambda\to\dcc$ uniformly on compact sets and $(G, \frac1\lambda d)$ is isometric to $(G, d_\lambda)$, and $(G, \dcc)$ is a Carnot group.
\end{corollary}

\subsection{Sub-Riemannian spaces as monotone limits of Riemannian spaces}

In this subsection, we show that every sub-Riemannian distance is a {\em monotone} limit of Riemannian distances. 
This fact has been well-known in nonholonomic geometry for the last 35 years, and it is probably due to V. Gershkovich. 


\begin{theorem}\label{thm_SR_monotone_limit} Every sub-Riemannian distance is an increasing limit of Riemannian distances. \end{theorem}
The argument is easy and well-known when
the sub-Riemannian distribution has constant rank.
 Here, we soon present a proof in the following general case:
As in Definition~\ref{rank-varying_structure}, let $E$ be a vector bundle over $M$ endowed with a scalar product $\langle\cdot, \cdot\rangle$ and let
$\sigma:E\to TM$
be a morphism of vector bundles. For each $p\in M$ and $v,v'\in T_p M$, set 
\begin{eqnarray}
\nonumber
\rho_p(v,v')&: =&\inf\{\langle u,u'\rangle: u,u'\in E_p, \sigma(u)=v, \sigma(u')=v'\}, \\ 
 \rho_p(v)&: =&\rho_p(v,v), \qquad \text{ and } \qquad
 N_p(v): = \sqrt{\rho_p(v)}\label{eq25may1717}
 .
 \end{eqnarray}
The sub-Riemannian distance $d_\rho$ associated with $\rho$ is 
given by \eqref{def_d_25may}.
The only extra assumption on $\rho$ is that the distance $d_\rho$ is finite and induces the manifold topology, as for example, in Definition~\ref{def_constant_rank_sub-Finsler} by Theorem~\ref{Chow2}.

\begin{proof}[Proof of Theorem~\ref{thm_SR_monotone_limit}]
Let $\rho$ and $d: =d_\rho$ 
as in \eqref{eq25may1717} and \eqref{def_d_25may}.
Notice that on each tangent space $T_pM$, with $p\in M$, the value of $\rho(\cdot, \cdot)$ is finite on a subspace, and on this it defines a scalar product. In particular, we have $\rho(v)=0$ only if $v=0$. We take some Riemannian tensor $g_{1}$ with the property that 
$g_{1}\leq \rho.$
Then, by recurrence, for each $m\in \N$, we consider $g_m$ to be a (smooth) Riemannian tensor with the property that, for every $p\in M$ and $v,w \in T_pM $, 
$$\max\{(g_{m-1})_p(v,w), \min\{ (1-2^{-m})\rho_p(v,w),m(g_{1})_p(v,w)\}\}\leq (g_{m})_p(v,w)\leq \rho_p(v,w) .$$
Obviously, we have that
$$g_{1}\leq g_m\leq g_{m+1}\leq\rho.$$
Then, for every absolutely continuous path $\gamma$, we have that
$$L_{g_m}(\gamma)\leq L_{\rho}(\gamma): =\int_0^1 \sqrt{\rho_{\gamma(t)}(\dot\gamma(t))}\dd t .$$
Thus, for every $p$ and $q$ in $M$, 
$$d_{g_m}(p,q)\leq d_\rho(p,q), $$
and therefore
$$\lim_{m\to\infty}d_{g_m}(p,q)\leq d_\rho(p,q) .$$
Assume, by contradiction, that, for some $p$ and $q$ in $M$, we have that
$$\lim_{m\to\infty}d_{g_m}(p,q)< d_\rho(p,q) .$$
Then there are curves $\gamma_m$ from $p$ to $q$ such that
$$\lim_{m\to\infty}L_{g_m}(\gamma_m)< d_\rho(p,q) .$$
Since 
$$L_{g_{1} }(\gamma_m)\leq L_{g_m}(\gamma_m),$$
we get a bound on the lengths $L_{g_{1} }(\gamma_m)$. Therefore, by Ascoli-Arzel\`a argument (Theorem~\ref{AATHM}), the sequence $\gamma_m$ converges to a curve $\gamma$ from $p$ to $q$.
We may assume that $\gamma$ is parametrized by arc length with respect to the distance of $g_1$.
Now, either $L_\rho( \gamma)$ is infinite, or it is finite.
Namely, either there is a positive-measure set $A\subset[0,L_{g_1}(\gamma)]$ such that
$$\rho_{\gamma(t)}(\dot\gamma(t))=\infty, \qquad\forall t\in A,$$
or, for almost every $t\in [0,L_{g_1}(\gamma)]$, the value
$\rho_{\gamma(t)}(\dot\gamma(t))$ is finite.
In the first case,
for all $t\in A$, 
$$(g_m)_{\gamma(t)} (\dot\gamma(t))\geq m \cdot (g_{1} )_{\gamma(t)} (\dot\gamma(t)).$$
From this, we have that
$$L_{g_m}(\gamma)\geq m L_{g_ {1}}(\gamma|_A)\to \infty, \quad \text{ as } m\to \infty.$$
We get a contradiction since by assumption $d_\rho(p,q)<\infty$.
In the second case,
for almost all $t$, for $m$ big enough, 
$$(1-2^{-m}) \rho_{\gamma(t)} (\dot\gamma(t))\leq (g_m)_{\gamma(t)} (\dot\gamma(t))\leq \rho_{\gamma(t)} (\dot\gamma(t)).$$
From this, we have that
$$L_{g_m}(\gamma)\to L_{\rho}(\gamma), \quad \text{ as } m\to \infty.$$
We get a contradiction since we have that 
 $d_\rho(p,q)\leq L_{\rho}(\gamma)$.
 \end{proof}


\section{Pansu Asymptotic Theorem}\label{sec_Pansu_asymp}
In this section, we explain and prove Pansu Theorem \ref{Pansu Asymptotic Theorem}.
Let $G$ be a Lie group with Lie algebra $\mathfrak{g}$. Let $\Delta$ be a left-invariant bracket-generating polarization corresponding to the subspace $\Delta_{1_G} \subset \mathfrak{g}$, as in \eqref{eq:induced:distribution}. Fix a norm $\norma{\cdot}$ on $\Delta_{1_G}$ and transfer it to a left-invariant norm on $\Delta$ as in \eqref{cont_var_norm}. Let $d_{sF}$ be the associated subFinsler metric on $\mathfrak{g}$ as in \eqref{dist_CC}. 
In addition, we assume that $G$ is simply connected and nilpotent. Consequently, the Lie algebra $\g$ is nilpotent, by Proposition~\ref{nilp_gp_nilp_alg}.

\subsection{Pansu asymptotic structure}\label{sec:Pansu asymptotic structure}

Given the nilpotent Lie algebra $\g$, we consider the associated Carnot algebra $\mathfrak{g}_\infty$ of $G$, as in Definition~\ref{Graded:algebra}. 
This is a stratifiable algebra, and $\mathfrak{g}/\lbrack\mathfrak{g}, \mathfrak{g}\rbrack$ is the first stratum of a stratification. 
Consider $G_\infty$ to be the nilpotent simply connected Lie group with Lie algebra $\g_\infty$, by recalling the existence by Birkhoff Theorem; see Proposition~\ref{prop:exists_group_algebra_nilpotent1}.
We shall put a Carnot metric on it.
Consider the projection $\pi: \mathfrak{g} \rightarrow \mathfrak{g}/\lbrack\mathfrak{g}, \mathfrak{g}\rbrack$ modulo $\lbrack\mathfrak{g}, \mathfrak{g}\rbrack$. Since $\Delta_{1_G}$ generates the whole Lie algebra and since $\mathfrak{g}/\lbrack\mathfrak{g}, \mathfrak{g}\rbrack$ is an abelian Lie algebra, we have $\pi(\Delta_{1_G}) = \mathfrak{g}/\lbrack\mathfrak{g}, \mathfrak{g}\rbrack$. We define a norm $\norma{\cdot}_\infty$, on $\mathfrak{g}/\lbrack\mathfrak{g}, \mathfrak{g}\rbrack$, which we name \textit{Pansu limit norm}, by imposing $B_{\norma{\cdot}_\infty}(0,1)=\pi(B_{\norma{\cdot}}(0,1)).$ 
 Let $d_\infty$ be the subFinsler distance associated with $(G_\infty, \mathfrak{g}/\lbrack\mathfrak{g}, \mathfrak{g}\rbrack, \norma{\cdot}_\infty).$ Thus $(G_\infty, d_\infty)$ is a Carnot group. We call $d_\infty$ the \textit{Pansu limit metric} or the \textit{asymptotic metric}.\label{def Pansu limit metric}\index{asymptotic! -- metric}\index{Pansu! -- limit metric}\index{Pansu! -- limit norm}

Pansu Theorem~\ref{Pansu Asymptotic Theorem} will follow from a more quantitative version as obtained in \cite{ Breuillard-LeDonne1}:
\begin{theorem}[\cite{Pansu-croissance, Breuillard-LeDonne1}; see Theorem~\ref{Pansu Asymptotic Theorem tris}]\label{Pansu Asymptotic Theorem bis}\index{Pansu! -- Asymptotic Theorem}\index{Theorem! Pansu Asymptotic --}
Let $(G, d)$ be a nilpotent simply connected Lie group equipped with a sub-Finsler left-invariant metric. Then, the associated Carnot group $G_\infty$ equipped with the Pansu limit metric $d_\infty$ 
is the asymptotic cone of $(G, d)$ and
there is 
a set-wise identification between $G$ and $G_\infty$ such that 
\begin{equation}\label{Pansu Asymptotic Theorem bis eq}
\Big| d(p,q)-d_\infty (p,q) \Big| = O\left(\max\{d(1,p), d(1,q)\}^{1- 1/s}\right), \qquad \text{ as } p,q\rightarrow \infty, \end{equation}
where $s$ denotes the nilpotency step of $G$.
\end{theorem}

\subsection{Structures of contracted metrics}

\subsubsection{Algebra-group(s) identification} \label{sec_BCH}

There is a natural way to identify the underlying vector space of $\g$ and the one of its associated Carnot algebra $\g_\infty$. This identification depends on the choice of 
a compatible linear grading on $\mathfrak{g}$, as defined in Definition~\ref{def_compatible_linear_grading}. 
The natural identification $\mathfrak{g}\simeq \mathfrak{g}_\infty$ has been discussed in Lemma~\ref{limit_bracket}.
Moreover, via the exponential maps, which in this case are global diffeomorphisms, by Theorem~\ref{CG-1.2.1}, we also have identification at the group level: 
\begin{equation}\label{identif18jun}
G \quad \underbrace{\simeq}_{\exp} \quad\g \overbrace{\simeq}^{\text{compatible grading}} \g_\infty \quad\underbrace{\simeq}_{\exp}\quad G_\infty.
\end{equation}
Hence, we may consider the Pansu limit metric $d_\infty$ as a metric on $G$ or on the vector space $\g$. 

\subsubsection{A family of products and metrics}

Consider the dilations $(\delta_\eps)_{\eps\in \R}$ relative to the compatible linear grading as in \eqref{def_dilation_relative_to_grading}, we stress that the family is defined also for nonpositive $\eps$ and it is polynomial in $\eps$. 
Using these maps, we shall modify the Lie group structure of $G$ and its metric into a one-parameter family.
Recall that we are identifying $G$ with $\g$, as in \eqref{identif18jun}.

We modify the Lie bracket of $\g$ as
\begin{equation}\label{def_bracket_Pansu}[X, Y]_\eps: = \delta_\eps[\delta_\eps^{-1} X, \delta_\eps^{-1} Y], \qquad \forall X, Y\in \g, \forall\eps\in \R\setminus\{0\}.
\end{equation}
We stress that for $X, Y\in \g$ the value $[X, Y]_\eps$, as defined for $\eps\neq0$ is a polynomial in $\eps$ and, by \eqref{formula-g-inf}, it extends to $\eps=0$ as
\begin{equation}
\label{eq: limit bracket asymp}
[X, Y]_0 : =\llbracket X, Y\rrbracket_\infty \stackrel{\stackrel{\rm def}{\eqref{formula-g-inf}}}{=} \lim_{\eps\to 0} [X, Y]_\eps,
\end{equation}
where $\llbracket X, Y\rrbracket_\infty $ denotes the Lie bracket on the associated Carnot algebra $\mathfrak{g}_\infty$.

For $\eps\in \R$, we denote by $\star_\eps$ the Dynkin product associated with the nilpotent Lie bracket $[\cdot, \cdot]_\eps$; see Definition~\ref{Dynkin_product}. The group $(\g, \star_0)$ is then the Carnot group with Lie algebra $\g_\infty$; see Proposition~\ref{prop:exists_group_algebra_nilpotent2}.

We also modify the metric $d=:\rho_1$ on $G$ with the
{\em contracted metrics} defined as\index{contracted metrics}
\begin{equation}\label{def_rho_23May} \rho_\eps (p,q) : = \begin{cases} |\eps| d(\delta_\eps^{-1} p, \delta_\eps^{-1} q)& \forall\eps\in \R\setminus\{0\} \\
d_\infty( p, q)& \text{ for } \eps=0
\end{cases}, \qquad \forall p,q\in G,
\end{equation} 
where $d_\infty$ is the Pansu limit metric from page \pageref{def Pansu limit metric}.
We point out that each distance $\rho_\eps$ is $\star_\eps$-left-invariant and, for $\eps\neq0$, it is the sub-Finsler distance induced by the following distribution and norm:
\begin{equation}\label{def_norm_eps2024}
\Delta_1^{(\eps)} : = \delta_\eps (\Delta_1) \qquad \text{ and }
\qquad \norm{v}^{(\eps)} : = |\eps| \cdot \norm{ \delta_\eps^{-1} v}, \quad \forall v\in \Delta_1^{(\eps)};
\end{equation}
see Exercise~\ref{ex3jun20242200}.
Therefore, for each $\eps\neq 0$, each map $\delta_\eps: (\g, \star_1, |\eps| d) \to (\g_\eps, \star_\eps, \rho_\eps)$ is a Lie-algebra isomorphism and an isometry.

\begin{remark}
Regarding the norm $ \norm{\cdot}^{(\eps)} $, we claim that, as soon as we fix a norm $\|\cdot\|_E$ on $\g$, there is a positive number $C$ such that 
\begin{equation}\label{eq3Jun0931}
\|u-\pi( u)\|_E \leq C \eps 
 \|u\|^{(\eps)}, \qquad\forall u\in \Delta_1^{(\eps)}, \forall \eps\in(0,1).
 \end{equation}
Indeed, we assume, without loss of generality, that $\|\cdot\|_E$ is Euclidean, that makes orthogonal the layers 
the compatible linear grading, and that 
$\|\cdot\|_E \leq \|\cdot\|_1$ on $\Delta_1$.
We write $u\in \Delta_1^{(\eps)}$ as $u=u_1+\cdots+u_s$ with $u_k \in V_k$.
Then, for each $k\in \{2, \ldots, s\}$ and $\eps\in(0,1)$,
we bound
\begin{align*}
 \eps^{-1} \|u_k\|_E \leq \eps^{1-k} \|u_k\|_E
 &\le \|u_1 + \eps^{-1} u_2 + \dots + \eps^{1-s} u_s\|_E \\
 &\le \| u_1 + \eps^{-1} u_2 + \dots + \eps^{1-s} u_s \|_1
 = \| \eps \delta_{\eps}^{-1} u \|_1 
 \stackrel{\rm{def}}{=} \|u\|^{(\eps)}.
\end{align*}
Therefore, we obtain \eqref{eq3Jun0931}:
\begin{equation*}
\|u-\pi( u)\|_E \le (\|u_2\|_E + \dots + \|u_s\|_E) \le (s-1) \eps \|u\|^{(\eps)}.
\end{equation*}
\end{remark}

All these sub-Finsler structures, $( \Delta^{(\eps)}, \norm{\cdot}^{(\eps)} )$, make the abelianization map a submetry:
\begin{proposition}\label{submetry for Pansu}
Under the identifications $V_1= G/[G,G]= G_\infty/[G_\infty,G_\infty]$, the following maps are submetries:
\begin{description}
\item[\ref{submetry for Pansu}.i.] The projection map $\pi|_{\Delta_{1_G}} \colon (\Delta_{1_G}, \norma{\cdot}) \rightarrow (V_1, \norma{\cdot}_\infty)$, modulo $[\g, \g]$.
\item[\ref{submetry for Pansu}.ii.] 
The projection map $\pi \colon (G, d) \rightarrow (V_1, \norma{\cdot}_\infty)$, modulo $[G,G]$.
\item[\ref{submetry for Pansu}.iii.] 
The projection map $\pi \colon (G_\infty, d_\infty ) \rightarrow (V_1, \norma{\cdot}_\infty)$, modulo $[G_\infty,G_\infty]$.
\item[\ref{submetry for Pansu}.iv.] 
The projection map $\pi : (\g, \rho_\eps) \to (V_1, \norm{\cdot}_{\infty})$, modulo $[\g, \g]$.
\end{description}
\end{proposition}
\begin{proof}
First, observe that the maps are well defined since the subspace $V_1$ is a complementary subspace of the commutator subalgebra, i.e.,
$ 
\g=[\g, \g]\oplus V_1
$. 
These maps are group homomorphisms.
The map \ref{submetry for Pansu}.i is a submetry by the way the norm $\norma{\cdot}_\infty$ has been defined.
It is clear that the maps \ref{submetry for Pansu}.ii and \ref{submetry for Pansu}.iii are submetries, for example, as a consequence of Proposition~\ref{prop_subFin_submetry_Lie}.
Regarding the map \ref{submetry for Pansu}.iv, the reason is again that the norm in the target is the push forward of the norm in the domain: 
\begin{eqnarray*}
\pi \left(B_{\norm{\cdot}^{(\eps)}}(0,1)\right) &=&
 \pi \left(\delta_\eps B_{\norm{\cdot}^{(1)}}(0,1/\eps)\right)\\
&=&\eps \pi \left(B_{\norm{\cdot}^{(1)}}(0,1/\eps)\right)\\
&=&\eps B_{\norm{\cdot}_{\infty}}(0,1/\eps)\\
&=& B_{\norm{\cdot}_{\infty}}(0,1),
\end{eqnarray*}
where we used the definition of $\norm{\cdot}^{(\eps)}$ and of $\norm{\cdot}_{\infty}$.
\end{proof}

%
%
%

\subsubsection{Guivarc'h apriori bound}

The distance function $
d_\infty(1, \cdot)$ on $G$ is 
a homogeneous quasi-norm with respect to the compatible linear grading, in the sense of Definition~\ref{def_homogeneous_quasi_norm}.
Therefore, we have an immediate consequence of Guivarc'h Theorem~\ref{thm_guiv}:
\begin{corollary}[Guivarc'h]\label{thm_guiv_bis}
 Let $G$ be a nilpotent simply connected Lie group equipped with a left-invariant sub-Finsler metric $d$. Let $d_\infty$ be the Pansu limit metric on $G$. Then there exists a constant $C >1$ such that
 \begin{align}\label{eq_thm_guiv0}
 \frac{1}{C}d_\infty(1,x) - C \leq d(1,x) \leq Cd_\infty(1,x) + C, \quad \forall x \in G.
 \end{align}
 in terms of the distances $\rho_\eps$, as in \eqref{def_rho_23May}, we have
 \begin{align}\label{eq_thm_guiv2}
 \frac{1}{C}\rho_0(1,x) - C \eps \leq \rho_\eps(1,x) \leq C \rho_0(1,x) + C\eps, \quad \forall x \in G, \forall \eps\geq0.
 \end{align}
\end{corollary}



%

\subsection{Proof of Pansu Asymptotic Theorem~\ref{Pansu Asymptotic Theorem}}
%
%
%
%
%
%
%
%

We are ready to state and prove a quantitative version of Pansu Asymptotic Theorem~\ref{Pansu Asymptotic Theorem}, which implies Theorem~\ref{Pansu Asymptotic Theorem bis}.
 \begin{theorem}[Quantitative Pansu Theorem; \cite{ Breuillard-LeDonne1}]\label{Pansu Asymptotic Theorem tris}\index{Pansu! -- Asymptotic Theorem}\index{Theorem! Pansu Asymptotic --}
Let $G$ be a simply connected $s$-step nilpotent Lie group equipped with a left-invariant sub-Finsler metric $\rho_1$.
With respect to some compatible linear grading, consider the Pansu limit metric $\rho_0$ on $G$ and the contracted metrics $\rho_\eps$ as in \eqref{def_rho_23May}. 
For every compact set $K\subseteq G$ there is a constant $C>0$ such that
 \begin{equation}\label{23May1034}
| \rho_\eps (p,q) - \rho_0(p,q)| \leq C \eps^{1/s}, \qquad \forall p,q\in K, \forall \eps\in (0,1).
\end{equation}
\end{theorem}
\begin{proof}
Without loss of generality, because of dilations, we assume $K= B_{\rho_0} (1_G,1)$. Fix $p,q\in K$.
We begin with {\bf a first inequality}: $\rho_0(p,q) \leq \rho_\eps (p,q) + C_1 \eps^{ 1/s}$, for some constant $C_1 \geq 0$ independent on $p$ and $q$.
 Let $C$ be the Guivarc'h constant from Corollary~\ref{thm_guiv_bis}.
 We have the uniform bound:
 \begin{equation}\label{eq22may1308tris}
 \rho_\eps (p,q) \leq \rho_\eps (p,1_G)+\rho_\eps (1_G,q)
 \stackrel{\eqref{eq_thm_guiv2}}{\leq}
 C\rho_0 (1_G,p) +C\eps + C\rho_0 (1_G,q) +C\eps < 4C, 
 \end{equation}
 where we used, in addition to the triangle inequality, the fact that $p$ and $q$ are in the unit $\rho_0$-ball.
 Let $\gamma_\eps : \lbrack 0,1 \rbrack \longrightarrow \mathfrak{g}$ be a $\rho_\eps$-geodesic parametrized by constant $\rho_\eps$-speed connecting $p$ to $q$.
 We check that, independently on $\eps\in (0,1)$, the curve $\gamma_\eps$ is inside the bounded set $B_{\rho_0}(1_G, 6 C^2 +C)$:
 \begin{eqnarray*}\label{eq22may1308bis}
\rho_0 (1_G, \gamma_\eps(t)) 
&\stackrel{\eqref{eq_thm_guiv2}}{\leq}&
C\rho_\eps (1_G, \gamma_\eps(t)) +C\eps 
\\&\leq& 
C ( \rho_\eps (1_G,p) +\rho_\eps (p, \gamma_\eps(t)))+C\eps 
\\
&\stackrel{\eqref{eq_thm_guiv2}}{\leq}&
C (C \rho_0 (1_G,p) + C\eps +\rho_\eps (p,q) )+C\eps
\\
&\stackrel{ \eqref{eq22may1308tris}}{<}&
6 C^2 +C.
 \end{eqnarray*} 
 
 Let $\gamma_0 \colon \lbrack 0,1 \rbrack \longrightarrow \mathfrak{g}$ be the left translation by $p$ of the multiplicative integral in $G_\infty$ of the curve $\pi \circ\gamma_\eps$, as in Proposition~\ref{prop development as projection}.
 Namely, if the curve $\gamma_\eps $ has the control $u\colon \lbrack 0,1 \rbrack \longrightarrow \mathfrak{g}$ for the group structure $\star_\eps$, then $\gamma_0$ is such that $\gamma_0(0)=p$ and it has the control $\pi\circ u$ for the group structure $\star_0$:
$$
\begin{cases} 
\dot\gamma_\eps(t) = (L^\eps_{\gamma_\eps(t)})_* u(t) : = \left. \frac{\dd}{\dd s} \gamma_\eps(t) \star_\eps (su(t)) \right|_{s=0} &, \\
\dot\gamma_0(t) = (L^0_{\gamma_0(t)})_* \pi(u(t)) : = \left. \frac{\dd}{\dd s} \gamma_0(t) \star_0 (s\pi(u(t))) \right|_{s=0} 
& .
\end{cases}
$$
In particular, the curves $\gamma_0$ and $ \pi\circ \gamma_\eps $ have the same length in their respective spaces: $(G_\infty, \rho_0)$ and $(V_1, ||\cdot||_\infty)$.
 Therefore, using in addition that $\pi : (\g, \rho_\eps) \to (V_1, \norm{\cdot}_{\infty})$ is 1-Lipschitz by Proposition~\ref{submetry for Pansu}, we have
 \begin{eqnarray}\nonumber
 L_{\rho_0}(\gamma_0) &=& L_{||\cdot||_\infty}(\pi\circ \gamma_\eps) 
 \\&
 {\leq}& L_{\rho_\eps}(\gamma_\eps)\nonumber
\\& =&\rho_\eps(p,q) \label{eq27may1308}\\
&\stackrel{ \eqref{eq22may1308tris}}{\leq}&
 4 C\nonumber
 \end{eqnarray} 
 and, also, 
 \begin{equation}
 \label{eq29may2024}
 \|u(t)\|^{(\eps)} \stackrel{ \eqref{eq22may1308tris}}{\leq} 4C, \qquad \text{ for a.e.~}t\in [0,1],
 \end{equation}
and
 \begin{eqnarray*}\label{eq22may1308A}
 \rho_0 (1_G, \gamma_0(t)) 
 &\leq & \rho_0 (1_G, p)+ \rho_0 (p, \gamma_0(t))\\
 &\leq & 1 + L_{\rho_0}(\gamma_0)\\
 &\leq & 1+ 4C.
 \end{eqnarray*} 
 We have established that the curves $\gamma_\eps $ and $\gamma_0$ lie in the bounded set $B_{\rho_0} (1_G,1+4C+6C^2)$.



Moreover, fixed an Euclidean norm $\norm{\cdot}_E$, from \eqref{eq3Jun0931}, together with 
\eqref{eq29may2024}, we obtain
\begin{equation}
\label{eq29May1540}
\|u-\pi\circ u\|_E \le (s-1)4 C\eps.
\end{equation}


We stress that $\gamma_\eps$ is an integral curve of the vector field $X^{u, \eps}$ as defined in \eqref{def_Xulambda_29May}. While the curve $\gamma_0$ is an integral curve of $X^{\pi\circ u,0}$.
As a consequence of Grönwall Lemma and the smoothness of the Lie structures on $\eps$, seen in Remark~\ref{rem29May1538}, together with \eqref{eq29may2024} and \eqref{eq29May1540}, we obtain that
$\|\gamma_0(1) - \gamma_\eps(1) \|_E \le \hat C \eps$.\index{Gr\"onwall Lemma} 
Consequently, by the Ball-Box Theorem for Carnot groups (Theorem~\ref{Ball-Box4Carnot}), we have for some constant $C'>0$
\begin{equation}\label{pansuBB}
\rho_0 (\gamma_0(1), \gamma_\eps(1)) \le C' \eps^{1/s}.
\end{equation}

Finally, we obtain the first bound:
\begin{eqnarray*}
 \rho_0 (p,q) 
 &\leq& \rho_0(p, \gamma_0(1)) + \rho_0(\gamma_0(1), \gamma_\eps(1)) \\ 
 &\stackrel{\eqref{pansuBB}}{\leq}& L_{\rho_0}(\gamma_0) + C' \eps^{1/s} \\
 &\stackrel{\eqref{eq27may1308}}{\leq }
 &
 \rho_\eps (p,q) + C' \eps^{1/s}. 
\end{eqnarray*}

\medskip
We consider the {\bf second inequality}: $\rho_\eps (p,q) \leq \rho_0 (p,q) + C_2 \eps^{1- 1/s}$, for some $C_2 \geq 0$. 
Let $\gamma_0 : \lbrack 0,1 \rbrack \longrightarrow \mathfrak{g}$ be a $\rho_0$-geodesic parametrized by constant $\rho_0$-speed connecting $p$ to $q$.

 Let $ \pi\circ \gamma_0: [0,1] \to (V_1, \norm{\cdot}_\infty)$ be the development of $\gamma_0$ and $\gamma_\eps: [0,1] \to G$ be a lift of $\gamma_0$ starting at $p$ that has the same length, which exists since the map $\pi: (G, \rho_\eps) \rightarrow (V_1, \norma{\cdot}_\infty)$ is a submetry; see Proposition~\ref{submetry for Pansu} and Corollary~\ref{lift pf geodesics via submetries}.
 Namely, if the curve $\gamma_0 $ has a control $v\colon \lbrack 0,1 \rbrack \longrightarrow V_1\subseteq\mathfrak{g}$ with respect to the group structure $\star_0$, then, the control of $\gamma_\eps$ with respect to the group structure $\star_\eps$ is a measurable function $ u\colon \lbrack 0,1 \rbrack \longrightarrow\Delta^{(\eps)}\subseteq \mathfrak{g}$ such that $\pi\circ u = v$ and $\norm{u}^{(\eps)}=\norm{v}$. 
 In particular, the curves have the same length in their respective spaces: $(G_\infty, \rho_0)$ and $(G, \rho_\eps)$.
 Therefore, we have
\begin{equation}\label{stessalunghezza28may}
 L_{\rho_\eps} ( \gamma_\eps ) 
=L_{\rho_0}(\gamma_0) =\rho_0(p,q)\leq 2.
\end{equation}
 and in particular 
$
 \|u\|^{(\eps)} = \norm{\pi \circ u}\le 2,$ almost everywhere.
Consequently, from \eqref{eq3Jun0931}, we obtain
\begin{equation}
\label{eq29May1545}
\|u-\pi\circ u\|_E = C'\eps,
\end{equation}
for some positive number $C'$ independent of $\eps$, where $\norm{\cdot}_E$ is an Euclidean norm that we fixed.

We claim that the two curves $\gamma_0$ and $\gamma_\eps$ stay within a bounded set independent of $\eps$. Indeed, regarding the first curve, we obviously have that $\gamma_0(t) \in B_{\rho_0}(1_G,2)$, since $p,q\in B_{\rho_0} (1_G,1)$ and $\gamma_0$ is a geodesic.
Regarding the second curve, we bound:
\begin{eqnarray*}
\rho_0 (1_G, \gamma_\eps(t) ) &\stackrel{\eqref{eq_thm_guiv2}}{\leq}
& C (\rho_\eps (1_G, \gamma_\eps(t)) +\eps )\\
&\leq & C (\rho_\eps (1_G, p)+\rho_\eps ( p, \gamma_\eps(t)) +\eps )\\
&\stackrel{\eqref{eq_thm_guiv2}}{\leq}& C (C\rho_0 (1_G, p)+C\eps + L_{\rho_\eps} ( \gamma_\eps ) +\eps )\\
&\stackrel{\eqref{stessalunghezza28may}}{\leq} & C(C+C +2 +1).
\end{eqnarray*}

We proceed similarly as in the first part of the proof: 
We notice that $\gamma_\eps$ is an integral curve of the vector field $X^{u, \eps}$ as defined in \eqref{def_Xulambda_29May}. While the curve $\gamma_0$ is an integral curve of $X^{\pi\circ u,0}$. 
 Again, we have a bound between the two controls: \eqref{eq29May1545}.
As a consequence of Grönwall Lemma and the smoothness of the Lie structures on $\eps$, as seen in Remark~\ref{rem29May1538}, we obtain that $\|\gamma_0(1) - \gamma_\eps(1) \|_E \le \hat C \eps$, for some constant $\hat C$ independent on $\eps\in (0,1)$. 
Consequently, since the map \eqref{eq6650a1ac} are Lipschitz on compact sets, by Lemma~\ref{lem664f4e91}, and the curves stay in a bounded set, for some constant $L$, we obtain 
 $$ \| \gamma_0(1)^{-1}\star_\eps \gamma_\eps(1) \|
 = \| \gamma_0(1)^{-1}\star_\eps \gamma_\eps(1) - \gamma_0(1)^{-1}\star_\eps \gamma_0(1) \|
 \leq L ( \norm{\gamma_\eps(1)- \gamma_0(1)} )\leq L(1+\hat C)\eps.
$$
Then, using the Ball-Box Theorem for Carnot groups (Theorem~\ref{Ball-Box4Carnot}), we have for some constant $L'>0$
\begin{equation}\label{pansuBB2}
\rho_0 (1_G, \gamma_0(1)^{-1}\star_\eps \gamma_\eps(1)) \le L' \eps^{1/s}.
\end{equation}

%
Consequently, for some constant $C''>0$, we bound the distance between the endpoints:
\begin{eqnarray}
\nonumber
\rho_\eps(\gamma_0(1), \gamma_\eps(1) )
&\leq &
\rho_\eps(1_G, \gamma_0(1)^{-1}\star_\eps \gamma_\eps(1) ) \\
\nonumber
&\stackrel{\eqref{eq_thm_guiv2}}{\leq}&
C\rho_0(1_G, \gamma_0(1)^{-1}\star_\eps \gamma_\eps(1) ) +C\eps \\
\label{eq29May2333}
&\stackrel{\eqref{pansuBB2}}{\leq}&
CL'\eps^{1/s} + C\eps \leq C''\eps^{1/s}.
\end{eqnarray}

 Finally, we can conclude 
 \begin{eqnarray*}
 \rho_\eps(p,q) & \leq & \rho_\eps(p, \gamma_\eps(1)) + \rho_\eps(\gamma_\eps(1), q) \\
 & {\leq} & L_{\rho_\eps}( \gamma_\eps) + \rho_\eps(\gamma_\eps(1), \gamma_0(1)) \\
 & \stackrel{\text{\eqref{eq29May2333}}}{\leq} & L_{\rho_\eps}( \gamma_\eps) + C'' \eps^{1 /s} \\
 &\stackrel{\eqref{stessalunghezza28may}}{=} &\rho_0(p,q) + C'' \eps^{1 /s} .
\end{eqnarray*}
\end{proof}

\section{Mitchell Tangent Theorem}
In this section, we explain and prove Mitchell Theorem~\ref{mitchellTheorem} for sub-Finsler Lie groups. Later, in Section~\ref{sec_Mitchell_mfds}, we will discuss the general case of manifolds.

\subsection{Carnot tangents of sub-Finsler Lie groups}
Let $G$ be a Lie group with Lie algebra $\mathfrak g$.
Let $(\Delta, \|\cdot\|)$ be a bracket-generating left-invariant sub-Finsler structure on $G$, with distance $d$.
Let $V_1: = \Delta_1$. 
The {\em osculating Carnot algebra} associated with the polarized Lie algebra $(\g,V_1)$ is the Lie algebra $\g_0$ given by the direct-sum decomposition\index{osculating! -- Carnot algebra} 
\begin{equation}\label{stratification_Mitchell}
\g_0 : = \bigoplus_{i=1}^\infty V_1^{( i)} /V_1^{( i -1)}, \qquad\text{ 
where }V_1^{( i)}: = \sum_{k=0}^{i-1}\ad_{V_1}^{k}{V_1}\subseteq \g,
\end{equation}
endowed with the unique Lie bracket $[\cdot, \cdot]_0$ that has the property that, if $X\in V_1^{( i)}$ and $Y\in V_1^{(j)}$, the bracket is defined, modulo $V_1^{( i+j-1)}$, as
\begin{equation}\label{Lie_Mitchell}
 \left[X+ V_1^{(i-1)},Y+ V_1^{(j-1)}\right]_0: = [X, Y] + V_1^{(i+j-1)},
 \end{equation}
which is well defined because $[V_1^{( i)}, V_1^{( j)}]\subseteq V_1^{( i+j)}$; see Exercise~\ref{ex_stratification_Mitchell}. 
The Lie algebra $\g_0$ is also known as {\em nilpotentization} of $\g$. But, be aware that if $\g$ is nilpotent but not Carnot, then $\g_0$ will not be isomorphic to $\g$ for any choice of $V_1$. Recall the example in Exercises~\ref{ex N51}.

Consider the simply connected polarized Lie group $(G_0, V_1)$ with Lie algebra $\g_0$.
We stress that $G_0$ is a stratified group with $V_1$ as the first strata.
In fact, the space $V_1$ is both a bracket-generating subspace of $\g$ and of $\g_0$.
If $V_1\subseteq \g$ is equipped with a norm $\norm{\cdot}$, then the pair
$(V_1, \norm{\cdot})$ can be seen either as a sub-Finsler structure on $G_0$ and as a sub-Finsler structure on $G$.
We consider $G_0$ equipped with the sub-Finsler distance $d_0$ induced by $(V_1, \norm{\cdot})$.
The Carnot group $(G_0, d_0)$ is called the {\em osculating Carnot group} of $(G, d)$.\index{osculating! -- Carnot group}
\index{nilpotentization|see {osculating Carnot group}}

Mitchell Theorem~\ref{mitchellTheorem} will follow from a more quantitative version, following the work of Bella{\"{\i}}che \cite[p.69]{bellaiche}.
\begin{theorem}[\cite{Mitchell, bellaiche}; see Theorem~\ref{mitchellTheorem tris}]\label{mitchellTheorem bis}\index{Mitchell Tangent Theorem}\index{Theorem! Mitchell Tangent --}
	Let $(G, d)$ be a Lie group equipped with a left-invariant sub-Finsler metric with respect to an $s$-step polarization.
Then, its osculating Carnot group
 $(G_0, d_0)$ 
 is the tangent metric space of $(G, d)$ at any point
and there is 
a local set-wise identification between $G$ and $G_0$ such that 
\begin{equation}\label{Mitchel Theorem bis eq}
\Big| d(p,q)-d_0 (p,q) \Big| = O\left(\max\{d(1,p), d(1,q)\}^{1+ 1/s}\right), \qquad \text{ as } p,q\rightarrow 1, \end{equation}
\end{theorem}

\subsection{Structure of dilated metrics}
There is a natural way to identify the underlying vector space of $\g$ and the one of its nilpotentization $\g_0$. This identification depends on the choice of a linear grading 
 $\mathfrak g = V_1\oplus\dots\oplus V_s$ 
 such that 
\begin{equation}\label{eq6655cc5c}
V_j \oplus \sum_{i=0}^{j-2}\ad_{V_1}^{i}{V_1} = \sum_{i=0}^{j-1}\ad_{V_1}^{i}{V_1} .
\end{equation}
Linear grading with property \eqref{eq6655cc5c} will be called {\em adapted linear grading}.\index{adapted! -- linear grading}
Via the exponential maps, we also have local identifications near the identity elements: 
\begin{equation}\label{iden7june1208}
G \stackrel{\rm loc}{\simeq}\g \simeq \g_0 \simeq G_0.
\end{equation}


Consider the one-parameter family of dilations $(\delta_\eps)_{\eps\in \R}$ relative to the linear grading as in \eqref{def_dilation_relative_to_grading}, which is polynomial in $\eps$. 
We modify the Lie bracket of $\g$ as 
\begin{equation}\label{def_bracket_Mitchel}
[X, Y]_\eps: = \delta_{1/\epsilon}[\delta_\epsilon X, \delta_\eps Y], \qquad \forall X, Y\in \g, \forall\eps\in \R\setminus\{0\}.
\end{equation}
The just defined bracket \eqref{def_bracket_Mitchel} has a different formula than the one defined in the proof of Pansu's Theorem, as in \eqref{def_bracket_Pansu}. In fact, while in \eqref{eq: limit bracket asymp}, we performed a large-scale asymptotic expansion, next, in \eqref{eq: limit bracket tangent}, we will expand infinitesimally. 

The bracket \eqref{def_bracket_Mitchel} satisfies the following formula: 
\begin{align*}
	[v,w]_\epsilon
	&= \left[ \sum_{i=1}^s v_i, \sum_{j=1}^s w_j, \right]_\epsilon
	= \sum_{i,j} \delta_{1/\epsilon} \epsilon^{i+j} [v_i,w_j] \\
	&= \sum_{i,j} \sum_{k\le i+j}\epsilon^{i+j-k} ([v_i,w_j])_k, 
\end{align*}
where we used the fact that if $v_i\in V_i$ and $w_j\in V_j$, 
then, thanks to~\eqref{eq6655cc5c}, $[v_i,w_j]\in 
[V_1^{( i)}, V_1^{( j)}]\subseteq V_1^{( i+j)}=
V_1\oplus\dots\oplus V_{i+j}$.
It follows that the Lie brackets $[\cdot, \cdot]_\epsilon$ are polynomial and so they extend to $\epsilon=0$:
\begin{equation}
\label{eq: limit bracket tangent}
[X, Y]_0 : = \lim_{\eps\to 0} [X, Y]_\eps, \qquad \forall X, Y\in \g.
\end{equation}
Observe that the identification \eqref{iden7june1208} gives a Lie algebra isomorphism between $(\g, [\cdot, \cdot]_0)$ and $(\g_0, [\cdot, \cdot]_0)$. 

Let $V$ be the vector space underlying $\g$. By Lemma~\ref{lem664f4e91}, the Dynkin formula gives a function $(\epsilon, X, Y)\mapsto x\star_\epsilon y$ that is analytic on a neighborhood of $\R\times \{0\}\times\{0\}$ in $\R\times V\times V$. The group $(\g, \star_0)$ is the Carnot group with Lie algebra $\g_0$; see Proposition~\ref{prop:exists_group_algebra_nilpotent2}.

Let $r_0>0$ small enough so that 
 the exponential map restricted to $B_0(0,r_0)\subset \g$
is an analytic diffeomorphism between $B_0(0,r_0)$ and its image.
Observe that for $\eps\in[-1,1]$, the map $\delta_\eps$ maps $B_0(0,r_0) $ into itself. Also, 
the map $\exp\circ\delta_\eps$ is the exponential map for the local Lie group $(\g, \star_\eps)$ and $\star_\epsilon$ is well defined and analytic on $B_0(0,r_0) \times B_0(0,r_0) $.

We modify the metric $d=:d_1$ on $G$ with the
{\em dilated metrics} defined as\index{dilated metrics}
\begin{equation}\label{def_d_23Maybis} 
d_\eps (p,q) : = \begin{cases} \frac{1}{|\eps|} d(\delta_\eps p, \delta_\eps q)& \forall\eps\in \R\setminus\{0\} \\
d_0( p, q)& \text{ for } \eps=0
\end{cases}, \qquad \forall p,q\in B_0(0,r_0).
\end{equation} 
We point out that each distance $d_\eps$ is $\star_\eps$-left-invariant and, for $\eps\neq0$, it is locally the sub-Finsler distance induced by the polarization $(V_1, \norm{\cdot})$ on $(\g, \star_\eps)$; see Exercise~\ref{ex3jun20242200}.
For each $\eps\neq 0$, each map $\delta_\eps: (\g_\eps, \star_\eps, d_\eps) \to
(\g, \star_1, \frac{1}{|\eps|} d) $ is a Lie-algebra isomorphism and an isometry.
 
 \begin{lemma}[Equicontinuity of distances - 1]\label{lem6655eb231}
 Fixed a Euclidean distance $\|\cdot\|_E$ on $\g$,
 there exists a constant $C$
 such that, for every $p\in B_0(0,r_0)$ and every $\eps\in[-1,1]$, one has
\begin{equation}\label{eq6655e75fb}
 d_\eps(1,p) \leq C \|p\|_E^{1/s}, 
\end{equation}
 and
 \begin{eqnarray}
 \label{eq66b229f1_i}
 d_\eps (1,p) \leq C d_0 (1,p)^{1/s}, \\
 \label{eq66b229f1_ii}
 d_0(1,p) \leq C d_\eps (1,p)^{1/s} .
 \end{eqnarray}
\end{lemma}
\begin{proof}
 The weak Ball-Box Theorem~\ref{corol_weak_BB} can be proved uniform in $\eps$, as seen in Exercise~\ref{ex66b22ab2}.
 Hence, there are 
 constants $C_1,C_2$ such that
 \(
 d_\eps(1,p) \le C_1 \|p\|^{1/s} \le C_2 d_{\eps'}(1,p)^{1/s} \), for all $\eps, \eps'\in [-1,1] $ and all $p \in B_0,(0,r_0)$.
 The equations of the lemma are special cases.
\end{proof}

Given an adapted linear grading, as in \eqref{eq6655cc5c}, of a polarized Lie algebra $\g$ and a norm $\norm{\cdot}$ on $\g$, we define the associated {\em boxes} as\index{box}
\begin{equation}\label{def_box_2024_bis}
\Bx(r) : = \left\{X_1+\cdots + X_n \;:\; X_j \in V_j, \;\norm{X_j}<r^j,
 j\in \{1, \dots, s \}
\right\}, \qquad \forall r\geq 0.
\end{equation}

\begin{theorem}[Ball-Box Theorem for Lie groups]\label{thm66b22d85}
 Let $G$ be a Lie group equipped with a sub-Finsler metric $d_1$ with respect to an $s$-step polarization.
With respect to some adapted linear grading, consider the osculating Carnot distance $d_0$ as in \eqref{def_d_23Maybis}, and boxes as in \eqref{def_box_2024_bis}.
 Then, there are a neighborhood $\Omega$ of $1$ in $G$ and a constant $C>0$ such that
 \begin{equation}\label{eq66b2324b}
 \frac1C d_0(1,p) \leq d_1(1,p) \leq C d_0(1,p), \qquad \forall p\in \Omega.
 \end{equation}
 In particular, there are constants $C>0$ and $R>0$ such that
 \begin{equation}\label{eq66b232e8}
 \exp(\Bx(r/C) )\subset B_{d_1}(1,r) \subset \exp(\Bx(Cr)), \qquad \forall r\in[0,R] .
 \end{equation}
\end{theorem}
\begin{proof}
 Notice that, by the Ball-Box Theorem~\ref{Ball-Box4Carnot} for Carnot groups, the estimates~\eqref{eq66b2324b} and~\eqref{eq66b232e8} are equivalent to
 \begin{equation}\label{eq66b23501}
 \exists C \, \forall \eps\in[-1,1]\qquad
 B_{d_0}(1, \eps/C) \subset B_{d_1}(1, \eps) \subset B_{d_0}(1,C\eps) .
 \end{equation}
 Let $d_\eps$ be the dilated metrics as in \eqref{def_d_23Maybis}.

 On the one hand, by~\eqref{eq66b229f1_ii},
 there is $ C>0$ such that, for all $\eps\in[-1,1]$,
 \[
 B_{d_\eps}(1,1) \subset B_{d_0}(1,C) .
 \]
 Since $B_{d_\eps}(1,1) = \delta_{1/\eps}(B_{d_1}(1, \eps))$ and $B_{d_0}(1,C) =\delta_{1/\eps}(B_{d_0}(1,C\eps) )$, we obtain
 \[
 B_{d_1}(1, \eps) \subset B_{d_0}(1,C\eps), 
 \]
 which is one of the two sides of~\eqref{eq66b23501}.

 On the other hand, by~\eqref{eq66b229f1_i},
 there is $ C>0$ such that, for all $\eps\in[-1,1]$,
 \[
 B_{d_0}(1,1/C) \subset B_{d_\eps}(1,1) .
 \]
 As before, after applying the dilation $\delta_\eps$, we obtain
 \[
 B_{d_0}(1, \eps/C) \subset B_{d_1}(1, \eps), 
 \]
 which is other side of~\eqref{eq66b23501}.
\end{proof}

Applying $\delta_\eps $ to \eqref{eq66b2324b} we immediately obtain another consequence:
\begin{corollary}[Equicontinuity of distances - 2]\label{cor6655eb231}
Under the same assumptions of Theorem~\ref{thm66b22d85},
 there exists a constant $C$
 such that, for every $p\in B_0(1,r_0)$ and every $\eps\in[-1,1]$,
\begin{equation}\label{eq6655d0e3b}
 \frac{1}{C} d_0 (1,p) \leq d_\eps(1,p) \leq Cd_0 (1,p) .
\end{equation} 
\end{corollary}

%
%
%
 
\subsubsection{A geodetically linearly connected subset}
For the proof of Mitchell's theorem, we shall need to consider a special set as follows.
For the value $r_0$ fixed above and the constant $C>1$ from Corollary~\ref{cor6655eb231},
 let $\Omega\subset \g$ be the set 
 \begin{equation}\label{def17Jun}
 \Omega: = B_0(0, r_1), \qquad \text{ with } r_1: =\frac{r_0}{4C^2}.
 \end{equation}
We claim that $\Omega$ is $d_\eps$-geodetically connected within $B_0(0, r_0)$, for all $\eps\in [-1,1]$, in the sense that for every pair $p,q\in \Omega$ there is a curve $\gamma_\eps$ from $p$ to $q$ valued in $B_0(0, r_0)$ such that $d_\eps(p,q) = L_{d_\eps}(\gamma_\eps)$.
For $\eps\neq 0 $, let $\tilde \gamma$ be a geodesic in $(G, d)$ between $\delta_\eps(p)$ and $\delta_\eps(q)$. Let $\gamma: = \delta_\eps^{-1}\circ \tilde \gamma$, which is a curve between $p$ and $q$. We claim that the curve $\gamma$ is well defined and valued in $B_0(0, r_0)$. Indeed, as long as $\gamma(t)$ is in $B_0(0, r_0)$, then in fact it is in $B_0(0,3 C^2 r_1)$, because we have
\begin{eqnarray*}
d_0 (1, \gamma(t) ) 
&\stackrel{\eqref{eq6655d0e3b}}{\le}& C d_\eps (1, \gamma(t) ) \\
&\leq & C( d_\eps (1, p) + L_{d_\eps}(\gamma ) )\\
&=& C( d_\eps (1, p) +d_\eps(p,q)) \\
&\leq & C( 2d_\eps (1, p) + d_\eps (1, q))\\ 
&\stackrel{\eqref{eq6655d0e3b}}{\le}&C^2( 2d_0 (1, p) + d_0 (1, q))\\ 
&<& 3 C^2 r_1< r_0
.
\end{eqnarray*}
Hence, the maximal time $t\in [0,1]$ for which $\gamma(t)\in B_0(0, r_0)$ is $t=1$.

For $\eps=0$, let $\gamma_0$ be a $d_0$-geodesic in $(G_0, d_0)$ between $p$ and $q$. Then, similarly, as above, we bound
$$d_0(1, \gamma_0(t) ) < 3 r_1< r_0.$$
Thus, the curve $\gamma_0$ is $B_0(0, r_0)$-valued.

\subsection{Proof of Mitchell Theorem~\ref{mitchellTheorem}}

We are ready to state and prove a quantitative version of Mitchell Theorem~\ref{mitchellTheorem}, which also implies Theorem~\ref{mitchellTheorem bis}.
\begin{theorem}[Quantitative Mitchell Theorem]\label{mitchellTheorem tris}\index{Mitchell Tangent Theorem}\index{Theorem! Mitchell Tangent --}
Let $G$ be a Lie group equipped with a sub-Finsler metric $d_1$ with respect to an $s$-step polarization.
With respect to some adapted linear grading, consider the osculating Carnot distance $d_0$ and the dilated metrics $d_\eps$ as in \eqref{def_d_23Maybis}.
	Then, there are a neighborhood $\Omega$ of $1$ in $G$ and a constant $C>0$ such that
	\begin{equation}\label{eq6655cdd8}
	| d_\eps (p,q) - d_0(p,q)| \leq C \eps^{1/s}, \qquad \forall p,q\in \Omega, \forall \eps\in [0,1].
	\end{equation}
\end{theorem}
 
\begin{proof}
We take as $\Omega$ the set constructed in \eqref{def17Jun}, which we saw is $d_\eps$-geodetically connected within $B_0(0,r_0)$, for all $\eps\in [-1,1]$.
Let $p,q\in \Omega$.

We begin with the {\bf first inequality}: $d_0(p,q) \leq d_\eps (p,q) + C_1 \eps^{ 1/s}$, for some constant $C_1 \geq 0$. 
Let $\gamma_\eps : \lbrack 0,1 \rbrack \longrightarrow B_0(0,r_0)$ be a $d_\eps$-geodesic parametrized by constant $d_\eps$-speed connecting $p$ to $q$. 
Let $u:\lbrack 0,1 \rbrack \longrightarrow V_1$ be the control of $\gamma_\eps$. 
Let $\gamma_0 : \lbrack 0,1 \rbrack \longrightarrow \mathfrak{g}$ 
be the $d_0$-rectifiable curve starting at $p$ with control $u$, with respect to $\star_0$. 
Namely, similarly, as in the proof of Pansu's theorem, we have 
$$
\begin{cases} 
\dot\gamma_\eps(t) = (L^\eps_{\gamma_\eps(t)})_* u(t) : = \left. \frac{\dd}{\dd s} \gamma_\eps(t) \star_\eps (su(t)) \right|_{s=0} &, \\
\dot\gamma_0(t) = (L^0_{\gamma_0(t)})_* u(t) : = \left. \frac{\dd}{\dd s} \gamma_0(t) \star_0 (su(t)) \right|_{s=0} 
& .
\end{cases}
$$

Notice that
\begin{eqnarray}
	L_{d_0}(\gamma_0) 
	&=& L_{d_\eps}(\gamma_\eps) \nonumber\\
	&=& d_\eps(p,q) \label{eq66573428} \\
	&\le& d_\eps(0,p) + d_\eps(0,q) \nonumber\\
	&\stackrel{\eqref{eq6655d0e3b}}{\le}& C(d_0(0,p)+d_0(0,q)) \nonumber\\
	&\le& 2C r_1 \nonumber
\end{eqnarray}
and thus $d_0(0, \gamma_0(t)) \le d_0(0,p) + L_{d_0}(\gamma_0) \le (2C+1) r_1 < r_0$,
that is, $\gamma_0([0,1])\subset B_0(0,r_0)$.

By Remark~\ref{rem29May1538}, we obtain $\|\gamma_0(1) - \gamma_\eps(1)\|_E \le C' \eps $. 
Consequently, by the Ball-Box Theorem for Carnot groups (Theorem~\ref{Ball-Box4Carnot}) we have for some constant $C''>0$
\begin{equation}\label{eq6655d196}
	d_0(\gamma_0(1), \gamma_\eps(1) ) \le C'' \eps^{1/s} .
\end{equation}

Finally, we obtain the first bound:
\begin{eqnarray*}
	d_0(p,q)
	&\leq& d_0(p, \gamma_0(1)) + d_0(\gamma_0(1), \gamma_\eps(1)) \\
	&\stackrel{\eqref{eq6655d196}}\leq& L_{d_0}(\gamma_0) + C'' \eps^{1/s} \\
	&\stackrel{\eqref{eq66573428}}{=}& d_\epsilon(p,q) + C'' \eps^{1/s} .
\end{eqnarray*}

We consider the {\bf second inequality}: $d_\eps (p,q) \leq d_0 (p,q) + C_2 \eps^{ 1/s}$, 
 for some $C_2 \geq 0$. 
Let $\gamma_0 : \lbrack 0,1 \rbrack \longrightarrow B_0(0,r_0)$ be a $d_0$-geodesic parametrized by constant $d_0$-speed connecting $p$ to $q$. 
Let $u:\lbrack 0,1 \rbrack \longrightarrow {V_1}$ be the control of $\gamma_0$, 
 with respect to $\star_0$.
 Let $\gamma_\eps : \lbrack 0,1 \rbrack \longrightarrow \mathfrak{g}$ 
be the $d_\eps$-rectifiable curve starting at $p$ with control $u$ with respect to $\star_\eps$. 
The error between $\gamma_\eps(1)$ and $q$ is estimated as
\begin{eqnarray*}
	d_\eps(\gamma_\eps(1),q)
	&= &d_\eps(1,(-\gamma_\eps(1))\star_\eps q) \\
	&\stackrel{\eqref{eq6655e75fb}}\le &C \|(-\gamma_\eps(1))\star_\eps q\|_E^{1/s} \\
	&= &C \|(-\gamma_\eps(1))\star_\eps q - (-q)\star_\eps q\|_E^{1/s} \\
	&\leq &C' \|-\gamma_\eps(1)+q\|_E^{1/s} \\
	&\stackrel{\rm Rem.~\ref{rem29May1538}}{\leq} &C'' \eps^{1/s} .
\end{eqnarray*}

Notice also that 
$d_\eps(p, \gamma_\eps(1))\leq L_{d_\eps}(\gamma_\eps) \leq L_{d_0}(\gamma_0) = d_0(p,q)$.
Therefore, we obtain
\begin{align*}
	d_\eps(p,q)
	&\leq d_\eps(p, \gamma_\eps(1)) + d_\eps(\gamma_\eps(1),q) \\
	&\leq d_0(p,q) + C'' \eps^{1/s}. 
\end{align*}
\end{proof}


\section{Mitchell's theorem for CC spaces}\label{sec_Mitchell_mfds}
Mitchell's result from the previous section holds more generally. In fact, 
Carnot groups are the tangent metric spaces of sub-Finsler manifolds at every regular point. 
Such a result, originally attributed to Mitchell, is quite technical and involved; a complete treatment is in \cite{bellaiche, jeancontrol}. 
In this section, we describe the Carnot group that appears as tangent at one such point.

\subsection{Osculating Carnot group, a.k.a. nilpotentization}
 Let $M$ be a manifold and let $\Delta $ be a bracket-generating distribution that is equi-regular 
in the sense of Section~\ref{sec:equiregular}.\index{equiregular! -- distribution} 
Denote by $s$ the step of $\Delta$.
Let 
$$\Delta=\Delta^{[1]}\subset\Delta^{[2]}\subset\ldots\subset\Delta^{[s]}=TM$$
be the flag of sub-bundles of $TM$ as in Definition~\ref{def:Delta_k}.
The simple but crucial fact is that 
\begin{equation}\label{Deltas-bracketed}
[\Delta^{[k]}, \Delta^{[l]}]\subseteq\Delta^{[k+l]}.\end{equation}
Equation \eqref{Deltas-bracketed} is obvious for $k=1$ and can be proved by induction using linearity and Jacobi identity:
\begin{eqnarray*}
[\Delta^{[k+1]}, \Delta^{[l]}]&=& \left[\Delta^{[k]} +[\Delta, \Delta^{[k]}], \Delta^{[l]}\right]\\
&=& [\Delta^{[k]}, \Delta^{[l]}] +\left[[\Delta, \Delta^{[k]}], \Delta^{[l]} \right]\\
&\subseteq&\Delta^{[k+l]}+ \left[[\Delta^{[k]}, \Delta^{[l]}], \Delta \right]+ \left[[\Delta^{[l]}, \Delta], \Delta^{[k]} \right]\\
&\subseteq&\Delta^{[k+l]}+[\Delta^{[k+l]}, \Delta]+[\Delta^{[l+1]}, \Delta^{[k]}]\\
&\subseteq&\Delta^{[k+l]}+\Delta^{[k+l+1]}+\Delta^{[k+l+1]}\\
&\subseteq&\Delta^{[k+l+1]},
\end{eqnarray*}
where we assumed that for given $k$ it has been proved for all $l$ and we showed for the value $k+1$ and every $l$.
This last argument is similar to the one for Exercises~\ref{ex:v2v2} and \ref{ex stratifications are gradings}.

Define $H_1: =\Delta$ and $H_j: =\Delta^{[j]}/\Delta^{[j-1]}$, for $j=2, \ldots,n$.
Still, the space $H_j$ is a bundle over $M$, but not a sub-bundle of the tangent bundle $TM$.
We obviously have the isomorphism
$$TM \simeq H_1\oplus H_2\oplus\ldots\oplus H_s.$$

With the aim of defining a Lie group that we will denote by $T_p(M, \Delta)$, for $p \in M$, we set 
\begin{equation}\label{def_Lie_nilp}
{\rm Lie}(T_p(M, \Delta)) : = \bigoplus_{j=1}^s \faktor{\Delta^{[j]}(p)}{\Delta^{[j-1]}(p)}.
\end{equation} We equip this set with the Lie bracket that has the property that for all $x\in \Delta^{[i]}(p) $ and $y\in \Delta^{[j]}(p)$, 
\begin{equation}\label{def_brack_nilp}
\left [x+\Delta^{[i-1]}(p), y+\Delta^{[j-1]}(p) \right ] : = [X, Y]_p + \Delta^{[i+j-1]}(p), 
\end{equation}
where $X\in\Gamma(\Delta^{[i]})$ with $X(p)=x$,
and $Y\in\Gamma(\Delta^{[j]})$ with $Y(p)=y$.
We show that the bracket \eqref{def_brack_nilp} is well defined: by symmetry, we just show the independence on the representative of $x$.
Let $\tilde X\in\Gamma(\Delta^{[i]})$ with $\tilde X(p)\in x+\Delta^{[i-1]}(p)$.
Write $X-\tilde X = \sum_{\ell=1}^k a^\ell f_\ell$ with $f_1, \dots,f_k$ a frame of $\Delta^{[i]}$
and $a^\ell(p)=0$ when $f_\ell\not\in\Gamma(\Delta^{[i-1]})$.
Then
\[
[X, Y] - [\tilde X, Y]
= [X-\tilde X, Y]
= \sum_{\ell=1}^k[a^\ell f_\ell, Y]
= \sum_{\ell=1}^ka^\ell [f_\ell, Y] - \sum_{\ell=1}^k(Ya^\ell) f_\ell, 
\]
which, when evaluated in $p$, gives
\[
[X, Y]_p - [\tilde X, Y]_p
\in [ \Gamma(\Delta^{[i-1]}_p), \Gamma(\Delta^{[j]}) ] + \Delta^{[i]}(p)
\subset \Delta^{[i+j-1]}(p) .
\]
Thus, the Lie bracket \eqref{def_brack_nilp} is well defined.

\begin{remark}
The Lie algebra ${\rm Lie}(T_p(M, \Delta))$ is stratified by
$ \left(\faktor{\Delta^{[j]}(p)}{\Delta^{[j-1]}(p)} \right)_{j\in\{1, \ldots,s\}}$.
The first layer is $\Delta_p$. There is a Carnot group with this stratified algebra by Proposition~\ref{prop:exists_group_algebra_nilpotent1}.
\end{remark}
\begin{definition}[Osculating Carnot group]\index{Carnot! osculating -- group}
When an equiregular polarized manifold $(M, \Delta)$ is equipped with a sub-Finsler norm $\|\cdot\|$, we consider the Carnot group $(T_p(M, \Delta), \Delta_p, \|\cdot\|_p)$ whose Lie algebra is 
${\rm Lie}(T_p(M, \Delta))$ as in \eqref{def_Lie_nilp},
and call it the \emph{osculating Carnot group}
\index{osculating! -- Carnot group}
of $(M, \Delta, \|\cdot\|)$ at $p$, and it is considered equipped with its Carnot-Carath\'eodory metric $d_0$.
\end{definition}

\subsection{Tangent metric spaces to equiregular CC spaces}
Mitchell's theorem for equiregular sub-Finsler manifolds is then the following.
\begin{theorem}[Mitchell]\label{thm Mitchell mnf}\index{Mitchell Tangent Theorem}\index{Theorem! Mitchell Tangent --}
 Let $(M, \Delta, \|\cdot\|)$ be a sub-Finsler manifold, with $\Delta$ equiregular.
 Let $p\in M$.
 Then, the tangent metric space of $(M, d)$ at $p$ is isometric to the osculating Carnot group $(T_p(M, \Delta), d_0)$.
\end{theorem}

We provided a full proof of this theorem for Lie groups; see Theorems~\ref{mitchellTheorem bis} and \ref{mitchellTheorem tris}. The general strategy is very similar. Again, one can see this result as a convergence of a one-parameter family of CC structures. We state the key result referring to \cite{jeancontrol} for a proof. With the terminology of Bella{\"{\i}}che and Jean, the next proposition says that exponential local parametrizations are examples of privileged coordinates. 
\begin{proposition}
 Let $M$ be a manifold equipped with an equiregular bracket-generating distribution $\Delta$.
 Fix $p\in M$.
 Let $X_1, \dots, X_n$ be a local frame of $TM$ around $p$ adapted to the flag $\Delta^{[1]}, \dots, \Delta^{[k]}, \dots$.
 Let $\Phi: U\subset\R^n\to M$ be the corresponding exponential local parametrization around $p$.
 For each $j\in\{1, \dots,n\}$ and $\eps\in(0,1)$, set
 \[
 X_j^{(\eps)} : = \eps^{{\rm deg}(j)} (\Phi\circ\delta_\eps)^* X_j, 
 \]
 where ${\rm deg}(j)$ is such that $X_j\in\Gamma(\Delta^{{\rm deg}(j)}) \setminus \Gamma(\Delta^{{\rm deg}(j)-1})$
 and $\delta_\eps:\R^n\to\R^n$ is the map
 $$\delta_\eps(x_1, \dots,x_n) : = (\eps x_1, \dots, \eps^{{\rm deg}(j)}x_j, \dots).$$
 Then, in the $C^\infty$-topology, each vector field $X_j^{(\eps)}$ converges to $X_j^{(0)}\in{\rm Vec}(\R^n)$ such that $X_1^{(0)}, \dots, X_n^{(0)}$ form a frame of $T(\R^n)$
 and generates a Lie algebra 
 that is isomorphic to ${\rm Lie}(T_p(M, \Delta))$.
\end{proposition}

\begin{remark}
 Differently from the Riemannian case, it is not true that sub-Riemannian manifolds are locally biLipschitz equivalent to their tangent metric spaces, even in the case when the tangent metric space is the same Carnot group at every point; see \cite{LOW}. It is, however, true for contact manifolds because of the Darboux Theorem. See also \cite{LeDonne_Young} for a measure-theoretic parametrization.
\end{remark}

\section{A general result on convergence of CC structures}\label{sec_general_limits}

Pansu's and Mitchell's theorem, together with the other examples that we have discussed, follow from a more general principle. In brief, when varying CC-bundle structures converge, the distances converge. 
The following theorem has been obtained in \cite[Appendix~C]{Antonelli_LeDonne_MR4645068}.
\begin{theorem}\label{thm:VIDEO}
	Let $\Lambda\subseteq\mathbb R$, 
	and let $\{(f_\lambda,N_\lambda)\}_{\lambda\in\Lambda}$ be a varying CC-bundle structure on a manifold $M$. Let $d_{\lambda}: =d_{(f_\lambda,N_\lambda)}$ for every $\lambda\in\Lambda$, as in \eqref{def_dist_bundles}. 
	Let $\lambda_0\in\Lambda$ be such that $f(\lambda_0,M\times\mathbb R^m)$ is a bracket-generating distribution and the metric space $(M, d_{\lambda_0})$ is boundedly compact. 
	Then $d_\lambda\to d_{\lambda_0}$ uniformly on compact sets of $M$ as $\lambda\to\lambda_0$. 
\end{theorem}

The proof of the previous theorem stands on the following crucial lemma, which should be compared with \eqref{eq_thm_guiv2} and \eqref{eq66b229f1_i}. 
\begin{lemma}[Equicontinuity of the distances]\label{LEMMACRUCIALE}
In the same assumptions of Theorem~\ref{thm:VIDEO}, 
	let $K\subseteq M$ be
	 compact set and 
	 $\rho$ 
	 Riemannian metric on $M$. 
	Then there exists a neighborhood $I_{\lambda_0}\subseteq\Lambda$ of $\lambda_0$, and $\beta$ homeomorphism of $[0,+\infty)$ such that
	$$
 d_\lambda(p,q)\leq \beta(\rho(p,q)), \qquad \text{for all $p,q\in K$ and $\lambda\in I_{\lambda_0}$}.
	$$
\end{lemma}

\begin{proof}[Sketch of the proof of Lemma~\ref{LEMMACRUCIALE}]
 One begins by showing that for all $x\in M$ and all $\eps>0$
 there exist $\delta>0$ and a neighborhood $I$ of $\lambda_0$ such that
 \begin{equation}\label{eq66718703}
 B_\rho(x, \delta) \subseteq B_{d_\lambda}(x, \eps),
 \qquad \forall \lambda\in I.
 \end{equation}
 As in the proof of Chow's Theorem~\ref{Chow}, we consider maps
 \[
 F_\lambda : (t_1, \dots,t_N)
 \mapsto
 \Phi^{t_N}_{X^\lambda_N} \circ \cdots \circ \Phi^{t_1}_{X^\lambda_1} (x)
 \]
 obtained as composition of flows starting at $x$ of vector fields in $f_\lambda$.
 Since $f_{\lambda_0}$ is bracket generating by assumption, from the proof of Chow's Theorem we have that
 $F_{\lambda_0}([0, \eps/N]^N)$ is a neighborhood of $x$.
 By smoothness of $F_\lambda$ in $\lambda$, the same is true for $\lambda$ in some neighborhood of $\lambda_0$.
 One infers~\eqref{eq66718703} and concludes by compactness arguments.
\end{proof}

\section{Finitely-generated groups of polynomial growth}\label{sec_largescale}

 Carnot-Carathéodory spaces appear when studying the asymptotic growth of finitely generated nilpotent groups, which essentially are the finitely generated groups of polynomial growth; see the theorems below.
 
 A group $\Gamma$ is {\em finitely generated}\index{finitely generated} if there is finite subset $S\subseteq\Gamma$, called {\em finite generating set}\index{generating! -- set} such that 
 $$\Gamma = \bigcup_{n\in \N} (S \cap S^{-1})^n.$$
In such a case, fixed a finite generating set $S\subset \Gamma$, then the distance function $d_S$ on $ \Gamma$ that is defined by
 $$\bar B_{d_S} (x, n ) = x (S\cap S^{-1})^n, \qquad \forall x\in \Gamma, \forall n\in \N$$
 is called the {\em word metric}\index{word metric} with respect to $S$

In addition, we say that a finitely generated group $\Gamma$ has {\em polynomial growth}\index{polynomial! -- growth} if for some (and consequently, for every) finite generating set $S$ one has
$$ {\rm Card} \left((S\cap S^{-1})^n\right)\leq c_1 n ^{c_2}, \qquad \forall n\in \N,$$
for some constants $c_1, c_2\in \R_+$, possibly depending on $S$. Here, the operator $ {\rm Card} $ denotes the cardinality of sets.

Sub-Finsler geometries will appear in the following two results:

\begin{theorem}[Wolf \cite{WolfMR0248688}, Bass \cite{MR0379672}, Pansu, \cite{Pansu-croissance}]\label{thm_WBP}\index{Theorem! Wolf-Bass-Pansu --}\index{Wolf-Bass-Pansu Theorem}\index{Theorem! Pansu --}\index{Pansu! -- Theorem}
 Let $\Gamma$ be a nilpotent finitely generated group, with a finite generating set $S\subset \Gamma$ and word distance $d_S$.
 \begin{description}
 \item[\eqref{thm_WBP}.i.] The asymptotic cone of $(\Gamma, d_S)$ exists and is a sub-Finsler Carnot group, whose Hausdorff dimension is some $Q\in \N$ and the $Q$-Hausdorff measure of the unit ball is some $v\in \R_+$;
 \item[\eqref{thm_WBP}.ii.] The group $\Gamma$ has polynomial growth, and in fact one has
 \begin{equation}\label{Pansu_expansion}
 \dfrac{{\rm Card} \left((S\cap S^{-1})^n\right)}{ v n ^Q}\longrightarrow 1, \qquad \text{as } n\to \infty.
 \end{equation}
 \end{description} 
 \end{theorem}
 
 The following theorem is a vice versa. 
 
 \begin{theorem}[Gromov's polynomial growth, \cite{Gromov-polygrowth}]\label{Gromov_Theorem}\index{Theorem! Gromov polynomial growth -- }\index{polynomial! -- growth}
 Let $\Gamma$ be a finitely generated group. 
If $\Gamma$ has a polynomial growth rate, then there is a subgroup $\Gamma_{\rm nil}\subset \Gamma$ that is finitely generated and nilpotent and the quotient $\Gamma / \Gamma_{\rm nil}$ is finite.
In particular, the asymptotic cone of $(\Gamma, d_S)$ is a sub-Finsler Carnot group for every finite generating set $S\subset \Gamma$ and word distance $d_S$.
\end{theorem}
The two above theorems are summarized by the following sentence: A finitely generated group has a polynomial growth rate if and only if it is virtually nilpotent. Actually, these are precisely the finitely generated groups that have Lie groups as asymptotic cones.

\subsection{A few comments on Wolf-Bass-Pansu's Theorem}
For the proof of Theorem~\ref{thm_WBP}, one can construct the asymptotic Carnot group explicitly.
Starting with a nilpotent finitely generated group $\Gamma$, equipped with a word distance $d_S$, we describe a collection of metric groups that are quasi-isometric to $(\Gamma, d_S)$, and actually, they are asymptotic to each other in the sense that, after a coarse identification, the ratio of the distances goes to 1 as the points go to infinity.
One of these metric spaces will be a nilpotent simply connected sub-Finsler Lie group. Hence, we will conclude, with the help of Pansu's Theorem~\ref{Pansu Asymptotic Theorem bis}, that this Lie group and, hence, $\Gamma$ are asymptotic to the associated Carnot group. Be aware that this final Carnot group may not be quasi-isometric to $\Gamma$.

After $\Gamma$, the next group in the collection that we consider is $\Gamma$ modulo its torsion elements 
$\mathrm{Tor}(\Gamma)$, as we discussed in Section~\ref{sec:torsion}.
One can prove that the set $\mathrm{Tor}(\Gamma)$ is a normal subgroup that is finite; see \cite[Theorem 9.17]{Macdonald}.
Let $\pi: \Gamma \to \Gamma/\mathrm{Tor}(\Gamma)$ be the quotient projection.
The group $\pi(\Gamma)= \Gamma/\mathrm{Tor}(\Gamma)$ is nilpotent, torsion-free, and finitely generated by $\pi(S)$.
Via the map $\pi$ the metric spaces $(\Gamma, d_S)$ and $( \pi(\Gamma), d_{\pi(S)})$ are $(1,C)$-quasi isometric.

By a fundamental result by Malcev \cite{Malcev}, called {\em Malcev completion},\index{Malcev! -- completion} there is a (unique) nilpotent simply connected Lie group $G$ that admits a discrete cocompact subgroup $\tilde \Gamma \subset G$ isomorphic to $\Gamma/\mathrm{Tor}(\Gamma)$; see \cite[Theorem 2.18]{Raghunathan}. 
Here, the cocompact means that the quotient $G/ \tilde \Gamma $ is compact.\index{cocompact}
Let us denote by $\phi : \Gamma/\mathrm{Tor}(\Gamma) \to \tilde \Gamma $ as isomorphism.

We equip the Lie group $G$ with the following metric: 
Let $\tilde S\subset G$ the finite set $\tilde S: = \phi (\pi(S))$. Then, 
the {\em Stoll distance} between $g_1, g_2 \in G$ ({\em relative to} $\tilde S$) is\index{Stoll metric}
$$d_{\rm Stoll}(g_1,g_2): =\inf\left\{|t_1|+\ldots+|t_n|\;:\; g_1^{-1}g_2=s_1^{t_1}\cdot \ldots \cdot s_n^{t_n}, n \in \N, s_1, \ldots,s_n \in \tilde S,
t_1, \ldots,t_n \in \R
\right\}$$
This distance function on $G$ is geodesic and left-invariant. Hence, it is a sub-Finsler metric by Berestovski's Theorem~\ref{Berestovskii}; see Exercise~\ref{ex:stollmetric}.
The Stoll distance is named after M. Stoll since he proved that if $G$ is step-two nilpotent, then on $\tilde \Gamma$, the difference between $d_{\rm Stoll}$ and the word metric $d_{\tilde S}$ is bounded; see \cite{stoll}. 
It is an open problem whether the same holds true in groups of steps of at least 3.
However, on $\tilde \Gamma$ the distances $d_{\rm Stoll}$ and $d_{\tilde S}$ are asymptotic, and also quantitatively, as shown in \cite{Breuillard-LeDonne1, GianellaThesis}:
$$ | d_{\tilde S}(p_1, p_2) -d_{\rm Stoll} (p_1, p_2) | \leq C \left(d_{\tilde S}(p_1, p_2)\right)^{1-\alpha}, \qquad \forall p_1, p_2\in \tilde \Gamma,$$
for some $C=C(\tilde S)$ and $\alpha=\alpha(G)>0$. 

After all the above steps, we obtained a quasi-isometry
$$\varphi : (\Gamma, d_S)\longrightarrow (G, d_{\rm Stoll}),$$
which is asymptotic: 
\begin{equation}\label{asympto_Stoll_S}
 | d_{ S}(p_1, p_2) -d_{\rm Stoll} (\varphi(p_1), \varphi(p_2)) | \leq C \left(d_{S}(p_1, p_2)\right)^{1-\alpha}, \qquad \forall p_1, p_2\in \Gamma, \end{equation}
for some $C=C(S, \phi)$ and $\alpha=\alpha(G)>0$.

Let $F$ be a compact fundamental domain for the action of $\Gamma$ on $G$.
Let $\vol$ be the Haar measure on $G$ such that 
$$\vol( F) = \vol (G/H) =1.$$
With this choice of normalization for the Haar measure, we have 
$$\vol(S^n F)={\rm Card}(S^n), \qquad \forall n\in \N.$$

Let $d_{\infty}$ be the Pansu limit metric associated with $d_{\rm Stoll} $, 
as in Section~\ref{sec:Pansu asymptotic structure}.
For the next observation, we use that $F$ is bounded, the asymptotic property \eqref{asympto_Stoll_S}, and Theorem~\ref{Pansu Asymptotic Theorem bis} applied to $d_{\rm Stoll} $.
There are constants $C, \alpha> 0$ such that 
$$ 	\begin{cases} S^n F \subseteq B_{d_{\infty}}(1,n + Cn^{1-\alpha } ),& 
\\ B_{d_{\infty}}(1,n)\subseteq S^{n+ Cn^{1-\alpha }} ) F,&
\end{cases}
\qquad \forall n \leq 1.$$
Recall the scaling property of $d_{\infty}$ from Proposition~\ref{prop_vol_Carnot}.
Consequently, we conclude the volume expansion:
$$ {\rm Card}(S^n) = \vol(B_{d_{\infty}}(1,1)) n^Q + O(n ^{Q-\alpha}), \qquad \text {as } n\to \infty. $$
We actually obtained a more quantitative expansion than \eqref{Pansu_expansion}.

\begin{example}[Standard generating set in the Heisenberg group]\label{example:standard generating set in Heisenberg}\index{Heisenberg! -- group} 
The Heisenberg group in the standard exponential coordinates as in Section~\ref{sec:Heisenberg} admits many discrete cocompact subgroups. In these coordinates, a standard one is $\Gamma: =\Z\times \Z\times \frac12\Z$, which
 is a nonabelian subgroup with respect to the product \eqref{Hprod}.
 The subset $S: =\{ (\pm 1,0,0), (0, \pm 1, 0)\}\subseteq \Gamma$ is a finite generating subset of $\Gamma$. Hence, it induces a word metric $d_S$ and $(\Gamma, d_S)$ has a polynomial volume growth.
 From what we have seen just above, the asymptotic geometry of $(\Gamma, d_S)$ is equal to the sub-Finsler geometry of the Heisenberg Lie group equipped with the $\ell_1$-Carnot metric. Such a metric has the standard horizontal bundle with the norm given by the convex hull of $\log(S )$. For a visual representation of the $\ell_1$-Carnot metric ball, we refer to Figure~\ref{fig:Heis_ell1}.
 
 \begin{figure}[h]
\centering
\subfigure[Unit ball for the Pansu limit] 
{
 \hspace{-0.2cm} \includegraphics[height=4.7cm]{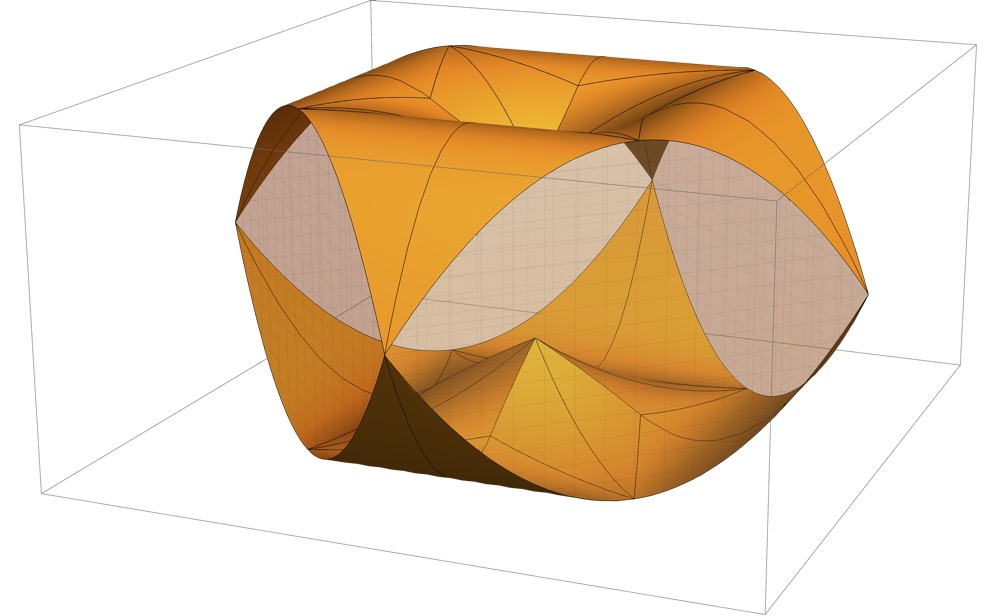}}
 \hspace{-0.5cm}
\subfigure[$S^7$ rescaled by $\delta^{-1}_7$] 
{
 \includegraphics[height=4.7cm]{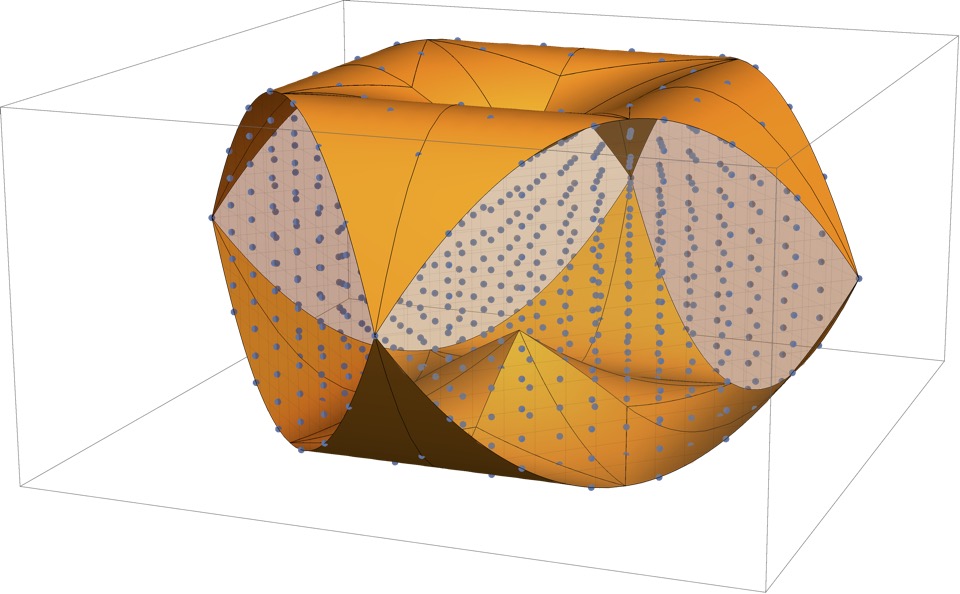}} 
 \caption{Balls for the case of the standard generating set in the Heisenberg group: $S: =\{ (\pm 1,0,0), (0, \pm 1, 0)\}$. The large balls with respect to the word distance, when rescaled by the intrinsic dilations, converge to the Carnot-Carathèodory ball.}
 \label{fig:Heis_ell1}
\end{figure} 
\end{example}

\subsection{A comment on Gromov's Theorem}
The full proof of Theorem \ref{Gromov_Theorem} requires some algebra that goes beyond the tasks of this book. However, we can spend some words explaining why there is a Carnot group among the possible blowdown spaces of $(\Gamma, d_S)$. 

Because $\Gamma$ is assumed to have polynomial growth, then there is $K\in \R_+$ such that the metric space $(\Gamma, d_S)$ is doubling with constant $K$. Consequently, for each $\eps>0$ the metric space $(\Gamma, \eps d_S)$ is $K$-doubling, as well.

Let $\eps_j\searrow 0$ be an infinitesimal sequence. We consider the Gromov-Hausdorff convergence of boundedly compact pointed metric spaces.
By a Gromov's argument, analogs to Ascoli-Arzel\'a results, up to passing to a subsequence, we may assume that the sequence of uniformly doubling pointed metric spaces $(\Gamma, \eps_j d_S,1)$ converges to some pointed complete metric space $(X, d_X, x_0)$, which still is $K$-doubling. In particular, the space $X$ is locally compact and has finite topological dimension.

Here are some more properties of $(X, d_X)$. 
As $(\Gamma, d_S)$ is $(1,1)$-quasi-geodesic (recall Definition~\ref{quasi-geodesic}), then, the metric space $(\Gamma, \eps d_S)$ is $(1, \eps)$-quasi-geodesic. Hence, the limit $(X, d_X)$ is geodesic. In particular, the space $X$ is connected and locally connected.

Since each $(\Gamma, \eps d_S)$ is isometrically homogeneous, then so is $(X, d_X)$. 
We checked all the conditions to apply Theorem~\ref{Montgomery-Zippin}, which was a consequence of the theory of Gleason-Montgomery-Yamabe-Zippin.
Consequently, the metric space $(X, d_X)$ has the structure of a Lie coset space. 
By Berestovskii's Theorem~\ref{Berestovskii}, the metric space $(X, d_X)$ is a sub-Finsler manifold. Therefore, by Mitchell Theorem~\ref{mitchellTheorem}, every tangent metric space of $(X, d_X)$ is a sub-Finsler Carnot group $G$.

As a general fact, the blowup of a blowdown is a blowdown. In fact, consider $\lambda_j\to \infty$ diverging so slow that 
$\eps_j: =\lambda_j\eps_j\searrow 0$.
Then, the sequence $(\Gamma, \tilde\eps_j d_S,1)$ converges to a tangent metric space of $(X, d_X,x_0)$, hence, to $G$. Thus, we proved that one of the blowups of $(\Gamma, d_S)$ is the sub-Finsler Carnot group $G$.

In Gromov's proof for Theorem~\ref{Gromov_Theorem}, one considers various possible actions of $\Gamma$ on the Carnot group $G$. We refer to \cite[Section~8]{Gromov-polygrowth} for the full argument.

\section{Exercises}
\begin{exercise}\label{ex8051255}
 Let $d_1, d_2$ be two left-invariant boundedly compact distances on a Lie group $G$ inducing the manifold topology.
 Then the increasing function $\xi:(0, \infty)\to(0, \infty)$ defined by
 \[
 \xi(r) = \diam_{d_1} \left(\overline B_{d_2}(1_G,r)\right)
 \]
 is such that $\xi(r)\to 0$, as $r\to 0$, and
 $d_1(p,q) \le \xi(d_2(p,q))$.
\end{exercise}

\begin{exercise}\label{ex66711a49}
 Let $X^{u, \lambda}(t,p)$ as in~\eqref{def_Xulambda_29May}.
 Let $K, \Lambda', C_1,C_2, \ell$ as in Remark~\ref{rem29May1538}.
 There exists $C>0$ such that, for all $\eps>0$, $u,v\in L^\infty([0,1];V)$ with $\|u-v\|_{\infty}< C_1 \eps$, $\|u\|_{\infty}$, $\|v\|_{\infty}<\ell$, all $a,b\in\Lambda'$ with $|a-b| < C_2 \eps$:
 \begin{description}
 \item[\eqref{ex66711a49}.i.] 
 $
 \|X^{u,a}(t,p) - X^{v,b}(t,p')\| \le C\|u\|_{\infty} + C\epsilon + C \|u\|_{\infty} \|p-p'\|,$ for all $p,p'\in K$;
 \item[\eqref{ex66711a49}.ii.] 
 if $\alpha$ and $\beta:[0,1]\to K$ are
 absolutely continuous integral curves of $X^{u,a}$ and $Y^{v,b}$, respectively, with $\alpha(0)=\beta(0)$, then
 $
 \|\alpha(1)-\beta(1)\| \leq C \|u\|_{\infty} \epsilon .
 $ \end{description}
\end{exercise}

\begin{exercise}
Let $X, Y, Z: =[X, Y]$ be the standard basis for the Lie algebra of the Heisenberg group $H$.
	For $n\in\N$, let $d_n$ be the Riemannian distance for which $X, Y, \frac1n Z$ form an orthonormal frame.
	The metric spaces $(H, \frac1n d_1)$ and $(H, d_n)$ are isometric.
\end{exercise}

\begin{exercise}
 	Let $M$ be a manifold and $\Delta\subset TM$ a bracket-generating subbundle.
	Let $(g_n)_{n\in\N}$ be a sequence of Riemannian metrics on $M$.
	Assume that 
	the orthogonal to $\Delta$ is the same for each $g_n$, that 
	$
	g_n|_\Delta = g_1|_\Delta$, 
	for all $ n\in\N$,
	and, for all $X\notin\Delta$, 
	$
	g_n(X, X) \to +\infty$,
	 as $ n\to\infty $.
	For all $p,q\in M$, we have
	$
	\lim_{n\to\infty} d_{g_n}(p,q) = \dcc(p,q)
	$
	where $\dcc$ is the sub-Riemannian distance associated with $(\Delta,g_1|_\Delta)$.\\
	{\it Hint.} See Example~\ref{teo:limit_dist_gruppi}.
\end{exercise}

\begin{exercise}
If $(M, d)$ is a CC space and $\lambda>0$, then $(M, \lambda d)$ is a CC space.
\end{exercise}

 \begin{exercise}
 Every sub-Riemannian Carnot group is the limit of some (left-invariant) Riemannian metrics on the same Lie group.
 \end{exercise}
\begin{exercise}
 	Let $G$ be an $s$-step stratified group, with a basis $X_1, \dots, X_n$ adapted to the stratification.
	For each $\lambda>0$, consider the dilations $\delta_\lambda$ and the Riemannian distance $d_\lambda$ for which 
		$
	X_1, \dots, \frac{\lambda^{{\rm deg}(X_j)}}\lambda X_j, \dots, \lambda^{s-1} X_n
	$
	are orthonormal. 
Then, the distance $\lambda d_1$ is associated with the Riemannian metric that makes $	\frac1\lambda X_1, \dots, \frac1\lambda X_n $	orthonormal.
	Every map $\delta_\lambda$ sends each vector field $\frac1\lambda X_j$ to the vectors field
	$
	\frac1\lambda \lambda^{{\rm deg}(X_j)} X_j .
	$
	For all $\lambda>0$ the metric space $(G, \lambda d_1)$ is isometric to $(G, d_\lambda)$ via the map $\delta_\lambda$. 	
 \end{exercise}

\begin{exercise}\label{ex3jun20242200}
Let $F: M_1\to M_2$ be a diffeomorphism between connected manifolds. Let $(\Delta, \norm{\cdot})$ be a sub-Finsler structure on $M_1$. Then the only metric on $M_2$ that makes $F$ an isometry comes from the sub-Finsler structure where $\norm{v}: = \norm {\dd F^{-1} v}$ for all $v\in \Delta_p: = (\dd F) \Delta_{F^{-1} (p)}$ and all $p\in M_2$. 
\end{exercise}

\begin{exercise}
In each $s$-step nilpotent simply connected Lie group, the Dynkin product for the BCH formula \eqref{Dynkin Formula}, has the form:\index{Dynkin product} 
\begin{align} \label{BCH}
x_1\star x_2= \log(\exp(x_1)\cdot \exp(x_2)) = x_1 + x_2 + \sum_{k=2}^s \sum_{q \in \{1,2\}^k}b_{k,q} \lbrack x_{q_1}, \ldots,x_{q_k} \rbrack,
\end{align}
where $b_{k,q}$ are universal real constants, as in Definition~\ref{Dynkin_product}, and here $\lbrack x_{q_1}, \ldots,x_{q_k} \rbrack$ denotes the left-iterated bracket: $\lbrack x_{q_1}, \ldots,x_{q_k} \rbrack: =\lbrack x_{q_1}, \lbrack\ldots, \lbrack x_{q_{k-1}},x_{q_k} \rbrack\ldots\rbrack\rbrack.$
\end{exercise}

\begin{exercise}\label{ex18june2311}
Suppose $G$ is a nilpotent simply connected Lie group.
Choose a basis $(e_1, \ldots,e_n)$ of $\g$ that is adapted to a compatible linear grading $\g=\oplus_i V_i$. 
Given an index $i \in \{1, \ldots, n\}$, let $d_i$ be the degree of $e_i$, namely the integer such that $e_i \in V_{d_i}$. For every multi-index $\alpha=(\alpha_1, \ldots, \alpha_n) \in \N^n$, we set $x^\alpha: =x_1^{\alpha_1}\cdot \ldots \cdot x_n^{\alpha_n}$ and $d_\alpha: =d_1 \alpha_1 + \ldots + d_n \alpha_n$. 
Let $\star$ be the Dynkin product on $\g$.
Then, for some constants $C_{\alpha, \beta} \in \R$, we have
$$(x\star y)_i=x_i + y_i + \sum_{\{\alpha, \beta \textnormal{ }|\textnormal{ } d_\alpha \geq 1, d_\beta \geq 1, d_\alpha+d_\beta \leq d_i\}} C_{\alpha, \beta} x^\alpha y^\beta.$$
For the Dynkin product $\star_0$ associated with the asymptotic Lie bracket $\llbracket \cdot, \cdot\rrbracket_\infty$ as in \eqref{eq: limit bracket asymp},
we have
$$(x\star_0 y)_i=x_i + y_i + \sum_{\{\alpha, \beta \textnormal{ }|\textnormal{ } d_\alpha \geq 1, d_\beta \geq 1, d_\alpha+d_\beta = d_i\}} C_{\alpha, \beta} x^\alpha y^\beta,$$
where we have chopped the terms with $d_\alpha+d_\beta<d_i$.
\end{exercise}

\begin{exercise} 
In the construction of the family of products on a nilpotent Lie algebra $\g$ as in the proof of Pansu's Theorem, 
for $\eps>0$, we denote by $\star_\eps$ the Dynkin product associated with the contracted brackets $[\cdot, \cdot]_\eps$ as in \eqref{def_bracket_Pansu}. We have
$$ g_1 \star_\eps g_2= \delta_\eps ( \delta_\eps^{-1} (g_1) \star \delta_\eps^{-1} (g_2)), \qquad \forall g_1, g_2\in \g.$$
\end{exercise}

\begin{exercise}
Let $G$ be a nilpotent simply connected subFinsler Lie group. Let $d_\infty$ its Pansu limit metric and $\rho_\eps$ the approximating distances, as in \eqref{def_rho_23May}.
Let $C>0$ be the Guivarc'h constant from \eqref{eq_thm_guiv2}.
Then, for all $R > 0$ and all $\epsilon\in (0,1)$, we have
$
B_{\rho_{\epsilon}}(1,R) \subseteq B_{d_\infty}(1,C(R+1)).$
\\ {\it Solution.} For $p \in B_{\rho_{\epsilon}}(0,R)$, one has $ d_\infty(1,p) = \epsilon d_\infty(1, \delta_\epsilon^{-1}(p))
 \leq C\epsilon d(1, \delta_\epsilon^{-1}(p)) + C\epsilon
 \leq CR+C\epsilon.$
 \end{exercise}

\begin{exercise}
Theorem~\ref{Pansu Asymptotic Theorem tris} implies
Theorem~\ref{Pansu Asymptotic Theorem bis}, 
and, in particular, 
\[\Big| d(1,g)-d_\infty (1,g) \Big| = O\left( d(1,g) ^{1- 1/s}\right), \qquad \text{ as } g\rightarrow \infty.\] 
{\it Hint.} As $d_\infty(1,g) \xrightarrow{}\infty$, take $\epsilon : = (d_\infty(1,g))^{-1}$ and $p : = \delta_\epsilon(g).$
 \end{exercise}
 
 \begin{exercise}[Pansu comparison theorem]
 Let $(G, d)$ be a nilpotent simply connected subFinsler Lie group.
 Let $d_\infty$ be the Pansu limit metric on $G$.
 Then
 $\left|\frac{d(1,g)}{d_\infty(1,g)} \right| \xrightarrow{} 1, $ as $ g \xrightarrow{}\infty$.
 \\{\it Hint.} Use Theorem~\ref{Pansu Asymptotic Theorem bis}
\end{exercise}

\begin{exercise}
 Using Exercise~\ref{ex66711a49} instead of Remark~\ref{rem29May1538}, upgrade equation \eqref{23May1034} to
 \[
 - C \epsilon^{1/s} \rho_0(p,q)^{1/s}
 \leq \rho_0(p,q) - \rho_\eps(p,q) \leq
 C \epsilon^{1/s} \rho_\epsilon(p,q)^{1/s} .
 \]
\end{exercise}


\begin{exercise}\label{ex_stratification_Mitchell}
Given a polarized Lie group $(G,V_1)$, with the notation from \eqref{stratification_Mitchell}, we have
$[V_1^{( i)}, V_1^{( j)}]\subseteq V_1^{( i+j)}$. Hence, the Lie bracket $[\cdot, \cdot]_0$, as in \eqref{Lie_Mitchell}, is well defined on $\g_0$ and \eqref{stratification_Mitchell} gives a stratification for $\g_0$.
\end{exercise}

\begin{exercise}\label{ex18june2311Mitchel}
Let $(G, \Delta)$ be a polarized Lie group.
Choose a basis $(e_1, \ldots,e_n)$ of $\g$ that is adapted to an 
adapted linear grading $\g=\oplus_i V_i$ as in \eqref{eq6655cc5c}.
Given an index $i \in \{1, \ldots, n\}$, let $d_i$ be the degree of $e_i$, namely the integer such that $e_i \in V_{d_i}$. For every multi-index $\alpha=(\alpha_1, \ldots, \alpha_n) \in \N^n$, we set $x^\alpha: =x_1^{\alpha_1}\cdot \ldots \cdot x_n^{\alpha_n}$ and $d_\alpha: =d_1 \alpha_1 + \ldots + d_n \alpha_n$. 
Let $\star$ be the Dynkin product on a neighborhood of 0 in $\g$.
Then, for some constants $C_{\alpha, \beta} \in \R$, we have
$$(x\star y)_i=x_i + y_i + \sum_{\{\alpha, \beta \textnormal{ }|\textnormal{ } d_\alpha \geq 1, d_\beta \geq 1, d_\alpha+d_\beta \geq d_i\}} C_{\alpha, \beta} x^\alpha y^\beta.$$
For the Dynkin product $\star_0$ associated with the osculating Carnot algebra Lie bracket $[\cdot, \cdot]_0$ as in \eqref{eq: limit bracket tangent},
we have
$$(x\star_0 y)_i=x_i + y_i + \sum_{\{\alpha, \beta \textnormal{ }|\textnormal{ } d_\alpha \geq 1, d_\beta \geq 1, d_\alpha+d_\beta = d_i\}} C_{\alpha, \beta} x^\alpha y^\beta,$$
where we have chopped the terms with $d_\alpha+d_\beta>d_i$.
\end{exercise}

 
\begin{exercise} 
In the construction of the family of products on a polarized Lie algebra $\g$ as in the proof of Mitchell's Theorem, for $\eps>0$, we denote by $\star_\eps$ the Dynkin product associated with the dilated brackets $[\cdot, \cdot]_\eps$ as in \eqref{def_bracket_Mitchel}. We have
$$ g_1 \star_\eps g_2= \delta_\eps^{-1} ( \delta_\eps (g_1) \star \delta_\eps (g_2)), \qquad \forall g_1, g_2\in \g.$$
\end{exercise}

 \begin{exercise}
 Theorem~\ref{mitchellTheorem tris}, implies Theorem~\ref{mitchellTheorem bis}, 
and, in particular, 
\[\Big| d(1,g)-d_0 (1,g) \Big| = O\left( d(1,g) ^{1+1/s}\right), \qquad \text{ as } g\rightarrow 1.\] 
{\it Hint.} 
As $d_0(1,g) \xrightarrow{}0$, take $\epsilon : = d_0(1,g)$ and $p : = \delta_\epsilon^{-1}(g).$
 \end{exercise}

\begin{exercise}\label{ex66b22ab2}
 Let $(G, \star_\eps)_{\eps\in[-1,1]}$ be a sequence of Lie group structures on a set $G$ as in Section~\ref{sec Analytically varying Lie structures}.
 Let $X^{1, \eps}, \dots, X^{k, \eps}$ be left-invariant vector fields on $(G, \star_\eps)$
 that depend smoothly on $\eps\in[-1,1]$.
 Assume that $X^{1, \eps}, \dots, X^{k, \eps}$ define a (bracket-generating) sub-Riemannian structure on $G$ of step $s$.
 Let $d_\eps$ be the sub-Riemannian metric.
 Fix a Riemannian metric $\rho$ on $G$.
 Then there exist a neighborhood $U$ of $1$ in $G$ and a constant $C>1$ such that
 \begin{equation*}
 \frac{1}{C} \rho \leq d_\eps \leq C\rho^{1/s}
 \qquad\text{on $U$, }\forall \eps\in[-1,1] .
 \end{equation*}
 {\it Hint:}
 Define a map $E_\eps$ as in Proposition~\ref{prop66b0c659}, which depends smoothly on $\eps$.
 Consider the map $\R^n\times[-1,1]\to G\times[-1,1]$, $(\mathbf t, \eps)\mapsto (E_\eps(\mathbf t), \eps)$.
 Proceed as in Corollary~\ref{corol_weak_BB}.
\end{exercise}

\begin{exercise}
 The tangent cone of each Carnot group $G$ is $G$ itself. In fact, dilations $\delta_\lambda$ provide isometries between $(G, \dcc )$ and $(G, \lambda \dcc )$.
\end{exercise}

\begin{exercise}
 Using Exercise~\ref{ex66711a49} instead of Remark~\ref{rem29May1538}, upgrade equation \eqref{eq6655cdd8} to
 \[
 - C \epsilon^{1/s} d_0(p,q)^{1/s}
 \leq d_0(p,q) - d_\eps(p,q) \leq
 C \epsilon^{1/s} d_\epsilon(p,q)^{1/s} .
 \]
\end{exercise}


\begin{exercise}
If $M$ is a Riemannian manifold, then 
$M$ is polarized by $TM$ and, for every $p\in M$, the 
osculating Carnot group of $M$ at $p$ is commutative.
\end{exercise}

\begin{exercise}
If $M$ is a contact $3$-manifold, then for every $p\in M$ the 
osculating Carnot group of $M$ at $p$ is the Heisenberg group.
\end{exercise}

\begin{exercise}
If $G$ is a Carnot group, then for every $p\in G$, the osculating Carnot group of $G$ is $G$ itself.
\end{exercise}

\begin{exercise}
Give an example of Lie group $G$ with a left-invariant bracket-generating distribution such that 
the osculating Carnot group of $G$ is not isomorphic to $G$.
\end{exercise}

\begin{exercise}\label{ex:stollmetric} The Stoll metric $d_{\rm Stoll}$ relative to a finite generating set $S\subset G$ coincides with the left-invariant sub-Finsler metric induced by the norm whose unit ball is the convex hull of $S$ in the Lie algebra $\Lie(G)$ and with left-invariant distribution induced by the subspace spanned by $S$ in $\Lie(G)$. 
\\{\it Hint:} Prove it directly or invoke Berestovski's theorem~\ref{Berestovskii}.
\end{exercise}

\chapter{Rank-one symmetric spaces}\label{ch_ROSS}

In Chapter \ref{ch_Heintze}, we will show that to every Riemannian symmetric space, one can associate a `visual boundary' that has the structure of a Carnot group. 
In this chapter, we review the classical notion of symmetric space. We will see each rank-one symmetric space of noncompact type has the structure of a metric Lie group. 

A \textit{Riemannian symmetric space} is a connected Riemannian manifold $M$ where for each point $p\in M$ there exists an isometry $\sigma_p$ of $M$ such that $\sigma_p(p)=p$ and the differential of $\sigma_p$ at $p$ is the multiplication by $-1$. 
Simple examples of symmetric spaces are round spheres, Euclidean spaces, and real-hyperbolic spaces. 
The \textit{rank} of a symmetric space is the largest dimension of a flat subspace in $M$, where a \textit{flat subspace of dimension $n$} in $M$ is the image of a local isometry from Euclidean $\mathbb{R}^n$ to $M$. 
For example, spheres and hyperbolic spaces have rank 1, whereas Euclidean $n$-space has rank $n$. 
A symmetric space is of \textit{noncompact type} if it is not the product of two symmetric spaces one of which is either compact or Euclidean. Hence, we shall exclude the spheres and the Euclidean spaces.
Symmetric spaces were first introduced and studied by \'Elie Cartan in 1926; see \cite{car1, car2}. 
In particular, he gave a complete description of these spaces by means of the classification of simple Lie algebras.

In this chapter, taking Cartan's description for granted, we first prove that every rank-one symmetric space of noncompact type admits a Lie group structure of a semidirect product with a precise formula for the left-invariant distance. 
The fact that such spaces admit semidirect-product structures has been known at least since Ernst Heintze's work in the 1970s; see \cite{Ernst}. 
However, the formula for the left-invariant distances cannot be traced in the literature. 
To study these spaces, we will need the following result: Let $M$ be a rank-one symmetric space of noncompact type, then $M$ is one of the following spaces, which we call $\mathbb{A}$\textit{-hyperbolic spaces} $\Hyp{A}^n$, with $n\in\mathbb{N}$: real-hyperbolic $n$-space $\Hyp{R}^n$, complex-hyperbolic $n$-space $\Hyp{C}^n$, quaternionic-hyperbolic $n$-space $\Hyp{H}^n$, or the octonionic plane $\Hyp{O}^2$. 
The proof of such last fact was indicated by Cartan but completely established in this form in the 1950s; see Arthur Besse's 1978 book \cite[Section~3.G]{besse} and see Heintze's 1974 paper \cite[Section 5]{Ernst} for a geometric proof. 

We shall introduce $\mathbb{A}$-hyperbolic spaces as metric spaces, where $\mathbb{A}$ is the set of the real, complex, or quaternionic numbers, giving a unified treatment of the subject. We will comment on the octonionic case in Remark~\ref{no_octo}.
For simplicity of notation, we shall restrict to the two-dimensional case, but all constructions generalize to higher dimensions without effort; see the exercise section. 
Following Felix Klein's construction, we shall describe the $\mathbb{A}$\textit{-hyperbolic space $\Hyp{A}^2$} as an open subset of the projectivization of the space $\mathbb{A}^{2,1}$ equipped with a Hermitian form of type $(2,1)$, as of Definition~\ref{def_A21}. 
We shall recall the distance function on $\Hyp{A}^2$, referring to Martin Bridson and Andr\'e H\"afliger's 1999 book \cite[Part II, Chapter 10]{Bridson_Haefliger}.

To recall the Lie group structure on each $\mathbb{A}\mathbf{H}^2$, in Section~\ref{sec: groups Heis on A}, we revise the continuous Heisenberg group $\mathcal{N}$ modeled on $\mathbb{A}$ and its intrinsic dilations; we follow \cite[Chapter~XII Section~1]{Stein:book}. 
We will double check that, after the identification of $\Hyp{\mathbb{A}}^2$ with $\mathcal{N}\rtimes\mathbb{R}_{+}$, the hyperbolic distance is invariant under left translations on $\mathcal{N}\rtimes\mathbb{R}_{+}$. 
We shall write explicitly the distance on the $\mathbb{A}$-hyperbolic space modeled as $\mathcal{N}\rtimes\mathbb{R}_{+}$ in terms of elementary functions of the coordinates. This whole chapter is devoted to the proof of the following summarizing result:
\begin{theorem}\label{thm_group_structure_ROSS2}
For every $\mathbb{A}\in\{\mathbb{R, \mathbb{C}, \mathbb{H}}\}$ the $\mathbb{A}$-hyperbolic space $\Hyp{\mathbb{A}}^2$ is isometric to the manifold $\mathbb{A}\times \imaginary(\mathbb{A})\times\mathbb{R}_{+}$ equipped with the group product
and the left-invariant distance $d$ given for all $(\xi, \mathbf{v};\lambda)$ and $(\xi', \mathbf{v}';\lambda')\in \mathbb{A}\times \imaginary(\mathbb{A})\times\mathbb{R}_{+}$ 
 by $$(\xi, \mathbf{v};\lambda)\cdot(\xi', \mathbf{v}';\lambda')=\left(\xi+\lambda\xi', \mathbf{v}+\lambda^{2}\mathbf{v}'+2\imaginary(\bar{\xi}\lambda \xi');\lambda\lambda'\right),
$$ and
$$ d(
(\xi, \mathbf{v}, \lambda),(\xi', \mathbf{v}', \lambda')
) = 
 \arccosh \sqrt { \left(1+ \frac { 
 |\xi-\xi'|^{2} + |\lambda-\lambda'|^2}{2\lambda\lambda'} \right)^2+ \left( \frac { |{ \mathbf{v} + \bar{\xi}\xi' - \mathbf{v'} -\bar{\xi}'\xi } |} {2\lambda\lambda'}\right)^2} 
.$$
\end{theorem}


The proof of Theorem~\ref{thm_group_structure_ROSS2} is spread in Theorem~\ref{theorem distance in horospherical coordinates} and Theorem~\ref{thm left-invariance hyp distance}, and their preliminaries.  See also the high-dimensional generalization in Exercise~\ref{final_ex}.

\begin{remark}\label{no_octo}
There is a remaining case: the octonionic-hyperbolic plane. It cannot be treated as described above due to the non-associativity of the octonions and, therefore, the impossibility of defining a notion of vector space over the octonions. 
For the basic ideas on how to deal with this case and build the octonionic-hyperbolic plane, see \cite{Pilastro_master_thesis}.  
\end{remark}

\section{Preliminary notions for rank-one symmetric spaces}

\subsection{Euclidean Hurwitz algebras}

The algebraic structure needed to construct rank-one symmetric spaces is exactly given by the Euclidean Hurwitz algebras, also referred to as the {\em normed division algebras over $\R$}.\index{Eucidean Hurwitz algebra}\index{normed! -- division algebra}\index{division algebra! normed --}
Even if they seem very general structures, there are only 4 examples of them; see Theorem~\ref{thm_Hurwitz}. Still, it is convenient to see them from a general viewpoint. 
\begin{definition}\label{Eucidean Hurwitz algebra}
A \emph{Eucidean Hurwitz algebra} \(\mathbb{A}\) is a not necessarily associative division algebra over \(\mathbb{R}\) with an identity \(1\), together with a multiplicative norm \(|\cdot|\). 
In order to clarify, a {\em(not necessarily associative) algebra} \(\mathbb{A}\) over $\R$ is a vector space over $\R$ that is equipped with a bilinear binary multiplication operation $\mathbb{A} \times \mathbb{A} \to \mathbb{A}$.
Whereas, an algebra $\mathbb{A}$ is said to be a {\em division algebra} if\index{division algebra}
\begin{equation}\label{division} \xi, \zeta \in \mathbb{A}, \, (\xi \zeta = 0) \implies (\xi = 0 \text{ or } \zeta = 0), \end{equation}
while a norm \(|\cdot|\) on an algebra $\mathbb{A}$ is said to be {\em multiplicative} if\index{multiplicative! -- norm}\index{norm! multiplicative --}
\begin{equation}\label{mult_norm} | \xi\zeta| = |\xi | \,|\zeta|, \qquad \forall \xi, \zeta \in \mathbb{A}.\end{equation}
\end{definition}

Every Euclidean Hurwitz algebra admits an anti-involution, called \emph{conjugation}. There is a general formula to define the conjugate in terms of the norm; see Exercise~\ref{exercise_norm_imaginary}.
In this sense, Hurwitz algebras may be thought of as a natural generalisation of the real numbers and the complex numbers.
We denote the \emph{conjugate of \(\xi\)} by \(\bar{\xi}\).

It is a result of Hurwitz that the only Euclidean Hurwitz algebras are the following four examples:
\begin{enumerate}
	\item The real numbers \(\mathbb{R}\), with the absolute value as norm.
	\item The complex numbers \(\mathbb{C}\), with Euclidean norm given by the natural identification with $\R^2$.
	\item The quaternionic numbers
\(\mathbb{H}\), with Euclidean norm given by the natural identification with $\R^4$.
	\item The
octonionic numbers \(\mathbb{O}\), with Euclidean norm given by the natural identification with $\R^8$.
\end{enumerate}

\begin{theorem}[Hurwitz \cite{Hurwitz1922}]\label{thm_Hurwitz}
	The only Euclidean Hurwitz algebras are $\mathbb{R}, \mathbb{C}, \mathbb{H}, $ and $\mathbb{O}$.
\end{theorem}

Given an element $\xi$ in a Euclidean Hurwitz algebra we define its {\em real part} as $\Re(\xi) := \frac{\xi+\bar \xi}{2}$
and its {\em imaginary part} as $\imaginary(\xi) := \xi - \Re(\xi) = \frac{\xi-\bar \xi}{2}$.\index{imaginary part, $\imaginary(\xi)$}\index{real part, $\Re(\xi)$}\index{$\imaginary(\xi)$}\index{$\Re(\xi)$}\index{purely imaginary numbers, $\imaginary( \mathbb{A})$}   
Notice that we follow the convention that the imaginary part is an element of the set $\imaginary( \mathbb{A})$ of {\em purely imaginary numbers}.

While we do not need to recall the real or complex numbers, we revise the quaternions and the octonions: 
The \emph{quaternions}, or {\em quaternionic numbers}, are a 4-dimensional associative algebra over $\R$ with basis \(\{ 1,i,j,k \}\), where \(1\) is central and \(i,j,k\) follow the rules: \[
ij=k, \quad jk=i, \quad ki=j,
\] and \[
i^{2} = j^{2} =k^{2} = -1.
\] The quaternions are typically denoted by \(\mathbb{H}\) in honour of its discovery by W. R. Hamilton.
On the quaternions, the \emph{conjugate} of an element 
\(\xi = \xi^{0} + \xi^{1}i + \xi^{2}j + \xi^{3}k\) is
\(\bar{\xi} = \xi^{0} - \xi^{1}i - \xi^{2}j - \xi^{3}k\), 
its real part is $ \Re(\xi):=\xi^{0}\in \R$, while the \emph{imaginary part} of \(\xi\) is 
$\imaginary(\xi):=\xi^{1}i + \xi^{2}j + \xi^{3}k \in \imaginary(\mathbb{H})\simeq \mathbb{R}^{3}$. 
On $\mathbb{H}$, we consider the Euclidean norm via the identification $\mathbb{H}\simeq\R^4$, which makes $\mathbb{H}$ into a Euclidean Hurwitz algebra; see Exercise~\ref{ex:quaternions_Hurwitz}.

The \emph{octonions}, or {\em octonionic numbers}, are a nonassociative Euclidean Hurwitz algebra of dimension 8 over \(\mathbb{R}\). This algebra is defined via the Cayley-Dickson construction \cite{Dickson1919} as pairs of quaternion numbers \((a,b), (c,d)\), together with the multiplication \[
(a,b)\cdot(c,d) := (ac - \bar{d}b,da + b \bar{c}), \qquad \forall a,b, c,d \in\mathbb{H}.
\]

Working with nonassociative algebras takes some time getting used to and we do not want to let it take too much space in this text. Therefore we restrict our discussion from now on to the associative Euclidean Hurwitz algebras. 
But the reader should keep in mind that with a little more work, the story extends to the octonion setting as well; see \cite{Pilastro_master_thesis}.

\subsection{Hermitian forms}

Notice that every associative Euclidean Hurwitz algebra is a ring. 
Let $\mathbb{A}$ be either $\mathbb{R}$, $\mathbb{C}$, or $\mathbb{H}$. Consider the right $\mathbb{A}$-module $\mathbb{A}^{3}$: for $\mathbf{z}=(z_1, z_2, z_3)\in\mathbb{A}^{3}$ and $\lambda \in \mathbb{A}$, we shall consider $\mathbf{z}\lambda:=(z_1\lambda, z_2\lambda, z_3\lambda)$.

\begin{definition}\index{Hermitian! -- form}\index{form! Hermitian --}
	A \emph{Hermitian form} over the right \(\mathbb{A}\)-module
\(\mathbb{A}^{3}\) is a 
map \[
\langle \cdot, \cdot \rangle: \mathbb{A}^{3}\times\mathbb{A}^{3} \to \mathbb{A},
\] such that for all \(\lambda, \mu \in\mathbb{A}\), and
\(\mathbf{z}, \mathbf{w}, \mathbf{z}', \mathbf{w}' \in \mathbb{A}^{3}\), 
\[
\langle \mathbf{z} +\mathbf{z}' , \mathbf{w}+ \mathbf{w}' \rangle = 
 \langle \mathbf{z}, \mathbf{w} \rangle + \langle \mathbf{z}, \mathbf{w}' \rangle + \langle \mathbf{z}', \mathbf{w} \rangle + \langle \mathbf{z}', \mathbf{w}' \rangle ,
\] 
\begin{equation}\label{eq_scalars_Hermitian}
\langle \mathbf{z}\lambda, \mathbf{w}\mu \rangle = \bar{\lambda} \langle \mathbf{z}, \mathbf{w} \rangle \mu,
\end{equation} and 
\begin{equation}\label{eq_sesqui}
\langle \mathbf{z}, \mathbf{w} \rangle = \overline{\langle \mathbf{w}, \mathbf{z} \rangle}.
\end{equation} 
\end{definition}
A Hermitian form is \emph{nondegenerate}, if \[
\left(\, \langle \mathbf{z}, \mathbf{y} \rangle = 0 , \, \forall \mathbf{y} \in \mathbb{A}^{3} \, \right)
\Longrightarrow \mathbf{z} = 0. \] 
The \emph{orthogonal complement} of \(\mathbf{z} \in \mathbb{A}^{3}\) is the set \[
\mathbf{z}^{\perp} := \{ \mathbf{y} \in \mathbb{A}^{3} | \langle \mathbf{z}, \mathbf{y} \rangle = 0 \}.
\] 
If $\langle \cdot, \cdot \rangle$ is nondegenerate, then
\(\mathbf{v} \oplus \mathbf{v}^{\perp} = \mathbb{A}^{3}\).
Every Hermitian form \(\langle \cdot, \cdot \rangle\) over the right
\(\mathbb{A}\)-module \(\mathbb{A}^{3}\) is given by a \(3\times 3\)
matrix \[
H = \begin{pmatrix}
a & b & c \\
d & e & f \\
g & h & i \\
\end{pmatrix}, \quad a,b,c,e,d,f,g,h,i \in \mathbb{A},
\] through
\(\langle \mathbf{z}, \mathbf{w} \rangle = \mathbf{z}^{*}H\mathbf{w}\).
Here \(\mathbf{z}^{*} = (\bar{z}_{1}, \bar{z}_{2}, \bar{z}_{3})\) is the
Hermitian conjugate of \(\mathbf{z}\).
Be aware that to use the matrix products, we should see $\mathbf{w} $ as a column vector.
In order to give a Hermitian form, such a matrix \(H\) satisfies \[
\begin{pmatrix}
a & b & c \\
d & e & f \\
g & h & i \\
\end{pmatrix} = 
\begin{pmatrix}
\bar{a} & \bar{d} & \bar{g} \\
\bar{b} & \bar{e} & \bar{h} \\
\bar{c} & \bar{f} & \bar{i} \\
\end{pmatrix}.
\] 
Equivalently said, a transformation is a Hermitian form if the matrix that represents it equals its Hermitian transpose, where given a matrix $M=(m_{ij})$ over a Euclidean Hurwitz algebra, the \textit{Hermitian transpose} of $M$ is $M^*:=(\conj{m}_{ji})$.\index{Hermitian! -- transpose} 

The classical spectral theory over $\mathbb{C}$ extends to associative Hurwitz algebras; see, for example, \cite{Farenick2003}.
In particular, the matrices representing Hermitian forms over $\mathbb{A}$ only have real eigenvalues.

\subsection{Hermitian forms of signature $(2,1)$}

\begin{definition}[The space \(\mathbb{A}^{2,1}\)]\label{def_A21}
A Hermitian form over \(\mathbb{A}^{3}\) has {\em signature} \((2,1)\) if it has exactly two strictly positive eigenvalues and one strictly negative eigenvalue. 
To stress when \(\mathbb{A}^{3}\) is considered equipped with a Hermitian form of signature \((2,1)\), we denote it by \(\mathbb{A}^{2,1}\).\index{\(\mathbb{A}^{2,1}\)}
Up to linear changes of variables, there is only one such a space \(\mathbb{A}^{2,1}\).
\end{definition}
\begin{lemma}\label{lemma_Hermitian_orthogonal}
Let \(\langle \cdot, \cdot \rangle\) be a Hermitian form of signature \((2,1)\) over the right \(\mathbb{A}\)-module \(\mathbb{A}^{3}\). 
Let \(\mathbf{z}\) be a vector in \(\mathbb{A}^{2,1}\) with \(\langle \mathbf{z}, \mathbf{z} \rangle < 0\), then \(\langle \cdot, \cdot \rangle\) is positive definite on \(\mathbf{z}^{\perp}\).
\end{lemma}

\begin{proof}
We begin by stressing that a Hermitian form of signature \((2,1)\) has a trivial kernel.
Assume that there exists a vector
\(\mathbf{y} \in \mathbf{z}^{\perp}\) such that
\(\langle \mathbf{y}, \mathbf{y} \rangle < 0\). 
We claim that \(\mathbf{z}\) and
\(\mathbf{y}\) are linearly independent over $\mathbb{A}$. Indeed, otherwise, noticing that 
$\mathbf{z}$ cannot be zero, we would have $\mathbf{y}= \mathbf{z} \lambda$ for some $\lambda \in \mathbb{A}$, but then 
$$ 0\stackrel{(\mathbf{y} \in \mathbf{z}^{\perp})}{=}
\langle \mathbf{z}, \mathbf{y} \rangle 
=\langle \mathbf{z}, \mathbf{ z \lambda} \rangle 
\stackrel{\eqref{eq_scalars_Hermitian}}{=}
\bar\lambda\langle \mathbf{z}, \mathbf{ z } \rangle, 
$$
and, because we are in a division algebra, \eqref{division}, either $\bar\lambda$ would be zero (but then this contradicts \(\langle \mathbf{y}, \mathbf{y} \rangle < 0\)) or $\langle \mathbf{z}, \mathbf{ z } \rangle=0$ (which contradicts the assumption \(\langle \mathbf{z}, \mathbf{z} \rangle < 0\)).
Consequently, the space
\(\operatorname{span}(\mathbf{z}, \mathbf{y})\) is a two-dimensional
subspace on which \(\langle \cdot, \cdot \rangle\) is negative-definite.
This is impossible since in \(\mathbb{A}^{2,1}\), we only have one negative eigenvalue by assumption.
\end{proof}
\begin{lemma}[Reverse Schwarz Inequality]\label{Reverse_Schwarz}
	Let \(\langle \cdot, \cdot \rangle\) be a Hermitian form
of signature \((2,1)\) over the right \(\mathbb{A}\)-module
\(\mathbb{A}^{3}\). For every \(\mathbf{z}, \mathbf{w} \in \mathbb{A}^{3}\) with
\(\langle \mathbf{z}, \mathbf{z} \rangle < 0\) and
\(\langle \mathbf{w}, \mathbf{w} \rangle< 0\), we have
\[\langle \mathbf{z}, \mathbf{w} \rangle \langle \mathbf{w}, \mathbf{z} \rangle \geq \langle \mathbf{z}, \mathbf{z} \rangle \langle \mathbf{w}, \mathbf{w} \rangle, \]
with equality if and only if \(\mathbf{z}\) and \(\mathbf{w}\) are
linearly dependent over $\mathbb{A}$.
\end{lemma}

\begin{proof}
	Suppose that \(\mathbf{z}\) and \(\mathbf{w}\) are linearly independent.
Since \(\langle \mathbf{w}, \mathbf{w} \rangle < 0\) and \(\langle \cdot, \cdot \rangle\) is positive definite on
\(\mathbf{z}^{\perp}\) by Lemma~\ref{lemma_Hermitian_orthogonal}, we infer that \(\mathbf{w} \notin \mathbf{z}^{\perp}\), 
i.e., \(\langle \mathbf{z}, \mathbf{w} \rangle \neq 0\).
Let
\(\lambda := -\langle \mathbf{z}, \mathbf{z} \rangle\langle \mathbf{z}, \mathbf{w} \rangle^{-1}\).
Since \(\mathbf{z}\) and \(\mathbf{w}\) are linearly independent,
\(\mathbf{z} + \mathbf{w}\lambda \neq 0\) and a quick calculation shows
that \(\mathbf{z} + \mathbf{w}\lambda \in \mathbf{z}^{\perp}\). 
Again by Lemma~\ref{lemma_Hermitian_orthogonal}, this implies that
\(\langle \mathbf{z} + \mathbf{w}\lambda, \mathbf{w}\lambda \rangle = \langle \mathbf{z} + \mathbf{w}\lambda, \mathbf{z}+\mathbf{w}\lambda \rangle > 0\).
Expanding the left-hand side gives \[
-\langle \mathbf{z}, \mathbf{z} \rangle + \langle \mathbf{z}, \mathbf{z} \rangle^{2}\langle \mathbf{w}, \mathbf{w} \rangle \langle \mathbf{w}, \mathbf{z} \rangle^{-1}\langle \mathbf{z}, \mathbf{w} \rangle^{-1} > 0.
\] Dividing by \(\langle \mathbf{z}, \mathbf{z} \rangle < 0\) and
rearranging, we obtain \[
\langle \mathbf{z}, \mathbf{w} \rangle \langle \mathbf{w}, \mathbf{z} \rangle > \langle \mathbf{z}, \mathbf{z} \rangle \langle \mathbf{w}, \mathbf{w} \rangle.
\] If \(\mathbf{z}\) and \(\mathbf{w}\) are linearly dependent, then
there exists \(\lambda \in \mathbb{A}\) such that
\(\mathbf{z} = \mathbf{w}\lambda\). Then the equality
\[\langle \mathbf{z}, \mathbf{w} \rangle \langle \mathbf{w}, \mathbf{z} \rangle = \langle \mathbf{z}, \mathbf{z} \rangle \langle \mathbf{w}, \mathbf{w} \rangle\]
 follows straightforwardly from the definition of Hermitian forms.
\end{proof}


\section{The $\mathbb{A}$-hyperbolic space $\Hyp{A}^2$}\label{sec A-hyperbolic space}

\subsection{Definition and properties}\index{hyperbolic! -- space}

Let \(\mathbb{A}\) be either \(\mathbb{R}\), \(\mathbb{C}\), or the quaternions \(\mathbb{H}\). Let \(\mathbb{A}^{2,1}\) be the right-module over \(\mathbb{A}\) of \(\mathbb{A}\)-dimension \(3\), equipped with some Hermitian form \(\langle \cdot, \cdot \rangle\) of signature \((2,1)\).
Define the subset
\(V_- := \{\mathbf{z} \in \mathbb{A}^{2,1} : \langle \mathbf{z}, \mathbf{z} \rangle < 0\}.\)
For every nonzero scalar \(\lambda\), we have that \(\mathbf{z}\lambda\) is in \(V_{-}\) if and only if so is \(\mathbf{z}\). We may therefore speak of negative $\mathbb{A}$-lines in \(\mathbb{A}^{2,1}\). 

The \emph{projective model of hyperbolic space} is the collection of negative lines in \(\mathbb{A}^{2,1}\). Formally, we have a projective map \(\mathbb{P}: \mathbb{A}^{2,1} \to \mathbb{P}(\mathbb{A}^{2,1})\),
where \(\mathbb{P}(\mathbb{A}^{2,1})\) is the projective space over \(\mathbb{A}^{2,1}\),
here $\mathbb{P}(\mathbf{z}) :=\{\mathbf{z} \lambda :\lambda\in \mathbb{A}\setminus\{0\}\} $, for $\mathbf{z}\in \mathbb{A}^{2,1}$.
Thus, we define the {\em hyperbolic space}
\(\Hyp{\mathbb{A}}^{2}\) to be \(\mathbb{P}(V_{-})\). 
In symbols, we write:
\begin{equation}\label{def_hyperbolic_space}
\Hyp{\mathbb{A}}^{2} := \{\mathbb{P}(\mathbf{z}) \in \mathbb{P}(\mathbb{A}^{2,1}) : \langle \mathbf{z}, \mathbf{z} \rangle < 0\}.\end{equation}
We denote points in \(\Hyp{\mathbb{A}}^{2}\) simply by \(z\) or \(w\).
Following \cite[Chapter II.10, p.302]{Bridson_Haefliger} for this projective model of hyperbolic space, the \emph{Bergman metric} on \(\Hyp{\mathbb{A}}^{2}\) is given by the distance function
\(d\) defined by the equation\footnote{We decided to follow the convention of not diving by 2 in the definition of the distance, as instead it is done by Parker \cite{Parker}. 
We instead follow, for example, Bridson-Haefliger in order to have a CAT($-1$) space.}
\begin{equation}\label{Bergman metric}
\left(\cosh\Big( {d(z,w)} \Big)\right)^2 = {\frac{\langle \mathbf{z}, \mathbf{w} \rangle \langle \mathbf{w}, \mathbf{z} \rangle}{\langle \mathbf{z}, \mathbf{z}\rangle \langle \mathbf{w}, \mathbf{w} \rangle }}, \end{equation}
where \(\mathbf{z}\) and \(\mathbf{w}\) are any vectors in
\(\mathbb{A}^{2,1}\) such that \(\mathbb{P}(\mathbf{z})=z\) and
\(\mathbb{P}(\mathbf{w})=w\); here $\cosh$ denotes the hyperbolic cosine \eqref{def sinh cosh}.
Equivalently, using the inverse hyperbolic cosine, we write
\[ {d(z,w)} :=\arccosh\sqrt { \frac{\langle \mathbf{z}, \mathbf{w} \rangle \langle \mathbf{w}, \mathbf{z} \rangle}{\langle \mathbf{z}, \mathbf{z}\rangle \langle \mathbf{w}, \mathbf{w} \rangle }}
= \arccosh\sqrt{ \frac{\langle \mathbf{z}, \mathbf{w} \rangle\overline{ \langle \mathbf{z}, \mathbf{w} \rangle}}{\langle \mathbf{z}, \mathbf{z}\rangle \langle \mathbf{w}, \mathbf{w} \rangle }}
= \arccosh\left(\frac{|\langle \mathbf{z}, \mathbf{w} \rangle |}{\sqrt{\langle \mathbf{z}, \mathbf{z}\rangle \langle \mathbf{w}, \mathbf{w} \rangle }}\right)
.\]
Thanks to the Reverse Schwarz Inequality (see Lemma~\ref{Reverse_Schwarz}), we have \[
\frac{\langle \mathbf{z}, \mathbf{w} \rangle \langle \mathbf{w}, \mathbf{z} \rangle}{\langle \mathbf{z}, \mathbf{z}\rangle \langle \mathbf{w}, \mathbf{w} \rangle } \geq 1, \qquad \forall\mathbf{z}, \mathbf{w}\in \mathbb{A}^{2,1}.
\] Hence, the equation makes sense: recall that $\cosh$ bijectively sends $[0,+\infty)$ onto $[1,+\infty)$.

We shall check that the formula does not depend on the representatives: If \(\mathbf{z}'\) and \(\mathbf{w}'\) are any other vectors in
\(\mathbb{A}^{2,1}\) such that \(\mathbb{P}(\mathbf{z}')=z\) and
\(\mathbb{P}(\mathbf{w}') = w\), then there exist nonzero scalars
\(\lambda\) and \(\mu\) in \(\mathbb{A}\) such that \(\mathbf{z}' = \mathbf{z}\lambda\)
and \(\mathbf{w}' = \mathbf{w}\mu\). Thus 
\begin{eqnarray*}
\frac{ |\langle \mathbf{z}', \mathbf{w}' \rangle |^2}{\langle \mathbf{z}', \mathbf{z}'\rangle \langle \mathbf{w}', \mathbf{w}' \rangle } 
&= &\frac{|\langle \mathbf{z}\lambda, \mathbf{w}\mu \rangle |^2}{\langle \mathbf{z}\lambda, \mathbf{z}\lambda\rangle \langle \mathbf{w}\mu, \mathbf{w}\mu \rangle } \\
&\stackrel{\eqref{eq_scalars_Hermitian}}{=} &\frac{|\bar{\lambda}\langle \mathbf{z}, \mathbf{w} \rangle \mu |^2}{|\lambda|^{2}\langle \mathbf{z}, \mathbf{z}\rangle |\mu|^{2} \langle \mathbf{w}, \mathbf{w} \rangle } \\
&\stackrel{\eqref{mult_norm}}{=} &\frac{ |\lambda|^{2}\cdot |\langle \mathbf{z}, \mathbf{w}\rangle|^2\cdot | \mu|^{2}}{|\lambda|^{2}\cdot | \mu|^{2}\cdot\langle \mathbf{z}, \mathbf{z}\rangle \langle \mathbf{w}, \mathbf{w} \rangle } \\
&=& \frac{|\langle \mathbf{z}, \mathbf{w}\rangle |^{2} }{ \langle \mathbf{z}, \mathbf{z}\rangle \langle \mathbf{w}, \mathbf{w} \rangle } .
\end{eqnarray*}
Here, we made explicit use of the fact that the norm on \(\mathbb{A}\) is multiplicative and that the norms are real numbers (and hence central elements).

\subsection{Models of $\mathbb{A}$-hyperbolic space}\label{sec models hyp spaces}

There are multiple models of real-hyperbolic space
\(\Hyp{\mathbb{R}}^{2}\). The most common ones are: 
\begin{enumerate}
	\item the \emph{Minkowski hyperboloid model} \cite[p.18]{Bridson_Haefliger}, 
	\item the \emph{Poincar\'e ball model} \cite[p.86]{Bridson_Haefliger},
	\item the \emph{Poincar\'e half-space model} \cite[p.90]{Bridson_Haefliger},
	\item the \emph{Klein ball model} \cite[p.83 and p.310]{Bridson_Haefliger},
	\item and the \emph{Siegel parabolic model} \cite[p.310]{Bridson_Haefliger}.
\end{enumerate}
The Poincar\'e ball model and the Poincar\'e half-space model are \emph{conformal}, in the sense that they are angle-preserving.
The Klein model and the Siegel domain model are \emph{projective} in the sense that hyperbolic straight lines and planes are the intersection of Euclidean straight lines and planes with the unit ball in \(\mathbb{R}^{2}\) and the interior of the paraboloid \(x_{1} = \frac{x_{2}^{2}}{2}\), respectively.
The Siegel domain model is not often mentioned in the real-hyperbolic case. 
It may be understood as the projective analog of the conformal Poincar\'e half-space model. 
The hyperboloid model, on the other hand, is particularly common in physics but neither conformal nor projective.

The projective models have straightforward generalizations to hyperbolic spaces modeled over each associative Hurwitz algebra.

Let \(\mathbb{A}\) be either \(\mathbb{R}\), \(\mathbb{C}\), or the quaternions \(\mathbb{H}\).
Let \(\mathbf{z}\) and \(\mathbf{w}\) be vectors in
\(\mathbb{A}^{2,1}\). The \emph{first Hermitian form} on
\(\mathbb{A}^{2,1}\) is defined as
\[\langle \mathbf{z}, \mathbf{w} \rangle_1 := \bar{z}_1 w_1 + \bar{z}_2 w_2 - \bar{z}_3 w_3.\]
It is associated with the Hermitian matrix
\[J_{1} := \begin{bmatrix}1 & 0 & 0 \\ 0 & 1 & 0 \\ 0 & 0 & -1\end{bmatrix}.\]
Every representative \(\mathbf{z}' \in \mathbb{A}^{2,1}\) of a point
\(z \in \Hyp{\mathbb{A}}^{2}\) satisfies
\(\langle \mathbf{z}', \mathbf{z}' \rangle_{1} < 0\), i.e., 
$0>\bar{z}_1 z_1 + \bar{z}_2 z_2 - \bar{z}_3 z_3=
|z_1|^2 + |z_2|^2 - | z_3|^2.
$
Therefore the third component of such a $z$ is nonzero. Hence, there exists exactly one representative \(\mathbf{z}\) of \(z\) of the form \[
\mathbf{z} = \begin{bmatrix}
z_{1} \\
z_{2} \\
1
\end{bmatrix} \simeq (z_1,z_2,1), \qquad \text{ for which } |z_1|^2 + |z_2|^2<1.
\] This gives an injection
\(\iota: \Hyp{\mathbb{A}}^{2} \hookrightarrow \mathbb{A}^{2,1}\)
as a section of
\(\mathbb{P}: V_{-} \rightarrow \Hyp{\mathbb{A}}^{2}\) and identifies the hyperbolic space \(\Hyp{\mathbb{A}}^{2}\) with Klein's model which setwise is the unit ball \(\mathbb{B}\) in \(\mathbb{A}^2\). We call such a \(\mathbf{z}\) the \emph{standard lift} of \(z\).

The \emph{second Hermitian form} on \(\mathbb{A}^{2,1}\) is defined as
\begin{equation}\label{def second Hermitian form}
\langle \mathbf{z}, \mathbf{w} \rangle_2 := \bar{z}_{1} w_3 + \bar{z}_{2} w_2 + \bar{z}_{3} w_1.\end{equation}
It is associated with the Hermitian matrix
\[J_{2} := \begin{bmatrix}0 & 0 & 1 \\ 0 & 1 & 0 \\ 1 & 0 & 0\end{bmatrix}.\]
Every representative \(\mathbf{z}' \in \mathbb{A}^{2,1}\) of a point
\(z \in \Hyp{\mathbb{A}}^{2}\) satisfies
\(\langle \mathbf{z}', \mathbf{z}' \rangle_{2} < 0\), i.e.,
$0> \bar{z}_{1} z_3 + \bar{z}_{2} z_2 + \bar{z}_{3} z_1= 2 \Re (\bar{z}_{1} z_3) +|z_2|^2.$ 
Therefore, in particular, the third component of such a $z$ is nonzero. 
Hence, there exists exactly one representative \(\mathbf{z}\) of \(z\) of the form \[
\mathbf{z} = \begin{bmatrix}
z_{1} \\
z_{2} \\
1
\end{bmatrix}.
\] This gives an injection
\(\iota: \Hyp{\mathbb{A}}^{2} \hookrightarrow \mathbb{A}^{2,1}\)
as a section of
\(\mathbb{P}: V_{-} \rightarrow \Hyp{\mathbb{A}}^{2}\) and
identifies each point $z\in \Hyp{\mathbb{A}}^{2}$ with the
point \((z_1,z_2) \in \mathbb{A}^2\) such that
\(2\Re(z_1) + |z_2|^2 < 0.\) This domain 
\begin{equation}\label{def_Siegel_domain}\mathbb{D}:=\{
(z_1,z_2) \in \mathbb{A}^2:
2\Re(z_1) + |z_2|^2 < 0 \}\end{equation}
is called the \emph{Siegel (parabolic) domain}.

A map that passes from one model of hyperbolic space to another is called a \emph{Cayley transform}. 
For example, a Cayley transform interchanging the first and second Hermitian forms is given by
\[\frac{1}{\sqrt{ 2 }} \begin{bmatrix} 1 & 0 & 1 \\ 0 & \sqrt{2} & 0 \\ 1 & 0 & -1\end{bmatrix}.\]
A Cayley transform is not unique since it may be pre- and post-composed by any unitary matrix preserving the relevant Hermitian form.

In the following we will denote by $\mathbf{z}$ both the standard lift $(z_{1},z_{2},1) \in \mathbb{A}^{2,1}$ of $z \in \Hyp{\mathbb{A}}^{2}$ and its identification with the point $(z_{1},z_{2})$ in the Siegel domain $\mathbb{D}$.

\subsection{A foliation of the Siegel domain}
We consider the hyperbolic space
\(\Hyp{\mathbb{A}}^{2}\), as defined in \eqref{def_hyperbolic_space}, as a subset of the projective space \(\mathbb{P}(\mathbb{A}^{2,1})\). Within this larger subset, the hyperbolic space
\(\Hyp{\mathbb{A}}^{2}\) has a topological boundary denoted by \(\partial\Hyp{\mathbb{A}}^{2}\). 

When we consider \(\mathbb{A}^{2,1}\) with respect to the second Hermitian form \(\langle \cdot, \cdot \rangle_{2}\), from \eqref{def second Hermitian form}, we distinguish two types of points in \(\partial\Hyp{\mathbb{A}}^{2}\): we say that \(z\in \partial\Hyp{\mathbb{A}}^{2}\)
is {\em finite}, if \(z\) admits a representative of the form \(\mathbf{z} = (z_1,z_2,1) \in \mathbb{A}^{2,1}\). Otherwise, we say that $z$ is {\em not finite}.

We point out that there exists exactly one point on the boundary that is not finite, and it is represented by $(1,0,0)$. Instead, the finite points in \(\partial\Hyp{\mathbb{A}}^{2}\) are those that are represented by a point in the boundary of the Siegel domain $\mathbb{D}$ as a subset of $\mathbb{A}^{2}$. 
Namely, a point $z\in \mathbb{P}(\mathbb{A}^{2,1})$ is a finite boundary point for \(\Hyp{\mathbb{A}}^{2}\) if (and only if) it has a representative \(\mathbf{z} = (z_1,z_2,1) \in \mathbb{A}^{2,1}\) such that 
$$z_1 + \bar{z}_{1} + |z_2|^2 = 0.$$
For this reason, we define the {\em parabolic boundary} of the hyperbolic space
\(\Hyp{\mathbb{A}}^{2}\) (modelled by the Siegel domain) as
 $$\partial_{\rm Par}\mathbb{D} := \{
(z_1,z_2) \in \mathbb{A}^2:
2\Re(z_1) + |z_2|^2 = 0 \}.$$ 
We stress that the topological boundary of \(\Hyp{\mathbb{A}}^{2}\) is the one-point compactification of $\partial_{\rm Par}\mathbb{D}.$ 
%
%
Setting ${\xi} := z_2/\sqrt{2}$, the condition for the parabolic boundary becomes $\Re(z_1) = -|{\xi}|^2$.
Hence, we may write $\mathbf{z} = (-|\xi|^2 - v, \sqrt{2}\xi)$ for a $v\in \imaginary(\mathbb{A})$. 
This shows that the parabolic boundary can be identified with \(\mathbb{A} \times \imaginary(\mathbb{A})\).
We shall extend this identification within the Siegel domain, seeing the boundary as a limiting leaf of a foliation of the domain.

Let \(\lambda \in \mathbb{R}_{+}\) and consider each point \(z \in \Hyp{\mathbb{A}}^{2}\) for which the standard lift \(\mathbf{z}\) satisfies \(\langle \mathbf{z}, \mathbf{z} \rangle_{2} = -2\lambda^2\). 
Define \[
H_{\lambda} := \{ z \in \Hyp{\mathbb{A}}^{2}\,|\, \langle \mathbf{z}, \mathbf{z} \rangle_{2} = -2\lambda^2 \}.
\]
The collection $\mathcal{F} = \{H_{\lambda}\}_{\lambda \in \mathbb{R}_{+}}$ is a foliation of \(\Hyp{\mathbb{A}}^{2}\).
A point $z \in \Hyp{A}^{2}$ lies in $H_{\lambda}$ if and only if its identification with $\mathbf{z} \in \mathbb{D}$, satisfies $2 \Re(z_{1}) = - |z_{2}|^{2} - 2\lambda^2$.
We write $z_{2}$ as $\sqrt{ 2 }\xi$, with $\xi\in \mathbb{A}$, and we get that $z_{1} =
\Re (z_{1})+ \imaginary(z_{1}) = -|\xi|^{2} - \lambda^2 - \mathbf{v}
$ for $\mathbf{v}:=
\imaginary(-z_{1})\in \imaginary(\mathbb{A})$.
Thus $ {z} \in \mathbb{D}$ corresponds to a point $(\xi, \mathbf{v}, \lambda) \in \mathbb{A} \times \imaginary(\mathbb{A}) \times \mathbb{R}_{+}$, via
$
(z_1, z_2) = (
- |\xi|^{2} - \lambda^2 - \mathbf{v} ,
\sqrt{ 2 } \xi 
).
$
This construction identifies the Siegel domain with the set 
\(\mathbb{A} \times \imaginary(\mathbb{A}) \times \mathbb{R}_{+}\). 
\begin{theorem}[Horospherical coordinates]\label{thm_horospherical coordinate}
	The Siegel domain $\mathbb{D}$, \eqref{def_Siegel_domain}, can be parametrized via the following diffeomorphism: 
\begin{equation}\label{def_horospherical_coordinates}
\varphi: \mathbb{A} \times \imaginary(\mathbb{A}) \times \mathbb{R}_{+} \to \mathbb{D}, \quad (\xi, \mathbf{v}, \lambda) \stackrel{\varphi}{\longmapsto} 
\left( - |\xi|^{2} - \lambda^2 - \mathbf{v}, \, \sqrt{ 2 } \xi \right).\end{equation}
\end{theorem}
\begin{proof}
	It follows from the above construction of \(\varphi\), that the map \[
(z_1, z_2) \mapsto
\left( \frac{z_{2}}{\sqrt{ 2 }}, \,-\imaginary(z_{1}), \sqrt{ -\frac{\langle \mathbf{z}, \, \mathbf{z} \rangle_{2}}{2} }\right),
\] where \(\mathbf{z} = (z_{1},z_{2},1)\), is the inverse of
\(\varphi\). Both \(\varphi\) and \(\varphi^{-1}\) are smooth.
Thus \(\varphi\) identifies
\(\mathbb{A} \times \imaginary(\mathbb{A}) \times \mathbb{R}_{+}\) with
\(\mathbb{D}\) as manifolds.
\end{proof}
The map $\varphi$ from Theorem~\ref{thm_horospherical coordinate} is called {\em horospherical coordinate map}. Next, we pull back the Bergman metric given by \eqref{Bergman metric} via such a map. In other words, the next result expresses the hyperbolic distance in horospherical coordinates.
\begin{theorem}[Hyperbolic distance in horospherical coordinates]\label{theorem distance in horospherical coordinates}\index{hyperbolic! -- distance}
The pullback \(\varphi^{*}d\) of the Bergman metric \(d\) to the manifold
\(\mathbb{A}\times\imaginary(\mathbb{A})\times \mathbb{R}_{+}\) by the horospherical coordinate map \(\varphi\) satisfies
\begin{equation}\label{formula hyperbolic distance in horospherical coordinates} \varphi^{*}d((\xi, \mathbf{v}, \lambda),(\xi', \mathbf{v}', \lambda')) 
= 
 \arccosh \sqrt { \left(1+ \frac { 
 |\xi-\xi'|^{2} + |\lambda-\lambda'|^2}{2\lambda\lambda'} \right)^2+ \left( \frac { |{ \mathbf{v} - \mathbf{v'} +2 \imaginary(\bar{\xi}\xi')|} } {2\lambda\lambda'}\right)^2} ,
\end{equation}
for all $(\xi, \mathbf{v}, \lambda)$ and $(\xi', \mathbf{v}', \lambda')\in \mathbb{A}\times\imaginary(\mathbb{A})\times \mathbb{R}_{+}$.
\end{theorem}

\begin{proof}
	It follows from the definition of the Bergman metric \eqref{Bergman metric} that \begin{eqnarray*}
\left(\cosh\left( {\varphi^{*}d((\xi, \mathbf{v}, \lambda),(\xi', \mathbf{v}', \lambda'))} \right)\right)^2 &=&\left(\cosh\left( {d(\varphi(\xi, \mathbf{v}, \lambda), \varphi(\xi', \mathbf{v}', \lambda'))} \right)\right)^2 \\
&=& \frac{|\langle \varphi(\xi, \mathbf{v}, \lambda),
\varphi(\xi', \mathbf{v}', \lambda') \rangle_{2} |^{2}}{\langle \varphi(\xi, \mathbf{v}, \lambda), \varphi(\xi, \mathbf{v}, \lambda) \rangle_{2} 
\langle \varphi(\xi', \mathbf{v}', \lambda'), \varphi(\xi', \mathbf{v}', \lambda') \rangle_{2} }.
\end{eqnarray*}
We start by explicitly computing
\( 
\langle \varphi(\xi, \mathbf{v}, \lambda), \varphi(\xi', \mathbf{v}', \lambda') \rangle_{2} 
\): using the definition of $\varphi$ from Theorem~\ref{thm_horospherical coordinate}, the definition of the second Hermitian form \eqref{def second Hermitian form}, and simple properties of elements in $\mathbb{A}$ (see Exercise~\ref{exercise_norm_imaginary}), we get
\[
\begin{split}
 \langle \varphi(\xi, \mathbf{v}, \lambda), \varphi(\xi', \mathbf{v}', \lambda') \rangle_{2} &= \langle (-|\xi|^{2} - \lambda^2 - \mathbf{v} ,
\sqrt{ 2 }\xi ,
1 )
, (-|\xi'|^{2} - \lambda'^2 - \mathbf{v}' ,
\sqrt{ 2 }\xi',
1) \rangle_{2}\\
&= \begin{bmatrix}-|\xi|^{2} - \lambda^2 + \mathbf{v}, \sqrt{ 2 }\bar{\xi},1\end{bmatrix}\begin{bmatrix}
0 & 0 & 1 \\
0 & 1 & 0 \\
1 & 0 & 0 \\
\end{bmatrix}\begin{bmatrix}
-|\xi'|^{2} - \lambda'^2 - \mathbf{v}' \\
\sqrt{ 2 }\xi \\
1 \\
\end{bmatrix} \\
&= \begin{bmatrix}
1, \sqrt{ 2 }\bar{\xi},
-|\xi|^{2} - \lambda^2 + \mathbf{v} \end{bmatrix} \begin{bmatrix}
-|\xi'|^{2} - \lambda'^2 - \mathbf{v'} \\
\sqrt{ 2 }\xi' \\
1 \\
\end{bmatrix} \\
&= 
- |\xi'|^{2} - \lambda'^2 - \mathbf{v'} 
+ 2\bar{\xi}\xi' 
-|\xi|^{2} - \lambda^2 + \mathbf{v} \\
& = 
- \bar{\xi}\xi + \bar{\xi}\xi' + \bar{\xi}'\xi - \bar{\xi}'\xi' + \bar{\xi}\xi' - \bar{\xi}'\xi - (\lambda^2+\lambda'^2) + (\mathbf{v} - \mathbf{v'}) 
\\
& = 
- |\xi-\xi'|^{2} +2 \imaginary(\bar{\xi}\xi') - (\lambda^2+\lambda'^2) + (\mathbf{v} - \mathbf{v'}) .
\end{split}
\] 
From this last calculation, we first see that
\(\langle \varphi(\xi, \mathbf{v}, \lambda), \varphi(\xi, \mathbf{v}, \lambda) \rangle_{2} = -2\lambda^2\)
and
\(\langle \varphi(\xi', \mathbf{v}', \lambda'), \varphi(\xi', \mathbf{v}', \lambda')\rangle_{2} = -2\lambda'^2\)
as we expect for horospherical coordinates. Then we see that the squared norm of
\(\langle \varphi(\xi, \mathbf{v}, \lambda), \varphi(\xi', \mathbf{v}', \lambda') \rangle_{2}\)
has value
$$ \left(
 |\xi-\xi'|^{2} + \lambda^2+\lambda'^2\right)^2 + |{ \mathbf{v} - \mathbf{v'} +2 \imaginary(\bar{\xi}\xi')}|^2 
 =
 \left(2\lambda \lambda'+
 |\xi-\xi'|^{2} + |\lambda-\lambda'|^2\right)^2 + |{ \mathbf{v} - \mathbf{v'} +2 \imaginary(\bar{\xi}\xi')}|^2 
 .$$
The claim of the theorem follows after normalization. See Exercise~\ref{ex other formulations distance ROSS} for other formulations for the distance function.
\end{proof}
Compare the formula in Theorem~\ref{theorem distance in horospherical coordinates} with the conformal Poincar\'e half-space model of real-hyperbolic space $\Hyp{R}^2$; see Exercise~\ref{ex distance on half space}.

\section{Group structures and their relations}\label{sec: groups Heis on A}\index{Heisenberg! -- group over $\mathbb A$} 

\subsection{The (first) Heisenberg group over $\mathbb{A}$}
The manifold 
\(\mathbb{A} \times \imaginary(\mathbb{A})\) can be
given the following Lie group structure,
\begin{equation}\label{Prod_Heisenberg_A}(\xi, \mathbf{v}) \cdot (\xi', \mathbf{v}') := \left( \xi + \xi', \mathbf{v} + \mathbf{v}' + 2\imaginary(\bar{\xi}\xi') \right), \qquad \forall (\xi, \mathbf{v}) , (\xi', \mathbf{v}')\in \mathbb{A} \times \imaginary(\mathbb{A}).\end{equation}
This construction turns \(\mathbb{A} \times \imaginary(\mathbb{A})\) (the parabolic boundary of the Siegel domain) into a Lie group, known as the \emph{Heisenberg group over} \(\mathbb{A}\) which we denote by
\(\mathcal{N}\) or \( \mathcal{N}^\mathbb{A}_{1}\).
\begin{proposition}
Formula \eqref{Prod_Heisenberg_A} is a group product on \(\mathbb{A} \times \imaginary(\mathbb{A})\) .
\end{proposition}
\begin{proof}
	The set is closed under the operation and for every
\((\xi, \mathbf{u}) \in \mathbb{A}\times \imaginary(\mathbb{A})\), we have \[
(\xi, \mathbf{u})\cdot(0,0) = (0,0)\cdot(\xi, \mathbf{u}) = (\xi, \mathbf{u}).
\] Thus, \((0,0)\) is the identity element. Moreover, for every
\((\zeta, \mathbf{v})\in\mathbb{A}\times \imaginary(\mathbb{A})\), we have \[
(\zeta, \mathbf{v})\cdot(-\zeta,-\mathbf{v}) = (\zeta-\zeta, \mathbf{v}-\mathbf{v}+2\imaginary(-|\zeta|^{2})) = (0,0) = (-\zeta,-\mathbf{v})\cdot(\zeta, \mathbf{v}).
\] Hence, the inverse element \((\zeta, \mathbf{v})^{-1}\) is given by
\((-\zeta,-\mathbf{v})\). Finally, associativity holds as well, since \[
\begin{split}
 ((\xi, \mathbf{u})\cdot(\zeta, \mathbf{v}))\cdot(\eta, \mathbf{w}) &= (\xi+\zeta, \mathbf{u}+\mathbf{v}+2\imaginary(\bar{\xi}\zeta))\cdot(\eta, \mathbf{w}) \\
 &= (\xi+\zeta+\eta, \mathbf{u}+\mathbf{v}+2\imaginary(\bar{\xi}\zeta) + \mathbf{w} + 2\imaginary((\overline{\xi+\zeta} )\eta)) \\
 &= (\xi+\zeta+\eta, \mathbf{u}+\mathbf{v}+\mathbf{w} + 2\imaginary(\bar{\xi}\zeta + \bar{\xi}\eta + \bar{\zeta}\eta)) \\
 &= (\xi+\zeta+\eta, \mathbf{u}+\mathbf{v}+\mathbf{w} + 2\imaginary(\bar{\zeta}\eta) + 2\imaginary(\bar{\xi}(\zeta + \eta))) \\
 &= (\xi, \mathbf{u})\cdot(\zeta+\eta, \mathbf{v}+\mathbf{w}+ 2\imaginary(\bar{\zeta}\eta)) \\
 &= (\xi, \mathbf{u})\cdot((\zeta, \mathbf{v})\cdot(\eta, \mathbf{w})).
\end{split}
\]
\end{proof}

It is possible to define a Lie group structure on \(\mathbb{A}^n \times \imaginary(\mathbb{A})\) similarly to \eqref{Prod_Heisenberg_A} and define the $n$-th Heisenberg group over $\mathbb{A}$; see Exercise~\ref{n_Heisenberg_A}.

\subsection{A group structure on \(\mathbb{A} \times \imaginary(\mathbb{A}) \times \mathbb{R}_{+}\)}

On the Heisenberg group \(\mathcal{N} := (\mathbb{A}\times\imaginary(\mathbb{A}), \cdot)\) with group product \eqref{Prod_Heisenberg_A}, 
we particularly identify a group
homomorphism
\(\delta: \mathbb{R}_{+} \to \operatorname{Aut}(\mathcal{N})\).
This shall endow the set
\(\mathbb{A}\times\imaginary(\mathbb{A})\times \mathbb{R}_{+}\) with a group structure given by the semidirect product \(\mathcal{N} \rtimes_{\delta} \mathbb{R}_{+}\).

The {\em one-parameter subgroup of standard Heisenberg homotheties} is the homomorphism
\(\delta: \mathbb{R}_{+} \to \operatorname{Aut}(\mathcal{N})\),
 where 
\begin{equation}\label{def delta on A ImA}\delta_{\lambda}(\xi, \mathbf{v}) := (\lambda \xi, \lambda^2\mathbf{v}), \qquad \forall \lambda \in \R_+, \forall(\xi, \mathbf{v})\in\mathbb{A}\times\imaginary(\mathbb{A}).
\end{equation}
We call \(\delta_{\lambda}\) the \emph{Heisenberg homothety of ratio}
\(\lambda\).\index{Heisenberg! -- homotheties} 

\begin{proposition}
	The map
\(\delta: \mathbb{R}_{+} \to \operatorname{Aut}(\mathcal{N})\) takes, indeed, values in $\operatorname{Aut}(\mathcal{N})$
and is a group homomorphism.
\end{proposition}
\begin{proof}
	For every \(\lambda \in \mathbb{R}_{+}\), we have that
\(\delta_{\lambda} \in \operatorname{Aut}(\mathcal{N})\), because
\(\delta_{\lambda}\) is bijective and for all
\((\xi, \mathbf{u}),(\zeta, \mathbf{v}) \in \mathcal{N}\) we have \[\begin{split}
\delta_{\lambda}((\xi, \mathbf{u})\cdot(\zeta, \mathbf{v})) &= \delta_{\lambda}((\xi+\zeta, \mathbf{u}+\mathbf{v}+2\imaginary(\bar{\xi}\zeta))) \\
&= (\lambda(\xi +\zeta), \lambda^{2}(\mathbf{u}+\mathbf{v}+2\imaginary(\bar{\xi}\zeta))) \\
&= (\lambda \xi + \lambda \zeta, \lambda^{2}\mathbf{u} + \lambda^{2}\mathbf{v} + 2\imaginary(\lambda \bar{\xi}\lambda \zeta)) \\
&= \delta_{\lambda}(\xi, \mathbf{u})\cdot \delta_{\lambda}(\zeta, \mathbf{v}).
\end{split}
\] Furthermore, the map \(\delta\) itself is a group homomorphism, since \[
\delta_{\lambda \lambda'} = \delta_{\lambda}\circ \delta_{\lambda'}, \qquad \forall \lambda \lambda'\in \R,
\]
which is easy to check; see also Execise~\ref{n_Heisenberg_homotheties}.
\end{proof}

\subsection{Left invariance of the pullback of the Bergman metric}\index{semi-direct product}\index{Bergman metric}
Recall from Theorem~\ref{thm_horospherical coordinate} that the hyperbolic space, in the form of the 
Siegel domain $\mathbb{D}$, can be modelled by
$\mathbb{A} \times \imaginary(\mathbb{A}) \times \mathbb{R}_{+} $. Next, we put a group structure on this manifold.
Recall that on \(\mathcal{N}:=\mathbb{A} \times \imaginary(\mathbb{A})\) we put the group product \eqref{Prod_Heisenberg_A} on which the group $\R_+$ acts by \eqref{def delta on A ImA}. Thus, we perform the semidirect group product \(\mathcal{N} \rtimes_{\delta}\mathbb{R}_{+}\); see Exercise~\ref{ex prod on semidirect A}, on which we put coordinates $(\xi, \mathbf{v};\lambda)\in \mathbb{A} \times \imaginary(\mathbb{A}) \times \mathbb{R}_{+}$. 
On this group, we consider the distance \(\varphi^{*}d\) from Theorem~\ref{theorem distance in horospherical coordinates}.
\begin{theorem}\label{thm left-invariance hyp distance}
	The pullback \(\varphi^{*}d\) of the
Bergman metric \(d\) to the Lie group
\(\mathcal{N} \rtimes_{\delta}\mathbb{R}_{+}\) by the horospherical coordinate map \(\varphi\) is left-invariant.
\end{theorem}
\begin{proof}
	In Theorem~\ref{theorem distance in horospherical coordinates}, we have shown that the distance \(\varphi^{*}d\) 
	between two points $(\xi, \mathbf{v};\lambda)$ and $(\xi', \mathbf{v}';\lambda') $ in $ \mathbb{A}\times\imaginary(\mathbb{A})\times \mathbb{R}_{+}=\mathcal{N} \rtimes_{\delta}\mathbb{R}_{+}$ satisfies
\[\Big(\cosh\left( \varphi^{*}d((\xi', \mathbf{v}';\lambda'),(\xi'', \mathbf{v}'';\lambda''))\right)\Big)^2 = 
\frac { \left(
 |\xi'-\xi''|^{2} + \lambda'^2+\lambda''^2\right)^2 + |{ \mathbf{v}' - \mathbf{v''} +2 \imaginary(\bar{\xi'}\xi'')}|^2} {4\lambda'^2\lambda''^2 } .
\]
Let
\((\xi, \mathbf{v}, \lambda) \in \mathcal{N}\rtimes \mathbb{R}_{+}\). 
Using the basic definition of semidirect product, as in Exercise~\ref{ex prod on semidirect A}, 
we get that the left translation is 
\begin{eqnarray*}
&&\hspace{-1cm}\Big(\cosh\left( \varphi^{*}d((\xi, \mathbf{v};\lambda)\cdot(\xi', \mathbf{v}';\lambda),(\xi, \mathbf{v};\lambda)\cdot(\xi'', \mathbf{v}'';\lambda'')) \right)\Big)^2 \\
&=& \Big(\cosh\left( \varphi^{*}d(
(\xi+\lambda\xi', \mathbf{v}+\lambda^2\mathbf{v}' + 2\imaginary(\bar{\xi}\lambda\xi');\lambda\lambda'),
(\xi+\lambda\xi'', \mathbf{v}+\lambda^2\mathbf{v}'' + 2\imaginary(\bar{\xi}\lambda\xi'');\lambda\lambda'')
 \right)\Big)^2 \\
&=&\frac { \left(
 |\xi+\lambda\xi'-(\xi +\lambda\xi'')|^{2} + \lambda^2\lambda'^2+\lambda^2\lambda''^2\right)^2} {4\lambda^2\lambda'^2\lambda^2\lambda''^2 } 
 \\
&& \hspace{2.2cm}
 +\frac{ |{ \mathbf{v}+\lambda^2\mathbf{v}'+2 \imaginary( \bar{\xi}\lambda\xi') 
 - (\mathbf{v}+\lambda^2\mathbf{v}'' + 2\imaginary(\bar{\xi}\lambda\xi''))
 +2 \imaginary(\overline{(\xi+\lambda\xi')}(\xi+\lambda\xi''))}|^2} {4\lambda^2\lambda'^2\lambda^2\lambda''^2 } \\
 &=&\frac { \left(
 | \xi'- \xi'')|^{2} + \lambda'^2+ \lambda''^2\right)^2} {4 \lambda'^2 \lambda''^2 } 
 \\
&& \hspace{3.5cm}
 +\frac{ |{ \mathbf{v}'-\mathbf{v}''+\frac{2}{\lambda^2} \imaginary\Big( \bar{\xi}\lambda\xi'
 - \bar{\xi}\lambda\xi''
 + \bar\xi\xi + \bar\xi\lambda\xi'' +\lambda\bar\xi' \xi+\lambda\bar\xi' \lambda\xi'')
 \Big)}|^2} {4 \lambda'^2\lambda''^2 } \\
&=&\frac { \left(
 |\xi'-\xi''|^{2} + \lambda'^2+\lambda''^2\right)^2 + |{ \mathbf{v}' - \mathbf{v''} +2 \imaginary(\bar{\xi'}\xi'')}|^2} {4\lambda'^2\lambda''^2 } \\
&=& \Big(\cosh\left( \varphi^{*}d((\xi', \mathbf{v}';\lambda'),(\xi'', \mathbf{v}'';\lambda'')) \right)\Big)^2,
\end{eqnarray*}
where we have used that 
$ \bar{\xi}\lambda\xi'
 + \bar\xi\xi +\lambda\bar\xi' \xi= -|\xi|^2+2\Re (\bar{\xi}\lambda\xi')$ has trivial imaginary part.
\end{proof}

 \section{Exercises}\label{ex:ch_ROSS}
\begin{exercise} Both the real numbers and the complex numbers, together with their standard norm, form a Euclidean Hurwitz algebra in the sense of Definition~\ref{Eucidean Hurwitz algebra}, which is associative and commutative.
\end{exercise}

\begin{exercise}\label{ex:quaternions_Hurwitz} The quaternions together with the norm $|z| = \sqrt{\bar{z}z}$ form a Euclidean Hurwitz algebra in the sense of Definition~\ref{Eucidean Hurwitz algebra}, which is associative but not commutative.
\end{exercise}

\begin{exercise}\label{exercise_norm_imaginary}
For every element $\xi$ in a Euclidean Hurwitz algebra, we have 
$|\xi| = \bar{\xi}\xi$, $ 2 \imaginary({\xi}) = \xi-\bar{\xi}$,
and
$\bar \xi= -\xi + |\xi+1|^2 -|\xi|^2-1$,
where the latter defines conjugate in terms of the norm.
\end{exercise}

 \begin{exercise} The Hermitian transpose of a product of two matrices over a Euclidean Hurwitz algebra is the product of the Hermitian transposes in the reverse order, that is, $(AB)^*=B^*A^*$. 
 \end{exercise}
 
 \begin{exercise}
Every eigenvalue of every Hermitian transformation is a real number.
\end{exercise}

 \begin{exercise}
 Because of \eqref{eq_sesqui}, for every Hermitian form we have $\QQ{x,x}\in\mathbb{R}$.
 \end{exercise}

\begin{exercise}\label{def hyperbolic functions}
For the definitions 
\begin{equation}\label{def sinh cosh}\index{hyperbolic! -- cosine}\index{hyperbolic! -- sine}
\sinh x:={\frac {e^{x}-e^{-x}}{2}} \qquad \text{ and } 
\qquad \cosh x:={\frac {e^{x}+e^{-x}}{2}},
\end{equation} prove 
the general formula 
\begin{equation}\label{alternative cosh} \ \arccosh{(x)}=\ln\left(x+\sqrt{x^2-1}\right), \qquad \forall x\ge 1, \end{equation} 
and
 the general {\em half-argument formula} 
$\left(\sinh \left({\frac {x}{2}}\right) \right)^2= {\frac {\cosh (x)-1}{2}} $.
\end{exercise}

\begin{exercise}[The $\mathbb{A}$-hyperbolic $n$-space $\Hyp{A}^n$]\label{def A Hyp n}
Let $\mathbb{A}\in\{\mathbb{R}, \mathbb{C}, \mathbb{H}\}$ and $n\in \N$.
Let $\mathbb{P}(\mathbb{A}^{n+1})$ be the \emph{$n$-dimensional $\mathbb{A}$-projective space}, defined as the quotient of $\mathbb{A}^{n+1}\setminus\{0\}$ by the equivalence relation that identifies $x=(x_1, \dotsc,x_{n+1})$ with $x\lambda=(x_1\lambda, \dotsc,x_{n+1}\lambda)$ for each $\lambda\in\mathbb{A}\setminus\{0\}$. 
The class of $x\in \mathbb{A}^{n+1}\setminus\{0\}$ is denoted by $\mathbb{P}(x)$. 
Fixing a Hermitian form $\QQ{\cdot, \cdot}$ of signature $(n,1)$, we define the \emph{$\mathbb{A}$-hyperbolic $n$-space} as the set 
\begin{equation}\label{def_hyperbolic_n}
\Hyp{A}^n:=\{\mathbb{P}(x)\in\mathbb{P}(\mathbb{A}^{n+1}):\QQ{x,x}<0\}
\end{equation} equipped with the distance $d$ such that \begin{align}\label{DIST}
 d(\mathbb{P}(x), \mathbb{P}(y)):=\arccosh\sqrt{ \frac{\QQ{x,y}\QQ{y,x}}{\QQ{x,x}\QQ{y,y}}}.
\end{align} 
%
%
%
The set $\Hyp{A}^n$ is well defined, and that the distance formula does not depend on the representative chosen.
\end{exercise}

\begin{exercise}[Alternative formulas for the hyperbolic distance in horospherical coordinates]\label{ex other formulations distance ROSS} Formula \eqref{formula hyperbolic distance in horospherical coordinates} has the equivalent formulations:
$$ 
 \arccosh \frac { \sqrt { \left(
 |\xi-\xi'|^{2} + \lambda^2+\lambda'^2\right)^2 + {| \mathbf{v} - \mathbf{v'} +2 \imaginary(\bar{\xi}\xi')|}^2}} {2\lambda\lambda' } 
$$
 and
 $$ 
 \arccosh \sqrt { \left(1+ \frac { 
 |\xi-\xi'|^{2} + |\lambda-\lambda'|^2}{2\lambda\lambda'} \right)^2+ \left( \frac { |{ \mathbf{v} + \bar{\xi}\xi' - \mathbf{v'} -\bar{\xi}'\xi }| } {2\lambda\lambda'}\right)^2}, $$
 from which one easily observes the symmetry.
\end{exercise}

\begin{exercise}Formula \eqref{formula hyperbolic distance in horospherical coordinates} gives distance between the neutral element $(0,0,1)$ and an arbitrary element 
$(\xi, \mathbf{v}, \lambda)$
 the value
\begin{equation*}
 \arccosh \frac { \sqrt { \left(
 |\xi |^{2} + \lambda^2+1\right)^2 + |{ \mathbf{v} }|^2}} {2\lambda } .
\end{equation*}
\end{exercise}

\begin{exercise}\label{ex distance on half space}
Consider the conformal Poincar\'e half-space model of real-hyperbolic space $\Hyp{R}^2$. 
The distance satisfies:
\begin{eqnarray*} d((\xi, \lambda), (\xi', \lambda') ) 
&=& \arccosh\left( \frac { 
 |\xi-\xi'|^{2} + \lambda^2+\lambda'^2 } {2 {\lambda\lambda'}} \right)\\
 &=& 2 {\rm arcsinh} \left({\frac {\|(\xi, \lambda) - (\xi', \lambda')\|}{2{\sqrt {\lambda \lambda'}}}}\right)
, \qquad \forall (\xi, \lambda) , (\xi', \lambda') \in \mathbb{R}\times\mathbb{R}_{+} ,
\end{eqnarray*}
where 
 $ \|\cdot\| $ is the standard Euclidean norm on $\R^2$. Check other formulas in \cite[Theorem~7.2.1]{Beardon}.
\\
{\it Hint}. First use Theorem~\ref{theorem distance in horospherical coordinates} and the fact that the imaginary part in $\R$ is trivial.
Then, 
recalling Exercise~\ref{def hyperbolic functions},
use the general half-argument formula. 
\end{exercise}

\begin{exercise}[$n$-th $\mathbb{A}$-Heisenberg group]\label{n_Heisenberg_A}\index{Heisenberg! -- group over $\mathbb A$} 
Let $\mathbb{A}\in\{\mathbb{R}, \mathbb{C}, \mathbb{H}\}$ and $n\in \N$. The 
{\em $n$-th $\mathbb{A}$-Heisenberg group} $\mathcal{N}^\mathbb{A}_{n}$ is the set $\mathbb{A}^{n}\times\imaginary(\mathbb{A})$
endowed with the multiplication law 
\begin{equation}\label{Prod_Heisenberg_A_n}(\xi, \mathbf{v}) \cdot (\xi', \mathbf{v}') := \left( \xi + \xi', \mathbf{v} + \mathbf{v}' + 2\imaginary(\bar{\xi}\cdot\xi') \right), \qquad \forall (\xi, \mathbf{v}) , (\xi', \mathbf{v}')\in \mathbb{A}^n \times \imaginary(\mathbb{A}), \end{equation}
where $\bar{\xi}\cdot \xi' = \bar{\xi}^t\xi'$ is the scalar product between the vector $\bar{\xi}$ and the vector $\xi'$.
This indeed defines a Lie group, whose inverse satisfies $(\xi, s)^{-1}= (-\xi, -s)$, for all $(\xi, s) \in \mathcal{N}^\mathbb{A}_{n}$.
The group $\mathcal{N}^\mathbb{A}_{n}$ is nilpotent.
\end{exercise}
\begin{exercise}[Heisenberg homotheties]\label{n_Heisenberg_homotheties} 
Let $\mathbb{A}\in\{\mathbb{R}, \mathbb{C}, \mathbb{H}\}$ and $n\in \N$.
On the group $\mathcal{N}^\mathbb{A}_{n}$ from Exercise~\ref{n_Heisenberg_A}, for each $\lambda\in\mathbb{R}$
 define \emph{Heisenberg homothety of ratio $\lambda$} by 
 \begin{equation}\label{def delta on A ImA n}\delta_{\lambda}(\xi, \mathbf{v}) := (\lambda \xi, \lambda^2\mathbf{v}), \qquad \forall(\xi, \mathbf{v})\in\mathbb{A}^n\times\imaginary(\mathbb{A}).
\end{equation}
The Heisenberg homothety satisfies the following properties:
\begin{description}
\item[\ref{n_Heisenberg_homotheties}.i.] $\delta_a((\xi,s)(v,t))=\delta_a(\xi,s)\delta_a(v,t)$, for all $a\in\mathbb{R}$, $\xi,v\in\mathbb{A}^{n},$ and $ s,t\in\imaginary(\mathbb{A})$;
\item[\ref{n_Heisenberg_homotheties}.ii.] $\delta_a^{-1}=\delta_{a^{-1}}$, for all $a\in\mathbb{R}\setminus\{0\}$;
\item[\ref{n_Heisenberg_homotheties}.iii.] $\delta_{a a'} = \delta_{a}\circ \delta_{a'}$, for all $a,a\in\mathbb{R}$ .
\end{description} 
\end{exercise} 

\begin{exercise}
Let $\mathbb{A}\in\{\mathbb{R}, \mathbb{C}, \mathbb{H}\}$ and $n\in \N$.
On $\mathbb{A}^{n+1}$ consider the forms given by \begin{equation}\label{Kmatrix0}\QQ{x,y}:=\QQ{x,y}_K=x^*Ky, \qquad \forall x,y\in\mathbb{A}^{n+1}\end{equation}
where $K$ is either
\begin{align}\label{Kmatrix}
J_1:=\begin{bmatrix}
\mathbb{I}_{n} & 0\\
 0 &- 1
\end{bmatrix}
, \, \text{ or }
\begin{bmatrix}
0 & 0 & -1\\
0 & \mathbb{I}_{n-1} & 0\\
-1 & 0 & 0
\end{bmatrix}
, \, \text{ or }\, J_2:=
\begin{bmatrix}
0 & 0 & 1\\
0 & \mathbb{I}_{n-1} & 0\\
1 & 0 & 0
\end{bmatrix}
.\end{align}
These forms are Hermitian forms of signature $(n,1)$. 
\end{exercise} 

\begin{exercise}[The group $O_\mathbb{A}(n,1)$]\label{group OK}
Let $\mathbb{A}\in\{\mathbb{R}, \mathbb{C}, \mathbb{H}\}$ and $n\in \N$. Let 
$\GL(n+1, \mathbb{A})$ be the group of invertible $(n+1,n+1)$ matrices with coefficients in $\mathbb{A}$.
Let $O_\mathbb{A}(n,1)$ denote the subgroup of $\GL(n+1, \mathbb{A})$ that preserves the form $\QQ{\cdot, \cdot}$ given by {(\ref{Kmatrix0})}, that is $$O_\mathbb{A}(n,1):=\{A\in \GL(n+1, \mathbb{A})\,:\, \QQ{Ax,Ay}=\QQ{x,y}, \ \forall x,y\in\mathbb{A}^{n,1}\}.$$
\begin{description}
\item[\ref{group OK}.i.] 
We have $A\in O_\mathbb{A}(n,1)\Leftrightarrow A^*KA=K$, for the $K$ in \eqref{Kmatrix}. 
%
\item[\ref{group OK}.ii.] The set $O_\mathbb{A}(n,1)$ with the matrix multiplication is a Lie group.
\end{description}
\end{exercise}

\begin{exercise}
The induced action of $O_\mathbb{A}(n,1)$ on the projective space $\mathbb{P}(\mathbb{A}^{n+1})$ preserves the subset $\Hyp{A}^n$, as defined in \eqref{def_hyperbolic_n}, and acts by isometries on $\Hyp{A}^n$.
\end{exercise}

\begin{exercise}[Two distinguished subgroups: $A$ and $N$]\label{def_A_N}
We denote by $A$ and $N$ the following subsets of the group $O_\mathbb{A}(n,1)$:
\begin{itemize}
\item $A$ denotes the 1-parameter set, formed by the elements 
$$A(a):=\begin{pmatrix}
e^a & 0 & 0\\
0 &I_{n-1}& 0\\
0&0&e^{-a}
\end{pmatrix}, \quad a\in\mathbb{R};$$
\item $N$ denotes the set of matrices of the form
\begin{align}\label{Nu}
\nu(M,M_{13}):=\begin{pmatrix}
1&M&M_{13}\\
0&I_{n-1}&M^*\\
0&0&1
\end{pmatrix},
\end{align}
where $M$ is a $(1,n-1)$-matrix with elements in $\mathbb{A}$ and $M_{13}$ is in $\mathbb{A}$ and satisfies $$|M|^2=M_{13}+\conj{M}_{13}.$$
\end{itemize}
\begin{description}
\item[\ref{def_A_N}.i.] 
For all $t\in\mathbb{R}$ and all $\nu(M,M_{13})\in N$, we have
$$\nu(M,M_{13}) A(t)=A(t)\nu(e^{-t}M,e^{-2t}M_{13}),$$ and 
$$A(t)\nu(M,M_{13})=\nu(e^{t}M,e^{2t}M_{13})A(t).$$
%

\item[\ref{def_A_N}.ii.] The sets $A$ and $N$ are subgroups of $O_\mathbb{A}(n,1)$;
\item[\ref{def_A_N}.iii.] The set $NA$ is a subgroup of $O_\mathbb{A}(n,1)$;
\item[\ref{def_A_N}.iv.] The group $N$ is normal in $NA$;
\item[\ref{def_A_N}.v.] The group $A$ is isomorphic to $\mathbb{R}$.
\item[\ref{def_A_N}.vi.] The group $N$ is isomorphic to the Heisenberg group $\mathcal{N}^\mathbb{A}_{n}$, defined in Exercise~\ref{n_Heisenberg_A}.
\item[\ref{def_A_N}.vii.] The group $NA$ acts simply transitively on $\Hyp{A}^n$.
\end{description}
\end{exercise}

\begin{exercise}\label{ex prod on semidirect A}
If on \(\mathcal{N} :=\mathbb{A} \times \imaginary(\mathbb{A})\) we consider the group product \eqref{Prod_Heisenberg_A} and $\R_+$ acts on $\mathcal{N} $ by \eqref{def delta on A ImA}, then the semidirect group product on \(\mathcal{N} \rtimes_{\delta}\mathbb{R}_{+}\) is
 \[
(\xi, \mathbf{v};\lambda)\cdot(\xi', \mathbf{v}';\lambda') = (\xi+\lambda\xi', \mathbf{v}+\lambda^2\mathbf{v}' + 2\imaginary(\bar{\xi}\cdot\lambda\xi');\lambda\lambda'), \qquad \forall
 (\xi, \mathbf{v};\lambda),(\xi, \mathbf{v}';\lambda')\in \mathbb{A}\times \imaginary(\mathbb{A})\times\mathbb{R}_{+}.
\] 
{\it Hint}. Use the definition \eqref{prod:semi}.

\end{exercise}

\begin{exercise} The Lie algebra of $\Hyp{\C}^2$, seen as $\mathbb{C}^{n-1}\times i\R\times\mathbb{R}_+$ as in Exercise~\ref{ex prod on semidirect A}, is isomorphic to the one in Exercise~\ref{ex sec curvature complex hyperbolic plane}.\end{exercise}

\begin{exercise}\label{ex hyperbolic as Riemannian}
Every $\mathbb{A}$-hyperbolic $n$-space is a Riemannian manifold whose sectional curvature is at most $-1$.
\\{\it Hint.} For some evidence on the curvature, see Exercise~\ref{ex sec curvature complex hyperbolic plane}.
\end{exercise}

\begin{exercise}[Final exercise]\label{final_ex}
We have the following generalization of Theorem~\ref{thm_group_structure_ROSS2} to arbitrary dimensions.
For $\mathbb{A}\in\{\mathbb{R, \mathbb{C}, \mathbb{H}}\}$ and $n\in\mathbb{N} $, the $\mathbb{A}$-hyperbolic $n$-space $\Hyp{A}^n$ is isometric to the manifold $\mathbb{A}^{n-1}\times \imaginary((\mathbb{A})\times\mathbb{R}_+$ equipped with the group product given by 
 \[
(\xi, \mathbf{v};\lambda)\cdot(\xi', \mathbf{v}';\lambda') = (\xi+\lambda\xi', \mathbf{v}+\lambda^2\mathbf{v}' + 2\imaginary(\bar{\xi}\cdot\lambda\xi');\lambda\lambda'), \qquad \forall
 (\xi, \mathbf{v};\lambda),(\xi, \mathbf{v}';\lambda')\in \mathbb{A}^{n-1}\times \imaginary(\mathbb{A})\times\mathbb{R}_{+}.
\] 
and the left-invariant distance $d$ such that
$$ d(1,(\xi, \mathbf{v}, \lambda)) 
= \arccosh \frac { \sqrt { \left(
 \norm{\xi}^{2} + \lambda^2+1\right)^2 + |{ \mathbf{v} }|^2}} {2\lambda } , \qquad \forall (\xi, \mathbf{v}, \lambda)\in \mathbb{A}^{n-1}\times\imaginary(\mathbb{A})\times \mathbb{R}_{+},
$$
where $\norm{\cdot}$ denotes the $\ell_2$-norm on $\mathbb{A}^{n-1}$.
\end{exercise}

\chapter{Heintze groups and their visual boundaries}\label{ch_Heintze}
 
In this chapter, we show that to every Riemannian symmetric space of rank-one and noncompact type one can associate a `visual boundary' that has the structure of a Carnot group. Visual boundaries are associated with spaces with negative sectional curvature. In fact, every homogeneous negatively curved manifold has the structure of a semidirect product of the form $N\rtimes \R$ for some positively graded nilpotent group $N$, which canonically represents the visual boundary.
Hence, these boundaries are equipped with the structure of metric Lie groups that admit dilations.

\section{CAT(-1) spaces and visual boundaries}

The $\mathbb{A}$-hyperbolic $n$-space $\Hyp{A}^n$ 
has been presented in Section~\ref{sec A-hyperbolic space} for $n=2$ and in Exercise~\ref{def A Hyp n} for arbitrary $n\in \N$.
This metric space can be seen as a Riemannian space whose sectional curvature is at most $-1$; see Exercise~\ref{ex hyperbolic as Riemannian}. From the metric viewpoint, we say that $\Hyp{A}^n$ is a $\mathrm{CAT}(-1)$ metric space in the sense of Definition~\ref{def_CAT}. An explicit proof of this last statement can be found in \cite[Part II, Chapter 10]{Bridson_Haefliger}. We shall proceed with applying the theory of CAT($-1$) spaces. In particular, we shall recall the notion of visual boundary.

There are several equivalent definitions of CAT($-1$).
We recall just the most standard one and refer to \cite{Bridson_Haefliger} for more.
 
\begin{figure}[h]
 \centering
 \begin{tikzpicture}[scale=2]
 \draw (0,0) .. controls (0.4,0.2) and (0.7,0.2) .. (1.5,0) 
 .. controls (0.6,0.3) and (0.6,0.5) .. (0.5,0.9)
 .. controls (0.4,0.4) and (0.3,0.2) .. (0,0)
 ;

 \draw (2,0) .. controls (2.5,0.1) and (3,0.1) .. (3.5,0) 
 .. controls (3.2,0.2) and (2.8, 0.5).. (2.5,0.9)
 .. controls (2.4,0.5) and (2.2, 0.2).. (2,0) 
 ;
 
 
 \coordinate (x) at (0.7,0.145);
 \coordinate (y) at (0.89,0.26);
 \coordinate (barx) at (2.64,0.074);
 \coordinate (bary) at (3,0.375);

 \draw[densely dashed, blue] (x) -- (y);
 \draw[densely dashed, blue] (barx) -- (bary);
 \fill (x) circle (0.5pt) node[below left] {$x$};
 \fill (y) circle (0.5pt) node[above right] {$y$};
 \fill (barx) circle (0.5pt) node[below left] {$\bar x$};
 \fill (bary) circle (0.5pt) node[above right] {$\bar y$};

 \node at (-0.2,0.8) {$M$};
 \node at (0,-0.15) {$q$};
 \node at (1.5,-0.15) {$p$};
 \node at (0.5,1) {$r$};
 
 \node at (3.6,0.8) {\(\Hyp{\mathbb{\R}}^{2}\)};
 \node at (2,-0.15) {$\bar q$};
 \node at (3.5,-0.15) {$\bar p $};
 \node at (2.5,1) {$\bar r$};
 
 \end{tikzpicture}
 \caption{Comparison of triangles for the CAT($-1$) condition:
 The triangle on the left is a triangle in the CAT($-1$) metric space $M$, and the triangle on the right is in the hyperbolic plane \(\Hyp{\mathbb{\R}}^{2}\). All marked distances are pairwise equal except the dashed ones. The CAT($-1$) condition requires that $d(x,y)\leq d(\bar x, \bar y)$. One loosely says that the triangle on the left is {\em thiner} than its {\em comparison triangle}, i.e., the one on the right.\index{comparison triangle}}
 \label{fig_CAT-1}
\end{figure}
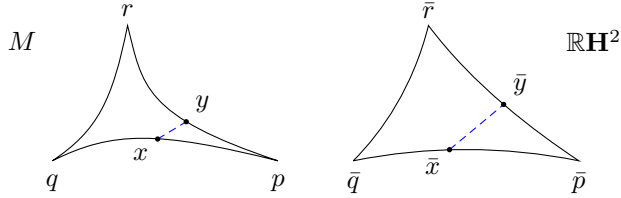

\begin{definition}[CAT($-1$) space]\label{def_CAT}\index{CAT$(-1)$}
A geodesic metric space $M$ is said to be {\rm CAT($-1$)}, or, {\em with generalized sectional curvature at most $-1$}, if the following comparison property holds (see Figure~\ref{fig_CAT-1}):
For all $p,q,r, x,y\in M$ and for all $\bar p, \bar q, \bar r, \bar x, \bar y\in \Hyp{\mathbb{\R}}^{2}$ such that
 \begin{equation*}
\begin{array}{ccccc}
d(p, x)+d(x,q) = d(p,q) ,& \quad& d(p,q)= d(\bar p, \bar q),& \quad & d(p, x)=d(\bar p, \bar x), \\
d(p,y)+d(y,r) = d(p,r), & &d(p,r)= d(\bar p, \bar r),& &d(x,q)=d(\bar x, \bar q), \\
&&d(r,q)= d(\bar r, \bar q),&&d(p,y)= d(\bar p, \bar y), \\
&&&& d(y,r)= d(\bar y, \bar r),
 \end{array}
 \end{equation*}
 we have $d(x,y)\leq d(\bar x, \bar y)$.
\end{definition}

\begin{definition}[Visual boundary]\index{visual! -- boundary}\index{geodesic! -- ray}
Let $M$ be a metric space.
A {\em geodesic ray} in $M$ is an isometric embedding $\gamma:[0, \infty)\to M$.
Given two geodesic rays $\gamma$ and $ \sigma$ in $M$, we write
$$\gamma\sim \sigma \stackrel{\rm def}{\Longleftrightarrow} \sup\{ d(\gamma(t), \sigma(t)) : t\in[0, \infty)\}\neq\infty, $$
and in this case, we say that 
$\gamma $ and $ \sigma $ are {\em asymptotic}. 
Fix $o\in M$. The {\em visual boundary of $M$ from} $o$ is defined as the set of geodesic rays starting at $o$, up to asymptotic equivalence:
$$\partial_{\rm vis} M :=\partial_{{\rm vis},o} M :=\{\gamma:[0, \infty)\to M\;:\; \gamma \text{ geodesic ray } , \gamma(0)=o \}/\sim.$$
\end{definition}

For rank-one symmetric spaces, it is convenient to consider visual boundaries where the visual point is at infinity. In this case, we fix a point in the visual boundary, and we consider geodesic lines originating from that point at infinity.
\begin{definition}[Parabolic visual boundary]\label{def_visual}\index{parabolic visual boundary}\index{visual! parabolic -- boundary}\index{geodesic! -- line}
Let $M$ be a metric space.
A {\em geodesic line} in $M$ is an isometric embedding $\gamma:\R\to M$.
Given two geodesic lines $\gamma$ and $ \sigma$ in $M$ (or two isometric embeddings $\gamma, \sigma:(-\infty,0]\to M$), 
we say that
\begin{equation}\label{eq_same_origin}
 \gamma\text{ and } \sigma \text{ have the {\em same origin} } \Longleftrightarrow 
\lim_{t\to - \infty} d(\gamma(t), \sigma(t))=0,
\end{equation}
while we say that they are {\em asymptotic} if the rays
$\gamma|_{[0, \infty)}$ and $ \sigma|_{[0, \infty)}$ are asymptotic.
Fix a geodesic line $\omega:\R\to M$, or just an isometric embedding $\omega:(-\infty,0]\to M$. 
The {\em parabolic visual boundary} of $M$ {\em from} $\omega$ is defined as: 
$$\partial_{\rm Par} M :=\partial_{{\rm Par}, \omega} M :=\{\gamma: \text{ geodesic line in } M \text{ with the same origin as }\omega \}/\sim.$$
We may refer to $\omega$ as the {\em viewpoint} of the visual boundary.\index{viewpoint}
\end{definition}

Visual boundaries can be naturally metrized, at least for CAT($-1$) spaces. Something similar holds for Gromov hyperbolic spaces, which we will not discuss here; see \cite[Proposition 3.3.3]{Buyalo_Schroeder_book}.
The construction of the parabolic boundary comes from the work of \cite[p.~383]{Hersonsky_Paulin}. 
Similar results on the metrization are also in the work of Hamenst\"adt in \cite{Hamenstadt89}, and for this reason, the visual distances, as we soon define, are also called {\em Hamenst\"adt distances}.\index{Hamenst\"adt! -- distance}

\begin{definition}[Visual distance]\label{def_visual_dist_gromov}\label{visual! -- distance}
Let $M$ be a metric space. Let 
 $\gamma$ and $\sigma$ 
be geodesic lines in $M$.
The \emph{Gromov product} 
between $\gamma$ and $\sigma$
 is defined as
 \begin{equation}\label{eq Gromov product}\index{Gromov product}
 (\,{\gamma}\,|\,{\sigma}\,)_\infty:=\lim_{t\to+\infty}\left(t- \frac{1}{2} d(\gamma(t), \sigma(t))\right), \end{equation}
when the limit exists. In this case, the \emph{visual distance} between $\gamma$ and $\sigma$
is defined as \begin{align}\label{Visdis}
\rho({\gamma},{\sigma}):= d_{\rm vis}({\gamma},{\sigma}):=e^{-({\gamma}|{\sigma})_\infty}
.
\end{align}
\end{definition}

The proof that when $M$ is a CAT($-1$) space, the limit in the Gromov product exists and that the visual distance as defined in (\ref{Visdis}) is a distance on $\partial_{{\rm Par}, \omega} M$ will be discussed in the next subsection. The main 
reference for the following presentation is \cite{Bourdon95}. 

\subsection{Hamenst\"adt distance}
We adopt Hamenst\"adt's viewpoint of horospherical distances via Busemann functions associated with geodesic lines. The aim is to prove that visual distances on parabolic boundaries of CAT($-1$) spaces satisfy the triangle inequality.

Let $M$ be a metric space. Fix $\omega: \R \to M$ a geodesic line in $M$. 
The {\em Busemann function} associated with $\omega$ is the function $b_\omega :M\to \R$ defined as the existing limit\index{Busemann function}
\begin{equation}\label{eq Busemann function} 
b_{\omega }(x):=\lim _{t\to -\infty }{\big (}d{\big (}\omega (t), x{\big )}+t{\big )}, \qquad \forall x\in M;
\end{equation}
see Exercise~\ref{ex_Busemann}.
In \eqref{eq Busemann function}, the chosen signs and the fact that the limit is to $-\infty$ may not coincide with the choices of other authors. The viewpoint that we take here is that the origin of $\omega $ is where $b_\omega$ equals $-\infty$, and it is where all the geodesic lines come from when we see the visual boundary from $\omega$. 
In fact, if $\gamma$ is a geodesic line in $M$ with the same origin as $\omega$, then 
\begin{eqnarray}
 b_\omega (\gamma(t)) &=& 
\lim _{s\to -\infty }{\big (}d{\big (}\omega (s), \gamma(t){\big )}+s{\big )} \nonumber \\
&\simeq& 
\lim _{s\to -\infty }{\big (}
d(\gamma (s), \gamma(t))
\pm
d(\omega (s), \gamma(s))
+s
{\big )} \nonumber \\
&=&
\lim _{s\to -\infty }{\big (}
t-s 
\pm
0+s
{\big )}=t, \label{eq18may1109}
\end{eqnarray}
where we used the symbol $\simeq$ to mean that we have upper and lower bounds according to the sign chosen as a consequence of triangle inequality.

Set $(x,y)_b := \frac12 (b(x) + b(y) - d(x,y))$ to be the 
 {\em Gromov product} between $x$ and $y$ with respect to a Busemann function $b:=b_\omega$, which is defined for $x,y\in M$, but actually would extends to $\partial_{{\rm Par}, \omega} M$. We call 
$\rho_b (\cdot, \cdot):= \exp( - (\cdot, \cdot)_b)$ the {\em Hamenstädt function}\index{Hamenst\"adt! -- function} with respect to $b$.
Notice that we have generalized \eqref{eq Gromov product} since if $\gamma$ and $\sigma$ have the same origin as $\omega$, then 
$$(\,{\gamma(t)}\,|\,{\sigma(t)}\,)_b = \frac12 (b(\gamma(t)) + b(\sigma(t)) - d(\gamma(t), \sigma(t)))
\stackrel{\eqref{eq18may1109}} {=} t -\frac12 d(\gamma(t), \sigma(t)).
$$

\begin{proposition}[Bourdon]\label{Hamenstädt gives a distance}\index{Bourdon}
Let $M$ be a CAT(-1) space. Fix $\omega$ a geodesic line and $b$ the Busemann function associated with $\omega$.
Then, the Hamenstädt function $\rho_b $ is a distance function on 
$\partial_{{\rm Par}, \omega} M$.
In fact, if for every $t\in \R$ we consider the function 
$$\alpha_t : \,(x,y) \in M\times M\longmapsto \alpha_t (x,y):= \frac{2}{e^t} \sinh \left( \frac{d(x,y)}{2}\right),$$ then we have: 
\begin{description}
\item[\ref{Hamenstädt gives a distance}.i]
As $ t $ goes to $\infty $, the function $\alpha_t $ converges to $ \rho_b $, in the sense that 
$$\lim_{t\to \infty} \alpha_t (\gamma(t), \sigma(t))= \rho_b (\gamma, \sigma) , \qquad \forall \gamma, \sigma\in \partial_{{\rm Par}, \omega} M.$$
\item[\ref{Hamenstädt gives a distance}.ii] 
The function $ \alpha_t$ is a distance on the set $b^{-1}(t)$.
\end{description}
\end{proposition}

\begin{proof}
The following argument is an adaptation of \cite[Théorème~2.5.1]{Bourdon95}. 
Referring to Exercise~\ref{ex_10May1246}, we begin by stressing that if $x, x' \in M$
and $z,z'\in \C$ are such that $\imaginary(z)= \imaginary(z')= e^{-t}$ and $d(z,z')= d(x, x')$, then 
\begin{equation}\label{alpha diversamente}\alpha_t (x, x') = |z-z'| .\end{equation}

{\it Proof of \ref{Hamenstädt gives a distance}.i.}
Take geodesic lines $\gamma $ and $ \gamma' $ that have the same origin as $\omega$.
We spell out the definitions:
\begin{eqnarray*}
 \alpha_t (\gamma(t) , \gamma'(t))
 &\stackrel{\rm def}{=}& 
 2 e^{-t} \sinh \left( \frac{d(\gamma(t) , \gamma'(t))}
 {2}\right) \\
 &\stackrel{\rm def}{=}& 
 2 e^{-t} \frac{1}{2} \left( e^{d(\gamma(t) , \gamma'(t))/2}
- e^{-d(\gamma(t) , \gamma'(t))/2}
 \right), 
\end{eqnarray*}
where we used the definition of $\alpha_t$ and the one of $\sinh$.
Then, since the function $e^{-d}$ is bounded by 1 and $\lim_{t\to +\infty} e^{-t}=0$, the above function has the same limit as
 \begin{eqnarray*}
e^{-t} e^{d(\gamma(t) , \gamma'(t))/2}
 &=& 
 \exp( - (t - d(\gamma(t) , \gamma'(t))/2)) \\
 &\to& \rho_b (\gamma, \gamma').
\end{eqnarray*}

{\it Proof of \ref{Hamenstädt gives a distance}.ii.}
We shall take $(0,e^{-t}) \in \C $ as comparison point for $\gamma(t)$, as $t\in \R$. 
Fix $x, x', x''$ points of the form 
$x=\gamma(t) , x'=\gamma'(t) , x''=\gamma''(t) $ such that $\gamma, \gamma', \gamma''$ are geodesic lines with the same origin as $\omega$. 
Take $z, z ',z'' \in \C$ so that $\imaginary(z)= \imaginary(z')=\imaginary(z'')= e^{-t}$,
$d(z, z ')= d(x, x ')$, 
$d(z ', z'')= d(x ', x'')$, and 
$\Re(z) \leq \Re(z') \leq \Re(z'') $; see Figure~\ref{ideal_triagles} for a visual representation.
(N.B. the triangles $z,z', \infty$ and $ z',z'', \infty$ are comparizon (ideal) triangles for $x, x', \omega$ and $ x', x'', \omega$, respectively.)
\begin{figure}[h]
 \centering
\begin{tikzpicture}
 \draw[dashed] (0,0) -- (0,5);
 \draw (0,0.4) -- (0,4.5);
 \draw[dashed] (2,0) -- (2,5);
 \draw (2,0.4) -- (2,4.5);
 \draw[dashed] (8,0) -- (8,5);
 \draw (8,0.4) -- (8,4.5);
 \draw[dashed] (-2,1) -- (9,1);
 \draw[dotted] (-2,0) -- (9,0);

 
 \draw (1,0) ++(45:1.4) arc (45:135:1.4);
 \draw (5,0) ++(18:3.16) arc (18:162.2:3.16);
 \draw (4,0) ++(14:4.123) arc (14:166:4.123);
 
 \node at (0,5.2) {$\infty$};
 \node at (2,5.2) {$\infty$};
 \node at (8,5.2) {$\infty$};
 \node at (1,3.2) {$\eta$};
 \filldraw[black] (0,1) circle (1pt);
 \filldraw[black] (2,1) circle (1pt);
 \filldraw[black] (8,1) circle (1pt);
 \node at (-0.2,0.7) {$z$};
	\node at (2.28,0.74) {$z'$};
	\node at ( 8.28,0.74) {$z''$};
 \node at (-1.2,1.2) {$e^{-t}$};
 \filldraw[black] (2,3.605) circle (1pt) node[anchor=west]{$z^*$}; 
\end{tikzpicture}
 \caption{The construction of two adjacent ideal triangles for the proof of Proposition~\ref{Hamenstädt gives a distance}.ii. The full curves are geodesic. The dashed horizontal line marks the height. The dotted horizontal line represents the parabolic visual boundary, which is the visual boundary minus a point, denoted by $\infty$.}
 \label{ideal_triagles}
\end{figure}
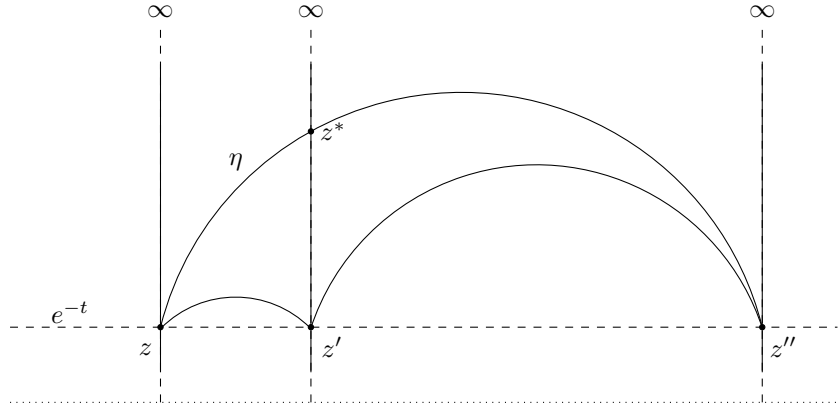

Let $\eta$ be the geodesic between $z$ and $z''$. Let $z^*$ be the point in $\eta$ with $Re(z^*)= Re(z')$. Let $x^*$ in $M$ be on the geodesic $\gamma'$ from $\omega$ to $x'$ with $d(x', x^*) = d(z',z^*)$.
Then, we bound using the CAT($-1$) comparison property, together with the fact that $z^*$ is along the geodesic $\eta$ and the triangle inequality:
\begin{eqnarray}
\nonumber
d(z,z'') &\stackrel{\eta \text{ geodesic}}{=}& d(z,z^*) + d(z^*,z'')\\
\nonumber
&\stackrel{\text{CAT}(-1)}{\geq}& d(x, x^*) + d(x^*, x'')\\
\label{23Apr24_1016}
&\stackrel{\text{triang. ineq.}}{\geq}& d(x, x'').
\end{eqnarray}
By monotonicity of $\sinh$ we deduce
\begin{eqnarray*}
\alpha_t(x, x'') &\stackrel{\rm def}{=} &\frac{2}{e^t} \sinh \left( \frac{d(x, x'')}{2}\right)\\
&\stackrel{\eqref{23Apr24_1016}}{\leq} & \frac{2}{e^t} \sinh \left( \frac{d(z,z'')}{2}\right) \\
&\stackrel{\eqref{alpha diversamente}}{=}& |z- z''|\\
&=& |z- z'|+ |z'- z''|\\
&\stackrel{\eqref{alpha diversamente}}{=}& 
\alpha_t(x, x') +\alpha_t(x', x'') ,
\end{eqnarray*}
where we used that by construction $\Re(z') $ is between $\Re(z) $ and $ \Re( z'')$.
Hence, the function $\alpha_t$ is a distance function.

 {\it End of proof of Proposition\ref{Hamenstädt gives a distance}} 
 Take $\gamma, \gamma',$ and $\gamma''\in \partial_{{\rm Par}, \omega} M$. Passing to the limit in the triangle inequality, given by \ref{Hamenstädt gives a distance}.ii,
 $$\alpha_t(\gamma (t ) , \gamma'' (t )) \leq \alpha_t(\gamma (t ), \gamma' (t )) +\alpha_t(\gamma' (t ), \gamma'' (t )), $$
 we obtain, by \ref{Hamenstädt gives a distance}.i, 
 $$\rho_b(\gamma , \gamma'' ) \leq \rho_b(\gamma , \gamma' ) +\rho_b(\gamma' , \gamma'' ). $$
Then, also the function $\rho_b $ satisfies the triangle inequality. 
\end{proof}


\section{The visual boundaries of $\mathbb{A}$-hyperbolic spaces}
\subsection{The boundaries of $\Hyp{A}^n$}\index{hyperbolic! -- space}
Let \(\mathbb{A}\) be either \(\mathbb{R}\), \(\mathbb{C}\), or the quaternions \(\mathbb{H}\). Let $n\in \N$.
Thanks to the way we defined the $\mathbb{A}$-hyperbolic $n$-space $\Hyp{A}^n$, in Section~\ref{sec A-hyperbolic space} and then in Exercise~\ref{def A Hyp n}, there is a natural way to represent its visual boundary. Denoting by $\mathbb{A}^{n,1}$ the space $\mathbb{A}^{n+1}$ equipped with a Hermitian form $\QQ{\cdot, \cdot}$ of signature $(n,1)$, 
 the set of points in the $\mathbb{A}$-hyperbolic $n$-space is
 $$\Hyp{A}^n=\{\mathbb{P}(x)\in\mathbb{P}(\mathbb{A}^{n,1})\,:\, \QQ{x, x}<0\},$$ where $\mathbb{P}(\mathbb{A}^{n,1})$ is the $\mathbb{A}$-projective space of $\mathbb{A}^{n,1}$. 
 
 As in Section~\ref{sec models hyp spaces}, we fix the Hermitian form \eqref{def second Hermitian form} so to have the Siegel domain $\mathbb{D}$, \eqref{def_Siegel_domain}, with the correspondence of Theorem~\ref{thm_horospherical coordinate}. The map $\varphi$ generalizes to arbitrary $n$. We summarize this correspondence with the following diagram:

\[
\begin{tikzcd}[column sep=1.5cm]
\mathbb{A}^{n-1} \times \imaginary(\mathbb{A}) \times \mathbb{R}_{+} 
\arrow[rr, bend right=20, "\hat\varphi"] 
\arrow[r, "{(\varphi,1)}"]
 &
 \mathbb{D}\times\{1\} 
\arrow[r, "\mathbb{P}"]
&
\Hyp{A}^n \subseteq \mathbb{P}(\mathbb{A}^{n,1}),
\end{tikzcd}
\]
defining the composition 
\begin{eqnarray}
\nonumber
\hat\varphi(\xi, \mathbf{v}, \lambda)&:=& \mathbb{P}(\varphi(\xi, \mathbf{v}, \lambda),1)
\\&\stackrel{\eqref{def_horospherical_coordinates}}{=}& \mathbb{P} 
\left( - |\xi|^{2} - \lambda^{2} - \mathbf{v}, \, \sqrt{ 2 } \xi ,1\right),
\qquad \forall (\xi, \mathbf{v}, \lambda) \in \mathbb{A}^{n-1} \times \imaginary(\mathbb{A}) \times \mathbb{R}_{+} . \label{definition hat phi}
\end{eqnarray}

In $ \mathbb{P}(\mathbb{A}^{n,1})$ we denote by $\omega$ the point $\omega:=\mathbb{P}(1, \mathbf{0},0)$, where here with $\mathbf{0}$ we indicate the zero in $\mathbb{A}^{n-1}$.
Obviously, the point $\omega$ belongs to the topological boundary of $\Hyp{A}^n$ within $\mathbb{P}(\mathbb{A}^{n,1})$.
We shall see that there is a correspondence between the topological boundary of $\Hyp{A}^n$ and the visual boundary of $\Hyp{A}^n$. Moreover, with this identification, the parabolic boundary of $\Hyp{A}^n$ from (a geodesic line representing) $\omega$ is in correspondence with $\mathbb{A}^{n-1} \times \imaginary(\mathbb{A}) \times \{0\}$.
With this purpose, we consider the following family of curves.

For every $\xi \in\mathbb{A}^{n-1}$ and $\mathbf{v}\in\imaginary(\mathbb{A})$, we consider the curve\index{$\gamma_{(\xi, \mathbf{v})}$} 
\begin{align}\label{visual geodesic} 
\gamma_{(\xi, \mathbf{v})}:\mathbb{R}&\longrightarrow\mathbb{A}^{n-1} \times \imaginary(\mathbb{A})\times\R_+
\\ t&\longmapsto\gamma_{(\xi, \mathbf{v})}(t):= 
(\xi, \mathbf{v}, e^{-t}).
 \nonumber\end{align} 
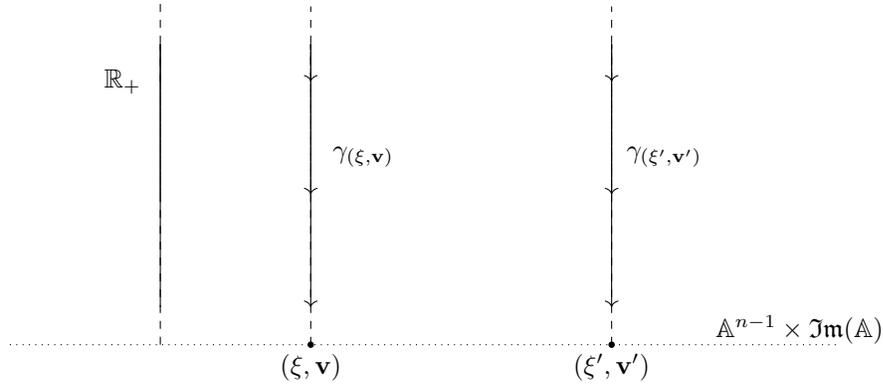
\begin{figure}[h]
 \centering
\begin{tikzpicture}
 \draw[dashed] (0,0) -- (0,4.55);
 \draw[dashed ] (0,2) -- (0,4);\draw[dashed ] (0,3.5) -- (0,4);
 \draw (0,0.5) -- (0,4);
 
 \draw[dashed] (2,0) -- (2,4.5);
 \draw[dashed,<-] (2,2) -- (2,4);\draw[dashed,<-] (2,3.5) -- (2,4);
 \draw[<-] (2,0.5) -- (2,4);

 \draw[dashed] (6,0) -- (6,4.5);
 \draw[dashed,<-] (6,2) -- (6,4);\draw[dashed,<-] (6,3.5) -- (6,4);
 \draw[<-] (6,0.5) -- (6,4);
 \draw[dotted] (-2,0) -- (9,0); 
 
 \node at (-.5,3.5) {$ \R_+$};
 \node at (2.7,2.5) {$\gamma_{(\xi, \mathbf{v})}$};
 \node at (6.7,2.5) {$\gamma_{(\xi', \mathbf{v}')}$};
 \node at (8.5,.2) {$\mathbb{A}^{n-1} \times \imaginary(\mathbb{A})$};
 \filldraw[black] (2,0) circle (1pt);
 \filldraw[black] (6,0) circle (1pt);
	\node at (2,-0.3) {$(\xi, \mathbf{v})$};
	\node at ( 6,-0.3) {$(\xi', \mathbf{v}')$};
\end{tikzpicture}
 \caption{The vertical geodesics in 
 the 
 $\mathbb{A}$-hyperbolic $n$-space $\Hyp{A}^n$ seen as $\mathbb{A}^{n-1} \times \imaginary(\mathbb{A})\times\R_+
$ in horospherical coordinates. Every point $(\xi, \mathbf{v})$ in $\mathbb{A}^{n-1} \times \imaginary(\mathbb{A})$ is associated with the geodesic line $\gamma_{(\xi, \mathbf{v})}:= 
(\xi, \mathbf{v}, e^{-t})$, whose height decreases when time increases.
 The dotted horizontal line is $\mathbb{A}^{n-1} \times \imaginary(\mathbb{A})$, which will represent the parabolic visual boundary.}
 \label{fig_geodesics_hyp}
\end{figure}
 \begin{lemma}\label{property of Geodesic in A hyperbolic}
The curves defined in \eqref{visual geodesic} satisfy the following properties:
\begin{description}
\item[\ref{property of Geodesic in A hyperbolic}.i.] each $\gamma_{(\xi, \mathbf{v})}$ is a geodesic;
\item[\ref{property of Geodesic in A hyperbolic}.ii.] every two $\gamma_{(\xi, \mathbf{v})}$ and $\gamma_{(\xi', \mathbf{v}')}$ have the same origin;
\item[\ref{property of Geodesic in A hyperbolic}.iii.] 
as $t\to -\infty$, the value $\hat \varphi(\gamma_{(\xi, \mathbf{v})}(t))$ tends to $\omega$
 within $\mathbb{P}(\mathbb{A}^{n,1})$;
\item[\ref{property of Geodesic in A hyperbolic}.iv.] 
as $t\to \infty$, the value $\gamma_{(\xi, \mathbf{v})}(t)$ tends to
$(\xi, \mathbf{v},0)$ in the topology of $\mathbb{A}^{n-1} \times \imaginary(\mathbb{A})\times[0,+\infty)$.
\end{description}
\end{lemma}

\begin{proof}
Fix $(\xi, \mathbf{v})\in \mathbb{A}^{n-1} \times \imaginary(\mathbb{A})$. 
To show that $\gamma_{(\xi, \mathbf{v})}$ is a geodesic, we directly use the formula \eqref{formula hyperbolic distance in horospherical coordinates} for the distance, for $s,t\in \R$:
\begin{eqnarray*}
\cosh d( \gamma_{(\xi, \mathbf{v})}(t), \gamma_{(\xi, \mathbf{v})}(s))
&\stackrel{\eqref{visual geodesic} } {=}&
\cosh d((\xi, \mathbf{v}, e^{-t}), (\xi, \mathbf{v}, e^{-s}))\\
&\stackrel{\eqref{formula hyperbolic distance in horospherical coordinates}}{=}
& \sqrt { \left(1+ \frac { 
 |e^{-t}-e^{-s}|^2}{2e^{-t-s}} \right)^2+ 0 } \\
&=&1+ \frac { 
 |e^{-t}-e^{-s}|^2}{2e^{-t-s}}\\
 &=&1 + \frac{e^{t+s}}{2} ( e^{-2t} - 2 e^{-t-s}+ e^{-2s})\\
 &=&1 + \frac{e^{ s- t }}{2} - 1 + \frac{e^{t-s}}{2} = \frac{e^{ s- t }+e^{t-s}}{2} \stackrel{\eqref{def sinh cosh}}{=} \cosh(s-t). 
\end{eqnarray*}
We deduce \ref{property of Geodesic in A hyperbolic}.i, thanks to the injectivity of $\cosh$ on the nonnegative real numbers.

Regarding \ref{property of Geodesic in A hyperbolic}.ii, 
according to the definition in \eqref{eq_same_origin} we need to prove that 
$$\lim_{t\to - \infty} d(\gamma_{(\xi, \mathbf{v})}(t), \gamma_{(\xi', \mathbf{v}')}(t))=0, \qquad \forall (\xi, \mathbf{v}) , (\xi', \mathbf{v}')\in \mathbb{A}^{n-1} \times \imaginary(\mathbb{A}).$$
We calculate: 
\begin{eqnarray}\nonumber
d(\gamma_{(\xi, \mathbf{v})}(t), \gamma_{(\xi', \mathbf{v}')}(t))&\stackrel{\eqref{visual geodesic} } {=} &
d((\xi, \mathbf{v}, e^{-t}), (\xi', \mathbf{v}', e^{-t}))\\
 &\stackrel{\eqref{formula hyperbolic distance in horospherical coordinates}}{=}& \arccosh \sqrt { \left(1+ \frac { 
 |\xi-\xi'|^{2} }{2e^{-2t}} \right)^2+ \left( \frac { |{ \mathbf{v} - \mathbf{v'} +2 \imaginary(\bar{\xi}\xi')|} } {2e^{-2t}}\right)^2}
 \label{calcolo_18may}
 \\
 &\stackrel{t\to -\infty}{\longrightarrow}& \arccosh (1) =0.\nonumber
 \end{eqnarray}

Regarding \ref{property of Geodesic in A hyperbolic}.iii, by the definition of $\hat \varphi$ in \eqref{definition hat phi}, we have, as $t\to -\infty$, that
\begin{eqnarray*}
\hat \varphi(\gamma_{(\xi, \mathbf{v})}(t))
& \stackrel{\eqref{visual geodesic}}{=}& \hat \varphi( \xi, \mathbf{v} , e^{-t})
\\&\stackrel{\eqref{definition hat phi}}{=}
& \mathbb{P} 
\left( - |\xi|^{2} - e^{-2t} - \mathbf{v}, \, \sqrt{ 2 } \xi ,1\right)
\\&=& \mathbb{P} 
\left( e^{2t}|\xi|^{2} + 1 + e^{2t}\mathbf{v}, \, -e^{2t} \sqrt{ 2 }\xi ,-e^{2t}\right) \longrightarrow 
\mathbb{P} (1, \mathbf{0},0) \stackrel{\rm def}{=} \omega.
\end{eqnarray*}
Regarding \ref{property of Geodesic in A hyperbolic}.iv, we clearly have
$\gamma_{(\xi, \mathbf{v})}(t)\stackrel{\rm def}{=}( \xi, \mathbf{v} , e^{-t})\to ( \xi, \mathbf{v} , 0)$, as $t\to \infty$.
\end{proof}

\begin{definition}\index{hyperbolic! boundary at infinity of -- space}
Let $\mathbb{A}\in\{\mathbb{R}, \mathbb{C}, \mathbb{H}\}$ and $n\in \N$.
We define the \emph{boundary at infinity of $\Hyp{A}^n$} as the set 
$$\Bhyp{A}^n:=\{\mathbb{P}(x)\in\mathbb{P}(\mathbb{A}^{n,1})\,:\, \QQ{x, x}=0\}.$$
\end{definition}

We arrived at the point where $\Bhyp{A}^n\setminus\{\omega\}$ can be identified with $\mathbb{A}^{n-1}\times\imaginary(\mathbb{A})\times\{0\}$, and hence with $\mathbb{A}^{n-1}\times\imaginary(\mathbb{A})$.
Recall that the latter may be endowed with the product from \eqref{Prod_Heisenberg_A_n}:
 $$(\xi, \mathbf{v};0) \cdot (\xi', \mathbf{v}';0) := \left( \xi + \xi', \mathbf{v} + \mathbf{v}' + 2\imaginary(\bar{\xi}\cdot\xi');0 \right), \qquad \forall (\xi, \mathbf{v}) , (\xi', \mathbf{v}')\in \mathbb{A}^{n-1} \times \imaginary(\mathbb{A}),$$
 which is actually isomorphic to $ \mathbb{A}^{n-1} \times \imaginary(\mathbb{A}) \times \{1\}$ inside
 $ \mathbb{A}^{n-1} \times \imaginary(\mathbb{A})\times \R_+=
 \mathcal{N}^\mathbb{A}_{n}\rtimes_{\delta}\mathbb{R}_{+}$, with its semidirect product; see Exercise~\ref{ex prod on semidirect A}.

 \subsection{The visual distance on the boundary of $\Hyp{A}^n$}

Our next aim is to explicitly write the visual distance for the $\mathbb{A}$-hyperbolic $n$-space $\Hyp{A}^n$ in the new coordinates.

\begin{definition}[Visual distance on the $n$-th $\mathbb{A}$-Heisenberg group]\index{Heisenberg! -- group over $\mathbb A$}\index{hyperbolic! visual distance for the -- space}
Let $\mathbb{A}\in\{\mathbb{R}, \mathbb{C}, \mathbb{H}\}$ and $n\in \N$.
For every pair $ (\xi, \mathbf{v})$ and $ (\xi', \mathbf{v}')\in \mathbb{A}^{n-1} \times \imaginary(\mathbb{A})$, we define 
the {\em visual distance} as
\begin{equation}\label{Heisenberg visual distance}
d\left((\xi, \mathbf{v}) , (\xi', \mathbf{v}') \right):= 
d_{\rm vis}\left( \gamma_{(\xi, \mathbf{v})}, \gamma_{(\xi', \mathbf{v}')}\right), 
\end{equation}
where $\gamma_{(\xi, \mathbf{v})}$ is defined in
\eqref{visual geodesic} and $d_{\rm vis}$ is defined in
\eqref{Visdis}.
\end{definition}
We finally prove that the parabolic visual boundaries of hyperbolic spaces are metric Lie groups:
\begin{theorem}
Let $\mathbb{A}\in\{\mathbb{R}, \mathbb{C}, \mathbb{H}\}$ and $n\in \N$.
On $\mathbb{A}^{n-1}\times\imaginary(\mathbb{A}) $, the visual distance reads as 
\begin{equation}\label{obtaining Koranyi}
d\left((\xi, \mathbf{v}) , (\xi', \mathbf{v}') \right)= 
\sqrt[4]{
\norm{ \xi - \xi'}^4 +\left| \mathbf{v} - \mathbf{v}' + 2\imaginary(\bar{\xi} \xi') \right|^2
}, \qquad \forall (\xi, \mathbf{v}) , (\xi', \mathbf{v}')\in \mathbb{A}^{n-1} \times \imaginary(\mathbb{A}),
\end{equation}
 it is left-invariant with respect to the group structure \eqref{Prod_Heisenberg_A_n}, and it is self-similar with respect to the dilations \eqref{def delta on A ImA n}.
\end{theorem}
\begin{proof}
We start by substituting the definitions \eqref{Heisenberg visual distance}, \eqref{Visdis}, and
\eqref{eq Gromov product}:
\begin{eqnarray}
\nonumber
d\left((\xi, \mathbf{v}) , (\xi', \mathbf{v}') \right)&\stackrel{\eqref{Heisenberg visual distance}}{=} &
d_{\rm vis}\left( \gamma_{(\xi, \mathbf{v})}, \gamma_{(\xi', \mathbf{v}')}\right) 
\\
\nonumber
&\stackrel{\eqref{Visdis}}{=}& e^{-
\left(
{\gamma_{(\xi, \mathbf{v})}}\Large|
{\gamma_{(\xi', \mathbf{v}')}}\right)_\infty
}\\
\label{141612022024}
&
\stackrel{\eqref{eq Gromov product}}{=}&\lim_{t\to\infty} \exp \left( \frac{1}{2} d(\gamma_{(\xi, \mathbf{v})}(t), \gamma_{(\xi', \mathbf{v}')}(t)) - t\right).
\end{eqnarray}
Therefore, we calculate the distance between the points along the curves, as defined in \eqref{visual geodesic}, as we did in \eqref{calcolo_18may} using \eqref{formula hyperbolic distance in horospherical coordinates}:
\begin{eqnarray*}
d(\gamma_{(\xi, \mathbf{v})}(t), \gamma_{(\xi', \mathbf{v}')}(t))
 &\stackrel{\eqref{calcolo_18may}}{=}&
 \arccosh \sqrt { \left(1+ \frac { 
 e^{2t}|\xi-\xi'|^{2} }{2} \right)^2+ \left( \frac { e^{2t}|{ \mathbf{v} - \mathbf{v'} +2 \imaginary(\bar{\xi}\xi')|} } {2}\right)^2} \\
 &=&
 \arccosh\left[\frac{ e^{2t}}{2}\sqrt { \left(2 e^{-2t}+ 
 |\xi-\xi'|^{2} \right)^2+ \left( |{ \mathbf{v} - \mathbf{v'} +2 \imaginary(\bar{\xi}\xi')|} \right)^2} \right].
\end{eqnarray*}
We denote by $\beta(t)$ the value 
\begin{equation}\label{auxiliary beta}
\beta(t):= \left(2 e^{-2t}+ 
 |\xi-\xi'|^{2} \right)^2+ \left( |{ \mathbf{v} - \mathbf{v'} +2 \imaginary(\bar{\xi}\xi')|} \right)^2,
 \end{equation}
so to continue, and using \eqref{alternative cosh}, 
\begin{eqnarray*}
d(\gamma_{(\xi, \mathbf{v})}(t), \gamma_{(\xi', \mathbf{v}')}(t))
&=&
 \arccosh\left(\frac{ e^{2t}}{2}\sqrt { \beta(t)} \right)
 \\ &\stackrel{\eqref{alternative cosh} 
}{=} &\ln\left(\frac{ e^{2t}}{2}\sqrt{\beta(t)}+\sqrt{ \frac{ e^{4t}}{4}\beta( t) -1}\right)\\
&=&\ln\left[e^{2t}\left( \frac{ \sqrt{\beta(t)}}{2}+ \sqrt{\frac{\beta(t)}{4}-e^{-4t}}\right)\right]\\
&=& 2t+\ln\left(\frac{ \sqrt{\beta(t)}}{2}+\sqrt{\frac{ \beta(t)}{4} - {e^{-4t}}}\right).
\end{eqnarray*}
Substituting the last equation in \eqref{141612022024}, we obtain
\begin{eqnarray*}\label{141712022024}
d\left((\xi, \mathbf{v}) , (\xi', \mathbf{v}') \right)&
=&
\lim_{t\to\infty} \exp \left( \frac{1}{2} 
 \ln\left(\frac{ \sqrt{\beta(t)}}{2}+\sqrt{\frac{ \beta(t)}{4} - {e^{-4t}}}\right)
 \right)\\
&
=&
\lim_{t\to\infty} \sqrt{
 \frac{ \sqrt{\beta(t)}}{2}+\sqrt{\frac{ \beta(t)}{4} - {e^{-4t}}}} 
\\
&
=& \sqrt[4]{ \beta(\infty)}\\
&=&\sqrt[4]{ \left( 
 |\xi-\xi'|^{2} \right)^2+ \left( |{ \mathbf{v} - \mathbf{v'} +2 \imaginary(\bar{\xi}\xi')|} \right)^2}
,
\end{eqnarray*}
where we took the limit of \eqref{auxiliary beta} as $t\to \infty$.
Therefore, we obtained \eqref{obtaining Koranyi}. The left-invariance and self-similarity are straightforward exercises.
\end{proof}

We refer to \eqref{obtaining Koranyi} as the {\em Cygan-Koranyi distance} on the $n$-th $\mathbb{A}$-Heisenberg group.\index{Cygan-Koranyi! -- distance}\index{Koranyi! -- distance}

\begin{proposition}\label{Visleft_triang}
The Cygan-Koranyi function on $\mathbb{A}^{n-1}\times\imaginary(\mathbb{A}) $, given by \eqref{obtaining Koranyi}, satisfies the triangle inequality; hence it is a distance function.
\end{proposition}
\begin{proof} 
The Cygan-Koranyi function is the visual parabolic distance function of the CAT($-1$) space $\Hyp{A}^n$. Therefore, by Bourdon's Proposition~\ref{Hamenstädt gives a distance}, it is a distance function. 

Still, we will give a more direct argument next. 
Because of left invariance, it is enough to prove the triangle inequality when the middle point is the identity element $(0, \mathbf{0})$. Namely, we need to prove that for every $(\xi, \mathbf{v}) , (\xi', \mathbf{v}')\in \mathbb{A}^{n-1}\times\imaginary(\mathbb{A}) $ we have
$$ 
d\left((\xi, \mathbf{v}) , (\xi', \mathbf{v}') \right) \leq 
d\left((\xi, \mathbf{v}) , (0, \mathbf{0}) \right)+
d\left((0, \mathbf{0}) , (\xi', \mathbf{v}') \right).$$
To do this, we rewrite the terms and bound using several times that the imaginary numbers are orthogonal to the real numbers:
\begin{eqnarray*}
d\left((\xi, \mathbf{v}) , (\xi', \mathbf{v}') \right)^2 
&=& \sqrt{
\norm{ \xi - \xi'}^4 +\left| \mathbf{v} - \mathbf{v}' + 2\imaginary(\bar{\xi} \xi') \right|^2}
\\&= & \left|{ \norm{ \xi - \xi'}^2 -\left( \mathbf{v} - \mathbf{v}' + 2\imaginary(\bar{\xi} \xi') \right) }\right|
\\&= &\left|{ \norm{ \xi }^2 - \bar\xi \xi' - \xi \bar\xi'+\norm{\xi'}^2 - \mathbf{v} + \mathbf{v}' - \bar{\xi} \xi'+ \xi\bar{ \xi}' }\right|
\\&= &\left|{ \norm{ \xi }^2 - \mathbf{v} - 2 \bar\xi \xi'+\norm{\xi'}^2 + \mathbf{v}' }\right|
\\&\leq &\left|{ \norm{ \xi }^2 - \mathbf{v} }\right|+\left|{ 2 \bar\xi \xi'}\right|+\left|{\norm{\xi'}^2 + \mathbf{v}' }\right|
\\&= &\sqrt{\norm{ \xi }^4 + \left|\mathbf{v} \right|^2 }+ 2 \norm{\xi}\norm{ \xi'}+\sqrt{\norm{ \xi' }^4 + \left|{\mathbf{v'} }\right|^2 }
\\&\leq &
\left(
d\left((\xi, \mathbf{v}) , (0, \mathbf{0}) \right)\right)^2+
2d\left((\xi, \mathbf{v}) , (0, \mathbf{0}) \right)
d\left((0, \mathbf{0}) , (\xi', \mathbf{v}') \right)
+\left(
d\left((0, \mathbf{0}) , (\xi', \mathbf{v}') \right)\right)^2
\\&= &
\left(
d\left((\xi, \mathbf{v}) , (0, \mathbf{0}) \right)+
d\left((0, \mathbf{0}) , (\xi', \mathbf{v}') \right)\right)^2,
\end{eqnarray*}
where in the second inequality, we used that the distance \eqref{visual distance from identity} is increasing in $|v|$.
\end{proof}


\section{Heintze groups}\index{semi-direct product}

A well-known motivation for the study of self-similar metric Lie groups in Sections~\ref{sec: self-similar} and~\ref{sec: self-similar2} is their appearance as parabolic visual boundaries of negatively curved homogeneous Riemannian manifolds.
More precisely, Heintze \cite{Heintze_74} showed that every simply connected negatively curved homogeneous Riemannian manifold is isometric to a Riemannian Lie group $ G$ that is a semidirect product \( N \rtimes_\alpha \R_+ \), where \( N \) is a simply connected nilpotent Lie group and at the Lie-algebra level, the element 1 in the Lie algebra \( \R \) of $\R_+$ acts on \( \mathfrak{n} \) by a derivation \( \alpha \) whose eigenvalues have strictly positive real parts.
The parabolic visual boundary of \( G \) may be identified with the Lie group \( N \) endowed with some \( \alpha \)-homogeneous distance in the sense of Definition~\ref{homogeneous metric group}, which is the visual distance.
It is important to remark that a quasi-isometric classification of Heintze groups, which at the moment is an unsolved problem, is equivalent to a quasi-symmetric classification of their parabolic boundaries,
which in turn reduces to a bi-Lipschitz classification of self-similar metric Lie groups (see \cite{Cornulier:qihlc, Kivioja_LeDonne_Nicolussi}, and references therein).

\subsection{Homogeneous groups and Heintze groups}
\begin{definition}\label{maar:homogeneous_group}\index{homogeneous! -- group}\index{group! homogeneous --} 
A \emph{homogeneous group} is a pair \( (N, \alpha) \) where \( N \) is a simply connected Lie group and \( \alpha \) is a derivation on the Lie algebra $\mathfrak n $ of \( N \) such that for each eigenvalue \( \lambda \) of \( \alpha \) it holds \( \Re(\lambda) > 0 \).
\end{definition}

If \( (N, \alpha) \) is a homogeneous group, then $N$ admits a contractive automorphism, hence \( \mathfrak{n} \) admits a positive grading; see Theorem~\ref{thm_Siebert}. In particular, the group $N$ is nilpotent; see Exercise~\ref{ex positive gradings give nilpotency}.

We already encountered homogeneous metric groups in Definition~\ref{homogeneous metric group}. However, in Chapter~\ref{ch_MetricGroups}, we were considering them already metrized with a self-similar distance. 
One of the punchlines of this chapter is to see such homogeneous groups as metric groups appearing as visual boundaries of negatively curved metric groups. 
We begin by introducing this latter type of groups.\index{metric! -- Lie group}\index{group! metric Lie --|see{metric Lie group}}\index{Lie group! metric --}

\begin{definition}\index{Heintze group}\index{group! Heintze --}  
A Lie group $G$ is called a {\em Heintze group} if it is of the form $N \rtimes_\alpha \R_+$, where $N$ is a simply connected Lie group, and $\alpha: \mathfrak n \to \mathfrak n $ is a derivation whose eigenvalues have positive real parts. In other words, the pair \( (N, \alpha) \) is a homogeneous group. 
Further, we say that a Heintze group $N \rtimes_\alpha \R_+$ (or, equivalently, the homogeneous group \( (N, \alpha) \)) is
	\begin{itemize}
\item \emph{purely real}, if the eigenvalues of \( \alpha \) are real numbers,\index{purely real}
\item \emph{of Carnot type} if it is purely real, if \( \alpha \) is diagonalizable over \( \mathbb{R} \), and if the eigenspace corresponding to the smallest of the eigenvalues Lie generates \( \mathfrak{n} \). Hence, the Lie group $N$ is a Carnot group, and some multiple of $\alpha$ is the derivation associated with the one-parameter subgroup of Carnot dilations from Definition~\ref{Dilations:groups}.\index{Carnot! -- group}\index{Carnot! -- type}\index{Carnot! -- dilation}
	\end{itemize}
\end{definition}


Notice that for every Heintze group $G:=N \rtimes_\alpha \R_+$, the nilpotent group $N$ represents the commutator subgroup $[G, G]$. Hence, every Heintze group is solvable; see Exercise~\ref{ex_def_solvable}.\index{commutator! -- subgroup}\index{solvable}

In this section, we discuss the following theorem from \cite{Heintze_74}.

\begin{theorem}[Heintze]\label{thm_Heintze_IFF}\index{isometrically homogeneous}\index{sectional curvature}\index{left-invariant! -- Riemannian metric}\index{Riemannian! --  metric}
 Every connected isometrically homogeneous Riemannian manifold of negative sectional curvature is isometric to a Heintze group metrized with a left-invariant Riemannian metric.
 
 Vice versa, every Heintze group may be metrized with a left-invariant Riemannian metric of negative sectional curvature.
\end{theorem}

The work of Heintze relies on earlier results by Wolf and by Kobayashi; see \cite{
WolfMR0163262, 
Kobayashi_MR0148015, Heintze_74}. 

\begin{theorem}[Heintze, after Wolf]\index{Wolf}
Every connected homogeneous Riemannian manifold of non-positive sectional curvature is isometric to a solvable Lie group metrized with a left-invariant Riemannian metric.
 \end{theorem}
 \begin{theorem}[Heintze]
Every connected solvable Lie group $G$ that admits a left-invariant Riemannian metric of negative sectional curvature is a Heintze group.
%
%
\end{theorem}
Obviously, the last two results give the first and harder part of Theorem~\ref{thm_Heintze_IFF}. We will instead discuss the second part of the statement of Theorem~\ref{thm_Heintze_IFF}.
\subsubsection{A discussion on metrics on some Heintze groups}

Not all left-invariant Riemannian metrics on Heintze groups are negatively curved.
For example, on the Heintze group of Carnot type given by the semidirect product between the Heisenberg group $\mathcal{N}^\mathbb{\C}_{1}$ and $\R_+$, there are left-invariant Riemannian metrics with some planes of positive sectional curvature. We provide an example of one such metric in Exercise~\ref{ex sec curvature positive plane}. In this example, the group is the same that represents the complex hyperbolic space $\Hyp{\C}^2$ of complex dimension 2.


\begin{example}[Heintze groups of Carnot type]\label{ex Heintze groups of Carnot type}
Let $N$ be a Carnot group. Let $G$ be the Lie group obtained by the semidirect product $\R_+\ltimes N$ where $\R_+$ acts on $N$ by Carnot dilations.
We now show that on $G$, there is a left-invariant Riemannian metric with negative sectional curvature.
\begin{proof}
We consider $\mathfrak n:=\Lie(N)$ as a positively graded vector space, denoting it by $V$, 
with grading $V=\bigoplus_{j=1}^s V_j$ and associated dilations $\delta_\eps$, for $\eps\in \R$.
Let $E_1, \dots, E_n$ be a basis of $V$ adapted to the direct-sum decomposition.
For every $j\in\{1, \dots,n\}$, let $w_j\in\{1, \dots,s\}$ be such that $E_j\in V_{w_j}$.
Let $\langle \cdot, \cdot \rangle$ be the scalar product that makes the basis $\{E_j\}_j$ orthonormal. 
Instead of changing the scalar product on $\mathfrak n$, we will change the Lie group structure on the vector space $V=\mathfrak n$, keeping the same scalar product but putting a different Lie bracket. We stress that the Lie brackets will be isomorphic to the one of $\mathfrak n$ for every value of a parameter $\eps$, except possibly for $\eps=0$.

Let $[\cdot, \cdot]_1$ be Lie brackets on $V$ seen as $\mathfrak n$.
Let $\{\alpha_{ij}^k\}_{ijk}$ be the structure constants of $\mathfrak n$, i.e.,
$
[E_i,E_j]_1 = \sum_{k=1}^n \alpha_{ij}^k E_k , $ for all $i,j\in \{1, \ldots, n\}$. 
For $\epsilon\in[0,1]$ define Lie brackets $[\cdot, \cdot]_\epsilon$ on $V$ by
\[
[X,Y]_\epsilon := \delta_\epsilon([X,Y]) ,
\]
i.e., the structure constants are
\[
[E_i,E_j]_\epsilon = \sum_{k=1}^n \epsilon^{w_k} \alpha_{ij}^k E_k, \qquad \forall i,j\in \{1, \ldots, n\} .
\]
Define $\mathfrak n_\epsilon$ as the vector space $V$ endowed with the Lie bracket $[\cdot, \cdot]_\epsilon$.
Whenever $\epsilon\neq0$, the Lie algebra $\mathfrak n_\epsilon$ is isomorphic to $\mathfrak n_0$ via the Carnot dilations $\delta_\eps$.
Notice that $[\cdot, \cdot]_\epsilon\to[\cdot, \cdot]_0 \equiv 0$ as $\epsilon\to0$. In other words, the Lie algebra $\mathfrak n_0$ is commutative.

Define the Lie algebra $\mathfrak g_\epsilon := \R\ltimes \mathfrak n_\epsilon$, 
with the Lie brackets $[\cdot, \cdot]_\epsilon$ extended from $V$ by
\[
[E_0,E_j]_\epsilon
:= w_j E_j , \qquad \forall j\in \{1, \ldots, s\},
\]
where $E_0:=(1, \bf{0})$, where $\bf{0}$ is the zero in $V=\mathfrak n_\epsilon$. In total, the structural constants of $\g_\eps$ with respect to the basis $E_0,E_1, \ldots, E_n$ are
\begin{equation}\label{struttura_eps}\left\{ \begin{array}{rcll}
\alpha_{ij}^k(\eps) & = & \epsilon^{w_k} \alpha_{ij}^k , & \forall i,j,k\neq 0 , \\
\alpha_{ij}^0(\eps) & = & 0, & \forall i,j \in \{0,1, \ldots,n\} , \\
\alpha_{0j}^k(\eps) & = & {w_j} \delta_{j}^k , & \forall j,k\neq 0 .
\end{array}\right.
\end{equation}
The map $(\id, \delta_\eps)$ is a Lie algebra isomomorphism between $\g_1$ and $\g_\eps$, for $\eps\neq 0$.

The sectional curvature, as in \eqref{def_sectional_curvature}, is
\begin{equation}\label{def_sectional_curvature2}
 {\rm Sec}_\epsilon(X,Y) := \frac{ \langle R_\epsilon(X,Y,Y), X \rangle }{\|X\|^2\|Y\|^2 - \langle X,Y \rangle^2 } ,
\end{equation}
where $R_\epsilon$ is the curvature tensor of $(\mathfrak g_\epsilon, \langle \cdot, \cdot \rangle)$. 
By Proposition~\ref{Prop_sectional_curvature_structural_constants}, each value $ \langle R_\epsilon(E_i,E_j,E_k) , E_h \rangle$ can be expressed in terms of the structure constants of $\mathfrak g_\epsilon$, which depend on $\epsilon$, with the formula:
\begin{equation}\label{eq66465fc5}
\begin{aligned}
 \langle R_\epsilon(E_i,E_j,E_k) , E_h \rangle
 = \sum_{\ell=0}^n \bigg[
 &\frac14 
 \left(\alpha_{jk}^\ell(\epsilon)-\alpha_{k\ell}^j(\epsilon)+\alpha_{\ell j}^k(\epsilon)\right)
 \left(\alpha_{i\ell}^h(\epsilon)-\alpha_{\ell h}^i(\epsilon)+\alpha_{hi}^\ell (\epsilon)\right) \\
 -&\frac14
 \left( \alpha_{ik}^\ell(\epsilon)- \alpha_{k\ell}^i(\epsilon)+\alpha_{\ell i}^k(\epsilon)\right)
 \left( \alpha_{j\ell}^h(\epsilon)- \alpha_{\ell h}^j(\epsilon)+\alpha_{hj}^\ell (\epsilon)\right) \\
 -&\frac12
 \alpha_{ij}^\ell (\epsilon)\left( \alpha_{\ell k}^h(\epsilon)-\alpha_{kh}^\ell(\epsilon)+ \alpha_{h\ell}^k(\epsilon)\right)
 \bigg].
\end{aligned}
\end{equation}

We focus for the moment on the limit case of $\eps=0$. In this case, the Lie algebra $\mathfrak g_0$ is the semidirect product of $\R$ with an abelian Lie algebra $\R^n$ that is graded by the degrees $w_1, \ldots, w_n$.
We claim that $(\mathfrak g_0, \langle \cdot, \cdot \rangle) $ has negative sectional curvature, that is
\begin{equation}\label{eq66465a42}
{\rm Sec}_0(X,Y)<0
\qquad\forall X,Y \in \R\times V \text{ linearly independent.}
\end{equation}
In fact, 
we need to show that 
\begin{equation}\label{eq66464feb}
 \langle R_0(X,Y,Y), X \rangle \le 0, \qquad \forall X,Y\in\mathfrak g,
\end{equation}
 with equality (if and) only if $X$ and $Y$ are linearly dependent.
Notice that, writing $X$ and $Y$ in components, by linearity we have
\[
\langle R_0(X,Y,Y), X \rangle \\
= \sum_{i,j,k,h=0}^n X^i Y^j Y^k X^h \langle R_0(E_i,E_j,E_k),E_h \rangle .
\]
So, we first distinguish three cases of indices in $\langle R_0(E_i,E_j,E_k),E_h \rangle$:
 the index 0 appears twice, it appears exactly once, or it does not appear.
Because of the symmetries of the curvature tensor, the three cases cover all possibilities.

{\bf Case 1}:
Suppose $i=h=0$ and $j,k\neq0$.
Then we claim that 
\[
\langle R_0(E_i,E_j,E_k),E_h \rangle
= - w_j^2 \delta_{jk} .
\]
Indeed, by the formula for the curvature tensor \eqref{eq66465fc5} in terms of the structural constants \eqref{struttura_eps}, we have:
\begin{eqnarray*}
 \langle R_0(E_i,E_j,E_k),E_h\rangle
&=&
\langle R_0(E_0,E_j,E_k),E_0 \rangle\\
&
 \stackrel{\eqref{eq66465fc5}}{=} 
 &\frac14 
 \left(\alpha_{jk}^0(0)-\alpha_{k0}^j(0)+\alpha_{0 j}^k(0)\right)
 \left(\alpha_{00}^0(0)-\alpha_{0 0}^0(0)+\alpha_{00}^0 (0)\right) \\
 &&-\frac14
 \left( \alpha_{0k}^0(0)- \alpha_{k0}^0(0)+\alpha_{0 0}^k(0)\right)
 \left( \alpha_{j0}^0(0)- \alpha_{0 0}^j(0)+\alpha_{0j}^0 (0)\right) \\
 &&-\frac12
 \alpha_{0j}^0 (0)\left( \alpha_{0 k}^0(0)-\alpha_{k0}^0(0)+ \alpha_{00}^k(0)\right)
 \\
 &&
 + \sum_{\ell=1}^n \bigg[
 \frac14 
 \left(\alpha_{jk}^\ell(0)-\alpha_{k\ell}^j(0)+\alpha_{\ell j}^k(0)\right)
 \left(\alpha_{0\ell}^0(0)-\alpha_{\ell 0}^0(0)+\alpha_{00}^\ell (0)\right) \\
 &&-\frac14
 \left( \alpha_{0k}^\ell(0)- \alpha_{k\ell}^0(0)+\alpha_{\ell 0}^k(0)\right)
 \left( \alpha_{j\ell}^0(0)- \alpha_{\ell 0}^j(0)+\alpha_{0j}^\ell (0)\right) \\
 &&-\frac12
 \alpha_{0j}^\ell (0)\left( \alpha_{\ell k}^0(0)-\alpha_{k0}^\ell(0)+ \alpha_{0\ell}^k(0)\right)
 \bigg] 
 \end{eqnarray*}
 \begin{eqnarray*} &
 \stackrel{\eqref{struttura_eps}}{=} 
 &
 \sum_{\ell=1}^n \bigg[ -\frac14
 \left( \alpha_{0k}^\ell(0) -\alpha_{0\ell }^k(0)\right)
 \left( \alpha_{0\ell }^j(0)+\alpha_{0j}^\ell (0)\right) \\
 &&\qquad -\frac12
 \alpha_{0j}^\ell (0)\left( \alpha_{0k}^\ell(0)+ \alpha_{0\ell}^k(0)\right)
 \bigg] \\ &
 \stackrel{\eqref{struttura_eps}}{=} 
 &
 - \sum_{\ell=1}^n 
 \frac12
 w_j \delta_{j}^\ell \left( w_k\delta_{k}^\ell+ w_k\delta_{\ell}^k\right)
 \\
 &=&- w_j^2 \delta_{j}^k.
\end{eqnarray*}

{\bf Case 2}:
Suppose $i=0$ and $j,k,h\neq0$.
Then, using \eqref{eq66465fc5}, similarly as in the previous case, one calculates:
\[
\langle R_0(E_i,E_j,E_k),E_h \rangle
= 0 .
\]

{\bf Case 3}:
Suppose $i,j,k,h\neq0$.
Then, using \eqref{eq66465fc5}, similarly as in the first case, one calculates: 
\[
\langle R_0(E_i,E_j,E_k),E_h \rangle
= w_iw_j (-\delta_{jk}\delta_{ih} + \delta_{ik}\delta_{jh}) .
\]

Using the symmetries of the curvature tensor, we see that in the general case, for $X, Y\in\mathfrak g$, we have
\begin{align*}
\langle R_0(X,Y,Y), X \rangle 
&= \sum_{i,j,k,h=0}^n X^i Y^j Y^k X^h \langle R_0(E_i,E_j,E_k),E_h \rangle \\
&= \sum_{j,h=1}^n X^0 Y^j Y^0 X^h \langle R_0(E_0,E_j,E_0),E_h \rangle \\
&\quad + \sum_{j,k=1}^n X^0 Y^j Y^k X^0 \langle R_0(E_0,E_j,E_k),E_0 \rangle \\
&\quad + \sum_{i,h=1}^n X^i Y^0 Y^0 X^h \langle R_0(E_i,E_0,E_0),E_h \rangle \\
&\quad + \sum_{i,k=1}^n X^i Y^0 Y^k X^0 \langle R_0(E_i,E_0,E_k),E_0 \rangle \\
&\quad + \sum_{i,j,k,h=1}^n X^i Y^j Y^k X^h \langle R_0(E_i,E_j,E_k),E_h \rangle 
\end{align*}
\begin{align*}
&= \sum_{j,h=1}^n X^0 Y^j Y^0 X^h (w_j^2\delta_{jh}) \\
&\quad + \sum_{j,k=1}^n X^0 Y^j Y^k X^0 (-w_j^2\delta_{jk}) \\
&\quad + \sum_{i,h=1}^n X^i Y^0 Y^0 X^h (-w_i^2\delta_{ih}) \\
&\quad + \sum_{i,k=1}^n X^i Y^0 Y^k X^0 (w_i^2\delta_{ik}) \\
&\quad + \sum_{i,j,k,h=1}^n X^i Y^j Y^k X^h w_iw_j (-\delta_{jk}\delta_{ih} + \delta_{ik}\delta_{jh}) 
\end{align*}
\begin{align*}
&= - \sum_{j=1}^n w_j^2 \left(X^0Y^j-Y^0X^j\right)^2 
 - \sum_{i,j=1}^n (X^i)^2 (Y^j)^2 w_iw_j 
 + \sum_{i,j=1}^n X^iY^jY^iX^j w_iw_j \\
&= - w_j^2 \left(X^0Y^j-Y^0X^j\right)^2 - \|A\|^2 \cdot \|B\|^2 + \langle A,B \rangle^2,
\end{align*}
where 
\[
A := \sum_{i=1}^n \sqrt{w_i}X^i E_i 
\quad\text{ and }\quad
B:= \sum_{j=1}^n \sqrt{w_j}Y^j E_j .
\]
Notice that $- \norm{A} \cdot \norm{B} + \langle A, B \rangle \le 0$, with equality if and only if $A$ and $B$ are linearly dependent.
We conclude that~\eqref{eq66464feb} holds for all $X,Y\in\mathfrak g_0$.
If equality holds, i.e., if $\langle R_0(X,Y,Y), X \rangle =0$, then $A$ and $B$ are linearly dependent, and thus so are $(X^1, \dots, X^n)$ and $(Y^1, \dots,Y^n)$, because $\{\sqrt{w_i} E_i\}_{i=1}^n$ is a basis of $V$.
Up to exchanging $X$ and $Y$, we can suppose $Y_j = r X_j$ for all $j$, with $r\in\R$.
Moreover, we also have $(X^0Y^j-Y^0X^j)^2=0$ for every $j$, that is, $(X^0 r-Y^0)^2(X^j)^2=0$.
Therefore, either $X$ and $Y$ are multiple of $E_0$, or $Y=rX$: in both cases, they are linearly dependent.
Thus, we obtained \eqref{eq66465a42}.

Next, we infer that there is $\epsilon_0>0$ such that 
\begin{equation}\label{eq664656ed}
{\rm Sec}_\epsilon(X,Y)<0
\quad\forall \epsilon\in[0, \epsilon_0]
\text{ and }\forall X,Y \in \R\times V \text{ linearly independent.}
\end{equation}
Indeed, 
thanks to~\eqref{eq66465a42} and the compactness of the Grassmannian of 2-planes in $\R\times V$,
we can just show that ${\rm Sec}_\epsilon(X,Y) \to {\rm Sec}_0(X,Y)$ as $\epsilon\to0$, locally uniformly in $X,Y$ linearly independent. 
But this is clear from~\eqref{def_sectional_curvature2} and~\eqref{eq66465fc5}.
\end{proof}
\end{example}

\subsection{Self-similar Lie groups as parabolic visual boundaries}\label{ssec:Heintze}\index{self-similar}\index{parabolic visual boundaries}\index{metric! -- Lie group}\index{Lie group! metric --}

We now explain how the parabolic boundary of a Heintze group is identified with a self-similar metric Lie group.

Let \( G:= N \rtimes \R_+ \) be a Heintze group equipped with a left-invariant Riemannian metric \( g \) such that the maximum of the sectional curvature is \( -1 \).
Up to replacing the $\R_+$ factor with a one-parameter subgroup orthogonal to $\Lie(N)$, we may assume that 
 \( N \times \{1\} \) is orthogonal to \( \{1_N\} \times \R_+\).
Or, more precisely, we have that \( T_{1_N}N \times \{0\}\) is perpendicular to \( \{0\} \times \R\) in $T_{1_G}G$.

We begin by parametrizing the boundary of $G$. In the following proposition, we say that the Lie group $\R_+$ acts {\em expansively}\index{expansively} on a Lie group $ N$ if every $\lambda\in(0,1)$ acts as a contractive automorphism in the sense defined at page \pageref{page_contractive}; see also Exercises~\ref{ex expansive action0} and \ref{ex expansive action}.\index{expansive action} 


\begin{proposition}\label{prop18may2205}
Let \( G:= N \rtimes \R_+ \) be a Heintze group, i.e., a semidirect product of a Lie group $N$ with the Lie group $\R_+$, such that $\R_+\acts N$ expansively.
Assume that $G$ is equipped with a left-invariant Riemannian metric such that \({\rm Sec}\leq -1 \) and such that 
 \( T_{1_N}N \times \{0\}\) is perpendicular to \( \{0\} \times \R\) in $T_{1_G}G$.
 Then there exists $\mu>0$ such that the map
\begin{equation}\label{def of gamma n}
n \in N \longmapsto \left( t\in \R\mapsto \gamma_n(t):= (n,e^{-\mu t})\right) \end{equation}
gives a continuous parametrization of the parabolic visual boundary $\partial_{{\rm Par}, \omega} G$ from $\omega$ where 
 \begin{equation}\label{def of omega}
 t \in \R \mapsto \omega(t):=(1_N, e^{-\mu t}).
\end{equation}
\end{proposition}

\begin{figure}[h]
 \centering
\begin{tikzpicture}
 \draw[dashed] (0,0) -- (0,4.55);
 \draw[dashed,<-] (0,2) -- (0,4);\draw[dashed,<-] (0,3.5) -- (0,4);
 \draw[<-] (0,0.5) -- (0,4);
 
 \draw[dashed] (2,0) -- (2,4.5);
 \draw[dashed,<-] (2,2) -- (2,4);\draw[dashed,<-] (2,3.5) -- (2,4);
 \draw[<-] (2,0.5) -- (2,4);

 \draw[dashed] (7,0) -- (7,4.5);
 \draw[dashed,<-] (7,2) -- (7,4);\draw[dashed,<-] (7,3.5) -- (7,4);
 \draw[<-] (7,0.5) -- (7,4);
 \draw[dotted] (-2,0) -- (9,0);

 \node at (-2,3.5) {$N \rtimes \R_+$};
 \node at (-0.5,2) {$\omega$};
 \node at (2.5,2) {$\gamma_n$};
 \node at (6.5,2) {$\gamma_m$};
 \node at (-2,.2) {$N$};
 \filldraw[black] (0,0) circle (1pt);
 \filldraw[black] (2,0) circle (1pt);
 \filldraw[black] (7,0) circle (1pt);
 \node at (-0,-0.3) {$1_N$};
	\node at (2,-0.3) {$n$};
	\node at ( 7,-0.3) {$m$};
\end{tikzpicture}
 \caption{The vertical geodesics in a Heitze group $N \rtimes \R_+$. Every point $n$ in $N$ is associated with the geodesic line $ \gamma_n(t):= (n,e^{-\mu t})$, whose hight decreases when time increases.
 The dotted horizontal line represents the parabolic visual boundary $N$, where the viewpoint is the geodesic line $ \omega(t):= (1_N,e^{-\mu t})$.
 }
 \label{fig_geodesics_Heintze}
\end{figure}
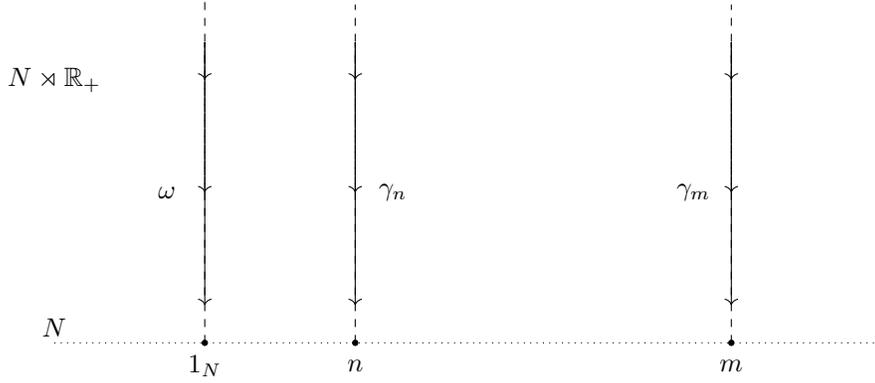
\begin{proof}
We consider the projection $N \rtimes \R_+ \to N \rtimes \R_+ /N\simeq \R_+$ modulo the normal subgroup $N$.
Since \( T_{1_N}N \times \{0\}\) is orthogonal to \( \{0\} \times \R\), then this projection is a submetry; see Corollary~\ref{prop:distance_quotient1}. 
Therefore, 
every injective curve with support in \(\{1_N\} \times \R_+ \) is length minimizing between every two of its points.
Therefore, the one-parameter subgroup \( s\in \R \mapsto (1_N,e^{ s}) \) is length minimizing, and, because the metric is left-invariant, it is parametrized at a constant speed, possibly different than 1.
We deduce that there exists \( \mu \in \R^+ \) such that $\omega$ as defined in \eqref{def of omega} is an isometric embedding.
By the left-invariance of the distance, we infer that every $\gamma_n$, as in \eqref{def of gamma n}, is a geodesic line. 

We denote by \( \lambda.n \) the action of $\lambda\in \R_+$ evaluated at \(n\in N \).
Using the left-invariance and the group laws of semidirect products, we express the distance between points on two such curves $\gamma_n$ and $\gamma_m$, for $n,m\in N$:
\begin{eqnarray}
\nonumber
 d\left( \gamma_n(t) , \gamma_m(t)\right )&\stackrel{\eqref{def of gamma n}}{=}&
 d\left( (n,e^{-\mu t} ),(m,e^{-\mu t} ) \right)\\
 \nonumber
&=&
 d\left ( (m,e^{-\mu t} )^{-1}\cdot(n,e^{-\mu t} ),1_G \right)
\nonumber\\
&\stackrel{{\rm Rem.\ref{rmk_semidir_grp}.v}}{=}&
 d \left( (e^{\mu t}.(m^{-1}),e^{\mu t} ) (n,e^{-\mu t} ),1_G \right)
\nonumber\\
&\stackrel{\eqref{prod:semi}}{=}&
 d \left( (e^{\mu t}.(m^{-1}) e^{\mu t} .(n), 1_\R) ,1_G\right )
 \nonumber\\
&=&
 d\left ( (e^{\mu t}.(m^{-1} n), 1_\R) ,1_G \right).\label{eq:calculation}
\end{eqnarray}

To explain that every $\gamma_n$ is an element in $\partial_{{\rm Par}, \omega} G$, we observe that it has the same origin as $\omega$: 
\begin{eqnarray*} 
d\left( \gamma_n(t), \omega(t) \right) &=&
d\left( \gamma_n(t), \gamma_{1_N}(t)\right )
\\
&\stackrel{\eqref{eq:calculation}}{=}&
 d\left ( (e^{\mu t}.( n), 1_\R) ,1_G \right)
 \\
&\stackrel{\rm Ex.\ref{ex expansive action0}}{\longrightarrow}&
 d\left ( (1_N, 1_\R) ,1_G\right )=
0 , \qquad \text{ as } t\to-\infty.
\end{eqnarray*}

We claim that every two curves \( \gamma_n \) and \( \gamma_{m} \) with $n, m\in N$ distinct, define two different points in $\partial_{{\rm Par}, \omega} G$, i.e., they are not asymptotic.
In fact, we show that 
\begin{equation}\label{diverging18may}
\lim_{t\to +\infty} d (\gamma_n(t), \gamma_m(t)) =\infty, \qquad \forall n, m\in N, n\neq m.
 \end{equation}
Indeed, since $m^{-1}n\neq 1_N$ and $\R_+\acts N$ is expansive, then, as $t\to +\infty$, then 
$e^{\mu t}.(m^{-1} n)$ leaves every compact set in $N$. Consequently, the curve
$(e^{\mu t}.(m^{-1} n), 1_\R)$ leaves every compact set in $G$, and so it diverges. Therefore, we infer 
\begin{equation*}
 d( \gamma_n(t) , \gamma_m(t) ) \stackrel{\eqref{eq:calculation}}{=}
 d \left( (e^{\mu t}.(m^{-1} n), 1_\R) ,1_G \right) \to \infty, \qquad \text{ as } t\to +\infty. 
 \end{equation*}
 Thus, we proved \eqref{diverging18may}.
  
We next show that the correspondence $n\in N\mapsto \gamma_n$ is surjective.
Consider \( \gamma \in \partial_{{\rm Par}, \omega} G \).
Take $n\in N$ and $t_0\in \R_+$ such that $\gamma(0)=
 (n, e^{t_0}) =\gamma_n(t_0)$.
 We shall show that $\gamma=\gamma_n$.
 Because both $\gamma$ and $\gamma_n$ have the same origin as $\omega$, then $\lim_{t\to -\infty} d (\gamma(t), \gamma_n(t) ) =0$.
Therefore, for very large $T>0$ the triangle $\gamma(0)$, $\gamma(-T)$, and $\gamma_n(-T) $ has side lengths that are $T$, $T+t_0$ and an infinitesimal length.
Because the metric space \( G \) is \( \mathrm{CAT}(-1) \) and \( T \) may be taken arbitrarily large, these triangles are degenerate, i.e., we have \( \gamma|_{(-\infty,0]} = \gamma_n|_{(-\infty,0]} \)
Since Riemannian geodesics do not branch, we conclude \( \gamma = \gamma_n \).

We check that the map in \eqref{def of gamma n} is continuous, where on $\partial_{{\rm Par}, \omega} G$ we consider the topology given by the visual distance from Definition~\ref{def_visual_dist_gromov}.
Namely, we check that if $m\to n$ in $N$, then $d_{\rm vis}(\gamma_m, \gamma_n)\to 0$. 
For this purpose, for $n,m\in N$ distinct, let $t^*\in \R$ such that 
\begin{equation*}
 d( \gamma_n(t^*) , \gamma_m(t^*) ) \stackrel{\eqref{eq:calculation}}{=}
 d \left( (e^{\mu t^*}.(m^{-1} n), 1_\R) ,1_G \right) = 1. 
 \end{equation*}
 We point out that, as $m\to n$, we have $t^*\to +\infty$.
By triangle inequality, we bound 
\begin{eqnarray*}
 d( \gamma_m(t) , \gamma_n(t) )
 &\leq &
 d( \gamma_m(t) , \gamma_m(t^*) )
 + d( \gamma_m(t^*) , \gamma_n(t^*) )
 + d( \gamma_n(t^*) , \gamma_n(t) )
 \\&=& 2 (t-t^*) + 1 
 \end{eqnarray*}
 and, consequently, 
 \begin{eqnarray*}
 d_{\rm vis}({\gamma_m},{\gamma_n})& \stackrel{\rm Def.\ref{def_visual_dist_gromov}}{=} &
 \exp\left(- \tfrac{1}{2} \lim_{t\to +\infty}(2t - d(\gamma_m(t), \gamma_n(t)))\right)
\\&\leq &
 \exp\left(-t^*+ \tfrac{1}{2} 
 \right)\to 0
 , \text{ as } m\to n.
 \end{eqnarray*}
 Thus, the parametrization is continuous -- it is actually a homeomorphism, but we shall motivate this fact later.
\end{proof}

Using the correspondence above, we view \( d_{\rm vis} \) as a distance function on \( N\).
\begin{definition}[Induced Hamenst\"adt distance function]\index{Hamenst\"adt! -- distance}
Let $N \rtimes \R_+$ be a Heintze group equipped with a left-invariant Riemannian metric such that \({\rm Sec}\leq -1 \). Up to changing the $\R_+$ factor, assume that 
 \( T_{1_N}N \times \{0\}\) is perpendicular to \( \{0\} \times \R\).
 The {\em induced Hamenst\"adt distance function} on $N$ is the function
 \begin{equation}\label{def_ham_18may}\hat d (n, m):= d_{\rm vis}(\gamma_n, \gamma_m), \qquad \forall n,m\in N,
 \end{equation}
 where $\gamma_n$ is the geodesic from \eqref{def of gamma n}.
\end{definition}

\begin{theorem}\index{admissible! -- distance function}
Given a Heintze group $N \rtimes \R_+$, the induced Hamenst\"adt distance function on $N$ is admissible, left-invariant, and self-similar. In fact, the group homomorphisms 
\begin{equation}\label{eq_dil_boundary_Heintze}
\delta_\lambda(n):= \lambda^\mu.n, \qquad\forall \lambda \in \R_+, \forall n\in N,
\end{equation}
are dilations for the distance, where $\mu$ is the normalizing factor from Proposition~\ref{prop18may2205}.
\end{theorem}
\begin{proof}
By \eqref{eq:calculation}, it is immediate that the visual distance is left-invariant on $N$.
For $s\in \R$, consider the action of $e^{\mu s}$ on $N$.
For $n,m\in N$, the visual distance changes as a translation in time for the geodesics as follows:
\begin{eqnarray*}
\hat d ( e^{\mu s}.n , e^{\mu s}.m ) &
\stackrel{\eqref{def_ham_18may}}{=}& d_{\rm vis}(\gamma_{e^{\mu s}.n}, \gamma_{e^{\mu s}.m})\\
&\stackrel{\rm Def.\ref{def_visual_dist_gromov}}{=}&\exp\left(-\lim_{t\to +\infty} \left(t - \tfrac{1}{2}d(\gamma_{e^{\mu s}.n}(t), \gamma_{e^{\mu s}.m}(t))\right)\right)
\\
&\stackrel{\eqref{eq:calculation}}{=}&
\exp\left(-\lim_{t\to +\infty} \left(t - \tfrac{1}{2}
 d ( (e^{\mu t}.((e^{\mu s}.m)^{-1} (e^{\mu s}.n)), 1_\R) ,1_G ) 
 \right)\right)
 \\
&=&\exp\left(s-\lim_{t\to +\infty} \left(t+s - \tfrac{1}{2}
 d ( (e^{\mu( t+s)}.( m^{-1} n), 1_\R) ,1_G ) 
 \right)\right)
 \\
&=&e^s\exp\left(-\lim_{t'\to +\infty} \left(t' - \tfrac{1}{2}
 d ( (e^{\mu t'}.( m^{-1} n), 1_\R) ,1_G ) 
 \right)\right)
 \\
 &\stackrel{\eqref{eq:calculation}}{=}&
 e^s \exp\left(-\lim_{t'\to +\infty} \left(t' - \tfrac{1}{2}d(\gamma_{ n}(t'), \gamma_{ m}(t'))\right)\right)
 \\
 &\stackrel{\rm Def.\ref{def_visual_dist_gromov}}{=}& e^s d_{\rm vis}(\gamma_{ n}, \gamma_{ m})
 \\
 &\stackrel{\eqref{def_ham_18may}}{=}&
 e^s \hat d ( n, m ).
 \end{eqnarray*}
 Setting $\lambda := e^s$ so \eqref{eq_dil_boundary_Heintze} gives $\delta_\lambda (n) = e^{\mu s}.n$, for $n\in N$, we deduce 
 $$\hat d(\delta_\lambda(n), \delta_\lambda(m)) = \lambda \hat d(n,m), \qquad \forall\lambda\in \R_+, \forall n,m\in N.$$

Regarding admissibility, it is a general property of left-invariant metrics on Lie groups admitting a one-parameter family $(\delta_\lambda)_\lambda$ of automorphisms each of which is a dilation of factor $\lambda$; see Exercise~\ref{ex_Enrico_Seba_ArkMat}.
However, we already proved in Proposition~\ref{prop18may2205} that the parametrization is continuous. Hence, the function $\hat d$ is continuous with respect to the topology. Consequently, using the dilations, it is easier to conclude that it induces the manifold topology, recall Exercise~\ref{property quasinorms}. 
\end{proof}


We conclude with a summary of this chapter.
Every Riemannian homogeneous negatively curved manifold has a notion of parabolic visual boundary, which has the structure of a metric Lie group admitting dilations.
All Carnot groups appear in this way, as boundaries of Heintze groups of Carnot type. In general, these Lie groups are precisely the nilpotent positively graded Lie groups.

\section{Exercise}
\begin{exercise}\label{ex_Busemann}
 Busemann functions associated with geodesic lines are well-defined (i.e., the limit exists) and are 1-Lipschitz. {\it Hint}. The limit is monotone. 
\end{exercise}

\begin{exercise}\label{ex_10May1246}
In the hyperbolic plane in the halfspace model, for every positive $t$ and $D$, we have
$d\left( (0,e^{-t}), ( \frac{2}{e^t} \sinh \left( \frac{D}{2}\right),e^{-t}) \right) = D $.
\\{\it Hint.} Use Exercise~\ref{ex distance on half space}.
\end{exercise}


\begin{exercise}\label{Visleft}
The visual distance defined as in {(\ref{Visdis})} is left-invariant, when $\Bhyp{A}^n$ is identified with $\mathbb{A}\times\imaginary(\mathbb{A}) $, with multiplication law \eqref{Prod_Heisenberg_A_n}.
\\{\it Hint.} Either use the left-invariance proved in Theorem~\ref{thm left-invariance hyp distance} or use formula \eqref{obtaining Koranyi}.
\end{exercise}

\begin{exercise}
The Cygan-Koranyi distance on $\mathbb{A}^{n-1}\times\imaginary(\mathbb{A}) $, given by \eqref{obtaining Koranyi}, between the point $(0, \mathbf 0) $ and a point $(\xi, \mathbf{v})\in \mathbb{A}^{n-1}\times\imaginary(\mathbb{A})$ is
\begin{equation}\label{visual distance from identity}
d ((0, \mathbf 0),(\xi, \mathbf{v}))=\sqrt[4]{ {\norm{\xi}^4} +|\mathbf{v}|^2}.
\end{equation}
This function is called {\em Cygan-Koranyi gauge}.\index{Koranyi! -- gauge}\index{Cygan-Koranyi! -- gauge}
\end{exercise}


\begin{exercise}
For $\mathbb{A}\in\{\mathbb{R, \mathbb{C}, \mathbb{H}}\}$ and $n\in\mathbb{N} $, equip $\mathbb{A}^{n-1}\times \imaginary((\mathbb{A})\times\mathbb{R}_+$ with the group product given in Exercise~\ref{final_ex}.
This Lie group is a Heintze group.
\end{exercise}

\begin{exercise}\label{ex:Gdelta3}\index{$G_z$}
Every group $G_z$, with $z\in \C\setminus i\R$, from Example~\ref{model_typeR} is a Heintze group.
\end{exercise}

\begin{exercise}\label{ex sec curvature positive plane}
Let $\g$ be the 4D Lie algebra with basis $X_1, \ldots, X_4$ and only nontrivial structural constants, as in \eqref{eq1655},
$$ c_{12}^3 =4, \quad c_{41}^1 = 1, \quad c_{42}^2 = 1, \quad c_{43}^3 = 2. $$
This Lie algebra is isomorphic to the one in Exercise~\ref{ex sec curvature complex hyperbolic plane}.
Consider a Lie group with such a $\g$ as Lie algebra, and on it, consider the left-invariant Riemannian structure for which $X_1, \ldots, X_4$ form an orthonormal basis.
We have the sectional curvature 
$	{\rm Sec}(X_1, X_3)= 2 $.
\\{\it Hint.} Apply \eqref{eq2130}.
\end{exercise}

\begin{exercise}
Let $N$ be a simply connected Lie group whose Lie algebra is positively graded. Let $G$ be the Lie group obtained by the semidirect product $\R_+\ltimes N$ where $\R_+$ acts on $N$ by the dilations induced by the grading.
On $G$, there is a left-invariant Riemannian metric with negative sectional curvature.
\\{\it Hint.} Adapt the proof of Example~\ref{ex Heintze groups of Carnot type}.
\end{exercise}

\begin{exercise}\label{ex expansive action0}\index{expansive action}
An action $\R_+\acts N$ by automorphisms is expansive if and only if 
$\lim_{t\to0}t.n = 1_N$, for all $n\in N$.
 \end{exercise}

\begin{exercise}\label{ex expansive action}
 If an action $\R_+\acts N$ is given by a derivation $A$, then the action is expansive if and only if the eigenvalues of $A$ have strictly positive real parts.
Equivalently, this is exactly when $N \rtimes_A \R_+ $ is a Heintze group.
\end{exercise}

\begin{exercise}\label{ex_Enrico_Seba_ArkMat}\index{admissible! -- distance function}
$\skull$ 
	Let $G$ be a Lie group.
	Let $\lambda\in (0, \infty)\mapsto\delta_\lambda\in\Aut(G)$ be a one-parameter group of Lie group automorphisms. 
	If there exists a left-invariant distance $d$ on $G$ with	$ d(\delta_\lambda x, \delta_\lambda y) = \lambda d(x,y)$, for all $x,y\in G $ and $\lambda>0 $,	then $d$ is admissible and $G$ is connected.
	\\{\it Hint.} Check
\cite[Theorem 1.1]{LeDonne_NicolussiGolo_metric}.
\end{exercise}

\newpage
\phantomsection
\addcontentsline{toc}{chapter}{Bibliography}
\bibliographystyle{alpha}
\bibliography{Book_chapters/general_bibliography,Book_chapters/bibliografia_applications}
\newpage
\phantomsection
\addcontentsline{toc}{chapter}{Index}
\printindex

\end{document}